\documentclass[12pt]{amsart}

\usepackage{graphicx}
\usepackage[margin=20truemm]{geometry}
\usepackage[hidelinks]{hyperref}
\usepackage[backend=biber,style=alphabetic]{biblatex}
\usepackage{amsmath, amssymb, amsthm}
\usepackage{mathrsfs}
\usepackage{mathtools}
\usepackage{empheq}
\usepackage{quiver}
\usepackage{tcolorbox}
\usepackage{cleveref}
\usepackage{enumitem}
\usepackage{here}

\tcbuselibrary{breakable, skins, theorems}

\newtcolorbox{simplebox}{
enhanced,
breakable=true,
colframe=black,
colback=white}

\newtcolorbox{titlebox}[1]{
enhanced,
breakable=true,
colframe=black,
colback=white,
top=4mm,
attach boxed title to top left={xshift=3mm,yshift=-3.5mm},
boxed title style={colframe=black,colback=white},
coltitle=black,
title=\textbf{#1}}

\theoremstyle{definition}
\newtheorem{df}{Definition}[section]
\crefname{df}{Definition}{Definitions}

\newtheorem{rem}[df]{Remark}
\crefname{rem}{Remark}{Remarks}

\newtheorem{nt}[df]{Notation}
\crefname{nt}{Notation}{Notations}

\newtheorem{as}[df]{Assumption}
\crefname{as}{Assumption}{Assumptions}

\theoremstyle{plain}
\newtheorem{lem}[df]{Lemma}
\crefname{lem}{Lemma}{Lemmas}

\newtheorem{prop}[df]{Proposition}
\crefname{prop}{Proposition}{Propositions}

\newtheorem{thm}[df]{Theorem}
\crefname{thm}{Theorem}{Theorems}

\newtheorem{cor}[df]{Corollary}
\crefname{cor}{Corollary}{Corollaries}

\mathtoolsset{showonlyrefs=true}

\newcommand{\id}{\mathrm{id}}
\newcommand{\op}{\mathrm{op}}

\newcommand{\can}{\mathrm{can}.}
\newcommand{\res}{\mathrm{res}.}
\newcommand{\inc}{\mathrm{inc}.}
\newcommand{\uhom}{\underline{\mathrm{Hom}}}
\newcommand{\colim}{\mathrm{colim}}
\newcommand{\mon}{\hookrightarrow}
\newcommand{\epi}{\twoheadrightarrow}

\newcommand{\of}{\circ}

\newcommand{\sub}{\subseteq}

\newcommand{\bus}{\supseteq}

\newcommand{\sm}{\setminus}
\newcommand{\set}[2]{\left\{~#1~\middle|~#2~\right\}}
\newcommand{\cat}[1]{\mathcal{#1}}
\newcommand{\ob}[1]{\left| \mathcal{#1} \right|}
\newcommand{\fk}[1]{\mathfrak{#1}}
\newcommand{\ub}[1]{\mathbf{#1}}

\newcommand{\N}{\mathbb{N}}
\newcommand{\Z}{\mathbb{Z}}
\newcommand{\Q}{\mathbb{Q}}
\newcommand{\R}{\mathbb{R}}

\newcommand{\ku}{\varnothing}

\newcommand{\To}{\Rightarrow}
\newcommand{\xto}[1]{\xrightarrow{#1}}
\newcommand{\xTo}[1]{\xRightarrow{#1}}

\newcommand{\ti}[1]{\textit{#1}}
\newcommand{\tup}[1]{\textup{#1}}
\newcommand{\tu}[1]{\underline{#1}}

\DeclareFieldFormat*{title}{\mkbibemph{#1}}
\DeclareFieldFormat*{journaltitle}{#1}
\DeclareFieldFormat*{booktitle}{#1}

\renewbibmacro{in:}{%
  \ifentrytype{article}
    {}
    {\printtext{\bibstring{in}\intitlepunct}}}

\DeclareFieldFormat[article]{volume}{\mkbibbold{#1}}
\DeclareFieldFormat[article]{number}{\bibstring{number}~#1}
\renewbibmacro*{journal+issuetitle}{%
  \usebibmacro{journal}%
  \setunit*{\addspace}%
  \iffieldundef{series}
    {}
    {\newunit
     \printfield{series}%
     \setunit{\addspace}}%
  \printfield{volume}%
  \setunit{\addspace}%
  \usebibmacro{issue+date}%
  \setunit{\addcolon\space}%
  \usebibmacro{issue}%
  \setunit{\addcomma\space}%
  \printfield{number}%
  \setunit{\addcomma\space}%
  \printfield{eid}%
  \newunit}

\DeclareFieldFormat[article]{pages}{#1}

\begin{document}

\title{Localization and Coalescence of Condensed Ringed Spaces}
\author{Naoto Fukutomi}
\date{\today}
\address{}
\email{}
\subjclass{}
\maketitle

\begin{abstract}
Gillam proved that the category of locally ringed spaces admits a fully faithful embedding into a certain category, which has a right adjoint that maps some simple objects to the spectra of rings. In this paper, we use condensed mathematics to define an analogous category and its coreflective full subcategories, and prove that one of their coreflections maps certain simple objects to the adic spectra of certain Huber pairs.
\end{abstract}

\tableofcontents

\section*{Introduction}

\subsection*{Condensed mathematics}

In algebra, certain algebraic objects carry topologies. For example, the field $\Q_p$ of $p$-adic numbers has the $p$-adic topology. If $I$ is an ideal of a commutative unital ring $R$, then one can define the $I$-adic topology on $R$. In the classical rigid analytic geometry, the affinoid algebras have topologies. However, the categories of such topological algebraic objects turn out to be not so desirable. We give two examples.
\begin{enumerate}
\item
The category of topological abelian groups is not an abelian category (\cite{Scholze:lectures}, 1. Lecture I). Consider the additive group $\R$ of the real numbers. With the usual topology $\R$ becomes a topological abelian group. On the other hand, $\R$ becomes topological abelian group if it is endowed with the discrete topology. The identity map
\begin{equation}
(\R, \text{discrete topology}) \to (\R, \text{usual topology})
\end{equation}
is both a monomorphism and an epimorphism in the category of topological abelian groups, but not an isomorphism.

\item
Colimits of topological algebraic objects are not so easy to handle. Consider the colimit $M = \colim \, M_i$ of some diagram in the category of topological abelian groups. The final topology $\cat{T}_0$ on $M$ induced by the canonical maps $M_i \to M$ is in general not compatible with the abelian group structure of $M$. Therefore the topology $\cat{T}$ of $M$ is in general coarser than $\cat{T}_0$, and it is not always the case that there is some explicit description of $\cat{T}$ . As a consequence, the notion of tensor products and localization is not always satisfactory.
\end{enumerate}

In order to overcome these difficulties, the idea of condensed mathmatics was introduced in \cite{Scholze:lectures}. This idea replaces the notion of topological abelian groups, topological rings, etc. with the notion of condensed abelian groups, condensed rings, etc. The category of condensed abelian groups does form an abelian category (\cref{prop:CAb is Grothendieck Abelian}). Moreover, these condensed objects tend to behave nicely with respect to colimits (\cref{prop:evaluation of CAb preserves limits filtered colimits and coproducts,prop:evaluation of CRing preserves limits and filtered colimits}, etc.). In particular, filtered colimits are easy to handle and this point is crucial in this paper (\cref{prop:compatibility of stalk and evaluation at S}, etc.).

One of the main features of topological abelian groups is the notion of completeness. The field $\Q_p$ of $p$-adic numbers and the affinoid algebras in the classical rigid analytic geometry are complete. If one has a topological abelian group which is not complete, then one often pass to its completion to obtain information. In condensed mathematics, the notion of comleteness is replaced by that of solidness, as defined in Definition 3.2.1 of \cite{Camargo:note} and in Definition 6.2.1 of \cite{Kedlaya:note}. In Definition 10.3.2 of \cite{Kedlaya:note}, the notion of $f$-coalescent modules is introduced, which can be seen as a generalization of solid condensed abelian groups. In this paper, we define the notion of $f/g$-coalescent modules (\cref{df:definition of coalescent modules}), which is a direct generalization of $f$-coalescent modules. As we shall see in \cref{cor:coalescent modules are closed under limits and colimits}, the category of $f/g$-coalescent modules turns out to be closed under colimits in the category of all condensed modules. This is a desirable property which is not shared by topological abelian groups: colimits of complete topological abelian groups need not be complete. 

Having said that, if the notion of condensed abelian groups/rings, etc. are totally different from that of topological abelian groups/rings, etc., then condensed mathematics is not so useful. In this paper, we will see that certain topological algebraic objects are described in terms of coalescence and localization (\cref{prop:coalescence of polynomial rings,prop:comparison of localizations of Huber rings}). In addition, although not written in this paper, it turns out that the class of first countable complete Hausdorff non-archimedean topological abelian groups is contained nicely in the category of condensed abelian groups.

\subsection*{Spectra and adic spectra}

Schemes are one of the main objects of algebraic geometry. They are obtained by gluing local models called affine schemes, which are constructed from commutative unital rings. For a commutative unital ring $R$, the corresponding affine scheme is called the spectrum of $R$. In the paper \cite{Gillam:paper}, this construction is proved to be a special case of a general construction called localization of ringed spaces. More precisely, Gillam proved the following. Let $\ub{PRS}$ be the category of ringed spaces $(X, \cat{O}_X)$ together with a set $M_x$ of prime ideals of the stalk $\cat{O}_{X,x}$ for each $x \in X$. Such system $M=(M_x)$ is called a prime system on $X$. Let $\ub{LRS}$ be the category of locally ringed spaces. Then we have a fully faithful functor $\cat{M} : \ub{LRS} \to \ub{PRS}$, $(X,\cat{O}_X) \mapsto (X,\cat{O}_X,\cat{M}_X)$, where $\cat{M}_{X,x}$ is the singleton consisting of the maximal ideal of the stalk $\cat{O}_{X,x}$ for $x \in X$.

\begin{thm} \tup{(\cite{Gillam:paper}, Theorem 2)}

The functor $\cat{M} : \ub{LRS} \to \ub{PRS}$ admits a right adjoint. For $X \in |\ub{PRS}|$, the image of $X$ under this right adjoint is called the localization of $X$, and is denoted by $X^{\mathrm{loc}}$.
\end{thm}

\begin{prop} \tup{(\cite{Gillam:paper}, Lemma 3)}

Let $A$ be a commutative unital ring. Let $(*,\cat{O})$ be the ringed space whose undelying topological space is the one point space $*$ and whose structure sheaf is determined by $\cat{O}(*) = A$. Let $\cat{T}$ be the prime system on $(*,\cat{O})$ such that $\cat{T}_*$ is equal to the set of all prime ideals of the stalk $\cat{O}_* = A$. Then the localization $(*,\cat{O},\cat{T})^{\mathrm{loc}}$ is equal to the spectrum of the commutative unital ring $A$. 
\end{prop}

On the other hand, the adic spaces, which are objects of non-archimedean analytic geometry, are obtained by gluing local models called adic spectra. These are constructed from (a certain class of) Huber pairs, which are topological algebraic objects. In this setting, one is interested in whether some results analogous to the results above hold. Let $\cat{C}$ be the category of triples $(X, \cat{O}_X, \cat{V}_X)$ consisting of a topological space $X$, a sheaf $\cat{O}_X$ of condensed rings on $X$ and a family $\cat{V}_X = (\cat{V}_{X,x})_{x \in X}$, where $\cat{V}_{X,x}$ is a set of continuous valuations on $\cat{O}_{X,x}$ for each $x \in X$. Let $\cat{C}_1$ be the full subcategory of $\cat{C}$ consisting of all objects $(X, \cat{O}_X, \cat{V}_X) \in |\cat{C}|$ such that the set $\cat{V}_{X,x}$ is a singleton for every $x \in X$. The unique element of $\cat{V}_{X,x}$ will be denoted by $|\cdot|_x$. 
Let $\cat{C}_l$ be the full subcategory of $\cat{C}_1$ consisting of all objects $(X, \cat{O}_X, \cat{V}_X) \in |\cat{C}_1|$ such that for every $x \in X$, the ring $\cat{O}_{X,x}(*)$ is a local ring and its maximal ideal is equal to the support of the valuation $|\cdot|_x$. Let $\cat{C}_f$ be the full subcategory of $\cat{C}_1$ consisting of all $(X, \cat{O}_X, \cat{V}_X) \in |\cat{C}_1|$ with the following property: For every open subset $U$ of $X$ and every $f, g \in \cat{O}_X(U)(*)$, the set
\begin{equation}
\set{x \in U}{|f_x|_x \leq |g_x|_x \neq 0}
\end{equation}
is an open subset of $X$. Let $\cat{C}_c$ be the full subcategory of $\cat{C}_f$ consisting of all $(X, \cat{O}_X, \cat{V}_X) \in |\cat{C}_f|$ such that for every $x \in X$ and every $f, g \in \cat{O}_{X,x}(*)$ with $|f|_x \leq |g|_x \neq 0$, the condensed $\cat{O}_{X,x}$-algebra $\cat{O}_{X,x}$ is $f/g$-coalescent. Then we prove the following proposition.

\begin{prop} \tup{(\cref{prop:intersections of C l C f C c are coreflective in C})}

Any intersection of these full subcategories $\cat{C}_1$, $\cat{C}_l$, $\cat{C}_f$, $\cat{C}_c$ is a coreflective full subcategory of $\cat{C}$.
\end{prop}

For any commutative unital topological ring $R$, let us write $\tu{R}$ for the condensed ring represented by $R$. Let $(A,A^+)$ be a Huber pair such that $A$ is complete Hausdorff. Let $X = \mathrm{Spa}(A,A^+)$ be the adic spectrum of $(A,A^+)$, with the structure presheaf $\cat{O}$. Let $\tu{\cat{O}}$ be the sheafification of the presheaf $\tu{\cat{O}}^{\mathrm{pre}}$ of condensed rings on $X$ defined by
\begin{equation}
\tu{\cat{O}}^{\mathrm{pre}}(U) := \tu{\cat{O}(U)} \quad (U \sub X \: \text{ open} ) .
\end{equation}
For $x \in X$, let $\tu{\cat{V}}_{x}$ be the singleton consisting of the unique continuous valuation on the stalk $\tu{\cat{O}}_x$ whose restriction to $A$ is equal to $x$. On the other hand, let $F$ be the sheaf of condensed rings on the one point topological space $\{*\}$ determined by $F(\{*\}) = \tu{A}$. Then $\cat{U} = \mathrm{Spa}(A,A^+)$ is a set of continuous valuations on the stalk $F_* = \tu{A}$ of $F$ at $* \in \{*\}$. Let $(\tilde{X}, \tilde{\cat{O}}, \tilde{\cat{V}})$ be the image of the object $(\{*\}, F, \cat{U})$ of $\cat{C}$ under the right adjoint of the inclusion $\cat{C}_l \cap \cat{C}_c \mon \cat{C}$.
\begin{thm} \tup{(\cref{thm:comparison of coreflection and adic spectrum})}

Suppose that one of the following holds.
\begin{enumerate}
\item $A$ is a strongly Noetherian Tate ring.
\item $A$ has a Noetherian ring of definition.
\end{enumerate}
Then we have a canonical isomorphism
\begin{equation}
(X, \tu{\cat{O}}, \tu{\cat{V}}) \simeq (\tilde{X}, \tilde{\cat{O}}, \tilde{\cat{V}})
\end{equation}
in the category $\cat{C}$.
\end{thm}

Here, we follow the idea of condensed mathematics and replace the sheaves of topologial rings by sheaves of condensed rings, in order to construct the category $\cat{C}$. Let us explain why this was successful. In the case of sheaves of topological rings, the stalk at a point is a ring, but in general does not have a very good topology. On the other hand, the stalks of sheaves of condensed rings are condensed rings. Moreover, these stalks have properties similar to the stalks of ordinary sheaves of rings (\cref{prop:testing sheaf isomorphisms by stalks,prop:testing equality of sheaf morphisms by stalks}, etc.) This is because stalks are defined as filtered colimits, and filtered colimits of condensed rings bahave especially well, as we have explained. Therefore we can rely on the notion of stalks, and this is one of the main reason of success in proving the theorem above.

\subsection*{Outline of this paper}

This paper is organized as follows. In \cref{sec:Review of condensed mathematics}, we review some basic definitions and facts about condensed mathematics. In addition, we explain some construction that is used in this paper. In \cref{sec:Valuations on condensed rings}, we study some properties of continuous valuations on condensed rings. The results in this section is used epecially in \cref{sec:Condensed ringed spaces}. In \cref{sec:Condensed sheaves}, we define sheaves of condensed rings and prove some basic results about them in detail. In general, the contents of this section are parallel to the basic theory of ordinary sheaves of rings. In \cref{sec:Condensed ringed spaces}, we define the category $\cat{C}$ and its full subcategories $\cat{C}_1 , \cat{C}_l , \cat{C}_f , \cat{C}_c$, and prove that any intersection of these full subcategories is coreflective in $\cat{C}$. In addition, we prove that these categories are complete and cocomplete. In \cref{sec:Relation to adic spectra}, we first review some basic definitions in topological algebra. After that, we study some relations between topological algebras and condensed mathematics. Using these results, we finally prove that the adic spectra of certain class of Huber pairs can be obtained as a coreflection along the inclusion functor $\cat{C}_l \cap \cat{C}_c \mon \cat{C}$.

\subsection*{Notations}

In this paper, we use the following notations.

\begin{enumerate}
\item Let $\cat{A}$ be a category.
\begin{enumerate}
\item
We write $\ob{A}$ for the class of all objects of $\cat{A}$.

\item
For each objects $A,B$ of $\cat{A}$, we write $\cat{A}(A,B)$ for the set of all morphisms $A \to B$ in $\cat{A}$.

\item
If $\cat{A}', \cat{A}''$ are full subcategories of $\cat{A}$, then we write $\cat{A}' \cap \cat{A}''$ for the full subcategory of $\cat{A}$ whose class of objects is equal to $|\cat{A}'| \cap |\cat{A}''|$.
\end{enumerate}

\item
\begin{enumerate}
\item $\ub{Set}$ denotes the category of sets and maps.

\item $\ub{Ab}$ denotes the category of abelian groups and homomorphisms.

\item $\ub{Ring}$ denotes the category of commutative unital rings and homomorphisms.

\item For any ring $R$, we write $\ub{Mod}_R$ for the category of $R$-modules and $R$-linear maps.

\item For any commutative unital ring $R$, we write $\ub{Alg}_R$ for the category of associative commutative unital $R$-algebras and homomorphisms.
\end{enumerate}

\item
\begin{enumerate}
\item $\ub{Top}$ denotes the category of topological spaces and continuous maps.

\item $\ub{TopAb}$ denotes the category of topological abelian groups and continuous homomorphisms.

\item $\ub{TopRing}$ denotes the category of commutative unital topological rings and continuous homomorphisms.

\item For any topological ring $R$, we write $\ub{TopMod}_R$ for the category of topological $R$-modules and continuous $R$-linear maps.

\item For any commutative unital topological ring $R$, we write $\ub{TopAlg}_R$ for the category of associative commutative unital topological $R$-algebras and continuous homomorphisms.
\end{enumerate}
\end{enumerate}

\subsection*{Acknowledgements}
The author is grateful to Atsushi Shiho for his support during the studies of the author.

\section{Review of condensed mathematics} \label{sec:Review of condensed mathematics}

In this section we review some basic definitions and facts concerning condensed mathematics. In \cite{Scholze:lectures}, the category of all profinite sets is used in the definition of condensed sets. However, this causes a set-theoretic size problem, since the category of all profinite sets is not essentially small. Here we follow \cite{Camargo:note} and \cite{Kedlaya:note} and we only consider light profinite sets.

\subsection{Light profinite sets}

\begin{df} (\cite{Camargo:note}, Proposition 2.1.1 and Definition 2.1.3)

A \ti{light profinite set} is a totally disconnected compact Hausdorff second countable topological space. We write $\ub{Prof}$ for the category of light profinite sets and continuous maps.
\end{df}

$\ub{Prof}$ is a full subcatgory of $\ub{Top}$.

\begin{prop} \label{prop:closure properties of light profinite sets} \,
\begin{enumerate}
\item Any closed subspace of a light profinite set is a light profinite set.

\item Any finite disjoint union of light profinite sets is a light profinite set.

\item $\ub{Prof}$ is closed under countable limits in $\ub{Top}$.
\end{enumerate}
\end{prop}

\begin{proof}
These assertions follow from the fact that the following classes of topological spaces are closed under taking closed subspaces, finite disjoint unions and countable limits.
\begin{itemize}
\item The class of totally disconnected spaces.
\item The class of compact Hausdorff spaces.
\item The class of second countable spaces.
\end{itemize}
\end{proof}

\begin{nt} \label{nt:one point topological space}
In this paper, we write $*$ for the topological space with a unique point. It is a light profinite set. The unique point of $*$ is also denoted by $*$.
\end{nt}

\begin{prop}
Let $X$ be a light profinite set. Then the cardinality $|X|$ of $X$ satisfies $|X| \leq |\R|$.
\end{prop}

\begin{proof}
This follows from the fact that any light profinite set admits a surjection from the Cantor set $\prod_{n \in \N} \{0,1\} \,$ (\cite{Camargo:note}, Proposition 2.1.9 ; \cite{Kedlaya:note}, Proposition 3.3.6).
\end{proof}

\begin{cor} \label{cor:Prof is essentially small}
The category $\ub{Prof}$ of light profinite sets and continuous maps is essentially small.
\end{cor}

\subsection{Condensed sets}

\begin{df} (\cite{Camargo:note}, Definition 2.2.1 ; \cite{Kedlaya:note}, Definition 4.1.1)
\begin{enumerate}
\item
Let $S$ be a light profinite set. A \ti{cover} of $S$ in $\ub{Prof}$ is a finite family $(S_i \to S)_i$ of morphisms in $\ub{Prof}$ such that the induced map
\[\begin{tikzcd}
	{\bigsqcup_{i} S_i} & S
	\arrow[from=1-1, to=1-2]
\end{tikzcd}\]
is surjective. We write $\textup{Cov}(S)$ for the class of all covers of $S$.

\item
One checks that the family $(\textup{Cov}(S))_{S \in |\ub{Prof}|}$ defines a pretopology on the category $\ub{Prof}$. We always endow $\ub{Prof}$ with this pretopology and consider it as a site.

\item The category $\ub{CSet}$ is defined to be the category of sheaves of sets on the site $\ub{Prof}$. Objects of $\ub{CSet}$ are called \ti{condensed sets}, and morphisms in $\ub{CSet}$ are called \ti{maps of condensed sets}.
\end{enumerate}
\end{df}

By \cref{cor:Prof is essentially small}, the set-theoretic size problem does not occur in this definition.

\begin{prop} \label{prop:monics and epics in CSet}
Let $X \xto{f} Y$ be a map of condensed sets.
\begin{enumerate}
\item $X \xto{f} Y$ is a monomorphism in $\ub{CSet}$ if and only if for every $S \in |\ub{Prof}|$, the map $X(S) \xto{f_S} Y(S)$ is injective.

\item $X \xto{f} Y$ is an epimorphism in $\ub{CSet}$ if and only if for every $S \in |\ub{Prof}|$ and $y \in Y(S)$, there exists a cover $(c_i : S_i \to S)_{i \in I}$ of $S$ in $\ub{Prof}$ such that $Y(c_i)(y) \in f_{S_i} \big( X(S_i) \big)$ for all $i \in I$.
\[\begin{tikzcd}
	& {Y(S)} \\
	{X(S_i)} & {Y(S_i)}
	\arrow["{Y(c_i)}", from=1-2, to=2-2]
	\arrow["{f_{S_i}}"', from=2-1, to=2-2]
\end{tikzcd}\]

\item
$X \xto{f} Y$ is an isomorphism in $\ub{CSet}$ if and only if $X \xto{f} Y$ is both a monomorphism and an epimorphism in $\ub{CSet}$.
\end{enumerate}
\end{prop}

\begin{proof}
This follows from \cite[\href{https://stacks.math.columbia.edu/tag/00WN}{Tag 00WN}]{stacks-project}.
\end{proof}

\begin{prop} \label{prop:limits and filtered colimits in CSet}
Let $\cat{P}$ be the category of presheaves of sets on the site $\ub{Prof}$. 
\begin{enumerate}
\item
The subcategory $\ub{CSet}$ of $\cat{P}$ is closed under limits and filtered colimits in $\cat{P}$. 

\item
Let $S$ be a light profinite set. Then the functor
\[\begin{tikzcd}
	{\ub{CSet}} & {\ub{Set},} & X & {X(S)}
	\arrow["{\mathrm{ev}_{S}}", from=1-1, to=1-2]
	\arrow[maps to, from=1-3, to=1-4]
\end{tikzcd}\]
preserves limits and filtered colimits.
\end{enumerate}
\end{prop}

\begin{proof}~
\begin{enumerate}
\item
Since $\ub{CSet}$ is a reflective subcategory of $\cat{P}$ (\cite{Borceux:cat3}, Theorem 3.3.12), the subcategory $\ub{CSet}$ is closed under limits in $\cat{P}$ (\cite{Borceux:cat1}, Proposition 3.5.3). We show that $\ub{CSet}$ is closed under filtered colimits in $\cat{P}$. Let $\cat{I}$ be a small filtered category and let $F : \cat{I} \to \ub{CSet}$ be a functor. Let $L \in |\cat{P}|$ be the colimit of the functor $\cat{I} \xto{F} \ub{CSet} \mon \cat{P}$. If $(S_{\lambda} \to S)_{\lambda \in \Lambda}$ is a cover of $S \in |\ub{Prof}|$ in $\ub{Prof}$, then we have an exact diagram
\[\begin{tikzcd}
	{F(i)(S)} & {\prod_{\lambda \in \Lambda} F(i)(S_\lambda)} & {\prod_{\lambda,\mu \in \lambda} F(i)(S_\lambda \times_S S_{\mu})}
	\arrow[from=1-1, to=1-2]
	\arrow[shift left, from=1-2, to=1-3]
	\arrow[shift right, from=1-2, to=1-3]
\end{tikzcd}\]
in $\ub{Set}$ for each $i \in \ob{I}$. Note that the index set $\Lambda$ is finite. Since filtered colimits commutes with finite limits in $\ub{Set}$ (\cite{Borceux:cat1}, Theorem 2.13.4), taking colimits over $\cat{I}$ leads to an exact diagram
\[\begin{tikzcd}
	{\colim_{i \in \cat{I}} \, (F(i)(S))} & {\prod_{\lambda \in \Lambda} \colim_{i \in \cat{I}} \, (F(i)(S_\lambda))} & {\prod_{\lambda,\mu \in \lambda} \colim_{i \in \cat{I}} \, (F(i)(S_\lambda \times_S S_{\mu}))}
	\arrow[from=1-1, to=1-2]
	\arrow[shift left, from=1-2, to=1-3]
	\arrow[shift right, from=1-2, to=1-3]
\end{tikzcd}\]
in $\ub{Set}$. This diagram is equal to
\[\begin{tikzcd}
	{L(S)} & {\prod_{\lambda \in \Lambda} L(S_\lambda)} & {\prod_{\lambda,\mu \in \lambda} L(S_\lambda \times_S S_{\mu})}
	\arrow[from=1-1, to=1-2]
	\arrow[shift left, from=1-2, to=1-3]
	\arrow[shift right, from=1-2, to=1-3]
\end{tikzcd}\]
by the dual of Theorem 2.15.2 of \cite{Borceux:cat1}. Therefore $L$ is a sheaf on the site $\ub{Prof}$. In other words, we have $F \in |\ub{CSet}|$.

\item
By (1), the inclusion functor $\ub{CSet} \mon \cat{P}$ preserves limits and filtered colimits. Moreover, the functor
\[\begin{tikzcd}
	{\cat{P}} & {\ub{Set},} & X & {X(S)}
	\arrow["{\mathrm{ev}_S}", from=1-1, to=1-2]
	\arrow[maps to, from=1-3, to=1-4]
\end{tikzcd}\]
preserves limits and colimits (\cite{Borceux:cat1}, Theorem 2.15.2 and its dual). Therefore the composition $\ub{CSet} \mon \cat{P} \xto{\mathrm{ev}_S} \ub{Set}$ preserves limits and filtered colimits. This composition is equal to the functor $\ub{CSet} \to \ub{Set}$, $X \mapsto X(S)$.
\end{enumerate}
\end{proof}

\subsection{Condensed abelian groups}

\begin{df}
The category $\ub{CAb}$ is defined to be the category of sheaves of abelian groups on the site $\ub{Prof}$. Objects of $\ub{CAb}$ are called \ti{condensed abelian groups}, and morphisms in $\ub{CAb}$ are called \ti{homomorphisms of condensed abelian groups}.
\end{df}

\begin{prop} \label{prop:free condensed abelian groups}
\tup{(\cite{Camargo:note}, Example 2.3.2, (1) ; \cite{Kedlaya:note}, Definition 5.1.5)}

The forgetful functor $\ub{CAb} \to \ub{CSet}$ admits a left adjoint. It is denoted by
\[\begin{tikzcd}
	{\ub{CSet}} & {\ub{CAb},} & X & {\Z X.}
	\arrow[from=1-1, to=1-2]
	\arrow[maps to, from=1-3, to=1-4]
\end{tikzcd}\]
\end{prop}

\begin{prop} \label{prop:CAb is Grothendieck Abelian}
\tup{(\cite{Camargo:note}, Theorem 2.3.3)}

The category $\ub{CAb}$ is a Grothendieck Abelian category.
\end{prop}

\begin{prop} \label{prop:evaluation of CAb preserves limits filtered colimits and coproducts}
Let $\cat{P}(\ub{Ab})$ be the category of presheaves of abelian groups on the site $\ub{Prof}$.
\begin{enumerate}
\item
The subcategory $\ub{CAb}$ of $\cat{P}(\ub{Ab})$ is closed under limits, filtered colimits and coproducts in $\cat{P}(\ub{Ab})$.

\item
Let $S$ be a light profinite set. Then the functor
\[\begin{tikzcd}
	{\ub{CAb}} & {\ub{Ab},} & M & {M(S)}
	\arrow["{\mathrm{ev}_{S}}", from=1-1, to=1-2]
	\arrow[maps to, from=1-3, to=1-4]
\end{tikzcd}\]
preserves limits, filtered colimits and coproducts.
\end{enumerate}
\end{prop}

\begin{proof}~
\begin{enumerate}
\item
Since $\ub{CAb}$ is a reflective subcategory of $\cat{P}(\ub{Ab})$, the subcategory $\ub{CAb}$ is closed under limits in $\cat{P}(\ub{Ab})$ (\cite{Borceux:cat1}, Proposition 3.5.3). Consider the category $\cat{P}$ of presheaves of sets on the site $\ub{Prof}$ and the forgetful functor $U : \cat{P}(\ub{Ab}) \to \cat{P}$. Then an object $M \in |\cat{P}(\ub{Ab})|$ is an object of  $\ub{CAb}$ if and only if $U(M)$ is an object of $\ub{CSet}$. Since $U : \cat{P}(\ub{Ab}) \to \cat{P}$ preserves filtered colimits, \cref{prop:limits and filtered colimits in CSet} immediately implies that $\ub{CAb}$ is closed under filtered colimits in $\cat{P}(\ub{Ab})$. Then the subcategory $\ub{CAb}$ is closed under finite direct sums and filtered colimits in $\cat{P}(\ub{Ab})$. Since coproducts in $\cat{P}(\ub{Ab})$ can be written as filtered colimits of finite direct sums in $\cat{P}(\ub{Ab})$, we conclude that $\ub{CAb}$ is closed under coproducts in $\cat{P}(\ub{Ab})$.

\item
By (1), the inclusion functor $\ub{CAb} \mon \cat{P}(\ub{Ab})$ preserves limits, filtered colimits and coproducts. Moreover, the functor
\[\begin{tikzcd}
	{\cat{P}(\ub{Ab})} & {\ub{Ab},} & M & {M(S)}
	\arrow["{\mathrm{ev}_S}", from=1-1, to=1-2]
	\arrow[maps to, from=1-3, to=1-4]
\end{tikzcd}\]
preserves limits and colimits (\cite{Borceux:cat1}, Theorem 2.15.2 and its dual). Therefore the composition $\ub{CAb} \mon \cat{P}(\ub{Ab}) \xto{\mathrm{ev}_S} \ub{Ab}$ preserves limits, filtered colimits and coproducts. This composition is equal to the functor $\ub{CAb} \to \ub{Ab}$, $M \mapsto M(S)$.
\end{enumerate}
\end{proof}

\begin{df} (\cite{Camargo:note}, Theorem 2.3.3 ; \cite{Kedlaya:note}, Definition 5.2.1)
\begin{enumerate}
\item
If $M,N$ are condensed abelian groups, define $M \otimes N$ to be the sheafification of the presheaf of abelian groups on $\ub{Prof}$
\[\begin{tikzcd}
	S & {M(S) \otimes N(S).}
	\arrow[maps to, from=1-1, to=1-2]
\end{tikzcd}\]
This defines a functor $\ub{CAb} \times \ub{CAb} \xto{\otimes} \ub{CAb}$.

\item
One checks that this functor $\otimes$ defines a closed symmetric monoidal structure on the category $\ub{CAb}$. We always endow $\ub{CAb}$ with this monoidal structure and consider it as a closed symmetric monoidal category. The internal hom functor is denoted by
\[\begin{tikzcd}
	{(M,N)} & {\uhom(M,N).}
	\arrow[maps to, from=1-1, to=1-2]
\end{tikzcd}\]
\end{enumerate}
\end{df}

\begin{df}
Let $M,N,Q$ be condensed abelian groups. A map $\beta : M \times N \to Q$ of condensed sets is called \ti{biadditive} if the map $\beta_S : M(S) \times N(S) \to Q(S)$ is biadditive for every $S \in |\ub{Prof}|$.
\end{df}

\begin{prop}
Let $M,N$ be condensed abelian groups. Let $M \otimes^{\mathrm{pre}} N$ be the presheaf $S \mapsto M(S) \otimes N(S)$ of abelian groups on $\ub{Prof}$. Let $\otimes^{\mathrm{pre}} : M \times N \to M \otimes^{\mathrm{pre}} N$ be the map of condensed sets defined by
\[\begin{tikzcd}
	{M(S) \times N(S)} & {M(S) \otimes N(S),} & {(x,y)} & {x \otimes y} & {(S \in |\ub{Prof}|) .}
	\arrow["{\otimes^{\mathrm{pre}}_S}", from=1-1, to=1-2]
	\arrow[maps to, from=1-3, to=1-4]
\end{tikzcd}\]
Let $\otimes : M \times N \to M \otimes N$ be the composition of $\otimes^{\mathrm{pre}} : M \times N \to M \otimes^{\mathrm{pre}} N$ and the sheafification $M \otimes^{\mathrm{pre}} N \to M \otimes N$. Then the map $\otimes : M \times N \to M \otimes N$ of condensed sets has the following universal property.
\begin{enumerate}
\item
The map $\otimes : M \times N \to M \otimes N$ of condensed sets is biadditive.

\item
If $Q$ is any condensed abelian group and if $\beta : M \times N \to Q$ is any biadditive map, then there exists a unique homomorphism $\gamma : M \otimes N \to Q$ of condensed abelian groups such that the diagram
\[\begin{tikzcd}
	{M \times N} & Q \\
	{M \otimes N}
	\arrow["\beta", from=1-1, to=1-2]
	\arrow["\otimes"', from=1-1, to=2-1]
	\arrow["\gamma"', from=2-1, to=1-2]
\end{tikzcd}\]
is commutative.
\end{enumerate}
\end{prop}

\begin{proof}
This follows from the universal property of the usual tensor products of abelian groups and that of the sheafification.
\end{proof}

\subsection{Condensed rings}

\begin{df}
The category $\ub{CRing}$ is defined to be the category of sheaves of commutative unital rings on the site $\ub{Prof}$. Objects of $\ub{CRing}$ are called \ti{condensed rings}, and morphisms in $\ub{CRing}$ are called \ti{homomorphisms of condensed rings}.
\end{df}

\begin{prop} \label{prop:evaluation of CRing preserves limits and filtered colimits}
Let $\cat{P}(\ub{Ring})$ be the category of presheaves of commutative unital rings on the site $\ub{Prof}$.
\begin{enumerate}
\item
The subcategory $\ub{CRing}$ of $\cat{P}(\ub{Ring})$ is closed under limits and filtered colimits in $\cat{P}(\ub{Ring})$.

\item
Let $S$ be a light profinite set. Then the functor
\[\begin{tikzcd}
	{\ub{CRing}} & {\ub{Ring},} & R & {R(S)}
	\arrow["{\mathrm{ev}_{S}}", from=1-1, to=1-2]
	\arrow[maps to, from=1-3, to=1-4]
\end{tikzcd}\]
preserves limits and filtered colimits.
\end{enumerate}
\end{prop}

\begin{proof}~
\begin{enumerate}
\item
Since $\ub{CRing}$ is a reflective subcategory of $\cat{P}(\ub{Ring})$, the subcategory $\ub{CRing}$ is closed under limits in $\cat{P}(\ub{Ring})$ (\cite{Borceux:cat1}, Proposition 3.5.3). Consider the category $\cat{P}$ of presheaves of sets on the site $\ub{Prof}$ and the forgetful functor $U : \cat{P}(\ub{Ring}) \to \cat{P}$. Then an object $R \in |\cat{P}(\ub{Ring})|$ is an object of  $\ub{CRing}$ if and only if $U(R)$ is an object of $\ub{CSet}$. Since $U : \cat{P}(\ub{Ring}) \to \cat{P}$ preserves filtered colimits, \cref{prop:limits and filtered colimits in CSet} immediately implies that $\ub{CRing}$ is closed under filtered colimits in $\cat{P}(\ub{Ring})$.

\item
By (1), the inclusion functor $\ub{CRing} \mon \cat{P}(\ub{Ring})$ preserves limits and filtered colimits. Moreover, the functor
\[\begin{tikzcd}
	{\cat{P}(\ub{Ring})} & {\ub{Ring},} & R & {R(S)}
	\arrow["{\mathrm{ev}_S}", from=1-1, to=1-2]
	\arrow[maps to, from=1-3, to=1-4]
\end{tikzcd}\]
preserves limits and colimits (\cite{Borceux:cat1}, Theorem 2.15.2 and its dual). Therefore the composition $\ub{CRing} \mon \cat{P}(\ub{Ring}) \xto{\mathrm{ev}_S} \ub{Ring}$ preserves limits and filtered colimits. This composition is equal to the functor $\ub{CRing} \to \ub{Ring}$, $R \mapsto R(S)$.
\end{enumerate}
\end{proof}

Recall that $*$ denotes the topological space with a unique point (\cref{nt:one point topological space}).

\begin{prop} \label{prop:universality of monoid algebras over condensed rings}
Let $R$ be a condensed ring. Let $M$ be a commutative monoid.
\begin{enumerate}
\item
The presheaf of commutative unital rings on $\ub{Prof}$
\[\begin{tikzcd}
	S & {R(S)[M]}
	\arrow[maps to, from=1-1, to=1-2]
\end{tikzcd}\]
is in fact a sheaf on $\ub{Prof}$. Therefore we obtain a condensed ring, which we write $R[M]$.

\item
Let $\iota : R \to R[M]$ be the homomorphism of condensed rings defined by
\[\begin{tikzcd}
	{R(S)} & {R(S)[M],} & r & r & {(S \in |\ub{Prof}|).}
	\arrow["{\iota_S}", from=1-1, to=1-2]
	\arrow[maps to, from=1-3, to=1-4]
\end{tikzcd}\]
This has the following universal property.
\begin{itemize}
\item
For every condensed ring $R'$, every homomorphism $\phi : R \to R'$ of condensed rings and every homomorphism $\sigma : M \to R'(*)$ of the monoid $M$ into the multiplicative monoid of the unital ring $R'(*)$, there exists a unique homomorphism $\psi : R[M] \to R'$ of condensed rings with the following properties. 
\begin{enumerate}
\item
The map
\[\begin{tikzcd}
	M & {R(*)[M] = R[M](*)} & {R'(*)}
	\arrow["\inc", hook, from=1-1, to=1-2]
	\arrow["{\psi_*}", from=1-2, to=1-3]
\end{tikzcd}\]
is equal to the map $\sigma : M \to R'(*)$.

\item
The diagram
\[\begin{tikzcd}
	R & {R'} \\
	{R[M]}
	\arrow["\phi", from=1-1, to=1-2]
	\arrow["\iota"', from=1-1, to=2-1]
	\arrow["\psi"', from=2-1, to=1-2]
\end{tikzcd}\]
is commutative.
\end{enumerate}
\end{itemize}
We call $\iota : R \to R[M]$ the canonical inclusion.
\end{enumerate}
\end{prop}

\begin{proof}~
\begin{enumerate}
\item
The underlying presheaf of abelian groups of the presheaf $S \mapsto R(S)[M]$ of commutative unital rings on $\ub{Prof}$ is equal to $\bigoplus_{m \in M} R$ in the category of presheaves of abelian groups on $\ub{Prof}$. By \cref{prop:evaluation of CAb preserves limits filtered colimits and coproducts}, this is in fact a sheaf of abelian groups on the site $\ub{Prof}$.

\item
Let $R'$ be a condensed ring. Let $\phi : R \to R'$ be a homomorphism of condensed rings. Let $\sigma : M \to R'(*)$ be a homomorphism of the monoid $M$ into the multiplicative monoid of the unital ring $R'(*)$. For each $S \in |\ub{Prof}|$, the universal property of the monoid ring $R(S)[M]$ implies that there exists a unique ring homomorphism $\psi_S : R(S)[M] \to R'(S)$ such that the diagrams
\[\begin{tikzcd}
	{R(S)} & {R'(S)} & M & {R'(*)} \\
	{R(S)[M]} && {R(S)[M]} & {R'(S)}
	\arrow["{\phi_S}", from=1-1, to=1-2]
	\arrow["{\iota_S}"', from=1-1, to=2-1]
	\arrow["\sigma", from=1-3, to=1-4]
	\arrow["\inc"', hook, from=1-3, to=2-3]
	\arrow["{R'(!_S)}", from=1-4, to=2-4]
	\arrow["{\psi_S}"', from=2-1, to=1-2]
	\arrow["{\psi_S}"', from=2-3, to=2-4]
\end{tikzcd}\]
are both commutative, where $!_S$ denotes the unique morphism $S \to *$ in $\ub{Prof}$. Then the family $\psi := (\psi_S)_{S \in |\ub{Prof}^{\op}|}$ is a homomorphism $R[M] \to R'$ of condensed rings, and has the properties (a) and (b). On the other hand, if another homomorphism $\rho : R[M] \to R'$ of condensed rings has the properties (a) and (b), then for every $S \in |\ub{Prof}|$, the diagrams
\[\begin{tikzcd}
	{R(S)} & {R'(S)} & M & {R(*)[M]} & {R'(*)} \\
	{R(S)[M]} &&& {R(S)[M]} & {R'(S)}
	\arrow["{\phi_S}", from=1-1, to=1-2]
	\arrow["{\iota_S}"', from=1-1, to=2-1]
	\arrow["\inc", hook, from=1-3, to=1-4]
	\arrow["\sigma", curve={height=-30pt}, from=1-3, to=1-5]
	\arrow["\inc"', hook, from=1-3, to=2-4]
	\arrow["{\rho_*}", from=1-4, to=1-5]
	\arrow["{R[M](!_S)}", from=1-4, to=2-4]
	\arrow["{R'(!_S)}", from=1-5, to=2-5]
	\arrow["{\rho_S}"', from=2-1, to=1-2]
	\arrow["{\rho_S}"', from=2-4, to=2-5]
\end{tikzcd}\]
are commutative. Therefore $\rho_S = \psi_S$ for every $S \in |\ub{Prof}|$, and hence $\rho = \psi$.
\end{enumerate}
\end{proof}

\begin{df}
Let $R$ be a condensed ring. Let $n \in \N$.
\begin{enumerate}
\item
We write
\begin{equation}
R[X_1, \ldots, X_n] := R[\N^n] .
\end{equation}
We call this the \ti{polynomial ring in} $n$ \ti{variables over} $R$.

\item
For $1 \leq i \leq n$, let $e_i$ be the element of $\N^n$ whose $i$-th entry is $1$ and whose other entries are $0$. For $S \in |\ub{Prof}|$, the image of $e_i$ in $R[X_1, \ldots, X_n](S)$ will be denoted by $X_i$.
\end{enumerate}
\end{df}

\subsection{Condensed modules}

\begin{df}
Let $R$ be a condensed ring. The category $\ub{CMod}_R$ is defined to be the category of $R$-module objects in $\ub{CSet}$. Objects of $\ub{CMod}_R$ are called \ti{condensed} $R$-\ti{modues}, and morphisms in $\ub{CMod}_R$ are called \ti{homomorphisms of condensed} $R$-\ti{modues}.
\end{df}

\begin{prop}
Let $R$ be a condensed ring.
\begin{enumerate}
\item \tup{(\cite{Kedlaya:note}, Definition 10.3.2)} The forgetful functor $\ub{CMod}_R \to \ub{CAb}$ admits a left adjoint
\[\begin{tikzcd}
	{\ub{CAb}} & {\ub{CMod}_R ,} & M & {R \otimes M .}
	\arrow[from=1-1, to=1-2]
	\arrow[maps to, from=1-3, to=1-4]
\end{tikzcd}\]

\item
The forgetful functor $\ub{CMod}_R \to \ub{CSet}$ admits a left adjoint. It is denoted by
\[\begin{tikzcd}
	{\ub{CSet}} & {\ub{CMod}_R ,} & X & {R X.}
	\arrow[from=1-1, to=1-2]
	\arrow[maps to, from=1-3, to=1-4]
\end{tikzcd}\]
\end{enumerate}
\end{prop}

\begin{proof}~
\begin{enumerate}
\item
This follows from the property of the usual tensor products of abelian groups and the universal property of the sheafification.

\item
This follows from (1) and \cref{prop:free condensed abelian groups}.
\end{enumerate}
\end{proof}

\begin{prop} \label{prop:CMod R to CAb preserves and reflects exact sequences}
Let $R$ be a condensed ring. Then the category $\ub{CMod}_R$ is a Grothendieck Abelian category. The forgetful functor $\ub{CMod}_R \to \ub{CAb}$ preserves and reflects exact sequences.
\end{prop}

\begin{proof}
By \cite[\href{https://stacks.math.columbia.edu/tag/03DA}{Tag 03DA}]{stacks-project} and \cite[\href{https://stacks.math.columbia.edu/tag/03DB}{Tag 03DB}]{stacks-project}, the category $\ub{CMod}_R$ is a cocomplete Abelian category in which filtered colimits are exact. Moreover, by \cref{cor:Prof is essentially small}, there is a set $\cat{S}$ of representatives of isomorphism classes of $\ub{Prof}$. Using adjunction, one checks that $\oplus_{S \in \cat{S}} \, R \tu{S}$ is a generator of $\ub{CMod}_R$, where $\tu{S}$ denotes the representable functor $\ub{Prof}^{\op} \to \ub{Set}$ represented by $S \in |\ub{Prof}|$ (see \cref{df:condensed sets represented by topological spaces} below). From \cite[\href{https://stacks.math.columbia.edu/tag/03D9}{Tag 03D9}]{stacks-project}, one deduces that the forgetful functor $\ub{CMod}_R \to \ub{CAb}$ preserves and reflects exact sequences. 
\end{proof}

\begin{prop} \label{prop:evaluation of CMod R preserves limits filtered colimits and coproducts}
Let $R$ be a condensed ring. Let $\ub{PMod}_R$ be the category of $R$-module objects in the category $\cat{P}$ of presheaves of sets on the site $\ub{Prof}$.
\begin{enumerate}
\item
The subcategory $\ub{CMod}_R$ of $\ub{PMod}_R$ is closed under limits, filtered colimits and coproducts in $\ub{PMod}_R$.

\item
Let $S$ be a light profinite set. Then the functor
\[\begin{tikzcd}
	{\ub{CMod}_R} & {\ub{Mod}_{R(S)},} & M & {M(S)}
	\arrow["{\mathrm{ev}_{S}}", from=1-1, to=1-2]
	\arrow[maps to, from=1-3, to=1-4]
\end{tikzcd}\]
preserves limits, filtered colimits and coproducts.
\end{enumerate}
\end{prop}

\begin{proof}~
\begin{enumerate}
\item
Consider the category $\cat{P}(\ub{Ab})$ of presheaves of abelian groups on the site $\ub{Prof}$ and the forgetful functor $U : \ub{PMod}_R \to \cat{P}(\ub{Ab})$. Then an object $M \in |\ub{PMod}_R|$ is an object of  $\ub{CMod}_R$ if and only if $U(M)$ is an object of $\ub{CAb}$. Since $U : \ub{PMod}_R \to \cat{P}(\ub{Ab})$ preserves limits and colimits (\cite[\href{https://stacks.math.columbia.edu/tag/03DB}{Tag 03DB}]{stacks-project}), \cref{prop:evaluation of CAb preserves limits filtered colimits and coproducts} immediately implies that $\ub{CMod}_R$ is closed under limits, filtered colimits and coproducts in $\ub{PMod}_R$. 

\item
By (1), the inclusion functor $\ub{CMod}_R \mon \ub{PMod}_R$ preserves limits, filtered colimits and coproducts. Moreover, the functor
\[\begin{tikzcd}
	{\ub{PMod}_R} & {\ub{Mod}_{R(S)},} & M & {M(S)}
	\arrow["{\mathrm{ev}_{S}}", from=1-1, to=1-2]
	\arrow[maps to, from=1-3, to=1-4]
\end{tikzcd}\]
preserves limits and colimits. Therefore the composition $\ub{CMod}_R \mon \ub{PMod}_R \xto{\mathrm{ev}_S} \ub{Mod}_{R(S)}$ preserves limits, filtered colimits and coproducts. This composition is equal to the functor $\ub{CMod}_R \to \ub{Mod}_{R(S)}$, $M \mapsto M(S)$.
\end{enumerate}
\end{proof}

\begin{df} (\cite[\href{https://stacks.math.columbia.edu/tag/03EK}{Tag 03EK}]{stacks-project}, \cite[\href{https://stacks.math.columbia.edu/tag/03EO}{Tag 03EO}]{stacks-project})

Let $R$ be a condensed ring.
\begin{enumerate}
\item
If $M,N$ are condensed $R$-modules, define $M \otimes_R N$ to be the sheafification of the presheaf on $\ub{Prof}$
\[\begin{tikzcd}
	S & {M(S) \otimes_{R(S)} N(S).}
	\arrow[maps to, from=1-1, to=1-2]
\end{tikzcd}\]
This defines a functor $\ub{CMod}_R \times \ub{CMod}_R \xto{\otimes_R} \ub{CMod}_R$.

\item 
One checks that this functor $\otimes_R$ defines a closed symmetric monoidal structure on the category $\ub{CMod}_R$. We always endow $\ub{CMod}_R$ with this monoidal structure and consider it as a closed symmetric monoidal category. The internal hom functor is denoted by
\[\begin{tikzcd}
	{(M,N)} & {\uhom_{R}(M,N).}
	\arrow[maps to, from=1-1, to=1-2]
\end{tikzcd}\]
\end{enumerate}
\end{df}

\begin{df}
Let $R$ be a condensed ring. Let $M,N,Q$ be condensed $R$-modules. A map $\beta : M \times N \to Q$ of condensed sets is called $R$\ti{-bilinear} if the map $\beta_S : M(S) \times N(S) \to Q(S)$ is $R(S)$-bilinear for every $S \in |\ub{Prof}|$.
\end{df}

\begin{prop} \label{prop:tensor product over R represents R bilinear maps}
Let $R$ be a condensed ring. Let $M,N$ be condensed $R$-modules. Let $M \otimes^{\mathrm{pre}}_R N$ be the presheaf $S \mapsto M(S) \otimes_{R(S)} N(S)$ on $\ub{Prof}$. Let $\otimes^{\mathrm{pre}}_R : M \times N \to M \otimes^{\mathrm{pre}}_R N$ be the map of condensed sets defined by
\[\begin{tikzcd}
	{M(S) \times N(S)} & {M(S) \otimes_{R(S)} N(S),} & {(x,y)} & {x \otimes y} & {(S \in |\ub{Prof}|) .}
	\arrow["{(\otimes^{\mathrm{pre}}_R)_S}", from=1-1, to=1-2]
	\arrow[maps to, from=1-3, to=1-4]
\end{tikzcd}\]
Let $\otimes_R : M \times N \to M \otimes_R N$ be the composition of $\otimes^{\mathrm{pre}}_R : M \times N \to M \otimes^{\mathrm{pre}}_R N$ and the sheafification $M \otimes^{\mathrm{pre}}_R N \to M \otimes_R N$. Then the map $\otimes_R : M \times N \to M \otimes_R N$ of condensed sets has the following universal property.
\begin{enumerate}
\item
The map $\otimes_R : M \times N \to M \otimes_R N$ of condensed sets is $R$-biinear.

\item
If $Q$ is any condensed $R$-module and if $\beta : M \times N \to Q$ is any $R$-bilinear map, then there exists a unique homomorphism $\gamma : M \otimes_R N \to Q$ of condensed $R$-modules such that the diagram
\[\begin{tikzcd}
	{M \times N} & Q \\
	{M \otimes_R N}
	\arrow["\beta", from=1-1, to=1-2]
	\arrow["{\otimes_R}"', from=1-1, to=2-1]
	\arrow["\gamma"', from=2-1, to=1-2]
\end{tikzcd}\]
is commutative.
\end{enumerate}
\end{prop}

\begin{proof}
This follows from the universal property of the usual tensor products of modules over commutative unital rings and that of the sheafification.
\end{proof}

Recall that $*$ denotes the topological space with a unique point (\cref{nt:one point topological space}).

\begin{rem} \label{rem:M(S) is an R(*) module}
Let $R$ be a condensed ring. If $S$ is any light profinite set, then we have a unique morphism $!_S : S \to *$ in $\ub{Prof}$, which induces a ring homomorphism $R(!_S) : R(*) \to R(S)$. 
\begin{enumerate}
\item
If $M$ is any condensed $R$-module, then the $R(S)$-module $M(S)$ can be considered as an $R(*)$-module, by restricting scalars via $R(!_S) : R(*) \to R(S)$. In this way we always consider $M(S)$ as an $R(*)$-module.

\item
If $M,N$ are condensed $R$-modules, then $\ub{CMod}_R (M,N)$ becomes an $R(*)$-module if we define the action of $f \in R(*)$ on $\alpha \in \ub{CMod}_R (M,N)$ by
\begin{align}
(f \cdot \alpha)_{S} := R(!_S)(f) \cdot (\alpha_S) && (S \in |\ub{Prof}|).
\end{align}
In this way we always consider $\ub{CMod}_R (M,N)$ as an $R(*)$-module.
\end{enumerate}
\end{rem}

\begin{prop}
Let $R$ be a condensed ring. Let $T \sub R(*)$ be a multiplicatively closed subset. Let $M$ be a condensed $R$-module. Then presheaf on $\ub{Prof}$
\[\begin{tikzcd}
	S & {T^{-1} M(S)}
	\arrow[maps to, from=1-1, to=1-2]
\end{tikzcd}\]
is in fact a sheaf on $\ub{Prof}$. Therefore we obtain a condensed $R$-module. We write $T^{-1} M$ for this condensed $R$-module and call it the \ti{localization} of $M$ by $T$.
\end{prop}

\begin{proof}
Let $(S_{\lambda} \to S)_{\lambda \in \Lambda}$ be a cover of $S \in |\ub{Prof}|$ in $\ub{Prof}$. Then we have an exact diagram
\[\begin{tikzcd}
	{M(S)} & {\prod_{\lambda \in \Lambda} M(S_{\lambda})} & {\prod_{\lambda, \mu \in \Lambda} M(S_{\lambda} \times_S S_{\mu})}
	\arrow[from=1-1, to=1-2]
	\arrow[shift left, from=1-2, to=1-3]
	\arrow[shift right, from=1-2, to=1-3]
\end{tikzcd}\]
in $\ub{Set}$. This diagram is also exact in $\ub{Mod}_{R(*)}$. Note that the index set $\Lambda$ is finite. By the flatness of $T^{-1} R(*)$ over $R(*)$, tensoring this diagram with $T^{-1} R(*)$ leads to an exact diagram
\[\begin{tikzcd}
	{T^{-1}M(S)} & {\prod_{\lambda \in \Lambda} T^{-1}M(S_{\lambda})} & {\prod_{\lambda, \mu \in \Lambda} T^{-1}M(S_{\lambda} \times_S S_{\mu})}
	\arrow[from=1-1, to=1-2]
	\arrow[shift left, from=1-2, to=1-3]
	\arrow[shift right, from=1-2, to=1-3]
\end{tikzcd}\]
in $\ub{Mod}_{R(*)}$. This diagram is also exact in $\ub{Set}$. Therefore $S \mapsto T^{-1} M(S)$ defines a sheaf on the site $\ub{Prof}$.
\end{proof}

\begin{nt}
Let $R$ be a condensed ring. Let $M$ be a condensed $R$-module.
\begin{enumerate}
\item
Given an element $g \in R(*)$, we write $M_g$ for the localization of $M$ by the multiplicatively closed subset $\set{g^{n}}{n \in \N}$ of $R(*)$ generated by $g$. We call $M_g$ the localization of $M$ by $g$.

\item
Given a prime ideal $\fk{p}$ of $R(*)$, we write $M_{\fk{p}}$ for the localization of $M$ by the multiplicatively closed subset $R(*) \sm \fk{p}$ of $R(*)$. We call $M_{\fk{p}}$ the localization of $M$ at $\fk{p}$.
\end{enumerate}
\end{nt}

\subsection{Condensed algebras}

\begin{df}
Let $R$ be a condensed ring. The category $\ub{CAlg}_R$ is defined to be the category of associative commutative unital $R$-algebra objects in $\ub{CSet}$. Objects of $\ub{CAlg}_R$ are called \ti{condensed} $R$-\ti{algebras}, and morphisms in $\ub{CAlg}_R$ are called \ti{homomorphisms of condensed} $R$-\ti{algebras}.
\end{df}

\begin{rem}
Let $R$ be a condensed ring. 
\begin{enumerate}
\item
The category $\ub{CAlg}_R$ is canonically equivalent to the under category $R/\ub{CRing}$.

\item
Given an object $A$ of $\ub{CAlg}_R$, the corresponding object $R \to A$ of $R/\ub{CRing}$ is called the \ti{structure homomorphism} of $A$.

\item
Conversely, given an object $R \to S$ of $R/\ub{CRing}$, we have a corresponding object of $\ub{CAlg}_R$ whose underlying condensed ring is equal to $S$. This condensed $R$-algebra is referred to as $S$ that is considered as a condensed $R$-algebra via $R \to S$.

\item
Unless otherwise stated explicitly, we consider $R$ as a condensed $R$-algebra via the identity $\id_R : R \to R$.
\end{enumerate}
\end{rem}

\begin{prop} \label{prop:evaluation of CAlg R preserves limits and filtered colimits}
Let $R$ be a condensed ring. Let $\ub{PAlg}_R$ be the category of associative commutative unital $R$-algebra objects in the category $\cat{P}$ of presheaves of sets on the site $\ub{Prof}$.
\begin{enumerate}
\item
The subcategory $\ub{CAlg}_R$ of $\ub{PAlg}_R$ is closed under limits and filtered colimits $\ub{PAlg}_R$.

\item
Let $S$ be a light profinite set. Then the functor
\[\begin{tikzcd}
	{\ub{CAlg}_R} & {\ub{Alg}_{R(S)},} & A & {A(S)}
	\arrow["{\mathrm{ev}_{S}}", from=1-1, to=1-2]
	\arrow[maps to, from=1-3, to=1-4]
\end{tikzcd}\]
preserves limits and filtered colimits.
\end{enumerate}
\end{prop}

\begin{proof}~
\begin{enumerate}
\item
Consider the category $\ub{PMod}_R$ of $R$-module objects in the category $\cat{P}$ and the forgetful functor $U : \ub{PAlg}_R \to \ub{PMod}_R$. Then an object $A \in |\ub{PAlg}_R|$ is an object of  $\ub{CAlg}_R$ if and only if $U(A)$ is an object of $\ub{CMod}_R$. Since $U : \ub{PAlg}_R \to \ub{PMod}_R$ preserves limits and filtered colimits, \cref{prop:evaluation of CMod R preserves limits filtered colimits and coproducts} immediately implies that $\ub{CAlg}_R$ is closed under limits and filtered colimits in $\ub{PAlg}_R$. 

\item
By (1), the inclusion functor $\ub{CAlg}_R \mon \ub{PAlg}_R$ preserves limits and  filtered colimits. Moreover, the functor
\[\begin{tikzcd}
	{\ub{PAlg}_R} & {\ub{Alg}_{R(S)},} & A & {A(S)}
	\arrow["{\mathrm{ev}_{S}}", from=1-1, to=1-2]
	\arrow[maps to, from=1-3, to=1-4]
\end{tikzcd}\]
preserves limits and colimits. Therefore the composition $\ub{CAlg}_R \mon \ub{PAlg}_R \xto{\mathrm{ev}_S} \ub{Alg}_{R(S)}$ preserves limits and filtered colimits. This composition is equal to the functor $\ub{CAlg}_R \to \ub{Alg}_{R(S)}$, $A \mapsto A(S)$.
\end{enumerate}
\end{proof}

\begin{prop} \label{prop:universality of localization of condensed algebra}
Let $R$ be a condensed ring. Let $T \sub R(*)$ be a multiplicatively closed subset. Let $A$ be a condensed $R$-algebra. Let $\sigma : R \to A$ be the structure homomorphism of $A$.
\begin{enumerate}
\item
The localization $T^{-1} A$ of the condensed $R$-module $A$ by $T$ has a natural structure of a condensed $R$-algebra. We always consider $T^{-1}A$ as a condensed $R$-algebra via this structure. 

\item
The map $\iota : A \to T^{-1} A$ of condensed sets defined by
\[\begin{tikzcd}
	{A(S)} & {(T^{-1}A)(S) = T^{-1}(A(S)) ,} & a & {\displaystyle \frac{a}{1}} & {(S \in |\ub{Prof}|)}
	\arrow["{\iota_S}", from=1-1, to=1-2]
	\arrow[maps to, from=1-3, to=1-4]
\end{tikzcd}\]
is a homomorphism of condensed $R$-algebras. We call $\iota : A \to T^{-1} A$ the canonical homomorphism.

\item
The canonical homomorphism $\iota : A \to T^{-1} A$ has the following universal property.
\begin{enumerate}
\item
For every $f \in T$, the element $\iota_* ( \sigma_* (f) ) \in (T^{-1}A)(*)$ is invertible in $(T^{-1}A)(*)$.

\item
Let $B$ be any condensed $R$-algebra. Let $\phi : A \to B$ be any homomorphism of condensed $R$-algebras. Suppose that the element $\phi_* (\sigma_* (f) ) \in B(*)$ is invertible in $B(*)$ for every $f \in T$. Then there exists a unique homomorphism $\psi : T^{-1}A \to B$ of condensed $R$-algebras such that the diagram
\[\begin{tikzcd}
	A & B \\
	{T^{-1}A}
	\arrow["\phi", from=1-1, to=1-2]
	\arrow["\iota"', from=1-1, to=2-1]
	\arrow["\psi"', from=2-1, to=1-2]
\end{tikzcd}\]
is commutative.
\end{enumerate}
\end{enumerate}
\end{prop}

\begin{proof}~
\begin{itemize}
\item[(1),(2)]
For every $S \in \ub{Prof}$, we have
\begin{equation}
(T^{-1}A)(S) = T^{-1} \big( A(S) \big)
\end{equation}
by definition. This has a natural structure of a commutative unital ring and the map
\[\begin{tikzcd}
	{A(S)} & {(T^{-1}A)(S) = T^{-1} \big( A(S) \big) ,} & a & {\displaystyle \frac{a}{1}}
	\arrow["{\iota_S}", from=1-1, to=1-2]
	\arrow[maps to, from=1-3, to=1-4]
\end{tikzcd}\]
is a ring homomorphism. Thus $T^{-1}A$ becomes a condensed ring and the map $\iota : A \to T^{-1}A$ of condensed sets is a homomorphism of condensed rings. If we consider $T^{-1}A$ as a condensed algebra via $R \xto{\sigma} A \xto{\iota} T^{-1}A$, then $T^{-1}A$ becomes a condensed $R$-algebra and the map $\iota : A \to T^{-1}A$ of condensed sets is a homomorphism of condensed $R$-algebras.

\item[(3)]
Since $(T^{-1}A)(*) = T^{-1} \big( A(*) \big)$ by definition, it is clear that the element $\iota_* ( \sigma_* (f) ) \in (T^{-1}A)(*)$ is invertible in $(T^{-1}A)(*)$ for every $f \in T$.

Let $B$ be any condensed $R$-algebra. Let $\phi : A \to B$ be any homomorphism of condensed $R$-algebras. Suppose that the element $\phi_* (\sigma_* (f) ) \in B(*)$ is invertible in $B(*)$ for every $f \in T$. For each $S \in |\ub{Prof}|$, write $!_S$ for the unique morphism $S \to *$ in $\ub{Prof}$. Then the following diagram is commutative.
\[\begin{tikzcd}
	{R(*)} & {A(*)} & {B(*)} \\
	{R(S)} & {A(S)} & {B(S)}
	\arrow["{\sigma_*}", from=1-1, to=1-2]
	\arrow["{R(!_S)}"', from=1-1, to=2-1]
	\arrow["{\phi_*}", from=1-2, to=1-3]
	\arrow["{A(!_S)}", from=1-2, to=2-2]
	\arrow["{B(!_S)}", from=1-3, to=2-3]
	\arrow["{\sigma_S}"', from=2-1, to=2-2]
	\arrow["{\phi_S}"', from=2-2, to=2-3]
\end{tikzcd}\]
Recall from \cref{rem:M(S) is an R(*) module} that $A(S)$ is considered as an $R(*)$-module via $R(*) \xto{R(!_S)} R(S) \xto{\sigma_S} A(S)$. For every $f \in T$, the element $\phi_S( \sigma_S( R(!_S)(f) ) ) \in B(S)$ is invertible in $B(S)$ since it is equal to $B(!_S)( \phi_*( \sigma_*(f) ) ) \in B(S)$ and the element $\phi_* (\sigma_* (f) ) \in B(*)$ is already invertible in $B(*)$ by assumption. Therefore, by the universality of the localization $T^{-1} \big( A(S) \big)$, there exists a unique ring homomorphism $\psi_S : T^{-1} \big( A(S) \big) \to B(S)$ such that the diagram
\[\begin{tikzcd}
	{A(S)} & {B(S)} \\
	{T^{-1} \big( A(S) \big)}
	\arrow["{\phi_S}", from=1-1, to=1-2]
	\arrow["{\iota_S}"', from=1-1, to=2-1]
	\arrow["{\psi_S}"', from=2-1, to=1-2]
\end{tikzcd}\]
is commutative. Then we obtain a homomoprhism $\psi :=(\psi_S)_{S \in |\ub{Prof}^{\op}|} : T^{-1}A \to B$ of condensed rings such that the diagram
\[\begin{tikzcd}
	A & B \\
	{T^{-1}A}
	\arrow["\phi", from=1-1, to=1-2]
	\arrow["\iota"', from=1-1, to=2-1]
	\arrow["\psi"', from=2-1, to=1-2]
\end{tikzcd}\]
is commutative. This homomorphism $\psi : T^{-1}A \to B$ is a homomorphism of condensed $R$-algebras since the diagram
\[\begin{tikzcd}
	R & A & B \\
	& {T^{-1}A}
	\arrow["\sigma", from=1-1, to=1-2]
	\arrow["{\text{str. hom. of } B}", curve={height=-24pt}, from=1-1, to=1-3]
	\arrow["{\text{str. hom. of } T^{-1}A}"', from=1-1, to=2-2]
	\arrow["\phi", from=1-2, to=1-3]
	\arrow["\iota"', from=1-2, to=2-2]
	\arrow["\psi"', from=2-2, to=1-3]
\end{tikzcd}\]
is commutative in $\ub{CRing}$, where str. hom. stands for structure homomorphism.

If $\rho : T^{-1}A \to B$ is another homomorphism of condensed $R$-algebras such that the diagram
\[\begin{tikzcd}
	A & B \\
	{T^{-1}A}
	\arrow["\phi", from=1-1, to=1-2]
	\arrow["\iota"', from=1-1, to=2-1]
	\arrow["\rho"', from=2-1, to=1-2]
\end{tikzcd}\]
is commutative, then the diagram
\[\begin{tikzcd}
	{A(S)} & {B(S)} \\
	{T^{-1} \big( A(S) \big)}
	\arrow["{\phi_S}", from=1-1, to=1-2]
	\arrow["{\iota_S}"', from=1-1, to=2-1]
	\arrow["{\rho_S}"', from=2-1, to=1-2]
\end{tikzcd}\]
is commutative in $\ub{Ring}$ for every $S \in |\ub{Prof}|$. By the definition of $\psi_S$, we conclude that $\rho_S = \psi_S$ for every $S \in |\ub{Prof}|$. Therefore $\rho = \psi$.
\end{itemize}
\end{proof}

\subsection{Condensed sets represented by topological spaces}

\begin{nt}
If $X,Y$ are topological spaces, we write $\mathrm{Cont}(X,Y)$ for the set of all continuous maps $X \to Y$.
\end{nt}

\begin{df} \label{df:condensed sets represented by topological spaces}
(\cite{Camargo:note}, Example 2.2.3, (2) ; \cite{Kedlaya:note}, Lemma 4.1.5)

Let $X$ be a topological space. The sheaf of sets on $\ub{Prof}$ defined by
\[\begin{tikzcd}
	S & {\mathrm{Cont}(S,X)}
	\arrow[maps to, from=1-1, to=1-2]
\end{tikzcd}\]
is called the condensed set \ti{represented by} $X$, and is denoted by $\tu{X}$. We obtain a functor $\ub{Top} \to \ub{CSet}$, $X \mapsto \tu{X}$.
\end{df}

Recall that $*$ denotes the topological space with a unique point. The unique point of $*$ is also denoted by $*$. (\cref{nt:one point topological space})

\begin{rem} \label{rem:identification of underbar X (*) with X}
Let $X$ be a topological space. Then we have a bijection
\[\begin{tikzcd}
	{\tu{X}(*) = \mathrm{Cont}(*,X)} & {X,} & t & {t(*) .}
	\arrow["\sim", from=1-1, to=1-2]
	\arrow[maps to, from=1-3, to=1-4]
\end{tikzcd}\]
Via this bijection, we always identify $\tu{X}(*)$ with the set $X$.
\end{rem}

\begin{prop} \label{prop:reflection along underbar construction}
\tup{(\cite{Camargo:note}, Proposition 2.2.4 ; \cite{Kedlaya:note}, Definition 4.3.1 and Proposition 4.3.2)}

Let $X$ be a condensed set. For each light profinite set $S$ and each $x \in S$, write $p_{S,x} : * \to S$ for the continuous map $* \mapsto x$. For each light profinite set $S$ and each $t \in X(S)$, we define a map $\xi_{S,t} : S \to X(*)$ by
\[\begin{tikzcd}
	{\xi_{S,t} : S} & {X(*),} & x & {X(p_{S,x})(t).}
	\arrow[from=1-1, to=1-2]
	\arrow[maps to, from=1-3, to=1-4]
\end{tikzcd}\]
Let $\cat{T}$ be the finest topology on $X(*)$ such that the maps $\xi_{S,t} : S \to X(*) \; (S \in |\ub{Prof}| ,\, t \in X(S))$ are all continuous when $X(*)$ is endowed with $\cat{T}$. Let us consider $X(*)$ as a topological space by giving it the topology $\cat{T}$. 

\begin{enumerate}
\item
For each light profinite set $S$, we define a map $\eta_S : X(S) \to \mathrm{Cont}(S, X(*))$ by
\[\begin{tikzcd}
	{\eta_S : X(S)} & {\mathrm{Cont}(S,X(*)),} & t & {\xi_{S,t} .}
	\arrow[from=1-1, to=1-2]
	\arrow[maps to, from=1-3, to=1-4]
\end{tikzcd}\]
Then the family $\eta := (\eta_S)_{S \in |\ub{Prof}|}$ is a map $X \to \tu{X(*)}$ of condensed sets.

\item
For any topological spaces $Y$ and any map $\alpha : X \to \tu{Y}$ of condensed sets, there exists a unique continuous map $f : X(*) \to Y$ such that the following diagram is commutative.
\[\begin{tikzcd}
	X & {\tu{Y}} \\
	{\tu{X(*)}}
	\arrow["\alpha", from=1-1, to=1-2]
	\arrow["\eta"', from=1-1, to=2-1]
	\arrow["{\tu{f}}"', from=2-1, to=1-2]
\end{tikzcd}\]
$f : X(*) \to Y$ is equal to the map
\[\begin{tikzcd}
	{X(*)} & {\tu{Y}(*) = Y .}
	\arrow["{\alpha_{*}}", from=1-1, to=1-2]
\end{tikzcd}\]
\end{enumerate}
\end{prop}

\begin{proof}~
\begin{enumerate}
\item
Let $c : S' \to S$ be an arbitrary morphism in $\ub{Prof}$. We prove that the following diagram is commutative.
\[\begin{tikzcd}
	{X(S)} & {\mathrm{Cont}(S, X(*))} \\
	{X(S')} & {\mathrm{Cont}(S', X(*))}
	\arrow["{\eta_S}", from=1-1, to=1-2]
	\arrow["{X(c)}"', from=1-1, to=2-1]
	\arrow["{\phi \mapsto \phi \of c}", from=1-2, to=2-2]
	\arrow["{\eta_{S'}}"', from=2-1, to=2-2]
\end{tikzcd}\]
Let $t \in X(S)$. We show that $\xi_{S,t} \of c = \xi_{S' ,\, X(c)(t)}$. For $x' \in S'$, we have
\begin{align}
\xi_{S,t} \of c \, (x')
& = X(p_{S,c(x')})(t) \: ; \\
\xi_{S' ,\, X(c)(t)} (x')
& = X(p_{S' ,\, x'})(X(c)(t)) \\
& = X(c \of p_{S' ,\, x'})(t).
\end{align}
However, both $p_{S,c(x')}$ and $c \of p_{S' ,\, x'}$ are equal to the map $ : * \to S$, $* \mapsto c(x')$. Therefore
\begin{equation}
\xi_{S,t} \of c \, (x')
= X(p_{S,c(x')})(t)
= X(c \of p_{S' ,\, x'})(t)
= \xi_{S' ,\, X(c)(t)} (x') .
\end{equation}
It follows that $\xi_{S,t} \of c = \xi_{S' ,\, X(c)(t)}$.

\item
Let $f : X(*) \to Y$ be the map $\alpha_* : X(*) \to \tu{Y}(*) = Y$. First prove the following claim $(\diamond)$.
\begin{itemize}
\item[$(\diamond)$]
For every $S \in |\ub{Prof}|$ and every $t \in X(s)$, the map $f \of \xi_{S,t} : S \to Y$ is equal to the map $\alpha_S (t)$.
\end{itemize}
Let $S \in |\ub{Prof}|$ and $t \in X(s)$. For $x \in S$, we have
\begin{equation}
f \of \xi_{S,t} \, (x) = \alpha_* \big( X(p_{S,x}) \, (t) \big) .
\end{equation}
The following diagram is commutative.
\[\begin{tikzcd}
	{X(S)} & {\tu{Y}(S)} \\
	{X(*)} & {\tu{Y}(*)}
	\arrow["{\alpha_S}", from=1-1, to=1-2]
	\arrow["{X(p_{S,x})}"', from=1-1, to=2-1]
	\arrow["{\tu{Y}(p_{S,x})}", from=1-2, to=2-2]
	\arrow["{\alpha_*}"', from=2-1, to=2-2]
\end{tikzcd}\]
Therefore
\begin{equation}
\alpha_* \big( X(p_{S,x}) \, (t) \big)
= \tu{Y}(p_{S,x}) \of \alpha_S \, (t)
= ( \alpha_S(t) ) \of p_{S,x} . 
\end{equation}
Under the identification $\tu{Y}(*) = Y$ made in \cref{rem:identification of underbar X (*) with X}, the element $( \alpha_S(t) ) \of p_{S,x} \in \tu{Y}(*)$ is identified with $( \alpha_S(t) ) \of p_{S,x} \, (*) \in Y$. However,
\begin{equation}
( \alpha_S(t) ) \of p_{S,x} \, (*) = \alpha_S (t)(x) .
\end{equation}
Consequently, the map $f \of \xi_{S,t} : S \to Y$ is equal to the map $x \mapsto \alpha_S(t) (x)$. In other words, $f \of \xi_{S,t} = \alpha_S(t)$. This completes the proof of the claim $(\diamond)$.

Next we prove that the map $f : X(*) \to Y$ is continuous. By the definition of the topology $\cat{T}$, it suffices to prove that the composition $f \of \xi_{S,t} : S \to Y$ is continuous for every $S \in |\ub{Prof}|$ and $t \in X(s)$. However, for every $S \in |\ub{Prof}|$ and $t \in X(s)$, we have
\begin{equation}
f \of \xi_{S,t} = \alpha_S(t) \in \tu{Y}(S) = \mathrm{Cont}(S,Y)
\end{equation}
by the claim $(\diamond)$. Thus the map $f \of \xi_{S,t} : S \to Y$ is continuous.

Thus we have a continuous map $f : X(*) \to Y$. For every $S \in |\ub{Prof}|$, the claim $(\diamond)$ shows that the following diagram is commutative.
\[\begin{tikzcd}
	{X(S)} & {\tu{Y}(S)} \\
	{\tu{X(*)}(S)}
	\arrow["{\alpha_S}", from=1-1, to=1-2]
	\arrow["{\eta_S}"', from=1-1, to=2-1]
	\arrow["{\tu{f}_S}"', from=2-1, to=1-2]
\end{tikzcd}\]
Therefore the diagram
\[\begin{tikzcd}
	X & {\tu{Y}} \\
	{\tu{X(*)}}
	\arrow["\alpha", from=1-1, to=1-2]
	\arrow["\eta"', from=1-1, to=2-1]
	\arrow["{\tu{f}}"', from=2-1, to=1-2]
\end{tikzcd}\]
is commutative. 

On the other hand, suppose that $g : X(*) \to Y$ is continuous map and that the diagram
\[\begin{tikzcd}
	X & {\tu{Y}} \\
	{\tu{X(*)}}
	\arrow["\alpha", from=1-1, to=1-2]
	\arrow["\eta"', from=1-1, to=2-1]
	\arrow["{\tu{g}}"', from=2-1, to=1-2]
\end{tikzcd}\]
is commutative. We prove that $g=f$. The following diagram is commutative.
\[\begin{tikzcd}
	{X(*)} & {\tu{Y}(*)} \\
	{\tu{X(*)}(*)}
	\arrow["{\alpha_*}", from=1-1, to=1-2]
	\arrow["{\eta_*}"', from=1-1, to=2-1]
	\arrow["{\tu{g}_*}"', from=2-1, to=1-2]
\end{tikzcd}\]
Therefore, for $t \in X(*)$, we have
\begin{equation}
\alpha_* (t) = \tu{g}_* \big( \eta_* (t) \big) = g \of \xi_{*,t} ,
\end{equation}
which shows
\begin{equation}
\alpha_* (t)(*) = g \big( \xi_{*,t} (*) \big) = g \big( X(p_{*,*})(t) \big) .
\end{equation}
On the other hand, the map $p_{*,*} : * \to *$ is equal to the identity $\id_* : * \to *$. Therefore
\begin{equation}
X(p_{*,*})(t) = \id_{X(*)} (t) = t .
\end{equation}
Consequently,
\begin{equation}
\alpha_* (t)(*) = g \big( X(p_{*,*})(t) \big) = g(t) .
\end{equation}
On the other hand, the map $f : X(*) \to Y$ was defined to be the map $\alpha_* : X(*) \to \tu{Y}(*) = Y$. Under the identification $\tu{Y}(*) = Y$ made in \cref{rem:identification of underbar X (*) with X}, the element $\alpha_* (t) \in \tu{Y}(*)$ is identified with $\alpha_* (t)(*) \in Y$. Therefore
\begin{equation}
f(t) = \alpha_* (t)(*) = g(t).
\end{equation}
Since this holds for any $t \in X(*)$, we conclude that $f=g$.
\end{enumerate}
\end{proof}

\begin{prop} \label{prop:functoriality of reflection along underbar construction}
Let $X,X'$ be condensed sets. For each light profinite set $S$ and each $x \in S$, write $p_{S,x} : * \to S$ for the continuous map $* \mapsto x$. For each light profinite set $S$ and each $t \in X(S)$ (resp.$\; t' \in X'(S)$), let $\xi_{S,t} : S \to X(*)$ (resp.$\; \xi'_{S,t'} : S \to X'(*)$) be the map $x \mapsto X(p_{S,x})(t)$ (resp.$\; x \mapsto X'(p_{S,x})(t')$). Let us consider $X(*)$ (resp.$\; X'(*)$) as a topological space by giving it the finest topology such that the maps $\xi_{S,t} : S \to X(*) \; (S \in |\ub{Prof}| ,\, t \in X(S))$ (resp.$\; \xi'_{S,t'} : S \to X'(*) \; (S \in |\ub{Prof}| ,\, t' \in X'(S))$) are all continuous. For each light profinite set $S$, let $\eta_S : X(S) \to \mathrm{Cont}(S, X(*))$ (resp.$\; \eta'_S : X'(S) \to \mathrm{Cont}(S, X'(*))$) be the map $t \mapsto \xi_{S,t}$ (resp.$\; t' \mapsto \xi'_{S,t'}$). By (1) of \cref{prop:reflection along underbar construction}, the family $\eta := (\eta_S)_{S \in |\ub{Prof}|}$ (resp.$\; \eta' := (\eta'_S)_{S \in |\ub{Prof}|}$) is a map $X \to \tu{X(*)}$ (resp.$\; X' \to \tu{X'(*)}$) of condensed sets.

If $f : X \to X'$ is a map of condensed sets, then the map $f_* : X(*) \to X'(*)$ is continuous and the following diagram is commutative.
\[\begin{tikzcd}
	X & {X'} \\
	{\tu{X(*)}} & {\tu{X'(*)}}
	\arrow["f", from=1-1, to=1-2]
	\arrow["\eta"', from=1-1, to=2-1]
	\arrow["{\eta'}", from=1-2, to=2-2]
	\arrow["{\tu{f_*}}"', from=2-1, to=2-2]
\end{tikzcd}\]
\end{prop}

\begin{proof}
By (2) of \cref{prop:reflection along underbar construction}, there exists a unique continuous map $g : X(*) \to X'(*)$ such that the following diagram is commutative.
\[\begin{tikzcd}
	X & {X'} \\
	{\tu{X(*)}} & {\tu{X'(*)}}
	\arrow["f", from=1-1, to=1-2]
	\arrow["\eta"', from=1-1, to=2-1]
	\arrow["{\eta'}", from=1-2, to=2-2]
	\arrow["{\tu{g}}"', from=2-1, to=2-2]
\end{tikzcd}\]
We prove that $g = f_*$. (2) of \cref{prop:reflection along underbar construction} shows that $g : X(*) \to X'(*)$ is equal to the map
\[\begin{tikzcd}
	{X(*)} & {X'(*)} & {\tu{X'(*)}(*) = X'(*).}
	\arrow["{f_*}", from=1-1, to=1-2]
	\arrow["{\eta'_*}", from=1-2, to=1-3]
\end{tikzcd}\]
Then it suffices to show that the map $\eta'_* : X'(*) \to \tu{X'(*)}(*) = X'(*)$ is equal to the identity. Let $t' \in X'(*)$. Under the identification $\tu{X'(*)}(*) = X'(*)$ made in \cref{rem:identification of underbar X (*) with X}, the element $\eta'_*(t') = \xi'_{*,t'} \in \tu{X'(*)}(*)$ is identified with $\xi'_{*,t'}(*) \in X'(*)$. On the other hand, we have
\begin{equation}
\xi'_{*,t'}(*) = X'(p_{*,*})(t') = t' ,
\end{equation}
since the map $p_{*,*} : * \to *$ is equal to the identity. It follows that the map $\eta'_* : X'(*) \to \tu{X'(*)}(*) = X'(*)$ is equal to the identity.
\end{proof}

\begin{cor} \label{cor:left adjoint of the underbar functor}
We have a functor
\[\begin{tikzcd}[ampersand replacement=\&,row sep=tiny]
	{\ub{CSet}} \& {\ub{Top}} \\
	X \& {\left( \begin{gathered}
	X(*) \text{ with the topology} \\
	\text{defined in \cref{prop:reflection along underbar construction}}
	\end{gathered} \right)} \& {(\text{on objects})} \\
	f \& {f_*} \& {(\text{on morphisms})}
	\arrow[from=1-1, to=1-2]
	\arrow[maps to, from=2-1, to=2-2]
	\arrow[maps to, from=3-1, to=3-2]
\end{tikzcd}\]
and this functor is left adjoint to the functor $\ub{Top} \to \ub{CSet}$, $X \mapsto \tu{X}$.
\end{cor}

\begin{rem} \label{rem:comparison of topologies of underbar X(*)}
Let $X$ be a topological space. Let $\cat{T}$ be the topology of $X$. Let $\cat{T}'$ be the topology on the set $\tu{X}(*) = X$ obtained by applying \cref{prop:reflection along underbar construction} to the condensed set $\tu{X}$. Then the topology $\cat{T}'$ is in general equal to or finer than the topology $\cat{T}$. If $X$ is first countable, then the topology $\cat{T}'$ coincides with the topology $\cat{T}$ (\cite{Kedlaya:note}, Proposition 4.4.2).
\end{rem}

\begin{rem} \,
\begin{enumerate}
\item
If $X$ is a topological abelian group (resp. commutative unital topological ring), then $\tu{X}$ has a natural structure of a condensed abelian group (resp. condensed ring). In this case we always consider $\tu{X}$ as a condensed abelian group (resp. condensed ring). Thus $X \mapsto \tu{X}$ defines a functor $\ub{TopAb} \to \ub{CAb}$ (resp. $\ub{TopRing} \to \ub{CRing}$).

\item
Let $R$ be a commutative unital topological ring. If $X$ is a topological $R$-module (resp. commutative unital topological $R$-algebra), then $\tu{X}$ has a natural structure of a condensed $\tu{R}$-module (resp. condensed $\tu{R}$-algebra). In this case we always consider $\tu{X}$ as a condensed $\tu{R}$-module (resp. condensed $\tu{R}$-algebra). Thus $X \mapsto \tu{X}$ defines a functor $\ub{TopMod}_R \to \ub{CMod}_{\tu{R}}$ (resp. $\ub{TopAlg}_R \to \ub{CAlg}_{\tu{R}}$).
\end{enumerate}
\end{rem}

\subsection{Coalescence}

\subsubsection{The sequence space}

\begin{nt} \label{nt:N infty} \,
\begin{enumerate}
\item (\cite{Kedlaya:note}, Example 3.2.3 and Remark 3.2.6) We write $\N_{\infty}$ for the one-point compactification of the discrete space $\N$. The point at infinity is denoted by $\infty$.

\item
Recall that $*$ denotes the topological space with a unique point. The unique point of $*$ is also denoted by $*$. (\cref{nt:one point topological space})
\end{enumerate}
\end{nt}

\begin{nt} \label{nt:translation and point map of N infty} \,
\begin{enumerate}
\item
We write $\tau$ for the continuous map
\[\begin{tikzcd}[ampersand replacement=\&]
	{\N_{\infty}} \& {\N_{\infty},} \& n \& {\left\{
	\begin{aligned}
	& n+1 & (n \in \N) \\ & \infty & (n=\infty)
	\end{aligned}
	\right. .}
	\arrow[from=1-1, to=1-2]
	\arrow[maps to, from=1-3, to=1-4]
\end{tikzcd}\]

\item For every $n \in \N_{\infty}$, we write $i_{n}$ for the continuous map
\[\begin{tikzcd}
	{*} & {\N_{\infty},} & {*} & {n.}
	\arrow[from=1-1, to=1-2]
	\arrow[maps to, from=1-3, to=1-4]
\end{tikzcd}\]
\end{enumerate}
\end{nt}

\begin{nt} \label{nt:relative translation and point map of N infty}
Let $S$ be a light profinite set.
\begin{enumerate}
\item
We write $\tau_S$ for the continuous map
\[\begin{tikzcd}
	{S \times \N_{\infty}} & {S \times \N_{\infty},} & {(x,n)} & {(x,\tau(n)).}
	\arrow[from=1-1, to=1-2]
	\arrow[maps to, from=1-3, to=1-4]
\end{tikzcd}\]

\item
For every $n \in \N_{\infty}$, we write $i_{S,n}$ for the continuous map
\[\begin{tikzcd}
	S & {S \times \N_{\infty},} & x & {(x,n).}
	\arrow[from=1-1, to=1-2]
	\arrow[maps to, from=1-3, to=1-4]
\end{tikzcd}\]
\end{enumerate}
\end{nt}

\begin{prop} \tup{(\cite{Camargo:note}, Example 2.1.8, (1))}

The topological space $\N_{\infty}$ is a light profinite set.
\end{prop}

\begin{proof}
Since the discrete space $\N$ is locally compact and Hausdorff, its one-point compactification $\N_{\infty}$ is compact Hausdorff (\cite{Bourbaki:top1}, Chapter I, \S 9.8, Theorem 4). Moreover, the set
\begin{equation}
\set{ \{n\} }{n \in \N} \cup \set{ \set{m \in \N}{m \geq n} \cup \{\infty\} }{n \in \N}
\end{equation}
is an open basis for the topology of $\N_{\infty}$. Since this set is countable, the topological space $\N_{\infty}$ is second countable. Let us show that the topological space $\N_{\infty}$ is totally disconnected. It suffices to show that for any distinct points $n,m$ of $\N_{\infty}$, there exists a subset $S$ of $\N_{\infty}$ which is both open and closed in $\N_{\infty}$ and satisfies $n \in S$ and $m \notin S$. Since $n \neq m$, we have either $n \neq \infty$ or $m \neq \infty$, and we may assume $n \neq \infty$. Then $S := \{n\}$ is both open and closed in $\N_{\infty}$ and satisfies $n \in S$ and $m \notin S$.
\end{proof}

\begin{df} \label{df:definition of sequence space}
Let $R$ be a condensed ring.
\begin{enumerate}
\item (\cite{Kedlaya:note}, Definition 6.1.1 and Definition 10.3.2)

The \ti{sequence space over} $R$ is defined to be the condensed $R$-module $P_R$ which is the cokernel of the homomorphism
\[\begin{tikzcd}
	{R \, \tu{*}} & {R \, \tu{\N_{\infty}}}
	\arrow["{R \, \tu{i_{\infty}}}", from=1-1, to=1-2]
\end{tikzcd}\]
in $\ub{CMod}_R$. We have a canonical epimorphism $R \, \tu{\N_{\infty}} \epi P_R$.

\item (\cite{Kedlaya:note}, Definition 6.1.1 and Definition 10.3.2)

The \ti{left shift} on $P_R$ is defined to be the homomorphism $P_R \xto{\sigma_R} P_R$ induced by $R \, \tu{\tau} : R \, \tu{\N_{\infty}} \to R \, \tu{\N_{\infty}}$. In other words, $\sigma_R$ is the unique homomorphism of condensed $R$-modules such that the diagram
\[\begin{tikzcd}
	{R \, \tu{\N_{\infty}}} & {P_R} \\
	{R \, \tu{\N_{\infty}}} & {P_R}
	\arrow[two heads, from=1-1, to=1-2]
	\arrow["{R \, \tu{\tau}}"', from=1-1, to=2-1]
	\arrow["{\sigma_R}", from=1-2, to=2-2]
	\arrow[two heads, from=2-1, to=2-2]
\end{tikzcd}\]
is commutative, where the horizontal maps are the canonical epimorphisms.

\item (cf. \cite{Kedlaya:note}, Definition 10.3.2)

For $f,g \in R(*)$, we define $\Delta_{f/g}$ to be the homomorphism of condensed $R$-modules
\[\begin{tikzcd}
	{g \cdot \id_{P_R} - f \cdot \sigma_R : P_R} & {P_R.}
	\arrow[from=1-1, to=1-2]
\end{tikzcd}\]
\end{enumerate}
\end{df}

\subsubsection{Coalescent modules}

\begin{lem} \label{lem:description of uhom R (PR M) (S)}
Let $R$ be a condensed ring.
\begin{enumerate}
\item \tup{(cf. \cite{Kedlaya:note}, Definition 5.3.1)}
The sequence
\[\begin{tikzcd}
	0 & {R \, \tu{*}} & {R \, \tu{\N_{\infty}}} & {P_R} & 0
	\arrow[from=1-1, to=1-2]
	\arrow["{R \, \tu{i_{\infty}}}", from=1-2, to=1-3]
	\arrow[two heads, from=1-3, to=1-4]
	\arrow[from=1-4, to=1-5]
\end{tikzcd}\]
is split exact in $\ub{CMod}_R$, where $R \, \tu{\N_{\infty}} \epi P_R$ denotes the canonical epimorphism.

\item
For every light profinite set $S$ and every condensed $R$-module $M$, there exist canonical isomorphisms of abelian groups
\begin{equation}
\uhom_R(P_R, M)(S) \simeq
\ub{CMod}_R \big( R \tu{S} \otimes_R P_R , M \big) \simeq
\ker M(i_{S,\infty})
\end{equation}
which are functorial in both $S \in |\ub{Prof}^{\op}|$ and $M \in |\ub{CMod}_R|$.
\end{enumerate}
\end{lem}

\begin{proof}~
\begin{enumerate}
\item
Let $!_{\N_{\infty}} : \N_{\infty} \to *$ be the unique morphism $\N_{\infty} \to *$ in $\ub{Prof}$. Then the diagram
\[\begin{tikzcd}
	{*} && {\N_{\infty}} \\
	& {*}
	\arrow["{i_{\infty}}", from=1-1, to=1-3]
	\arrow["{\id_*}"', from=1-1, to=2-2]
	\arrow["{!_{\N_{\infty}}}", from=1-3, to=2-2]
\end{tikzcd}\]
is commutative in $\ub{Prof}$. Therefore
\[\begin{tikzcd}
	{R \, \tu{*}} && {R \, \tu{\N_{\infty}}} \\
	& {R \, \tu{*}}
	\arrow["{R \, \tu{i_{\infty}}}", from=1-1, to=1-3]
	\arrow["{\id_{R \, \tu{*}}}"', from=1-1, to=2-2]
	\arrow["{R \, \tu{!_{\N_{\infty}}}}", from=1-3, to=2-2]
\end{tikzcd}\]
is commutative in $\ub{CMod}_R$.

\item
By Yoneda lemma, we have canonical bijections
\begin{align}
\uhom_R(P_R, M)(S) & \simeq
\ub{CSet} \big( \tu{S}, \uhom_R(P_R, M) \big) \: , \\
M(S \times \N_{\infty}) & \simeq
\ub{CSet} \big( \tu{S \times \N_{\infty}}, M \big) \: , \\
M(S) & \simeq
\ub{CSet} \big( \tu{S}, M \big) .
\end{align}
By adjunction, we have canonical bijections
\begin{align}
\ub{CSet} \big( \tu{S}, \uhom_R(P_R, M) \big) & \simeq
\ub{CMod}_R \big( R \tu{S}, \uhom_R(P_R, M) \big) \simeq
\ub{CMod}_R \big( R \tu{S} \otimes_R P_R , M \big) \: , \\
\ub{CSet} \big( \tu{S \times \N_{\infty}}, M \big) & \simeq
\ub{CMod}_R \big( R \, \tu{S \times \N_{\infty}}, M \big) \: , \\
\ub{CSet} \big( \tu{S}, M \big) & \simeq
\ub{CMod}_R \big( R \tu{S}, M \big) .
\end{align}
By the construction of the tensor product $\otimes_R$, one can check that there exist canonical isomorphisms
\begin{align}
R \, \tu{S \times \N_{\infty}} & \simeq 
R \tu{S} \otimes_R R \, \tu{\N_{\infty}} \: , \\
R \tu{S} \simeq R \, \tu{S \times *} & \simeq
R \tu{S} \otimes_R R \, \tu{*} 
\end{align}
and that the diagram
\[\begin{tikzcd}
	{R \tu{S}} && {R \, \tu{S \times \N_{\infty}}} \\
	{R \tu{S} \otimes_R R \, \tu{*}} && {R \tu{S} \otimes_R R \, \tu{\N_{\infty}}}
	\arrow["{R \, \tu{i_{S,\infty}}}", from=1-1, to=1-3]
	\arrow["{\text{\rotatebox{90}{$\sim$}}}", no head, from=2-1, to=1-1]
	\arrow["{R \, \tu{\id_S} \otimes_R R \, \tu{i_{\infty}}}"', from=2-1, to=2-3]
	\arrow["{\text{\rotatebox{90}{$\sim$}}}", no head, from=2-3, to=1-3]
\end{tikzcd}\]
is commutative. Therefore we have canonical bijections
\begin{align}
\ub{CMod}_R \big( R \, \tu{S \times \N_{\infty}}, M \big) & \simeq
\ub{CMod}_R \big( R \tu{S} \otimes_R R \, \tu{\N_{\infty}} , M \big) \: , \\
\ub{CMod}_R \big( R \tu{S}, M \big) & \simeq
\ub{CMod}_R \big( R \tu{S} \otimes_R R \, \tu{*} , M \big)
\end{align}
and the diagram
\[\begin{tikzcd}
	{\ub{CMod}_R \big( R \tu{S}, M \big)} && {\ub{CMod}_R \big( R \, \tu{S \times \N_{\infty}}, M \big)} \\
	{\ub{CMod}_R \big( R \tu{S} \otimes_R R \, \tu{*} , M)} && {\ub{CMod}_R \big( R \tu{S} \otimes_R R \, \tu{\N_{\infty}} , M)}
	\arrow["{\big( R \, \tu{i_{S,\infty}} \big)^*}"', from=1-3, to=1-1]
	\arrow["{\text{\rotatebox{90}{$\sim$}}}", no head, from=2-1, to=1-1]
	\arrow["{\text{\rotatebox{90}{$\sim$}}}", no head, from=2-3, to=1-3]
	\arrow["{\big( R \, \tu{\id_S} \otimes_R R \, \tu{i_{\infty}} \big)^*}", from=2-3, to=2-1]
\end{tikzcd}\]
is commutative. Consequently, we have canonical bijections
\begin{align}
\uhom_R(P_R, M)(S) & \simeq
\ub{CMod}_R \big( R \tu{S} \otimes_R P_R , M \big) \: , \\
M(S \times \N_{\infty}) & \simeq
\ub{CMod}_R \big( R \tu{S} \otimes_R R \, \tu{\N_{\infty}} , M \big) \: , \\
M(S) & \simeq
\ub{CMod}_R \big( R \tu{S} \otimes_R R \, \tu{*} , M \big) .
\end{align}
One can check that these are isomorphisms of abelian groups. Moreover, the diagram
\[\begin{tikzcd}
	{M(S)} && {M(S \times \N_{\infty})} \\
	{\ub{CMod}_R \big( R \tu{S} \otimes_R R \, \tu{*} , M)} && {\ub{CMod}_R \big( R \tu{S} \otimes_R R \, \tu{\N_{\infty}} , M)}
	\arrow["{M(i_{S,\infty})}"', from=1-3, to=1-1]
	\arrow["{\text{\rotatebox{90}{$\sim$}}}", no head, from=2-1, to=1-1]
	\arrow["{\text{\rotatebox{90}{$\sim$}}}", no head, from=2-3, to=1-3]
	\arrow["{\big( R \, \tu{\id_S} \otimes_R R \, \tu{i_{\infty}} \big)^*}", from=2-3, to=2-1]
\end{tikzcd}\]
is commutative. On the other hand, since the sequence
\[\begin{tikzcd}
	0 & {R \, \tu{*}} & {R \, \tu{\N_{\infty}}} & {P_R} & 0
	\arrow[from=1-1, to=1-2]
	\arrow["{R \, \tu{i_{\infty}}}", from=1-2, to=1-3]
	\arrow["\can", two heads, from=1-3, to=1-4]
	\arrow[from=1-4, to=1-5]
\end{tikzcd}\]
is split exact in $\ub{CMod}_R$ by (1), the sequence
\[\begin{tikzcd}[column sep=huge]
	0 & {\ub{CMod}_R \big( R \tu{S} \otimes R \, \tu{*} , M)} & {\ub{CMod}_R \big( R \tu{S} \otimes R \, \tu{\N_{\infty}} , M)} \\
	& {} & {\ub{CMod}_R \big( R \tu{S} \otimes P_R , M)} & 0
	\arrow[from=1-2, to=1-1]
	\arrow["{\big( R \, \tu{\id_S} \otimes_R R \, \tu{i_{\infty}} \big)^*}"', from=1-3, to=1-2]
	\arrow[from=2-3, to=2-2]
	\arrow[from=2-4, to=2-3]
\end{tikzcd}\]
is exact in $\ub{Ab}$. Therefore we obtain a canonical isomorphism of abelian groups
\begin{equation}
\ub{CMod}_R \big( R \tu{S} \otimes_R P_R , M \big) \simeq \ker M(i_{S,\infty}) .
\end{equation}
From the construction, one can check that the isomorphisms of abelian groups
\begin{equation}
\uhom_R(P_R, M)(S) \simeq
\ub{CMod}_R \big( R \tu{S} \otimes_R P_R , M \big) \simeq
\ker M(i_{S,\infty})
\end{equation}
we have just constructed are functorial in both $S \in |\ub{Prof}^{\op}|$ and $M \in |\ub{CMod}_R|$.
\end{enumerate}
\end{proof}

\begin{prop} \label{prop:characterization of coalescence}
Let $R$ be a condensed ring. Let $f,g \in R(*)$. Let $M$ be a condensed $R$-module. Then the following conditions are equivalent.
\begin{enumerate}
\item
The homomorphism
\[\begin{tikzcd}
	{\uhom_R(P_R, M)} && {\uhom_R(P_R, M)}
	\arrow["{(\Delta_{f/g})^{*}}", from=1-1, to=1-3]
\end{tikzcd}\]
is an isomorphism of condensed $R$-modules.

\item For every light profinite set $S$ and for every homomorphism $R \, \tu{S} \otimes_R P_R \xto{\phi} M$ of condensed $R$-modules, there exists a unique homomorphism $R \, \tu{S} \otimes_R P_R \xto{\psi} M$ of condensed $R$-modules such that the diagram
\[\begin{tikzcd}
	{R \, \tu{S} \otimes_R P_R} & M \\
	{R \, \tu{S} \otimes_R P_R}
	\arrow["\phi", from=1-1, to=1-2]
	\arrow["{\id_{R \tu{S}} \, \otimes_R \, \Delta_{f/g}}"', from=1-1, to=2-1]
	\arrow["\psi"', from=2-1, to=1-2]
\end{tikzcd}\]
is commutative.

\item For every light profinite set $S$, the map
\[\begin{tikzcd}
	{M(S \times \N_{\infty})} && {M(S \times \N_{\infty})}
	\arrow["{g \cdot \id - f \cdot M(\tau_S)}", from=1-1, to=1-3]
\end{tikzcd}\]
restricts to a bijection
\[\begin{tikzcd}
	{\ker M(i_{S,\infty})} && {\ker M(i_{S,\infty}).}
	\arrow["\sim", from=1-1, to=1-3]
\end{tikzcd}\]

\end{enumerate}
\end{prop}

\begin{proof}
By \cref{lem:description of uhom R (PR M) (S)}, we have canonical bijections
\begin{equation}
\uhom_R(P_R, M)(S) \simeq
\ub{CMod}_R \big( R \tu{S} \otimes_R P_R , M \big) \simeq
\ker M(i_{S,\infty})
\end{equation}
for every light profinite set $S$. Moreover, from the construction of these bijections given in the proof of \cref{lem:description of uhom R (PR M) (S)}, one can check that the diagram
\[\begin{tikzcd}
	{\uhom_R(P_R, M)(S)} & {\ub{CMod}_R \big( R \tu{S} \otimes_R P_R , M \big)} & {\ker M(i_{S,\infty})} \\
	{\uhom_R(P_R, M)(S)} & {\ub{CMod}_R \big( R \tu{S} \otimes_R P_R , M \big)} & {\ker M(i_{S,\infty})}
	\arrow["\sim", no head, from=1-1, to=1-2]
	\arrow["{\big( (\Delta_{f/g})^* \big)_S}"', from=1-1, to=2-1]
	\arrow["\sim", no head, from=1-2, to=1-3]
	\arrow["{\big( \id_{R \tu{S}} \, \otimes_R \, \Delta_{f/g} \big)^*}"', from=1-2, to=2-2]
	\arrow["{g \cdot \id - f \cdot M(\tau_S)}", from=1-3, to=2-3]
	\arrow["\sim"', no head, from=2-1, to=2-2]
	\arrow["\sim"', no head, from=2-2, to=2-3]
\end{tikzcd}\]
is commutative for every light profinite set $S$. Therefore the equivalence of (1), (2) and (3) follows.
\end{proof}

\begin{df} \label{df:definition of coalescent modules}
(cf. \cite{Kedlaya:note}, Definition 10.3.2)

Let $R$ be a condensed ring.
\begin{enumerate}
\item
Let $f,g \in R(*)$. A condensed $R$-module $M$ is said to be $f/g$-\ti{coalescent} if $M$ satisfies the equivalent conditions of \cref{prop:characterization of coalescence}. 

\item
Let $T \sub R(*) \times R(*)$. A condensed $R$-module $M$ is said to be $T$-\ti{coalescent} if $M$ is $f/g$-coalescent for every $(f,g) \in T$.

\item
We write $\ub{CMod}_{R \approx T}$ for the full subcategory of $\ub{CMod}_R$ consisting of all $T$-coalescent condensed $R$-modules.
\end{enumerate}
\end{df}

\begin{df}
Let $R$ be a condensed ring.
\begin{enumerate}
\item
Let $T \sub R(*) \times R(*)$. A condensed $R$-algebra $A$ is said to be $T$-\ti{coalescent} if the underlying condensed $R$-module of $A$ is $T$-coalescent.

\item
We write $\ub{CAlg}_{R \approx T}$ for the full subcategory of $\ub{CAlg}_R$ consisting of all $T$-coalescent condensed $R$-algebras.
\end{enumerate}
\end{df}

\begin{rem} Let $R$ be a condensed ring. Let $f,g \in R(*)$.
\begin{enumerate}
\item
If $u \in R(*)$ is an invertible element of $R(*)$, then a condensed $R$-module $M$ is $f/g$-coalescent if and only if it is $uf/ug$-coalescent.

\item (cf. \cite{Kedlaya:note}, Remark 10.3.6)

Let $R'$ be another condensed ring. Let $\phi : R \to R'$ be a homomorphism of condensed rings. If $N$ is any condensed $R'$-module, then we obtain a condensed $R$-module $N_{/R}$ from $N$ by restricting scalars via $\phi : R \to R'$. By the characterization (3) of \cref{prop:characterization of coalescence}, the condensed $R$-module $N_{/R}$ is $f/g$-coalescent if and only if the condensed $R'$-module $N$ is $\phi_*(f) / \phi_*(g)$-coalescent.
\end{enumerate}
\end{rem}

\begin{rem} Let $R$ be a condensed ring.
\begin{enumerate}
\item
The notion of $f/g$-coalescent condensed $R$-modules given in \cref{df:definition of coalescent modules} is a direct generalization of the notion of $f$-coalescent modules defined in Definition 10.3.2 of \cite{Kedlaya:note}. In fact, if $f \in R(*)$, then a condensed $R$-module $M$ is $f$-coalescent in the sense of Definition 10.3.2 of \cite{Kedlaya:note} if and only if the condensed $R$-module $M$ is $f/1$-coalescent in the sense of \cref{df:definition of coalescent modules}.

\item
A condensed $R$-module $M$ is $1/1$-coalescent if and only if the underlying condensed abelian group of $M$ is \ti{solid} in the sense of Definition 3.2.1 of \cite{Camargo:note} and Definition 6.2.1 of \cite{Kedlaya:note}.
\end{enumerate}
\end{rem}

\subsubsection{The coalescence functor}

\begin{prop} \tup{(cf. \cite{Kedlaya:note}, Proposition 10.3.4)}

Let $R$ be a condensed ring. Let $T \sub R(*) \times R(*)$. Then the inclusion functor $\ub{CMod}_{R \approx T} \mon \ub{CMod}_R$ has a left adjoint.
\end{prop}

\begin{proof}
We use Theorem 5.4.7 of \cite{Borceux:cat1} to prove this proposition. First of all, the category $\ub{CMod}_R$ is cocomplete by \cref{prop:CMod R to CAb preserves and reflects exact sequences}. Next, by \cref{cor:Prof is essentially small}, there is a set $\cat{S}$ of representatives of isomorphism classes of $\ub{Prof}$. Let $\Sigma$ be the set of all morphisms $\id_{R \tu{S}} \otimes_R \, \Delta_{f/g} : R \tu{S} \otimes_R P_R \to R \tu{S} \otimes_R P_R$ in $\ub{CMod}_R$, where $S \in \cat{S}$ and $(f,g) \in T$. Note that $\Sigma$ is a \ti{set} since both $\cat{S}$ and $T$ are sets. The characterization (2) of \cref{prop:characterization of coalescence} shows that the subcategory $\ub{CMod}_{R \approx T}$ of $\ub{CMod}_R$ is equal to the orthogonal subcategory of $\ub{CMod}_R$ determined by the set of morphisms $\Sigma$ in the sense of Definition 5.4.5 of \cite{Borceux:cat1}. Then Theorem 5.4.7 of \cite{Borceux:cat1} shows that the category $\ub{CMod}_{R \approx T}$ is a reflective subcategory of $\ub{CMod}_R$ provided that every object $M$ of $\ub{CMod}_R$ is presentable in the sense that the representable functor $\ub{CMod}_R \to \ub{Set}$, $N \mapsto \ub{CMod}_R(M,N)$ preserves $\alpha$-filtered colimits for some regular cardinal $\alpha$ (\cite{Borceux:cat2}, Definition 5.1.1).

We prove that every object of $\ub{CMod}_R$ is presentable in the sense we have just described above. First we show that $R \tu{S}$ is a presentable object of $\ub{CMod}_R$ for every $S \in \cat{S}$. If $S \in \cat{S}$, then the representable functor $\ub{CMod}_R \to \ub{Set}$, $N \mapsto \ub{CMod}_R(R \tu{S},N)$ is isomorphic to the functor $\ub{CMod}_R \to \ub{Set}$, $N \mapsto N(S)$ by adjunction and the Yoneda lemma. This functor is the composition of the functor $\mathrm{ev}_S : \ub{CMod}_R \to \ub{Mod}_{R(S)}, N \mapsto N(S)$ and the forgetful functor $\ub{Mod}_{R(S)} \to \ub{Set}$. By \cref{prop:evaluation of CMod R preserves limits filtered colimits and coproducts}, the functor $\mathrm{ev}_S : \ub{CMod}_R \to \ub{Mod}_{R(S)}, N \mapsto N(S)$ preserves filtered colimits. The forgetful functor $\ub{Mod}_{R(S)} \to \ub{Set}$ also preserves filtered colimits. Therefore we conclude that the representable functor $\ub{CMod}_R \to \ub{Set}$, $N \mapsto \ub{CMod}_R(R \tu{S},N)$ preserves filtered colimits. This shows that $R \tu{S}$ is a presentable object of $\ub{CMod}_R$. 

Next we claim that every object $M$ of $\ub{CMod}_R$ admits an epimorphism from a coproduct of some set-indexed family of objects which belong to the set $\set{R \tu{S}}{S \in \cat{S}}$. Indeed, consider the set
\begin{equation}
I := \set{(S,f)}{S \in \cat{S} \: ; \: f \in \ub{CMod}_R(R \tu{S},M)} .
\end{equation}
Note that $I$ is a \ti{set} since $\cat{S}$ is a set. We have a unique morphism $\gamma : (\bigoplus_{(S,f) \in I} \, R \tu{S}) \to M$ in $\ub{CMod}_R$ such that the diagram
\[\begin{tikzcd}
	{\displaystyle \bigoplus_{(S,f) \in I} R \tu{S}} & M \\
	{R \tu{S}}
	\arrow["\gamma", from=1-1, to=1-2]
	\arrow["\can", from=2-1, to=1-1]
	\arrow["f"', from=2-1, to=1-2]
\end{tikzcd}\]
is commutative for every $(S,f) \in I$. Using adjunction and the Yoneda lemma, it immediately follows that the morphism $\gamma : (\bigoplus_{(S,f) \in I} \, R \tu{S}) \to M$ is an epimorphism in $\ub{CMod}_R$.

Thus we have proved that every object of $\ub{CMod}_R$ admits an epimorphism from a coproduct of some set-indexed family of objects which belong to the set $\set{R \tu{S}}{S \in \cat{S}}$. Then we conclude that for every object $M$ of $\ub{CMod}_R$, there exists an exact sequence of the form
\[\begin{tikzcd}
	{\displaystyle \bigoplus_{j \in J} R \tu{S'_j}} & {\displaystyle \bigoplus_{i \in I} R \tu{S_i}} & M & 0
	\arrow[from=1-1, to=1-2]
	\arrow[from=1-2, to=1-3]
	\arrow[from=1-3, to=1-4]
\end{tikzcd}\]
where $(S_i)_{i \in I}$ and $(S'_j)_{j \in J}$ are both set-indexed families of objects belonging to $\cat{S}$. Since $R \tu{S}$ is a presentable object of $\ub{CMod}_R$ for every $S \in \cat{S}$, Proposition 5.1.4 of \cite{Borceux:cat2} shows that $M$ is a presentable object of $\ub{CMod}_R$. This completes the proof.
\end{proof}

\begin{df} (cf. \cite{Kedlaya:note}, Proposition 10.3.4)

Let $R$ be a condensed ring. Let $T \sub R(*) \times R(*)$.
\begin{enumerate}
\item
The left adjoint of the inclusion functor $\ub{CMod}_{R \approx T} \mon \ub{CMod}_R$ is called the $T$-\ti{coalscence}. It is denoted by $\ub{CMod}_R \to \ub{CMod}_{R \approx T}$, $M \mapsto M_{\approx T}$.

\item
If $M$ is a condensed $R$-module, then the canonical homomorphism $M \to M_{\approx T}$ is called the $T$-\ti{coalescence of} $M$.
\end{enumerate}
\end{df}

\begin{prop} \label{prop:coalescence and tensor and internal hom}
Let $R$ be a condensed ring. Let $T \sub R(*) \times R(*)$.
\begin{enumerate}
\item \tup{(cf. \cite{Kedlaya:note}, Proposition 10.3.3)}
Let $M,N$ be condensed $R$-modules. If $N$ is $T$-coalescent, then $\uhom_R (M,N)$ is $T$-coalescent.

\item
Let $M \xto{\iota_1} M_{\approx T}$, $N \xto{\iota_2} N_{\approx T}$, $M \otimes_R N \xto{\iota_3} (M \otimes_R N)_{\approx T}$ and $M_{\approx T} \otimes_R N_{\approx T} \xto{\iota_4} (M_{\approx T} \otimes_R N_{\approx T})_{\approx T}$ be the $T$-coalescence of $M,N,M \otimes_R N$ and $M_{\approx T} \otimes_R N_{\approx T}$ respectively. Then there exists a unique homomorphism $(M \otimes_R N)_{\approx T} \xto{\theta} (M_{\approx T} \otimes_R N_{\approx T})_{\approx T}$ of condensed $R$-modules such that the diagram 
\[\begin{tikzcd}
	{M \otimes_R N} & {M_{\approx T} \otimes_R N_{\approx T}} \\
	{(M \otimes_R N)_{\approx T}} & {(M_{\approx T} \otimes_R N_{\approx T})_{\approx T}}
	\arrow["{\iota_1 \, \otimes_R \, \iota_2}", from=1-1, to=1-2]
	\arrow["{\iota_3}"', from=1-1, to=2-1]
	\arrow["{\iota_4}", from=1-2, to=2-2]
	\arrow["\theta"', from=2-1, to=2-2]
\end{tikzcd}\]
is commutative. Moreover, the homomorphism $(M \otimes_R N)_{\approx T} \xto{\theta} (M_{\approx T} \otimes_R N_{\approx T})_{\approx T}$ is an isomorphism of condensed $R$-modules.
\end{enumerate}
\end{prop}

\begin{proof}~
\begin{enumerate}
\item
We use the characterization (1) of \cref{prop:characterization of coalescence}. By adjunction, we have a canonical isomorphism of condensed $R$-modules
\begin{align}
\uhom_R \big( P_R, \uhom_R (M,N) \big)
& \simeq \uhom_R \big( P_R \, \otimes_R \, M ,\, N \big) \\
& \simeq \uhom_R \big( M \, \otimes_R \, P_R ,\, N \big) \\
& \simeq \uhom_R \big( M, \uhom_R (P_R,N) \big)
\end{align}
and the diagram
\[\begin{tikzcd}
	{\uhom_R \big( P_R, \uhom_R (M,N) \big)} & {\uhom_R \big( M, \uhom_R (P_R,N) \big)} \\
	{\uhom_R \big( P_R, \uhom_R (M,N) \big)} & {\uhom_R \big( M, \uhom_R (P_R,N) \big)}
	\arrow["\sim", no head, from=1-1, to=1-2]
	\arrow["{\uhom_R \big( \Delta_{f/g}, \, \id_{\uhom_R (M,N)} \big)}"', from=1-1, to=2-1]
	\arrow["{\uhom_R \big( \id_M, \, \uhom_R (\Delta_{f/g},\id_N) \big)}", from=1-2, to=2-2]
	\arrow["\sim"', no head, from=2-1, to=2-2]
\end{tikzcd}\]
is commutative for every $f,g \in R(*)$. Using the characterization (1) of \cref{prop:characterization of coalescence}, we conclude that if the condensed $R$-module $N$ is $T$-coalescent, then the condensed $R$-module $\uhom_R (M,N)$ is $T$-coalescent.

\item
We have a homomorphism $M \otimes_R N \xto{\iota_1 \, \otimes_R \, \iota_2} M_{\approx T} \otimes_R N_{\approx T} \xto{\iota_4} (M_{\approx T} \otimes_R N_{\approx T})_{\approx T}$ of condensed $R$-modules and the condensed $R$-module $(M_{\approx T} \otimes_R N_{\approx T})_{\approx T}$ is $T$-coalescent. By the universality of the $T$-coalescence $M \otimes_R N \xto{\iota_3} (M \otimes_R N)_{\approx T}$ of $M \otimes_R N$, there exists a unique homomorphism $(M \otimes_R N)_{\approx T} \xto{\theta} (M_{\approx T} \otimes_R N_{\approx T})_{\approx T}$ of condensed $R$-modules such that the diagram 
\[\begin{tikzcd}
	{M \otimes_R N} & {M_{\approx T} \otimes_R N_{\approx T}} \\
	{(M \otimes_R N)_{\approx T}} & {(M_{\approx T} \otimes_R N_{\approx T})_{\approx T}}
	\arrow["{\iota_1 \, \otimes_R \, \iota_2}", from=1-1, to=1-2]
	\arrow["{\iota_3}"', from=1-1, to=2-1]
	\arrow["{\iota_4}", from=1-2, to=2-2]
	\arrow["\theta"', from=2-1, to=2-2]
\end{tikzcd}\]
is commutative. We show that this homomorphism $(M \otimes_R N)_{\approx T} \xto{\theta} (M_{\approx T} \otimes_R N_{\approx T})_{\approx T}$ is an isomorphism of condensed $R$-modules.

Let $Q$ be an arbitrary $T$-coalescent condensed $R$-module. By the universality of the $T$-coalescence $M \otimes_R N \xto{\iota_3} (M \otimes_R N)_{\approx T}$ and $M_{\approx T} \otimes_R N_{\approx T} \xto{\iota_4} (M_{\approx T} \otimes_R N_{\approx T})_{\approx T}$, the maps
\begin{align}
\ub{CMod}_{R \approx T} \big( (M \otimes_R N)_{\approx T}, Q \big)
& = \ub{CMod}_R \big( (M \otimes_R N)_{\approx T}, Q \big) \\
& \xto{(\iota_3)^*} \ub{CMod}_R \big( M \otimes_R N, Q \big) \: ; \\
\ub{CMod}_{R \approx T} \big( (M_{\approx T} \otimes_R N_{\approx T})_{\approx T}, Q \big)
& = \ub{CMod}_R \big( (M_{\approx T} \otimes_R N_{\approx T})_{\approx T}, Q \big) \\
& \xto{(\iota_4)^*} \ub{CMod}_R \big( M_{\approx T} \otimes_R N_{\approx T}, Q \big)
\end{align}
are bijective. By adjunction, we have a canonical bijection
\begin{equation}
\ub{CMod}_R \big( M_{\approx T} \otimes_R N_{\approx T}, Q \big)
\simeq \ub{CMod}_R \big( M_{\approx T}, \uhom_R( N_{\approx T}, Q ) \big) .
\end{equation}
By (1), the condensed $R$-module $\uhom_R( N_{\approx T}, Q )$ is $T$-coalescent. Then the universality of the $T$-coalescence $M \xto{\iota_1} M_{\approx T}$ of $M$ shows that the map
\begin{equation}
\ub{CMod}_R \big( M_{\approx T}, \uhom_R( N_{\approx T}, Q ) \big)
\xto{(\iota_1)^*} \ub{CMod}_R \big( M, \uhom_R( N_{\approx T}, Q ) \big)
\end{equation}
is bijective. By adjunction, we have a canonical bijection
\begin{align}
\ub{CMod}_R \big( M, \uhom_R( N_{\approx T}, Q ) \big)
& \simeq \ub{CMod}_R \big( M \otimes_R N_{\approx T}, Q \big) \\
& \simeq \ub{CMod}_R \big( N_{\approx T} \otimes_R M, Q \big) \\
& \simeq \ub{CMod}_R \big( N_{\approx T}, \uhom_R(M,Q) \big) .
\end{align}
By (1), the condensed $R$-module $\uhom_R(M,Q)$ is $T$-coalescent. Then the universality of the $T$-coalescence $N \xto{\iota_2} N_{\approx T}$ of $N$ shows the map
\begin{equation}
\ub{CMod}_R \big( N_{\approx T}, \uhom_R(M,Q) \big)
\xto{(\iota_2)^*} \ub{CMod}_R \big( N, \uhom_R(M,Q) \big)
\end{equation}
is bijective. By adjunction, we have canonical bijections
\begin{align}
\ub{CMod}_R \big( N, \uhom_R(M,Q) \big) \simeq
\ub{CMod}_R \big( N \otimes_R M, Q \big) \simeq
\ub{CMod}_R \big( M \otimes_R N, Q \big) .
\end{align}
Consequently, we obtain a bijection
\begin{equation}
\ub{CMod}_{R \approx T} \big( (M \otimes_R N)_{\approx T}, Q \big) \simeq
\ub{CMod}_{R \approx T} \big( (M_{\approx T} \otimes_R N_{\approx T})_{\approx T}, Q \big) ,
\end{equation}
and the construction shows that this map is equal to the map
\begin{equation}
\ub{CMod}_{R \approx T} \big( (M \otimes_R N)_{\approx T}, Q \big) \xto{\theta^*}
\ub{CMod}_{R \approx T} \big( (M_{\approx T} \otimes_R N_{\approx T})_{\approx T}, Q \big).
\end{equation}
Since $Q$ is an arbitrary $T$-coalescent condensed $R$-module, the Yoneda lemma implies that the homomorphism $(M \otimes_R N)_{\approx T} \xto{\theta} (M_{\approx T} \otimes_R N_{\approx T})_{\approx T}$ is an isomorphism in $\ub{CMod}_{R \approx T}$. Therefore $(M \otimes_R N)_{\approx T} \xto{\theta} (M_{\approx T} \otimes_R N_{\approx T})_{\approx T}$ is an isomorphism in $\ub{CMod}_R$.
\end{enumerate}
\end{proof}

\begin{cor} \label{cor:extension of R bilinear maps along coalescence}
Let $R$ be a condensed ring. Let $T \sub R(*) \times R(*)$. Let $M,N$ be condensed $R$-modules. Let $Q$ be a $T$-coalescent $R$-module. If $M \times N \xto{\beta} Q$ is any $R$-bilinear map, then there exists a unique $R$-bilinear map $M_{\approx T} \times N_{\approx T} \xto{\gamma} Q$ such that the diagram
\[\begin{tikzcd}
	{M \times N} & Q \\
	{M_{\approx T} \times N_{\approx T}}
	\arrow["\beta", from=1-1, to=1-2]
	\arrow["{\iota_1 \times \iota_2}"', from=1-1, to=2-1]
	\arrow["\gamma"', from=2-1, to=1-2]
\end{tikzcd}\]
is commutative, where $M \xto{\iota_1} M_{\approx T}$ and $N \xto{\iota_2} N_{\approx T}$ are the $T$-coalescence of $M,N$ respectively.
\end{cor}

\begin{proof}
For each condensed $R$-modules $U,V$ and $W$, let us write $\mathrm{Bil}_R(U,V;W)$ for the set of $R$-bilinear maps $U \times V \to W$. By \cref{prop:tensor product over R represents R bilinear maps}, the maps
\[\begin{tikzcd}[row sep=tiny]
	{\mathrm{Bil}_R(M,N;Q)} & {\ub{CMod}_R(M \otimes_R N,Q) \: ;} \\
	{\mathrm{Bil}_R(M_{\approx T}, N_{\approx T};Q)} & {\ub{CMod}_R(M_{\approx T} \otimes_R N_{\approx T},Q)}
	\arrow["{(\otimes_R)^*}"', from=1-2, to=1-1]
	\arrow["{(\otimes_R)^*}"', from=2-2, to=2-1]
\end{tikzcd}\]
are bijective and the diagram
\[\begin{tikzcd}
	{M \times N} & {M \otimes_R N} \\
	{M_{\approx T} \times N_{\approx T}} & {M_{\approx T} \otimes_R N_{\approx T}}
	\arrow["{\otimes_R}", from=1-1, to=1-2]
	\arrow["{\iota_1 \times \iota_2}"', from=1-1, to=2-1]
	\arrow["{\iota_1 \, \otimes_R \, \iota_2}", from=1-2, to=2-2]
	\arrow["{\otimes_R}"', from=2-1, to=2-2]
\end{tikzcd}\]
is commutative. On the other hand, let $M \otimes_R N \xto{\iota_3} (M \otimes_R N)_{\approx T}$ and $M_{\approx T} \otimes_R N_{\approx T} \xto{\iota_4} (M_{\approx T} \otimes_R N_{\approx T})_{\approx T}$ be the $T$-coalescence of $M \otimes_R N$ and $M_{\approx T} \otimes_R N_{\approx T}$ respectively. Since $Q$ is a $T$-coalescent $R$-module, the universality of $T$-coalescence shows that the maps
\[\begin{tikzcd}[row sep=tiny]
	{\ub{CMod}_R(M \otimes_R N,Q)} & {\ub{CMod}_R \big( (M \otimes_R N)_{\approx T},Q \big) \: ;} \\
	{\ub{CMod}_R(M_{\approx T} \otimes_R N_{\approx T},Q)} & {\ub{CMod}_R \big( (M_{\approx T} \otimes_R N_{\approx T})_{\approx T},Q \big)}
	\arrow["{(\iota_3)^*}"', from=1-2, to=1-1]
	\arrow["{(\iota_4)^*}"', from=2-2, to=2-1]
\end{tikzcd}\]
are bijective. Moreover, \cref{prop:coalescence and tensor and internal hom} shows that there exists an isomorphism $(M \otimes_R N)_{\approx T} \xto{\theta} (M_{\approx T} \otimes_R N_{\approx T})_{\approx T}$ of condensed $R$-modules such that the diagram 
\[\begin{tikzcd}
	{M \otimes_R N} & {(M \otimes_R N)_{\approx T}} \\
	{M_{\approx T} \otimes_R N_{\approx T}} & {(M_{\approx T} \otimes_R N_{\approx T})_{\approx T}}
	\arrow["{\iota_3}", from=1-1, to=1-2]
	\arrow["{\iota_1 \, \otimes_R \, \iota_2}"', from=1-1, to=2-1]
	\arrow["\theta", from=1-2, to=2-2]
	\arrow["{\iota_4}"', from=2-1, to=2-2]
\end{tikzcd}\]
is commutative. Consequently, the diagram
\[\begin{tikzcd}
	{\mathrm{Bil}_R(M,N;Q)} & {\ub{CMod}_R(M \otimes_R N,Q)} & {\ub{CMod}_R \big( (M \otimes_R N)_{\approx T},Q \big)} \\
	{\mathrm{Bil}_R(M_{\approx T}, N_{\approx T};Q)} & {\ub{CMod}_R(M_{\approx T} \otimes_R N_{\approx T},Q)} & {\ub{CMod}_R \big( (M_{\approx T} \otimes_R N_{\approx T})_{\approx T},Q \big)}
	\arrow["{(\otimes_R)^*}"', from=1-2, to=1-1]
	\arrow["{(\iota_3)^*}"', from=1-3, to=1-2]
	\arrow["{(\iota_1 \times \iota_2)^*}", from=2-1, to=1-1]
	\arrow["{(\iota_1 \, \otimes_R \, \iota_2)^*}", from=2-2, to=1-2]
	\arrow["{(\otimes_R)^*}", from=2-2, to=2-1]
	\arrow["{\theta^*}"', from=2-3, to=1-3]
	\arrow["{(\iota_4)^*}", from=2-3, to=2-2]
\end{tikzcd}\]
is commutative, and the right vertical map and the horizontal maps are bijective. Therefore the left vertical map
\[\begin{tikzcd}
	{\mathrm{Bil}_R(M_{\approx T}, N_{\approx T};Q)} && {\mathrm{Bil}_R(M,N;Q)}
	\arrow["{(\iota_1 \times \iota_2)^*}", from=1-1, to=1-3]
\end{tikzcd}\]
is also bijective.
\end{proof}

\begin{prop} \label{prop:R algebra structure on the coalescence of an R algebra}
Let $R$ be a condensed ring. Let $T \sub R(*) \times R(*)$. Let $A$ be a condensed $R$-algebra. Let $A \xto{\iota} A_{\approx T}$ be the $T$-coalesence of the underlying condensed $R$-module of $A$.
\begin{enumerate}
\item
The condensed $R$-module $A_{\approx T}$ has a unique structure of a condensed $R$-algebra such that the $T$-colescence $A \xto{\iota} A_{\approx T}$ is a homomorphism of condensed $R$-algebras.

\item
Let us give $A_{\approx T}$ the structure of a condensed $R$-algebra defined in (1). For every $T$-coalescent condensed $R$-algebra $B$ and every homomorphism $A \xto{\phi} B$ of condensed $R$-algebras, there exists a unique homomorphism $A_{\approx T} \xto{\psi} B$ of condensed $R$-algebras such that the diagram
\[\begin{tikzcd}
	A & B \\
	{A_{\approx T}}
	\arrow["\phi", from=1-1, to=1-2]
	\arrow["\iota"', from=1-1, to=2-1]
	\arrow["\psi"', from=2-1, to=1-2]
\end{tikzcd}\]
is commutative.
\end{enumerate}
\end{prop}

If $A$ is a condensed $R$-algebra, we always consider $A_{\approx T}$ as a condensed $R$-algebra by giving it the unique structure of a condensed $R$-algebra such that the $T$-colescence $A \to A_{\approx T}$ is a homomorphism of condensed $R$-algebras.

\begin{proof}[Proof of \cref{prop:R algebra structure on the coalescence of an R algebra}]
Let $\mu : A \times A \to A$ be the multiplication map of $A$.
\begin{enumerate}
\item
Since the map $A \times A \xto{\mu} A \xto{\iota} A_{\approx T}$ is $R$-bilinear and the condensed $R$-module $A_{\approx T}$ is $T$-coalescent, \cref{cor:extension of R bilinear maps along coalescence} shows that there exists a unique $R$-bilinear map $\nu : A_{\approx T} \times A_{\approx T} \to A_{\approx T}$ such that the diagram \[\begin{tikzcd}
	{A \times A} & A \\
	{A_{\approx T} \times A_{\approx T}} & {A_{\approx T}}
	\arrow["\mu", from=1-1, to=1-2]
	\arrow["{\iota \times \iota}"', from=1-1, to=2-1]
	\arrow["\iota", from=1-2, to=2-2]
	\arrow["\nu"', from=2-1, to=2-2]
\end{tikzcd}\]
is commutative. Using the uniqueness assertion of \cref{cor:extension of R bilinear maps along coalescence}, one checks that this $R$-bilinear map $\nu : A_{\approx T} \times A_{\approx T} \to A_{\approx T}$ defines a structure of a condensed $R$-algebra on the condensed $R$-module $A_{\approx T}$ such that the $R$-bilinear map $\nu : A_{\approx T} \times A_{\approx T} \to A_{\approx T}$ is the multiplication map and the map $\iota : A \to A_{\approx T}$ is a homomorphism of condensed $R$-algebras.

Suppose that there is a another structure of a condensed $R$-algebra on the condensed $R$-module $A_{\approx T}$ such that $A \xto{\iota} A_{\approx T}$ is a homomorphism of condensed $R$-algebras. Let $\nu' : A_{\approx T} \times A_{\approx T} \to A_{\approx T}$ be the multiplication map with respect to this structure. Then the diagram
\[\begin{tikzcd}
	{A \times A} & A \\
	{A_{\approx T} \times A_{\approx T}} & {A_{\approx T}}
	\arrow["\mu", from=1-1, to=1-2]
	\arrow["{\iota \times \iota}"', from=1-1, to=2-1]
	\arrow["\iota", from=1-2, to=2-2]
	\arrow["{\nu'}"', from=2-1, to=2-2]
\end{tikzcd}\]
is commutative. By the definition of $\nu$, we have $\nu' = \nu$. This proves the uniqueness of a structure of a condensed $R$-algebra on the condensed $R$-module $A_{\approx T}$ such that $A \xto{\iota} A_{\approx T}$ is a homomorphism of condensed $R$-algebras.

\item
By the universality of the $T$-coalescence $A \xto{\iota} A_{\approx T}$ of the condensed $R$-module $A$, there exists a unique homomorphism $A_{\approx T} \xto{\psi} B$ of condensed $R$-modules such that the diagram
\[\begin{tikzcd}
	A & B \\
	{A_{\approx T}}
	\arrow["\phi", from=1-1, to=1-2]
	\arrow["\iota"', from=1-1, to=2-1]
	\arrow["\psi"', from=2-1, to=1-2]
\end{tikzcd}\]
is commutative. We prove that this homomorphism $A_{\approx T} \xto{\psi} B$ is in fact a homomorphism of condensed $R$-algebras. For every light profinite set $S$, we have $\psi_S(1) = \psi_S (\iota_S(1) ) = \phi_S(1) = 1$. If we write $\lambda : B \times B \to B$ for the multipliction map of $B$, then the following diagram is commutative.
\[\begin{tikzcd}
	{A_{\approx T} \times A_{\approx T}} & {B \times B} \\
	{A \times A} & A & B \\
	{A_{\approx T} \times A_{\approx T}} & {A_{\approx T}}
	\arrow["{\psi \times \psi}", from=1-1, to=1-2]
	\arrow["\lambda", from=1-2, to=2-3]
	\arrow["{\iota \times \iota}", from=2-1, to=1-1]
	\arrow["{\phi \times \phi}"', from=2-1, to=1-2]
	\arrow["\mu"', from=2-1, to=2-2]
	\arrow["{\iota \times \iota}"', from=2-1, to=3-1]
	\arrow["\phi"', from=2-2, to=2-3]
	\arrow["\iota"', from=2-2, to=3-2]
	\arrow["\nu"', from=3-1, to=3-2]
	\arrow["\psi"', from=3-2, to=2-3]
\end{tikzcd}\]
Since both the maps $A_{\approx T} \times A_{\approx T} \xto{\psi \times \psi} B \times B \xto{\lambda} B$ and $A_{\approx T} \times A_{\approx T} \xto{\nu} A_{\approx T} \xto{\psi} B$ are $R$-bilinear and the underlying condensed $R$-module of $B$ is $T$-coalescent, \cref{cor:extension of R bilinear maps along coalescence} shows that the following diagram is commutative.
\[\begin{tikzcd}
	{A_{\approx T} \times A_{\approx T}} & {B \times B} \\
	&& B \\
	{A_{\approx T} \times A_{\approx T}} & {A_{\approx T}}
	\arrow["{\psi \times \psi}", from=1-1, to=1-2]
	\arrow[equals, from=1-1, to=3-1]
	\arrow["\lambda", from=1-2, to=2-3]
	\arrow["\nu"', from=3-1, to=3-2]
	\arrow["\psi"', from=3-2, to=2-3]
\end{tikzcd}\]
Therefore $\psi_S( a \cdot b) = \psi_S(a) \cdot \psi_S(b)$ for every light profinite set $S$ and every elements $a,b \in A_{\approx T}(S)$. This completes the proof that $A_{\approx T} \xto{\psi} B$ is  a homomorphism of condensed $R$-algebras.
\end{enumerate}
\end{proof}

\begin{cor}
Let $R$ be a condensed ring. Let $T \sub R(*) \times R(*)$. Then the inclusion functor $\ub{CAlg}_{R \approx T} \mon \ub{CAlg}_R$ has a left adjoint.
\end{cor}

\begin{df}
Let $R$ be a condensed ring. Let $T \sub R(*) \times R(*)$.
\begin{enumerate}
\item
The left adjoint of the inclusion functor $\ub{CAlg}_{R \approx T} \mon \ub{CAlg}_R$ is called the $T$-\ti{coalescence}. It is denoted by $\ub{CAlg}_R \to \ub{CAlg}_{R \approx T}$, $A \mapsto A_{\approx T}$.

\item
If $A$ is a condensed $R$-algebra, then the canonical homomorphism $A \to A_{\approx T}$ is called the $T$-\ti{coalescence of} $A$. By \cref{prop:R algebra structure on the coalescence of an R algebra}, it is equal to the $T$-coalescence of the underlying condensed $R$-module of $A$.
\end{enumerate}
\end{df}

\subsubsection{Internal projectivity of the sequence space}

\begin{prop} \label{prop:P R is internally projective} \tup{(cf. \cite{Kedlaya:note} , Definition 10.3.2)}

Let $R$ be a condensed ring. The the functor
\[\begin{tikzcd}
	{\ub{CMod}_R} & {\ub{CMod}_R,} & M & {\uhom_R(P_R,M)}
	\arrow[from=1-1, to=1-2]
	\arrow[maps to, from=1-3, to=1-4]
\end{tikzcd}\]
preserves limits and colimits.
\end{prop}

\begin{proof}
Let us write $H : \ub{CMod}_R \to \ub{CMod}_R$ for the functor $M \mapsto \uhom_R(P_R,M)$. Since it is right adjoint to the functor $\ub{CMod}_R \to \ub{CMod}_R$, $N \mapsto N \otimes_R P_R$, the functor $H : \ub{CMod}_R \to \ub{CMod}_R$ preserves limits.

Next we show that the functor $H : \ub{CMod}_R \to \ub{CMod}_R$ preserves colimits. By the dual of Proposition 2.9.2 of \cite{Borceux:cat1}, it suffices to show that the functor $H : \ub{CMod}_R \to \ub{CMod}_R$ preserves coproducts and cokernels.

First we show that the functor $H : \ub{CMod}_R \to \ub{CMod}_R$ preserves coproducts. Let $(M_{\lambda})_{\lambda \in \Lambda}$ be a family of condensed $R$-modules. Let $M := \bigoplus_{\lambda \in \Lambda} M_{\lambda}$. By \cref{lem:description of uhom R (PR M) (S)}, for each light profinite set $S$, we have canonical isomorphisms of abelian groups
\begin{align}
\uhom_R(P_R, M)(S) & \simeq \ker M(i_{S,\infty}) \: ; \\
\uhom_R(P_R, M_{\lambda})(S) & \simeq \ker M_{\lambda}(i_{S,\infty}) ,
\end{align}
which are functorial in $S \in |\ub{Prof}^{\op}|$ and compatible with the canonical morphisms $M_{\lambda} \to M$ $(\lambda \in \Lambda)$. From the construction of these isomorphisms given in the proof of \cref{lem:description of uhom R (PR M) (S)}, one can check that these isomorphisms are in fact isomorphisms of $R(S)$-modules. On the other hand, \cref{prop:evaluation of CMod R preserves limits filtered colimits and coproducts} shows that we have canonical isomorphisms of $R(S)$-modules
\begin{align}
\ker M(i_{S,\infty})
& \simeq \ker \left( \bigoplus_{\lambda \in \Lambda} M_{\lambda}(i_{S,\infty}) \, : \,
\bigoplus_{\lambda \in \Lambda} M_{\lambda}(S \times \N_{\infty}) \to
\bigoplus_{\lambda \in \Lambda} M_{\lambda}(S) \right) \\
& \simeq \bigoplus_{\lambda \in \Lambda} \ker \Big( M_{\lambda}(i_{S,\infty}) \Big) \: ; \\
\left( \bigoplus_{\lambda \in \Lambda} \uhom_R(P_R, M_{\lambda}) \right)(S)
& \simeq \bigoplus_{\lambda \in \Lambda} \Big( \uhom_R(P_R, M_{\lambda}) (S) \Big).
\end{align}
One can check that these isomorphisms are functorial in $S \in |\ub{Prof}^{\op}|$ and compatible with the canonical morphisms $M_{\lambda} \to M$ $(\lambda \in \Lambda)$. Consequently, we obtain an isomorphism of $R(S)$-modules
\begin{equation}
\uhom_R(P_R, M)(S) \simeq \left( \bigoplus_{\lambda \in \Lambda} \uhom_R(P_R, M_{\lambda}) \right)(S) ,
\end{equation}
which is functorial in $S \in |\ub{Prof}^{\op}|$ and compatible with the canonical morphisms $M_{\lambda} \to M$ $(\lambda \in \Lambda)$. This proves that
\begin{equation}
\uhom_R(P_R, M) \simeq \left( \bigoplus_{\lambda \in \Lambda} \uhom_R(P_R, M_{\lambda}) \right) .
\end{equation}
Therefore the functor $H : \ub{CMod}_R \to \ub{CMod}_R$ preserves coproducts.

Next we show that the functor $H : \ub{CMod}_R \to \ub{CMod}_R$ preserves cokernels. Let $P$ be the cokernel of the homomorphism $\Z \, \tu{i_{\infty}} : \Z \, \tu{*} \to \Z \, \tu{\N_{\infty}}$ in $\ub{CAb}$. Note that $P_R \simeq R \otimes P$ by definition. The functor
\[\begin{tikzcd}
	{\ub{CAb}} & {\ub{CAb},} & M & {\uhom(P,M)}
	\arrow[from=1-1, to=1-2]
	\arrow[maps to, from=1-3, to=1-4]
\end{tikzcd}\]
preserves cokernels (\cite{Camargo:note}, Theorem 2.3.3, proof of (3) ; \cite{Kedlaya:note}, Proposition 5.4.8). On the other hand, by adjunction, we have canonical isomorphisms of condensed abelian groups
\begin{equation}
\uhom(P,M) \simeq \uhom_R(R \otimes P,M) \simeq \uhom_R(P_R,M)
\end{equation}
which are functorial in $M \in |\ub{CMod}_R|$. Moreover, the forgetful functor $\ub{CMod}_R \to \ub{CAb}$ preserves and reflects cokernels by \cref{prop:CMod R to CAb preserves and reflects exact sequences}. Therefore we conclude that the functor $H : \ub{CMod}_R \to \ub{CMod}_R$, $M \mapsto \uhom_R(P_R,M)$ preserves cokernels.
\end{proof}

\begin{cor} \label{cor:coalescent modules are closed under limits and colimits}
Let $R$ be a condensed ring. Let $T \sub R(*) \times R(*)$. 
\begin{enumerate}
\item
$\ub{CMod}_{R \approx T}$ is closed under limits and colimits in $\ub{CMod}_R$.

\item
$\ub{CAlg}_{R \approx T}$ is closed under limits and filtered colimits in $\ub{CAlg}_R$.
\end{enumerate}
\end{cor}

\begin{proof}~
\begin{enumerate}
\item
This follows immediately from \cref{prop:P R is internally projective} and the characterization (1) of \cref{prop:characterization of coalescence}.

\item
This follows from (1) and the fact that the forgetful functor $\ub{CAlg}_R \to \ub{CMod}_R$ preserves limits and filtered colimits.
\end{enumerate}
\end{proof}

\begin{cor} \label{cor:CModR approx T is Grothendieck Abelian}
Let $R$ be a condensed ring. Let $T \sub R(*) \times R(*)$. Then the category $\ub{CMod}_{R \approx T}$ is a Grothendieck Abelian category. The inclusion functor $\ub{CMod}_{R \approx T} \mon \ub{CMod}_R$ preserves limits and colimits. 
\end{cor}

\begin{proof}
By \cref{prop:CMod R to CAb preserves and reflects exact sequences}, the category $\ub{CMod}_R$ is a bicomplete Abelian category in which filtered colimits are exact. Since the subcategory $\ub{CMod}_{R \approx T}$ is closed under limits and colimits in $\ub{CMod}_R$ by \cref{cor:coalescent modules are closed under limits and colimits}, the category $\ub{CMod}_{R \approx T}$ is also a bicomplete Abelian category in which filtered colimits are exact. In addition, the inclusion functor $\ub{CMod}_{R \approx T} \mon \ub{CMod}_R$ preserves limits and colimits. Furthermore, the category $\ub{CMod}_R$ has a generator $G$ by \cref{prop:CMod R to CAb preserves and reflects exact sequences}. By adjunction, the $T$-coalescence $G_{\approx T}$ of $G$ is a generator of the category $\ub{CMod}_{R \approx T}$.
\end{proof}

\subsubsection{Coalescence and other constructions}

\begin{prop} \label{prop:localization and monoid algebra of coalescent algebras}
Let $R$ be a condensed ring. Let $T \sub R(*) \times R(*)$.
\begin{enumerate}
\item
If $M$ is a $T$-coalescent condensed $R$-module and if $U$ is a multiplicatively closed subset of $R(*)$, then $U^{-1} M$ is a $T$-coalescent condensed $R$-module.

\item
If $A$ is a $T$-coalescent condensed $R$-algebra and $M$ is a commutative monoid, then $A[M]$ is a $T$-coalescent condensed $R$-algebra.
\end{enumerate}
\end{prop}

\begin{proof}~
\begin{enumerate}
\item
We use the characterization (3) of \cref{prop:characterization of coalescence}. Let $M$ be a $T$-coalescent condensed $R$-module. Let $U$ be a multiplicatively closed subset of $R(*)$. For every $(f,g) \in T$ and every light profinite set $S$, the homomorphism of $R(*)$-modules
\[\begin{tikzcd}
	{M(S \times \N_{\infty})} && {M(S \times \N_{\infty})}
	\arrow["{g \cdot \id - f \cdot M(\tau_S)}", from=1-1, to=1-3]
\end{tikzcd}\]
restricts to an isomorphism of $R(*)$-modules
\[\begin{tikzcd}
	{\ker M(i_{S,\infty})} && {\ker M(i_{S,\infty}) .}
	\arrow["{\delta}", from=1-1, to=1-3]
\end{tikzcd}\]
Since $U^{-1}R(*)$ is flat over $R(*)$, we have
\begin{equation}
\ker \big( U^{-1}M(i_{S,\infty}) \big)
= U^{-1} \big( \ker M(i_{S,\infty}) \big) ,
\end{equation}
and the homomorphism of $R(*)$-modules
\[\begin{tikzcd}[column sep=large]
	{U^{-1}M(S \times \N_{\infty})} && {U^{-1}M(S \times \N_{\infty})}
	\arrow["{g \cdot \id - f \cdot (U^{-1}M(\tau_S))}", from=1-1, to=1-3]
\end{tikzcd}\]
restricts to the homomorphism of $R(*)$-modules
\[\begin{tikzcd}
	{\ker \big( U^{-1}M(i_{S,\infty}) \big) = U^{-1} \big( \ker M(i_{S,\infty}) \big)} & {U^{-1} \big( \ker M(i_{S,\infty}) \big) = \ker \big( U^{-1}M(i_{S,\infty}) \big) .}
	\arrow["{U^{-1}\delta}", from=1-1, to=1-2]
\end{tikzcd}\]
This is an isomorphism of $R(*)$-modules since the homomorphism $\delta$ is an isomorphism of $R(*)$-modules. This shows that $U^{-1} M$ is a $T$-coalescent condensed $R$-module.

\item
This immediately follows from \cref{cor:coalescent modules are closed under limits and colimits} and the fact that the underlying condensed $R$-module of $A[M]$ is equal to $\bigoplus_{m \in M} A$.
\end{enumerate}
\end{proof}

\begin{prop} \label{prop:localization and coalescence commute}
Let $R$ be a condensed ring. Let $T \sub R(*) \times R(*)$. Let $U \sub R(*)$ be a multiplicatively closed subset of $R(*)$. Let $A$ be a condensed $R$-algebra. Let $R \xto{\sigma} A$ be the structure homomorphism of $A$. Let $A \xto{\iota_1} A_{\approx T}$ and $U^{-1}A \xto{\iota_2} (U^{-1}A)_{\approx T}$ be the $T$-coalescence of the condensed $R$-algebras $A$ and $U^{-1}A$ respectively. Let $A_{\approx T} \xto{\lambda_1} U^{-1}(A_{\approx T})$ and $A \xto{\lambda_2} U^{-1}A$ be the canonical homomorphism into the localization. Then the compositions $\lambda_1 \of \iota_1 : A \to U^{-1}(A_{\approx T})$ and $\iota_2 \of \lambda_2 : A \to (U^{-1}A)_{\approx T}$ have the following universal property.

\begin{enumerate}
\item
The condensed $R$-algebras $U^{-1}(A_{\approx T})$ and $(U^{-1}A)_{\approx T}$ are $T$-coalescent. For every $f \in U$, the element $(\lambda_1)_* (\iota_1)_* \sigma_* (f) \in U^{-1}(A_{\approx T})(*)$ is invertible in $U^{-1}(A_{\approx T})(*)$ and the element $(\iota_2)_* (\lambda_2)_* \sigma_* (f) \in (U^{-1}A)_{\approx T}(*)$ is invertible in $(U^{-1}A)_{\approx T}(*)$.

\item
Let $B$ be any $T$-coalescent condensed $R$-algebra. Let $\phi : A \to B$ be any homomorphism of condensed $R$-algebras. Suppose that the element $\phi_* \sigma_* (f) \in B(*)$ is invertible in $B(*)$ for every $f \in U$. Then there exist unique homomorphisms of condensed $R$-algebras $\psi_1 : U^{-1}(A_{\approx T}) \to B$ and $\psi_2 : (U^{-1}A)_{\approx T} \to B$ such that the diagrams
\[\begin{tikzcd}
	A & B & A & B \\
	{U^{-1}(A_{\approx T})} && {(U^{-1}A)_{\approx T}}
	\arrow["\phi", from=1-1, to=1-2]
	\arrow["{\lambda_1 \of \iota_1}"', from=1-1, to=2-1]
	\arrow["\phi", from=1-3, to=1-4]
	\arrow["{\iota_2 \of \lambda_2}"', from=1-3, to=2-3]
	\arrow["{\psi_1}"', from=2-1, to=1-2]
	\arrow["{\psi_2}"', from=2-3, to=1-4]
\end{tikzcd}\]
are commutative.
\end{enumerate}
In particular, there exists a unique isomorphism $\theta : U^{-1}(A_{\approx T}) \to (U^{-1}A)_{\approx T}$ of condensed $R$-algebras such that the diagram
\[\begin{tikzcd}
	A & {(U^{-1}A)_{\approx T}} \\
	{U^{-1}(A_{\approx T})}
	\arrow["{\iota_2 \of \lambda_2}", from=1-1, to=1-2]
	\arrow["{\lambda_1 \of \iota_1}"', from=1-1, to=2-1]
	\arrow["\theta"', from=2-1, to=1-2]
\end{tikzcd}\]
is commutative.
\end{prop}

\begin{proof}
By \cref{prop:localization and monoid algebra of coalescent algebras}, the condensed $R$-algebra $U^{-1}(A_{\approx T})$ is $T$-coalescent. The condensed $R$-algebra $(U^{-1}A)_{\approx T}$ is $T$-coalescent by definition. For every $f \in U$, the element $(\lambda_1)_* (\iota_1)_* \sigma_* (f) \in U^{-1}(A_{\approx T})(*)$ is invertible in $U^{-1}(A_{\approx T})(*)$ by \cref{prop:universality of localization of condensed algebra}. Moreover, \cref{prop:universality of localization of condensed algebra} also shows that the element $(\lambda_2)_* \sigma_* (f) \in (U^{-1}A)(*)$ is invertible in $(U^{-1}A)(*)$. Therefore the element $(\iota_2)_* (\lambda_2)_* \sigma_* (f) \in (U^{-1}A)_{\approx T}(*)$ is invertible in $(U^{-1}A)_{\approx T}(*)$. This proves the property $(1)$. The property (2) is an immediate consequence of the universality of the localization (\cref{prop:universality of localization of condensed algebra}) and the universality of the $T$-coalescence of condensed $R$-algebras.
\end{proof}

\begin{prop} \label{prop:localization and taking monoid algebra commute}
Let $R$ be a condensed ring. Let $T \sub R(*) \times R(*)$. Let $M$ be a commutative monoid. Let $A$ be a condensed $R$-algebra. Let $R \xto{\sigma} A$ be the structure homomorphism of $A$. Let $A \xto{\iota_1} A_{\approx T}$ and $A[M] \xto{\iota_2} (A[M])_{\approx T}$ be the $T$-coalescence of the condensed $R$-algebras $A$ and $A[M]$ respectively. Let $A_{\approx T} \xto{\lambda_1} (A_{\approx T})[M]$ and $A \xto{\lambda_2} A[M]$ be the canonical inclusion. Then the compositions $\lambda_1 \of \iota_1 : A \to (A_{\approx T})[M]$ and $\iota_2 \of \lambda_2 : A \to (A[M])_{\approx T}$ have the following universal property.
\begin{enumerate}
\item
The condensed $R$-algebras $(A_{\approx T})[M]$ and $(A[M])_{\approx T}$ are $T$-coalescent.

\item
For every $T$-coalescent condensed $R$-algebra $B$, every homomorphism $\phi : A \to B$ of condensed $R$-algebras and every homomorphism $\tau : M \to B(*)$ of the monoid $M$ into the multiplicative monoid of the unital ring $B(*)$, there exist unique homomorphisms of condensed $R$-algebras $\psi_1 : (A_{\approx T})[M] \to B$ and $\psi_2 : (A[M])_{\approx T} \to B$ with the following properties. 
\begin{enumerate}
\item
The maps
\[\begin{tikzcd}
	M & {(A_{\approx T})(*)[M] = (A_{\approx T})[M](*)} && {B(*) \: ;} \\
	M & {A(*)[M] = A[M](*)} & {(A[M])_{\approx T}(*)} & {B(*)}
	\arrow["\inc", hook, from=1-1, to=1-2]
	\arrow["{(\psi_1)_*}", from=1-2, to=1-4]
	\arrow["\inc", hook, from=2-1, to=2-2]
	\arrow["{(\iota_2)_*}", from=2-2, to=2-3]
	\arrow["{(\psi_2)_*}", from=2-3, to=2-4]
\end{tikzcd}\]
are both equal to the map $\tau : M \to B(*)$.

\item
The diagrams
\[\begin{tikzcd}
	A & B & A & B \\
	{(A_{\approx T})[M]} && {(A[M])_{\approx T}}
	\arrow["\phi", from=1-1, to=1-2]
	\arrow["{\lambda_1 \of \iota_1}"', from=1-1, to=2-1]
	\arrow["\phi", from=1-3, to=1-4]
	\arrow["{\iota_2 \of \lambda_2}"', from=1-3, to=2-3]
	\arrow["{\psi_1}"', from=2-1, to=1-2]
	\arrow["{\psi_2}"', from=2-3, to=1-4]
\end{tikzcd}\]
are commutative.
\end{enumerate}
\end{enumerate}
In particular, there exists a unique isomorphism $\theta : (A_{\approx T})[M] \to (A[M])_{\approx T}$ of condensed $R$-algebras such that the diagram
\[\begin{tikzcd}
	A & {(A[M])_{\approx T}} \\
	{(A_{\approx T})[M]}
	\arrow["{\iota_2 \of \lambda_2}", from=1-1, to=1-2]
	\arrow["{\lambda_1 \of \iota_1}"', from=1-1, to=2-1]
	\arrow["\theta"', from=2-1, to=1-2]
\end{tikzcd}\]
is commutative in $\ub{CAlg}_R$ and the diagram
\[\begin{tikzcd}
	M & {A(*)[M] = A[M](*)} \\
	{(A_{\approx T})(*)[M] = (A_{\approx T})[M](*)} & {(A[M])_{\approx T}(*)}
	\arrow["\inc", hook, from=1-1, to=1-2]
	\arrow["\inc"', hook, from=1-1, to=2-1]
	\arrow["{(\iota_2)_*}", from=1-2, to=2-2]
	\arrow["{\theta_*}"', from=2-1, to=2-2]
\end{tikzcd}\]
is commutative in the category of commutative monoids.
\end{prop}

\begin{proof}
By \cref{prop:localization and monoid algebra of coalescent algebras}, the condensed $R$-algebra $(A_{\approx T})[M]$ is $T$-coalescent. The condensed $R$-algebra $(A[M])_{\approx T}$ is $T$-coalescent by definition. This proves the property (1). The property (2) is an immediate consequence of \cref{prop:universality of monoid algebras over condensed rings}.
\end{proof}

\section{Valuations on condensed rings} \label{sec:Valuations on condensed rings}

In this section, we study the notion of continuous valuations on condensed rings, which is a natural generalization of continuous valuations on commutative unital topological rings (\cite{Morel:note}, Definition II.2.1.1 ; \cite{Wedhorn:note}, Definition 7.7). This section mainly serves as a preperation for \cref{sec:Condensed ringed spaces}.

\begin{nt}
Let $v$ be a valuation on a commutative unital ring $R$. Let $\Gamma$ be the value group of $v$.
\begin{enumerate}
\item
In this paper, we always use the multiplicative notation for the value group $\Gamma$.

\item
If $f \in R$, then $|f|_v$ denotes the image of $f$ in $\Gamma \cup \{0\}$.

\item
$\mathrm{Supp}(v)$ denotes the prime ideal $\set{f \in R}{|f|_v = 0}$ of $R$. It is called the support of $v$.
\end{enumerate}
\end{nt}

\subsection{Definition}

\begin{df} \;
\begin{enumerate}
\item Let $R$ be a condensed ring. A \ti{valuation} on $R$ is a valuation on the ring $R(*)$. 

\item
Two valuations on a condensed ring $R$ is called \ti{equivalent} if they are equivalent as valuations on $R(*)$. As usual, two equivalent valuations are considered to be equal.

\item Let $\phi : R \to S$ be a homomorphism of condensed rings. Then we have a ring homomorphism $\phi_* : R(*) \to S(*)$. Therefore if $v$ is a valuation on $S$, then we have a valuation $v \of \phi_*$ on $R$. This is called the \ti{inverse image} of $v$ under $\phi$, and is denoted by $\phi^{-1}(v)$.
\end{enumerate}
\end{df}

\subsection{Continuous valuations}

\begin{prop} \label{prop:continuous valuations}
Let $R$ be a condensed ring. Let $v$ be a valuation on $R$. Let $K$ be the valued field of $v$. Let $\phi : R(*) \to K$ be the canonical map. We consider $K$ as a topological ring by giving it the valuation topology associated to $v$. Then the following conditions are equivalent.
\begin{enumerate}
\item
The map $\phi : R(*) \to K$ is continuous if $R(*)$ is given the topology defined in \cref{prop:reflection along underbar construction}.

\item
There exists a homomorphism $\psi : R \to \tu{K}$ of condensed rings such that $\psi_* = \phi$.
\end{enumerate}
Moreover, if these hold, then the homomorphism $\psi$ described in (2) is unique.
\end{prop}

\begin{df} \label{df:definition of continuous valuation on CRing}
We use the notation given in \cref{prop:continuous valuations}.
\begin{enumerate}
\item
The valuation $v$ on $R$ is called \ti{continuous} if it satisfies the equivalent conditions of \cref{prop:continuous valuations}.

\item
If $v$ is continuous, then the homomorphism $\psi : R \to \tu{K}$ of condensed rings described in \cref{prop:continuous valuations} is called the \ti{canonical homomorphism}.
\end{enumerate}
\end{df}

\begin{proof}[Proof of \cref{prop:continuous valuations}]
For each light profinite set $S$ and each $x \in S$, write $p_{S,x} : * \to S$ for the continuous map $* \mapsto x$. For each light profinite set $S$ and each $t \in R(S)$, let $\xi_{S,t} : S \to R(*)$ be the map
\[\begin{tikzcd}
	{\xi_{S,t} : S} & {R(*),} & x & {R(p_{S,x})(t).}
	\arrow[from=1-1, to=1-2]
	\arrow[maps to, from=1-3, to=1-4]
\end{tikzcd}\]
Then the topology on $R(*)$ defined in \cref{prop:reflection along underbar construction} is the finest topology $\cat{T}$ on $R(*)$ such that the maps $\xi_{S,t} : S \to R(*) \; (S \in |\ub{Prof}| ,\, t \in R(S))$ are all continuous when $R(*)$ is endowed with $\cat{T}$. Let us consider $R(*)$ as a topological space by giving it the topology $\cat{T}$. Then \cref{prop:reflection along underbar construction} shows that the following hold.
\begin{enumerate}
\item
For each light profinite set $S$, let $\eta_S : R(S) \to \mathrm{Cont}(S, R(*))$ be the map
\[\begin{tikzcd}
	{\eta_S : R(S)} & {\mathrm{Cont}(S,R(*)),} & t & {\xi_{S,t} .}
	\arrow[from=1-1, to=1-2]
	\arrow[maps to, from=1-3, to=1-4]
\end{tikzcd}\]
Then the family $\eta := (\eta_S)_{S \in |\ub{Prof}|}$ is a map $R \to \tu{R(*)}$ of condensed sets.

\item
The following maps are bijective and inverse to each other.
\[\begin{tikzcd}[row sep=tiny]
	{\mathrm{Cont}(R(*),K)} & {\ub{CSet}(R, \tu{K})} \\
	f & {\tu{f} \of \eta} \\
	{( \alpha_* : R(*) \to \tu{K}(*) = K )} & \alpha
	\arrow["\sim", tail reversed, from=1-1, to=1-2]
	\arrow["\Gamma", maps to, from=2-1, to=2-2]
	\arrow["\Delta"', maps to, from=3-2, to=3-1]
\end{tikzcd}\]
\end{enumerate}

Therefore if there exists a homomorphism $\psi : R \to \tu{K}$ of condensed rings such that $\psi_* : R(*) \to \tu{K}(*) = K$ is equal to $\phi$, then $\phi = \psi_* = \Delta(\psi)$ is an element of $\mathrm{Cont}(R(*),K)$. In other words, the map $\phi : R(*) \to K$ is continuous. Furthermore, if another homomorphism $\rho : R \to \tu{K}$ of condensed rings satisfies $\rho_* = \phi$, then
\begin{equation}
\Delta(\rho) = \rho_* = \phi = \psi_* = \Delta(\psi) .
\end{equation}
Since the map $\Delta$ is injective, we conclude that $\rho = \psi$. This proves the uniqueness of $\psi$.

Conversely, suppose that the map $\phi : R(*) \to K$ is continuous. We prove that there exists a homomorphism $\psi : R \to \tu{K}$ of condensed rings such that $\psi_* = \phi$. Let us define $\psi := \Gamma(\phi)$. This is a map $R \to \tu{K}$ of condensed sets and satisfies $\psi_* = \Delta(\psi) = \phi$. Therefore it remains to prove that $\psi : R \to \tu{K}$ is a homomorphism of condensed rings. Let $S \in |\ub{Prof}|$. We show that the map $\psi_S : R(S) \to \tu{K}(S)$ is a ring homomorphism. Since $\psi = \Gamma(\phi) = \tu{\phi} \of \eta$, the map $\psi_S$ is equal to the composition
\[\begin{tikzcd}
	{R(S)} & {\mathrm{Cont}(S,R(*))} && {\mathrm{Cont}(S,K).}
	\arrow["{\eta_S}", from=1-1, to=1-2]
	\arrow["{f \mapsto \phi \of f}", from=1-2, to=1-4]
\end{tikzcd}\]
Let $t_1, t_2 \in R(S)$. For each $x \in S$, the map $R(p_{S,x}) : R(S) \to R(*)$ is a ring homomorphism. Therefore
\begin{align}
\xi_{S, (t_1 + t_2)}(x) & = R(p_{S,x})(t_1 + t_2)
= R(p_{S,x})(t_1) + R(p_{S,x})(t_2) = \xi_{S,t_1}(x) + \xi_{S,t_2}(x) \: ; \\
\xi_{S, (t_1 \cdot t_2)}(x) & = R(p_{S,x})(t_1 \cdot t_2)
= R(p_{S,x})(t_1) \cdot R(p_{S,x})(t_2) = \xi_{S,t_1}(x) \cdot \xi_{S,t_2}(x) \: ; \\
\xi_{S,1}(x) & = R(p_{S,x})(1) = 1 .
\end{align}
It follows that
\begin{align}
\psi_S(t_1 + t_2)(x)
& = \Big( \phi \of \big( \eta_S (t_1 + t_2) \big) \Big) (x)
   = \phi \big( \xi_{S, (t_1 + t_2)}(x) \big) \\
& = \phi \big( \xi_{S,t_1}(x) + \xi_{S,t_2}(x) \big)
   = \phi \big( \xi_{S,t_1}(x) \big) + \phi \big( \xi_{S,t_2}(x) \big) \\
& = \Big( \phi \of \big( \eta_S (t_1) \big) \Big) (x) + \Big( \phi \of \big( \eta_S (t_2) \big) \Big) (x)
   = \psi_S(t_1)(x) + \psi_S(t_2)(x) \: ; \\
\psi_S(t_1 \cdot t_2)(x)
& = \Big( \phi \of \big( \eta_S (t_1 \cdot t_2) \big) \Big) (x)
   = \phi \big( \xi_{S, (t_1 \cdot t_2)}(x) \big) \\
& = \phi \big( \xi_{S,t_1}(x) \cdot \xi_{S,t_2}(x) \big)
   = \phi \big( \xi_{S,t_1}(x) \big) \cdot \phi \big( \xi_{S,t_2}(x) \big) \\
& = \Big( \phi \of \big( \eta_S (t_1) \big) \Big) (x) \cdot \Big( \phi \of \big( \eta_S (t_2) \big) \Big) (x)
   = \psi_S(t_1)(x) \cdot \psi_S(t_2)(x) \: ; \\
\psi_S(1)(x)
& = \Big( \phi \of \big( \eta_S (1) \big) \Big) (x)
   = \phi \big( \xi_{S,1}(x) \big)
   = \phi (1) = 1.
\end{align}
Since these hold for every $x \in S$, we conclude that
\begin{align}
\psi_S(t_1 + t_2) & = \psi_S(t_1) + \psi_S(t_2) \: ; \\
\psi_S(t_1 \cdot t_2) & = \psi_S(t_1) \cdot \psi_S(t_2) \: ; \\
\psi_S(1) & = 1 .
\end{align}
This completes the proof that the map $\psi_S : R(S) \to \tu{K}(S)$ is a ring homomorphism.
\end{proof}

\begin{rem} \label{rem:comparison of continuous valuations on R and on underbar R}
Let $R$ be a commutative unital topological ring. By \cref{rem:comparison of topologies of underbar X(*)}, the following hold.
\begin{enumerate}
\item
If $v$ is a continuous valuation on the topological ring $R$ in the sense of Definition II.2.1.1 of \cite{Morel:note} and Definition 7.7 of \cite{Wedhorn:note}, then $v$ is a continuous valuation on the condensed ring $\tu{R}$ in the sense of \cref{df:definition of continuous valuation on CRing}. 

\item
If $R$ is first countable, then the notion of continuous valuations on the topological ring $R$ coincides with the notion of continuous valuations on the condensed ring $\tu{R}$.
\end{enumerate} 
\end{rem}

\begin{prop}
Let $\rho : R \to S$ be a homomorphism of condensed rings. If $w$ is a continuous valuation on $S$, then $\rho^{-1}(w)$ is a continuous valuation on $R$.
\end{prop}

\begin{proof}
Let us write $v := \rho^{-1}(w)$. Let $K,L$ be the valued field of $v,w$ respectively. Let $\phi : R(*) \to K$, $\psi : S(*) \to L$ be the canonical maps. Then we have a canonical extension $\iota : K \mon L$ of valued fields such that the following diagram is commutative in $\ub{Ring}$.
\[\begin{tikzcd}
	{R(*)} & K \\
	{S(*)} & L
	\arrow["\phi", from=1-1, to=1-2]
	\arrow["{\rho_*}"', from=1-1, to=2-1]
	\arrow["\iota", hook, from=1-2, to=2-2]
	\arrow["\psi"', from=2-1, to=2-2]
\end{tikzcd}\]
Let us consider $R(*), S(*)$ as topological spaces by giving them the topologies defined in \cref{prop:reflection along underbar construction}. Let us endow $L$ with the valuation topology associated to $w$. Let $\cat{T}$ be the initial topology on $K$ induced by the map $\iota : K \mon L$. Then the map $\psi : S(*) \to L$ is continuous. By \cref{prop:functoriality of reflection along underbar construction}, the map $\rho_* : R(*) \to S(*)$ is continuous. Therefore if $K$ is endowed with the topology $\cat{T}$, the map $\phi : R(*) \to K$ is continuous. Since the valuation topology on $K$ associated to $v$ is coarser than the topology $\cat{T}$, it follows that the map $\phi : R(*) \to K$ is continuous when $K$ is endowed with the valuation topology associated to $v$.
\end{proof}

\subsection{Valued condensed rings}

\begin{df} \label{df:valued condensed ring} \;
\begin{enumerate}
\item
A \ti{valued condensed ring} is a pair $(R,v)$ consisting of a condensed ring $R$ and a continuous valuation $v$ on $R$.

\item
Let $(R,v),(S,w)$ be valued condensed rings. A \ti{homomorphism} $(R,v) \to (S,w)$ \ti{of valued condensed rings} is a homomorphism $\phi : R \to S$ of condensed rings such that $\phi^{-1}(w) = v$.

\item
We obtain the category $\ub{VCRing}$ of valued condensed rings and homomorphisms between them.
\end{enumerate}
\end{df}

\begin{prop} \label{prop:extension of continuous valuation to filtered colimits}
Let $\cat{I}$ be a small filtered category. Let $D : \cat{I} \to \ub{VCRing}$ be a functor. For each $i \in \ob{I}$, let us write $D(i) = (R_i,v_i)$. For each morphism $\alpha$ in $\cat{I}$, let us write $D(\alpha) = \phi_{\alpha}$. Let $F : \ub{VCRing} \to \ub{CRing}$ be the forgetful functor $(R,v) \mapsto R$. Let $(R, (R_i \xto{\rho_i} R)_{i \in \ob{I}})$ be the colimit of $F \of D : \cat{I} \to \ub{CRing}$. Suppose that for each morphism $\alpha : i \to j$ in $\cat{I}$, the map $(\phi_{\alpha})_* : R_i(*) \to R_j(*)$ induces an isomorphism of the value group of $v_i$ onto that of $v_j$. Then there exists a unique continuous valuation $v$ on $R$ such that $\rho_i^{-1}(v) = v_i$ for every $i \in \ob{I}$.
\end{prop}

\begin{proof}
The forgetful functor $U : \ub{CRing} \to \ub{CSet}$ preserves filtered colimits. The functor $E : \ub{CSet} \to \ub{Set}$, $X \mapsto X(*)$ preserves filtered colimits by \cref{prop:limits and filtered colimits in CSet}. Therefore $(R(*), (R_i(*) \xto{(\rho_i)_*} R(*))_{i \in \ob{I}})$ is the colimit of the functor $E \of U \of F \of D : \cat{I} \to \ub{Set}$. On the other hand, the assumption implies that the value groups of the valuations $v_i$ $(i \in \ob{I})$ are all the same. Let us write $\Gamma$ for this common value group. Then we have the map $|\cdot|_{v_i} : R_i(*) \to \Gamma \cup \{0\}$ for each $i \in \ob{I}$, and the diagram
\[\begin{tikzcd}
	{R_i(*)} & {\Gamma \cup \{0\}} \\
	{R_j(*)}
	\arrow["{|\cdot|_{v_i}}", from=1-1, to=1-2]
	\arrow["{(\phi_{\alpha})_*}"', from=1-1, to=2-1]
	\arrow["{|\cdot|_{v_j}}"', from=2-1, to=1-2]
\end{tikzcd}\]
is commutative for each morphism $\alpha : i \to j$ in $\cat{I}$. Therefore there exist a unique map $|\cdot|_v : R(*) \to \Gamma \cup \{0\}$ such that the diagram 
\[\begin{tikzcd}
	{R(*)} & {\Gamma \cup \{0\}} \\
	{R_i(*)}
	\arrow["{|\cdot|_v}", from=1-1, to=1-2]
	\arrow["{(\rho_i)_*}", from=2-1, to=1-1]
	\arrow["{|\cdot|_{v_i}}"', from=2-1, to=1-2]
\end{tikzcd}\]
is commutative for every $i \in\ob{I}$. Using the properties of filtered colimits in $\ub{Set}$, one checks that this map $|\cdot|_v : R(*) \to \Gamma \cup \{0\}$ defines a unique valuation $v$ on $R$ such that $\rho_i^{-1}(v) = v_i$ for every $i \in \ob{I}$. 

We claim that $v$ is a continuous valuation on $R$. Let us consider $R(*)$ and $R_i(*)$ $(i \in \ob{I})$ as topological spaces by giving them the topologies defined in \cref{prop:reflection along underbar construction}. Let $K$ be the valued field of $v$, which is given the valuation tpology asociated to $v$. We prove that the canonical homomorphism $\psi : R(*) \to K$ is continuous. The continuity of $\psi$ is equivalent to the condition that the valuation topology on $R(*)$ associated to $v$ is coarser than the topology of $R(*)$. Thus it suffices to show that for every $x \in R(*)$ and every $\gamma \in \Gamma$, the subset $\set{y \in R(*)}{|y-x|_v < \gamma}$ is open in $R(*)$. Let $x \in R(*)$ and $\gamma \in \Gamma$. Write $A := \set{y \in R(*)}{|y-x|_v < \gamma}$. We show that $A$ is an open subset of $R(*)$. Consider the functor $V : \ub{CSet} \to \ub{Top}$ defined in \cref{cor:left adjoint of the underbar functor}. Since this functor is a left adjoint, it preserves colimits. Moreover, the forgetful functor $U : \ub{CRing} \to \ub{CSet}$ preserves filtered colimits. It follows that $(R(*), (R_i(*) \xto{(\rho_i)_*} R(*))_{i \in \ob{I}})$ is the colimit of the functor $V \of U \of F \of D : \cat{I} \to \ub{Top}$. Then the construction of colimits in $\ub{Top}$ shows that $A$ is open in $R(*)$ if and only if $(\rho_i)_*^{-1}(A)$ is an open subset of $R_i(*)$ for every $i \in \ob{I}$. Therefore we fix $i \in \ob{I}$ arbitrarily and prove that $(\rho_i)_*^{-1}(A)$ is an open subset of $R_i(*)$. The construction of filtered colimits in $\ub{Top}$ shows that there exists an object $j \in \ob{I}$ and an element $x_j \in R_j(*)$ such that $(\rho_j)_*(x_j) = x$. Since the category $\cat{I}$ is filtered, there exists an object $k \in \ob{I}$ and morphisms $\alpha : i \to k$ and $\beta : j \to k$ in $\cat{I}$. Then the diagrams
\[\begin{tikzcd}
	{R_i(*)} \\
	{R_k(*)} & {R(*)} & {R(*)} & {\Gamma \cup \{0\}} \\
	{R_j(*)} && {R_k(*)}
	\arrow["{(\phi_{\alpha})_*}"', from=1-1, to=2-1]
	\arrow["{(\rho_i)_*}", from=1-1, to=2-2]
	\arrow["{(\rho_k)_*}", from=2-1, to=2-2]
	\arrow["{|\cdot|_v}", from=2-3, to=2-4]
	\arrow["{(\phi_{\beta})_*}", from=3-1, to=2-1]
	\arrow["{(\rho_j)_*}"', from=3-1, to=2-2]
	\arrow["{(\rho_k)_*}", from=3-3, to=2-3]
	\arrow["{|\cdot|_{v_k}}"', from=3-3, to=2-4]
\end{tikzcd}\]
are commutative. Therefore $x= (\rho_j)_*(x_j) = (\rho_k)_*(\phi_{\beta})_*(x_j)$ and 
\begin{align}
(\rho_k)_*^{-1}(A)
& = \set{y \in R_k(*)}{|(\rho_k)_*(y) - (\rho_k)_*(\phi_{\beta})_*(x_j)|_v < \gamma} \\
& = \set{y \in R_k(*)}{|y - (\phi_{\beta})_*(x_j)|_{v_k} < \gamma}.
\end{align}
This is an open subset of $R_k(*)$ since $v_k$ is a continuous valuation on $R_k$. Therefore $(\rho_i)_*^{-1}(A) = (\phi_{\alpha})_*^{-1} (\rho_k)_*^{-1}(A)$ is an open subset of $R_i(*)$.
\end{proof}

\subsection{Localization of valued condensed rings}

\begin{df} \label{df:local valued ring}
The category $\ub{VCRing}_l$ is defined to be the full subcategory of $\ub{VCRing}$ consisting of all $(R,v) \in |\ub{VCRing}|$ with the following property: The ring $R(*)$ is a local ring and its maximal ideal is equal to the support of $v$.
\end{df}

\begin{lem} \label{lem:extension of continuous valuation along localization}
Let $(R,v)$ be a valued condensed ring. Let $T$ be a multiplicatively closed subset of $R(*)$. Let $\iota : R \to T^{-1}R$ be the canonical homomorphism of the condensed $R$-algebra $R$ into its localization $T^{-1}R$ by $T$. Suppose that $|g|_v \neq 0$ for every $g \in T$. Then there exists a unique continuous valuation $w$ on $T^{-1}R$ such that $\iota^{-1}(w) = v$.
\end{lem}

\begin{proof}
The uniqueness follows immediately from the fact that $(T^{-1}R)(*) = T^{-1}(R(*))$. We prove the existence. Let $K$ be the valued field of the valuation $v$ on $R$. We consider $K$ as a topological ring by giving it the valuation topology associated to $v$. Let $\tilde{v}$ be the valuation on $K$ which is the canonical extension of $v$. Then $\tilde{v}$ is a continuous valuation on $\tu{K}$. Let $\psi : R \to \tu{K}$ be the canonical homomorphism. For every $g \in T$, we have $|g|_v \neq 0$, and hence $\psi_*(g) \neq 0$. Thus $\psi_*(g) \in K$ is invertible in $K$ for every $g \in T$. Consequently, \cref{prop:universality of localization of condensed algebra} shows that there exists a unique hmomorphism $\rho : T^{-1} R \to \tu{K}$ of condensed rings such that the following diagram is commutative.
\[\begin{tikzcd}
	R & {\tu{K}} \\
	{T^{-1}R}
	\arrow["\psi", from=1-1, to=1-2]
	\arrow["\iota"', from=1-1, to=2-1]
	\arrow["\rho"', from=2-1, to=1-2]
\end{tikzcd}\]
Then $w := \rho^{-1}(\tilde{v})$ is a continuous valuation on $T^{-1}R$ and satisfies $\iota^{-1}(w) = \iota^{-1} \big( \rho^{-1}(\tilde{v}) \big) = \psi^{-1}(\tilde{v}) = v$.
\end{proof}

\begin{prop} \label{prop:construction of localization of vcrings}
Let $(R,v)$ be a valued condensed ring. Let $\fk{p}$ be the support of $v$. Let $\iota : R \to R_{\fk{p}}$ be the canonical homomorphism of $R$ into the localization $R_{\fk{p}}$.  Then the following hold.
\begin{enumerate}
\item
There exists a unique continuous valuation $\tilde{v}$ on $R_{\fk{p}}$ such that $\iota^{-1}(\tilde{v}) = v$.

\item
The valued condensed ring $(R_{\fk{p}}, \tilde{v})$ is an object of $\ub{VCRing}_l$.

\item
The pair $( (R_{\fk{p}}, \tilde{v}) , (R,v) \xto{\iota} (R_{\fk{p}}, \tilde{v}) )$ is a reflection of $(R,v) \in |\ub{VCRing}|$ along the inclusion functor $\ub{VCRing}_l \mon \ub{VCRing}$.
\end{enumerate}
\end{prop}

\begin{proof} ~
\begin{enumerate}
\item
This immediately follows from \cref{lem:extension of continuous valuation along localization}.

\item
$(R_{\fk{p}})(*) = (R(*))_{\fk{p}}$ is a local ring. Moreover, since $\iota^{-1}(\tilde{v}) = v$, we have $(\iota_*)^{-1} (\mathrm{Supp}(\tilde{v})) = \mathrm{Supp}(v) = \fk{p}$. Therefore $\mathrm{Supp}(\tilde{v}) = \fk{p} \cdot (R(*))_{\fk{p}}$ is the maximal ideal of $(R(*))_{\fk{p}} = (R_{\fk{p}})(*)$.

\item
Suppose that $(S,w)$ is an object of $\ub{VCRing}_l$ and that $\phi : (R,v) \to (S,w)$ is a homomorphism of valued condensed rings. We prove that there exists a unique homomorphism $\psi : (R_{\fk{p}}, \tilde{v}) \to (S,w)$ of valued condensed rings such that the following diagram is commutative.
\[\begin{tikzcd}
	{(R,v)} & {(S,w)} \\
	{(R_{\fk{p}},\tilde{v})}
	\arrow["\phi", from=1-1, to=1-2]
	\arrow["\iota"', from=1-1, to=2-1]
	\arrow["\psi"', from=2-1, to=1-2]
\end{tikzcd}\]

If $g \in R(*) \sm \fk{p}$, then $|g|_v \neq 0$. Then $|\phi_*(g)|_w \neq 0$ since $\phi^{-1}(w) = v$. Thus $\phi_*(g) \in S(*) \sm \mathrm{Supp}(w)$. Since $(S,w)$ is an object of $\ub{VCRing}_l$, the ring $S(*)$ is a local ring and its maximal ideal is equal to $\mathrm{Supp}(w)$. Therefore $S(*) \sm \mathrm{Supp}(w)$ is equal to the set of invertible elements of $S(*)$. Consequently, $\phi_*(g) \in S(*)$ is invertible in $S(*)$. Then \cref{prop:universality of localization of condensed algebra} shows that there exists a unique homomorphism $\psi : R_{\fk{p}} \to S$ of condensed rings such that the following diagram is commutative.
\[\begin{tikzcd}
	R & S \\
	{R_{\fk{p}}}
	\arrow["\phi", from=1-1, to=1-2]
	\arrow["\iota"', from=1-1, to=2-1]
	\arrow["\psi"', from=2-1, to=1-2]
\end{tikzcd}\]
Then $\iota^{-1} \big( \psi^{-1}(w) \big) = \phi^{-1}(w) = v$. By (1), we have $\psi^{-1}(w) = \tilde{v}$. Therefore $\psi$ is a homomorphism $(R_{\fk{p}}, \tilde{v}) \to (S,w)$ of valued condensed rings. This completes the proof.
\end{enumerate}
\end{proof}

\begin{cor}
$\ub{VCRing}_l$ is a reflective full subcategory of $\ub{VCRing}$.
\end{cor}

\begin{df} \label{df:localization of valued condensed ring}
Let $(R,v)$ be a valued condensed ring. Let $\big( (\tilde{R}, \tilde{v}) , (R,v) \xto{\iota} (\tilde{R}, \tilde{v}) \big)$ be the reflection of $(R,v) \in |\ub{VCRing}|$ along the inclusion functor $\ub{VCRing}_l \mon \ub{VCRing}$. Then $(\tilde{R}, \tilde{v}) \in |\ub{VCRing}_l|$ is called the \ti{localization of the valued condensed ring} $(R,v)$. The homomorphism $(R,v) \xto{\iota} (\tilde{R}, \tilde{v})$ of valued condensed rings is also called the localization of the valued condensed ring $(R,v)$.
\end{df}

\subsection{Coalescence of valued condensed rings}

\begin{df} \label{df:coalescent valued ring}
The category $\ub{VCRing}_c$ is defined to be the full subcategory of $\ub{VCRing}$ consisting of all $(R,v) \in |\ub{VCRing}|$ with the following property: For every $f, g \in R(*)$ with $|f|_v \leq |g|_v \neq 0$, the condensed $R$-algebra $R$ is $f/g$-coalescent.
\end{df}

\begin{lem} \label{lem:restriction of valuations on fields}
Let $L$ be a field. Let $w$ be a valuation on $L$. We consider $L$ as a topological ring by giving it the valuation topology associated to $w$. Let $\iota : L \to \hat{L}$ be the Hausdorff completion of $L$. Let $K$ be a subfield of $L$. Let $v$ be the valuation on $K$ which is the restriction of $w$ to $K$.
\begin{enumerate}
\item
The closure $\overline{\iota(K)}$ of $\iota(K)$ in $\hat{L}$ is a subfield of $\hat{L}$.

\item
Suppose that $\overline{\iota(K)} = \hat{L}$. Then the following hold.
\begin{enumerate}
\item
$K$ is dense in $L$.

\item
The value group of $v$ coincides with the value group of $w$.

\item
The valuation topology on $K$ associated to $v$ is equal to the subspace topology on $K$ induced by the topology of $L$.
\end{enumerate}
\end{enumerate}
\end{lem}

\begin{proof} ~
\begin{enumerate}
\item
$\hat{L}$ is a field and $\iota : L \to \hat{L}$ is a field extension. Therefore $\iota(K)$ is a subfield of $\hat{L}$. Since $\hat{L}$ is a topological field, the closure $\overline{\iota(K)}$ is also a subfield of $\hat{L}$.

\item
\begin{enumerate}
\item
The closure $\overline{K}$ of $K$ in $L$ is equal to $\iota^{-1} ( \overline{\iota(K)} )$. Therefore $\overline{K} = \iota^{-1} ( \overline{\iota(K)} ) = \iota^{-1} (\hat{L}) = L$.

\item
The value group $\Gamma$ of $v$ is a subgroup of the value group $\Delta$ of $w$. If $\delta \in \Delta$, then $\set{x \in L}{|x|_w = \delta}$ is a nonempty open subset of $L$. Since $K$ is dense in $L$ by (a), this set intersects with $K$. In other words, there exists an $x \in K$ such that $|x|_w = \delta$. It follows that $\delta \in \Gamma$. Consequently $\Gamma = \Delta$.

\item
This follows from (b) and the definition of the valuation topology.
\end{enumerate}
\end{enumerate}
\end{proof}

\begin{lem} \label{lem:extension of continuous valuation along coalescence}
Let $(R,v)$ be a valued condensed ring. Let $T$ be a subset of $R(*) \times R(*)$. Let $\iota : R \to R_{\approx T}$ be the $T$-coalescence of the condensed $R$-algebra $R$. Suppose that $|f|_v \leq |g|_v \neq 0$ for every $(f,g) \in T$. Then there exists a unique continuous valuation $w$ on $R_{\approx T}$ such that $\iota^{-1}(w) = v$. Moreover, the map $\iota_* : R(*) \to R_{\approx T}(*)$ induces an isomorphism of the value group of $v$ onto that of $w$.
\end{lem}

\begin{proof}
Let $K$ be the valued field of $v$. We consider $K$ as a topological ring by giving it the valuation topology associated to $v$. Let $\psi : R \to \tu{K}$ be the canonical homomorphism. Let $i : K \to \hat{K}$ be the Hausdorff completion of $K$. We write $\hat{v}$ for the valuation on $\hat{K}$ which is the canonical extension of $v$ to $\hat{K}$. Then $\hat{v}$ is a continuous valuation on $\tu{\hat{K}}$. On the other hand, we have the following homomorphism of condensed rings.
\[\begin{tikzcd}
	R & {\tu{K}} & {\tu{\hat{K}}}
	\arrow["\psi", from=1-1, to=1-2]
	\arrow["{\tu{i}}", from=1-2, to=1-3]
\end{tikzcd}\]
If $(f,g) \in T$, then $|f|_v \leq |g|_v \neq 0$. Therefore $\psi_*(g)$ is invertible in $K$ and $\psi_*(f) \cdot \psi_*(g)^{-1} \in K$ is in the valuation ring of $v$. Hence $i(\psi_*(g))$ is invertible in $\hat{K}$ and $i(\psi_*(f)) \cdot i(\psi_*(g))^{-1} \in \hat{K}$ is in the valuation ring of $\hat{v}$. Therefore $i(\psi_*(f)) \cdot i(\psi_*(g))^{-1}$ is a power-bounded element of $\hat{K}$. Then \cref{prop:power-boundedness and coalescence} below shows that $\tu{\hat{K}}$ is a $\left( i(\psi_*(f)) \cdot i(\psi_*(g))^{-1} \right) / 1$-coalescent condensed $\tu{\hat{K}}$-module. Since $i(\psi_*(g))$ is invertible in $\hat{K}$, the condensed $\tu{\hat{K}}$-module $\tu{\hat{K}}$ is $i(\psi_*(f)) / i(\psi_*(g))$-coalescent. It follows that $\tu{\hat{K}}$ is an $f/g$-coalescent condensed $R$-algebra via $R \xto{\psi} \tu{K} \xto{\tu{i}} \tu{\hat{K}}$. Since this holds for every $(f,g) \in T$, we conclude that $\tu{\hat{K}}$ is a $T$-coalescent condensed $R$-algebra via $R \xto{\psi} \tu{K} \xto{\tu{i}} \tu{\hat{K}}$. Therefore there exists a unique homomorphism $\rho : R_{\approx T} \to \tu{\hat{K}}$ of condensed rings such that the following diagram is commutative.
\[\begin{tikzcd}
	R & {\tu{K}} \\
	{R_{\approx T}} & {\tu{\hat{K}}}
	\arrow["\psi", from=1-1, to=1-2]
	\arrow["\iota"', from=1-1, to=2-1]
	\arrow["{\tu{i}}", from=1-2, to=2-2]
	\arrow["\rho"', from=2-1, to=2-2]
\end{tikzcd}\]
Then $w := \rho^{-1}(\hat{v})$ is a continuous valuation on $R_{\approx T}$ and satisfies $\iota^{-1}(w) = v$. Moreover, since the value group of $v$ and the value group of $\hat{v}$ coincide, the map $\iota_* : R(*) \to R_{\approx T}(*)$ necessarily induces an isomorphism of the value group of $v$ onto the value group of $w$.

Suppose that $w'$ is another continuous valuation on $R_{\approx T}$ such that $\iota^{-1}(w') = v$. We prove that $w' = w$. Let $L$ be the valued field of $w'$. We consider $L$ as a topological ring by giving it the valuation topology associated to $w'$. Let $\tau : R_{\approx T} \to \tu{L}$ be the canonical homomorphism. Let $j : L \to \hat{L}$ be the Hausdorff completion of $L$. We write $\hat{w}'$ for the valuation on $\hat{L}$ which is the canonical extension of $w'$ to $\hat{L}$. Since $\iota^{-1}(w') = v$, there exists a canonical extension $k : K \to L$ of valued fields such that the following diagram is commutative in $\ub{Ring}$.
\[\begin{tikzcd}
	{R(*)} & K \\
	{R_{\approx T}(*)} & L
	\arrow["{\psi_*}", from=1-1, to=1-2]
	\arrow["{\iota_*}"', from=1-1, to=2-1]
	\arrow["k", from=1-2, to=2-2]
	\arrow["{\tau_*}"', from=2-1, to=2-2]
\end{tikzcd}\]
Let $K'$ be the field $K$ endowed with the initial topology induced by the map $k : K \to L$. $K'$ is a topological field. We write $k' : K' \to L$ for the map $k$. By the definition of valuation topology, the map $K' \to K$, $x \mapsto x$ is continuous. We write $c : K' \to K$ for this map. Furthermore, let us consider $R(*)$ and $R_{\approx T}(*)$ as topological spaces by giving them the topologies defined in \cref{prop:reflection along underbar construction}. Since $w'$ is a continuous valuation on $R_{\approx T}$, the map $\tau_* : R_{\approx T}(*) \to L$ is continuous. By \cref{prop:functoriality of reflection along underbar construction}, the map $\iota_* : R(*) \to R_{\approx T}(*)$ is continuous. Therefore the map $\psi_* : R(*) \to K'$ is continuous. To avoid confusion, we write $\tilde{\psi_*} : R(*) \to K'$ for this map $\psi_*$. Let $\eta : R \to \tu{R(*)}$ be the map of condensed sets obtained by applying \cref{prop:reflection along underbar construction} to the condensed set $R$. We write $\psi' : R \to \tu{K'}$ for the map of condensed sets $R \xto {\eta} \tu{R(*)} \xto{\tu{\tilde{\psi_*}} \;} \tu{K'}$. Then (2) of \cref{prop:reflection along underbar construction} shows that $\tilde{\psi_*} = \psi'_* : R(*) \to K'$. Therefore
\begin{equation}
(\tu{c} \of \psi')_* = \tu{c}_* \of \psi'_* = c \of \tilde{\psi_*} = \psi_* .
\end{equation}
By (2) of \cref{prop:reflection along underbar construction}, we conclude that $\tu{c} \of \psi' = \psi$. In other words, the following diagram is commutative in $\ub{CSet}$.
\[\begin{tikzcd}
	& {\tu{K}} \\
	R & {\tu{K'}}
	\arrow["\psi", from=2-1, to=1-2]
	\arrow["{\psi'}"', from=2-1, to=2-2]
	\arrow["{\tu{c}}"', from=2-2, to=1-2]
\end{tikzcd}\]
$\tu{c} : \tu{K'} \to \tu{K}$ defines an isomorphism of the condensed ring $\tu{K'}$ onto a condensed subring of $\tu{K}$. Moreover, $\psi : R \to \tu{K}$ is a homomorphism of condensed rings. Therefore we conclude that $\psi' : R \to \tu{K'}$ is a homomorphism of condensed rings. Furthermore, the following diagram is commutative.
\[\begin{tikzcd}
	{R(*)} & {K'} \\
	{R_{\approx T}(*)} & L
	\arrow["{\psi'_* = \tilde{\psi_*}}", from=1-1, to=1-2]
	\arrow["{\iota_*}"', from=1-1, to=2-1]
	\arrow["{k'}", from=1-2, to=2-2]
	\arrow["{\tau_*}"', from=2-1, to=2-2]
\end{tikzcd}\]
By (2) of \cref{prop:reflection along underbar construction}, we conclude that the following diagram is commutative in $\ub{CSet}$, and hence in $\ub{CRing}$.
\[\begin{tikzcd}
	R & {\tu{K'}} \\
	{R_{\approx T}} & {\tu{L}}
	\arrow["{\psi'}", from=1-1, to=1-2]
	\arrow["\iota"', from=1-1, to=2-1]
	\arrow["{\tu{k'}}", from=1-2, to=2-2]
	\arrow["\tau"', from=2-1, to=2-2]
\end{tikzcd}\]

Let $\hat{K}'$ be the closure of $j \big( k'(K') \big)$ in $\hat{L}$. By \cref{lem:restriction of valuations on fields}, $\hat{K}'$ is a subfield of $\hat{L}$. Let $\hat{k}' : \hat{K}' \to \hat{L}$ be the inclusion. Let $i' : K' \to \hat{K}'$ be the restriction of $K' \xto{k'} L \xto{j} \hat{L}$. Then $i' : K' \to \hat{K}'$ is the Hausdorff completion of $K'$, and the following diagram is commutative in $\ub{CRing}$.
\[\begin{tikzcd}
	R & {\tu{K'}} & {\tu{\hat{K}'}} \\
	{R_{\approx T}} & {\tu{L}} & {\tu{\hat{L}}}
	\arrow["{\psi'}", from=1-1, to=1-2]
	\arrow["\iota"', from=1-1, to=2-1]
	\arrow["{\tu{i'}}", from=1-2, to=1-3]
	\arrow["{\tu{k'}}", from=1-2, to=2-2]
	\arrow["{\tu{\hat{k}'}}", from=1-3, to=2-3]
	\arrow["\tau"', from=2-1, to=2-2]
	\arrow["{\tu{j}}"', from=2-2, to=2-3]
\end{tikzcd}\]
We claim that $\tu{\hat{L}}$ is a $T$-coalescent condensed $R$-algebra via $R \xto{\iota} R_{\approx T} \xto{\tau} \tu{L} \xto{\tu{j} \;} \tu{\hat{L}}$ and that $\tu{\hat{K}'}$ is a $T$-coalescent condensed $R$-algebra via $R \xto{\psi'} \tu{K'} \xto{\tu{i'} \;} \tu{\hat{K}'}$. Let $(f,g) \in T$. Then $|f|_v \leq |g|_v \neq 0$. Since the valuation $\hat{w}'$ extends $v$ in the sense that $\hat{w}' \of j \of \tau_* \of \iota_* = v$, we have
\begin{equation}
| (j \of \tau_* \of \iota_*)(f) |_{\hat{w}'} \leq
| (j \of \tau_* \of \iota_*)(g) |_{\hat{w}'} \neq 0 .
\end{equation}
Therefore
\begin{equation}
(j \of \tau_* \of \iota_*)(g) = (\hat{k}' \of i' \of \psi'_*)(g) \, \in \, \hat{L}
\end{equation}
is invertible in $\hat{L}$ and the element
\begin{equation}
(j \of \tau_* \of \iota_*)(f) \cdot (j \of \tau_* \of \iota_*)(g)^{-1}
= (\hat{k}' \of i' \of \psi'_*)(f) \cdot (\hat{k}' \of i' \of \psi'_*)(g)^{-1}
\, \in \, \hat{L}
\end{equation}
is in the valuation ring of $\hat{w}'$, which shows that this is a power-bounded element of $\hat{L}$. Then since $\hat{K}'$ is a topological subfield of $\hat{L}$,
\begin{equation}
(i' \of \psi'_*)(g) \, \in \, \hat{K}'
\end{equation}
is invertible in $\hat{K}'$ and
\begin{equation}
(i' \of \psi'_*)(f) \cdot (i' \of \psi'_*)(g)^{-1}
\, \in \, \hat{K}'
\end{equation}
is a power-bounded element of $\hat{K}'$. Then \cref{prop:power-boundedness and coalescence} below shows that $\tu{\hat{L}}$ is a $\big( (j \of \tau_* \of \iota_*)(f) \cdot (j \of \tau_* \of \iota_*)(g)^{-1} \big) /1$-coalescent condensed $\tu{\hat{L}}$-module and that $\tu{\hat{K}'}$ is a $\big( (i' \of \psi'_*)(f) \cdot (i' \of \psi'_*)(g)^{-1} \big) /1$-coalescent condensed $\tu{\hat{K}'}$-module. Since $(j \of \tau_* \of \iota_*)(g)$ is invertible in $\hat{L}$, it follows that $\tu{\hat{L}}$ is a $\big( (j \of \tau_* \of \iota_*)(f) \big) / \big( (j \of \tau_* \of \iota_*)(g) \big)$-coalescent condensed $\tu{\hat{L}}$-module. This shows that $\tu{\hat{L}}$ is a $f/g$-coalescent condensed $R$-algebra via $R \xto{\iota} R_{\approx T} \xto{\tau} \tu{L} \xto{\tu{j} \;} \tu{\hat{L}}$. Similarly, since $(i' \of \psi'_*)(g)$ is invertible in $\hat{K}'$, it follows that $\tu{\hat{K}'}$ is a $\big( (i' \of \psi'_*)(f) \big) / \big( (i' \of \psi'_*)(g) \big)$-coalescent condensed $\tu{\hat{K}'}$-module. This shows that $\tu{\hat{K}'}$ is a $f/g$-coalescent condensed $R$-algebra via $R \xto{\psi'} \tu{K'} \xto{\tu{i'} \;} \tu{\hat{K}'}$. This completes the proof that $\tu{\hat{L}}$ is a $T$-coalescent condensed $R$-algebra via $R \xto{\iota} R_{\approx T} \xto{\tau} \tu{L} \xto{\tu{j} \;} \tu{\hat{L}}$ and that $\tu{\hat{K}'}$ is a $T$-coalescent condensed $R$-algebra via $R \xto{\psi'} \tu{K'} \xto{\tu{i'} \;} \tu{\hat{K}'}$.

Since $\tu{\hat{K}'}$ is a $T$-coalescent condensed $R$-algebra via $R \xto{\psi'} \tu{K'} \xto{\tu{i'} \;} \tu{\hat{K}'}$, there exists a unique homomorphism $\sigma : R_{\approx T} \to \tu{\hat{K}'}$ of condensed rings such that the following diagram is commutative.
\[\begin{tikzcd}
	R & {\tu{K'}} & {\tu{\hat{K}'}} \\
	{R_{\approx T}}
	\arrow["{\psi'}", from=1-1, to=1-2]
	\arrow["\iota"', from=1-1, to=2-1]
	\arrow["{\tu{i'}}", from=1-2, to=1-3]
	\arrow["\sigma"', from=2-1, to=1-3]
\end{tikzcd}\]
Then the following diagram is commutative in $\ub{CRing}$.
\[\begin{tikzcd}
	{R_{\approx T}} \\
	R & {\tu{K'}} & {\tu{\hat{K}'}} \\
	{R_{\approx T}} & {\tu{L}} & {\tu{\hat{L}}}
	\arrow["\sigma", from=1-1, to=2-3]
	\arrow["\iota", from=2-1, to=1-1]
	\arrow["{\psi'}"', from=2-1, to=2-2]
	\arrow["\iota"', from=2-1, to=3-1]
	\arrow["{\tu{i'}}"', from=2-2, to=2-3]
	\arrow["{\tu{k'}}", from=2-2, to=3-2]
	\arrow["{\tu{\hat{k}'}}", from=2-3, to=3-3]
	\arrow["\tau"', from=3-1, to=3-2]
	\arrow["{\tu{j}}"', from=3-2, to=3-3]
\end{tikzcd}\]
Since $\tu{\hat{L}}$ is a $T$-coalescent condensed $R$-algebra via $R \xto{\iota} R_{\approx T} \xto{\tau} \tu{L} \xto{\tu{j} \;} \tu{\hat{L}}$, the universal property of the $T$-coalescence $\iota : R \to R_{\approx T}$ shows that the following diagram is commutative in $\ub{CRing}$.
\[\begin{tikzcd}
	&& {\tu{\hat{K}'}} \\
	{R_{\approx T}} & {\tu{L}} & {\tu{\hat{L}}}
	\arrow["{\tu{\hat{k}'}}", from=1-3, to=2-3]
	\arrow["\sigma", from=2-1, to=1-3]
	\arrow["\tau"', from=2-1, to=2-2]
	\arrow["{\tu{j}}"', from=2-2, to=2-3]
\end{tikzcd}\]
Thus we obtain the following commutative diagram in $\ub{CRing}$.
\[\begin{tikzcd}
	&& {\hat{K}'} \\
	{R_{\approx T}(*)} & L & {\hat{L}}
	\arrow["{\hat{k}'}", from=1-3, to=2-3]
	\arrow["{\sigma_*}", from=2-1, to=1-3]
	\arrow["{\tau_*}"', from=2-1, to=2-2]
	\arrow["j"', from=2-2, to=2-3]
\end{tikzcd}\]
It follows that $\hat{K}'$ contains the image of the map $R_{\approx T}(*) \xto{\tau_*} L \xto{j} \hat{L}$. However, the field $L$ is generated by the image of the map $R_{\approx T}(*) \xto{\tau_*} L$ and $\hat{K}'$ is a subfield of $\hat{L}$. Therefore $\hat{K}'$ contains the image of $L \xto{j} \hat{L}$. Since $L \xto{j} \hat{L}$ is the Hausdorff completion of $L$, the image $j(L)$ is dense in $\hat{L}$. By definition, $\hat{K}'$ is closed in $\hat{L}$. Therefore we conclude that $\hat{K}' = \hat{L}$.

Recall that $\hat{K}'$ was defined to be the closure of $j(k'(K')) = j(k(K))$ in $\hat{L}$. Since $\hat{K}' = \hat{L}$, (2)(c) of \cref{lem:restriction of valuations on fields} shows that the valuation topology on $K$ associated $v$ is equal to the initial topology on $K$ induced by the map $k : K \to L$. In other words, we have $K = K'$ as topological rings. Then we have $\hat{K} = \hat{K}' = \hat{L}$ as topological rings, and the following diagram is commutative in $\ub{CRing}$.
\[\begin{tikzcd}
	R && {\tu{K}} & {\tu{\hat{K}}} \\
	&& {\tu{K'}} & {\tu{\hat{K}'}} \\
	{R_{\approx T}} && {\tu{L}} & {\tu{\hat{L}}}
	\arrow["\psi", from=1-1, to=1-3]
	\arrow["{\psi'}"', from=1-1, to=2-3]
	\arrow["\iota"', from=1-1, to=3-1]
	\arrow["{\tu{i}}", from=1-3, to=1-4]
	\arrow["{\tu{c}}"', equals, from=1-3, to=2-3]
	\arrow[equals, from=1-4, to=2-4]
	\arrow["{i'}", from=2-3, to=2-4]
	\arrow["{\tu{k'}}"', from=2-3, to=3-3]
	\arrow[equals, from=2-4, to=3-4]
	\arrow["\tau"', from=3-1, to=3-3]
	\arrow["{\tu{j}}"', from=3-3, to=3-4]
\end{tikzcd}\]
By the uniqueness of $\rho : R_{\approx T} \to \tu{\hat{K}}$, the homomorphim $R_{\approx T} \xto{\tau} \tu{L} \xto{\tu{j} \;} \tu{\hat{L}}$ must be equal to $\rho : R_{\approx T} \to \tu{\hat{K}}$. Moreover, both $\hat{v}$ and $\hat{w}'$ are valuations on $\tu{\hat{K}} = \tu{\hat{L}}$, and the following equality holds.
\begin{equation}
\psi^{-1} \big( \big( \tu{i} \big)^{-1} (\hat{w}') \big)
= \iota^{-1} \big( \tau^{-1} \big( \big( \tu{j} \big)^{-1} (\hat{w}') \big) \big) = v
= \psi^{-1} \big( \big( \tu{i} \big)^{-1} (\hat{v}) \big) .
\end{equation}
Since the field $K$ is generated by the image of the map $\psi_* : R(*) \to K$, we conclude that
\begin{equation}
\big( \tu{i} \big)^{-1} (\hat{w}') = \big( \tu{i} \big)^{-1} (\hat{v}) .
\end{equation}
Since both $\hat{v}$ and $\hat{w}'$ are continuous valuations on $\hat{K} = \hat{L}$ and since the map $i : K \to \hat{K}$ has a dense image, we conclude that
\begin{equation}
\hat{w'} = \hat{v} .
\end{equation}
Consequently,
\begin{equation}
w' =  \tau^{-1} \big( \big( \tu{j} \big)^{-1} (\hat{w}') \big)
= \rho^{-1}(\hat{w}') = \rho^{-1}(\hat{v}) = w.
\end{equation}
This completes the proof.
\end{proof}

\begin{prop} \label{prop:construction of coalescence of vcrings}
Let $(R,v)$ be a valued condensed ring. Let $\N$ be the set of natural numbers ordered by the usual order, which we consider as a small category. By induction, we define a functor 
\[\begin{tikzcd}[row sep=tiny]
	\N & {\ub{Ring}} \\
	n & {(R_n,v_n)} \\
	{(n \leq m)} & {\tau_{n,m}}
	\arrow["D", from=1-1, to=1-2]
	\arrow[maps to, from=2-1, to=2-2]
	\arrow[maps to, from=3-1, to=3-2]
\end{tikzcd}\]
as follows.
\begin{itemize}
\item
We define $(R_0,v_0) := (R,v)$.

\item
Let $n \in \N$ with $n \geq 1$. Suppose we have defined $(R_{n-1}, v_{n-1}) \in |\ub{VCRing}|$. Let $T_{n-1}$ be the set of all $(f,g) \in R_{n-1}(*) \times R_{n-1}(*)$ such that $|f|_{v_{n-1}} \leq |g|_{v_{n-1}} \neq 0$. Let $\tau_{n-1,n} : R_{n-1} \to R_n$ be the $T_{n-1}$-coalescence of the condensed $R_{n-1}$-algebra $R_{n-1}$. By \cref{lem:extension of continuous valuation along coalescence}, there exists a unique continuous valuation $v_n$ on $R_n$ such that $(\tau_{n-1,n})^{-1}(v_n) = v_{n-1}$. This completes the definition of $(R_n,v_n)$ and $\tau_{n-1,n}$.
\end{itemize}
Let $F : \ub{VCRing} \to \ub{CRing}$ be the forgetful functor $(R,v) \mapsto R$. Let $(R_{\infty}, (R_n \xto{\rho_n} R_{\infty})_{n \in \N})$ be the colimit of the functor $F \of D : \N \to \ub{CRing}$. Then the following hold.
\begin{enumerate}
\item
There exists a unique continuous valuation $v_{\infty}$ on $R_{\infty}$ such that $(\rho_0)^{-1}(v_{\infty}) = v$. Moreover, this valuation $v_{\infty}$ satisfies $(\rho_n)^{-1}(v_{\infty}) = v_n$ for every $n \in \N$.

\item
The valued condensed ring $(R_{\infty}, v_{\infty})$ is an object of $\ub{VCRing}_c$.

\item
The pair $( (R_{\infty}, v_{\infty}) , (R,v) \xto{\rho_0} (R_{\infty}, v_{\infty}) )$ is a reflection of $(R,v) \in |\ub{VCRing}|$ along the inclusion functor $\ub{VCRing}_c \mon \ub{VCRing}$.
\end{enumerate}
\end{prop}

\begin{proof} ~
\begin{enumerate}
\item
\cref{lem:extension of continuous valuation along coalescence} and induction shows that for every $n,m \in \N$ with $n \leq m$, the map $(\tau_{n,m})_* : R_n(*) \to R_m(*)$ induces an isomorphism of the value group of the valuation $v_n$ onto that of $v_m$. Then \cref{prop:extension of continuous valuation to filtered colimits} shows that there exists a unique continuous valuation $v_{\infty}$ on $R_{\infty}$ such that $(\rho_n)^{-1}(v_{\infty}) = v_n$ for every $n \in \N$.

Suppose that $v'$ is a continuous valuation on $R_{\infty}$ such that $(\rho_0)^{-1}(v') = v$. We prove that $v' = v_{\infty}$. By the definition of $v_{\infty}$, it suffices to show that $(\rho_n)^{-1}(v') = v_n$ holds for every $n \in \N$. We use induction on $n$. The case $n=0$ follows directly from the assumption. Let $n \in \N$ with $n \geq 1$. Suppose that $(\rho_{n-1})^{-1}(v') = v_{n-1}$. Then the commutativity of the diagram
\[\begin{tikzcd}
	{R_n} & {R_{\infty}} \\
	{R_{n-1}}
	\arrow["{\rho_n}", from=1-1, to=1-2]
	\arrow["{\tau_{n-1,n}}", from=2-1, to=1-1]
	\arrow["{\rho_{n-1}}"', from=2-1, to=1-2]
\end{tikzcd}\]
shows that $(\tau_{n-1,n})^{-1} \big( (\rho_n)^{-1}(v') \big) = (\rho_{n-1})^{-1}(v') = v_{n-1}$. By the definition of $v_n$, we conclude that $(\rho_n)^{-1}(v') = v_n$.

\item
Let $f,g \in R_{\infty}(*)$ and suppose $|f|_{v_{\infty}} \leq |g|_{v_{\infty}} \neq 0$. Since the category $\N$ is filtered, \cref{prop:evaluation of CRing preserves limits and filtered colimits} shows that $(R_{\infty}(*), ( R_n(*) \xto{(\rho_n)_*} R_{\infty}(*) )_{n \in \N})$ is the colimit of the functor
\[\begin{tikzcd}[row sep=tiny]
	\N & {\ub{Ring}} \\
	n & {R_n(*)} \\
	{(n \leq m)} & {(\tau_{n,m})_*}
	\arrow[from=1-1, to=1-2]
	\arrow[maps to, from=2-1, to=2-2]
	\arrow[maps to, from=3-1, to=3-2]
\end{tikzcd}\]
Since the category $\N$ is filtered, there exist an $N \in \N$ and $f_N, g_N \in R_N(*)$ such that $(\rho_N)_*(f_N) = f$ and $(\rho_N)_*(g_N) = g$. If $n \in \N$ and $n \geq N+1$, then we have 
\begin{align}
(\rho_{n-1})_* ( (\tau_{N,n-1})_*(f_N) ) & = (\rho_N)_*(f_N) = f \: ; \\
(\rho_{n-1})_* ( (\tau_{N,n-1})_*(g_N) ) & = (\rho_N)_*(g_N) = g .
\end{align}
Moreover, (1) shows that 
\begin{equation}
(\rho_{n-1})^{-1}(v_{\infty}) = v_{n-1} .
\end{equation}
Since $|f|_{v_{\infty}} \leq |g|_{v_{\infty}} \neq 0$ by assumption, we have
\begin{equation}
| (\tau_{N,n-1})_*(f_N) |_{v_{n-1}} \leq | (\tau_{N,n-1})_*(g_N) |_{v_{n-1}} \neq 0 .
\end{equation}
Then the construction of $D$ shows that $R_n$ is a $(\tau_{N,n-1})_*(f_N) / (\tau_{N,n-1})_*(g_N)$-coalescent condensed $R_{n-1}$-algebra via $\tau_{n-1,n} : R_{n-1} \to R_n$. It follows that $R_n$ is an $f_N / g_N$-coalescent condensed $R_N$-algebra via $\tau_{N,n} : R_N \to R_n$. On the other hand, if $R_{\infty}$ is considered as a condensed $R_N$-algebra via $\rho_N : R_N \to R_{\infty}$, then $R_{\infty}$ is the colimit of the functor
\[\begin{tikzcd}[ampersand replacement=\&,row sep=tiny]
	{\N_{\geq N+1}} \& {\ub{CAlg}_{R_N}} \\
	n \& {\left( \begin{gathered}
	R_n \text{ viewed as a condensed} \\
	R_N \text{-algebra via } \tau_{N,n} : R_N \to R_n
	\end{gathered} \right)} \\
	{(n \leq m)} \& {\tau_{n,m}}
	\arrow[from=1-1, to=1-2]
	\arrow[maps to, from=2-1, to=2-2]
	\arrow[maps to, from=3-1, to=3-2]
\end{tikzcd}\]
where $\N_{\geq N+1}$ denotes the full subcategory of $\N$ consisting of all $n \in \N$ with $n \geq N+1$. By \cref{cor:coalescent modules are closed under limits and colimits}, we conclude that $R_{\infty}$ is an $f_N / g_N$-coalescent condensed $R_N$-algebra via $\rho_N : R_N \to R_{\infty}$. Since $(\rho_N)_*(f_N) = f$ and $(\rho_N)_*(g_N) = g$, it follows that $R_{\infty}$ is an $f/g$-coalescent condensed $R_{\infty}$-algebra.

\item
Suppose that $(S,w)$ is an object of $\ub{VCRing}_c$ and that $\phi : (R,v) \to (S,w)$ is a homomorphism of valued condensed rings. We prove that there exists a unique homomorphism $\phi_{\infty} : (R_{\infty}, v_{\infty}) \to (S,w)$ of valued condensed rings such that the following diagram is commutative.
\[\begin{tikzcd}
	{(R_{\infty},v_{\infty})} & {(S,w)} \\
	{(R,v)}
	\arrow["{\phi_{\infty}}", from=1-1, to=1-2]
	\arrow["{\rho_0}", from=2-1, to=1-1]
	\arrow["\phi"', from=2-1, to=1-2]
\end{tikzcd}\]

By induction, we construct a family $( (R_n, v_n) \xto{\phi_n} (S,w) )_{n \in \N}$ of homomorphisms of valued condensed rings as follows.
\begin{itemize}
\item
We define $\phi_0 := \phi : (R,v) \to (S,w)$. 

\item
Let $n \in \N$ with $n \geq 1$. Suppose we have defined the homomorphism $\phi_{n-1} : (R_{n-1}, v_{n-1}) \to (S,w)$ of valued condensed rings. If $f,g \in R_{n-1}(*)$ and $|f|_{v_{n-1}} \leq |g|_{v_{n-1}} \neq 0$, then $|(\phi_{n-1})_*(f)|_w \leq |(\phi_{n-1})_*(f)|_w \neq 0$ since $(\phi_{n-1})^{-1}(w) = v_{n-1}$. Since $(S,w)$ is an object of $\ub{VCRing}_c$ by assumption, the condensed $S$-algebra $S$ is $(\phi_{n-1})_*(f) / (\phi_{n-1})_*(g)$-coalescent. It follows that $S$ is a $f/g$-coalescent $R_{n-1}$-algebra via $\phi_{n-1} : R_{n-1} \to S$. By construction of $D$, there exists a unique homomorphism $\phi_n : R_n \to S$ of condensed rings such that the following diagram is commutative.
\[\begin{tikzcd}
	{R_n} & S \\
	{R_{n-1}}
	\arrow["{\phi_n}", from=1-1, to=1-2]
	\arrow["{\tau_{n-1,n}}", from=2-1, to=1-1]
	\arrow["{\phi_{n-1}}"', from=2-1, to=1-2]
\end{tikzcd}\]
Then we have $(\tau_{n-1,n})^{-1} \big( (\phi_n)^{-1}(w) \big) = (\phi_{n-1})^{-1}(w) = v_{n-1}$. By the definition of $v_n$, we have $(\phi_n)^{-1}(w) = v_n$. Therefore $\phi_n$ is a homomorphism $(R_n,v_n) \to (S,w)$ of valued condensed rings. This completes the definition of $\phi_n$.
\end{itemize}

By construction, $(S, (R_n \xto{\phi_n} S)_{n \in \N} )$ is a cocone on the functor $F \of D : \N \to \ub{CRing}$. Since $(R_{\infty}, (R_n \xto{\rho_n} R_{\infty})_{n \in \N})$ is the colimit of the functor $F \of D : \N \to \ub{CRing}$, there exists a unique homomorphism $\phi_{\infty} : R_{\infty} \to S$ of condensed rings such that the following diagram is commutative for every $n \in \N$.
\[\begin{tikzcd}
	{R_{\infty}} & S \\
	{R_n}
	\arrow["{\phi_{\infty}}", from=1-1, to=1-2]
	\arrow["{\rho_n}", from=2-1, to=1-1]
	\arrow["{\phi_n}"', from=2-1, to=1-2]
\end{tikzcd}\]
In particular, the following diagram is commutative.
\[\begin{tikzcd}
	{R_{\infty}} & S \\
	R
	\arrow["{\phi_{\infty}}", from=1-1, to=1-2]
	\arrow["{\rho_0}", from=2-1, to=1-1]
	\arrow["\phi"', from=2-1, to=1-2]
\end{tikzcd}\]
Then we have $(\rho_0)^{-1} \big( (\phi_{\infty})^{-1}(w) \big) = (\phi)^{-1}(w) = v$. By (1), we conclude that $(\phi_{\infty})^{-1}(w) = v_{\infty}$. Therefore $\phi_{\infty}$ is a homomorphism $(R_{\infty}, v_{\infty}) \to (S,w)$ of valued condensed rings.

On the other hand, suppose $\psi : (R_{\infty}, v_{\infty}) \to (S,w)$ is another homomorphism of valued condensed rings such that the following diagram is commutative. 
\[\begin{tikzcd}
	{(R_{\infty},v_{\infty})} & {(S,w)} \\
	{(R,v)}
	\arrow["\psi", from=1-1, to=1-2]
	\arrow["{\rho_0}", from=2-1, to=1-1]
	\arrow["\phi"', from=2-1, to=1-2]
\end{tikzcd}\]
We show that $\psi = \phi_{\infty}$. By the definition of $\phi_{\infty}$, it suffices to show that the following diagram is commutaive for every $n \in \N$.
\[\begin{tikzcd}
	{R_{\infty}} & S \\
	{R_n}
	\arrow["\psi", from=1-1, to=1-2]
	\arrow["{\rho_n}", from=2-1, to=1-1]
	\arrow["{\phi_n}"', from=2-1, to=1-2]
\end{tikzcd}\]
We use induction on $n$. The case $n=0$ follows directly from the assumption. Let $n \in \N$ with $n \geq 1$. Suppose that the diagram
\[\begin{tikzcd}
	{R_{\infty}} & S \\
	{R_{n-1}}
	\arrow["\psi", from=1-1, to=1-2]
	\arrow["{\rho_{n-1}}", from=2-1, to=1-1]
	\arrow["{\phi_{n-1}}"', from=2-1, to=1-2]
\end{tikzcd}\]
is commutative. Then the following diagram is commutative.
\[\begin{tikzcd}
	{R_n} & {R_{\infty}} & S \\
	{R_{n-1}}
	\arrow["{\rho_n}", from=1-1, to=1-2]
	\arrow["\psi", from=1-2, to=1-3]
	\arrow["{\tau_{n-1,n}}", from=2-1, to=1-1]
	\arrow["{\rho_{n-1}}"'{pos=0.7}, from=2-1, to=1-2]
	\arrow["{\phi_{n-1}}"', curve={height=12pt}, from=2-1, to=1-3]
\end{tikzcd}\]
By the definition of $\phi_n : R_n \to S$, we have $\psi \of \rho_n = \phi_n$. In other words, the following diagram is commutaive.
\[\begin{tikzcd}
	{R_{\infty}} & S \\
	{R_n}
	\arrow["\psi", from=1-1, to=1-2]
	\arrow["{\rho_n}", from=2-1, to=1-1]
	\arrow["{\phi_n}"', from=2-1, to=1-2]
\end{tikzcd}\]
This completes the proof.
\end{enumerate}
\end{proof}

\begin{cor}
$\ub{VCRing}_c$ is a reflective full subcategory of $\ub{VCRing}$.
\end{cor}

\begin{df} \label{df:coalescence of valued condensed ring}
Let $(R,v)$ be a valued condensed ring. Let $\big( (\tilde{R}, \tilde{v}) , (R,v) \xto{\iota} (\tilde{R}, \tilde{v}) \big)$ be the reflection of $(R,v) \in |\ub{VCRing}|$ along the inclusion functor $\ub{VCRing}_c \mon \ub{VCRing}$. Then $(\tilde{R}, \tilde{v}) \in |\ub{VCRing}_c|$ is called the \ti{coalescence of the valued condensed ring} $(R,v)$. The homomorphism $(R,v) \xto{\iota} (\tilde{R}, \tilde{v})$ of valued condensed rings is also called the coalescence of the valued condensed ring $(R,v)$.
\end{df}

\subsection{Localization of objects of $\ub{VCRing}_c$}

\begin{prop} \label{prop:localization of objects of vcring c is in vcring c}
Let $(R,v)$ be an object of $\ub{VCRing}_c$. Let $(R,v) \xto{\iota} (\tilde{R}, \tilde{v})$ be the localization of the valued condensed ring $(R,v)$. Then $(\tilde{R}, \tilde{v})$ is an object of $\ub{VCRing}_c$.
\end{prop}

\begin{proof}
Let $\fk{p}$ be the support of $v$. By \cref{prop:construction of localization of vcrings}, we may assume that $\tilde{R} = R_{\fk{p}}$ and that the homomorphism $R \xto{\iota} \tilde{R}$ is equal to the canonical homomorphism $R \to R_{\fk{p}}$ into the localization $R_{\fk{p}}$. Let $f,g \in R_{\fk{p}}(*)$ with $|f|_{\tilde{v}} \leq |g|_{\tilde{v}} \neq 0$. We show that the condensed $R_{\fk{p}}$-algebra $R_{\fk{p}}$ is $f/g$-coalescent. Since $R_{\fk{p}}(*) = \big( R(*) \big)_{\fk{p}}$, there exist elements $f_0, g_0 \in R(*)$ and $h_0 \in R(*) \sm \fk{p}$ such that $f = \iota_*(f_0) \cdot \iota_*(h_0)^{-1}$ and $g = \iota_*(g_0) \cdot \iota_*(h_0)^{-1}$. Then $|\iota_*(h_0)|_{\tilde{v}} \neq 0$ and the assumption $|f|_{\tilde{v}} \leq |g|_{\tilde{v}} \neq 0$ implies that $|\iota_*(f_0)|_{\tilde{v}} \leq |\iota_*(g_0)|_{\tilde{v}} \neq 0$. Since $v = \iota^{-1}(\tilde{v})$, we have $|f_0|_v \leq |g_0|_v \neq 0$. Since $(R,v)$ is an object of $\ub{VCRing}_c$, the condensed $R$-algebra $R$ is $f_0/g_0$-coalescent. Then \cref{prop:localization and monoid algebra of coalescent algebras} shows that $R_{\fk{p}}$ is an $f_0/g_0$-coalescent condensed $R$-algebra via $R \xto{\iota} R_{\fk{p}}$. It follows that the condensed $R_{\fk{p}}$-algebra $R_{\fk{p}}$ is $\iota_*(f_0)/\iota_*(g_0)$-coalescent. Since the element $\iota_*(h_0) \in R_{\fk{p}}(*)$ is invertible in $R_{\fk{p}}(*)$, we conclude that the condensed $R_{\fk{p}}$-algebra $R_{\fk{p}}$ is $\big( \iota_*(f_0) \cdot \iota_*(h_0)^{-1} \big)/\big( \iota_*(g_0) \cdot \iota_*(h_0)^{-1} \big)$-coalescent. Therefore the condensed $R_{\fk{p}}$-algebra $R_{\fk{p}}$ is $f/g$-coalescent.
\end{proof}

\begin{cor} \label{cor:coalescence followed by localization of vcrings}
Let $(R,v)$ be a valued condensed ring. Let $(R,v) \xto{\iota} (R_c, v_c)$ be the coalescence of the valued condensed ring $(R,v)$. Let $(R_c, v_c) \xto{\lambda} (R_{c,l}, v_{c,l})$ be the localization of the valued condensed ring $(R_c, v_c)$. Then the pair $( (R_{c,l}, v_{c,l}), (R,v) \xto{\lambda \of \iota} (R_{c,l}, v_{c,l}) )$ is the reflection of $(R,v) \in |\ub{VCRing}|$ along the inclusion functor $\ub{VCRing}_l \cap \ub{VCRing}_c \mon \ub{VCRing}$.
\end{cor}

\begin{proof}
By \cref{prop:localization of objects of vcring c is in vcring c}, the valued condensed ring $(R_{c,l}, v_{c,l})$ is an object of $\ub{VCRing}_l \cap \ub{VCRing}_c$. On the other hand, for each object $(S,w) \in |\ub{VCRing}_l \cap \ub{VCRing}_c|$, the map
\[\begin{tikzcd}
	{\ub{VCRing} \Big( (R_c, v_c), (S,w) \Big)} & {\ub{VCRing} \Big( (R,v), (S,w) \Big) ,} & \psi & {\psi \of \iota}
	\arrow[from=1-1, to=1-2]
	\arrow[maps to, from=1-3, to=1-4]
\end{tikzcd}\]
is bijective by the universality of the coalescence $(R,v) \xto{\iota} (R_c, v_c)$ of the valued condensed ring $(R,v)$. Furthermore, the map
\[\begin{tikzcd}
	{\ub{VCRing} \Big( (R_{c,l}, v_{c,l}), (S,w) \Big)} & {\ub{VCRing} \Big( (R_c, v_c), (S,w) \Big) ,} & \rho & {\rho \of \lambda}
	\arrow[from=1-1, to=1-2]
	\arrow[maps to, from=1-3, to=1-4]
\end{tikzcd}\]
is bijective by the universality of the localization $(R_c, v_c) \xto{\lambda} (R_{c,l}, v_{c,l})$ of the valued condensed ring $(R_c, v_c)$. Therefore the map
\[\begin{tikzcd}
	{\ub{VCRing} \Big( (R_{c,l}, v_{c,l}), (S,w) \Big)} & {\ub{VCRing} \Big( (R,v), (S,w) \Big) ,} & \rho & {\rho \of \lambda \of \iota}
	\arrow[from=1-1, to=1-2]
	\arrow[maps to, from=1-3, to=1-4]
\end{tikzcd}\]
is bijective. This completes the proof.
\end{proof}

\begin{cor} \;
\begin{enumerate}
\item
$\ub{VCRing}_l \cap \ub{VCRing}_c$ is a reflective full subcategory of $\ub{VCRing}_c$.

\item
$\ub{VCRing}_l \cap \ub{VCRing}_c$ is a reflective full subcategory of $\ub{VCRing}$.
\end{enumerate}
\end{cor}

\begin{proof}~
\begin{enumerate}
\item
By \cref{prop:localization of objects of vcring c is in vcring c}, the left adjoint $\ub{VCRing} \to \ub{VCRing}_l$ of the inclusion functor $\ub{VCRing}_l \mon \ub{VCRing}$ restricts a functor $\ub{VCRing}_c \to \ub{VCRing}_l \cap \ub{VCRing}_c$. This functor is left adjoint to the inclusion functor $\ub{VCRing}_l \cap \ub{VCRing}_c \to \ub{VCRing}_c$.

\item
This follows from \cref{cor:coalescence followed by localization of vcrings}.
\end{enumerate}
\end{proof}

\section{Condensed sheaves} \label{sec:Condensed sheaves}

In this section, we study the notion of sheaves of condensed rings, which is an analogue of sheaves of commutative unital rings in the context of condensed mathematics. Although one can consider sheaves of condensed sets, condensed abelian groups, condensed modules, and so on, we concentrate on sheaves of condensed rings here because we do not need sheaves of objects other than condensed rings in this paper.

\subsection{Definition}

\begin{df}
Let $X$ be a topological space. Let $\cat{T}$ be the set of all open subsets in $X$ ordered by inclusion. We consider $\cat{T}$ as a small category.
\begin{enumerate}
\item
The category $\ub{CPSh}(X)$ is defined to be the category of functors $\cat{T}^{\op} \to \ub{CRing}$ and natural transformation between them. Objects of $\ub{CPSh}(X)$ are called \ti{presheaves of condensed rings} on $X$, and morphisms in $\ub{CPSh}(X)$ are called \ti{morphisms of presheaves of condensed rings} on $X$.

\item
Let $F$ be a presheaf of condensed rings on $X$. Let $U,V$ be open subsets of $X$ such that $V \sub U$. Then we have a homomorphism of condensed rings $F(V \sub U) : F(U) \to F(V)$, which we call the \ti{restriction} from $F(U)$ to $F(V)$. We often write just $F(U) \xto{\res} F(V)$ for this.
\end{enumerate}
\end{df}

\begin{df}
Let $X$ be a topological space. Let $\cat{T}$ be the set of all open subsets in $X$ ordered by inclusion. We consider $\cat{T}$ as a small category.
\begin{enumerate}
\item
Let $F$ be a presheaf of condensed rings on $X$. Let $U$ be an open subset of $X$. Let $(U_i)_{i \in I}$ be an open covering of $U$. Consider the following homomorphisms of condensed rings.
\begin{enumerate}
\item
The unique homomorphism $\alpha : F(U) \to \prod_{i \in I} F(U_i)$ of condensed rings such that the diagram
\[\begin{tikzcd}
	{F(U)} & {\displaystyle \prod_{i \in I} F(U_i)} \\
	& {F(U_i)}
	\arrow["\alpha", from=1-1, to=1-2]
	\arrow["\res"', from=1-1, to=2-2]
	\arrow["{i\text{-th projection}}", from=1-2, to=2-2]
\end{tikzcd}\]
is commutative for all $i \in I$.

\item
For each $k \in \{1,2\}$, the unique homomorphism $\beta_k : \prod_{i \in I} F(U_i) \to \prod_{(i_1,i_2) \in I \times I} F(U_{i_1} \cap U_{i_2})$ of condensed rings such that the diagram
\[\begin{tikzcd}
	{\displaystyle \prod_{i \in I} F(U_i)} & {\displaystyle
	\prod_{(i_1,i_2) \in I \times I} F(U_{i_1} \cap U_{i_2})} \\
	{F(U_{i_k})} & {F(U_{i_1} \cap U_{i_2})}
	\arrow["{\beta_k}", from=1-1, to=1-2]
	\arrow["{i_k \text{-th projection}}"', from=1-1, to=2-1]
	\arrow["{(i_1 , i_2)\text{-th projection}}", from=1-2, to=2-2]
	\arrow["\res"', from=2-1, to=2-2]
\end{tikzcd}\]
is commutative for all $(i_1,i_2) \in I \times I$.
\end{enumerate}

The \ti{sheaf condition with respect to the open covering} $(U_i)_{i \in I}$ \ti{of} $U$ is the condition that the diagram
\[\begin{tikzcd}
	{F(U)} & {\prod_{i \in I} F(U_i)} & {\prod_{(i_1,i_2) \in I\times I} F(U_{i_1} \cap U_{i_2})}
	\arrow["\alpha", from=1-1, to=1-2]
	\arrow["{\beta_2}"', shift right, from=1-2, to=1-3]
	\arrow["{\beta_1}", shift left, from=1-2, to=1-3]
\end{tikzcd}\]
be exact in $\ub{CRing}$, i.e., the homomorphism $\alpha$ be an equalizer of the homomorphisms $\beta_1 , \beta_2$ in $\ub{CRing}$.

\item
The category $\ub{CSh}(X)$ is defined to be the full subcategory of $\ub{CPSh}(X)$ consisting of all presheaves of condensed rings on $X$ which satisfy the sheaf condition with respect to all open coverings of all open subsets of $X$. Objects of $\ub{CSh}(X)$ are called \ti{sheaves of condensed rings} on $X$, and morphisms in $\ub{CSh}(X)$ are called \ti{morphisms of sheaves of condensed rings} on $X$.
\end{enumerate}
\end{df}

\begin{nt} \label{nt:passing from condesed sheaves to ordinary sheaves}
Let $X$ be a topological space. Let $\cat{T}$ be the set of all open subsets in $X$ ordered by inclusion, which we consider as a small category. Let $S$ be a light profinite set.
\begin{enumerate}
\item
Let $F$ be a presheaf of condensed rings on $X$. We write $F^S$ for the presheaf of commutative unital rings on $X$ defined as follows.
\[\begin{tikzcd}[row sep=tiny]
	{\cat{T}^{\op}} & {\ub{Ring}} \\
	U & {F(U)(S)} & {(\text{on objects})} \\
	{(V \sub U)} & {F(V \sub U)_S} & {(\text{on morphisms})}
	\arrow["{F^S}", from=1-1, to=1-2]
	\arrow[maps to, from=2-1, to=2-2]
	\arrow[maps to, from=3-1, to=3-2]
\end{tikzcd}\]

\item
Let $\alpha : F \to G$ be a morphism of presheaves of condensed rings on $X$. We write $\alpha^S : F^S \to G^S$ for the morphism of presheaves of commutative unital rings on $X$ defined as follows.
\begin{equation}
\alpha^S := \left( F(U)(S) \xto{\alpha_{U,S}} G(U)(S) \right)_{U \in |\cat{T}^{\op}|}
\end{equation}

\item
We obtain a functor
\[\begin{tikzcd}[row sep=tiny]
	{\ub{CPSh}(X)} & {\ub{PSh}(X)} \\
	F & {F^S} & {(\text{on objects})} \\
	\alpha & {\alpha^S} & {(\text{on morphisms}),}
	\arrow[from=1-1, to=1-2]
	\arrow[maps to, from=2-1, to=2-2]
	\arrow[maps to, from=3-1, to=3-2]
\end{tikzcd}\]
where $\ub{PSh}(X)$ denotes the category of presheaves of commutative unital rings on $X$. 
\end{enumerate}
\end{nt}

\begin{rem} \label{rem:sheaf property in terms of F S}
Let $X$ be a topological space. Let $F$ be a presheaf of condensed rings on $X$. Then $F$ is a sheaf of condensed rings on $X$ if and only if $F^S$ is a sheaf of commutative unital rings on $X$ for every light profinite set $S$.
\end{rem}

\begin{nt}
Let $X$ be a topological space. Let $F$ be a presheaf of condensed rings on $X$. If $U,V$ are open subsets of $X$ with $V \sub U$ and $c: S' \to S$ is a morphism in $\ub{Prof}$, then the following diagram is commutative.
\[\begin{tikzcd}
	{F(U)(S)} && {F(U)(S')} \\
	{F(V)(S)} && {F(V)(S')}
	\arrow["{F(U)(c)}", from=1-1, to=1-3]
	\arrow["{F(V \sub U)_S}"', from=1-1, to=2-1]
	\arrow["{F(V \sub U)_{S'}}", from=1-3, to=2-3]
	\arrow["{F(V)(c)}"', from=2-1, to=2-3]
\end{tikzcd}\]
\begin{enumerate}
\item 
We define $F(V \sub U, c) : F(U)(S) \to F(V)(S')$ by
\begin{equation}
F(V \sub U, c) := F(V \sub U)_{S'} \of F(U)(c) = F(V)(c) \of F(V \sub U)_S .
\end{equation}

\item
For $t \in F(U)(S)$, we write
\begin{equation}
t|_{V,c} := F(V \sub U, c) (t) .
\end{equation}

\item
If $c$ is equal to the identity $S \to S$, then we also write
\begin{equation}
t|_V := F(V \sub U, \id_S) (t) = F(V \sub U)_S \, (t) .
\end{equation}
for $t \in F(U)(S)$.
\end{enumerate}
\end{nt}

\subsection{Subsheaves}

\begin{df}
Let $X$ be a topological space. Let $F$ be a presheaf of condensed rings on $X$.
\begin{enumerate}
\item
A \ti{subpresheaf of condensed rings} of $F$ is a presheaf $G$ of condensed rings on $X$ which satisfies the following conditions.
\begin{enumerate}
\item
For every open subset $U$ of $X$, the condensed ring $G(U)$ is a condensed subring of $F(U)$.

\item
For every open subsets $U,V$ of $X$ with $V \sub U$, the following diagram is commutative.
\[\begin{tikzcd}[column sep=large]
	{G(U)} & {F(U)} \\
	{G(V)} & {F(V)}
	\arrow["{\text{inclusion}}", hook, from=1-1, to=1-2]
	\arrow["{G(V \sub U)}"', from=1-1, to=2-1]
	\arrow["{F(V \sub U)}", from=1-2, to=2-2]
	\arrow["{\text{inclusion}}"', hook, from=2-1, to=2-2]
\end{tikzcd}\]
\end{enumerate}
In this case, the family of inclusions $G(U) \mon F(U) \; (U \sub X \text{ open})$ defines a morphism $G \to F$ of presheaves of condensed rings on $X$, which we call the \ti{inclusion} of $G$ into $F$.

\item
A \ti{subsheaf of condensed rings} of $F$ is a subpresheaf $G$ of condensed rings of $F$ such that $G$ itself is a sheaf of condensed rings on $X$. 
\end{enumerate}
\end{df}

\begin{rem}
Let $X$ be a topological space. Let $F$ be a presheaf of condensed rings on $X$. To give a subpresheaf $G$ of condensed rings of $F$ is equivalent to give a family of condensed subrings $G(U) \sub F(U)$, where $U$ runs through all open subsets of $X$, with the following property: For every open subsets $U,V$ of $X$ with $V \sub U$ and every light profinite set $S$, we have
\begin{equation}
F(V \sub U)_S \Big( G(U)(S) \Big) \sub G(V)(S) .
\end{equation}
\end{rem}

\subsection{Stalks}

\begin{df}
Let $X$ be a topological space. Let $x \in X$. Let $\cat{T}(x)$ be the set of all open neighbourhoods of $x$ in $X$, ordered by inclusion. We consider $\cat{T}(x)$ as a small category.
\begin{enumerate}
\item
Let $F$ be a presheaf of condensed rings on $X$. The \ti{stalk} $F_x$ of $F$ at $x \in X$ is defined to be the colimit
\begin{equation}
\underset{U \in \cat{T}(x)^{\op}}{\colim} \; F(U) .
\end{equation}
Here the colimit is taken in $\ub{CRing}$.

\item
Let $\alpha : F \to G$ be a morphism of presheaves of condensed rings on $X$. Then the homomorphisms of condensed rings
\[\begin{tikzcd}
	{F(U)} & {G(U)} & {( U \in |\cat{T}(x)^{\op}| )}
	\arrow["{\alpha_U}", from=1-1, to=1-2]
\end{tikzcd}\]
induce a homomorphism of condensed rings
\[\begin{tikzcd}
	{F_x} & {G_x ,}
	\arrow[from=1-1, to=1-2]
\end{tikzcd}\]
which we write $\alpha_x$.

\item
We obtain a functor
\[\begin{tikzcd}
	{\ub{CPSh}(X)} & {\ub{CRing},} & F & {F_x .}
	\arrow[from=1-1, to=1-2]
	\arrow[maps to, from=1-3, to=1-4]
\end{tikzcd}\]
\end{enumerate}
\end{df}

\begin{nt}
Let $X$ be a topological space. Let $U$ be an open subset of $X$. Let $F$ be a presheaf of condensed rings on $X$. If $x \in U$, we have a canonical homomorphism $c : F(U) \to F_x$ of condensed rings. If $S$ is a light profinite set, then we have a homomorphism $c_S : F(U)(S) \to F_x(S)$ of rings. For $t \in F(U)(S)$, the image of $t$ under $c_S : F(U)(S) \to F_x(S)$ is denoted by $t_x$.
\end{nt}

\begin{prop} \label{prop:compatibility of stalk and evaluation at S}
Let $X$ be a topological space. Let $x \in X$. Let $\cat{T}(x)$ be the set of all open neighbourhoods of $x$ in $X$, ordered by inclusion. We consider $\cat{T}(x)$ as a small category. Let $S$ be a light profinite set.
\begin{enumerate}
\item
Let $F$ be a presheaf of condensed rings on $X$. For each $U \in \cat{T}(x)$, write $c_U : F(U) \to F_x$ for the canonical homomorphism to the stalk. Then we have
\begin{equation}
F_x (S) = \underset{U \in \cat{T}(x)^{\op}}{\colim} \; F(U)(S) = (F^S)_x ,
\end{equation}
where the colimit is taken in $\ub{Ring}$. Moreover, for each $U \in \cat{T}(x)$, the canonical map $F^S(U) \to (F^S)_x$ coincides with the map $(c_U)_S : F(U)(S) \to F_x(S)$.

\item
Let $\alpha : F \to G$ be a morphism of presheaves of condensed rings on $X$. Then the map $(\alpha^S)_x : (F^S)_x \to (G^S)_x$ coincides with the map $(\alpha_x)_S : F_x(S) \to G_x(S)$.
\end{enumerate}
\end{prop}

\begin{proof}
Since the category $\cat{T}(x)^{\op}$ is filtered, these assertions immediately follow from \cref{prop:evaluation of CRing preserves limits and filtered colimits}.
\end{proof}

\begin{prop} \label{prop:testing equality of sheaf morphisms by stalks}
Let $X$ be a topological space. Let $\alpha, \beta : F \to G$ be two morphisms of sheaves of condensed rings on $X$. Then $\alpha = \beta$ if and only if $\alpha_x = \beta_x$ for every $x \in X$.
\end{prop}

\begin{proof}
Suppose that $\alpha_x = \beta_x$ for every $x \in X$. We prove that $\alpha = \beta$. Let $S$ be any light profinite set. If $x \in X$, then \cref{prop:compatibility of stalk and evaluation at S} shows that the maps $(\alpha^S)_x , (\beta^S)_x : (F^S)_x \to (G^S)_x$ coincide with the maps $(\alpha_x)_S , (\beta_x)_S : F_x(S) \to G_x(S)$, respectively. Since $\alpha_x = \beta_x$ by hypothesis, we conclude that $(\alpha^S)_x = (\beta^S)_x$. Thus we have proved that $(\alpha^S)_x = (\beta^S)_x : (F^S)_x \to (G^S)_x$ for every $x \in X$. On the other hand, by \cref{rem:sheaf property in terms of F S}, $F^S, G^S$ are sheaves of commutative unital rings on $X$. Then the ordinary sheaf theory shows that $\alpha^S = \beta^S : F^S \to G^S$. Consequently,
\begin{equation}
\alpha_{U,S} = (\alpha^S)_U = (\beta^S)_U = \beta_{U,S}
\end{equation}
for every open subset $U$ of $X$ and every light profinite set $S$. It follows that $\alpha = \beta$.
\end{proof}

\begin{prop} \label{prop:testing sheaf isomorphisms by stalks}
Let $X$ be a topological space. Let $\alpha : F \to G$ be a morphism of sheaves of condensed rings on $X$. Then $\alpha : F \to G$ is an isomorphism in $\ub{CSh}(X)$ if and only if $\alpha_x : F_x \to G_x$ is an isomorphism in $\ub{CRing}$ for every $x \in X$.
\end{prop}

\begin{proof}
Suppose that $\alpha : F \to G$ is an isomorphism in $\ub{CSh}(X)$. If $x \in X$, then we have a functor
\[\begin{tikzcd}[row sep=tiny]
	{\ub{CSh}(X)} & {\ub{CRing}} \\
	H & {H_x} & {(\text{on objects})} \\
	\beta & {\beta_x} & {(\text{on morphisms})}
	\arrow[from=1-1, to=1-2]
	\arrow[maps to, from=2-1, to=2-2]
	\arrow[maps to, from=3-1, to=3-2]
\end{tikzcd}\]
It follows that $\alpha_x : F_x \to G_x$ is an isomorphism in $\ub{CRing}$.

Conversely, suppose that $\alpha_x : F_x \to G_x$ is an isomorphism in $\ub{CRing}$ for every $x \in X$. Let $S$ be any light profinite set. If $x \in X$, then \cref{prop:compatibility of stalk and evaluation at S} shows that the map $(\alpha^S)_x : (F^S)_x \to (G^S)_x$ coincides with the map $(\alpha_x)_S : F_x(S) \to G_x(S)$. Since $\alpha_x : F_x \to G_x$ is an isomorphism in $\ub{CRing}$ by assumption, we conclude that $(\alpha^S)_x : (F^S)_x \to (G^S)_x$ is an isomorphism in $\ub{Ring}$. This holds for every $x \in X$. On the other hand, by \cref{rem:sheaf property in terms of F S}, $F^S, G^S$ are sheaves of commutative unital rings on $X$. Then the ordinary sheaf theory shows that $\alpha^S : F^S \to G^S$ is an isomorphism of sheaves of commutative unital rings on $X$. Consequently, for every light profinite set $S$ and every open subset $U$ of $X$, the map $(\alpha^S)_U : F^S(U) \to G^S(U)$ is an isomorphism in $\ub{Ring}$. However, the map $(\alpha^S)_U : F^S(U) \to G^S(U)$ coincides with the map $\alpha_{U,S} : F(U)(S) \to G(U)(S)$. It follows that $\alpha : F \to G$ is an isomorphism in $\ub{CSh}(X)$.
\end{proof}

\subsection{Construction of condensed sheaves}

$\\[2mm]$ In this subsection, we explain a method for constructing a sheaf of condensed rings.

\subsubsection{Settings}

$\\[2mm]$ Throughout this subsection, we fix the following notation.

\begin{nt} \;
\begin{enumerate}
\item $X$ is a topological space. We write $\cat{T}$ for the set of all open subsets in $X$, ordered by inclusion. We consider $\cat{T}$ as a small category.

\item $(I ,\leq)$ is directed set. We also consider $I$ as a small category.

\item $D : I \to \cat{T}^{\op}$ is a functor.

\item $R : I \to \ub{CRing}$ is a functor.
\end{enumerate}
\end{nt}

\begin{nt} Let $x \in X$.
\begin{enumerate}
\item We write $\cat{T}(x)$ for the full subcategory of $\cat{T}$ consisting of all open neighbourhoods of $x$ in $X$.

\item We write $I(x)$ for the full subcategory of $I$ consisting of all $i \in I$ such that $x \in D(i)$.
\end{enumerate}
\end{nt}

Throughout this subsection, we make the following assumption.

\begin{as} \label{as:assumption for constructing condensed sheaves} \;
\begin{enumerate}
\item
The set $\set{D(i)}{i \in I}$ is a basis for the topology of $X$.

\item
For every $i,j \in I$ and every $x \in D(i) \cap D(j)$, there exists a $k \in I$ such that $i \leq k$, $\; j \leq k$ and $x \in D(k)$.
\end{enumerate}
\end{as}

\subsubsection{Statement}

\begin{prop} \label{prop:constructing condensed sheaves}
Under the \cref{as:assumption for constructing condensed sheaves}, there exist a sheaf $G$ of condensed rings on $X$ and a natural transformation $\beta : R \To G \of D$ of functors $I \to \ub{CRing}$ with the following property.
\begin{enumerate}
\item
For each $x \in X$, the homomorphisms of condensed rings
\[\begin{tikzcd}
	{R(i)} & {G(D(i))} & {(i \in I(x))}
	\arrow["{\beta_i}", from=1-1, to=1-2]
\end{tikzcd}\]
induce an isomorphism of condensed rings
\[\begin{tikzcd}
	{\underset{i \in I(x)}{\colim} \; R(i)} & {G_x .}
	\arrow["\sim", from=1-1, to=1-2]
\end{tikzcd}\]

\item
For every sheaf $H$ of condensed rings on $X$ and a natural transformation $\gamma : R \To H \of D$ of functors $I \to \ub{CRing}$, there exists a unique morphism $\delta : G \to H$ of sheaves of condensed rings on $X$ such that the diagram
\[\begin{tikzcd}
	R & {H \of D} \\
	{G \of D}
	\arrow["\gamma", from=1-1, to=1-2]
	\arrow["\beta"', from=1-1, to=2-1]
	\arrow["{\delta * \id_D}"', from=2-1, to=1-2]
\end{tikzcd}\]
is commutative.
\end{enumerate}
\end{prop}

The proof of \cref{prop:constructing condensed sheaves} is the content of the remainder of this subsection.

\subsubsection{Construction of $G$ and $\beta$}

\begin{nt}
Let $i,j \in I$ with $i \leq j$. Let $c : S' \to S$ be a morphism in $\ub{Prof}$. Then the following diagram is commutative.
\[\begin{tikzcd}
	{R(i)(S)} & {R(i)(S')} \\
	{R(j)(S)} & {R(j)(S')}
	\arrow["{R(i)(c)}", from=1-1, to=1-2]
	\arrow["{R(i \leq j)_S}"', from=1-1, to=2-1]
	\arrow["{R(i \leq j)_{S'}}", from=1-2, to=2-2]
	\arrow["{R(j)(c)}"', from=2-1, to=2-2]
\end{tikzcd}\]
We define $R(i \leq j, c) : R(i)(S) \to R(j)(S')$ by
\begin{equation}
R(i \leq j, c) := R(i \leq j)_{S'} \of R(i)(c) = R(j)(c) \of R(i \leq j)_S .
\end{equation}
\end{nt}

\begin{lem} \label{lem:I(x) is filtered}
For each $x \in X$, the preordered set $(I(x), \leq)$ is a nonempty directed set.
\end{lem}

\begin{proof} Suppose that $x \in X$. Then (1) of \cref{as:assumption for constructing condensed sheaves} shows that there exists at least one $i \in I$ such that $x \in D(i)$. Then $i \in I(x)$, which shows that $I(x) \neq \ku$.

Suppose $i,j \in I(x)$. Then $x \in D(i) \cap D(j)$. By (2) of \cref{as:assumption for constructing condensed sheaves}, there exists a $k \in I$ such that $i \leq k$, $\; j \leq k$ and $x \in D(k)$. Then $k \in I(x)$, and we have $i \leq k$ and $j \leq k$.
\end{proof}

\begin{nt} For each $x \in X$, we write
\begin{equation}
R_x := \underset{i \in I(x)}{\colim} \; R(i) ,
\end{equation}
where the colimit is taken in $\ub{CRing}$. For each $i \in I(x)$, we write $\rho_{x,i} : R(i) \to R_x$ for the canonical homomorphism.
\end{nt}

\begin{nt} \;
\begin{enumerate}
\item
For every open subset $U$ of $X$, we define
\begin{equation}
F(U) := \prod_{x \in U} R_x ,
\end{equation}
where the product is taken in $\ub{CRing}$. For each $x \in U$, we write $\pi_{U,x} : F(U) \to R_x$ for the $x$-th projection.

\item
If $U,V$ are open subsets of $X$ such that $V \sub U$, then we define $F(V \sub U) : F(U) \to F(V)$ to be the unique homomorphism of condensed rings such that the diagram
\[\begin{tikzcd}
	{F(U)} \\
	{F(V)} & {R_x}
	\arrow["{F(V \sub U)}"', from=1-1, to=2-1]
	\arrow["{\pi_{U,x}}", from=1-1, to=2-2]
	\arrow["{\pi_{V,x}}"', from=2-1, to=2-2]
\end{tikzcd}\]
is commutative for all $x \in V$.

\item Then we obtain a sheaf $F$ of condensed rings on $X$.
\end{enumerate}
\end{nt}

\begin{prop} \;
\begin{enumerate}
\item
For every $i \in I$, there exists a unique homomorphism $\alpha_i : R(i) \to F(D(i))$ of condensed rings such that the diagram
\[\begin{tikzcd}
	{R(i)} & {F(D(i))} \\
	& {R_x}
	\arrow["{\alpha_i}", from=1-1, to=1-2]
	\arrow["{\rho_{x,i}}"', from=1-1, to=2-2]
	\arrow["{\pi_{D(i), x}}", from=1-2, to=2-2]
\end{tikzcd}\]
is commutative for every $x \in D(i)$. 

\item
The family $\alpha := (\alpha_i)_{i \in I}$ is a natural transformation $R \To F \of D$ of functors $I \to \ub{CRing}$.
\end{enumerate}
\end{prop}

\begin{proof} $\\$
\begin{enumerate}
\item
Let $i \in I$. For each $x \in D(i)$, we have $i \in I(x)$. Therefore we have a canonical homomorphism $\rho_{x,i} : R(i) \to R_x$. Since $F(D(i)) = \prod_{x \in D(i)} R_x$, there exists a unique homomorphism $\alpha_i : R(i) \to F(D(i))$ of condensed rings such that the diagram
\[\begin{tikzcd}
	{R(i)} & {F(D(i))} \\
	& {R_x}
	\arrow["{\alpha_i}", from=1-1, to=1-2]
	\arrow["{\rho_{x,i}}"', from=1-1, to=2-2]
	\arrow["{\pi_{D(i), x}}", from=1-2, to=2-2]
\end{tikzcd}\]
is commutative for every $x \in D(i)$. 

\item
Suppose $i,j \in I$ and that $i \leq j$. Then $D(j) \sub D(i)$. For each $x \in D(j)$, we have $x \in D(j)$ and $x \in D(i)$, and therefore $i,j \in I(x)$. Then the following diagram is commutative.
\[\begin{tikzcd}
	{R(i)} &&&& {F(D(i))} \\
	&& {R_x} && {F(D(j))} \\
	{R(j)} &&&& {F(D(j))}
	\arrow["{\alpha_i}", from=1-1, to=1-5]
	\arrow["{\rho_{x,i}}"', from=1-1, to=2-3]
	\arrow["{R(i \leq j)}"', from=1-1, to=3-1]
	\arrow["{\pi_{D(i),x}}"'{pos=0.7}, curve={height=6pt}, from=1-5, to=2-3]
	\arrow["{F(D(j) \sub D(i))}", from=1-5, to=2-5]
	\arrow["{\pi_{D(j),x}}"'{pos=0.4}, from=2-5, to=2-3]
	\arrow["{\rho_{x,j}}", from=3-1, to=2-3]
	\arrow["{\alpha_j}"', from=3-1, to=3-5]
	\arrow["{\pi_{D(j),x}}", from=3-5, to=2-3]
\end{tikzcd}\]
Since $F(D(j)) = \prod_{x \in D(j)} R_x$, we conclude that the following diagram is commutative.
\[\begin{tikzcd}
	{R(i)} & {F(D(i))} \\
	{R(j)} & {F(D(j))}
	\arrow["{\alpha_i}", from=1-1, to=1-2]
	\arrow["{R(i \leq j)}"', from=1-1, to=2-1]
	\arrow["{F(D(j) \sub D(i))}", from=1-2, to=2-2]
	\arrow["{\alpha_j}"', from=2-1, to=2-2]
\end{tikzcd}\]
This proves that the family $\alpha := (\alpha_i)_{i \in I}$ is a natural transformation $R \To F \of D$ of functors $I \to \ub{CRing}$.
\end{enumerate}
\end{proof}

\begin{nt} \label{nt:definition of G(U)(S)}
For every open subset $U$ of $X$ and every light profinite set $S$, we define $G(U)(S)$ to be the subset of $F(U)(S)$ consisting of all $t \in F(U)(S)$ which satisfies the following condition: For every $x \in U$, there exist an $i \in I$ and a cover $(c_{\lambda} : S_{\lambda} \to S)_{\lambda \in \Lambda}$ of $S$ in $\ub{Prof}$ such that $x \in D(i) \sub U$ and $t|_{D(i),c_{\lambda}} \in \alpha_{i, S_{\lambda}}(R(i)(S_{\lambda}))$ for every $\lambda \in \Lambda$.
\end{nt}

\begin{lem} \label{lem:choosing i and (c lambda) uniformly}
Let $U$ be an open subset of $X$. Let $x \in U$. Let $S$ be a light profinite set. Suppose that $t_1 ,\ldots, t_n \in G(U)(S)$, where $n \in \N$ and $n \geq 1$. Then there exist an $i \in I$ and a cover $(c_{\lambda} : S_{\lambda} \to S)_{\lambda \in \Lambda}$ of $S$ in $\ub{Prof}$ such that $x \in D(i) \sub U$ and $t_p |_{D(i),c_{\lambda}} \in \alpha_{i, S_{\lambda}}(R(i)(S_{\lambda}))$ for every $\lambda \in \Lambda$ and every $1 \leq p \leq n$.
\end{lem}

\begin{proof}
We use the induction on $n$. The case $n=1$ follows from the very definition of $G(U)(S)$. Suppose $n \geq 2$ and that the lemma is true for $n-1$. We apply this assumption to $t_1 ,\ldots, t_{n-1} \in G(U)(S)$. We obtain an $i \in I$ and a cover $(c_{\lambda} : S_{\lambda} \to S)_{\lambda \in \Lambda}$ of $S$ in $\ub{Prof}$ such that $x \in D(i) \sub U$ and $t_p |_{D(i),c_{\lambda}} \in \alpha_{i, S_{\lambda}}(R(i)(S_{\lambda}))$ for every $\lambda \in \Lambda$ and every $1 \leq p \leq n-1$. On the other hand, since $t_n \in G(U)(S)$, there exist an $j \in I$ and a cover $(d_{\mu} : S'_{\mu} \to S)_{\mu \in M}$ of $S$ in $\ub{Prof}$ such that $x \in D(j) \sub U$ and $t_n |_{D(j),d_{\mu}} \in \alpha_{j, S'_{\mu}}(R(j)(S'_{\mu}))$ for every $\mu \in M$.

First we note that $x \in D(i) \cap D(j)$. By (2) of \cref{as:assumption for constructing condensed sheaves}, there exists a $k \in I$ such that $i \leq k$, $\; j \leq k$ and $x \in D(k)$. Then $x \in D(k) \sub D(i) \sub U$.

Next, for each $(\lambda, \mu) \in \Lambda \times M$, we define $S''_{\lambda \mu} := S_{\lambda} \times_{S} S'_{\mu}$, with projections $p_{\lambda \mu} : S''_{\lambda \mu} \to S_{\lambda}$ and $q_{\lambda \mu} : S''_{\lambda \mu} \to S'_{\mu}$. This fiber product is taken in $\ub{Top}$. In other words, the following diagram is a pullback square in $\ub{Top}$.
\[\begin{tikzcd}
	{S''_{\lambda \mu}} & {S'_{\mu}} \\
	{S_{\lambda}} & S
	\arrow["{q_{\lambda \mu}}", from=1-1, to=1-2]
	\arrow["{p_{\lambda \mu}}"', from=1-1, to=2-1]
	\arrow["{d_{\mu}}", from=1-2, to=2-2]
	\arrow["{c_{\lambda}}"', from=2-1, to=2-2]
\end{tikzcd}\]
By \cref{prop:closure properties of light profinite sets}, $S''_{\lambda \mu}$ is a light profinite set. We define $e_{\lambda \mu} : S''_{\lambda \mu} \to S$ to be the continuous map $c_{\lambda} \of p_{\lambda \mu} = d_{\mu} \of q_{\lambda \mu}$. Then we obtain a cover $(e_{\lambda \mu} : S''_{\lambda \mu} \to S)_{(\lambda, \mu) \in \Lambda \times M}$ of $S$ in $\ub{Prof}$. We prove that $t_p |_{D(k), e_{\lambda \mu}} \in \alpha_{k, S''_{\lambda \mu}}(R(k)(S''_{\lambda \mu}))$ for every $(\lambda, \mu) \in \Lambda \times M$ and every $1 \leq p \leq n$.

Suppose $p \in \N$ and $1 \leq p \leq n-1$. For every $(\lambda, \mu) \in \Lambda \times M$, the following diagram is commutative.
\[\begin{tikzcd}
	{R(i)(S_{\lambda})} & {F(D(i))(S_{\lambda})} && {F(U)(S)} \\
	{R(k)(S''_{\lambda \mu})} & {F(D(k))(S''_{\lambda \mu})}
	\arrow["{\alpha_{i,S_{\lambda}}}", from=1-1, to=1-2]
	\arrow["{R(i \leq k ,\, p_{\lambda \mu})}"', from=1-1, to=2-1]
	\arrow["{F(D(k) \sub D(i) ,\, p_{\lambda \mu})}", from=1-2, to=2-2]
	\arrow["{F(D(i) \sub U ,\, c_{\lambda})}"', from=1-4, to=1-2]
	\arrow["{F(D(k) \sub U ,\, e_{\lambda \mu})}", curve={height=-18pt}, from=1-4, to=2-2]
	\arrow["{\alpha_{k,S''_{\lambda \mu}}}"', from=2-1, to=2-2]
\end{tikzcd}\]
Since $t_p |_{D(i),c_{\lambda}} \in \alpha_{i, S_{\lambda}}(R(i)(S_{\lambda}))$, it follows that $t_p |_{D(k), e_{\lambda \mu}} = \big( t_p |_{D(i),c_{\lambda}} \big) |_{D(k), p_{\lambda \mu}} \in \alpha_{k, S''_{\lambda \mu}}(R(k)(S''_{\lambda \mu}))$.

Suppose $p=n$. For every $(\lambda, \mu) \in \Lambda \times M$, the following diagram is commutative.
\[\begin{tikzcd}
	{R(j)(S'_{\mu})} & {F(D(j))(S'_{\mu})} && {F(U)(S)} \\
	{R(k)(S''_{\lambda \mu})} & {F(D(k))(S''_{\lambda \mu})}
	\arrow["{\alpha_{j,S'_{\mu}}}", from=1-1, to=1-2]
	\arrow["{R(j \leq k ,\, q_{\lambda \mu})}"', from=1-1, to=2-1]
	\arrow["{F(D(k) \sub D(j) ,\, q_{\lambda \mu})}", from=1-2, to=2-2]
	\arrow["{F(D(j) \sub U ,\, d_{\mu})}"', from=1-4, to=1-2]
	\arrow["{F(D(k) \sub U ,\, e_{\lambda \mu})}", curve={height=-18pt}, from=1-4, to=2-2]
	\arrow["{\alpha_{k,S''_{\lambda \mu}}}"', from=2-1, to=2-2]
\end{tikzcd}\]
Since $t_n |_{D(j),d_{\mu}} \in \alpha_{j, S'_{\mu}}(R(j)(S'_{\mu}))$, it follows that $t_n |_{D(k), e_{\lambda \mu}} = \big( t_n |_{D(j),d_{\mu}} \big) |_{D(k), q_{\lambda \mu}} \in \alpha_{k, S''_{\lambda \mu}}(R(k)(S''_{\lambda \mu}))$. This completes the proof.
\end{proof}

\begin{lem} \label{lem:stability of G under restriction}
Let $U,V$ be open subsets of $X$ such that $V \sub U$. Let $c : S' \to S$ be a morphism in $\ub{Prof}$. Suppose that $t \in G(U)(S)$. Then $F(V \sub U, c)(t) \in G(V)(S')$.
\end{lem}

\begin{proof}
For simplicity, we write $t' := F(V \sub U, c)(t)$. This is an element of $F(V)(S')$.

Suppose $x \in V$. Then $x \in U$. Since $t \in G(U)(S)$, there exist an $i \in I$ and a cover $(c_{\lambda} : S_{\lambda} \to S)_{\lambda \in \Lambda}$ of $S$ in $\ub{Prof}$ such that $x \in D(i) \sub U$ and $t|_{D(i),c_{\lambda}} \in \alpha_{i, S_{\lambda}}(R(i)(S_{\lambda}))$ for every $\lambda \in \Lambda$. On the other hand, (1) of \cref{as:assumption for constructing condensed sheaves} shows that there exists a $j \in I$ such that $x \in D(j) \sub V$.

First we note that $x \in D(i) \cap D(j)$. By (2) of \cref{as:assumption for constructing condensed sheaves}, there exists a $k \in I$ such that $i \leq k$, $\; j \leq k$ and $x \in D(k)$. Then $x \in D(k) \sub D(j) \sub V$.

Next, for each $\lambda \in \Lambda$, we define $S'_{\lambda} := S' \times_{S} S_{\lambda}$, with projections $p_{\lambda} : S'_{\lambda} \to S'$ and $q_{\lambda} : S'_{\lambda} \to S_{\lambda}$. This fiber product is taken in $\ub{Top}$. In other words, the following diagram is a pullback square in $\ub{Top}$.
\[\begin{tikzcd}
	{S'_{\lambda}} & {S_{\lambda}} \\
	{S'} & S
	\arrow["{q_{\lambda}}", from=1-1, to=1-2]
	\arrow["{p_{\lambda}}"', from=1-1, to=2-1]
	\arrow["{c_{\lambda}}", from=1-2, to=2-2]
	\arrow["c"', from=2-1, to=2-2]
\end{tikzcd}\]
By \cref{prop:closure properties of light profinite sets}, $S'_{\lambda}$ is a light profinite set. Then we obtain a cover $(p_{\lambda} : S'_{\lambda} \to S')_{\lambda \in \Lambda}$ of $S'$ in $\ub{Prof}$. 

We prove that $t'|_{D(k), p_{\lambda}} \in \alpha_{k, S'_{\lambda}}(R(k)(S'_{\lambda}))$ for every $\lambda \in \Lambda$. If $\lambda \in \Lambda$, then the following diagram is commutative.
\[\begin{tikzcd}
	{R(i)(S_{\lambda})} & {F(D(i))(S_{\lambda})} && {F(U)(S)} \\
	{R(k)(S'_{\lambda})} & {F(D(k))(S'_{\lambda})} && {F(V)(S')}
	\arrow["{\alpha_{i,S_{\lambda}}}", from=1-1, to=1-2]
	\arrow["{R(i \leq k ,\, q_{\lambda})}"', from=1-1, to=2-1]
	\arrow["{F(D(k) \sub D(i) ,\, q_{\lambda})}", from=1-2, to=2-2]
	\arrow["{F(D(i) \sub U ,\, c_{\lambda})}"', from=1-4, to=1-2]
	\arrow["{F(V \sub U ,\, c)}", from=1-4, to=2-4]
	\arrow["{\alpha_{k,S'_{\lambda}}}"', from=2-1, to=2-2]
	\arrow["{F(D(k) \sub V ,\, p_{\lambda})}", from=2-4, to=2-2]
\end{tikzcd}\]
Then we see that
\begin{align}
t'|_{D(k), p_{\lambda}}
& = F(D(k) \sub V, p_{\lambda}) \of F(V \sub U, c) \; (t) \\
& = F(D(k) \sub D(i), q_{\lambda}) \of F(D(i) \sub U, c_{\lambda}) \; (t) \\
& = F(D(k) \sub D(i), q_{\lambda}) (t|_{D(i),c_{\lambda}}) .
\end{align}
Moreover, since $t|_{D(i),c_{\lambda}} \in \alpha_{i, S_{\lambda}}(R(i)(S_{\lambda}))$, it follows that $F(D(k) \sub D(i), q_{\lambda}) (t|_{D(i),c_{\lambda}}) \in \alpha_{k, S'_{\lambda}}(R(k)(S'_{\lambda}))$. Therefore 
\begin{equation}
t'|_{D(k), p_{\lambda}}
= F(D(k) \sub D(i), q_{\lambda}) (t|_{D(i),c_{\lambda}}) 
\in \alpha_{k, S'_{\lambda}}(R(k)(S'_{\lambda})) .
\end{equation}
This completes the proof.
\end{proof}

\begin{lem}
For every open subset $U$ of $X$ and every light profinite set $S$, the subset $G(U)(S)$ is a subring of $F(U)(S)$.
\end{lem}

\begin{proof}
First we prove that $1 \in G(U)(S)$. Let $x \in U$. (1) of \cref{as:assumption for constructing condensed sheaves} shows that there exists an $i \in I$ such that $x \in D(i) \sub U$. $(\id_S : S \to S)$ is a cover of $S$ in $\ub{Prof}$. Furthermore, we have $1|_{D(i), \id_S} = 1 \in \alpha_{i,S}(R(i)(S))$, since $\alpha_{i,S} : R(i)(S) \to F(D(i))(S)$ is a homomorphism of unital rings. Thus $1 \in G(U)(S)$.

Next suppose that $t_1, t_2 \in G(U)(S)$. We prove that $t_1 - t_2, \, t_1 \cdot t_2 \in G(U)(S)$. Let $x \in U$. By \cref{lem:choosing i and (c lambda) uniformly}, there exist an $i \in I$ and a cover $(c_{\lambda} : S_{\lambda} \to S)_{\lambda \in \Lambda}$ of $S$ in $\ub{Prof}$ such that $x \in D(i) \sub U$ and $t_p |_{D(i),c_{\lambda}} \in \alpha_{i, S_{\lambda}}(R(i)(S_{\lambda}))$ for every $\lambda \in \Lambda$ and every $p \in \{ 1,2 \}$. Then, for every $\lambda \in \Lambda$, we have
\begin{align}
(t_1 - t_2)|_{D(i),c_{\lambda}}
& = (t_1 |_{D(i),c_{\lambda}}) - (t_2 |_{D(i),c_{\lambda}})
\in \alpha_{i, S_{\lambda}}(R(i)(S_{\lambda})) \: ; \\
(t_1 \cdot t_2)|_{D(i),c_{\lambda}}
& = (t_1 |_{D(i),c_{\lambda}}) \cdot (t_2 |_{D(i),c_{\lambda}})
\in \alpha_{i, S_{\lambda}}(R(i)(S_{\lambda})) ,
\end{align}
since $\alpha_{i,S} : R(i)(S) \to F(D(i))(S)$ is a homomorphism of rings. Thus $t_1 - t_2, \, t_1 \cdot t_2 \in G(U)(S)$.
\end{proof}

\begin{lem}
For every open subset $U$ of $X$, the family of subrings
\begin{equation}
G(U)(S) \sub F(U)(S) \quad (S \in |\ub{Prof}|)
\end{equation}
defines a condensed subring of $F(U)$. We write $G(U)$ for this.
\end{lem}

\begin{proof}
By \cref{lem:stability of G under restriction}, the family of subrings $G(U)(S) \sub F(U)(S) \; (S \in |\ub{Prof}|)$ defines a subpresheaf of commuatative unital rings of the sheaf $F(U)$ of commuatative unital rings on the site $\ub{Prof}$. Then it suffices to prove the following.
\begin{itemize}
\item
Let $S \in |\ub{Prof}|$. Let $(c_{\lambda} : S_{\lambda} \to S)_{\lambda \in \Lambda}$ be a cover of $S$ in $\ub{Prof}$. If $t \in F(U)(S)$ satisfies $F(U)(c_{\lambda})(t) \in G(U)(S_{\lambda})$ for every $\lambda \in \Lambda$, then $t \in G(U)(S)$.
\end{itemize}

Let us prove this assertion. Let $S \in |\ub{Prof}|$. Let $(c_{\lambda} : S_{\lambda} \to S)_{\lambda \in \Lambda}$ be a cover of $S$ in $\ub{Prof}$. Suppose that $t \in F(U)(S)$ satisfies $F(U)(c_{\lambda})(t) \in G(U)(S_{\lambda})$ for every $\lambda \in \Lambda$. We prove $t \in G(U)(S)$.

Let $x \in U$. For each $\lambda$, the relation $F(U)(c_{\lambda})(t) \in G(U)(S_{\lambda})$ implies that there exist an $i_{\lambda} \in I$ and a cover $(d_{\lambda \mu} : S_{\lambda \mu} \to S_{\lambda})_{\mu \in M_{\lambda}}$ of $S_{\lambda}$ in $\ub{Prof}$ such that $x \in D(i_{\lambda}) \sub U$ and $\big( F(U)(c_{\lambda})(t) \big) |_{D(i_{\lambda}), d_{\lambda \mu}} \in \alpha_{i_{\lambda}, S_{\lambda \mu}}(R(i_{\lambda})(S_{\lambda \mu}))$ for every $\mu \in M_{\lambda}$. On the other hand, (1) of \cref{as:assumption for constructing condensed sheaves} shows that there exists a $j \in I$ such that $x \in D(j) \sub U$.

First we note that $i_{\lambda} \in I(x)$ for every $\lambda \in \Lambda$ and that $j \in I(x)$. By \cref{lem:I(x) is filtered}, there exists a $k \in I(x)$ such that $j \leq k$ and $i_{\lambda} \leq k$ for every $\lambda \in \Lambda$. Then $x \in D(k) \sub D(j) \sub U$.

Next we note that $(c_{\lambda} \of d_{\lambda \mu} : S_{\lambda \mu} \to S)_{\mu \in M_{\lambda} ,\, \lambda \in \Lambda}$ is a cover of $S$ in $\ub{Prof}$.

We show that $t|_{D(k) ,\, c_{\lambda} \of d_{\lambda \mu}} \in \alpha_{k, S_{\lambda \mu}}(R(k)(S_{\lambda \mu}))$ for every $\mu \in M_{\lambda}$, $\lambda \in \Lambda$. This will imply that $t \in G(U)(S)$. If $\lambda \in \Lambda$ and $\mu \in M_{\lambda}$, then the following diagram is commutative.
\[\begin{tikzcd}
	{R(i_{\lambda})(S_{\lambda \mu})} & {F(D(i_{\lambda}))(S_{\lambda \mu})} && {F(U)(S_{\lambda})} \\
	{R(k)(S_{\lambda \mu})} & {F(D(k))(S_{\lambda \mu})} && {F(U)(S)}
	\arrow["{\alpha_{i_{\lambda},S_{\lambda \mu}}}", from=1-1, to=1-2]
	\arrow["{R(i_{\lambda} \leq k)_{S_{\lambda \mu}}}"', from=1-1, to=2-1]
	\arrow["{F(D(k) \sub D(i_{\lambda}))_{S_{\lambda \mu}}}", from=1-2, to=2-2]
	\arrow["{F(D(i_{\lambda}) \sub U ,\, d_{\lambda \mu})}"', from=1-4, to=1-2]
	\arrow["{\alpha_{k,S_{\lambda \mu}}}"', from=2-1, to=2-2]
	\arrow["{F(U)(c_{\lambda})}"', from=2-4, to=1-4]
	\arrow["{F(D(k) \sub U ,\, c_{\lambda} \of d_{\lambda \mu})}", from=2-4, to=2-2]
\end{tikzcd}\]
Then we see that
\begin{align}
t|_{D(k) ,\, c_{\lambda} \of d_{\lambda \mu}}
& = F(D(k) \sub U ,\, c_{\lambda} \of d_{\lambda \mu})(t) \\
& = F(D(k) \sub D(i_{\lambda}))_{S_{\lambda \mu}} \of F(D(i_{\lambda}) \sub U, d_{\lambda \mu}) \of F(U)(c_{\lambda}) \; (t) \\
& = F(D(k) \sub D(i_{\lambda}))_{S_{\lambda \mu}} \Big( \big( F(U)(c_{\lambda})(t) \big) |_{D(i_{\lambda}), d_{\lambda \mu}} \Big) .
\end{align}
Moreover, since $\big( F(U)(c_{\lambda})(t) \big) |_{D(i_{\lambda}), d_{\lambda \mu}} \in \alpha_{i_{\lambda}, S_{\lambda \mu}}(R(i_{\lambda})(S_{\lambda \mu}))$, it follows that
\begin{equation}
F(D(k) \sub D(i_{\lambda}))_{S_{\lambda \mu}} \Big( \big( F(U)(c_{\lambda})(t) \big) |_{D(i_{\lambda}), d_{\lambda \mu}} \Big) \in \alpha_{k, S_{\lambda \mu}}(R(k)(S_{\lambda \mu})) .
\end{equation}
Therefore
\begin{equation}
t|_{D(k) ,\, c_{\lambda} \of d_{\lambda \mu}}
= F(D(k) \sub D(i_{\lambda}))_{S_{\lambda \mu}} \Big( \big( F(U)(c_{\lambda})(t) \big) |_{D(i_{\lambda}), d_{\lambda \mu}} \Big)
\in \alpha_{k, S_{\lambda \mu}}(R(k)(S_{\lambda \mu})) .
\end{equation}
\end{proof}

\begin{lem}
The family of condensed subrings
\begin{equation}
G(U) \sub F(U) \quad (U \in |\cat{T}^{\op}|)
\end{equation}
defines a subsheaf of condensed rings of $F$. We write $G$ for this.
\end{lem}

\begin{proof}
By \cref{lem:stability of G under restriction}, the family of condensed subrings $G(U) \sub F(U) \; (U \in \cat{T}^{\op})$ defines a subpresheaf of condensed rings of $F$. Let us write $G$ for this subpresheaf of $F$. According to \cref{rem:sheaf property in terms of F S}, in order to prove that $G$ is a sub\ti{sheaf} of condensed rings of $F$, it suffices to verify the following assertion.
\begin{itemize}
\item
For each light profinite set $S$, the presheaf $G^S$ of commutative unital rings on $X$ is a sub\ti{sheaf} of the sheaf $F^S$ of commutative unital rings on $X$.
\end{itemize}
This is equivalent to the following assertion.
\begin{itemize}
\item
Let $S$ be a light profinite set. Let $U$ be an open subset of $X$. Let $(U_{\lambda})_{\lambda \in \Lambda}$ be an open covering of $U$ in $X$. If $t \in F(U)(S)$ satisfies $F(U_{\lambda} \sub U)(\id_S)(t) \in G(U_{\lambda})(S)$ for every $\lambda \in \Lambda$, then $t \in G(U)(S)$.
\end{itemize}
However, this immediately follows from the definition of $G(U)(S) \; (U \in \cat{T} ,\, S \in |\ub{Prof}|)$ given in \cref{nt:definition of G(U)(S)}.
\end{proof}

\begin{nt} \;
\begin{enumerate}
\item
We write $\iota : G \mon F$ for the inclusion of the subsheaf $G$ into $F$.

\item
By the definition of $G$, we have a unique natural transformation $\beta : R \To G \of D$ of functors $I \to \ub{CRing}$ such that the diagram
\[\begin{tikzcd}
	R & {F \of D} \\
	& {G \of D}
	\arrow["\alpha", from=1-1, to=1-2]
	\arrow["\beta"', from=1-1, to=2-2]
	\arrow["{\iota * \id_D}"', hook, from=2-2, to=1-2]
\end{tikzcd}\]
is commutative.
\end{enumerate}
\end{nt}

\subsubsection{Proof of (1) of \cref{prop:constructing condensed sheaves}}

\begin{lem} \label{lem:description of stalk in terms of D(i)}
Let $L$ be any sheaf of condensed rings on $X$. For each $x \in X$, the homomorphisms of condensed rings
\[\begin{tikzcd}
	{L(D(i))} & {L_x} & {(i \in I(x))}
	\arrow["\can", from=1-1, to=1-2]
\end{tikzcd}\]
induce an isomorphism of condensed rings
\[\begin{tikzcd}
	{\underset{i \in I(x)}{\colim} \; L(D(i))} & {L_x .}
	\arrow["\sim", from=1-1, to=1-2]
\end{tikzcd}\]
\end{lem}

\begin{proof}
(1) of \cref{as:assumption for constructing condensed sheaves} implies that for each $U \in \cat{T}(x)$, there exists an $i \in I(x)$ such that $D(i) \sub U$. (2) of \cref{as:assumption for constructing condensed sheaves} implies that for each $U \in \cat{T}(x)$ and each $i,j \in I(x)$ with $D(i) \sub U$ and $D(j) \sub U$, there exists a $k \in I(x)$ such that $i \leq k$ and $j \leq k$. Then the dual of Proposition 2.11.2 of \cite{Borceux:cat1} shows that the functor $D : I(x) \to \cat{T}(x)^{\op}$, $i \mapsto D(i)$ is cofinal. Therefore $\big( L_x , (L(D(i)) \xto{\can} L_x)_{i \in I(x)} \big)$ is the colimit of the functor $I(x) \xto{D} \cat{T}(x)^{\op} \xto{L} \ub{CRing}$.
\end{proof}

\begin{nt}
For each $x \in X$, we define a homomorphism $\sigma_x : F_x \to R_x$ of condensed rings as follows. For each $i \in I(x)$, we have a homomorphism of condensed rings
\[\begin{tikzcd}
	{F(D(i))} & {R_x .}
	\arrow["{\pi_{D(i),x}}", from=1-1, to=1-2]
\end{tikzcd}\]
If $i,j \in I(x)$ and $i \leq j$, the the following diagram is commutative.
\[\begin{tikzcd}
	{F(D(i))} \\
	{F(D(j))} & {R_x}
	\arrow["{F(D(j) \sub D(i))}"', from=1-1, to=2-1]
	\arrow["{\pi_{D(i),x}}", from=1-1, to=2-2]
	\arrow["{\pi_{D(j),x}}"', from=2-1, to=2-2]
\end{tikzcd}\]
By \cref{lem:description of stalk in terms of D(i)}, there exists a unique homomorphism $\sigma_x : F_x \to R_x$ of condensed rings such that the following diagram is commutative for all $i \in I(x)$.
\[\begin{tikzcd}
	{F(D(i))} & {R_x} \\
	{F_x}
	\arrow["{\pi_{D(i),x}}", from=1-1, to=1-2]
	\arrow["\can"', from=1-1, to=2-1]
	\arrow["{\sigma_x}"', from=2-1, to=1-2]
\end{tikzcd}\]
\end{nt}

\begin{prop} \label{prop:beta induces isomorphisms on stalks}
For each $x \in X$, the homomorphisms of condensed rings
\[\begin{tikzcd}
	{R(i)} & {G(D(i))} & {(i \in I(x))}
	\arrow["{\beta_i}", from=1-1, to=1-2]
\end{tikzcd}\]
induce an isomorphism of condensed rings
\[\begin{tikzcd}
	{\beta_x : R_x} & {G_x .}
	\arrow["\sim", from=1-1, to=1-2]
\end{tikzcd}\]
\end{prop}

\begin{proof}
For each $i \in I(x)$, we have a homomorphism of condensed rings
\[\begin{tikzcd}
	{R(i)} & {G(D(i))} & {G_x .}
	\arrow["{\beta_i}", from=1-1, to=1-2]
	\arrow["\can", from=1-2, to=1-3]
\end{tikzcd}\]
If $i,j \in I(x)$ and $i \leq j$, then the following diagram is commutative.
\[\begin{tikzcd}
	{R(i)} && {G(D(i))} \\
	{R(j)} && {G(D(j))} & {G_x}
	\arrow["{\beta_i}", from=1-1, to=1-3]
	\arrow["{R(i \leq j)}"', from=1-1, to=2-1]
	\arrow["{G(D(j) \sub D(i))}"', from=1-3, to=2-3]
	\arrow["\can", from=1-3, to=2-4]
	\arrow["{\beta_j}"', from=2-1, to=2-3]
	\arrow["\can"', from=2-3, to=2-4]
\end{tikzcd}\]
Therefore there exists a unique homomorphism $\beta_x : R_x \to G_x$ of condensed rings such that the following diagram is commutative for all $i \in I(x)$.
\[\begin{tikzcd}
	{R(i)} & {G(D(i))} \\
	{R_x} & {G_x}
	\arrow["{\beta_i}", from=1-1, to=1-2]
	\arrow["{\rho_{x,i}}"', from=1-1, to=2-1]
	\arrow["\can", from=1-2, to=2-2]
	\arrow["{\beta_x}"', from=2-1, to=2-2]
\end{tikzcd}\]

We prove that $\beta_x : R_x \to G_x$ is an isomorphism of condensed rings. In order to prove this assertion, it suffices to show that $\beta_x : R_x \to G_x$ is an isomorphism when it is considered as a morphism in $\ub{CSet}$. By \cref{prop:monics and epics in CSet}, it is enough to show that $\beta_x : R_x \to G_x$ is both a monomorphism and an epimorphism when it is considered as a morphism in $\ub{CSet}$.

First we show that $\beta_x : R_x \to G_x$ is a monomorphism when considered as a morphism in $\ub{CSet}$. For every $i \in I(x)$, the following diagram is commutative. 
\[\begin{tikzcd}
	&& {R(i)} \\
	{R(i)} & {G(D(i))} & {F(D(i))} & {R_x} \\
	{R_x} & {G_x} & {F_x}
	\arrow["{\alpha_i}"', from=1-3, to=2-3]
	\arrow["{\rho_{x,i}}", from=1-3, to=2-4]
	\arrow[equals, from=2-1, to=1-3]
	\arrow["{\beta_i}"', from=2-1, to=2-2]
	\arrow["{\rho_{x,i}}"', from=2-1, to=3-1]
	\arrow["{\iota_{D(i)}}"', from=2-2, to=2-3]
	\arrow["\can", from=2-2, to=3-2]
	\arrow["{\pi_{D(i),x}}"', from=2-3, to=2-4]
	\arrow["\can"', from=2-3, to=3-3]
	\arrow["{\beta_x}"', from=3-1, to=3-2]
	\arrow["{\iota_x}"', from=3-2, to=3-3]
	\arrow["{\sigma_x}"', curve={height=12pt}, from=3-3, to=2-4]
\end{tikzcd}\]
Then we have $\sigma_x \of \iota_x \of \beta_x \of \rho_{x,i} = \rho_{x,i}$ for all $i \in I(x)$. Since $R_x = \colim_{i \in I(x)} \; R(i)$, we conclude that $\sigma_x \of \iota_x \of \beta_x = \id_{R_x}$. It follows that $\beta_x$ is a monomorphism when considered as a morphism in $\ub{CSet}$.

Next we show that $\beta_x : R_x \to G_x$ is an epimorphism when considered as a morphism in $\ub{CSet}$. We use \cref{prop:monics and epics in CSet}. Let $S$ be a light profinite set. Let $s \in G_x(S)$. By \cref{prop:compatibility of stalk and evaluation at S}, we have 
\begin{equation}
G_x (S) = \underset{U \in \cat{T}(x)^{\op}}{\colim} \; G(U)(S) ,
\end{equation}
where the colimit is taken in $\ub{Ring}$. Since the category $\cat{T}(x)^{\op}$ is filtered, there exist a $U \in |\cat{T}(x)^{\op}|$ and a $t \in G(U)(S)$ such that $t_x = s$. Since $t \in G(U)(S)$, there exist an $i \in I$ and a cover $(c_{\lambda} : S_{\lambda} \to S)_{\lambda \in \Lambda}$ of $S$ in $\ub{Prof}$ such that $x \in D(i) \sub U$ and $t|_{D(i),c_{\lambda}} \in \alpha_{i, S_{\lambda}}(R(i)(S_{\lambda}))$ for every $\lambda \in \Lambda$. For each $\lambda \in \Lambda$, the following diagram is commutative.
\[\begin{tikzcd}
	{R(i)(S_{\lambda})} & {G(D(i))(S_{\lambda})} && {G(U)(S)} \\
	{R_x(S_{\lambda})} & {G_x(S_{\lambda})} && {G_x(S)}
	\arrow["{\beta_{i,S_{\lambda}}}", from=1-1, to=1-2]
	\arrow["{(\rho_{x,i})_{S_{\lambda}}}"', from=1-1, to=2-1]
	\arrow["\can", from=1-2, to=2-2]
	\arrow["{G(D(i) \sub U ,\,c_{\lambda})}"', from=1-4, to=1-2]
	\arrow["\can", from=1-4, to=2-4]
	\arrow["{(\beta_x)_{S_{\lambda}}}"', from=2-1, to=2-2]
	\arrow["{G_x(c_{\lambda})}", from=2-4, to=2-2]
\end{tikzcd}\]
Then we have
\begin{equation}
G_x(c_{\lambda})(s) = G_x(c_{\lambda})(t_x) = \big( t|_{D(i),c_{\lambda}} \big)_x .
\end{equation}
Moreover, since we have 
\begin{equation}
t|_{D(i),c_{\lambda}}
\in \alpha_{i, S_{\lambda}}(R(i)(S_{\lambda}))
= \beta_{i, S_{\lambda}}(R(i)(S_{\lambda})) ,
\end{equation}
it follows that
\begin{equation}
\big( t|_{D(i),c_{\lambda}} \big)_x \in (\beta_x)_{S_{\lambda}} \Big( R_x(S_{\lambda}) \Big) .
\end{equation}
Consequently,
\begin{equation}
G_x(c_{\lambda})(s)
= \big( t|_{D(i),c_{\lambda}} \big)_x
\in (\beta_x)_{S_{\lambda}} \Big( R_x(S_{\lambda}) \Big) .
\end{equation}
By \cref{prop:monics and epics in CSet}, we conclude that $\beta_x : R_x \to G_x$ is an epimorphism when considered as a morphism in $\ub{CSet}$.
\end{proof}

\subsubsection{Proof of (2) of \cref{prop:constructing condensed sheaves}: the uniqueness part}

\begin{nt} In the remainder of this subsection, we fix the following notation.
\begin{enumerate}
\item $H$ is a sheaf of condensed rings on $X$.

\item $\gamma : R \To H \of D$ is a natural transformation of functors $I \to \ub{CRing}$.
\end{enumerate}
\end{nt}

\begin{nt}
Let $x \in X$. We define a homomorphism $\gamma_x : R_x \to H_x$ of condensed rings as follows. For each $i \in I(x)$, we have a homomorphism of condensed rings
\[\begin{tikzcd}
	{R(i)} & {H(D(i))} & {H_x .}
	\arrow["{\gamma_i}", from=1-1, to=1-2]
	\arrow["\can", from=1-2, to=1-3]
\end{tikzcd}\]
If $i,j \in I(x)$ and $i \leq j$, then the following diagram is commutative.
\[\begin{tikzcd}
	{R(i)} && {H(D(i))} \\
	{R(j)} && {H(D(j))} & {H_x}
	\arrow["{\gamma_i}", from=1-1, to=1-3]
	\arrow["{R(i \leq j)}"', from=1-1, to=2-1]
	\arrow["{H(D(j) \sub D(i))}"', from=1-3, to=2-3]
	\arrow["\can", from=1-3, to=2-4]
	\arrow["{\gamma_j}"', from=2-1, to=2-3]
	\arrow["\can"', from=2-3, to=2-4]
\end{tikzcd}\]
Therefore there exists a unique homomorphism $\gamma_x : R_x \to H_x$ of condensed rings such that the following diagram is commutative for all $i \in I(x)$.
\[\begin{tikzcd}
	{R(i)} & {H(D(i))} \\
	{R_x} & {H_x}
	\arrow["{\gamma_i}", from=1-1, to=1-2]
	\arrow["{\rho_{x,i}}"', from=1-1, to=2-1]
	\arrow["\can", from=1-2, to=2-2]
	\arrow["{\gamma_x}"', from=2-1, to=2-2]
\end{tikzcd}\]
\end{nt}

\begin{prop}
Suppose that $\delta, \epsilon : G \to H$ are morphisms of sheaves of condensed rings on $X$ such that the following diagrams are commutative.
\[\begin{tikzcd}
	R & {H \of D} & R & {H \of D} \\
	{G \of D} && {G \of D}
	\arrow["\gamma", from=1-1, to=1-2]
	\arrow["\beta"', from=1-1, to=2-1]
	\arrow["\gamma", from=1-3, to=1-4]
	\arrow["\beta"', from=1-3, to=2-3]
	\arrow["{\delta * \id_{D}}"', from=2-1, to=1-2]
	\arrow["{\epsilon * \id_{D}}"', from=2-3, to=1-4]
\end{tikzcd}\]
Then we have $\delta = \epsilon$. 
\end{prop}

\begin{proof}
Let $x \in X$. For each $i \in I(x)$, the following diagrams are commutative.
\[\begin{tikzcd}
	{R_x} && {R(i)} && {R_x} \\
	& {G(D(i))} && {H(D(i))} \\
	& {G_x} && {H_x}
	\arrow["{\beta_x}"', curve={height=12pt}, from=1-1, to=3-2]
	\arrow["{\rho_{x,i}}"', from=1-3, to=1-1]
	\arrow["{\rho_{x,i}}", from=1-3, to=1-5]
	\arrow["{\beta_i}"', from=1-3, to=2-2]
	\arrow["{\gamma_i}", from=1-3, to=2-4]
	\arrow["{\gamma_x}", curve={height=-12pt}, from=1-5, to=3-4]
	\arrow["{\delta_{D(i)}}"', from=2-2, to=2-4]
	\arrow["\can", from=2-2, to=3-2]
	\arrow["\can"', from=2-4, to=3-4]
	\arrow["{\delta_x}"', from=3-2, to=3-4]
\end{tikzcd}\]

\[\begin{tikzcd}
	{R_x} && {R(i)} && {R_x} \\
	& {G(D(i))} && {H(D(i))} \\
	& {G_x} && {H_x}
	\arrow["{\beta_x}"', curve={height=12pt}, from=1-1, to=3-2]
	\arrow["{\rho_{x,i}}"', from=1-3, to=1-1]
	\arrow["{\rho_{x,i}}", from=1-3, to=1-5]
	\arrow["{\beta_i}"', from=1-3, to=2-2]
	\arrow["{\gamma_i}", from=1-3, to=2-4]
	\arrow["{\gamma_x}", curve={height=-12pt}, from=1-5, to=3-4]
	\arrow["{\epsilon_{D(i)}}"', from=2-2, to=2-4]
	\arrow["\can", from=2-2, to=3-2]
	\arrow["\can"', from=2-4, to=3-4]
	\arrow["{\epsilon_x}"', from=3-2, to=3-4]
\end{tikzcd}\]
Since $R_x = \colim_{i \in I(x)} \, R(i)$, we conclude that the following diagrams are commutative.
\[\begin{tikzcd}
	{R_x} & {H_x} & {R_x} & {H_x} \\
	{G_x} && {G_x}
	\arrow["{\gamma_x}", from=1-1, to=1-2]
	\arrow["{\beta_x}"', from=1-1, to=2-1]
	\arrow["{\gamma_x}", from=1-3, to=1-4]
	\arrow["{\beta_x}"', from=1-3, to=2-3]
	\arrow["{\delta_x}"', from=2-1, to=1-2]
	\arrow["{\epsilon_x}"', from=2-3, to=1-4]
\end{tikzcd}\]
By \cref{prop:beta induces isomorphisms on stalks}, the homomorphism $\beta_x : R_x \to G_x$ is an isomorphism of condensed rings. It follows that $\delta_x = \epsilon_x$. Thus we have proved that $\delta_x = \epsilon_x$ for every $x \in X$. Then \cref{prop:testing equality of sheaf morphisms by stalks} shows that $\delta = \epsilon$.
\end{proof}

\subsubsection{Proof of (2) of \cref{prop:constructing condensed sheaves}: the existence part}

\begin{nt} \;
\begin{enumerate}
\item
For every open subset $U$ of $X$, we define
\begin{equation}
\tilde{H}(U) := \prod_{x \in U} H_x ,
\end{equation}
where the product is taken in $\ub{CRing}$. For each $x \in U$, we write $\pi'_{U,x} : \tilde{H}(U) \to H_x$ for the $x$-th projection.

\item
If $U,V$ are open subsets of $X$ such that $V \sub U$, then we define $\tilde{H}(V \sub U) : \tilde{H}(U) \to \tilde{H}(V)$ to be the unique homomorphism of condensed rings such that the diagram
\[\begin{tikzcd}
	{\tilde{H}(U)} \\
	{\tilde{H}(V)} & {H_x}
	\arrow["{\tilde{H}(V \sub U)}"', from=1-1, to=2-1]
	\arrow["{\pi'_{U,x}}", from=1-1, to=2-2]
	\arrow["{\pi'_{V,x}}"', from=2-1, to=2-2]
\end{tikzcd}\]
is commutative for all $x \in V$. 

\item Then we obtain a sheaf $\tilde{H}$ of condensed rings on $X$.
\end{enumerate}
\end{nt}

\begin{nt} \;
\begin{enumerate}
\item
For every open subset $U$ of $X$, let $\tilde{\gamma}_U : F(U) \to \tilde{H}(U)$ be the unique homomorphism of condensed rings such that the diagram
\[\begin{tikzcd}
	{F(U)} & {\tilde{H}(U)} \\
	{R_x} & {H_x}
	\arrow["{\tilde{\gamma}_U}", from=1-1, to=1-2]
	\arrow["{\pi_{U,x}}"', from=1-1, to=2-1]
	\arrow["{\pi'_{U,x}}", from=1-2, to=2-2]
	\arrow["{\gamma_x}"', from=2-1, to=2-2]
\end{tikzcd}\]
is commutative for all $x \in U$.

\item
Then the family $\tilde{\gamma} := (\tilde{\gamma}_U)_{U \in |\cat{T}^{\op}|}$ is a morphism $F \to \tilde{H}$ of sheaves of condensed rings on $X$.
\end{enumerate}
\end{nt}

\begin{nt} \;
\begin{enumerate}
\item
For every open subset $U$ of $X$, let $\theta_{U} : H(U) \to \tilde{H}(U)$ be the unique homomorphism of condensed rings such that the diagram
\[\begin{tikzcd}
	{H(U)} & {\tilde{H}(U)} \\
	& {H_x}
	\arrow["{\theta_U}", from=1-1, to=1-2]
	\arrow["\can"', from=1-1, to=2-2]
	\arrow["{\pi'_{U,x}}", from=1-2, to=2-2]
\end{tikzcd}\]
is commutative for all $x \in U$.

\item
Then the family $\theta := (\theta_U)_{U \in |\cat{T}^{\op}|}$ is a morphism $H \to \tilde{H}$ of sheaves of condensed rings on $X$.
\end{enumerate}
\end{nt}

\begin{nt}
Let $x \in X$. We define a homomorphism $\tau_x : \tilde{H}_x \to H_x$ of condensed rings as follows. For each open neighbourhood $U$ of $x$ in $X$, we have a homomorphism $\pi'_{U,x} : \tilde{H}(U) \to H_x$ of condensed rings. If $U,V$ are open neighbourhoods of $x$ in $X$ such that $V \sub U$, then the following diagram is commutative.
\[\begin{tikzcd}
	{\tilde{H}(U)} \\
	{\tilde{H}(V)} & {H_x}
	\arrow["{\tilde{H}(V \sub U)}"', from=1-1, to=2-1]
	\arrow["{\pi'_{U,x}}", from=1-1, to=2-2]
	\arrow["{\pi'_{V,x}}"', from=2-1, to=2-2]
\end{tikzcd}\]
Therefore there exists a unique homomorphism $\tau_x : \tilde{H}_x \to H_x$ of condensed rings such that the following diagram is commutative for all open neighbourhoods $U$ of $x$ in $X$.
\[\begin{tikzcd}
	{\tilde{H}(U)} \\
	{\tilde{H}_x} & {H_x}
	\arrow["\can"', from=1-1, to=2-1]
	\arrow["{\pi'_{U,x}}", from=1-1, to=2-2]
	\arrow["{\tau_x}"', from=2-1, to=2-2]
\end{tikzcd}\]
\end{nt}

\begin{lem} \label{lem:tilde gamma of alpha equals theta of gamma}
The following diagram is commutative.
\[\begin{tikzcd}
	R & {H \of D} \\
	{F \of D} & {\tilde{H} \of D}
	\arrow["\gamma", from=1-1, to=1-2]
	\arrow["\alpha"', from=1-1, to=2-1]
	\arrow["{\theta * \id_D}", from=1-2, to=2-2]
	\arrow["{\tilde{\gamma} * \id_D}"', from=2-1, to=2-2]
\end{tikzcd}\]
\end{lem}

\begin{proof}
Let $i \in I$. For every $x \in D(i)$, the following diagram is commutative.
\[\begin{tikzcd}
	{H(D(i))} && {\tilde{H}(D(i))} \\
	{R(i)} & {R_x} & {H_x} \\
	{F(D(i))} && {\tilde{H}(D(i))}
	\arrow["{\theta_{D(i)}}", from=1-1, to=1-3]
	\arrow["\can", from=1-1, to=2-3]
	\arrow["{\pi'_{D(i),x}}", from=1-3, to=2-3]
	\arrow["{\gamma_i}", from=2-1, to=1-1]
	\arrow["{\rho_{x,i}}", from=2-1, to=2-2]
	\arrow["{\alpha_i}"', from=2-1, to=3-1]
	\arrow["{\gamma_x}"', from=2-2, to=2-3]
	\arrow["{\pi_{D(i),x}}"', from=3-1, to=2-2]
	\arrow["{\tilde{\gamma}_{D(i)}}"', from=3-1, to=3-3]
	\arrow["{\pi'_{D(i),x}}"', from=3-3, to=2-3]
\end{tikzcd}\]
Since $\tilde{H}(D(i)) = \prod_{x \in D(i)} H_x$, it follows that the following diagram is commutative.
\[\begin{tikzcd}
	{R(i)} & {H(D(i))} \\
	{F(D(i))} & {\tilde{H}(D(i))}
	\arrow["{\gamma_i}", from=1-1, to=1-2]
	\arrow["{\alpha_i}"', from=1-1, to=2-1]
	\arrow["{\theta_{D(i)}}", from=1-2, to=2-2]
	\arrow["{\tilde{\gamma}_{D(i)}}"', from=2-1, to=2-2]
\end{tikzcd}\]
\end{proof}

\begin{lem} \label{lem:theta x has a retraction}
For each $x \in X$, the following diagram is commutative.
\[\begin{tikzcd}
	{H_x} & {H_x} \\
	{\tilde{H}_x}
	\arrow["{\id_{H_x}}", from=1-1, to=1-2]
	\arrow["{\theta_x}"', from=1-1, to=2-1]
	\arrow["{\tau_x}"', from=2-1, to=1-2]
\end{tikzcd}\]
\end{lem}

\begin{proof}
For every open neighbourhood $U$ of $x$ in $X$, the following diagram is commutative.
\[\begin{tikzcd}
	{H_x} \\
	{H(U)} & {\tilde{H}(U)} & {H_x} \\
	{H_x} & {\tilde{H}_x}
	\arrow["{\id_{H_x}}", curve={height=-12pt}, from=1-1, to=2-3]
	\arrow["\can", from=2-1, to=1-1]
	\arrow["{\theta_U}"', from=2-1, to=2-2]
	\arrow["\can"', from=2-1, to=3-1]
	\arrow["{\pi'_{U,x}}"', from=2-2, to=2-3]
	\arrow["\can"', from=2-2, to=3-2]
	\arrow["{\theta_x}"', from=3-1, to=3-2]
	\arrow["{\tau_x}"', curve={height=12pt}, from=3-2, to=2-3]
\end{tikzcd}\]
Since $H_x = \colim_{U \in \cat{T}(x)^{\op}} \; H(U)$, we conclude that $\tau_x \of \theta_x = \id_{H_x}$.
\end{proof}

\begin{lem} \label{lem:injectivity of theta U S}
Let $U$ be an open subset of $X$. Let $S$ be a light profinite set. Then the map $\theta_{U,S} : H(U)(S) \to \tilde{H}(U)(S)$ is injective.
\end{lem}

\begin{proof}
Let $x \in X$. \cref{prop:compatibility of stalk and evaluation at S} shows that the map $(\theta^S)_x : (H^S)_x \to (\tilde{H}^S)_x$ coincides with the map $(\theta_x)_S : H_x(S) \to \tilde{H}_x(S)$. However, \cref{lem:theta x has a retraction} shows that
\begin{equation}
(\tau_x)_S \of (\theta_x)_S = (\id_{H_x})_S = \id_{H_x(S)}.
\end{equation}
Therefore $(\theta_x)_S : H_x(S) \to \tilde{H}_x(S)$ is injective. It follows that $(\theta^S)_x : (H^S)_x \to (\tilde{H}^S)_x$ is injective. This holds for every $x \in X$. On the other hand, by \cref{rem:sheaf property in terms of F S}, $H^S, \tilde{H}^S$ are sheaves of commutative unital rings on $X$. Then the ordinary sheaf theory shows that for every open subset $U$ of $X$, the map $(\theta^S)_U : H^S(U) \to \tilde{H}^S(U)$ is injective. However, this map coincides with the map $\theta_{U,S} : H(U)(S) \to \tilde{H}(U)(S)$.
\end{proof}

\begin{lem} \label{lem:image of G under tilde gamma is in H}
Let $U$ be an open subset of $X$. Let $S$ be a light profinite set. Then we have
\begin{equation}
\tilde{\gamma}_{U,S} \Big( G(U)(S) \Big) \sub \theta_{U,S} \Big( H(U)(S) \Big) .
\end{equation}
\end{lem}

\begin{proof}
Let $t \in G(U)(S)$. We show that $\tilde{\gamma}_{U,S}(t) \in \theta_{U,S} \big( H(U)(S) \big)$. For every $x \in U$, there exist an $i \in I$ and a cover $(c_{\lambda} : S_{\lambda} \to S)_{\lambda \in \Lambda}$ of $S$ in $\ub{Prof}$ such that $x \in D(i) \sub U$ and $t|_{D(i),c_{\lambda}} \in \alpha_{i, S_{\lambda}}(R(i)(S_{\lambda}))$ for every $\lambda \in \Lambda$. For each $\lambda \in \Lambda$, \cref{lem:tilde gamma of alpha equals theta of gamma} shows that the following diagram is commutative.
\[\begin{tikzcd}
	{R(i)(S_{\lambda})} & {H(D(i))(S_{\lambda})} \\
	{F(D(i))(S_{\lambda})} & {\tilde{H}(D(i))(S_{\lambda})} \\
	{F(U)(S)} & {\tilde{H}(U)(S)}
	\arrow["{\gamma_{i,S_{\lambda}}}", from=1-1, to=1-2]
	\arrow["{\alpha_{i,S_{\lambda}}}"', from=1-1, to=2-1]
	\arrow["{\theta_{D(i),S_{\lambda}}}", from=1-2, to=2-2]
	\arrow["{\tilde{\gamma}_{D(i),S_{\lambda}}}"', from=2-1, to=2-2]
	\arrow["{F(D(i) \sub U ,\, c_{\lambda})}", from=3-1, to=2-1]
	\arrow["{\tilde{\gamma}_{U,S}}"', from=3-1, to=3-2]
	\arrow["{\tilde{H}(D(i) \sub U ,\, c_{\lambda})}"', from=3-2, to=2-2]
\end{tikzcd}\]
It follows that
\begin{align}
\tilde{H}(D(i))(c_{\lambda}) \Big( \tilde{H}(D(i) \sub U)_S \big( \tilde{\gamma}_{U,S}(t) \big) \Big)
& = \tilde{H}(D(i) \sub U, c_{\lambda}) \of \tilde{\gamma}_{U,S} \, (t) \\
& = \tilde{\gamma}_{D(i), S_{\lambda}} \of F(D(i) \sub U, c_{\lambda}) \, (t) \\
& = \tilde{\gamma}_{D(i), S_{\lambda}} \big( t|_{D(i),c_{\lambda}} \big) .
\end{align}
Since $t|_{D(i),c_{\lambda}} \in \alpha_{i, S_{\lambda}}(R(i)(S_{\lambda}))$, we have
\begin{equation}
\tilde{\gamma}_{D(i), S_{\lambda}} \big( t|_{D(i),c_{\lambda}} \big) 
\in \theta_{D(i), S_{\lambda}} \big( H(D(i))(S_{\lambda}) \big) .
\end{equation}
Therefore
\begin{equation}
\tilde{H}(D(i))(c_{\lambda}) \Big( \tilde{H}(D(i) \sub U)_S \big( \tilde{\gamma}_{U,S}(t) \big) \Big)
= \tilde{\gamma}_{D(i), S_{\lambda}} \big( t|_{D(i),c_{\lambda}} \big)
\in \theta_{D(i), S_{\lambda}} \big( H(D(i))(S_{\lambda}) \big) .
\end{equation}
However, \cref{lem:injectivity of theta U S} implies that the family of subrings
\begin{equation}
\theta_{D(i), S'} \big( H(D(i))(S') \big) \sub \tilde{H}(D(i))(S') \quad (S' \in |\ub{Prof}|)
\end{equation}
defines a condensed subring of $\tilde{H}(D(i))$. Since $(c_{\lambda} : S_{\lambda} \to S)_{\lambda \in \Lambda}$ is a cover of $S$ in $\ub{Prof}$, we conclude that
\begin{equation}
\tilde{H}(D(i) \sub U)_S \big( \tilde{\gamma}_{U,S}(t) \big)
\in \theta_{D(i), S} \big( H(D(i))(S) \big) .
\end{equation}

Thus we have proved the following assertion.
\begin{itemize}
\item[$(\dag)$]
For every $x \in U$, there exists an $i \in I$ such that $x \in D(i) \sub U$ and 
\begin{equation}
\tilde{H}(D(i) \sub U)_S \big( \tilde{\gamma}_{U,S}(t) \big)
\in \theta_{D(i), S} \big( H(D(i))(S) \big) .
\end{equation}
\end{itemize}
If we use \cref{nt:passing from condesed sheaves to ordinary sheaves}, this assertion can be expressed as follows.
\begin{itemize}
\item[$(\ddag)$]
For every $x \in U$, there exists an $i \in I$ such that $x \in D(i) \sub U$ and 
\begin{equation}
\tilde{H}^S (D(i) \sub U) \big( \tilde{\gamma}_{U,S}(t) \big)
\in (\theta^S)_{D(i)} \big( H^S (D(i)) \big) .
\end{equation}
\end{itemize}
On the other hand, by \cref{rem:sheaf property in terms of F S}, $H^S, \tilde{H}^S$ are sheaves of commutative unital rings on $X$. Moreover, \cref{lem:injectivity of theta U S} shows that the morphism $\theta^S : H^S \to \tilde{H}^S$ of sheaves of commutative unital rings on $X$ has the property that the map $(\theta^S)_V : H^S(V) \to \tilde{H}^S(V)$ is injective for every open subset $V$ of $X$. Hence the family of subrings
\begin{equation}
(\theta^S)_V \big( H^S(V) \big)
\sub \tilde{H}^S(V) \quad (V \in |\cat{T}^{\op}|)
\end{equation}
defines a subsheaf of commutative unital rings of $\tilde{H}^S$ on $X$. Consequently, $(\ddag)$ implies that
\begin{equation}
\tilde{\gamma}_{U,S}(t) \in (\theta^S)_U \big( H^S(U) \big) .
\end{equation}
Thus
\begin{equation}
\tilde{\gamma}_{U,S}(t)
\in (\theta^S)_U \big( H^S(U) \big)
= \theta_{U,S} \Big( H(U)(S) \Big) .
\end{equation}
\end{proof}

\begin{nt} \;
\begin{enumerate}
\item
Let $U$ be an open subset of $X$. Let $S$ be a light profinite set. By \cref{lem:injectivity of theta U S} and \cref{lem:image of G under tilde gamma is in H}, there exists a unique ring homomorphism $\delta_{U,S} : G(U)(S) \to H(U)(S)$ such that the following diagram is commutative.
\[\begin{tikzcd}
	{F(U)(S)} & {\tilde{H}(U)(S)} \\
	{G(U)(S)} & {H(U)(S)}
	\arrow["{\tilde{\gamma}_{U,S}}", from=1-1, to=1-2]
	\arrow["{\iota_{U,S}}", hook', from=2-1, to=1-1]
	\arrow["{\delta_{U,S}}"', from=2-1, to=2-2]
	\arrow["{\theta_{U,S}}"', hook', from=2-2, to=1-2]
\end{tikzcd}\]

\item
For every open subset $U$ of $X$, the family $\delta_U := (\delta_{U,S})_{S \in |\ub{Prof}^{\op}|}$ is a homomorphism $G(U) \to H(U)$ of condensed rings.

\item
The family $\delta := (\delta_U)_{U \in |\cat{T}^{\op}|}$ is a morphism $G \to H$ of sheaves of condensed rings on $X$.
\end{enumerate}
\end{nt}

\begin{prop}
The following diagram is commutative.
\[\begin{tikzcd}
	R & {H \of D} \\
	{G \of D}
	\arrow["\gamma", from=1-1, to=1-2]
	\arrow["\beta"', from=1-1, to=2-1]
	\arrow["{\delta * \id_D}"', from=2-1, to=1-2]
\end{tikzcd}\]
\end{prop}

\begin{proof}
Let $i \in I$. For every $x \in D(x)$, the following diagram is commutative.
\[\begin{tikzcd}
	& {R(i)} & {H(D(i))} & {\tilde{H}(D(i))} \\
	{G(D(i))} & {F(D(i))} & {R_x} \\
	{H(D(i))} & {\tilde{H}(D(i))} && {H_x}
	\arrow["{\gamma_i}", from=1-2, to=1-3]
	\arrow["{\beta_i}"', from=1-2, to=2-1]
	\arrow["{\alpha_i}"', from=1-2, to=2-2]
	\arrow["{\rho_{x,i}}", from=1-2, to=2-3]
	\arrow["{\theta_{D(i)}}", from=1-3, to=1-4]
	\arrow["\can", from=1-3, to=3-4]
	\arrow["{\pi'_{D(i),x}}", from=1-4, to=3-4]
	\arrow["{\iota_{D(i)}}"', from=2-1, to=2-2]
	\arrow["{\delta_{D(i)}}"', from=2-1, to=3-1]
	\arrow["{\pi_{D(i),x}}"', from=2-2, to=2-3]
	\arrow["{\tilde{\gamma}_{D(i)}}", from=2-2, to=3-2]
	\arrow["{\gamma_x}"', from=2-3, to=3-4]
	\arrow["{\theta_{D(i)}}"', from=3-1, to=3-2]
	\arrow["{\pi'_{D(i),x}}"', from=3-2, to=3-4]
\end{tikzcd}\]
Since $\tilde{H}(D(i)) = \prod_{x \in D(i)} H_x$, we conclude that the following diagram is commutative.
\[\begin{tikzcd}
	{R(i)} && {H(D(i))} \\
	{G(D(i))} \\
	{H(D(i))} && {\tilde{H}(D(i))}
	\arrow["{\gamma_i}", from=1-1, to=1-3]
	\arrow["{\beta_i}"', from=1-1, to=2-1]
	\arrow["{\theta_{D(i)}}", from=1-3, to=3-3]
	\arrow["{\delta_{D(i)}}"', from=2-1, to=3-1]
	\arrow["{\theta_{D(i)}}"', from=3-1, to=3-3]
\end{tikzcd}\]
By \cref{lem:injectivity of theta U S}, the homomorphism $\theta_{D(i)} : H(D(i)) \to \tilde{H}(D(i))$ is a monomorphism in $\ub{CRing}$. Therefore we conclude that the following diagram is commutative.
\[\begin{tikzcd}
	{R(i)} & {H(D(i))} \\
	{G(D(i))}
	\arrow["{\gamma_i}", from=1-1, to=1-2]
	\arrow["{\beta_i}"', from=1-1, to=2-1]
	\arrow["{\delta_{D(i)}}"', from=2-1, to=1-2]
\end{tikzcd}\]
\end{proof}

\subsection{Sheafification}

\begin{prop} \label{prop:construction of sheafification of condensed sheaves}
Let $X$ be a topological space. Let $F$ be a presheaf of condensed rings on $X$. Then there exist a sheaf $G$ of condensed rings on $X$ and a morphism $\eta : F \to G$ of presheaves of condensed rings on $X$ with the following property.
\begin{enumerate}
\item For every $x \in X$, the homomorphism $\eta_x : F_x \to G_x$ is an isomorphism of condensed rings.

\item For every sheaf $H$ of condensed rings on $X$ and every morphism $\gamma : F \to H$ of presheaves of condensed rings on $X$, there exists a unique morphism $\delta : G \to H$ of sheaves of condensed rings on $X$ such that the diagram
\[\begin{tikzcd}
	F & H \\
	G
	\arrow["\gamma", from=1-1, to=1-2]
	\arrow["\eta"', from=1-1, to=2-1]
	\arrow["\delta"', from=2-1, to=1-2]
\end{tikzcd}\]
is commutative.
\end{enumerate}
\end{prop}

\begin{proof}
We use \cref{prop:constructing condensed sheaves}. Let $\cat{T}$ be the set of all open subsets of $X$ ordered by inclusion, which we consider as a small category. Then $\cat{T}^{\op}$ is a directed set. Let $D : \cat{T}^{\op} \to \cat{T}^{\op}$ be the identity functor. Then \cref{as:assumption for constructing condensed sheaves} is trivially satisfied. Moreover, we have a functor $F : \cat{T}^{\op} \to \ub{CRing}$. In this situation, \cref{prop:constructing condensed sheaves} precisely states that there exist a sheaf $G$ of condensed rings on $X$ and a morphism $\eta : F \to G$ of presheaves of condensed rings on $X$ having the properties (1) and (2).
\end{proof}

\begin{cor} \label{cor:existence of sheafification functor of condensed presheaves}
Let $X$ be a topological space. 
\begin{enumerate}
\item
The category $\ub{CSh}(X)$ is a reflective full subcategory of $\ub{CPSh}(X)$.

\item
Let $F \in |\ub{CPSh}(X)|$. Let $(G', F \xto{\eta'} G')$ be any reflection of $F$ along the inclusion functor $\ub{CSh}(X) \mon \ub{CPSh}(X)$. Then the homomorphism $\eta'_x : F_x \to G'_x$ is an isomorphism of condensed rings for every $x \in X$.
\end{enumerate}
\end{cor}

\begin{proof}~
\begin{enumerate}
\item
\cref{prop:construction of sheafification of condensed sheaves} implies that every $F \in |\ub{CPSh}(X)|$ has a reflection along the inclusion functor $\ub{CSh}(X) \mon \ub{CPSh}(X)$. This shows (1).

\item
By \cref{prop:construction of sheafification of condensed sheaves}, there exists a reflection $(G, F \xto{\eta} G)$ of $F$ along the inclusion functor $\ub{CSh}(X) \mon \ub{CPSh}(X)$ with the property that the homomorphism $\eta_x : F_x \to G_x$ is an isomorphism of condensed rings for every $x \in X$. By the uniqueness of reflection, there exists a unique isomorphism $\alpha : G \to G'$ in $\ub{CSh}(X)$ such that the following diagram is commutative.
\[\begin{tikzcd}
	F & G \\
	& {G'}
	\arrow["\eta", from=1-1, to=1-2]
	\arrow["{\eta'}"', from=1-1, to=2-2]
	\arrow["\alpha", from=1-2, to=2-2]
\end{tikzcd}\]
Then the following diagram is commutative for all $x \in X$.
\[\begin{tikzcd}
	{F_x} & G \\
	& {G'}
	\arrow["{\eta_x}", from=1-1, to=1-2]
	\arrow["{\eta'_x}"', from=1-1, to=2-2]
	\arrow["{\alpha_x}", from=1-2, to=2-2]
\end{tikzcd}\]
Since $\eta_x : F_x \to G_x$ and $\alpha_x : G_x \to G'_x$ are isomorphisms of condensed rings, it follows that $\eta'_x : F_x \to G'_x$ is an isomorphism of condensed rings.
\end{enumerate}
\end{proof}

\begin{df}
Let $X$ be a topological space.
\begin{enumerate}
\item
The left adjoint of the inclusion functor $\ub{CSh}(X) \mon \ub{CPSh}(X)$ is called the \ti{sheafification functor}.

\item
Let $F$ be a presheaf of condensed rings on $X$. Let $(G, F \xto{\eta} G)$ be the reflection of $F$ along the inclusion functor $\ub{CSh}(X) \mon \ub{CPSh}(X)$. Then the sheaf $G$ of condensed rings on $X$ is called the \ti{sheafification} of $F$. The morphism $\eta : F \to G$ of presheaves of condensed rings on $X$ is also called the sheafification of $F$.
\end{enumerate}
\end{df}

\begin{cor} \label{cor:CSh(X) is complete and cocomplete}
Let $X$ be a topological space. Then the category $\ub{CSh}(X)$ is complete and cocomplete.
\end{cor}

\begin{proof}
Since the category $\ub{CRing}$ is complete and cocomplete, the category $\ub{CPSh}(X)$ is also complete and cocomplete by Theorem 2.15.2 of \cite{Borceux:cat1}. By \cref{cor:existence of sheafification functor of condensed presheaves}, the category $\ub{CSh}(X)$ is a reflective full subcategory of $\ub{CPSh}(X)$. Then the result follows from Proposition 3.5.3 and Proposition 3.5.4 of \cite{Borceux:cat1}.
\end{proof}

\subsection{Direct images}

\begin{df}
Let $X,Y$ be topological spaces. Let $f : X \to Y$ be a continuous map.
\begin{enumerate}
\item
Let $F$ be a sheaf of condensed ring on $X$. For each open subset $V$ of $Y$, define
\begin{equation}
f_* F (V) := F(f^{-1}(V)) .
\end{equation}
We obtain a sheaf $f_* F$ of condensed rings on $Y$, which we call the \ti{direct image} of $F$ under $f$.

\item
We obtain a functor $f_* : \ub{CSh}(X) \to \ub{CSh}(Y)$, $F \mapsto f_* F$.
\end{enumerate}
\end{df}

\begin{nt}
Let $X,Y$ be topological spaces. Let $f : X \to Y$ be a continuous map. Let $F$ be a sheaf of condensed ring on $X$. Let $G$ be a sheaf of condensed ring on $Y$. Let $\alpha : G \to f_* F$ be a morphism of sheaves of condensed rings on $Y$. For each $x \in X$, the homomorphisms of condensed rings
\[\begin{tikzcd}[ampersand replacement=\&]
	{G(V)} \& {F(f^{-1}(V))} \& { \left(
	\begin{gathered}
	V \text{ is an open} \\ \text{neighbourhood of } f(x)
	\end{gathered}
	\right) }
	\arrow["{\alpha_V}", from=1-1, to=1-2]
\end{tikzcd}\]
induce a homomorphism of condensed rings
\[\begin{tikzcd}
	{G_{f(x)}} & {F_x .}
	\arrow[from=1-1, to=1-2]
\end{tikzcd}\]
By abuse of notation, we write $\alpha_x$ for this homomorphism.
\end{nt}

\begin{prop} \label{prop:testing equality of sheaf morphisms to direct images by stalks}
Let $X,Y$ be topological spaces. Let $f : X \to Y$ be a continuous map. Let $F$ be a sheaf of condensed ring on $X$. Let $G$ be a sheaf of condensed ring on $Y$. Let $\alpha, \beta : G \to f_* F$ be two morphisms of sheaves of condensed rings on $Y$. Then $\alpha = \beta$ if and only if $\alpha_x = \beta_x : G_{f(x)} \to F_x$ for every $x \in X$.
\end{prop}

\begin{proof}
Suppose that $\alpha_x = \beta_x : G_{f(x)} \to F_x$ for every $x \in X$. We prove that $\alpha = \beta$. Let $S$ be any light profinite set. Then we have $(f_*F)^S = f_*(F^S)$ by definition. Thus we have two morphisms $\alpha^S, \beta^S : G^S \to f_*(F^S)$ of presheaves of commutative unital rings on $Y$. If $x \in X$, then (1) of \cref{prop:compatibility of stalk and evaluation at S} shows that the maps $(\alpha^S)_x , (\beta^S)_x : (G^S)_{f(x)} \to (F^S)_x$ coincide with the maps $(\alpha_x)_S , (\beta_x)_S : G_{f(x)}(S) \to F_x(S)$, respectively. Since $\alpha_x = \beta_x : G_{f(x)} \to F_x$ by hypothesis, we conclude that $(\alpha^S)_x = (\beta^S)_x : (G^S)_{f(x)} \to (F^S)_x$. This holds for every $x \in X$. On the other hand, by \cref{rem:sheaf property in terms of F S}, $F^S, G^S$ are sheaves of commutative unital rings on $X,Y$ respectively. Then the ordinary sheaf theory shows that $\alpha^S = \beta^S : G^S \to f_*(F^S)$. Consequently,
\begin{equation}
\alpha_{V,S} = (\alpha^S)_V = (\beta^S)_V = \beta_{V,S}
\end{equation}
for every open subset $V$ of $Y$ and every light profinite set $S$. It follows that $\alpha = \beta$.
\end{proof}

\subsection{Inverse images}

\begin{prop} \label{prop:construction of inverse images of condensed sheaves}
Let $X,Y$ be a topological space. Let $f : X \to Y$ be a continuous map. Let $G$ be a sheaf of condensed rings on $Y$. Then there exists a sheaf $F$ of condensed rings on $X$ and a morphism $\eta : G \to f_* F$ of sheaves of condensed rings on $Y$ with the following property.
\begin{enumerate}
\item For every $x \in X$, the homomorphism $\eta_x : G_{f(x)} \to F_x$ is an isomorphism of condensed rings.

\item For every sheaf $E$ of condensed rings on $X$ and every morphism $\gamma : G \to f_* E$ of sheaves of condensed rings on $Y$, there exists a unique morphism $\delta : F \to E$ of sheaves of condensed rings on $X$ such that the diagram
\[\begin{tikzcd}
	G & {f_* E} \\
	{f_* F}
	\arrow["\gamma", from=1-1, to=1-2]
	\arrow["\eta"', from=1-1, to=2-1]
	\arrow["{f_* \delta}"', from=2-1, to=1-2]
\end{tikzcd}\]
is commutative.
\end{enumerate}
\end{prop}

\begin{proof}
We use \cref{prop:constructing condensed sheaves}. Let $\cat{T}$ (resp.$\; \cat{S}$) be the set of all open subsets of $X$ (resp.$\; Y$) ordered by inclusion, which we consider as a small category. Let $(I, \leq)$ be the directed set defined as follows.
\begin{itemize}
\item
We define $I$ to be the set of all pairs $(U,V)$ consisting of an open subset $U$ of $X$ and an open subset $V$ of $Y$ such that $f(U) \sub V$.

\item
We define the relation $\leq$ on $I$ by declaring that for $(U,V) ,\, (U',V') \in I$, the relation $(U,V) \leq (U',V')$ holds if and only if both $U \bus U'$ and $V \bus V'$ hold.
\end{itemize}
We also consider the directed set $I$ as a small category. Let $D : I \to \cat{T}^{\op}$ be the functor $(U,V) \mapsto U$. Then one checks easily that \cref{as:assumption for constructing condensed sheaves} is satisfied. Let $B : I \to \cat{S}^{\op}$ be the functor $(U,V) \mapsto V$, and define the functor $R : I \to \ub{CRing}$ to be the composition $I \xto{B} \cat{S}^{\op} \xto{G} \ub{CRing}$. Then \cref{prop:constructing condensed sheaves} shows that there exists a sheaf $F$ of condensed rings on $X$ and a natural transformation $\beta : R \To F \of D$ with the following property.
\begin{enumerate}
\item
For each $x \in X$, let us write $I(x)$ for the full subcategory of $I$ consisting of all $i \in I$ such that $x \in D(i)$. Then the homomorphisms of condensed rings
\[\begin{tikzcd}
	{R(i)} & {F(D(i))} & {(i \in I(x))}
	\arrow["{\beta_i}", from=1-1, to=1-2]
\end{tikzcd}\]
induce an isomorphism of condensed rings
\[\begin{tikzcd}
	{\underset{i \in I(x)}{\colim} \, R(i)} & {F_x .}
	\arrow["{\beta_x}", from=1-1, to=1-2]
\end{tikzcd}\]

\item
For any sheaf $E$ of condensed rings on $X$ and any natural transformation $\gamma : R \To E \of D$ of functors $I \to \ub{CRing}$, there exists a unique morphism $\delta : F \to E$ of sheaves of condensed rings on $X$ such that the following diagram is commutative.
\[\begin{tikzcd}
	R & {E \of D} \\
	{F \of D}
	\arrow["\gamma", from=1-1, to=1-2]
	\arrow["\beta"', from=1-1, to=2-1]
	\arrow["{\delta * \id_D}"', from=2-1, to=1-2]
\end{tikzcd}\]
\end{enumerate}
Let $C : \cat{S}^{\op} \to I$ be the functor $V \mapsto (f^{-1}(V),V)$. Then the functor $B \of C : \cat{S}^{\op} \to \cat{S}^{\op}$ is equal to the identity, and the functor $D \of C : \cat{S}^{\op} \to \cat{T}^{\op}$ is equal to the functor $V \mapsto f^{-1}(V)$. Consequently, $\eta := \beta * \id_C : R \of C \to F \of D \of C$ is a morphism $G \to f_* F$ of sheaves of condensed rings on $Y$.

Let $x \in X$. We claim that the homomorphism $\eta_x : G_{f(x)} \to F_x$ is an isomorphism of condensed rings. Let $\cat{S}(f(x))$ be the full subcategory of $\cat{S}$ consisting of all open neighbourhoods of $f(x)$ in $Y$. Then we have a functor
\[\begin{tikzcd}
	{B : I(x)} & {\cat{S}(f(x))^{\op},} & {(U,V)} & {V .}
	\arrow[from=1-1, to=1-2]
	\arrow[maps to, from=1-3, to=1-4]
\end{tikzcd}\]
Using the dual of Proposition 2.11.2 of \cite{Borceux:cat1}, one checks easily that this functor is cofinal. Therefore
\begin{equation}
\left( G_{f(x)} ,\, \left( R(i) = G(B(i)) \xto{\can} G_{f(x)} \right)_{i \in |I(x)|} \right)
\end{equation}
is the colimit of the functor $R : I(x) \xto{B} \cat{S}(f(x))^{\op} \xto{G} \ub{CRing}$. Consequently, the isomorphism $\beta_x : \colim_{i \in I(x)} \, R(i) \to F_x$ is an isomorphism $G_{f(x)} \to F_x$ of condensed rings with the property that the following diagram is commutative for all $i \in I(x)$.
\[\begin{tikzcd}
	{G(B(i))} & {F(D(i))} \\
	{G_{f(x)}} & {F_x}
	\arrow["{\beta_i}", from=1-1, to=1-2]
	\arrow["\can"', from=1-1, to=2-1]
	\arrow["\can", from=1-2, to=2-2]
	\arrow["{\beta_x}"', from=2-1, to=2-2]
\end{tikzcd}\]
In particular, for every $V \in |\cat{S}(f(x))^{\op}|$, the following diagram is commutative.
\[\begin{tikzcd}
	{G(V)} & {F(f^{-1}(V))} \\
	{G(B(C(V)))} & {F(D(C(V)))} \\
	{G_{f(x)}} & {F_x}
	\arrow["{\eta_V}", from=1-1, to=1-2]
	\arrow[equals, from=1-1, to=2-1]
	\arrow[equals, from=1-2, to=2-2]
	\arrow["{\beta_{C(V)}}", from=2-1, to=2-2]
	\arrow["\can"', from=2-1, to=3-1]
	\arrow["\can", from=2-2, to=3-2]
	\arrow["{\beta_x}"', from=3-1, to=3-2]
\end{tikzcd}\]
Then the definition of $\eta_x$ shows that $\eta_x$ is equal to $\beta_x$. Therefore $\eta_x : G_{f(x)} \to F_x$ is an isomorphism of condensed rings.

Next suppose that $E$ is a sheaf of condensed rings on $X$ and that $\gamma : G \to f_* E$ is a morphism of sheaves of condensed rings on $Y$. We prove that there exists a unique morphism $\delta : F \to E$ of sheaves of condensed rings on $X$ such that the diagram
\[\begin{tikzcd}
	G & {f_* E} \\
	{f_* F}
	\arrow["\gamma", from=1-1, to=1-2]
	\arrow["\eta"', from=1-1, to=2-1]
	\arrow["{f_* \delta}"', from=2-1, to=1-2]
\end{tikzcd}\]
is commutative. For each $i=(U,V) \in I$, let $\tilde{\gamma}_i : R(i) \to E(D(i))$ be the homomorphism of condensed rings
\[\begin{tikzcd}
	{R(i) = G(V)} & {E(f^{-1}(V))} && {E(U) = E(D(i)) .}
	\arrow["{\gamma_V}", from=1-1, to=1-2]
	\arrow["{E(U \sub f^{-1}(V))}", from=1-2, to=1-4]
\end{tikzcd}\]
Then the family $\tilde{\gamma} := (\tilde{\gamma}_i)_{i \in I}$ is a natural transformation $R \To E \of D$ of functors $I \to \ub{CRing}$. Therefore there exists a unique morphism $\delta : F \to E$ of sheaves of condensed rings on $X$ such that the diagram
\[\begin{tikzcd}
	R & {E \of D} \\
	{F \of D}
	\arrow["{\tilde{\gamma}}", from=1-1, to=1-2]
	\arrow["\beta"', from=1-1, to=2-1]
	\arrow["{\delta * \id_D}"', from=2-1, to=1-2]
\end{tikzcd}\]
is commutative. Then the following diagram is commutative for every open subset $V$ of $Y$.
\[\begin{tikzcd}
	{G(V)} \\
	& {R(C(V))} & {E(D(C(V)))} & {E(f^{-1}(V))} \\
	& {F(D(C(V)))} \\
	& {F(f^{-1}(V))}
	\arrow[equals, from=1-1, to=2-2]
	\arrow["{\gamma_V}", curve={height=-12pt}, from=1-1, to=2-4]
	\arrow["{\eta_V}"', curve={height=18pt}, from=1-1, to=4-2]
	\arrow["{\tilde{\gamma}_{C(V)}}", from=2-2, to=2-3]
	\arrow["{\beta_{C(V)}}"', from=2-2, to=3-2]
	\arrow[equals, from=2-4, to=2-3]
	\arrow["{\delta_{D(C(V))}}"', from=3-2, to=2-3]
	\arrow["{\delta_{f^{-1}(V)}}"', from=4-2, to=2-4]
	\arrow[equals, from=4-2, to=3-2]
\end{tikzcd}\]
It follows that the diagram
\[\begin{tikzcd}
	G & {f_* E} \\
	{f_* F}
	\arrow["\gamma", from=1-1, to=1-2]
	\arrow["\eta"', from=1-1, to=2-1]
	\arrow["{f_* \delta}"', from=2-1, to=1-2]
\end{tikzcd}\]
is commutative. On the other hand, suppose that $\epsilon : F \to E$ is another morphism of sheaves of condensed rings on $X$ such that the diagram
\[\begin{tikzcd}
	G & {f_* E} \\
	{f_* F}
	\arrow["\gamma", from=1-1, to=1-2]
	\arrow["\eta"', from=1-1, to=2-1]
	\arrow["{f_* \epsilon}"', from=2-1, to=1-2]
\end{tikzcd}\]
is commutative. If $i=(U,V) \in I$, then the following diagram is commutative.
\[\begin{tikzcd}
	{G(V)} & {E(f^{-1}(V))} && {E(U)} \\
	{F(f^{-1}(V))} && {F(U)}
	\arrow["{\gamma_V}", from=1-1, to=1-2]
	\arrow["{\eta_V}"', from=1-1, to=2-1]
	\arrow["{E(U \sub f^{-1}(V))}", from=1-2, to=1-4]
	\arrow["{\epsilon_{f^{-1}(V)}}"', from=2-1, to=1-2]
	\arrow["{F(U \sub f^{-1}(V))}"', from=2-1, to=2-3]
	\arrow["{\epsilon_U}"', from=2-3, to=1-4]
\end{tikzcd}\]
However, the following diagram is commutative.
\[\begin{tikzcd}
	{G(V)} &&& {G(V)} \\
	& {R(C(V))} && {R(i)} \\
	& {F(D(C(V)))} && {F(D(i))} \\
	{F(f^{-1}(V))} &&& {F(U)}
	\arrow["{\id_{G(V)}}", from=1-1, to=1-4]
	\arrow[equals, from=1-1, to=2-2]
	\arrow["{\eta_V}"', from=1-1, to=4-1]
	\arrow[equals, from=1-4, to=2-4]
	\arrow["{R(C(V) \leq i)}", from=2-2, to=2-4]
	\arrow["{\beta_{C(V)}}"', from=2-2, to=3-2]
	\arrow["{\beta_i}", from=2-4, to=3-4]
	\arrow["{F(D(i) \sub D(C(V)))}"', from=3-2, to=3-4]
	\arrow[equals, from=3-4, to=4-4]
	\arrow[equals, from=4-1, to=3-2]
	\arrow["{F(U \sub f^{-1}(V))}"', from=4-1, to=4-4]
\end{tikzcd}\]
Therefore $G(V) \xto{\eta_V} F(f^{-1})(V)) \xto{F(U \sub f^{-1}(V))} F(U)$ is equal to $R(i) \xto{\beta_i} F(D(i))$. Moreover, since $D(i) = U$, the homomorphism $F(U) \xto{\epsilon_U} E(U)$ is equal to $F(D(i)) \xto{\epsilon_{D(i)}} E(D(i))$. The homomorphism $G(V) \xto{\gamma_V} E(f^{-1}(V)) \xto{E(U \sub f^{-1}(V))} E(U)$ is equal to $R(i) \xto{\tilde{\gamma}_i} E(D(i))$ by definition. Consequently, the following diagram is commutative.
\[\begin{tikzcd}
	{R(i)} & {E(D(i))} \\
	{F(D(i))}
	\arrow["{\tilde{\gamma}_i}", from=1-1, to=1-2]
	\arrow["{\beta_i}"', from=1-1, to=2-1]
	\arrow["{\epsilon_{D(i)}}"', from=2-1, to=1-2]
\end{tikzcd}\]
This diagram is commutative for every $i \in I$. Therefore the diagram
\[\begin{tikzcd}
	R & {E \of D} \\
	{F \of D}
	\arrow["{\tilde{\gamma}}", from=1-1, to=1-2]
	\arrow["\beta"', from=1-1, to=2-1]
	\arrow["{\epsilon * \id_D}"', from=2-1, to=1-2]
\end{tikzcd}\]
is commutative. Then the uniqueness of $\delta$ shows that $\epsilon = \delta$.
\end{proof}

\begin{cor} \label{cor:existence of inverse image functor}
Let $X,Y$ be topological spaces. Let $f : X \to Y$ be a continuous map. 
\begin{enumerate}
\item
The functor $f_* : \ub{CSh}(X) \to \ub{CSh}(Y)$ has a left adjoint.

\item
Let $G \in |\ub{CSh}(Y)|$. Let $(F', G \xto{\eta'} f_* F')$ be any reflection of $G$ along the functor $f_* : \ub{CSh}(X) \to \ub{CSh}(Y)$. Then the homomorphism $\eta'_x : G_{f(x)} \to F'_x$ is an isomorphism of condensed rings for every $x \in X$.
\end{enumerate}
\end{cor}

\begin{proof}~
\begin{enumerate}
\item
\cref{prop:construction of inverse images of condensed sheaves} implies that every $G \in |\ub{CSh}(Y)|$ has a reflection along the functor $f_* : \ub{CSh}(X) \to \ub{CSh}(Y)$. This shows (1).

\item
By \cref{prop:construction of inverse images of condensed sheaves}, there exists a reflection $(F ,\, G \xto{\eta} f_* F)$ of $G$ along the functor $f_* : \ub{CSh}(X) \to \ub{CSh}(Y)$ with the property that the homomorphism $\eta_x : G_{f(x)} \to F_x$ is an isomorphism of condensed rings for every $x \in X$. By the uniqueness of reflection, there exists a unique isomorphism $\alpha : F \to F'$ in $\ub{CSh}(X)$ such that the following diagram is commutative.
\[\begin{tikzcd}
	G & {f_* F} \\
	& {f_* F'}
	\arrow["\eta", from=1-1, to=1-2]
	\arrow["{\eta'}"', from=1-1, to=2-2]
	\arrow["{f_* \alpha}", from=1-2, to=2-2]
\end{tikzcd}\]
It follows that the following diagram is commutative for all $x \in X$.
\[\begin{tikzcd}
	{G_{f(x)}} & {F_x} \\
	& {F'_x}
	\arrow["{\eta_x}", from=1-1, to=1-2]
	\arrow["{\eta'_x}"', from=1-1, to=2-2]
	\arrow["{\alpha_x}", from=1-2, to=2-2]
\end{tikzcd}\]
Since $\eta_x : G_{f(x)} \to F_x$ and $\alpha_x : F_x \to F'_x$ are isomorphisms of condensed rings, it follows that $\eta'_x : G_{f(x)} \to F'_x$ is an isomorphism of condensed rings.
\end{enumerate}
\end{proof}

\begin{df}
Let $X,Y$ be topological spaces. Let $f : X \to Y$ be a continuous map.
\begin{enumerate}
\item
The left adjoint of the functor $f_* : \ub{CSh}(X) \to \ub{CSh}(Y)$ is denoted by $f^{-1} : \ub{CSh}(Y) \to \ub{CSh}(X)$.

\item
If $G$ is a sheaf of condensed rings on $Y$, then the sheaf $f^{-1} G$ of condensed rings on $X$ is called the \ti{inverse image} of $G$ under $f$.
\end{enumerate}
\end{df}

\subsection{Inverse images and stalks}

\begin{prop} \label{prop:taking global section over * is an equivalence}
Let $\{*\}$ be the topological space with a unique point $*$. Then the functor
\[\begin{tikzcd}[row sep=tiny]
	{\ub{CSh}(\{*\})} & {\ub{CRing}} \\
	F & {F(\{*\})} & {(\text{on objects})} \\
	\alpha & {\alpha_{\{*\}}} & {(\text{on morphisms})}
	\arrow[from=1-1, to=1-2]
	\arrow[maps to, from=2-1, to=2-2]
	\arrow[maps to, from=3-1, to=3-2]
\end{tikzcd}\]
is an equivalence of categories.
\end{prop}

\begin{proof}
This immediately follows from the fact that a presheaf $F$ of condensed rings on $\{*\}$ is a sheaf of condensed rings on $\{*\}$ if and only if $F(\ku) = 0$.
\end{proof}

\begin{prop} \label{prop:description of stalks via inverse images}
Let $X$ be a topological space. Let $x \in X$. Let $\{*\}$ be the topological space with a unique point $*$. Let $j : \{*\} \to X$ be the continuous map $* \mapsto x$. 
Let $\Gamma : \ub{CSh}(\{*\}) \to \ub{CRing}$ be the functor
\[\begin{tikzcd}[row sep=tiny]
	{\ub{CSh}(\{*\})} & {\ub{CRing}} \\
	F & {F(\{*\})} & {(\text{on objects})} \\
	\alpha & {\alpha_{\{*\}}} & {(\text{on morphisms})}
	\arrow["\Gamma", from=1-1, to=1-2]
	\arrow[maps to, from=2-1, to=2-2]
	\arrow[maps to, from=3-1, to=3-2]
\end{tikzcd}\]
Then the functor
\[\begin{tikzcd}[row sep=tiny]
	{\ub{CSh}(X)} & {\ub{CRing}} \\
	F & {F_x} & {(\text{on objects})} \\
	\alpha & {\alpha_x} & {(\text{on morphisms})}
	\arrow["\Sigma", from=1-1, to=1-2]
	\arrow[maps to, from=2-1, to=2-2]
	\arrow[maps to, from=3-1, to=3-2]
\end{tikzcd}\]
is isomorphic to the composition
\[\begin{tikzcd}
	{\ub{CSh}(X)} & {\ub{CSh}(\{*\})} & {\ub{CRing} .}
	\arrow["{j^{-1}}", from=1-1, to=1-2]
	\arrow["\Gamma", from=1-2, to=1-3]
\end{tikzcd}\]
\end{prop}

\begin{proof}
By defintion, the functor $j^{-1} : \ub{CSh}(X) \to \ub{CSh}(\{*\})$ is left adjoint to the functor $j_* : \ub{CSh}(\{*\}) \to \ub{CSh}(X)$. Let $\eta : \id_{\ub{CSh}(X)} \To j_* \of j^{-1}$ be the unit of of this adjunction. (2) of \cref{cor:existence of inverse image functor} shows that for every $F \in |\ub{CSh}(X)|$, the homomorphism $\left( \eta_F \right)_* : F_x \to \left( j^{-1}F \right)_{*}$ is an isomorphism of condensed rings. On the other hand, since $*$ is the unique point of $\{*\}$, we have 
\begin{equation}
\left( j^{-1}F \right)_{*} =  \left( j^{-1}F \right) (\{*\}) .
\end{equation}
Therefore $\left( \eta_F \right)_*$ is an isomorphism $F_x \to \left( j^{-1}F \right) (\{*\})$ of condensed rings. Furthermore, if $\alpha : F \to G$ is a morphism in $\ub{CSh}(X)$, the following diagram is commutative in $\ub{CSh}(X)$.
\[\begin{tikzcd}
	F && G \\
	{j_* j^{-1} F} && {j_* j^{-1} G}
	\arrow["\alpha", from=1-1, to=1-3]
	\arrow["{\eta_F}"', from=1-1, to=2-1]
	\arrow["{\eta_G}", from=1-3, to=2-3]
	\arrow["{j_* j^{-1} \alpha}"', from=2-1, to=2-3]
\end{tikzcd}\]
It follows that the following diagram is commutative in $\ub{CRing}$. 
\[\begin{tikzcd}
	{F_x} && {G_x} \\
	{\left( j^{-1}F \right)_{*}} && {\left( j^{-1}G \right)_{*}}
	\arrow["{\alpha_x}", from=1-1, to=1-3]
	\arrow["{\left( \eta_F \right)_*}"', from=1-1, to=2-1]
	\arrow["{\left( \eta_G \right)_*}", from=1-3, to=2-3]
	\arrow["{\left( j^{-1}\alpha \right)_{*}}"', from=2-1, to=2-3]
\end{tikzcd}\]
On the other hand, since $*$ is the unique point of $\{*\}$, we have
\begin{equation}
\left( j^{-1}\alpha \right)_{*} =  \left( j^{-1}\alpha \right)_{\{*\}} .
\end{equation}
Consequently, the following diagram is commutative.
\[\begin{tikzcd}
	{F_x} && {G_x} \\
	{\left( j^{-1}F \right) (\{*\})} && {\left( j^{-1}G \right) (\{*\})}
	\arrow["{\alpha_x}", from=1-1, to=1-3]
	\arrow["{\left( \eta_F \right)_*}"', from=1-1, to=2-1]
	\arrow["{\left( \eta_G \right)_*}", from=1-3, to=2-3]
	\arrow["{\left( j^{-1}\alpha \right)_{\{*\}}}"', from=2-1, to=2-3]
\end{tikzcd}\]
It follows that the family
\[\begin{tikzcd}
	{F_x} & {\left( j^{-1}F \right) (\{*\})} & {(F \in |\ub{CSh}(X)|)}
	\arrow["{\left( \eta_F \right)_*}", from=1-1, to=1-2]
\end{tikzcd}\]
is an isomorphism $\Sigma \xTo{\sim} \Gamma \of j^{-1}$ of functors $\ub{CSh}(X) \to \ub{CRing}$.
\end{proof}

\begin{cor} \label{cor:right adjoint of the stalk functor}
Let $X$ be a topological space. Let $x \in X$. Then the functor
\[\begin{tikzcd}[row sep=tiny]
	{\ub{CSh}(X)} & {\ub{CRing}} \\
	F & {F_x} & {(\text{on objects})} \\
	\alpha & {\alpha_x} & {(\text{on morphisms})}
	\arrow["\Sigma", from=1-1, to=1-2]
	\arrow[maps to, from=2-1, to=2-2]
	\arrow[maps to, from=3-1, to=3-2]
\end{tikzcd}\]
has a right adjoint. In particular, it preserves colimits.
\end{cor}

\begin{proof}
Let $\{*\}$ be the topological space with a unique point $*$. Let $j : \{*\} \to X$ be the continuous map $* \mapsto x$. By \cref{prop:taking global section over * is an equivalence}, the functor
\[\begin{tikzcd}[row sep=tiny]
	{\ub{CSh}(\{*\})} & {\ub{CRing}} \\
	F & {F(\{*\})} & {(\text{on objects})} \\
	\alpha & {\alpha_{\{*\}}} & {(\text{on morphisms})}
	\arrow["\Gamma", from=1-1, to=1-2]
	\arrow[maps to, from=2-1, to=2-2]
	\arrow[maps to, from=3-1, to=3-2]
\end{tikzcd}\]
is an equivalence. Moreover, \cref{prop:description of stalks via inverse images} shows that the functor $\Sigma : \ub{CSh}(X) \to \ub{CRing}$ is isomorphic to the composition
\[\begin{tikzcd}
	{\ub{CSh}(X)} & {\ub{CSh}(\{*\})} & {\ub{CRing} .}
	\arrow["{j^{-1}}", from=1-1, to=1-2]
	\arrow["\Gamma", from=1-2, to=1-3]
\end{tikzcd}\]
Since the functor $j^{-1} : \ub{CSh}(X) \to \ub{CSh}(\{*\})$ has a right adjoint $j_* : \ub{CSh}(\{*\}) \to \ub{CSh}(X)$, it follows that $\Sigma : \ub{CSh}(X) \to \ub{CRing}$ also has a right adjoint.
\end{proof}

\subsection{Stalk criteria}

\begin{prop} \label{prop:testing invertibility by stalk}
Let $X$ be a topological space. Let $F$ be a sheaf of condensed rings on $X$. Let $U$ be an open subset of $X$. Let $S$ be a light profinite set. Let $f \in F(U)(S)$. Then the following conditions are equivalent.
\begin{enumerate}
\item
$f \in F(U)(S)$ is an invertible element of $F(U)(S)$.

\item
For every $x \in U$, the element $f_x \in F_x(S)$ is an invertible element of $F_x(S)$.
\end{enumerate}
\end{prop}

\begin{proof}
Consider the presheaf $F^S$ of commutative unital rings on $X$ (\cref{nt:passing from condesed sheaves to ordinary sheaves}). This is in fact a sheaf of commutative unital rings on $X$ (\cref{rem:sheaf property in terms of F S}). Moreover, for each $x \in U$, we have
\begin{equation}
(F^S)_x = F_x (S)
\end{equation}
by \cref{prop:compatibility of stalk and evaluation at S}. Therefore the assertion follows from the ordinary sheaf theory. 
\end{proof}

\begin{prop} \label{prop:testing coalescence by stalks}
Let $X$ be a topological space. Let $F$ be a sheaf of condensed rings on $X$. Let $U$ be an open subset of $X$. Let $f,g \in F(U)(*)$. Then the following conditions are equivalent.
\begin{enumerate}
\item
For every open subset $V$ of $X$ with $V \sub U$, the condensed ring $F(V)$ becomes an $f/g$-coalescent condensed $F(U)$-algebra via the restriction $F(V \sub U) : F(U) \to F(V)$.

\item
For every $x \in U$, the condensed ring $F_x$ becomes an $f/g$-coalescent condensed $F(U)$-algebra via the canonical homomorphism $F(U) \to F_x$.
\end{enumerate}
\end{prop}

\begin{proof}
Let $\cat{T}$ be the set of all open subsets of $X$ contained in $U$, ordered by inclusion. We consider $\cat{T}$ as a small category. For every open subset $V$ of $X$ with $V \sub U$, we consider $F(V)$ as a condensed $F(U)$-algebra via the restriction $F(V \sub U) : F(U) \to F(V)$. For every $x \in U$, we consider $F_x$ as a condensed $F(U)$-algebra via the canonical homomorphism $F(U) \to F_x$. Then the following hold.
\begin{itemize}
\item
We have a functor
\[\begin{tikzcd}[row sep=tiny]
	{\cat{T}^{\op}} & {\ub{CMod}_{F(U)}} \\
	V & {F(V)} & {(\text{on objects})} \\
	{(W \sub V)} & {F(W \sub V)} & {(\text{on morphisms}).}
	\arrow["G", from=1-1, to=1-2]
	\arrow[maps to, from=2-1, to=2-2]
	\arrow[maps to, from=3-1, to=3-2]
\end{tikzcd}\]

\item
The forgetful functor $\ub{CAlg}_{F(U)} \to \ub{CRing}$ reflects limits and filtered colimits. The forgetful functor $\ub{CAlg}_{F(U)} \to \ub{CMod}_{F(U)}$ preserves limits and filtered colimits. 

\item
Using the notation given in \cref{df:definition of sequence space}, consider the functor
\[\begin{tikzcd}
	{D : \ub{CMod}_{F(U)}} & {\ub{CMod}_{F(U)} ,} & M & {\uhom_{F(U)}(P_{F(U)},M)}
	\arrow[from=1-1, to=1-2]
	\arrow[maps to, from=1-3, to=1-4]
\end{tikzcd}\]
and the natural transformation $\delta : D \To D$ defined by
\begin{equation}
\delta := \left( \uhom_{F(U)}(\Delta_{f/g}, \id_M) \right)_{M \in |\ub{CMod}_{F(U)}|} .
\end{equation}
By \cref{prop:P R is internally projective}, the functor $D : \ub{CMod}_{F(U)} \to \ub{CMod}_{F(U)}$ preserves limits and colimits. 

\item
For each light profinite set $S$, \cref{prop:evaluation of CMod R preserves limits filtered colimits and coproducts} shows that the functor
\[\begin{tikzcd}
	{\mathrm{ev}_S : \ub{CMod}_{F(U)}} & {\ub{Mod}_{F(U)(S)} ,} & M & {M(S)}
	\arrow[from=1-1, to=1-2]
	\arrow[maps to, from=1-3, to=1-4]
\end{tikzcd}\]
preserves limits and filtered colimits. 
\end{itemize}
From these facts, one concludes that the following hold.
\begin{itemize}
\item
For each light profinite set $S$, the functor $H_S := \mathrm{ev}_{S} \of D \of G : \cat{T}^{\op} \to \ub{Mod}_{F(U)(S)}$ is a sheaf of $F(U)(S)$-modules on the topological space $U$. We have a morphism $\epsilon_S := \id_{\mathrm{ev}_{S}} * \delta * \id_{G} : H_S \to H_S$ of sheaves of $F(U)(S)$-modules on $U$.

\item
For each light profinite set $S$ and each $x \in U$, the stalk $(H_S)_x$ is equal to 
\begin{equation}
\mathrm{ev}_S ( D(F_x) ) = \uhom_{F(U)}(P_{F(U)}, F_x)(S) ,
\end{equation}
and the map $(\epsilon_S)_x : (H_S)_x \to (H_S)_x$ is equal to the map 
\begin{equation}
\uhom_{F(U)}(\Delta_{f/g}, \id_{F_x})_S :
\uhom_{F(U)}(P_{F(U)}, F_x)(S) \to \uhom_{F(U)}(P_{F(U)}, F_x)(S) .
\end{equation}
\end{itemize}
By the characterization (1) of \cref{prop:characterization of coalescence}, we conclude that the condition (2) of this proposition holds if and only if the map $(\epsilon_S)_x : (H_S)_x \to (H_S)_x$ is an isomorphism of $F(U)(S)$-modules for every $x \in U$ and every light profinite set $S$. Since $H_S$ is a sheaf of $F(U)(S)$-modules on $U$ for every light profinite $S$, the ordinary sheaf theory shows that the map $(\epsilon_S)_x : (H_S)_x \to (H_S)_x$ is an isomorphism of $F(U)(S)$-modules for every $x \in U$ and every light profinite set $S$ if and only if $\epsilon_S : H_S \to H_S$ is an isomorphism of sheaves of $F(U)(S)$-modules on $U$ for every light profinite set $S$. On the other hand, the definition shows that for each light profinite set $S$ and each open subset $V$ of $X$ with $V \sub U$ we have 
\begin{equation}
H_S(V) = \uhom_{F(U)}(P_{F(U)}, F(V))(S) ,
\end{equation}
and the map $(\epsilon_S)_V : H_S(V) \to H_S(V)$ is equal to the map 
\begin{equation}
\uhom_{F(U)}(\Delta_{f/g}, \id_{F(V)})_S :
\uhom_{F(U)}(P_{F(U)}, F(V))(S) \to \uhom_{F(U)}(P_{F(U)}, F(V))(S) .
\end{equation}
By the characterization (1) of \cref{prop:characterization of coalescence}, we conclude that $\epsilon_S : H_S \to H_S$ is an isomorphism of sheaves of $F(U)(S)$-modules on $U$ for every light profinite set $S$ if and only if the condition (1) of this proposition holds. This completes the proof.
\end{proof}

\section{Condensed ringed spaces} \label{sec:Condensed ringed spaces}

In this section, we first give the definition of the categories $\cat{D}$ and $\cat{C}$, together with methods for constructing limits and colimits in these categories. The category $\cat{D}$ is an analogue of the category of ringed spaces in the setting of condensed mathematics. On the other hand, the category $\cat{C}$ corresponds to the category $\ub{PRS}$ of primed ringed spaces in Definition 3 of \cite{Gillam:paper}. After that, we introduce the categories $\cat{C}_1$, $\cat{C}_l$, $\cat{C}_f$ and $\cat{C}_c$, all of which are proved to be coreflective full subcategories of $\cat{C}$.

\subsection{The category $\cat{D}$}

\subsubsection{Definition}

\begin{df}
We define the category $\cat{D}$ as follows.
\begin{enumerate}
\item
An object of $\cat{D}$ is a pair $(X, \cat{O}_X)$ consisting of a topological space $X$ and a sheaf $\cat{O}_X$ of condensed rings on $X$.

\item
A morphism $(X, \cat{O}_X) \to (Y, \cat{O}_Y)$ in $\cat{D}$ is a pair $(f,f^{\#})$ consisting of a continuous map $f : X \to Y$ and a morphism $f^{\#} : \cat{O}_Y \to f_* \cat{O}_X$ of sheaves of condensed rings on $Y$.

\item
The composition of morphisms $(f,f^{\#}) : (X, \cat{O}_X) \to (Y, \cat{O}_Y)$ and $(g,g^{\#}) : (Y, \cat{O}_Y) \to (Z, \cat{O}_Z)$ in $\cat{D}$ is defined to be $(g \of f ,\, (g_* f^{\#}) \of g^{\#} )$.
\end{enumerate}
\end{df}

\begin{rem}
We have a forgetful functor
\[\begin{tikzcd}[row sep=tiny]
	{\cat{D}} & {\ub{Top}} \\
	{(X, \cat{O}_X)} & X & {(\text{on objects})} \\
	{(f, f^{\#})} & f & {(\text{on morphisms})}
	\arrow[from=1-1, to=1-2]
	\arrow[maps to, from=2-1, to=2-2]
	\arrow[maps to, from=3-1, to=3-2]
\end{tikzcd}\]
In this subsection, we write $F$ for this forgetful functor.
\end{rem}

\subsubsection{Completeness}

\begin{prop} \label{prop:construction of limits in the category D}
Let $\cat{I}$ be a small category. Let $D : \cat{I} \to \cat{D}$ be a functor. For each $i \in |\cat{I}|$, let us write $(X_i , \cat{O}_i) := D(i)$. For each morphism $\alpha$ in $\cat{I}$, let us write $(f_{\alpha}, f^{\#}_{\alpha}) := D(\alpha)$. Let
\begin{equation}
\left( X, (X \xto{\pi_i} X_i)_{i \in |\cat{I}|} \right)
\end{equation}
be the limit of the functor $F \of D : \cat{I} \to \ub{Top}$.
\begin{enumerate}
\item
For each $i \in |\cat{I}|$, let $\tilde{\cat{O}}_i$ be the inverse image of $\cat{O}_i$ under $\pi_i : X \to X_i$. Let us write $\rho_i : \cat{O}_i \to \pi_{i,*} \tilde{\cat{O}}_i$ for the canonical morphism of sheaves of condensed rings on $X_i$. Then, for each morphism $\alpha : i \to j$ in $\cat{I}$, there exists a unique morphism $\tilde{f}^{\#}_{\alpha} : \tilde{\cat{O}}_j \to \tilde{\cat{O}}_i$ of sheaves of condensed rings on $X$ such that the following diagram is commutative in $\ub{CSh}(X_j)$.
\[\begin{tikzcd}
	{\pi_{j,*} \tilde{\cat{O}}_i} & {f_{\alpha,*} \cat{O}_i} \\
	{\pi_{j,*} \tilde{\cat{O}}_j} & {\cat{O}_j}
	\arrow["{f_{\alpha,*} \rho_i}"', from=1-2, to=1-1]
	\arrow["{\pi_{j,*} \tilde{f}^{\#}_{\alpha}}", from=2-1, to=1-1]
	\arrow["{f^{\#}_{\alpha}}"', from=2-2, to=1-2]
	\arrow["{\rho_j}", from=2-2, to=2-1]
\end{tikzcd}\]

\item
\[\begin{tikzcd}[row sep=tiny]
	{\cat{I}^{\op}} & {\ub{CSh}(X)} \\
	{i^{\op}} & {\tilde{\cat{O}}_i} & {(\text{on objects})} \\
	{\alpha^{\op}} & {\tilde{f}^{\#}_{\alpha}} & {(\text{on morphisms})}
	\arrow["{\tilde{D}}", from=1-1, to=1-2]
	\arrow[maps to, from=2-1, to=2-2]
	\arrow[maps to, from=3-1, to=3-2]
\end{tikzcd}\]
defines a functor.

\item
By \cref{cor:CSh(X) is complete and cocomplete}, the category $\ub{CSh}(X)$ is cocomplete. Let
\begin{equation}
\left(\cat{O} , ( \tilde{\cat{O}}_i \xto{\sigma_i} \cat{O} )_{i \in |\cat{I}^{\op}|} \right)
\end{equation}
be the colimit of the functor $\tilde{D} : \cat{I}^{\op} \to \ub{CSh}(X)$. For each $i \in |\cat{I}|$, let us write $\pi^{\#}_i : \cat{O}_i \to \pi_{i,*} \cat{O}$ for the composition $\cat{O}_i \xto{\rho_i} \pi_{i,*} \tilde{\cat{O}}_i \xto{\pi_{i,*} \sigma_i} \pi_{i,*} \cat{O}$ in $\ub{CSh}(X_i)$. Then
\begin{equation}
\left( (X,\cat{O}),
\left( (X,\cat{O}) \xto{(\pi_i, \pi^{\#}_i)} (X_i , \cat{O}_i) \right)_{i \in |\cat{I}|} \right)
\end{equation}
is the limit of the functor $D : \cat{I} \to \cat{D}$.

\item
Let $x \in X$. For each $i \in |\cat{I}|$, let us write $x_i := \pi_i(x) \in X_i$. Then
\[\begin{tikzcd}[row sep=tiny]
	{\cat{I}^{\op}} & {\ub{CRing}} \\
	{i^{\op}} & {\cat{O}_{i, x_i}} & {(\text{on objects})} \\
	{(i \xto{\alpha} j)^{\op}} & {\left( \cat{O}_{j, x_j} \xto{f^{\#}_{\alpha,x_i}} \cat{O}_{i,x_i} \right)} & {(\text{on morphisms})}
	\arrow["{D_x}", from=1-1, to=1-2]
	\arrow[maps to, from=2-1, to=2-2]
	\arrow[maps to, from=3-1, to=3-2]
\end{tikzcd}\]
defines a functor, and
\begin{equation}
\left( \cat{O}_x,
\left( \cat{O}_{i, x_i} \xto{\pi^{\#}_{i,x}} \cat{O}_x \right)_{i \in |\cat{I}^{\op}|} \right)
\end{equation}
is the colimit of this functor $D_x : \cat{I}^{\op} \to \ub{CRing}$.
\end{enumerate}
\end{prop}

\begin{cor} \label{cor:completeness of the category D}
The category $\cat{D}$ is complete. The forgetful functor $\cat{D} \to \ub{Top}$ preserves limits.
\end{cor}

\begin{proof}[Proof of \cref{prop:construction of limits in the category D}]
$\\$
\begin{enumerate}
\item
Since the diagram
\[\begin{tikzcd}
	& {X_i} \\
	X & {X_j}
	\arrow["{f_{\alpha}}", from=1-2, to=2-2]
	\arrow["{\pi_i}", from=2-1, to=1-2]
	\arrow["{\pi_j}"', from=2-1, to=2-2]
\end{tikzcd}\]
is commutative in $\ub{Top}$, we have
\begin{equation}
f_{\alpha,*} \, \pi_{i,*} \tilde{\cat{O}}_i = \pi_{j,*} \tilde{\cat{O}}_i .
\end{equation}
Therefore we have a morphism of sheaves of condensed rings on $X_j$
\[\begin{tikzcd}
	{\cat{O}_j} & {f_{\alpha,*} \cat{O}_i} & {f_{\alpha,*} \, \pi_{i,*} \tilde{\cat{O}}_i = \pi_{j,*} \tilde{\cat{O}}_i .}
	\arrow["{f^{\#}_{\alpha}}", from=1-1, to=1-2]
	\arrow["{f_{\alpha,*} \, \rho_i}", from=1-2, to=1-3]
\end{tikzcd}\]
Since $\tilde{\cat{O}}_j$ is the inverse image of $\cat{O}_j$ under $\pi_j : X \to X_j$, the assertion of (1) follows.

\item
Suppose that $i \xto{\alpha} j \xto{\beta} k$ are morphisms in $\cat{I}$. Then the following diagram is commutative in $\ub{Top}$.
\[\begin{tikzcd}
	X & {X_j} & {X_i} \\
	& {X_k}
	\arrow["{\pi_j}", from=1-1, to=1-2]
	\arrow["{\pi_k}"', from=1-1, to=2-2]
	\arrow["{f_{\beta}}", from=1-2, to=2-2]
	\arrow["{f_{\alpha}}"', from=1-3, to=1-2]
	\arrow["{f_{\beta \of \alpha}}", from=1-3, to=2-2]
\end{tikzcd}\]
It follows that the following diagram is commutative in $\ub{CSh}(X_k)$.
\[\begin{tikzcd}
	{\pi_{k,*} \tilde{\cat{O}}_i} && {f_{\beta \of \alpha,*} \cat{O}_i} \\
	{\pi_{k,*} \tilde{\cat{O}}_j} && {f_{\beta,*} \cat{O}_j} \\
	{\pi_{k,*} \tilde{\cat{O}}_k} && {\cat{O}_k}
	\arrow["{f_{\beta \of \alpha,*} \rho_i}"', from=1-3, to=1-1]
	\arrow["{\pi_{k,*} \tilde{f}^{\#}_{\alpha}}", from=2-1, to=1-1]
	\arrow["{f_{\beta,*} f^{\#}_{\alpha}}", from=2-3, to=1-3]
	\arrow["{f_{\beta,*} \rho_j}", from=2-3, to=2-1]
	\arrow["{\pi_{k,*} \tilde{f}^{\#}_{\beta}}", from=3-1, to=2-1]
	\arrow["{f^{\#}_{\beta \of \alpha}}"', curve={height=30pt}, from=3-3, to=1-3]
	\arrow["{f^{\#}_{\beta}}", from=3-3, to=2-3]
	\arrow["{\rho_k}", from=3-3, to=3-1]
\end{tikzcd}\]
By the definition of $\tilde{f}^{\#}_{\beta \of \alpha}$, we conclude that $\tilde{f}^{\#}_{\beta \of \alpha} = \tilde{f}^{\#}_{\alpha} \of \tilde{f}^{\#}_{\beta}$. 

If $i \in |\cat{I}|$, then the following diagram is commutative in $\ub{CSh}(X_i)$.
\[\begin{tikzcd}
	{\pi_{i,*} \tilde{\cat{O}}_i} && {f_{\id_i,*} \cat{O}_i} \\
	{\pi_{i,*} \tilde{\cat{O}}_i} && {\cat{O}_i}
	\arrow["{f_{\id_i,*} \rho_i \,=\, \rho_i}"', from=1-3, to=1-1]
	\arrow["{\pi_{i,*} \id_{\tilde{\cat{O}}_i}}", from=2-1, to=1-1]
	\arrow["{f^{\#}_{\id_i} \,=\, \id_{\cat{O}_i}}"', from=2-3, to=1-3]
	\arrow["{\rho_i}", from=2-3, to=2-1]
\end{tikzcd}\]
By the definition of $\tilde{f}^{\#}_{\id_i}$, we conclude that $\tilde{f}^{\#}_{\id_i} = \id_{\tilde{\cat{O}}_i}$.

\item
First of all, we prove that
\begin{equation}
\left( (X,\cat{O}),
\left( (X,\cat{O}) \xto{(\pi_i, \pi^{\#}_i)} (X_i , \cat{O}_i) \right)_{i \in |\cat{I}|} \right)
\end{equation}
is a cone on the functor $D : \cat{I} \to \cat{D}$. Let $\alpha : i \to j$ be any morphism in $\cat{I}$. Then the following diagram is commutative in $\ub{Top}$.
\[\begin{tikzcd}
	X & {X_i} \\
	& {X_j}
	\arrow["{\pi_i}", from=1-1, to=1-2]
	\arrow["{\pi_j}"', from=1-1, to=2-2]
	\arrow["{f_{\alpha}}", from=1-2, to=2-2]
\end{tikzcd}\]
Moreover, the following diagram is commutative in $\ub{CSh}(X)$. 
\[\begin{tikzcd}
	{\cat{O}} & {\tilde{\cat{O}}_i} \\
	& {\tilde{\cat{O}}_j}
	\arrow["{\sigma_i}"', from=1-2, to=1-1]
	\arrow["{\sigma_j}", from=2-2, to=1-1]
	\arrow["{\tilde{f}^{\#}_{\alpha}}"', from=2-2, to=1-2]
\end{tikzcd}\]
It follows that the following diagram is commutative in $\ub{CSh}(X_j)$.
\[\begin{tikzcd}
	{\pi_{j,*} \cat{O}} & {\pi_{j,*} \tilde{\cat{O}}_i} & {f_{\alpha,*} \cat{O}_i} \\
	{\pi_{j,*} \cat{O}} & {\pi_{j,*} \tilde{\cat{O}}_j} & {\cat{O}_j}
	\arrow["{\pi_{j,*} \sigma_i}"', from=1-2, to=1-1]
	\arrow["{f_{\alpha,*} \pi^{\#}_i}"', curve={height=30pt}, from=1-3, to=1-1]
	\arrow["{f_{\alpha,*} \rho_i}"', from=1-3, to=1-2]
	\arrow[equals, from=2-1, to=1-1]
	\arrow["{\pi_{j,*} \tilde{f}^{\#}_{\alpha}}"', from=2-2, to=1-2]
	\arrow["{\pi_{j,*} \sigma_j}", from=2-2, to=2-1]
	\arrow["{f^{\#}_{\alpha}}"', from=2-3, to=1-3]
	\arrow["{\pi^{\#}_j}", curve={height=-30pt}, from=2-3, to=2-1]
	\arrow["{\rho_j}", from=2-3, to=2-2]
\end{tikzcd}\]
Therefore the following diagram is commutative in $\cat{D}$.
\[\begin{tikzcd}
	{(X,\cat{O})} & {(X_i,\cat{O}_i)} \\
	& {(X_j,\cat{O}_j)}
	\arrow["{(\pi_i,\pi^{\#}_i)}", from=1-1, to=1-2]
	\arrow["{(\pi_j,\pi^{\#}_j)}"', from=1-1, to=2-2]
	\arrow["{(f_{\alpha}, f^{\#}_{\alpha})}", from=1-2, to=2-2]
\end{tikzcd}\]
This shows that
\begin{equation}
\left( (X,\cat{O}),
\left( (X,\cat{O}) \xto{(\pi_i, \pi^{\#}_i)} (X_i , \cat{O}_i) \right)_{i \in |\cat{I}|} \right)
\end{equation}
is a cone on the functor $D : \cat{I} \to \cat{D}$.

Next suppose that
\begin{equation}
\left( (Y,\cat{O}_Y),
\left( (Y,\cat{O}_Y) \xto{(\mu_i, \mu^{\#}_i)} (X_i , \cat{O}_i) \right)_{i \in |\cat{I}|} \right)
\end{equation}
is a cone on $D : \cat{I} \to \cat{D}$. We prove that there exists a unique morphism $(\nu, \nu^{\#}) : (Y,\cat{O}_Y) \to (X, \cat{O})$ in $\cat{D}$ such that the following diagram is commutative for every $i \in |\cat{I}|$. 
\[\begin{tikzcd}
	{(Y, \cat{O}_Y)} & {(X,\cat{O})} \\
	& {(X_i,\cat{O}_i)}
	\arrow["{(\nu,\nu^{\#})}", from=1-1, to=1-2]
	\arrow["{(\mu,\mu^{\#}_i)}"', from=1-1, to=2-2]
	\arrow["{(\pi_i,\pi^{\#}_i)}", from=1-2, to=2-2]
\end{tikzcd}\]

First of all, we have a cone
\begin{equation}
\left( Y, ( Y \xto{\mu_i} X_i )_{i \in |\cat{I}|} \right)
\end{equation}
on the functor $F \of D : \cat{I} \to \ub{Top}$. Since
\begin{equation}
\left( X, (X \xto{\pi_i} X_i)_{i \in |\cat{I}|} \right)
\end{equation}
is the limit of the functor $F \of D : \cat{I} \to \ub{Top}$, there exists a unique contnuous map $\nu : Y \to X$ such that the following diagram is commutative in $\ub{Top}$ for every $i \in |\cat{I}|$.
\[\begin{tikzcd}
	Y & X \\
	& {X_i}
	\arrow["\nu", from=1-1, to=1-2]
	\arrow["{\mu_i}"', from=1-1, to=2-2]
	\arrow["{\pi_i}", from=1-2, to=2-2]
\end{tikzcd}\]
Then, for every $i \in \ob{I}$, we have a morphism of sheaves of condensed rings on $X_i$
\[\begin{tikzcd}
	{\cat{O}_i} & {\mu_{i,*} \cat{O}_Y = \pi_{i,*} \, \nu_* \cat{O}_Y .}
	\arrow["{\mu^{\#}_i}", from=1-1, to=1-2]
\end{tikzcd}\]
Since $\tilde{\cat{O}}_i$ is the inverse image of $\cat{O}_i$ under $\pi_i : X \to X_i$, there exists a unique morphism $\tilde{\mu}^{\#}_i : \tilde{\cat{O}}_i \to \nu_* \cat{O}_Y$ of sheaves of condensed rings on $X$ such that the following diagram is commutative.
\[\begin{tikzcd}
	{\pi_{i,*} \nu_* \cat{O}_Y} \\
	{\pi_{i,*} \tilde{\cat{O}}_i} & {\cat{O}_i}
	\arrow["{\pi_{i,*} \tilde{\mu}^{\#}_i}", from=2-1, to=1-1]
	\arrow["{\mu^{\#}_i}"', from=2-2, to=1-1]
	\arrow["{\rho_i}", from=2-2, to=2-1]
\end{tikzcd}\]
If $i \xto{\alpha} j$ is a morphism in $\cat{I}$, then the following diagram is commutative in $\ub{Top}$.
\[\begin{tikzcd}
	& {X_i} \\
	X & {X_j}
	\arrow["{f_{\alpha}}", from=1-2, to=2-2]
	\arrow["{\pi_i}", from=2-1, to=1-2]
	\arrow["{\pi_j}"', from=2-1, to=2-2]
\end{tikzcd}\]
It follows that the following diagram is commutative in $\ub{CSh}(X_j)$.
\[\begin{tikzcd}
	{\pi_{j,*} \nu_* \cat{O}_Y} & {\pi_{j,*} \nu_* \cat{O}_Y} \\
	{\pi_{j,*} \tilde{\cat{O}}_i} & {f_{\alpha,*} \cat{O}_i} \\
	{\pi_{j,*} \tilde{\cat{O}}_j} & {\cat{O}_j}
	\arrow[equals, from=1-2, to=1-1]
	\arrow["{\pi_{j,*} \tilde{\mu}^{\#}_i}", from=2-1, to=1-1]
	\arrow["{f_{\alpha,*} \mu^{\#}_i}", from=2-2, to=1-2]
	\arrow["{f_{\alpha,*} \rho_i}", from=2-2, to=2-1]
	\arrow["{\pi_{j,*} \tilde{f}^{\#}_{\alpha}}", from=3-1, to=2-1]
	\arrow["{\mu^{\#}_j}"', curve={height=30pt}, from=3-2, to=1-2]
	\arrow["{f^{\#}_{\alpha}}", from=3-2, to=2-2]
	\arrow["{\rho_j}", from=3-2, to=3-1]
\end{tikzcd}\]
By the definition of $\tilde{\mu}^{\#}_j$, we conclude that $\tilde{\mu}^{\#}_j = \tilde{\mu}^{\#}_i \of \tilde{f}^{\#}_{\alpha}$. Thus
\begin{equation}
\left( \nu_* \cat{O}_Y , 
( \tilde{\cat{O}}_i \xto{\tilde{\mu}^{\#}_i} \nu_* \cat{O}_Y )_{i \in |\cat{I}^{\op}|} \right)
\end{equation}
is a cocone on the functor $\tilde{D} : \cat{I}^{\op} \to \ub{CSh}(X)$. Since
\begin{equation}
\left(\cat{O} , ( \tilde{\cat{O}}_i \xto{\sigma_i} \cat{O} )_{i \in |\cat{I}^{\op}|} \right)
\end{equation}
is the colimit of the functor $\tilde{D} : \cat{I}^{\op} \to \ub{CSh}(X)$, there exists a unique morphism $\nu^{\#} : \cat{O} \to \nu_* \cat{O}_Y$ of sheaves of condensed rings on $X$ such that the following diagram is commutative in $\ub{CSh}(X)$ for every $i \in \ob{I}$.
\[\begin{tikzcd}
	{\nu_* \cat{O}_Y} & {\cat{O}} \\
	& {\tilde{\cat{O}}_i}
	\arrow["{\nu^{\#}}"', from=1-2, to=1-1]
	\arrow["{\tilde{\mu}^{\#}_i}", from=2-2, to=1-1]
	\arrow["{\sigma_i}"', from=2-2, to=1-2]
\end{tikzcd}\]
Then the following diagram is commutative in $\ub{CSh}(X_i)$ for every $i \in \ob{I}$.
\[\begin{tikzcd}
	& {\pi_{i,*} \cat{O}} \\
	{\pi_{i,*}  \nu_* \cat{O}_Y} & {\pi_{i,*} \tilde{\cat{O}}_i} \\
	& {\cat{O}_i}
	\arrow["{\pi_{i,*} \nu^{\#}}"', curve={height=12pt}, from=1-2, to=2-1]
	\arrow["{\pi_{i,*} \sigma_i}", from=2-2, to=1-2]
	\arrow["{\pi_{i,*} \tilde{\mu}^{\#}_i}", from=2-2, to=2-1]
	\arrow["{\pi^{\#}_i}"', curve={height=30pt}, from=3-2, to=1-2]
	\arrow["{\mu^{\#}_i}", curve={height=-12pt}, from=3-2, to=2-1]
	\arrow["{\rho_i}", from=3-2, to=2-2]
\end{tikzcd}\]
Thus we obtain a morphism $(\nu, \nu^{\#}) : (Y,\cat{O}_Y) \to (X, \cat{O})$ in $\cat{D}$ such that the following diagram is commutative in $\cat{D}$ for every $i \in |\cat{I}|$. 
\[\begin{tikzcd}
	{(Y, \cat{O}_Y)} & {(X,\cat{O})} \\
	& {(X_i,\cat{O}_i)}
	\arrow["{(\nu,\nu^{\#})}", from=1-1, to=1-2]
	\arrow["{(\mu,\mu^{\#}_i)}"', from=1-1, to=2-2]
	\arrow["{(\pi_i,\pi^{\#}_i)}", from=1-2, to=2-2]
\end{tikzcd}\]

On the other hand, suppose that $(\lambda, \lambda^{\#}) : (Y,\cat{O}_Y) \to (X, \cat{O})$ is another morphism in $\cat{D}$ such that the following diagram is commutative in $\cat{D}$ for every $i \in |\cat{I}|$. 
\[\begin{tikzcd}
	{(Y, \cat{O}_Y)} & {(X,\cat{O})} \\
	& {(X_i,\cat{O}_i)}
	\arrow["{(\lambda, \lambda^{\#})}", from=1-1, to=1-2]
	\arrow["{(\mu,\mu^{\#}_i)}"', from=1-1, to=2-2]
	\arrow["{(\pi_i,\pi^{\#}_i)}", from=1-2, to=2-2]
\end{tikzcd}\]
Then the following diagram is commutative in $\ub{Top}$ for every $i \in \ob{I}$.
\[\begin{tikzcd}
	Y & X \\
	& {X_i}
	\arrow["\lambda", from=1-1, to=1-2]
	\arrow["{\mu_i}"', from=1-1, to=2-2]
	\arrow["{\pi_i}", from=1-2, to=2-2]
\end{tikzcd}\]
By the definition of $\nu$, we conclude that $\lambda = \nu$. Then, for every $i \in \ob{I}$, the following diagram is commutative in $\ub{CSh}(X_i)$.
\[\begin{tikzcd}
	{\pi_{i,*} \cat{O}} & {\pi_{i,*}  \nu_* \cat{O}_Y} \\
	{\pi_{i,*} \tilde{\cat{O}}_i} & {\cat{O}_i}
	\arrow["{\pi_{i,*} \lambda^{\#}}", from=1-1, to=1-2]
	\arrow["{\pi_{i,*} \sigma_i}", from=2-1, to=1-1]
	\arrow["{\pi^{\#}_i}", from=2-2, to=1-1]
	\arrow["{\mu^{\#}_i}"', from=2-2, to=1-2]
	\arrow["{\rho_i}", from=2-2, to=2-1]
\end{tikzcd}\]
By the definition of $\tilde{\mu}^{\#}_i$, we conclude that $\lambda^{\#} \of \sigma_i = \tilde{\mu}^{\#}_i$. In other words, the following diagram is commutative in $\ub{CSh}(X)$.
\[\begin{tikzcd}
	{\nu_* \cat{O}_Y} & {\cat{O}} \\
	& {\tilde{\cat{O}}_i}
	\arrow["{\lambda^{\#}}"', from=1-2, to=1-1]
	\arrow["{\tilde{\mu}^{\#}_i}", from=2-2, to=1-1]
	\arrow["{\sigma_i}"', from=2-2, to=1-2]
\end{tikzcd}\]
This holds for every $i \in \ob{I}$. By the definition of $\nu^{\#}$, we conclude that $\lambda^{\#} = \nu^{\#}$. Consequently, we have $(\lambda, \lambda^{\#}) = (\nu, \nu^{\#})$. This completes the proof of the assertion (3).

\item
Since 
\begin{equation}
\left(\cat{O} , ( \tilde{\cat{O}}_i \xto{\sigma_i} \cat{O} )_{i \in |\cat{I}^{\op}|} \right)
\end{equation}
is the colimit of the functor $\tilde{D} : \cat{I}^{\op} \to \ub{CSh}(X)$, \cref{cor:right adjoint of the stalk functor} implies that
\begin{equation}
\left(\cat{O}_x , ( \tilde{\cat{O}}_{i,x} \xto{\sigma_{i,x}} \cat{O}_x )_{i \in |\cat{I}^{\op}|} \right)
\end{equation}
is the colimit of the following functor $\tilde{D}_x : \cat{I}^{\op} \to \ub{CRing}$.
\[\begin{tikzcd}[row sep=tiny]
	{\cat{I}^{\op}} & {\ub{CRing}} \\
	{i^{\op}} & {\tilde{\cat{O}}_{i,x}} & {(\text{on objects})} \\
	{\alpha^{\op}} & {\tilde{f}^{\#}_{\alpha,x}} & {(\text{on morphisms})}
	\arrow["{\tilde{D}_x}", from=1-1, to=1-2]
	\arrow[maps to, from=2-1, to=2-2]
	\arrow[maps to, from=3-1, to=3-2]
\end{tikzcd}\]

On the other hand, \cref{cor:existence of inverse image functor} shows that for every $i \in \ob{I}$, the homomorphism $\rho_{i,x} : \cat{O}_{i,x_i} \to \tilde{\cat{O}}_{i,x}$ is an isomorphism of condensed rings. Moreover, for every morphism $i \xto{\alpha} j$ in $\cat{I}$, the commutativity of the diagram 
\[\begin{tikzcd}
	{\pi_{j,*} \tilde{\cat{O}}_i} && {f_{\alpha,*} \cat{O}_i} \\
	{\pi_{j,*} \tilde{\cat{O}}_j} && {\cat{O}_j}
	\arrow["{f_{\alpha,*} \rho_i}"', from=1-3, to=1-1]
	\arrow["{\pi_{j,*} \tilde{f}^{\#}_{\alpha}}", from=2-1, to=1-1]
	\arrow["{f^{\#}_{\alpha}}"', from=2-3, to=1-3]
	\arrow["{\rho_j}", from=2-3, to=2-1]
\end{tikzcd}\]
shows that the following diagram is commutative in $\ub{CRing}$.
\[\begin{tikzcd}
	{\tilde{\cat{O}}_{i,x}} && {\cat{O}_{i,x_i}} \\
	{\tilde{\cat{O}}_{j,x}} && {\cat{O}_{j,x_j}}
	\arrow["{\rho_{i,x}}"', from=1-3, to=1-1]
	\arrow["{\tilde{f}^{\#}_{\alpha,x}}", from=2-1, to=1-1]
	\arrow["{f^{\#}_{\alpha,x_i}}"', from=2-3, to=1-3]
	\arrow["{\rho_{j,x}}", from=2-3, to=2-1]
\end{tikzcd}\]
It follows that 
\[\begin{tikzcd}[row sep=tiny]
	{\cat{I}^{\op}} & {\ub{CRing}} \\
	{i^{\op}} & {\cat{O}_{i, x_i}} & {(\text{on objects})} \\
	{\alpha^{\op}} & {f^{\#}_{\alpha,x_i}} & {(\text{on morphisms})}
	\arrow["{D_x}", from=1-1, to=1-2]
	\arrow[maps to, from=2-1, to=2-2]
	\arrow[maps to, from=3-1, to=3-2]
\end{tikzcd}\]
defines a functor $D_x : \cat{I}^{\op} \to \ub{CRing}$ which is isomorphic to $\tilde{D}_x : \cat{I}^{\op} \to \ub{CRing}$ via the isomorphism
\begin{equation}
\left( \cat{O}_{i,x_i} \xto{\rho_{i,x}} \tilde{\cat{O}}_{i,x} \right)_{i \in |\cat{I}^{\op}|} .
\end{equation}
Consequently,
\begin{equation}
\left(\cat{O}_x , 
( \cat{O}_{i,x_i} \xto{\rho_{i,x}} \tilde{\cat{O}}_{i,x} \xto{\sigma_{i,x}} \cat{O}_x )
_{i \in |\cat{I}^{\op}|} \right)
\end{equation}
is the colimit of the functor $D_x : \cat{I}^{\op} \to \ub{CRing}$. Moreover, for every $i \in \ob{I}$, the commutativity of the diagram
\[\begin{tikzcd}
	{\cat{O}_i} & {\pi_{i,*} \tilde{\cat{O}}_i} \\
	& {\pi_{i,*} \cat{O}}
	\arrow["{\rho_i}", from=1-1, to=1-2]
	\arrow["{\pi^{\#}_i}"', from=1-1, to=2-2]
	\arrow["{\pi_{i,*} \sigma_i}", from=1-2, to=2-2]
\end{tikzcd}\]
shows that the following diagram is commutative in $\ub{CRing}$.
\[\begin{tikzcd}
	{\cat{O}_{i,x_i}} & {\tilde{\cat{O}}_{i,x}} \\
	& {\cat{O}_x}
	\arrow["{\rho_{i,x}}", from=1-1, to=1-2]
	\arrow["{\pi^{\#}_{i,x}}"', from=1-1, to=2-2]
	\arrow["{\sigma_{i,x}}", from=1-2, to=2-2]
\end{tikzcd}\]
Therefore
\begin{equation}
\left(\cat{O}_x , ( \cat{O}_{i,x_i} \xto{\pi^{\#}_{i,x}} \cat{O}_x )_{i \in |\cat{I}^{\op}|} \right)
\end{equation}
is the colimit of the functor $D_x : \cat{I}^{\op} \to \ub{CRing}$.
\end{enumerate}
\end{proof}

\subsubsection{Cocompleteness}

\begin{prop} \label{prop:construction of colimits in the category D}
Let $\cat{I}$ be a small category. Let $D : \cat{I} \to \cat{D}$ be a functor. For each $i \in |\cat{I}|$, let us write $(X_i , \cat{O}_i) := D(i)$. For each morphism $\alpha$ in $\cat{I}$, let us write $(f_{\alpha}, f^{\#}_{\alpha}) := D(\alpha)$. Let
\begin{equation}
\left( X, (X_i \xto{\pi_i} X)_{i \in |\cat{I}|} \right)
\end{equation}
be the colimit of the functor $F \of D : \cat{I} \to \ub{Top}$.
\begin{enumerate}
\item
For each $i \in |\cat{I}|$, let $\tilde{\cat{O}}_i$ be the direct image of $\cat{O}_i$ under $\pi_i : X_i \to X$. For each morphism $\alpha : i \to j$ in $\cat{I}$, let us write $\tilde{f}^{\#}_{\alpha} : \tilde{\cat{O}}_j \to \tilde{\cat{O}}_i$ for the morphism $\pi_{j,*} f^{\#}_{\alpha} : \pi_{j,*} \cat{O}_j \to \pi_{i,*} \cat{O}_i$ of sheaves of condensed rings on $X$. Then
\[\begin{tikzcd}[row sep=tiny]
	{\cat{I}^{\op}} & {\ub{CSh}(X)} \\
	{i^{\op}} & {\tilde{\cat{O}}_i} & {(\text{on objects})} \\
	{\alpha^{\op}} & {\tilde{f}^{\#}_{\alpha}} & {(\text{on morphisms})}
	\arrow["{\tilde{D}}", from=1-1, to=1-2]
	\arrow[maps to, from=2-1, to=2-2]
	\arrow[maps to, from=3-1, to=3-2]
\end{tikzcd}\]
defines a functor.

\item
By \cref{cor:CSh(X) is complete and cocomplete}, the category $\ub{CSh}(X)$ is complete. Let
\begin{equation}
\left(\cat{O} , ( \cat{O} \xto{\pi^{\#}_i} \tilde{\cat{O}}_i )_{i \in |\cat{I}^{\op}|} \right)
\end{equation}
be the limit of the functor $\tilde{D} : \cat{I}^{\op} \to \ub{CSh}(X)$. Then
\begin{equation}
\left( (X,\cat{O}),
\left( (X_i , \cat{O}_i) \xto{(\pi_i, \pi^{\#}_i)} (X,\cat{O}) \right)_{i \in |\cat{I}|} \right)
\end{equation}
is the colimit of the functor $D : \cat{I} \to \cat{D}$.
\end{enumerate}
\end{prop}

\begin{cor} \label{cor:cocompleteness of the category D}
The category $\cat{D}$ is cocomplete. The forgetful functor $\cat{D} \to \ub{Top}$ preserves colimits.
\end{cor}

\begin{proof}[Proof of \cref{prop:construction of colimits in the category D}]
$\\$
\begin{enumerate}
\item
If $\alpha : i \to j$ is a morphism in $\cat{I}$, then the diagram
\[\begin{tikzcd}
	{X_i} \\
	{X_j} & X
	\arrow["{f_{\alpha}}"', from=1-1, to=2-1]
	\arrow["{\pi_i}", from=1-1, to=2-2]
	\arrow["{\pi_j}"', from=2-1, to=2-2]
\end{tikzcd}\]
is commutative in $\ub{Top}$. Thus $\pi_{j,*} f_{\alpha,*} \cat{O}_i = \pi_{i,*} \cat{O}_i$, and we have a morphism of sheaves of condensed rings on $X$
\[\begin{tikzcd}
	{\tilde{f}^{\#}_{\alpha} = \pi_{j,*} f^{\#}_{\alpha} : \pi_{j,*} \cat{O}_j} & {\pi_{j,*} f_{\alpha,*} \cat{O}_i = \pi_{i,*} \cat{O}_i .}
	\arrow[from=1-1, to=1-2]
\end{tikzcd}\]

Suppose that $i \xto{\alpha} j \xto{\beta} k$ are morphisms in $\cat{I}$. Then the following diagram is commutative in $\ub{CSh}(X_k)$.
\[\begin{tikzcd}
	{f_{\beta \of \alpha, *} \cat{O}_i} & {f_{\beta,*} \cat{O}_j} \\
	& {\cat{O}_k}
	\arrow["{f_{\beta,*} f^{\#}_{\alpha}}"', from=1-2, to=1-1]
	\arrow["{f^{\#}_{\beta \of \alpha}}", from=2-2, to=1-1]
	\arrow["{f^{\#}_{\beta}}"', from=2-2, to=1-2]
\end{tikzcd}\]
Applying the functor $\pi_{k,*} : \ub{CSh}(X_k) \to \ub{CSh}(X)$, we conclude that the diagram 
\[\begin{tikzcd}
	{\pi_{i,*} \cat{O}_i} & {\pi_{j,*} \cat{O}_j} \\
	& {\pi_{k,*} \cat{O}_k}
	\arrow["{\pi_{j,*} f^{\#}_{\alpha}}"', from=1-2, to=1-1]
	\arrow["{\pi_{k,*} f^{\#}_{\beta \of \alpha}}", from=2-2, to=1-1]
	\arrow["{\pi_{k,*} f^{\#}_{\beta}}"', from=2-2, to=1-2]
\end{tikzcd}\]
is commutative in $\ub{CSh}(X)$, which shows that the diagram
\[\begin{tikzcd}
	{\tilde{\cat{O}}_i} & {\tilde{\cat{O}}_j} \\
	& {\tilde{\cat{O}}_k}
	\arrow["{\tilde{f}^{\#}_{\alpha}}"', from=1-2, to=1-1]
	\arrow["{\tilde{f}^{\#}_{\beta \of \alpha}}", from=2-2, to=1-1]
	\arrow["{\tilde{f}^{\#}_{\beta}}"', from=2-2, to=1-2]
\end{tikzcd}\]
is commutative in $\ub{CSh}(X)$.

If $i \in |\cat{I}|$, then
\[\begin{tikzcd}
	{f^{\#}_{\id_i} = \id_{\cat{O}_i} : \cat{O}_i} & {f_{\id_i,*} \cat{O}_i = \cat{O}_i .}
	\arrow[from=1-1, to=1-2]
\end{tikzcd}\]
Applying the functor $\pi_{i,*} : \ub{CSh}(X_i) \to \ub{CSh}(X)$, we conclude that $\pi_{i,*} f^{\#}_{\id_i} = \id_{\pi_{i,*} \cat{O}_i}$. Therefore
\begin{equation}
\tilde{f}^{\#}_{\id_i} = \pi_{i,*} f^{\#}_{\id_i}
= \id_{\pi_{i,*} \cat{O}_i} = \id_{\tilde{\cat{O}}_i} .
\end{equation}
This completes the proof of the assertion (1).

\item
First of all, we prove that
\begin{equation}
\left( (X,\cat{O}),
\left( (X_i , \cat{O}_i) \xto{(\pi_i, \pi^{\#}_i)} (X,\cat{O}) \right)_{i \in |\cat{I}|} \right)
\end{equation}
is a cocone on the functor $D : \cat{I} \to \cat{D}$. Let $\alpha : i \to j$ be any morphism in $\cat{I}$. Then the following diagram is commutative in $\ub{Top}$.
\[\begin{tikzcd}
	{X_i} & X \\
	{X_j}
	\arrow["{\pi_i}", from=1-1, to=1-2]
	\arrow["{f_{\alpha}}"', from=1-1, to=2-1]
	\arrow["{\pi_j}"', from=2-1, to=1-2]
\end{tikzcd}\]
Moreover, the following diagram is commutative in $\ub{CSh}(X)$.
\[\begin{tikzcd}
	{\pi_{i,*} \cat{O}_i} & {\tilde{\cat{O}}_i} & {\cat{O}} \\
	{\pi_{j,*} \cat{O}_j} & {\tilde{\cat{O}}_j}
	\arrow[equals, from=1-2, to=1-1]
	\arrow["{\pi^{\#}_i}"', from=1-3, to=1-2]
	\arrow["{\pi^{\#}_j}", from=1-3, to=2-2]
	\arrow["{\pi_{j,*} f^{\#}_{\alpha}}", from=2-1, to=1-1]
	\arrow["{\tilde{f}^{\#}_{\alpha}}", from=2-2, to=1-2]
	\arrow[equals, from=2-2, to=2-1]
\end{tikzcd}\]
Therefore the following diagram is commutative in $\cat{D}$.
\[\begin{tikzcd}
	{(X_i,\cat{O}_i)} & {(X,\cat{O})} \\
	{(X_j,\cat{O}_j)}
	\arrow["{(\pi_i,\pi^{\#}_i)}", from=1-1, to=1-2]
	\arrow["{(f_{\alpha},f^{\#}_{\alpha})}"', from=1-1, to=2-1]
	\arrow["{(\pi_j,\pi^{\#}_j)}"', from=2-1, to=1-2]
\end{tikzcd}\]
This shows that
\begin{equation}
\left( (X,\cat{O}),
\left( (X_i , \cat{O}_i) \xto{(\pi_i, \pi^{\#}_i)} (X,\cat{O}) \right)_{i \in |\cat{I}|} \right)
\end{equation}
is a cocone on the functor $D : \cat{I} \to \cat{D}$.

Next suppose that
\begin{equation}
\left( (Y,\cat{O}_Y),
\left( (X_i , \cat{O}_i) \xto{(\mu_i, \mu^{\#}_i)} (Y,\cat{O}_Y) \right)_{i \in |\cat{I}|} \right)
\end{equation}
is a cocone on $D : \cat{I} \to \cat{D}$. We prove that there exists a unique morphism $(\nu, \nu^{\#}) : (X, \cat{O}) \to (Y,\cat{O}_Y)$ in $\cat{D}$ such that the following diagram is commutative for every $i \in |\cat{I}|$. 
\[\begin{tikzcd}
	{(X,\cat{O})} & {(Y,\cat{O}_Y)} \\
	{(X_i,\cat{O}_i)}
	\arrow["{(\nu,\nu^{\#})}", from=1-1, to=1-2]
	\arrow["{(\pi_i,\pi^{\#}_i)}", from=2-1, to=1-1]
	\arrow["{(\mu_i,\mu^{\#}_i)}"', from=2-1, to=1-2]
\end{tikzcd}\]

First of all, we have a cocone
\begin{equation}
\left( Y, ( X_i \xto{\mu_i} Y )_{i \in |\cat{I}|} \right)
\end{equation}
on the functor $F \of D : \cat{I} \to \ub{Top}$. Since
\begin{equation}
\left( X, (X_i \xto{\pi_i} X)_{i \in |\cat{I}|} \right)
\end{equation}
is the colimit of the functor $F \of D : \cat{I} \to \ub{Top}$. there exists a unique contnuous map $\nu : X \to Y$ such that the following diagram is commutative in $\ub{Top}$ for every $i \in |\cat{I}|$.
\[\begin{tikzcd}
	X & Y \\
	{X_i}
	\arrow["\nu", from=1-1, to=1-2]
	\arrow["{\pi_i}", from=2-1, to=1-1]
	\arrow["{\mu_i}"', from=2-1, to=1-2]
\end{tikzcd}\]
Then, for every $i \in \ob{I}$, we have a morphism of sheaves of condensed rings on $Y$
\[\begin{tikzcd}
	{\cat{O}_Y} & {\mu_{i,*} \cat{O}_i = \nu_* \, \pi_{i,*} \cat{O}_i = \nu_* \tilde{\cat{O}}_i .}
	\arrow["{\mu^{\#}_i}", from=1-1, to=1-2]
\end{tikzcd}\]
Moreover, for any morphism $\alpha : i \to j$ in $\cat{I}$, the following diagram is commutative in $\ub{CSh}(Y)$.
\[\begin{tikzcd}
	{\nu_* \tilde{\cat{O}}_i} & {\nu_* \pi_{i,*} \cat{O}_i} & {\mu_{i,*} \cat{O}_i} & {\cat{O}_Y} \\
	{\nu_* \tilde{\cat{O}}_j} & {\nu_* \pi_{j,*} \cat{O}_j} & {\mu_{j,*} \cat{O}_j}
	\arrow[equals, from=1-2, to=1-1]
	\arrow[equals, from=1-3, to=1-2]
	\arrow["{\mu^{\#}_i}"', from=1-4, to=1-3]
	\arrow["{\mu^{\#}_i}", from=1-4, to=2-3]
	\arrow["{\nu_* \tilde{f}^{\#}_{\alpha}}", from=2-1, to=1-1]
	\arrow["{\nu_* \pi_{j,*} f^{\#}_{\alpha}}", from=2-2, to=1-2]
	\arrow[equals, from=2-2, to=2-1]
	\arrow["{\mu_{j,*} f^{\#}_{\alpha}}", from=2-3, to=1-3]
	\arrow[equals, from=2-3, to=2-2]
\end{tikzcd}\]
On the other hand,
\begin{equation}
\left(\cat{O} , ( \cat{O} \xto{\pi^{\#}_i} \tilde{\cat{O}}_i )_{i \in |\cat{I}^{\op}|} \right)
\end{equation}
is the limit of the functor $\tilde{D} : \cat{I}^{\op} \to \ub{CSh}(X)$. The functor $\nu_* : \ub{CSh}(X) \to \ub{CSh}(Y)$ preserves limits since it has a left adjoint $\nu^{-1} : \ub{CSh}(Y) \to \ub{CSh}(X)$. Therefore 
\begin{equation}
\left(\nu_* \cat{O} , 
( \nu_* \cat{O} \xto{\nu_* \pi^{\#}_i} \nu_* \tilde{\cat{O}}_i )_{i \in |\cat{I}^{\op}|} \right)
\end{equation}
is the limit of the functor $\nu_* \of \tilde{D} : \cat{I}^{\op} \to \ub{CSh}(Y)$. Consequently, there exists a unique morphism $\nu^{\#} : \cat{O}_Y \to \nu_* \cat{O}$ of sheaves of condensed rings on $Y$ such that the following diagram is commutative in $\ub{CSh}(Y)$ for every $i \in \ob{I}$.
\[\begin{tikzcd}
	{\nu_* \cat{O}} & {\cat{O}_Y} \\
	{\nu_* \tilde{\cat{O}}_i} & {\mu_{i,*} \cat{O}_i}
	\arrow["{\nu_* \pi^{\#}_i}"', from=1-1, to=2-1]
	\arrow["{\nu^{\#}}"', from=1-2, to=1-1]
	\arrow["{\mu^{\#}_i}", from=1-2, to=2-2]
	\arrow[equals, from=2-2, to=2-1]
\end{tikzcd}\]
Thus we obtain a morphism $(\nu, \nu^{\#}) : (X, \cat{O}) \to (Y,\cat{O}_Y)$ in $\cat{D}$ such that the following diagram is commutative in $\cat{D}$ for every $i \in |\cat{I}|$. 
\[\begin{tikzcd}
	{(X,\cat{O})} & {(Y,\cat{O}_Y)} \\
	{(X_i,\cat{O}_i)}
	\arrow["{(\nu,\nu^{\#})}", from=1-1, to=1-2]
	\arrow["{(\pi_i,\pi^{\#}_i)}", from=2-1, to=1-1]
	\arrow["{(\mu_i,\mu^{\#}_i)}"', from=2-1, to=1-2]
\end{tikzcd}\]

On the other hand, suppose that $(\lambda, \lambda^{\#}) : (X, \cat{O}) \to (Y,\cat{O}_Y)$ is another morphism in $\cat{D}$ such that the following diagram is commutative in $\cat{D}$ for every $i \in |\cat{I}|$. 
\[\begin{tikzcd}
	{(X,\cat{O})} & {(Y,\cat{O}_Y)} \\
	{(X_i,\cat{O}_i)}
	\arrow["{(\lambda, \lambda^{\#})}", from=1-1, to=1-2]
	\arrow["{(\pi_i,\pi^{\#}_i)}", from=2-1, to=1-1]
	\arrow["{(\mu_i,\mu^{\#}_i)}"', from=2-1, to=1-2]
\end{tikzcd}\]
Then the following diagram is commutative in $\ub{Top}$ for every $i \in \ob{I}$.
\[\begin{tikzcd}
	X & Y \\
	{X_i}
	\arrow["\lambda", from=1-1, to=1-2]
	\arrow["{\pi_i}", from=2-1, to=1-1]
	\arrow["{\mu_i}"', from=2-1, to=1-2]
\end{tikzcd}\]
By the definition of $\nu$, we conclude that $\lambda = \nu$. Then, for every $i \in \ob{I}$, the following diagram is commutative in $\ub{CSh}(Y)$.
\[\begin{tikzcd}
	{\nu_* \cat{O}} & {\lambda_* \cat{O}} & {\cat{O}_Y} \\
	{\nu_* \tilde{\cat{O}}_i} & {\mu_{i,*} \cat{O}_i}
	\arrow["{\nu_* \pi^{\#}_i}"', from=1-1, to=2-1]
	\arrow[equals, from=1-2, to=1-1]
	\arrow["{\lambda_* \pi^{\#}_i}"', from=1-2, to=2-2]
	\arrow["{\lambda^{\#}}"', from=1-3, to=1-2]
	\arrow["{\mu^{\#}_i}", from=1-3, to=2-2]
	\arrow[equals, from=2-2, to=2-1]
\end{tikzcd}\]
By the definition of $\nu^{\#}$, we conclude that $\lambda^{\#} = \nu^{\#}$. Thus $(\lambda, \lambda^{\#}) = (\nu, \nu^{\#})$. This completes the proof of the assertion (2).
\end{enumerate}
\end{proof}

\subsection{The category $\cat{C}$}

\subsubsection{Definition}

\begin{df}
We define the category $\cat{C}$ as follows.
\begin{enumerate}
\item
An obejct of $\cat{C}$ is a triple $(X, \cat{O}_X, \cat{V}_X)$ consisting of the following objects.
\begin{enumerate}
\item $(X, \cat{O}_X)$ is an object of $\cat{D}$.
\item $\cat{V}_X$ is a family $(\cat{V}_{X,x})_{x \in X}$, where $\cat{V}_{X,x}$ is a set of continuous valuations on $\cat{O}_{X,x}$ for each $x \in X$.
\end{enumerate}

\item
A morphism $(X, \cat{O}_X, \cat{V}_X) \to (Y, \cat{O}_Y, \cat{V}_Y)$ in $\cat{C}$ is a morphism $(f,f^{\#}) : (X, \cat{O}_X) \to (Y, \cat{O}_Y)$ in $\cat{D}$ which satisfies the following condition: For every $x \in X$ and every $v \in \cat{V}_{X,x}$, the homomorphism $f^{\#}_x : \cat{O}_{Y,f(x)} \to \cat{O}_{X,x}$ of condensed rings satisfies $(f^{\#}_x)^{-1}(v) \in \cat{V}_{Y, f(x)}$.

\item
The composition of morphisms in $\cat{C}$ is the same as those in $\cat{D}$. In other words, the composition of morphisms $(f,f^{\#}) : (X, \cat{O}_X, \cat{V}_X) \to (Y, \cat{O}_Y, \cat{V}_Y)$ and $(g,g^{\#}) : (Y, \cat{O}_Y, \cat{V}_Y) \to (Z, \cat{O}_Z, \cat{V}_Z)$ in $\cat{C}$ is defined to be $(g \of f ,\, (g_* f^{\#}) \of g^{\#} )$.
\end{enumerate}
\end{df}

\begin{rem}
We have a forgetful functor
\[\begin{tikzcd}[row sep=tiny]
	{\cat{C}} & {\cat{D}} \\
	{(X, \cat{O}_X, \cat{V}_X)} & {(X, \cat{O}_X)} & {(\text{on objects})} \\
	{(f, f^{\#})} & {(f, f^{\#})} & {(\text{on morphisms})}
	\arrow[from=1-1, to=1-2]
	\arrow[maps to, from=2-1, to=2-2]
	\arrow[maps to, from=3-1, to=3-2]
\end{tikzcd}\]
In this subsection, we write $G$ for this forgetful functor.
\end{rem}

\subsubsection{Completeness}

\begin{prop} \label{prop:construction of limits in the category C}
Let $\cat{I}$ be a small category. Let $D : \cat{I} \to \cat{C}$ be a functor. For each $i \in |\cat{I}|$, let us write $(X_i , \cat{O}_i, \cat{V}_i) := D(i)$. For each morphism $\alpha$ in $\cat{I}$, let us write $(f_{\alpha}, f^{\#}_{\alpha}) := D(\alpha)$. Let
\begin{equation}
\left( (X,\cat{O}),
\left( (X,\cat{O}) \xto{(\pi_i, \pi^{\#}_i)} (X_i , \cat{O}_i) \right)_{i \in |\cat{I}|} \right)
\end{equation}
be the limit of the functor $G \of D : \cat{I} \to \cat{D}$. For each $x \in X$, let us define $\cat{V}_x$ to be the set of all continous valuations $v$ on $\cat{O}_x$ with the following property: For every $i \in |\cat{I}|$, the homomorphism $\pi^{\#}_{i,x} : \cat{O}_{i, \pi_i(x)} \to \cat{O}_x$ of condensed rings satisfies $( \pi^{\#}_{i,x} )^{-1} (v) \in \cat{V}_{i, \pi_i(x)}$. Let us write $\cat{V} := (\cat{V}_x)_{x \in X}$. Then
\begin{equation}
\left( (X,\cat{O}, \cat{V}),
\left( (X,\cat{O}, \cat{V}) \xto{(\pi_i, \pi^{\#}_i)} (X_i , \cat{O}_i, \cat{V}_i) \right)_{i \in |\cat{I}|} \right)
\end{equation}
is the limit of the functor $D : \cat{I} \to \cat{C}$. 
\end{prop}

\begin{proof}
It immediately follows from the definition that $(\pi_i, \pi^{\#}_i)$ is a morphism $(X,\cat{O}, \cat{V}) \to (X_i , \cat{O}_i, \cat{V}_i)$ in $\cat{C}$ for every $i \in \ob{I}$. Therefore
\begin{equation}
\left( (X,\cat{O}, \cat{V}),
\left( (X,\cat{O}, \cat{V}) \xto{(\pi_i, \pi^{\#}_i)} (X_i , \cat{O}_i, \cat{V}_i) \right)_{i \in |\cat{I}|} \right)
\end{equation}
is a cone on the functor $D : \cat{I} \to \cat{C}$. On the other hand, suppose that
\begin{equation}
\left( (Y,\cat{O}_Y, \cat{V}_Y),
\left( (Y,\cat{O}_Y, \cat{V}_Y) \xto{(\rho_i, \rho^{\#}_i)} (X_i , \cat{O}_i, \cat{V}_i) \right)_{i \in |\cat{I}|} \right)
\end{equation}
is another cone on the functor $D : \cat{I} \to \cat{C}$. Then
\begin{equation}
\left( (Y,\cat{O}_Y),
\left( (Y,\cat{O}_Y) \xto{(\rho_i, \rho^{\#}_i)} (X_i , \cat{O}_i) \right)_{i \in |\cat{I}|} \right)
\end{equation}
is a cone on the functor $G \of D : \cat{I} \to \cat{D}$. Since
\begin{equation}
\left( (X,\cat{O}),
\left( (X,\cat{O}) \xto{(\pi_i, \pi^{\#}_i)} (X_i , \cat{O}_i) \right)_{i \in |\cat{I}|} \right)
\end{equation}
is the limit of the functor $G \of D : \cat{I} \to \cat{D}$, there exists a unique morphism $(\sigma, \sigma^{\#}) : (Y,\cat{O}_Y) \to (X,\cat{O})$ in $\cat{D}$ such that the following diagram is commutative in $\cat{D}$ for every $i \in \ob{I}$.
\[\begin{tikzcd}
	{(Y,\cat{O}_Y)} & {(X,\cat{O})} \\
	& {(X_i, \cat{O}_i)}
	\arrow["{(\sigma,\sigma^{\#})}", from=1-1, to=1-2]
	\arrow["{(\rho_i,\rho^{\#}_i)}"', from=1-1, to=2-2]
	\arrow["{(\pi_i, \pi^{\#}_i)}", from=1-2, to=2-2]
\end{tikzcd}\]
We show that $(\sigma, \sigma^{\#})$ is a morphism $(Y,\cat{O}_Y, \cat{V}_Y) \to (X,\cat{O}, \cat{V})$ in $\cat{C}$. Let $y \in Y$ and $w \in \cat{V}_Y$. For every $i \in \ob{I}$, the following diagram is commutative in $\ub{CRing}$.
\[\begin{tikzcd}
	{\cat{O}_{Y,y}} & {\cat{O}_{\sigma(y)}} \\
	& {\cat{O}_{i,\rho_i(y)}}
	\arrow["{\sigma^{\#}_{y}}"', from=1-2, to=1-1]
	\arrow["{\rho^{\#}_{i,y}}", from=2-2, to=1-1]
	\arrow["{\pi^{\#}_{i,\sigma(y)}}"', from=2-2, to=1-2]
\end{tikzcd}\]
Therefore
\begin{equation}
\left( \pi^{\#}_{i,\sigma(y)} \right)^{-1} 
\left( \left( \sigma^{\#}_y \right) ^{-1} (w) \right)
= \left( \rho^{\#}_{i,y} \right)^{-1} (w) .
\end{equation}
Since $(\rho_i, \rho^{\#}_i)$ is a morphism $(Y,\cat{O}_Y, \cat{V}_Y) \to (X_i , \cat{O}_i, \cat{V}_i)$ in $\cat{C}$, we have $( \rho^{\#}_{i,y} )^{-1} (w) \in \cat{V}_{i,\rho_i(y)}$. Consequently,
\begin{equation}
\left( \pi^{\#}_{i,\sigma(y)} \right)^{-1} 
\left( \left( \sigma^{\#}_y \right) ^{-1} (w) \right)
= \left( \rho^{\#}_{i,y} \right)^{-1} (w)
\, \in \, \cat{V}_{i,\rho_i(y)} = \cat{V}_{i, \pi_i(\sigma(y))} .
\end{equation}
Since this holds for every $i \in \ob{I}$, we conclude that
\begin{equation}
\left( \sigma^{\#}_y \right) ^{-1} (w) \in \cat{V}_{\sigma(y)} .
\end{equation}
by the definition of $\cat{V}_{\sigma(y)}$. This shows that $(\sigma, \sigma^{\#})$ is a morphism $(Y,\cat{O}_Y, \cat{V}_Y) \to (X,\cat{O}, \cat{V})$ in $\cat{C}$. This completes the proof.
\end{proof}

\begin{cor} \label{cor:completeness of the category C}
The category $\cat{C}$ is complete. The forgetful functor $\cat{C} \to \cat{D}$ preserves limits.
\end{cor}

\begin{proof}
This follows from \cref{prop:construction of limits in the category C} and \cref{cor:completeness of the category D}.
\end{proof}

\subsubsection{Cocompleteness}

\begin{prop} \label{prop:construction of colimits in the category C}
Let $\cat{I}$ be a small category. Let $D : \cat{I} \to \cat{C}$ be a functor. For each $i \in |\cat{I}|$, let us write $(X_i , \cat{O}_i, \cat{V}_i) := D(i)$. For each morphism $\alpha$ in $\cat{I}$, let us write $(f_{\alpha}, f^{\#}_{\alpha}) := D(\alpha)$. Let
\begin{equation}
\left( (X,\cat{O}),
\left( (X_i , \cat{O}_i) \xto{(\pi_i, \pi^{\#}_i)} (X,\cat{O}) \right)_{i \in |\cat{I}|} \right)
\end{equation}
be the colimit of the functor $G \of D : \cat{I} \to \cat{D}$. For each $x \in X$, let us define $\cat{V}_x$ to be the set of all continous valuations $v$ on $\cat{O}_x$ with the following property: There exist an $i \in |\cat{I}|$, an $x_i \in X_i$ and a $v_i \in \cat{V}_{i,x_i}$ such that $\pi_i(x_i) = x$ and the homomorphism $\pi^{\#}_{i,x_i} : \cat{O}_x \to \cat{O}_{i,x_i}$ of condensed rings satisfies $( \pi^{\#}_{i,x_i} )^{-1} (v_i) = v$. Let us write $\cat{V} := (\cat{V}_x)_{x \in X}$. Then
\begin{equation}
\left( (X,\cat{O}, \cat{V}),
\left( (X_i , \cat{O}_i, \cat{V}_i) \xto{(\pi_i, \pi^{\#}_i)} (X,\cat{O}, \cat{V}) \right)_{i \in |\cat{I}|} \right)
\end{equation}
is the colimit of the functor $D : \cat{I} \to \cat{C}$. 
\end{prop}

\begin{proof}
It immediately follows from the definition that $(\pi_i, \pi^{\#}_i)$ is a morphism $(X_i , \cat{O}_i, \cat{V}_i) \to (X,\cat{O}, \cat{V})$ in $\cat{C}$ for every $i \in \ob{I}$. Therefore
\begin{equation}
\left( (X,\cat{O}, \cat{V}),
\left( (X_i , \cat{O}_i, \cat{V}_i) \xto{(\pi_i, \pi^{\#}_i)} (X,\cat{O}, \cat{V}) \right)_{i \in |\cat{I}|} \right)
\end{equation}
is a cocone on the functor $D : \cat{I} \to \cat{C}$. On the other hand, suppose that
\begin{equation}
\left( (Y,\cat{O}_Y, \cat{V}_Y),
\left( (X_i , \cat{O}_i, \cat{V}_i) \xto{(\rho_i, \rho^{\#}_i)} (Y,\cat{O}_Y, \cat{V}_Y) \right)_{i \in |\cat{I}|} \right)
\end{equation}
is another cocone on the functor $D : \cat{I} \to \cat{C}$. Then
\begin{equation}
\left( (Y,\cat{O}_Y),
\left( (X_i , \cat{O}_i) \xto{(\rho_i, \rho^{\#}_i)} (Y,\cat{O}_Y) \right)_{i \in |\cat{I}|} \right)
\end{equation}
is a cocone on the functor $G \of D : \cat{I} \to \cat{D}$. Since
\begin{equation}
\left( (X,\cat{O}),
\left( (X_i , \cat{O}_i) \xto{(\pi_i, \pi^{\#}_i)} (X,\cat{O}) \right)_{i \in |\cat{I}|} \right)
\end{equation}
is the colimit of the functor $G \of D : \cat{I} \to \cat{D}$, there exists a unique morphism $(\sigma, \sigma^{\#}) : (X,\cat{O}) \to (Y,\cat{O}_Y)$ in $\cat{D}$ such that the following diagram is commutative in $\cat{D}$ for every $i \in \ob{I}$.
\[\begin{tikzcd}
	{(X,\cat{O})} & {(Y,\cat{O}_Y)} \\
	{(X_i, \cat{O}_i)}
	\arrow["{(\sigma,\sigma^{\#})}", from=1-1, to=1-2]
	\arrow["{(\pi_i, \pi^{\#}_i)}", from=2-1, to=1-1]
	\arrow["{(\rho_i,\rho^{\#}_i)}"', from=2-1, to=1-2]
\end{tikzcd}\]
We show that $(\sigma, \sigma^{\#})$ is a morphism $(X,\cat{O}, \cat{V}) \to (Y,\cat{O}_Y, \cat{V}_Y)$ in $\cat{C}$. Let $x \in X$ and $v \in \cat{V}_x$. By the definition of $\cat{V}_x$, there exist an $i \in |\cat{I}|$, an $x_i \in X_i$ and a $v_i \in \cat{V}_{i,x_i}$ such that $\pi_i(x_i) = x$ and the homomorphism $\pi^{\#}_{i,x_i} : \cat{O}_x \to \cat{O}_{i,x_i}$ of condensed rings satisfies $( \pi^{\#}_{i,x_i} )^{-1} (v_i) = v$. On the other hand, the following diagram is commutative in $\ub{CRing}$.
\[\begin{tikzcd}
	{\cat{O}_{x}} & {\cat{O}_{Y,\rho_i(x_i)}} \\
	{\cat{O}_{i,x_i}}
	\arrow["{\pi^{\#}_{i,x_i}}"', from=1-1, to=2-1]
	\arrow["{\sigma^{\#}_{x}}"', from=1-2, to=1-1]
	\arrow["{\rho^{\#}_{i,x_i}}", from=1-2, to=2-1]
\end{tikzcd}\]
Therefore
\begin{equation}
\left( \sigma^{\#}_x \right)^{-1} (v)
= \left( \sigma^{\#}_x \right)^{-1} \left( \left( \pi^{\#}_{i,x_i} \right)^{-1} (v_i) \right)
= \left( \rho^{\#}_{i,x_i} \right)^{-1} (v_i) .
\end{equation}
Since $(\rho_i, \rho^{\#}_i)$ is a morphism $(X_i , \cat{O}_i, \cat{V}_i) \to (Y,\cat{O}_Y, \cat{V}_Y)$ in $\cat{C}$, we have $( \rho^{\#}_{i,x_i} )^{-1} (v_i) \in \cat{V}_{Y,\rho_i(x_i)}$. Consequently,
\begin{equation}
\left( \sigma^{\#}_x \right)^{-1} (v)
= \left( \rho^{\#}_{i,x_i} \right)^{-1} (v_i) 
\, \in \, \cat{V}_{Y,\rho_i(x_i)} = \cat{V}_{Y,\sigma(x)} .
\end{equation}
This shows that $(\sigma, \sigma^{\#})$ is a morphism $(X,\cat{O}, \cat{V}) \to (Y,\cat{O}_Y, \cat{V}_Y)$ in $\cat{C}$. This completes the proof.
\end{proof}

\begin{cor} \label{cor:cocompleteness of the category C}
The category $\cat{C}$ is cocomplete. The forgetful functor $\cat{C} \to \cat{D}$ preserves colimits.
\end{cor}

\begin{proof}
This follows from \cref{prop:construction of colimits in the category C} and \cref{cor:cocompleteness of the category D}.
\end{proof}

\subsection{The category $\cat{C}_1$}

\subsubsection{Definition}

\begin{df}
The category $\cat{C}_1$ is defined to be the full subcategory of $\cat{C}$ consisting of all $(X, \cat{O}_X, \cat{V}_X) \in \ob{C}$ such that for every $x \in X$, the set $\cat{V}_{X,x}$ is a singleton.
\end{df}

\begin{nt}
Let $(X, \cat{O}_X, \cat{V}_X) \in |\cat{C}_1|$. If $x \in X$, then we write $|\cdot|_x$ for the continuous valuation on $\cat{O}_{X,x}$ which is the unique element of $\cat{V}_{X,x}$. Then the pair $(\cat{O}_{X,x}, |\cdot|_x)$ is a valued condensed ring in the sense of \cref{df:valued condensed ring}.
\end{nt}

\subsubsection{Coreflection}

\begin{prop} \label{prop:construction of coreflection in C 1}
Let $(X, \cat{O}, \cat{V})$ be an object of $\cat{C}$. Then there exists a coreflection
\begin{equation}
\left( (\tilde{X}, \tilde{\cat{O}}, \tilde{\cat{V}}) ,\,  
(\tilde{X}, \tilde{\cat{O}}, \tilde{\cat{V}}) \xto{(\pi, \pi^{\#})} (X, \cat{O}, \cat{V})
\right)
\end{equation}
of $(X, \cat{O}, \cat{V})$ along the inclusion functor $\cat{C}_1 \mon \cat{C}$ which has the following property.
\begin{enumerate}
\item
The underlying set of $\tilde{X}$ is the set of all pairs $(x,v)$ consisting of a point $x \in X$ and a valuation $v \in \cat{V}_x$.

\item
The map $\tilde{X} \xto{\pi} X$ is equal to the map $(x,v) \mapsto x$.

\item
The topology of $\tilde{X}$ is equal to the initial topology induced by the topology of $X$ under the map $\tilde{X} \xto{\pi} X$.

\item
For every $(x,v) \in \tilde{X}$, the homomorphism $\pi^{\#}_{(x,v)} : \cat{O}_x \to \tilde{\cat{O}}_{(x,v)}$ is an isomorphism of condensed rings. If $\tilde{v}$ denotes the unique element of $\tilde{\cat{V}}_{(x,v)}$, then $\left( \pi^{\#}_{(x,v)} \right)^{-1}(\tilde{v}) = v$.
\end{enumerate}
\end{prop}

\subsubsection{Proof of \cref{prop:construction of coreflection in C 1}}
\label{sec:Proof of construction of coreflection in C 1}

\begin{proof}
Let $(X, \cat{O}, \cat{V})$ be an object of $\cat{C}$. Let us define $\tilde{X} := \set{(x,v)}{x \in X ,\; v \in \cat{V}_x}$. Let $\pi : \tilde{X} \to X$ be the map $(x,v) \to x$. Let us endow $\tilde{X}$ with the initial topology induced by the topology of $X$ under the map $\pi : \tilde{X} \to X$. Thus a subset of $\tilde{X}$ is open if and only if it is of the form $\pi^{-1}(U)$ for some open subset $U$ of $X$. Then the map $\pi : \tilde{X} \to X$ is continuous. Let $\tilde{\cat{O}}$ be the inverse image $\pi^{-1} \cat{O}$ of $\cat{O}$ under $\pi$. Let us write $\pi^{\#} : \cat{O} \to \pi_* \tilde{\cat{O}}$ for the canonical morphism of sheaves of condensed rings on $X$. Then \cref{cor:existence of inverse image functor} shows that for each $(x,v) \in \tilde{X}$, the homomorphism $\pi^{\#}_{(x,v)} : \cat{O}_x \to \tilde{\cat{O}}_{(x,v)}$ is an isomorphism of condensed rings. Therefore there exists a unique valuation $\tilde{v}_{(x,v)}$ on $\tilde{\cat{O}}_{(x,v)}$ such that $\left( \pi^{\#}_{(x,v)} \right)^{-1} ( \tilde{v}_{(x,v)} ) = v$. $\tilde{v}_{(x,v)}$ is a continuous valuation on $\tilde{\cat{O}}_{(x,v)}$ since $v$ is a continuous valuation on $\cat{O}_x$. We define $\tilde{\cat{V}}_{(x,v)} := \{ \tilde{v}_{(x,v)} \}$. Then $(\tilde{X}, \tilde{\cat{O}}, \tilde{\cat{V}})$ is an object of $\cat{C}_1$, and $(\pi, \pi^{\#})$ is a morphism $(\tilde{X}, \tilde{\cat{O}}, \tilde{\cat{V}}) \to (X, \cat{O}, \cat{V})$ in $\cat{C}$.

Let us prove that 
\begin{equation}
\left( (\tilde{X}, \tilde{\cat{O}}, \tilde{\cat{V}}) ,\,  
(\tilde{X}, \tilde{\cat{O}}, \tilde{\cat{V}}) \xto{(\pi, \pi^{\#})} (X, \cat{O}, \cat{V})
\right)
\end{equation}
is a coreflection of $(X, \cat{O}, \cat{V})$ along the inclusion functor $\cat{C}_1 \mon \cat{C}$. Suppose that $(Y, \cat{O}_Y, \cat{V}_Y)$ is an object of $\cat{C}_1$ and that $(\rho, \rho^{\#}) : (Y, \cat{O}_Y, \cat{V}_Y) \to (X, \cat{O}, \cat{V})$ is a morphism in $\cat{C}$. We prove that there exists a unique morphism $(\sigma, \sigma^{\#}) : (Y, \cat{O}_Y, \cat{V}_Y) \to (\tilde{X}, \tilde{\cat{O}}, \tilde{\cat{V}})$ in $\cat{C}$ such that the following diagram is commutative.
\[\begin{tikzcd}
	{(Y, \cat{O}_Y, \cat{V}_Y)} & {(\tilde{X}, \tilde{\cat{O}}, \tilde{\cat{V}})} \\
	& {(X, \cat{O}, \cat{V})}
	\arrow["{(\sigma,\sigma^{\#})}", from=1-1, to=1-2]
	\arrow["{(\rho,\rho^{\#})}"', from=1-1, to=2-2]
	\arrow["{(\pi, \pi^{\#})}", from=1-2, to=2-2]
\end{tikzcd}\]

Let $y \in Y$. Then we have a homomorphism $\rho^{\#}_y : \cat{O}_{\rho(y)} \to \cat{O}_{Y,y}$ of condensed rings. Let us write $w_y$ for the valuation on $\cat{O}_{Y,y}$ which is the unique element of $\cat{V}_{Y,y}$. Then $\left( \rho^{\#}_y \right) ^{-1} (w_y) \in \cat{V}_{\rho(y)}$ since $(\rho, \rho^{\#}) : (Y, \cat{O}_Y, \cat{V}_Y) \to (X, \cat{O}, \cat{V})$ is a morphism in $\cat{C}$. Consequently, we have an element
\begin{equation}
\sigma(y) := \left( \rho(y) ,\, \left( \rho^{\#}_y \right) ^{-1} (w_y) \right) \in \tilde{X} .
\end{equation}
Thus we obtain a map $\sigma :Y \to \tilde{X}$, $y \mapsto \sigma(y)$. The following diagram is commutative.
\[\begin{tikzcd}
	Y & {\tilde{X}} \\
	& X
	\arrow["\sigma", from=1-1, to=1-2]
	\arrow["\rho"', from=1-1, to=2-2]
	\arrow["\pi", from=1-2, to=2-2]
\end{tikzcd}\]
Since $\tilde{X}$ has the initial topology induced by the topology of $X$ under the map $\pi : \tilde{X} \to X$ and since $\pi \of \sigma = \rho : Y \to X$ is a continuous map, we conclude that $\sigma : Y \to \tilde{X}$ is a continuous map. Then the equality
\begin{equation}
\rho_* \cat{O}_Y = \pi_* \sigma_* \cat{O}_Y
\end{equation}
holds. Therefore the morphism $\rho^{\#} : \cat{O} \to \rho_* \cat{O}_Y$ of sheaves of condensed rings on $X$ is a morphism $\cat{O} \to \pi_* \sigma_* \cat{O}_Y$. Since $\tilde{\cat{O}}$ is the inverse image of $\cat{O}$ under $\pi$, there exists a unique morphism $\sigma^{\#} : \tilde{\cat{O}} \to \sigma_* \cat{O}_Y$ of sheaves of condensed rings on $\tilde{X}$ such that the following diagram is commutative.
\[\begin{tikzcd}
	{\pi_* \sigma_* \cat{O}_Y} & {\pi_* \tilde{\cat{O}}} \\
	& {\cat{O}}
	\arrow["{\pi_* \sigma^{\#}}"', from=1-2, to=1-1]
	\arrow["{\rho^{\#}}", from=2-2, to=1-1]
	\arrow["{\pi^{\#}}"', from=2-2, to=1-2]
\end{tikzcd}\]
Thus we have a morphism $(\sigma, \sigma^{\#}) : (Y, \cat{O}_Y) \to (\tilde{X}, \tilde{\cat{O}})$ in $\cat{D}$ such that the following diagram is commutative.
\[\begin{tikzcd}
	{(Y, \cat{O}_Y)} & {(\tilde{X}, \tilde{\cat{O}})} \\
	& {(X, \cat{O})}
	\arrow["{(\sigma,\sigma^{\#})}", from=1-1, to=1-2]
	\arrow["{(\rho,\rho^{\#})}"', from=1-1, to=2-2]
	\arrow["{(\pi, \pi^{\#})}", from=1-2, to=2-2]
\end{tikzcd}\]

Let us prove that $(\sigma, \sigma^{\#})$ is a morphism $(Y, \cat{O}_Y, \cat{V}_Y) \to (\tilde{X}, \tilde{\cat{O}}, \tilde{\cat{V}})$ in $\cat{C}$. For every $y \in Y$, the following diagram is commutative.
\[\begin{tikzcd}
	{\cat{O}_{Y,y}} & {\tilde{\cat{O}}_{\sigma(y)}} \\
	& {\cat{O}_{\rho(y)}}
	\arrow["{\sigma^{\#}_{y}}"', from=1-2, to=1-1]
	\arrow["{\rho^{\#}_y}", from=2-2, to=1-1]
	\arrow["{\pi^{\#}_{\sigma(y)}}"', from=2-2, to=1-2]
\end{tikzcd}\]
Then we have
\begin{equation}
\left( \pi^{\#}_{\sigma(y)} \right)^{-1} \left( \left( \sigma^{\#}_y \right)^{-1} (w_y) \right)
= \left( \rho^{\#}_y \right) ^{-1} (w_y) .
\end{equation}
Since $\tilde{v}_{\sigma(y)}$ is the unique valuation on $\tilde{\cat{O}}_{\sigma(y)}$ such that $\left( \pi^{\#}_{\sigma(y)} \right)^{-1} ( \tilde{v}_{\sigma(y)} ) = \left( \rho^{\#}_y \right) ^{-1} (w_y)$, it follows that
\begin{equation}
\left( \sigma^{\#}_y \right)^{-1} (w_y) = \tilde{v}_{\sigma(y)} .
\end{equation}
Consequently $\left( \sigma^{\#}_y \right)^{-1} (w_y) = \tilde{v}_{\sigma(y)} \in \tilde{\cat{V}}_{\sigma(y)}$. This proves that $(\sigma, \sigma^{\#}) : (Y, \cat{O}_Y, \cat{V}_Y) \to (\tilde{X}, \tilde{\cat{O}}, \tilde{\cat{V}})$ is a morphism in $\cat{C}$.

Since the diagram
\[\begin{tikzcd}
	{(Y, \cat{O}_Y)} & {(\tilde{X}, \tilde{\cat{O}})} \\
	& {(X, \cat{O})}
	\arrow["{(\sigma,\sigma^{\#})}", from=1-1, to=1-2]
	\arrow["{(\rho,\rho^{\#})}"', from=1-1, to=2-2]
	\arrow["{(\pi, \pi^{\#})}", from=1-2, to=2-2]
\end{tikzcd}\]
is commutative in $\cat{D}$, the diagram
\[\begin{tikzcd}
	{(Y, \cat{O}_Y, \cat{V}_Y)} & {(\tilde{X}, \tilde{\cat{O}}, \tilde{\cat{V}})} \\
	& {(X, \cat{O}, \cat{V})}
	\arrow["{(\sigma,\sigma^{\#})}", from=1-1, to=1-2]
	\arrow["{(\rho,\rho^{\#})}"', from=1-1, to=2-2]
	\arrow["{(\pi, \pi^{\#})}", from=1-2, to=2-2]
\end{tikzcd}\]
is commutative in $\cat{C}$.

On the other hand, suppose that $(\tau, \tau^{\#}) : (Y, \cat{O}_Y, \cat{V}_Y) \to (\tilde{X}, \tilde{\cat{O}}, \tilde{\cat{V}})$ is another morphism in $\cat{C}$ such that the following diagram is commutative.
\[\begin{tikzcd}
	{(Y, \cat{O}_Y, \cat{V}_Y)} & {(\tilde{X}, \tilde{\cat{O}}, \tilde{\cat{V}})} \\
	& {(X, \cat{O}, \cat{V})}
	\arrow["{(\tau,\tau^{\#})}", from=1-1, to=1-2]
	\arrow["{(\rho,\rho^{\#})}"', from=1-1, to=2-2]
	\arrow["{(\pi, \pi^{\#})}", from=1-2, to=2-2]
\end{tikzcd}\]
We prove that $(\tau, \tau^{\#}) = (\sigma, \sigma^{\#})$. Let $y \in Y$. Write $\tau(y) = (x,v) \in \tilde{X}$. Since the diagram
\[\begin{tikzcd}
	Y & {\tilde{X}} \\
	& X
	\arrow["\tau", from=1-1, to=1-2]
	\arrow["\rho"', from=1-1, to=2-2]
	\arrow["\pi", from=1-2, to=2-2]
\end{tikzcd}\]
is commutative, we have $x = \pi( \tau(y) ) = \rho(y)$. Furthermore, since $(\tau, \tau^{\#}) : (Y, \cat{O}_Y, \cat{V}_Y) \to (\tilde{X}, \tilde{\cat{O}}, \tilde{\cat{V}})$ is a morphism in $\cat{C}$, we have
\begin{equation}
\left( \tau^{\#}_y \right)^{-1}(w_y) \in \tilde{\cat{V}}_{(x,v)} .
\end{equation}
Since $\tilde{\cat{V}}_{(x,v)} = \{ \tilde{v}_{(x,v)} \}$, it follows that
\begin{equation}
\left( \tau^{\#}_y \right)^{-1}(w_y) = \tilde{v}_{(x,v)} .
\end{equation}
Therefore
\begin{equation}
\left( \pi^{\#}_{(x,v)} \right)^{-1} \left( \left( \tau^{\#}_y \right)^{-1}(w_y) \right)
= \left( \pi^{\#}_{(x,v)} \right)^{-1} ( \tilde{v}_{(x,v)} )
= v.
\end{equation}
On the other hand, the following diagram is commutative.
\[\begin{tikzcd}
	{\cat{O}_{Y,y}} & {\tilde{\cat{O}}_{(x,v)}} \\
	& {\cat{O}_{\rho(y)}}
	\arrow["{\tau^{\#}_{y}}"', from=1-2, to=1-1]
	\arrow["{\rho^{\#}_y}", from=2-2, to=1-1]
	\arrow["{\pi^{\#}_{(x,v)}}"', from=2-2, to=1-2]
\end{tikzcd}\]
Therefore
\begin{equation}
\left( \pi^{\#}_{(x,v)} \right)^{-1} \left( \left( \tau^{\#}_y \right)^{-1}(w_y) \right)
= \left( \rho^{\#}_y \right)^{-1} ( w_y ).
\end{equation}
Consequently, we have
\begin{equation}
v = 
\left( \pi^{\#}_{(x,v)} \right)^{-1} \left( \left( \tau^{\#}_y \right)^{-1}(w_y) \right)
= \left( \rho^{\#}_y \right)^{-1} ( w_y ).
\end{equation}
Thus
\begin{equation}
\tau(y) = (x,v) = \left( \rho(y) ,\, \left( \rho^{\#}_y \right)^{-1} ( w_y ) \right) = \sigma(y).
\end{equation}
This shows that $\tau=\sigma : Y \to \tilde{X}$. Then the morphism $\tau^{\#} : \tilde{\cat{O}} \to \tau_* \cat{O}_Y$ of sheaves of condensed rings on $\tilde{X}$ is a morphism $\tilde{\cat{O}} \to \sigma_* \cat{O}_Y$. Moreover, the following diagram is commutative.
\[\begin{tikzcd}
	{\pi_* \sigma_* \cat{O}_Y} & {\pi_* \tilde{\cat{O}}} \\
	& {\cat{O}}
	\arrow["{\pi_* \tau^{\#}}"', from=1-2, to=1-1]
	\arrow["{\rho^{\#}}", from=2-2, to=1-1]
	\arrow["{\pi^{\#}}"', from=2-2, to=1-2]
\end{tikzcd}\]
Then the uniqueness of $\sigma^{\#} : \tilde{\cat{O}} \to \sigma_* \cat{O}_Y$ shows that $\tau^{\#} = \sigma^{\#}$. Thus $(\tau, \tau^{\#}) = (\sigma, \sigma^{\#})$. This completes the proof of \cref{prop:construction of coreflection in C 1}.
\end{proof}

\subsection{The category $\cat{C}_l$}

\subsubsection{Definition}

\begin{df}
The category $\cat{C}_l$ is defined to be the full subcategory of $\cat{C}_1$ consisting of all $(X, \cat{O}_X, \cat{V}_X) \in |\cat{C}_1|$ with the following property: For every $x \in X$, the ring $\cat{O}_{X,x}(*)$ is a local ring, whose maximal ideal is equal to the support of the unique element of $\cat{V}_{X,x}$. In other words, the valued condensed ring $(\cat{O}_{X,x}, |\cdot|_x)$ is an object of $\ub{VCRing}_l$ (\cref{df:local valued ring}) for every $x \in X$.
\end{df}

\subsubsection{Coreflection}

\begin{prop} \label{prop:construction of coreflection in C l}
Let $(X, \cat{O}, \cat{V})$ be an object of $\cat{C}_1$. Then there exists a coreflection
\begin{equation}
\left( (\tilde{X}, \tilde{\cat{O}}, \tilde{\cat{V}}) ,\,  
(\tilde{X}, \tilde{\cat{O}}, \tilde{\cat{V}}) \xto{(\pi, \pi^{\#})} (X, \cat{O}, \cat{V})
\right)
\end{equation}
of $(X, \cat{O}, \cat{V})$ along the inclusion functor $\cat{C}_l \mon \cat{C}_1$ which has the following properties.
\begin{enumerate}
\item
The underlying set of $\tilde{X}$ is equal to the underlying set of $X$.

\item
The map $\tilde{X} \xto{\pi} X$ is equal to the identity map $x \mapsto x$.

\item
Let $I$ be the set of all pairs $(U;g)$ consisting of an open subset $U$ of $X$ and an element $g \in \cat{O}(U)(*)$. For $(U;g) \in I$, write $D(U;g)$ for the set of all $x \in U$ such that $|g_x|_x \neq 0$. Then the set
\begin{equation}
\set{D(U;g)}{(U;g) \in I}
\end{equation}
is a basis for the topology of $\tilde{X}$.

\item
For every $x \in \tilde{X}$, the homomorphism $\pi^{\#}_x : (\cat{O}_x, |\cdot|_x) \to (\tilde{\cat{O}}_x, |\cdot|_x)$ of valued condensed rings is the localization of the valued condensed ring $(\cat{O}_x, |\cdot|_x)$ in the sense of \cref{df:localization of valued condensed ring}. In particular, the homomorphism $\pi^{\#}_x : \cat{O}_x \to \tilde{\cat{O}}_x$ of condensed rings induces an isomorphism of condensed rings
\begin{equation}
(\cat{O}_x)_{\mathrm{Supp}(|\cdot|_x)} \simeq \tilde{\cat{O}}_x .
\end{equation}
\end{enumerate}
\end{prop}

\subsubsection{Proof of \cref{prop:construction of coreflection in C l}}
\label{sec:Proof of construction of coreflection in C l}

\begin{proof}
Let $(X, \cat{O}, \cat{V})$ be an object of $\cat{C}_1$. For each $x \in X$, let us write $v_x$ for the valuation on $\cat{O}_x$ which is the unique element of $\cat{V}_x$. Let $I$ be the set of all pairs $(U;g)$ consisting of an open subset $U$ of $X$ and an element $g \in \cat{O}(U)(*)$. For $(U;g) \in I$, write 
\begin{equation}
D(U;g) := \set{x \in U}{|g_x|_{v_x} \neq 0} .
\end{equation}

\begin{lem} \label{lem:properties of D(U g)} \;
\begin{enumerate}
\item
Suppose that $(U;g) ,\, (U';g') \in I$. Let us write
\begin{align}
U'' & := U \cap U' \: ; \\
g'' & := (g|_{U''}) \cdot (g'|_{U''}) \in \cat{O}(U'')(*) .
\end{align}
Then we have $(U'';g'') \in I$ and
\begin{equation}
D(U;g) \cap D(U';g') = D(U'';g'') .
\end{equation}

\item
For every open subset $U$ of $X$, we have
\begin{equation}
U = D(U;1) .
\end{equation}
In particular, we have $X = D(X;1)$.

\end{enumerate}
\end{lem}

\begin{proof}~
\begin{enumerate}
\item
$U''$ is an open subset of $X$ which is contained in both $U$ and $U'$. Therefore $g|_{U''}, g'|_{U''}$ are well-defined elements of $\cat{O}(U'')(*)$. Thus $g'' = (g|_{U''}) \cdot (g'|_{U''})$ is an element $\cat{O}(U'')(*)$ and we have $(U'';g'') \in I$. Then
\begin{align}
D(U;g) \cap D(U';g')
& = \set{x \in U}{|g_x|_{v_x} \neq 0} \cap \set{x \in U'}{|g'_x|_{v_x} \neq 0} \\
& = \set{x \in U \cap U'}{|g_x \cdot g_x'|_{v_x} \neq 0} \\
& = \set{x \in U''}{|g''_x|_{v_x} \neq 0} = D(U'';g'') .
\end{align}

\item
Let $U$ be an open subset of $X$. For every $x \in U$, we have $|1_x|_{v_x} = 1 \neq 0$. Therefore $U \sub D(U;1)$. Since $D(U;1) \sub U$ by defintion, we have $U = D(U;1)$.
\end{enumerate}
\end{proof}

By \cref{lem:properties of D(U g)}, there exists a unique topology on the underlying set of $X$ for which the set
\begin{equation}
\set{D(U;g)}{(U;g) \in I}
\end{equation}
is an open basis. Let $\tilde{X}$ be the set $X$ endowed with this topology. Let $\pi : \tilde{X} \to X$ be the identity map $x \mapsto x$. Then (2) of \cref{lem:properties of D(U g)} shows that $\pi : \tilde{X} \to X$ is continuous.

We use \cref{prop:constructing condensed sheaves} in order to construct a sheaf $\tilde{\cat{O}}$ of condensed rings on $\tilde{X}$. Let $\tilde{\cat{T}}$ be the set of all open sets of $\tilde{X}$ ordered by inclusion, which we consider as a small category. We define a preorder $\leq$ on the set $I$ by declaring that for $(U;g) ,\, (U';g') \in I$, the relation $(U,g) \leq (U';g')$ holds if and only if all of the following conditions are satisfied. 
\begin{enumerate}
\item $U \bus U'$.

\item
The image of $g|_{U'} \in \cat{O}(U')(*)$ in $\cat{O}(U')(*)_{g'}$ is invertible in $\cat{O}(U')(*)_{g'}$.

\item $D(U;g) \bus D(U';g')$.
\end{enumerate}
We consider the preordered set $I$ as a small category. Then (1) of  \cref{lem:properties of D(U g)} shows that $(I, \leq)$ is a directed set and that the functor
\[\begin{tikzcd}
	{D:I} & {\tilde{\cat{T}}^{\op},} & {(U;g)} & {D(U;g)}
	\arrow[from=1-1, to=1-2]
	\arrow[maps to, from=1-3, to=1-4]
\end{tikzcd}\]
satisfies (2) of \cref{as:assumption for constructing condensed sheaves}.
Moreover, the definition of the topology of $\tilde{X}$ shows that (1) of \cref{as:assumption for constructing condensed sheaves} is satisfied.

Next we define a functor $R : I \to \ub{CRing}$ as follows.
\begin{enumerate}
\item
For $i = (U;g) \in I$, we define $R(i)$ to be the localization $\cat{O}(U)_g$ of the condensed $\cat{O}(U)$-algebra $\cat{O}(U)$ by $g \in \cat{O}(U)(*)$.

\item
Suppose $i = (U;g) \in I$ and $j = (U';g') \in I$ satisfy $i \leq j$. Then $U \bus U'$ and the image of $g|_{U'} \in \cat{O}(U')(*)$ in $\cat{O}(U')(*)_{g'} = \cat{O}(U')_{g'}(*)$ is invertible in $\cat{O}(U')_{g'}(*)$. By \cref{prop:universality of localization of condensed algebra}, there exists a unique homomorphism $R(i \leq j) : \cat{O}(U)_g \to \cat{O}(U')_{g'}$ of condensed rings such that the following diagram is commutative.
\[\begin{tikzcd}
	{\cat{O}(U)} & {\cat{O}(U)_g} \\
	{\cat{O}(U')} & {\cat{O}(U')_{g'}}
	\arrow["\can", from=1-1, to=1-2]
	\arrow["{\cat{O}(U' \sub U)}"', from=1-1, to=2-1]
	\arrow["{R(i \leq j)}", from=1-2, to=2-2]
	\arrow["\can"', from=2-1, to=2-2]
\end{tikzcd}\]
This is the definition of $R(i \leq j) : R(i) \to R(j)$.
\end{enumerate}

By \cref{prop:constructing condensed sheaves}, there exist a sheaf $\tilde{\cat{O}}$ of condensed rings on $\tilde{X}$ and a natural transformation $\beta : R \To \tilde{\cat{O}} \of D$ of functors $I \to \ub{CRing}$ with the following properties.
\begin{enumerate}
\item
For each $x \in \tilde{X}$, let us write $I(x)$ for the full subcategory of $I$ consisting of all $i \in I$ such that $x \in D(i)$. Then the homomorphisms of condensed rings
\[\begin{tikzcd}
	{R(i)} & {\tilde{\cat{O}}(D(i))} & {(i \in I(x))}
	\arrow["{\beta_i}", from=1-1, to=1-2]
\end{tikzcd}\]
induce an isomorphism of condensed rings
\[\begin{tikzcd}
	{\underset{i \in I(x)}{\colim} \, R(i)} & {\tilde{\cat{O}}_x .}
	\arrow["{\beta_x}", from=1-1, to=1-2]
\end{tikzcd}\]

\item
For every sheaf $H$ of condensed rings on $\tilde{X}$ and a natural transformation $\gamma : R \To H \of D$ of functors $I \to \ub{CRing}$, there exists a unique morphism $\delta : \tilde{\cat{O}} \to H$ of sheaves of condensed rings on $\tilde{X}$ such that the diagram
\[\begin{tikzcd}
	R & {H \of D} \\
	{\tilde{\cat{O}} \of D}
	\arrow["\gamma", from=1-1, to=1-2]
	\arrow["\beta"', from=1-1, to=2-1]
	\arrow["{\delta * \id_D}"', from=2-1, to=1-2]
\end{tikzcd}\]
is commutative.
\end{enumerate}

Let $\cat{T}$ be the set of all open subsets of $X$ ordered by inclusion, which we consider as a small category. Consider the functor
\[\begin{tikzcd}
	{C : \cat{T}^{\op}} & {I,} & U & {(U;1).}
	\arrow[from=1-1, to=1-2]
	\arrow[maps to, from=1-3, to=1-4]
\end{tikzcd}\]
Then we have $\cat{O} = R \of C : \cat{T}^{\op} \to \ub{CRing}$ by definition of $R$. Moreover, (2) of \cref{lem:properties of D(U g)} shows that $D \of C : \cat{T}^{\op} \to \tilde{\cat{T}}^{\op}$ is equal to the functor $U \mapsto \pi^{-1}(U)$. Therefore the natural transformation $\beta * \id_C : R \of C \to \tilde{\cat{O}} \of D \of C$ of functors $\cat{T}^{\op} \to \ub{CRing}$ is a morphism $\cat{O} \to \pi_* \tilde{\cat{O}}$ of sheaves of condensed rings on $X$. We define $\pi^{\#} := \beta * \id_C : \cat{O} \to \pi_* \tilde{\cat{O}}$.

\begin{lem} \label{lem:stalk of tilde O is localization at support}
Let $x \in \tilde{X}$.
\begin{enumerate}
\item
For each $i = (U;g) \in I(x)$, there exists a unique homomorphism $\epsilon_{x,i} : R(i) \to (\cat{O}_{x})_{\mathrm{Supp}(v_x)}$ of condensed rings such that the following diagram is commutative.
\[\begin{tikzcd}
	{\cat{O}(U)} & {\cat{O}(U)_g} & {R(i)} \\
	{\cat{O}_x} && {(\cat{O}_{x})_{\mathrm{Supp}(v_x)}}
	\arrow["\can", from=1-1, to=1-2]
	\arrow["\can"', from=1-1, to=2-1]
	\arrow[equals, from=1-2, to=1-3]
	\arrow["{\epsilon_{x,i}}", from=1-3, to=2-3]
	\arrow["\can"', from=2-1, to=2-3]
\end{tikzcd}\]

\item
The homomorphisms of condensed rings
\[\begin{tikzcd}
	{R(i)} & {(\cat{O}_{x})_{\mathrm{Supp}(v_x)}} & {(i \in I(x))}
	\arrow["{\epsilon_{x,i}}", from=1-1, to=1-2]
\end{tikzcd}\]
induce an isomorphism of condensed rings
\[\begin{tikzcd}
	{\underset{i \in I(x)}{\colim} \, R(i)} & {(\cat{O}_{x})_{\mathrm{Supp}(v_x)} .}
	\arrow["{\epsilon_x}", from=1-1, to=1-2]
\end{tikzcd}\]

\item
The following diagram is commutative.
\[\begin{tikzcd}
	{\cat{O}_x} & {\tilde{\cat{O}}_x} \\
	{(\cat{O}_{x})_{\mathrm{Supp}(v_x)}} & {\underset{i \in I(x)}{\colim} \, R(i)}
	\arrow["{\pi^{\#}_x}", from=1-1, to=1-2]
	\arrow["\can"', from=1-1, to=2-1]
	\arrow["{\epsilon_x^{-1}}"', from=2-1, to=2-2]
	\arrow["{\beta_x}"', from=2-2, to=1-2]
\end{tikzcd}\]
In other words, the homomorphism $\pi^{\#}_x : \cat{O}_x \to \tilde{\cat{O}}_x$ of condensed rings induces an isomorphism of condensed rings
\[\begin{tikzcd}
	{\beta_x \of \epsilon_x^{-1} : (\cat{O}_{x})_{\mathrm{Supp}(v_x)}} & {\tilde{\cat{O}}_x .}
	\arrow["\sim", from=1-1, to=1-2]
\end{tikzcd}\]

\end{enumerate}
\end{lem}

\begin{proof}~
\begin{enumerate}
\item
Let $i = (U;g) \in I(x)$. Then $x \in D(U;g)$ and therefore $|g_x|_{v_x} \neq 0$. In other words, we have $g_x \in \cat{O}_x(*) \sm \mathrm{Supp}(v_x)$. Therefore the image of $g_x \in \cat{O}_x(*)$ in $(\cat{O}_x(*))_{\mathrm{Supp}(v_x)} = (\cat{O}_x)_{\mathrm{Supp}(v_x)}(*)$ is invertible in $(\cat{O}_x)_{\mathrm{Supp}(v_x)}(*)$. Then the assertion $(1)$ follows from \cref{prop:universality of localization of condensed algebra}.

\item
We prove that the pair
\begin{equation}
\left( (\cat{O}_x)_{\mathrm{Supp}(v_x)} \; , \;
\left( R(i) \xto{\epsilon_{x,i}} (\cat{O}_{x})_{\mathrm{Supp}(v_x)} )\right)_{i \in |I(x)|}
\right)
\end{equation}
is a colimit of the functor $I(x) \mon I \xto{R} \ub{CRing}$.

First suppose that $i = (U;g) \in I(x)$ and $j = (U';g') \in I(x)$ satisfy $i \leq j$. The the following diagram is commutative.
\[\begin{tikzcd}
	{\cat{O}(U)} & {\cat{O}(U)_g} & {R(i)} \\
	{\cat{O}(U')} & {\cat{O}(U')_{g'}} & {R(j)} \\
	{\cat{O}_x} && {(\cat{O}_x)_{\mathrm{Supp}(v_x)}}
	\arrow["\can", from=1-1, to=1-2]
	\arrow["{\cat{O}(U' \sub U)}", from=1-1, to=2-1]
	\arrow["\can"', curve={height=24pt}, from=1-1, to=3-1]
	\arrow[equals, from=1-2, to=1-3]
	\arrow["{R(i \leq j)}", from=1-3, to=2-3]
	\arrow["\can", from=2-1, to=2-2]
	\arrow["\can", from=2-1, to=3-1]
	\arrow[equals, from=2-2, to=2-3]
	\arrow["{\epsilon_{x,j}}", from=2-3, to=3-3]
	\arrow["\can", from=3-1, to=3-3]
\end{tikzcd}\]
Then the definition of $\epsilon_{x,i}$ shows that $\epsilon_{x,j} \of R(i \leq j) = \epsilon_{x,i}$. This shows that the pair
\begin{equation}
\left( (\cat{O}_x)_{\mathrm{Supp}(v_x)} \; , \;
\left( R(i) \xto{\epsilon_{x,i}} (\cat{O}_{x})_{\mathrm{Supp}(v_x)} )\right)_{i \in |I(x)|}
\right)
\end{equation}
is a cocone on the functor $I(x) \mon I \xto{R} \ub{CRing}$.

On the other hand, suppose that
\begin{equation}
\left( A \; , \;
\left( R(i) \xto{\zeta_i} A \right)_{i \in |I(x)|}
\right)
\end{equation}
is another cocone on the functor $I(x) \mon I \xto{R} \ub{CRing}$. We show that there exists a unique homomorphism $\theta : (\cat{O}_x)_{\mathrm{Supp}(v_x)} \to A$ of condensed rings such that the diagram
\[\begin{tikzcd}
	{(\cat{O}_x)_{\mathrm{Supp}(v_x)}} & A \\
	{R(i)}
	\arrow["\theta", from=1-1, to=1-2]
	\arrow["{\epsilon_{x,i}}", from=2-1, to=1-1]
	\arrow["{\zeta_i}"', from=2-1, to=1-2]
\end{tikzcd}\]
is commutative for every $i \in I(x)$.

Let $\cat{T}(x)$ be the full subcategory of $\cat{T}$ consisting of all open neighbourhoods of $x$ in $X$. If $U,V$ are open neighbourhoods of $x$ in $X$ such that $V \sub U$, then the following diagram is commutative.
\[\begin{tikzcd}
	{\cat{O}(U)} & {R(C(U))} && A \\
	{\cat{O}(V)} & {R(C(V))}
	\arrow[equals, from=1-1, to=1-2]
	\arrow["{\cat{O}(V \sub U)}"', from=1-1, to=2-1]
	\arrow["{\zeta_{C(U)}}", from=1-2, to=1-4]
	\arrow["{R(C(U) \leq C(V))}", from=1-2, to=2-2]
	\arrow[equals, from=2-1, to=2-2]
	\arrow["{\zeta_{C(V)}}"', curve={height=12pt}, from=2-2, to=1-4]
\end{tikzcd}\]
Therefore the pair
\begin{equation}
\left( A \; , \;
\left( \cat{O}(U) \xto{\zeta_{C(U)}} A \right)_{U \in |\cat{T}(x)^{\op}|} \right)
\end{equation}
is a cocone on the functor $\cat{T}(x)^{\op} \mon \cat{T}^{\op} \xto{\cat{O}} \ub{CRing}$. Consequently, there exists a unique homomorphism $\tilde{\zeta} : \cat{O}_x \to A$ of condensed rings such that the diagram
\[\begin{tikzcd}
	{\cat{O}(U)} & {R(C(U))} \\
	{\cat{O}_x} & A
	\arrow[equals, from=1-1, to=1-2]
	\arrow["\can"', from=1-1, to=2-1]
	\arrow["{\zeta_{C(U)}}", from=1-2, to=2-2]
	\arrow["{\tilde{\zeta}}"', from=2-1, to=2-2]
\end{tikzcd}\]
is commutative for every open neighbourhood $U$ of $x$ in $X$. Next suppose that $t \in \cat{O}_x(*) \sm \mathrm{Supp}(v_x)$. By \cref{prop:compatibility of stalk and evaluation at S}, we have
\begin{equation}
\cat{O}_x (*) = \underset{U \in \cat{T}(x)^{\op}}{\colim} \; \cat{O}(U)(*) ,
\end{equation}
where the colimit is taken in $\ub{Ring}$. Then the construction of filtered colimits in $\ub{Ring}$ shows that there exist an open neighbourhood $U$ of $x$ in $X$ and an element $g \in \cat{O}(U)(*)$ such that $g_x = t$. Then $i := (U;g) \in I$. Furthermore, we have $|g_x|_{v_x} = |t|_{v_x} \neq 0$ since $t \in \cat{O}_x(*) \sm \mathrm{Supp}(v_x)$. Therefore $x \in D(U;g)$ and hence $i \in I(x)$. Then the diagram
\[\begin{tikzcd}
	{\cat{O}_x} \\
	{\cat{O}(U)} && {R(C(U))} & A \\
	{\cat{O}(U)} & {\cat{O}(U)_g} & {R(i)}
	\arrow["{\tilde{\zeta}}", curve={height=-12pt}, from=1-1, to=2-4]
	\arrow["\can", from=2-1, to=1-1]
	\arrow[equals, from=2-1, to=2-3]
	\arrow["{\zeta_{C(U)}}"', from=2-3, to=2-4]
	\arrow["{R(C(U) \leq i)}"', from=2-3, to=3-3]
	\arrow[equals, from=3-1, to=2-1]
	\arrow["\can"', from=3-1, to=3-2]
	\arrow[equals, from=3-2, to=3-3]
	\arrow["{\zeta_i}"', curve={height=12pt}, from=3-3, to=2-4]
\end{tikzcd}\]
is commutative. Since the image $g/1$ of $g$ in $\cat{O}(U)_g(*) = \cat{O}(U)(*)_g$ is invertible in $\cat{O}(U)_g(*)$, it follows that the element $(\tilde{\zeta})_*(t) = (\tilde{\zeta})_*(g_x) = (\zeta_i)_*(g/1) \in A(*)$ is invertible in $A(*)$. This holds for every $t \in \cat{O}_x(*) \sm \mathrm{Supp}(v_x)$. Then \cref{prop:universality of localization of condensed algebra} shows that there exists a unique homomorphism $\theta : (\cat{O}_x)_{\mathrm{Supp}(v_x)} \to A$ of condensed rings such that the diagram
\[\begin{tikzcd}
	{(\cat{O}_x)_{\mathrm{Supp}(v_x)}} & A \\
	{\cat{O}_x}
	\arrow["\theta", from=1-1, to=1-2]
	\arrow["\can", from=2-1, to=1-1]
	\arrow["{\tilde{\zeta}}"', from=2-1, to=1-2]
\end{tikzcd}\]
is commutative. Then, for each $i = (U;g) \in I$, the diagram
\[\begin{tikzcd}
	{R(i)} & {(\cat{O}_x)_{\mathrm{Supp}(v_x)}} \\
	{\cat{O}(U)_g} & {\cat{O}_x} && A \\
	{\cat{O}(U)} && {R(C(U))} \\
	{\cat{O}(U)} & {\cat{O}(U)_g} & {R(i)}
	\arrow["{\epsilon_{x,i}}", from=1-1, to=1-2]
	\arrow["\theta", from=1-2, to=2-4]
	\arrow[equals, from=2-1, to=1-1]
	\arrow["\can", from=2-2, to=1-2]
	\arrow["{\tilde{\zeta}}"', from=2-2, to=2-4]
	\arrow["\can", from=3-1, to=2-1]
	\arrow["\can", from=3-1, to=2-2]
	\arrow[equals, from=3-1, to=3-3]
	\arrow[equals, from=3-1, to=4-1]
	\arrow["{\zeta_{C(U)}}"{pos=0.4}, from=3-3, to=2-4]
	\arrow["{R(C(U) \leq i)}"', from=3-3, to=4-3]
	\arrow["\can"', from=4-1, to=4-2]
	\arrow[equals, from=4-2, to=4-3]
	\arrow["{\zeta_i}"', curve={height=18pt}, from=4-3, to=2-4]
\end{tikzcd}\]
is commutative. By \cref{prop:universality of localization of condensed algebra}, we conclude that the diagram
\[\begin{tikzcd}
	{(\cat{O}_x)_{\mathrm{Supp}(v_x)}} & A \\
	{R(i)}
	\arrow["\theta", from=1-1, to=1-2]
	\arrow["{\epsilon_{x,i}}", from=2-1, to=1-1]
	\arrow["{\zeta_i}"', from=2-1, to=1-2]
\end{tikzcd}\]
is commutative. On the other hand, suppose that $\theta' : (\cat{O}_x)_{\mathrm{Supp}(v_x)} \to A$ is another homomorphism of condensed rings such that the diagram
\[\begin{tikzcd}
	{(\cat{O}_x)_{\mathrm{Supp}(v_x)}} & A \\
	{R(i)}
	\arrow["{\theta'}", from=1-1, to=1-2]
	\arrow["{\epsilon_{x,i}}", from=2-1, to=1-1]
	\arrow["{\zeta_i}"', from=2-1, to=1-2]
\end{tikzcd}\]
is commutative for every $i \in I(x)$. Then, for each open neighbourhood $U$ of $x$ in $X$, the diagram
\[\begin{tikzcd}
	{\cat{O}(U)} & {R(C(U))} \\
	{\cat{O}_x} & {(\cat{O}_x)_{\mathrm{Supp}(v_x)}} & A
	\arrow[equals, from=1-1, to=1-2]
	\arrow["\can"', from=1-1, to=2-1]
	\arrow["{\epsilon_{x,C(U)}}"', from=1-2, to=2-2]
	\arrow["{\zeta_{C(U)}}", from=1-2, to=2-3]
	\arrow["\can"', from=2-1, to=2-2]
	\arrow["{\theta'}"', from=2-2, to=2-3]
\end{tikzcd}\]
is commutative. Then the definition of $\tilde{\zeta}$ shows that the composition $\cat{O}_x \xto{\can} (\cat{O}_x)_{\mathrm{Supp}(v_x)} \xto{\theta'} A$ is equal to $\tilde{\zeta} : \cat{O}_x \to A$. In other words, the diagram 
\[\begin{tikzcd}
	{(\cat{O}_x)_{\mathrm{Supp}(v_x)}} & A \\
	{\cat{O}_x}
	\arrow["{\theta'}", from=1-1, to=1-2]
	\arrow["\can"', from=2-1, to=1-1]
	\arrow["{\tilde{\zeta}}"', from=2-1, to=1-2]
\end{tikzcd}\]
is commutative. Then the definition of $\theta$ shows that $\theta' = \theta$. This completes the proof of the assertion (2).

\item
For each open neighbourhood $U$ of $x$ in $X$, the following diagram is commutative.
\[\begin{tikzcd}
	& {\cat{O}_x} &&& {\tilde{\cat{O}}_x} \\
	{\cat{O}_x} & {\cat{O}(U)} & {\tilde{\cat{O}}(\pi^{-1}(U))} & {\tilde{\cat{O}}(D(C(U)))} \\
	& {R(C(U))} && {R(C(U))} \\
	{(\cat{O}_x)_{\mathrm{Supp}(v_x)}} &&&& {\underset{i \in I(x)}{\colim} \, R(i)}
	\arrow["{\pi^{\#}_x}", from=1-2, to=1-5]
	\arrow["\can"', from=2-1, to=4-1]
	\arrow["\can", from=2-2, to=1-2]
	\arrow["\can"', from=2-2, to=2-1]
	\arrow["{\pi^{\#}_U}", from=2-2, to=2-3]
	\arrow[equals, from=2-2, to=3-2]
	\arrow[equals, from=2-3, to=2-4]
	\arrow["\can", from=2-4, to=1-5]
	\arrow[equals, from=3-2, to=3-4]
	\arrow["{\epsilon_{x,C(U)}}", from=3-2, to=4-1]
	\arrow["{\beta_{C(U)}}"', from=3-4, to=2-4]
	\arrow["\can"', from=3-4, to=4-5]
	\arrow["{\epsilon_x^{-1}}"', from=4-1, to=4-5]
	\arrow["{\beta_x}"', from=4-5, to=1-5]
\end{tikzcd}\]
Therefore the diagram
\[\begin{tikzcd}
	{\cat{O}_x} & {\tilde{\cat{O}}_x} \\
	{(\cat{O}_{x})_{\mathrm{Supp}(v_x)}} & {\underset{i \in I(x)}{\colim} \, R(i)}
	\arrow["{\pi^{\#}_x}", from=1-1, to=1-2]
	\arrow["\can"', from=1-1, to=2-1]
	\arrow["{\epsilon_x^{-1}}"', from=2-1, to=2-2]
	\arrow["{\beta_x}"', from=2-2, to=1-2]
\end{tikzcd}\]
is commutative.
\end{enumerate}
\end{proof}

Let $x \in \tilde{X}$. By \cref{prop:construction of localization of vcrings} and (3) of \cref{lem:stalk of tilde O is localization at support}, there exists a unique continuous valuation $\tilde{v}_x$ on the condensed ring $\tilde{\cat{O}}_x$ such that $\left( \pi^{\#}_x \right)^{-1} (\tilde{v}_x) = v_x$. Moreover, the pair
\begin{equation}
\left( (\tilde{\cat{O}}_x , \tilde{v}_x) , \,
(\cat{O}_x , v_x) \xto{\pi^{\#}_x} (\tilde{\cat{O}}_x , \tilde{v}_x) \right)
\end{equation}
is the localization of the valued condensed ring $(\cat{O}_x , v_x)$ in the sense of \cref{df:localization of valued condensed ring}. We define $\tilde{\cat{V}}_x := \{ \tilde{v}_x \}$. Then we obtain a family $\tilde{\cat{V}} := (\tilde{\cat{V}}_x)_{x \in \tilde{X}}$. The triple $(\tilde{X}, \tilde{\cat{O}}, \tilde{\cat{V}})$ is an object of $\cat{C}_l$ and $(\pi, \pi^{\#})$ is a morphism $(\tilde{X}, \tilde{\cat{O}}, \tilde{\cat{V}}) \to (X, \cat{O}, \cat{V})$ in $\cat{C}_1$.

Let us prove that 
\begin{equation}
\left( (\tilde{X}, \tilde{\cat{O}}, \tilde{\cat{V}}) ,\,  
(\tilde{X}, \tilde{\cat{O}}, \tilde{\cat{V}}) \xto{(\pi, \pi^{\#})} (X, \cat{O}, \cat{V})
\right)
\end{equation}
is a coreflection of $(X, \cat{O}, \cat{V})$ along the inclusion functor $\cat{C}_l \mon \cat{C}_1$. Suppose that $(Y, \cat{O}_Y, \cat{V}_Y)$ is an object of $\cat{C}_l$ and that $(\rho, \rho^{\#}) : (Y, \cat{O}_Y, \cat{V}_Y) \to (X, \cat{O}, \cat{V})$ is a morphism in $\cat{C}_1$. We prove that there exists a unique morphism $(\sigma, \sigma^{\#}) : (Y, \cat{O}_Y, \cat{V}_Y) \to (\tilde{X}, \tilde{\cat{O}}, \tilde{\cat{V}})$ in $\cat{C}_l$ such that the following diagram is commutative.
\[\begin{tikzcd}
	{(Y, \cat{O}_Y, \cat{V}_Y)} & {(\tilde{X}, \tilde{\cat{O}}, \tilde{\cat{V}})} \\
	& {(X, \cat{O}, \cat{V})}
	\arrow["{(\sigma,\sigma^{\#})}", from=1-1, to=1-2]
	\arrow["{(\rho,\rho^{\#})}"', from=1-1, to=2-2]
	\arrow["{(\pi, \pi^{\#})}", from=1-2, to=2-2]
\end{tikzcd}\]

For each $y \in Y$, let us write $w_y$ for the valuation on $\cat{O}_{Y,y}$ which is the unique element of $\cat{V}_{Y,y}$. Let $\sigma :Y \to \tilde{X}$ be the map $y \mapsto \rho(y)$. We claim that this map $\sigma :Y \to \tilde{X}$ is continuous. Since the set $\set{D(U;g)}{(U;g) \in I}$ is a basis for the topology of $\tilde{X}$, it suffices to show that for every $(U;g) \in I$, the set $\sigma^{-1} \big( D(U;g) \big)$ is an open subset of $Y$. Let $(U;g) \in I$. Then we have a homomorphism $\rho^{\#}_U : \cat{O}(U) \to \cat{O}_Y(\rho^{-1}(U))$ of condensed rings. Thus for each $h \in \cat{O}(U)(*)$, we have an element $( \rho^{\#}_U )_* (h) \in \cat{O}_Y ( \rho^{-1}(U) )(*)$.

\begin{lem} \label{lem:descriotion of sigma inverse D(U g)}
The following equality holds.
\begin{equation}
\sigma^{-1} \big( D(U;g) \big)
= \set{y \in \rho^{-1}(U)}{ \Big| \Big( ( \rho^{\#}_U )_* (g) \Big)_y \Big|_{w_y} \neq 0 } .
\end{equation}
\end{lem}

\begin{proof}
We have
\begin{equation}
\sigma^{-1} \big( D(U;g) \big)
= \set{y \in Y}{ \rho(y) \in U \, ; \, \left| g_{\rho(y)} \right|_{v_{\rho(y)}} \neq 0 } .
\end{equation}
For each $y \in Y$, we have a homomorphism $\rho^{\#}_y : \cat{O}_{\rho(y)} \to \cat{O}_{Y,y}$ of condensed rings. Since $(\rho, \rho^{\#}) : (Y, \cat{O}_Y, \cat{V}_Y) \to (X, \cat{O}, \cat{V})$ is a morphism in $\cat{C}_1$, we have
\begin{equation}
v_{\rho(y)} = \left( \rho^{\#}_y \right) ^{-1} (w_y) .
\end{equation}
Therefore
\begin{align}
& \set{y \in Y}{ \rho(y) \in U \, ; \,
\left| g_{\rho(y)} \right|_{v_{\rho(y)}} \neq 0 } \\
= & \set{y \in Y}{ \rho(y) \in U \, ; \,
\left| g_{\rho(y)} \right|_{\left( \rho^{\#}_y \right) ^{-1} (w_y)} \neq 0 } \\
= & \set{y \in \rho^{-1}(U)}{
\left| (\rho^{\#}_y)_* ( g_{\rho(y)} ) \right|_{w_y} \neq 0 } .
\end{align}
For each $y \in \rho^{-1}(U)$, the following diagram is commutative.
\[\begin{tikzcd}
	{\cat{O}(U)(*)} & {\cat{O}_Y(\rho^{-1}(U))(*)} \\
	{\cat{O}_{\rho(y)}(*)} & {\cat{O}_{Y,y}(*)}
	\arrow["{(\rho^{\#}_U)_*}", from=1-1, to=1-2]
	\arrow["\can"', from=1-1, to=2-1]
	\arrow["\can", from=1-2, to=2-2]
	\arrow["{(\rho^{\#}_y)_*}"', from=2-1, to=2-2]
\end{tikzcd}\]
It follows that 
\begin{align}
& \set{y \in \rho^{-1}(U)}{
\left| (\rho^{\#}_y)_* ( g_{\rho(y)} ) \right|_{w_y} \neq 0 } \\
= & \set{y \in \rho^{-1}(U)}{
\Big| \Big( ( \rho^{\#}_U )_* (g) \Big)_y \Big|_{w_y} \neq 0 } .
\end{align}
Thus
\begin{align}
\sigma^{-1} \big( D(U;g) \big)
& = \set{y \in Y}{
\rho(y) \in U \, ; \, \left| g_{\rho(y)} \right|_{v_{\rho(y)}} \neq 0 } \\
& = \set{y \in \rho^{-1}(U)}{
\left| (\rho^{\#}_y)_* ( g_{\rho(y)} ) \right|_{w_y} \neq 0 } \\
& = \set{y \in \rho^{-1}(U)}{
\Big| \Big( ( \rho^{\#}_U )_* (g) \Big)_y \Big|_{w_y} \neq 0 } .
\end{align}
\end{proof}

\begin{lem} \label{lem:nonvanishing locus is open in objects of C l}
Let $V$ be an open subset of $Y$. Let $h \in \cat{O}_Y(V)(*)$. Then the set
\begin{equation}
V' := \set{y \in V}{ |h_y|_{w_y} \neq 0}
\end{equation}
is an open subset of $Y$. Moreover, the element $h|_{V'} \in \cat{O}_Y(V')(*)$ is invertible in $\cat{O}_Y(V')(*)$.
\end{lem}

\begin{proof}
Since $(Y, \cat{O}_Y, \cat{V}_Y)$ is an object of $\cat{C}_l$, the ring $\cat{O}_{Y,y} (*)$ is a local ring and its maximal ideal is equal to the support of the valuation $w_y$ for every $y \in V$. Therefore
\begin{equation}
\set{y \in V}{ |h_y|_{w_y} \neq 0} =
\set{y \in V}{ h_y \in \cat{O}_{Y,y} (*) \text{ is invertible in } \cat{O}_{Y,y} (*)} .
\end{equation}
On the other hand, consider the presheaf $\left( \cat{O}_Y \right)^*$ of commutative unital rings on $Y$ (\cref{nt:passing from condesed sheaves to ordinary sheaves}). For $y \in V$, \cref{prop:compatibility of stalk and evaluation at S} shows that $\cat{O}_{Y,y} (*) = \left( \left( \cat{O}_Y \right)^* \right)_y$ and that the canonical map $\left( \cat{O}_Y \right)^*(V) \to \left( \left( \cat{O}_Y \right)^* \right)_y$ coincides with the canonical map $\cat{O}_Y(V)(*) \to \cat{O}_{Y,y} (*)$. Therefore
\begin{align}
& \set{y \in V}{ h_y \in \cat{O}_{Y,y} (*) \text{ is invertible in } \cat{O}_{Y,y} (*)} \\
= & \set{y \in V}{ h_y \in \left( \left( \cat{O}_Y \right)^* \right)_y \text{ is invertible in } \left( \left( \cat{O}_Y \right)^* \right)_y} .
\end{align}
By \cref{rem:sheaf property in terms of F S}, the presheaf $\left( \cat{O}_Y \right)^*$ is in fact a sheaf of commutative unital rings on $Y$. Then the ordinary sheaf theory shows that the set
\begin{equation}
V' = \set{y \in V}{ h_y \in \left( \left( \cat{O}_Y \right)^* \right)_y \text{ is invertible in } \left( \left( \cat{O}_Y \right)^* \right)_y}
\end{equation}
is an open subset of $Y$, and the element $h|_{V'} \in \left( \cat{O}_Y \right)^*(V')$ is invertible in $\left( \cat{O}_Y \right)^*(V')$. This completes the proof.
\end{proof}

\cref{lem:descriotion of sigma inverse D(U g)} shows that $\sigma^{-1} \big( D(U;g) \big)$ is equal to the set
\begin{equation}
\set{y \in \rho^{-1}(U)}{
\Big| \Big( ( \rho^{\#}_U )_* (g) \Big)_y \Big|_{w_y} \neq 0 } ,
\end{equation}
and \cref{lem:nonvanishing locus is open in objects of C l} shows that this set is an open subset of $Y$. This completes the proof that $\sigma : Y \to \tilde{X}$ is a continuous map.

Since $\sigma : Y \to \tilde{X}$ is a continuous map, we have the sheaf $\sigma_* \cat{O}_Y$ of condensed rings on $\tilde{X}$. We construct a natural transformation $\gamma : R \To \left( \sigma_* \cat{O}_Y \right) \of D$ of functors $I \to \ub{CRing}$. Let $i = (U;g) \in I$. Then \cref{lem:descriotion of sigma inverse D(U g)} shows that
\begin{equation}
\sigma^{-1} (D(i)) = 
\set{y \in \rho^{-1}(U)}{
\Big| \Big( ( \rho^{\#}_U )_* (g) \Big)_y \Big|_{w_y} \neq 0 } ,
\end{equation}
and \cref{lem:nonvanishing locus is open in objects of C l} shows that the element
\begin{equation}
\left. ( \rho^{\#}_U )_* (g) \right|_{\sigma^{-1} (D(i))} \, \in \,
\cat{O}_Y \big( \sigma^{-1} (D(i)) \big) (*)
\end{equation}
is invertible in $\cat{O}_Y \big( \sigma^{-1} (D(i)) \big) (*)$. Then \cref{prop:universality of localization of condensed algebra} shows that there exists a unique homomorphism $\gamma_i : \cat{O}(U)_g \to \cat{O}_Y \big( \sigma^{-1} (D(i)) \big)$ of condensed rings such that the diagram
\[\begin{tikzcd}
	{\cat{O}(U)} & {\cat{O}_Y(\rho^{-1}(U))} \\
	{\cat{O}(U)_g} & {\cat{O}_Y \big( \sigma^{-1} (D(i)) \big)}
	\arrow["{\rho^{\#}_U}", from=1-1, to=1-2]
	\arrow["\can"', from=1-1, to=2-1]
	\arrow["{\cat{O}_Y \big( \sigma^{-1} (D(i)) \, \sub \, \rho^{-1}(U) \big)}", from=1-2, to=2-2]
	\arrow["{\gamma_i}"', from=2-1, to=2-2]
\end{tikzcd}\]
is commutative. If $i = (U;g) \in I$ and $j = (U';g') \in I$ satisfy $i \leq j$, then the following diagram is commutative.
\[\begin{tikzcd}
	& {\cat{O}(U)_g} \\
	{\cat{O}(U)_g} & {\cat{O}(U)} & {\cat{O}_Y(\rho^{-1}(U))} &&& {\cat{O}_Y \big( \sigma^{-1}(D(i)) \big)} \\
	& {\cat{O}(U')} & {\cat{O}_Y(\rho^{-1}(U'))} &&& {\cat{O}_Y \big( \sigma^{-1}(D(j)) \big)} \\
	& {\cat{O}(U')_{g'}}
	\arrow["{\gamma_i}", curve={height=-24pt}, from=1-2, to=2-6]
	\arrow["{R(i \leq j)}"', curve={height=18pt}, from=2-1, to=4-2]
	\arrow["\can", from=2-2, to=1-2]
	\arrow["\can", from=2-2, to=2-1]
	\arrow["{\rho^{\#}_U}", from=2-2, to=2-3]
	\arrow["{\cat{O}(U' \sub U)}"', from=2-2, to=3-2]
	\arrow["{\cat{O}_Y \big( \sigma^{-1}(D(i)) \, \sub \, \rho^{-1}(U) \big)}", from=2-3, to=2-6]
	\arrow["{\cat{O}_Y(\rho^{-1}(U') \sub \rho^{-1}(U))}", from=2-3, to=3-3]
	\arrow["{\cat{O}_Y \big( \sigma^{-1}(D(j)) \, \sub \, \sigma^{-1}(D(i)) \big)}", from=2-6, to=3-6]
	\arrow["{\rho^{\#}_U{'}}", from=3-2, to=3-3]
	\arrow["\can"', from=3-2, to=4-2]
	\arrow["{\cat{O}_Y \big( \sigma^{-1}(D(j)) \, \sub \, \rho^{-1}(U') \big)}"', from=3-3, to=3-6]
	\arrow["{\gamma_j}"', curve={height=24pt}, from=4-2, to=3-6]
\end{tikzcd}\]
By \cref{prop:universality of localization of condensed algebra}, we conclude that the following diagram is commutative.
\[\begin{tikzcd}
	{\cat{O}(U)_g} & {\cat{O}_Y \big( \sigma^{-1}(D(i)) \big)} \\
	{\cat{O}(U')_{g'}} & {\cat{O}_Y \big( \sigma^{-1}(D(j)) \big)}
	\arrow["{\gamma_i}", from=1-1, to=1-2]
	\arrow["{R(i \leq j)}"', from=1-1, to=2-1]
	\arrow["{\cat{O}_Y \big( \sigma^{-1}(D(j)) \, \sub \, \sigma^{-1}(D(i)) \big)}", from=1-2, to=2-2]
	\arrow["{\gamma_j}"', from=2-1, to=2-2]
\end{tikzcd}\]
Consequently, the family $\gamma := (\gamma_i)_{i \in I}$ is a natural transformation $R \To \left( \sigma_* \cat{O}_Y \right) \of D$ of functors $I \to \ub{CRing}$. Note that the definition shows that the following diagram is commutative.
\[\begin{tikzcd}
	{\cat{O}} && {\pi_* \sigma_* \cat{O}_Y} \\
	{R \of C} && {\left( \sigma_* \cat{O}_Y \right) \of D \of C}
	\arrow["{\rho^{\#}}", from=1-1, to=1-3]
	\arrow[equals, from=2-1, to=1-1]
	\arrow["{\gamma * \id_C}"', from=2-1, to=2-3]
	\arrow[equals, from=2-3, to=1-3]
\end{tikzcd}\]

By the definition of $\tilde{\cat{O}}$, there exists a unique morphism $\sigma^{\#} : \tilde{\cat{O}} \to \sigma_* \cat{O}_Y$ of sheaves of condensed rings on $\tilde{X}$ such that the diagram
\[\begin{tikzcd}
	R & {\left( \sigma_* \cat{O}_Y \right) \of D} \\
	{\tilde{\cat{O}} \of D}
	\arrow["\gamma", from=1-1, to=1-2]
	\arrow["\beta"', from=1-1, to=2-1]
	\arrow["{\sigma^{\#} * \id_D}"', from=2-1, to=1-2]
\end{tikzcd}\]
is commutative. Thus we have a morphism $(\sigma, \sigma^{\#}) : (Y, \cat{O}_Y) \to (\tilde{X}, \tilde{\cat{O}})$ in the category $\cat{D}$. Moreover, the following diagrams are commutative. 
\[\begin{tikzcd}
	&&&&&& {\pi_* \tilde{\cat{O}}} \\
	Y & {\tilde{X}} & {\pi_* \sigma_* \cat{O}_Y} & {\left( \sigma_* \cat{O}_Y \right) \of D \of C} && {\tilde{\cat{O}} \of D \of C} \\
	& X &&&& {R \of C} \\
	&&&&& {\cat{O}}
	\arrow["{\pi_* \sigma^{\#}}"', curve={height=18pt}, from=1-7, to=2-3]
	\arrow["\sigma", from=2-1, to=2-2]
	\arrow["\rho"', from=2-1, to=3-2]
	\arrow["\pi", from=2-2, to=3-2]
	\arrow[equals, from=2-4, to=2-3]
	\arrow[equals, from=2-6, to=1-7]
	\arrow["{\sigma^{\#} * \id_D * \id_C}"', from=2-6, to=2-4]
	\arrow["{\gamma * \id_C}", from=3-6, to=2-4]
	\arrow["{\beta * \id_C}"', from=3-6, to=2-6]
	\arrow[equals, from=3-6, to=4-6]
	\arrow["{\pi^{\#}}"', curve={height=18pt}, from=4-6, to=1-7]
	\arrow["{\rho^{\#}}", curve={height=-12pt}, from=4-6, to=2-3]
\end{tikzcd}\]
Therefore the diagram
\[\begin{tikzcd}
	{(Y,\cat{O}_Y)} & {(\tilde{X},\tilde{\cat{O}})} \\
	& {(X,\cat{O})}
	\arrow["{(\sigma,\sigma^{\#})}", from=1-1, to=1-2]
	\arrow["{(\rho,\rho^{\#})}"', from=1-1, to=2-2]
	\arrow["{(\pi,\pi^{\#})}", from=1-2, to=2-2]
\end{tikzcd}\]
is commutative in $\cat{D}$.

Let us prove that $(\sigma, \sigma^{\#})$ is a morphism $(Y, \cat{O}_Y, \cat{V}_Y) \to (\tilde{X}, \tilde{\cat{O}}, \tilde{\cat{V}})$ in $\cat{C}_l$. For every $y \in Y$, the following diagram is commutative.
\[\begin{tikzcd}
	{\cat{O}_{Y,y}} & {\tilde{\cat{O}}_{\sigma(y)}} \\
	& {\cat{O}_{\rho(y)}}
	\arrow["{\sigma^{\#}_{y}}"', from=1-2, to=1-1]
	\arrow["{\rho^{\#}_y}", from=2-2, to=1-1]
	\arrow["{\pi^{\#}_{\sigma(y)}}"', from=2-2, to=1-2]
\end{tikzcd}\]
Then we have
\begin{equation}
\left( \pi^{\#}_{\sigma(y)} \right)^{-1} \left( \left( \sigma^{\#}_y \right)^{-1} (w_y) \right)
= \left( \rho^{\#}_y \right) ^{-1} (w_y) .
\end{equation}
On the other hand, since $(\rho, \rho^{\#}) : (Y, \cat{O}_Y, \cat{V}_Y) \to (X, \cat{O}, \cat{V})$ is a morphism in $\cat{C}_1$, we have
\begin{equation}
\left( \rho^{\#}_y \right) ^{-1} (w_y) = v_{\rho(y)} .
\end{equation}
Therefore
\begin{equation}
\left( \pi^{\#}_{\sigma(y)} \right)^{-1} \left( \left( \sigma^{\#}_y \right)^{-1} (w_y) \right)
= \left( \rho^{\#}_y \right) ^{-1} (w_y) 
= v_{\rho(y)} = v_{\sigma(y)} .
\end{equation}
Since $\tilde{v}_{\sigma(y)}$ is the unique continuous valuation on $\tilde{\cat{O}}_{\sigma(y)}$ such that $\left( \pi^{\#}_{\sigma(y)} \right)^{-1} ( \tilde{v}_{\sigma(y)} ) = v_{\sigma(y)}$, it follows that
\begin{equation}
\left( \sigma^{\#}_y \right)^{-1} (w_y) = \tilde{v}_{\sigma(y)} .
\end{equation}
This shows that $(\sigma, \sigma^{\#}) : (Y, \cat{O}_Y, \cat{V}_Y) \to (\tilde{X}, \tilde{\cat{O}}, \tilde{\cat{V}})$ is a morphism in $\cat{C}_l$.

Since the diagram
\[\begin{tikzcd}
	{(Y, \cat{O}_Y)} & {(\tilde{X}, \tilde{\cat{O}})} \\
	& {(X, \cat{O})}
	\arrow["{(\sigma,\sigma^{\#})}", from=1-1, to=1-2]
	\arrow["{(\rho,\rho^{\#})}"', from=1-1, to=2-2]
	\arrow["{(\pi, \pi^{\#})}", from=1-2, to=2-2]
\end{tikzcd}\]
is commutative in $\cat{D}$, the diagram
\[\begin{tikzcd}
	{(Y, \cat{O}_Y, \cat{V}_Y)} & {(\tilde{X}, \tilde{\cat{O}}, \tilde{\cat{V}})} \\
	& {(X, \cat{O}, \cat{V})}
	\arrow["{(\sigma,\sigma^{\#})}", from=1-1, to=1-2]
	\arrow["{(\rho,\rho^{\#})}"', from=1-1, to=2-2]
	\arrow["{(\pi, \pi^{\#})}", from=1-2, to=2-2]
\end{tikzcd}\]
is commutative in $\cat{C}_1$.

On the other hand, suppose that $(\tau, \tau^{\#}) : (Y, \cat{O}_Y, \cat{V}_Y) \to (\tilde{X}, \tilde{\cat{O}}, \tilde{\cat{V}})$ is another morphism in $\cat{C}_l$ such that the following diagram is commutative.
\[\begin{tikzcd}
	{(Y, \cat{O}_Y, \cat{V}_Y)} & {(\tilde{X}, \tilde{\cat{O}}, \tilde{\cat{V}})} \\
	& {(X, \cat{O}, \cat{V})}
	\arrow["{(\tau,\tau^{\#})}", from=1-1, to=1-2]
	\arrow["{(\rho,\rho^{\#})}"', from=1-1, to=2-2]
	\arrow["{(\pi, \pi^{\#})}", from=1-2, to=2-2]
\end{tikzcd}\]
We prove that $(\tau, \tau^{\#}) = (\sigma, \sigma^{\#})$. The diagram
\[\begin{tikzcd}
	Y & {\tilde{X}} \\
	& X
	\arrow["\tau", from=1-1, to=1-2]
	\arrow["\rho"', from=1-1, to=2-2]
	\arrow["\pi", from=1-2, to=2-2]
\end{tikzcd}\]
is commutative. Since the map $\pi : \tilde{X} \to X$ is the identity map $x \mapsto x$, the map $\tau : Y \to \tilde{X}$ must be equal to the map $y \mapsto \rho(y)$. On the other hand, the map $\sigma : Y \to \tilde{X}$ was defined to be the map $y \mapsto \rho(y)$. Therefore $\tau = \sigma$. Then the morphism $\tau^{\#} : \tilde{\cat{O}} \to \tau_* \cat{O}_Y$ of sheaves of condensed rings on $\tilde{X}$ is a morphism $\tilde{\cat{O}} \to \sigma_* \cat{O}_Y$. Moreover, for each $i = (U;g) \in I$, the following diagram is commutative.
\[\begin{tikzcd}
	& {\cat{O}_Y \big( \sigma^{-1} (D(i)) \big)} &&& {\tilde{\cat{O}}(D(i))} \\
	& {\cat{O}_Y(\rho^{-1}(U))} && {\tilde{\cat{O}}(\pi^{-1}(U))} & {\tilde{\cat{O}}(D(C(U)))} \\
	{\cat{O}(U)_g} &&& {\cat{O}(U)} & {R(C(U))} \\
	&&& {\cat{O}(U)} && {\cat{O}(U)_g}
	\arrow["{\tau^{\#}_{D(i)}}"', from=1-5, to=1-2]
	\arrow["{\cat{O}_Y \big( \sigma^{-1} (D(i)) \, \sub \, \rho^{-1}(U) \big)}"', from=2-2, to=1-2]
	\arrow["{\tau^{\#}_{\pi^{-1}(U)}}"{pos=0.3}, from=2-4, to=2-2]
	\arrow[equals, from=2-4, to=2-5]
	\arrow["{\tilde{\cat{O}}(D(i) \sub D(C(U)))}", from=2-5, to=1-5]
	\arrow["{\gamma_i}", curve={height=-30pt}, from=3-1, to=1-2]
	\arrow["{\rho^{\#}_U}", from=3-4, to=2-2]
	\arrow["{\pi^{\#}_U}"', from=3-4, to=2-4]
	\arrow["\can", from=3-4, to=3-1]
	\arrow[equals, from=3-4, to=3-5]
	\arrow[equals, from=3-4, to=4-4]
	\arrow["{\beta_{C(U)}}", from=3-5, to=2-5]
	\arrow["{R(C(U) \leq i)}"'{pos=0.3}, from=3-5, to=4-6]
	\arrow["\can"', from=4-4, to=4-6]
	\arrow["{\beta_i}"', curve={height=30pt}, from=4-6, to=1-5]
\end{tikzcd}\]
By \cref{prop:universality of localization of condensed algebra}, we conclude that the folloing diagram is commutative.
\[\begin{tikzcd}
	{\cat{O}_Y \big( \sigma^{-1} (D(i)) \big)} && {\tilde{\cat{O}}(D(i))} \\
	&& {\cat{O}(U)_g}
	\arrow["{\tau^{\#}_{D(i)}}"', from=1-3, to=1-1]
	\arrow["{\gamma_i}", from=2-3, to=1-1]
	\arrow["{\beta_i}"', from=2-3, to=1-3]
\end{tikzcd}\]
Since this holds for every $i = (U;g) \in I$, it follows that the following diagram is commutative.
\[\begin{tikzcd}
	R & {\left( \sigma_* \cat{O}_Y \right) \of D} \\
	{\tilde{\cat{O}} \of D}
	\arrow["\gamma", from=1-1, to=1-2]
	\arrow["\beta"', from=1-1, to=2-1]
	\arrow["{\tau^{\#} * \id_D}"', from=2-1, to=1-2]
\end{tikzcd}\]
By the definition of $\sigma^{\#}$, we conclude that $\tau^{\#} = \sigma^{\#}$. This completes the proof of \cref{prop:construction of coreflection in C l}.
\end{proof}

\subsection{The category $\cat{C}_f$}

\subsubsection{Definition}

\begin{df}
The category $\cat{C}_f$ is defined to be the full subcategory of $\cat{C}_1$ consisting of all $(X, \cat{O}_X, \cat{V}_X) \in |\cat{C}_1|$ with the following property: For every open subset $U$ of $X$ and every $f, g \in \cat{O}_X(U)(*)$, the set
\begin{equation}
\set{x \in U}{|f_x|_x \leq |g_x|_x \neq 0}
\end{equation}
is an open subset of $X$.
\end{df}

\subsubsection{Coreflection}

\begin{prop} \label{prop:construction of coreflection in C f}
Let $(X, \cat{O}, \cat{V})$ be an object of $\cat{C}_1$. Then there exists a coreflection
\begin{equation}
\left( (\tilde{X}, \tilde{\cat{O}}, \tilde{\cat{V}}) ,\,  
(\tilde{X}, \tilde{\cat{O}}, \tilde{\cat{V}}) \xto{(\pi, \pi^{\#})} (X, \cat{O}, \cat{V})
\right)
\end{equation}
of $(X, \cat{O}, \cat{V})$ along the inclusion functor $\cat{C}_f \mon \cat{C}_1$ which has the following property.
\begin{enumerate}
\item
The underlying set of $\tilde{X}$ is equal to the underlying set of $X$.

\item
The map $\tilde{X} \xto{\pi} X$ is equal to the identity map $x \mapsto x$.

\item
Let $I$ be the set of all triples $(U;T,g)$ consisting of an open subset $U$ of $X$, a nonempty finite subset $T$ of $\cat{O}(U)(*)$ and an element $g \in \cat{O}(U)(*)$. For $(U;T,g) \in I$, write $D(U;T,g)$ for the set of all $x \in U$ such that $|f_x|_x \leq |g_x|_x \neq 0$ for every $f \in T$. Then the set
\begin{equation}
\set{D(U;T,g)}{(U;T,g) \in I}
\end{equation}
is a basis for the topology of $\tilde{X}$.

\item
For every $x \in \tilde{X}$, the homomorphism $\pi^{\#}_x : (\cat{O}_x, |\cdot|_x) \to (\tilde{\cat{O}}_x, |\cdot|_x)$ is an isomorphism of valued condensed rings.
\end{enumerate}
\end{prop}

\subsubsection{Proof of \cref{prop:construction of coreflection in C f}}
\label{sec:Proof of construction of coreflection in C f}

\begin{proof}
Let $(X, \cat{O}, \cat{V})$ be an object of $\cat{C}_1$. For each $x \in X$, let us write $v_x$ for the valuation on $\cat{O}_x$ which is the unique element of $\cat{V}_x$. Let $I$ be the set of all triples $(U;T,g)$ consisting of an open subset $U$ of $X$, a nonempty finite subset $T$ of $\cat{O}(U)(*)$ and an element $g \in \cat{O}(U)(*)$. For each $(U;T,g) \in I$, write
\begin{equation}
D(U;T,g) := \set{x \in U}{ |f_x|_{v_x} \leq |g_x|_{v_x} \neq 0 \text{ for every } f \in T} .
\end{equation}

\begin{lem} \label{lem:properties of D(U T g)} \;
\begin{enumerate}
\item
Suppose that $(U;T,g) ,\, (U';T',g') \in I$. Let us write
\begin{align}
U'' & := U \cap U' \: ; \\
T'' & := \set{ (f|_{U''}) \cdot (g'|_{U''}) \in \cat{O}(U'')(*)}{f \in T} \cup
\set{ (f'|_{U''}) \cdot (g|_{U''}) \in \cat{O}(U'')(*)}{f' \in T'} \: ; \\
g'' & := (g|_{U''}) \cdot (g'|_{U''}) \in \cat{O}(U'')(*) .
\end{align}
Then we have $(U'';T'',g'') \in I$ and
\begin{equation}
D(U;T,g) \cap D(U';T',g') = D(U'';T'',g'') .
\end{equation}

\item
For every open subset $U$ of $X$, we have
\begin{equation}
U = D(U; \{1\}, 1) .
\end{equation}
In particular, we have $X = D(X; \{1\}, 1)$.

\end{enumerate}
\end{lem}

\begin{proof}~
\begin{enumerate}
\item
$U''$ is an open subset of $X$ which is contained in both $U$ and $U'$. Therefore $f|_{U''}, f'|_{U''}$ $(f \in T ,\; f' \in T')$ and $g|_{U''}, g'|_{U''}$ are well-defined elements of $\cat{O}(U'')(*)$. Thus $g'' = (g|_{U''}) \cdot (g'|_{U''})$ is an element $\cat{O}(U'')(*)$ and $T'' = \set{ (f|_{U''}) \cdot (g'|_{U''}) }{f \in T} \cup \set{ (f'|_{U''}) \cdot (g|_{U''}) }{f' \in T'}$ is a subset of $\cat{O}(U'')(*)$. Moreover, the set $T''$ is nonempty and finite since the sets $T,T'$ are nonempty finite sets. Consequently we have $(U'';T'',g'') \in I$. We prove that $D(U;T,g) \cap D(U';T',g') = D(U'';T'',g'')$. First suppose $x \in D(U'';T'',g'')$. Then $x \in U'' = U \cap U'$. For every $f \in T$, we have $|f_x \cdot g'_x|_{v_x} \leq |g_x \cdot g'_x|_{v_x} \neq 0$. Then we have $|g'_x|_{v_x} \neq 0$ and hence $|f_x|_{v_x} \leq |g_x|_{v_x} \neq 0$. Therefore $x \in D(U;T,g)$. Similarly, for every $f' \in T'$, we have $|f'_x \cdot g_x|_{v_x} \leq |g_x \cdot g'_x|_{v_x} \neq 0$. Then we have $|g_x|_{v_x} \neq 0$ and hence $|f'_x|_{v_x} \leq |g'_x|_{v_x} \neq 0$. Therefore $x \in D(U';T',g')$. Thus $x \in D(U;T,g) \cap D(U';T',g')$. Conversely, suppose that $x \in D(U;T,g) \cap D(U';T',g')$. Then $x \in U$ and $x \in U'$. Therefore $x \in U''$. Moreover, if $f \in T$, then $|f_x|_{v_x} \leq |g_x|_{v_x} \neq 0$. Since $T'$ is nonempty, there exists some element $f' \in T'$. Then $|f'_x|_{v_x} \leq |g'_x|_{v_x} \neq 0$. In particular, we have $|g'_x|_{v_x} \neq 0$. Therefore $|f_x \cdot g'_x|_{v_x} \leq |g_x \cdot g'_x|_{v_x} \neq 0$. Similarly, if $f' \in T'$, then $|f'_x|_{v_x} \leq |g'_x|_{v_x} \neq 0$. Since $T$ is nonempty, there exists some element $f \in T$. Then $|f_x|_{v_x} \leq |g_x|_{v_x} \neq 0$. In particular, we have $|g_x|_{v_x} \neq 0$. Therefore $|f'_x \cdot g_x|_{v_x} \leq |g_x \cdot g'_x|_{v_x} \neq 0$. Consequently we have $x \in D(U'';T'',g'')$.

\item
Let $U$ be an open subset of $X$. For every $x \in U$, we have $|1_x|_{v_x} \leq |1_x|_{v_x} = 1 \neq 0$. Therefore $U \sub D(U; \{1\}, 1)$. Since $D(U; \{1\}, 1) \sub U$ by defintion, we have $U = D(U; \{1\}, 1)$.
\end{enumerate}
\end{proof}

By \cref{lem:properties of D(U T g)}, there exists a unique topology $\cat{T}$ on the underlying set of $X$ for which the set
\begin{equation}
\set{D(U;T,g)}{(U;T,g) \in I}
\end{equation}
is an open basis. Let $\tilde{X}$ be the set $X$ endowed with this topology. Let $\pi : \tilde{X} \to X$ be the identity map $x \mapsto x$. Then (2) of \cref{lem:properties of D(U T g)} shows that $\pi : \tilde{X} \to X$ is continuous. Let $\tilde{\cat{O}}$ be the inverse image $\pi^{-1} \cat{O}$ of $\cat{O}$ under $\pi$. Let us write $\pi^{\#} : \cat{O} \to \pi_* \tilde{\cat{O}}$ for the canonical morphism of sheaves of condensed rings on $X$. Then \cref{cor:existence of inverse image functor} shows that for each $x \in \tilde{X}$, the homomorphism $\pi^{\#}_x : \cat{O}_x \to \tilde{\cat{O}}_x$ is an isomorphism of condensed rings. Therefore there exists a unique valuation $\tilde{v}_x$ on $\tilde{\cat{O}}_x$ such that $\left( \pi^{\#}_x \right)^{-1} ( \tilde{v}_x ) = v_x$. Then $\tilde{v}_x$ is a continuous valuation on $\tilde{\cat{O}}_x$ since $v_x$ is a continuous valuation on $\cat{O}_x$. We define $\tilde{\cat{V}}_x := \{ \tilde{v}_x \}$. Then we obtain a family $\tilde{\cat{V}} := (\tilde{\cat{V}}_x)_{x \in \tilde{X}}$. Then $(\tilde{X}, \tilde{\cat{O}}, \tilde{\cat{V}})$ is an object of $\cat{C}_1$, and $(\pi, \pi^{\#})$ is a morphism $(\tilde{X}, \tilde{\cat{O}}, \tilde{\cat{V}}) \to (X, \cat{O}, \cat{V})$ in $\cat{C}$.

Let us prove that $(\tilde{X}, \tilde{\cat{O}}, \tilde{\cat{V}})$ is an object of $\cat{C}_f$. Let $\tilde{U}$ be an open subset of $\tilde{X}$. Let $\tilde{f},\tilde{g} \in \tilde{\cat{O}}(\tilde{U})(*)$. Write
\begin{equation}
\tilde{V} :=
\set{ x \in \tilde{U} }{
| \tilde{f}_x |_{\tilde{v}_x} \leq
| \tilde{g}_x |_{\tilde{v}_x} \neq 0 } .
\end{equation}
We prove that $\tilde{V}$ is an open subset of $\tilde{X}$. For this it suffices to show that for every $x \in \tilde{U}$, there exists an open neighbourhood $\tilde{W}$ of $x$ in $\tilde{X}$ such that $\tilde{V} \cap \tilde{W}$ is open in $\tilde{X}$. Therefore let us fix $x \in \tilde{U}$ arbitrarily. Since the homomorphism $\pi^{\#}_x : \cat{O}_x \to \tilde{\cat{O}}_x$ is an isomorphism of condensed rings, there exist elements $s,t \in \cat{O}_x(*)$ such that $(\pi^{\#}_x)_*(s) = \tilde{f}_x$ and $(\pi^{\#}_x)_*(t) = \tilde{g}_x$. By \cref{prop:compatibility of stalk and evaluation at S}, we have
\begin{equation}
\cat{O}_x (*) = \underset{x \in U \sub X \text{ open}}{\colim} \; \cat{O}(U)(*) ,
\end{equation}
where the colimit is taken in $\ub{Ring}$. Then the construction of filtered colimits in $\ub{Ring}$ shows that there exist an open neighbourhood $U$ of $x$ in $X$ and elements $f,g \in \cat{O}(U)(*)$ such that $f_x = s$ and $g_x = t$. Then we have an element $(U; \{f\}, g) \in I$. On the other hand, the commutativity of the diagram
\[\begin{tikzcd}
	{\cat{O}(U)} & {\tilde{\cat{O}}(\pi^{-1}(U))} \\
	{\cat{O}_x} & {\tilde{\cat{O}}_x}
	\arrow["{\pi^{\#}_U}", from=1-1, to=1-2]
	\arrow["\can"', from=1-1, to=2-1]
	\arrow["\can", from=1-2, to=2-2]
	\arrow["{\pi^{\#}_x}"', from=2-1, to=2-2]
\end{tikzcd}\]
shows that
\begin{align}
\big( (\pi^{\#}_U)_*(f) \big)_x = (\pi^{\#}_x)_*(f_x) = (\pi^{\#}_x)_*(s) = \tilde{f}_x \: ; \\
\big( (\pi^{\#}_U)_*(g) \big)_x = (\pi^{\#}_x)_*(g_x) = (\pi^{\#}_x)_*(t) = \tilde{g}_x .
\end{align}
By \cref{prop:compatibility of stalk and evaluation at S}, we have
\begin{equation}
\tilde{\cat{O}}_x (*) =
\underset{x \in \tilde{W} \sub \tilde{X} \text{ open}}{\colim} \; 
\tilde{\cat{O}}(\tilde{W})(*) ,
\end{equation}
where the colimit is taken in $\ub{Ring}$. Then the construction of filtered colimits in $\ub{Ring}$ shows that there exists an open neighbourhood $\tilde{W}$ of $x$ in $\tilde{X}$ such that $\tilde{W} \sub \pi^{-1}(U) \cap \tilde{U}$ and 
\begin{equation}
(\pi^{\#}_U)_*(f)|_{\tilde{W}} = \tilde{f}|_{\tilde{W}} \quad \text{ and } \quad
(\pi^{\#}_U)_*(g)|_{\tilde{W}} = \tilde{g}|_{\tilde{W}} .
\end{equation}
We claim that $\tilde{V} \cap \tilde{W} = D(U; \{f\}, g) \cap \tilde{W}$. Indeed, for $y \in \tilde{W}$, we have
\begin{equation}
\big( (\pi^{\#}_U)_*(f) \big)_y = \tilde{f}_y \quad \text{ and } \quad
\big( (\pi^{\#}_U)_*(g) \big)_y = \tilde{g}_y 
\end{equation}
since $(\pi^{\#}_U)_*(f)|_{\tilde{W}} = \tilde{f}|_{\tilde{W}}$ and $(\pi^{\#}_U)_*(g)|_{\tilde{W}} = \tilde{g}|_{\tilde{W}}$. Moreover, the commutativity of the diagram
\[\begin{tikzcd}
	{\cat{O}(U)} & {\tilde{\cat{O}}(\pi^{-1}(U))} \\
	{\cat{O}_y} & {\tilde{\cat{O}}_y}
	\arrow["{\pi^{\#}_U}", from=1-1, to=1-2]
	\arrow["\can"', from=1-1, to=2-1]
	\arrow["\can", from=1-2, to=2-2]
	\arrow["{\pi^{\#}_y}"', from=2-1, to=2-2]
\end{tikzcd}\]
shows that
\begin{equation}
\big( (\pi^{\#}_U)_*(f) \big)_y = (\pi^{\#}_y)_*(f_y) \quad \text{ and } \quad
\big( (\pi^{\#}_U)_*(g) \big)_y = (\pi^{\#}_y)_*(g_y) .
\end{equation}
Consequently, we have $\tilde{f}_y = (\pi^{\#}_y)_*(f_y)$ and $\tilde{g}_y = (\pi^{\#}_y)_*(g_y)$. On the other hand, the definition of $\tilde{v}_y$ shows that $(\pi^{\#}_y)^{-1}(\tilde{v}_y) = v_y$. Therefore the condition
\begin{equation}
| \tilde{f}_y |_{\tilde{v}_y} \leq | \tilde{g}_y |_{\tilde{v}_y} \neq 0
\end{equation}
is equivalent to the condition
\begin{equation}
| f_y |_{v_y} \leq | g_y |_{v_y} \neq 0 .
\end{equation}
Therefore we have $y \in \tilde{V}$ if and only if $y \in D(U; \{f\}, g)$. This shows that $\tilde{V} \cap \tilde{W} = D(U; \{f\}, g) \cap \tilde{W}$. Thus we have proved that there exists an open neighbourhood $\tilde{W}$ of $x$ in $\tilde{X}$ such that $\tilde{V} \cap \tilde{W} = D(U; \{f\}, g) \cap \tilde{W}$. Since both $D(U; \{f\}, g)$ and $\tilde{W}$ are open subsets of $\tilde{X}$, we conclude that $\tilde{V} \cap \tilde{W}$ is open in $\tilde{X}$. This completes the proof that $(\tilde{X}, \tilde{\cat{O}}, \tilde{\cat{V}})$ is an object of $\cat{C}_f$.

Next suppose that $(Y, \cat{O}_Y, \cat{V}_Y)$ is an object of $\cat{C}_f$ and that $(\rho, \rho^{\#}) : (Y, \cat{O}_Y, \cat{V}_Y) \to (X, \cat{O}, \cat{V})$ is a morphism in $\cat{C}_1$. We prove that there exists a unique morphism $(\sigma, \sigma^{\#}) : (Y, \cat{O}_Y, \cat{V}_Y) \to (\tilde{X}, \tilde{\cat{O}}, \tilde{\cat{V}})$ in $\cat{C}_f$ such that the following diagram is commutative.
\[\begin{tikzcd}
	{(Y, \cat{O}_Y, \cat{V}_Y)} & {(\tilde{X}, \tilde{\cat{O}}, \tilde{\cat{V}})} \\
	& {(X, \cat{O}, \cat{V})}
	\arrow["{(\sigma,\sigma^{\#})}", from=1-1, to=1-2]
	\arrow["{(\rho,\rho^{\#})}"', from=1-1, to=2-2]
	\arrow["{(\pi, \pi^{\#})}", from=1-2, to=2-2]
\end{tikzcd}\]

For each $y \in Y$, let us write $w_y$ for the valuation on $\cat{O}_{Y,y}$ which is the unique element of $\cat{V}_{Y,y}$. Let $\sigma :Y \to \tilde{X}$ be the map $y \mapsto \rho(y)$. We claim that this map $\sigma :Y \to \tilde{X}$ is continuous. Since the set $\set{D(U;T,g)}{(U;T,g) \in I}$ is a basis for the topology of $\tilde{X}$, it suffices to show that for every $(U;T,g) \in I$, the set $\sigma^{-1} \big( D(U;T,g) \big)$ is an open subset of $Y$. Let $(U;T,g) \in I$. Then we have a homomorphism $\rho^{\#}_U : \cat{O}(U) \to \cat{O}_Y(\rho^{-1}(U))$ of condensed rings. Thus for each $h \in \cat{O}(U)(*)$, we have an element $( \rho^{\#}_U )_* (h) \in \cat{O}_Y ( \rho^{-1}(U) )(*)$.

\begin{lem} \label{lem:descriotion of sigma inverse D(U T g)}
The following equality holds.
\begin{equation}
\sigma^{-1} \big( D(U;T,g) \big)
= \bigcap_{f \in T} \set{y \in \rho^{-1}(U)}{
\Big| \Big( ( \rho^{\#}_U )_* (f) \Big)_y \Big|_{w_y} \leq
\Big| \Big( ( \rho^{\#}_U )_* (g) \Big)_y \Big|_{w_y} \neq 0 } .
\end{equation}
\end{lem}

\begin{proof}
We have
\begin{equation}
\sigma^{-1} \big( D(U;T,g) \big)
= \set{y \in Y}{ \rho(y) \in U \, ; \,
\left| f_{\rho(y)} \right|_{v_{\rho(y)}} \leq 
\left| g_{\rho(y)} \right|_{v_{\rho(y)}} \neq 0
\text{ for every } f \in T } .
\end{equation}
For each $y \in Y$, we have a homomorphism $\rho^{\#}_y : \cat{O}_{\rho(y)} \to \cat{O}_{Y,y}$ of condensed rings. Since $(\rho, \rho^{\#}) : (Y, \cat{O}_Y, \cat{V}_Y) \to (X, \cat{O}, \cat{V})$ is a morphism in $\cat{C}_1$, we have
\begin{equation}
v_{\rho(y)} = \left( \rho^{\#}_y \right) ^{-1} (w_y) .
\end{equation}
Therefore
\begin{align}
& \set{y \in Y}{ \rho(y) \in U \, ; \,
\left| f_{\rho(y)} \right|_{v_{\rho(y)}} \leq 
\left| g_{\rho(y)} \right|_{v_{\rho(y)}} \neq 0
\text{ for every } f \in T } \\
= & \set{y \in Y}{ \rho(y) \in U \, ; \,
\left| f_{\rho(y)} \right|_{\left( \rho^{\#}_y \right) ^{-1} (w_y)} \leq 
\left| g_{\rho(y)} \right|_{\left( \rho^{\#}_y \right) ^{-1} (w_y)} \neq 0
\text{ for every } f \in T } \\
= & \set{y \in \rho^{-1}(U)}{
\left| (\rho^{\#}_y)_* ( f_{\rho(y)} ) \right|_{w_y} \leq 
\left| (\rho^{\#}_y)_* ( g_{\rho(y)} ) \right|_{w_y} \neq 0
\text{ for every } f \in T } .
\end{align}
For each $y \in \rho^{-1}(U)$, the following diagram is commutative.
\[\begin{tikzcd}
	{\cat{O}(U)(*)} & {\cat{O}_Y(\rho^{-1}(U))(*)} \\
	{\cat{O}_{\rho(y)}(*)} & {\cat{O}_{Y,y}(*)}
	\arrow["{(\rho^{\#}_U)_*}", from=1-1, to=1-2]
	\arrow["\can"', from=1-1, to=2-1]
	\arrow["\can", from=1-2, to=2-2]
	\arrow["{(\rho^{\#}_y)_*}"', from=2-1, to=2-2]
\end{tikzcd}\]
It follows that 
\begin{align}
& \set{y \in \rho^{-1}(U)}{
\left| (\rho^{\#}_y)_* ( f_{\rho(y)} ) \right|_{w_y} \leq 
\left| (\rho^{\#}_y)_* ( g_{\rho(y)} ) \right|_{w_y} \neq 0
\text{ for every } f \in T } \\
= & \set{y \in \rho^{-1}(U)}{
\Big| \Big( ( \rho^{\#}_U )_* (f) \Big)_y \Big|_{w_y} \leq
\Big| \Big( ( \rho^{\#}_U )_* (g) \Big)_y \Big|_{w_y} \neq 0
\; \text{ for every } f \in T } .
\end{align}
Thus
\begin{align}
\sigma^{-1} \big( D(U;T,g) \big)
& = \set{y \in Y}{ \rho(y) \in U \, ; \,
\left| f_{\rho(y)} \right|_{v_{\rho(y)}} \leq 
\left| g_{\rho(y)} \right|_{v_{\rho(y)}} \neq 0
\text{ for every } f \in T } \\
& = \set{y \in \rho^{-1}(U)}{
\left| (\rho^{\#}_y)_* ( f_{\rho(y)} ) \right|_{w_y} \leq 
\left| (\rho^{\#}_y)_* ( g_{\rho(y)} ) \right|_{w_y} \neq 0
\text{ for every } f \in T } \\
& = \set{y \in \rho^{-1}(U)}{
\Big| \Big( ( \rho^{\#}_U )_* (f) \Big)_y \Big|_{w_y} \leq
\Big| \Big( ( \rho^{\#}_U )_* (g) \Big)_y \Big|_{w_y} \neq 0
\; \text{ for every } f \in T } \\
& = \bigcap_{f \in T} \set{y \in \rho^{-1}(U)}{
\Big| \Big( ( \rho^{\#}_U )_* (f) \Big)_y \Big|_{w_y} \leq
\Big| \Big( ( \rho^{\#}_U )_* (g) \Big)_y \Big|_{w_y} \neq 0 } .
\end{align}
\end{proof}

Since $(Y, \cat{O}_Y, \cat{V}_Y)$ is an object of $\cat{C}_f$, the set
\begin{equation}
\set{y \in \rho^{-1}(U)}{
\Big| \Big( ( \rho^{\#}_U )_* (f) \Big)_y \Big|_{w_y} \leq
\Big| \Big( ( \rho^{\#}_U )_* (g) \Big)_y \Big|_{w_y} \neq 0 } 
\end{equation}
is an open subset of $Y$ for every $f \in T$. Since $T$ is a finite set, \cref{lem:descriotion of sigma inverse D(U T g)} shows that $\sigma^{-1} \big( D(U;T,g) \big)$ is equal to a finite intersection of such open subsets of $Y$. It follows that $\sigma^{-1} \big( D(U;T,g) \big)$ is an open subset of $Y$. This completes the proof that $\sigma : Y \to \tilde{X}$ is a continuous map. 

The following diagram is commutative.
\[\begin{tikzcd}
	Y & {\tilde{X}} \\
	& X
	\arrow["\sigma", from=1-1, to=1-2]
	\arrow["\rho"', from=1-1, to=2-2]
	\arrow["\pi", from=1-2, to=2-2]
\end{tikzcd}\]
Then the equality
\begin{equation}
\rho_* \cat{O}_Y = \pi_* \sigma_* \cat{O}_Y
\end{equation}
holds. Therefore the morphism $\rho^{\#} : \cat{O} \to \rho_* \cat{O}_Y$ of sheaves of condensed rings on $X$ is a morphism $\cat{O} \to \pi_* \sigma_* \cat{O}_Y$. Since $\tilde{\cat{O}}$ is the inverse image of $\cat{O}$ under $\pi$, there exists a unique morphism $\sigma^{\#} : \tilde{\cat{O}} \to \sigma_* \cat{O}_Y$ of sheaves of condensed rings on $\tilde{X}$ such that the following diagram is commutative.
\[\begin{tikzcd}
	{\pi_* \sigma_* \cat{O}_Y} & {\pi_* \tilde{\cat{O}}} \\
	& {\cat{O}}
	\arrow["{\pi_* \sigma^{\#}}"', from=1-2, to=1-1]
	\arrow["{\rho^{\#}}", from=2-2, to=1-1]
	\arrow["{\pi^{\#}}"', from=2-2, to=1-2]
\end{tikzcd}\]
Thus we have a morphism $(\sigma, \sigma^{\#}) : (Y, \cat{O}_Y) \to (\tilde{X}, \tilde{\cat{O}})$ in $\cat{D}$ such that the following diagram is commutative.
\[\begin{tikzcd}
	{(Y, \cat{O}_Y)} & {(\tilde{X}, \tilde{\cat{O}})} \\
	& {(X, \cat{O})}
	\arrow["{(\sigma,\sigma^{\#})}", from=1-1, to=1-2]
	\arrow["{(\rho,\rho^{\#})}"', from=1-1, to=2-2]
	\arrow["{(\pi, \pi^{\#})}", from=1-2, to=2-2]
\end{tikzcd}\]

Let us prove that $(\sigma, \sigma^{\#})$ is a morphism $(Y, \cat{O}_Y, \cat{V}_Y) \to (\tilde{X}, \tilde{\cat{O}}, \tilde{\cat{V}})$ in $\cat{C}_1$. For every $y \in Y$, the following diagram is commutative.
\[\begin{tikzcd}
	{\cat{O}_{Y,y}} & {\tilde{\cat{O}}_{\sigma(y)}} \\
	& {\cat{O}_{\rho(y)}}
	\arrow["{\sigma^{\#}_{y}}"', from=1-2, to=1-1]
	\arrow["{\rho^{\#}_y}", from=2-2, to=1-1]
	\arrow["{\pi^{\#}_{\sigma(y)}}"', from=2-2, to=1-2]
\end{tikzcd}\]
Then we have
\begin{equation}
\left( \pi^{\#}_{\sigma(y)} \right)^{-1} \left( \left( \sigma^{\#}_y \right)^{-1} (w_y) \right)
= \left( \rho^{\#}_y \right) ^{-1} (w_y) .
\end{equation}
On the other hand, since $(\rho, \rho^{\#}) : (Y, \cat{O}_Y, \cat{V}_Y) \to (X, \cat{O}, \cat{V})$ is a morphism in $\cat{C}_1$, we have
\begin{equation}
\left( \rho^{\#}_y \right) ^{-1} (w_y) = v_{\rho(y)} .
\end{equation}
Therefore
\begin{equation}
\left( \pi^{\#}_{\sigma(y)} \right)^{-1} \left( \left( \sigma^{\#}_y \right)^{-1} (w_y) \right)
= \left( \rho^{\#}_y \right) ^{-1} (w_y) 
= v_{\rho(y)} = v_{\sigma(y)} .
\end{equation}
Since $\tilde{v}_{\sigma(y)}$ is the unique valuation on $\tilde{\cat{O}}_{\sigma(y)}$ such that $\left( \pi^{\#}_{\sigma(y)} \right)^{-1} ( \tilde{v}_{\sigma(y)} ) = v_{\sigma(y)}$, it follows that
\begin{equation}
\left( \sigma^{\#}_y \right)^{-1} (w_y) = \tilde{v}_{\sigma(y)} .
\end{equation}
This shows that $(\sigma, \sigma^{\#}) : (Y, \cat{O}_Y, \cat{V}_Y) \to (\tilde{X}, \tilde{\cat{O}}, \tilde{\cat{V}})$ is a morphism in $\cat{C}_1$.

Since the diagram
\[\begin{tikzcd}
	{(Y, \cat{O}_Y)} & {(\tilde{X}, \tilde{\cat{O}})} \\
	& {(X, \cat{O})}
	\arrow["{(\sigma,\sigma^{\#})}", from=1-1, to=1-2]
	\arrow["{(\rho,\rho^{\#})}"', from=1-1, to=2-2]
	\arrow["{(\pi, \pi^{\#})}", from=1-2, to=2-2]
\end{tikzcd}\]
is commutative in $\cat{D}$, the diagram
\[\begin{tikzcd}
	{(Y, \cat{O}_Y, \cat{V}_Y)} & {(\tilde{X}, \tilde{\cat{O}}, \tilde{\cat{V}})} \\
	& {(X, \cat{O}, \cat{V})}
	\arrow["{(\sigma,\sigma^{\#})}", from=1-1, to=1-2]
	\arrow["{(\rho,\rho^{\#})}"', from=1-1, to=2-2]
	\arrow["{(\pi, \pi^{\#})}", from=1-2, to=2-2]
\end{tikzcd}\]
is commutative in $\cat{C}_1$.

On the other hand, suppose that $(\tau, \tau^{\#}) : (Y, \cat{O}_Y, \cat{V}_Y) \to (\tilde{X}, \tilde{\cat{O}}, \tilde{\cat{V}})$ is another morphism in $\cat{C}_f$ such that the following diagram is commutative.
\[\begin{tikzcd}
	{(Y, \cat{O}_Y, \cat{V}_Y)} & {(\tilde{X}, \tilde{\cat{O}}, \tilde{\cat{V}})} \\
	& {(X, \cat{O}, \cat{V})}
	\arrow["{(\tau,\tau^{\#})}", from=1-1, to=1-2]
	\arrow["{(\rho,\rho^{\#})}"', from=1-1, to=2-2]
	\arrow["{(\pi, \pi^{\#})}", from=1-2, to=2-2]
\end{tikzcd}\]
We prove that $(\tau, \tau^{\#}) = (\sigma, \sigma^{\#})$. The diagram
\[\begin{tikzcd}
	Y & {\tilde{X}} \\
	& X
	\arrow["\tau", from=1-1, to=1-2]
	\arrow["\rho"', from=1-1, to=2-2]
	\arrow["\pi", from=1-2, to=2-2]
\end{tikzcd}\]
is commutative. Since the map $\pi : \tilde{X} \to X$ is the identity map $x \mapsto x$, the map $\tau : Y \to \tilde{X}$ must be equal to the map $y \mapsto \rho(y)$. However, the map $\sigma : Y \to \tilde{X}$ was defined to be the map $y \mapsto \rho(y)$. Therefore $\tau = \sigma$. Then the morphism $\tau^{\#} : \tilde{\cat{O}} \to \tau_* \cat{O}_Y$ of sheaves of condensed rings on $\tilde{X}$ is a morphism $\tilde{\cat{O}} \to \sigma_* \cat{O}_Y$. Moreover, the following diagram is commutative.
\[\begin{tikzcd}
	{\pi_* \sigma_* \cat{O}_Y} & {\pi_* \tilde{\cat{O}}} \\
	& {\cat{O}}
	\arrow["{\pi_* \tau^{\#}}"', from=1-2, to=1-1]
	\arrow["{\rho^{\#}}", from=2-2, to=1-1]
	\arrow["{\pi^{\#}}"', from=2-2, to=1-2]
\end{tikzcd}\]
Then the uniqueness of $\sigma^{\#} : \tilde{\cat{O}} \to \sigma_* \cat{O}_Y$ shows that $\tau^{\#} = \sigma^{\#}$. Thus $(\tau, \tau^{\#}) = (\sigma, \sigma^{\#})$. This completes the proof of \cref{prop:construction of coreflection in C f}.
\end{proof}

\subsection{The category $\cat{C}_c$}

\subsubsection{Definition}

\begin{df}
The category $\cat{C}_c$ is defined to be the full subcategory of $\cat{C}_f$ consisting of all $(X, \cat{O}_X, \cat{V}_X) \in |\cat{C}_f|$ with the following property: For every $x \in X$ and every $f, g \in \cat{O}_{X,x}(*)$ with $|f|_x \leq |g|_x \neq 0$, the condensed $\cat{O}_{X,x}$-algebra $\cat{O}_{X,x}$ is $f/g$-coalescent. In other words, the valued condensed ring $(\cat{O}_{X,x}, |\cdot|_x)$ is an object of $\ub{VCRing}_c$ (\cref{df:coalescent valued ring}) for every $x \in X$.
\end{df}

\subsubsection{Coreflection}

\begin{prop} \label{prop:construction of coreflection in C c}
Let $(X, \cat{O}, \cat{V})$ be an object of $\cat{C}_f$. Then there exists a coreflection
\begin{equation}
\left( (\tilde{X}, \tilde{\cat{O}}, \tilde{\cat{V}}) ,\,  
(\tilde{X}, \tilde{\cat{O}}, \tilde{\cat{V}}) \xto{(\pi, \pi^{\#})} (X, \cat{O}, \cat{V})
\right)
\end{equation}
of $(X, \cat{O}, \cat{V})$ along the inclusion functor $\cat{C}_c \mon \cat{C}_f$ which has the following property.
\begin{enumerate}
\item
The underlying set of $\tilde{X}$ is equal to the underlying set of $X$.

\item
The map $\tilde{X} \xto{\pi} X$ is equal to the identity map $x \mapsto x$.

\item
For every $x \in \tilde{X}$, the homomorphism $\pi^{\#}_x : (\cat{O}_x, |\cdot|_x) \to (\tilde{\cat{O}}_x, |\cdot|_x)$ of valued condensed rings is the coalescence of the valued condensed ring $(\cat{O}_x, |\cdot|_x)$ in the sense of \cref{df:coalescence of valued condensed ring}.
\end{enumerate}
\end{prop}

\subsubsection{Proof of \cref{prop:construction of coreflection in C c}}
\label{sec:Proof of construction of coreflection in C c}

\begin{lem} \label{lem:induction step for construction of coreflection in C c}
Let $(X, \cat{O}, \cat{V})$ be an object of $\cat{C}_f$. Then there exists an object $(X', \cat{O}', \cat{V}')$ of $\cat{C}_f$ and a morphism $(\beta, \beta^{\#}) : (X', \cat{O}', \cat{V}') \to (X, \cat{O}, \cat{V})$ in $\cat{C}_f$ with the following properties.
\begin{enumerate}
\item
The underlying set of $X'$ is equal to the underlying set of $X$.

\item
The map $\beta : X' \to X$ is equal to the identity map $x \mapsto x$.

\item
For every $x \in X$, let us define
\begin{equation}
T_x := \set{(f,g) \in \cat{O}_x(*) \times \cat{O}_x(*)}{|f|_x \leq |g|_x \neq 0}.
\end{equation}
Then the homomorphism $\beta^{\#}_x : \cat{O}_x \to \cat{O}'_x$ of condensed rings is the $T_x$-coalescence of the condensed $\cat{O}_x$-algebra $\cat{O}_x$.

\item
Suppose that $(Y, \cat{O}_Y, \cat{V}_Y)$ is an object of $\cat{C}_c$ and that $(\rho, \rho^{\#}) : (Y, \cat{O}_Y, \cat{V}_Y) \to (X, \cat{O}, \cat{V})$ is a morphism in $\cat{C}_f$. Then there exists a unique morphism $(\rho', {\rho'}^{\#}) : (Y, \cat{O}_Y, \cat{V}_Y) \to (X', \cat{O}', \cat{V}')$ in $\cat{C}_f$ such that the following diagram is commutative.
\[\begin{tikzcd}
	{(Y, \cat{O}_Y, \cat{V}_Y)} & {(X', \cat{O}', \cat{V}')} \\
	& {(X, \cat{O}, \cat{V})}
	\arrow["{(\rho', {\rho'}^{\#})}", from=1-1, to=1-2]
	\arrow["{(\rho,\rho^{\#})}"', from=1-1, to=2-2]
	\arrow["{(\beta, \beta^{\#})}", from=1-2, to=2-2]
\end{tikzcd}\]
\end{enumerate}
\end{lem}

\begin{proof}
For each $x \in X$, let us write $v_x$ for the unique element of $\cat{V}_x$. For each open subset $U$ of $X$, let us define
\begin{equation}
T(U) := \set{(f,g) \in \cat{O}(U)(*) \times \cat{O}(U)(*)}{
|f_x|_{v_x} \leq |g_x|_{v_x} \neq 0 \text{ for every } x \in U }.
\end{equation}
Let $\iota_U : \cat{O}(U) \to \cat{O}(U)_{\approx T(U)}$ be the $T(U)$-coalescence of the condensed $\cat{O}(U)$-algebra $\cat{O}(U)$. If $U,V$ are open subsets of $X$ with $V \sub U$ and if $(f,g) \in T(U)$, then we have $(f|_V, g|_V) \in T(V)$. Therefore $\cat{O}(V)_{\approx T(V)}$ is a $T(U)$-coalescent condensed $\cat{O}(U)$-algebra via $\cat{O}(U) \xto{\cat{O}(V \sub U)} \cat{O}(V) \xto{\iota_V} \cat{O}(V)_{\approx T(V)}$. Then the universality of the $T(U)$-coalescence $\iota_U : \cat{O}(U) \to \cat{O}(U)_{\approx T(U)}$ of the condensed $\cat{O}(U)$-algebra $\cat{O}(U)$ shows that there exists a unique homomorphism $\phi_{U,V} : \cat{O}(U)_{\approx T(U)} \to \cat{O}(V)_{\approx T(V)}$ of condensed rings such that the diagram
\[\begin{tikzcd}
	{\cat{O}(U)} & {\cat{O}(U)_{\approx T(U)}} \\
	{\cat{O}(V)} & {\cat{O}(V)_{\approx T(V)}}
	\arrow["{\iota_U}", from=1-1, to=1-2]
	\arrow["{\cat{O}(V \sub U)}"', from=1-1, to=2-1]
	\arrow["{\phi_{U,V}}", from=1-2, to=2-2]
	\arrow["{\iota_V}"', from=2-1, to=2-2]
\end{tikzcd}\]
is commutative. We define a presheaf $F$ of condensed rings on $X$ by
\begin{align}
F(U) := \cat{O}(U)_{\approx T(U)} \quad (U \sub X \text{ open}) \: ; &&
F(V \sub U) := \phi_{U,V} \quad (V \sub U \sub X \text{ open}) .
\end{align}
Then the family $\iota := (\iota_U)_{U \sub X \text{ open}}$ is a morphism $\cat{O} \to F$ of presheaves of condensed rings on $X$. Let $\eta : F \to G$ be the sheafification of the presheaf $F$ of condensed rings on $X$. Let us write $\theta : \cat{O} \to G$ for the composition $\cat{O} \xto{\iota} F \xto{\eta} G$. Then we obtain an object $(X,G)$ of the category $\cat{D}$, and the pair $(\id_X, \theta)$ is a morphism $(X,G) \to (X,\cat{O})$ in $\cat{D}$.

Let $x \in X$. We claim that the homomorphism $\theta_x : \cat{O}_x \to G_x$ of condensed rings is the $T_x$-coalescence of the condensed $\cat{O}_x$-algebra $\cat{O}_x$. The homomorphism $\theta_x : \cat{O}_x \to G_x$ is equal to the composition $\cat{O}_x \xto{\iota_x} F_x \xto{\eta_x} G_x$. Since $F_x \xto{\eta_x} G_x$ is an isomorphism of condensed rings by \cref{cor:existence of sheafification functor of condensed presheaves}, it suffices to show that the homomorphism $\iota_x : \cat{O}_x \to F_x$ of condensed rings is the $T_x$-coalescence of the condensed $\cat{O}_x$-algebra $\cat{O}_x$. First we show that $F_x$ is a $T_x$-coalescent condensed $\cat{O}_x$-algebra via $\iota_x : \cat{O}_x \to F_x$. Let $(f,g) \in T_x$. Since
\begin{equation}
\cat{O}_x (*) = \underset{x \in U \sub X \text{ open}}{\colim} \; \cat{O}(U)(*)
\end{equation}
by \cref{prop:compatibility of stalk and evaluation at S}, there exist an open neighbourhood $U$ of $x$ in $X$ and an element $\tilde{f}, \tilde{g} \in \cat{O}(U)(*)$ such that $\tilde{f}_x = f$ and $\tilde{g}_x = g$. Since $(X, \cat{O}, \cat{V})$ be an object of $\cat{C}_f$, the set
\begin{equation}
V := \set{x' \in U}{|\tilde{f}_{x'}|_{v_{x'}} \leq |\tilde{g}_{x'}|_{v_{x'}} \neq 0}
\end{equation}
is an open subset of $X$. Moreover, since $(\tilde{f}_x, \tilde{g}_x) = (f,g) \in T_x$, we have $x \in V$. Therefore $V$ is an open neighbourhood of $x$ in $X$. On the other hand, if $W$ is any open neighbourhood of $x$ such that $W \sub V$, then we have $(\tilde{f}|_W, \tilde{g}|_W) \in T(W)$. Therefore $F(W) = \cat{O}(W)_{\approx T(W)}$ is an $(\tilde{f}|_W) / (\tilde{g}|_W)$-coalescent condensed $\cat{O}(W)$-algebra via $\iota_W : \cat{O}(W) \to \cat{O}(W)_{\approx T(W)} = F(W)$. Then $F(W)$ is an $\tilde{f} / \tilde{g}$-coalescent condensed $\cat{O}(U)$-algebra via $\cat{O}(U) \xto{\cat{O}(W \sub U)} \cat{O}(W) \xto{\iota_W} F(W)$. Let us consider $F(W)$ as a condensed $\cat{O}(U)$-algebra via $\cat{O}(U) \xto{\cat{O}(W \sub U)} \cat{O}(W) \xto{\iota_W} F(W)$. Furthermore, let us consider $F_x$ as a condensed $\cat{O}(U)$-algebra via $\cat{O}(U) \xto{\can} \cat{O}_x \xto{\iota_x} F_x$. Then we can write
\begin{equation}
F_x = \underset{x \in W \sub V \text{ open}}{\colim} \, F(W)
\end{equation}
in the category $\ub{CAlg}_{\cat{O}(U)}$. By \cref{cor:coalescent modules are closed under limits and colimits}, we conclude that $F_x$ is an $\tilde{f} / \tilde{g}$-coalescent condensed $\cat{O}(U)$-algebra via $\cat{O}(U) \xto{\can} \cat{O}_x \xto{\iota_x} F_x$. Since $\tilde{f}_x = f$ and $\tilde{g}_x = g$, it follows that $F_x$ is an $f/g$-coalescent condensed $\cat{O}_x$-algebra via $\iota_x : \cat{O}_x \to F_x$. Since this holds for any $(f,g) \in T_x$, we conclude that $F_x$ is a $T_x$-coalescent condensed $\cat{O}_x$-algebra via $\iota_x : \cat{O}_x \to F_x$.

Next suppose that $A$ is any $T_x$-coalescent condensed $\cat{O}_x$-algebra and that $\psi : \cat{O}_x \to A$ is any homomorphism of condensed $\cat{O}_x$-algebras. We prove that there exists a unique homomorphism $\tilde{\psi} : F_x \to A$ of condensed $\cat{O}_x$-algebras such that the diagram
\[\begin{tikzcd}
	{\cat{O}_x} & A \\
	{F_x}
	\arrow["\psi", from=1-1, to=1-2]
	\arrow["{\iota_x}"', from=1-1, to=2-1]
	\arrow["{\tilde{\psi}}"', from=2-1, to=1-2]
\end{tikzcd}\]
is commutative. Let $U$ be any open neighbourhood of $x$ in $X$. For every $(f,g) \in T(U)$, we have $(f_x,g_x) \in T_x$ by definition. Therefore $A$ is a $T(U)$-coalescent condensed $\cat{O}(U)$-algebra via $\cat{O}(U) \xto{\can} \cat{O}_x \xto{\psi} A$. Then the universality of the $T(U)$-coalescence $\iota_U : \cat{O}(U) \to \cat{O}(U)_{\approx T(U)}$ of the condensed $\cat{O}(U)$-algebra $\cat{O}(U)$ shows that there exists a unique homomorphism $\psi_U : \cat{O}(U)_{\approx T(U)} \to A$ of condensed rings such that the diagram
\[\begin{tikzcd}
	{\cat{O}(U)} & {\cat{O}(U)_{\approx T(U)}} \\
	{\cat{O}_x} & A
	\arrow["{\iota_U}", from=1-1, to=1-2]
	\arrow["\can"', from=1-1, to=2-1]
	\arrow["{\psi_U}", from=1-2, to=2-2]
	\arrow["\psi"', from=2-1, to=2-2]
\end{tikzcd}\]
is commutative. If $U,V$ are open neighbourhoods of $x$ in $X$ such that $V \sub U$, then the commutativity of the diagram
\[\begin{tikzcd}
	{\cat{O}(U)} & {\cat{O}(U)_{\approx T(U)}} \\
	{\cat{O}(V)} & {\cat{O}(V)_{\approx T(V)}} \\
	{\cat{O}_x} & A
	\arrow["{\iota_U}", from=1-1, to=1-2]
	\arrow["{\cat{O}(U \sub V)}", from=1-1, to=2-1]
	\arrow["\can"', curve={height=24pt}, from=1-1, to=3-1]
	\arrow["{\phi_{U,V}}", from=1-2, to=2-2]
	\arrow["{\iota_V}"', from=2-1, to=2-2]
	\arrow["\can", from=2-1, to=3-1]
	\arrow["{\psi_V}", from=2-2, to=3-2]
	\arrow["\psi"', from=3-1, to=3-2]
\end{tikzcd}\]
shows that $\psi_V \of \phi_{U,V} = \psi_U$. Therefore there exists a unique homomorphism $\tilde{\psi} : F_x \to A$ of condensed rings such that the diagram
\[\begin{tikzcd}
	{F(U)} & A \\
	{F_x}
	\arrow["{\psi_U}", from=1-1, to=1-2]
	\arrow["\can"', from=1-1, to=2-1]
	\arrow["{\tilde{\psi}}"', from=2-1, to=1-2]
\end{tikzcd}\]
is commutative for every open neighbourhood $U$ of $x$ in $X$. Then the diagram
\[\begin{tikzcd}
	{\cat{O}_x} \\
	{\cat{O}(U)} & {F(U)} & A \\
	{\cat{O}_x} & {F_x}
	\arrow["\psi", curve={height=-12pt}, from=1-1, to=2-3]
	\arrow["\can", from=2-1, to=1-1]
	\arrow["{\iota_U}", from=2-1, to=2-2]
	\arrow["\can"', from=2-1, to=3-1]
	\arrow["{\psi_U}", from=2-2, to=2-3]
	\arrow["\can"', from=2-2, to=3-2]
	\arrow["{\iota_x}", from=3-1, to=3-2]
	\arrow["{\tilde{\psi}}"', from=3-2, to=2-3]
\end{tikzcd}\]
is commutative for every open neighbourhood $U$ of $x$ in $X$. It follows that the diagram
\[\begin{tikzcd}
	{\cat{O}_x} & A \\
	{F_x}
	\arrow["\psi", from=1-1, to=1-2]
	\arrow["{\iota_x}"', from=1-1, to=2-1]
	\arrow["{\tilde{\psi}}"', from=2-1, to=1-2]
\end{tikzcd}\]
is commutative. Then $\tilde{\psi} : F_x \to A$ is a homomorphism of condensed $\cat{O}_x$-algebras. On the other hand, suppose that $\tilde{\psi}' : F_x \to A$ is another homomorphism of condensed $\cat{O}_x$-algebras such that the diagram
\[\begin{tikzcd}
	{\cat{O}_x} & A \\
	{F_x}
	\arrow["\psi", from=1-1, to=1-2]
	\arrow["{\iota_x}"', from=1-1, to=2-1]
	\arrow["{\tilde{\psi}'}"', from=2-1, to=1-2]
\end{tikzcd}\]
is commutative. Then, for every open neighbourhood $U$ of $x$ in $X$, the diagram
\[\begin{tikzcd}
	{\cat{O}(U)} & {F(U)} \\
	& {F_x} \\
	{\cat{O}_x} & A
	\arrow["{\iota_U}", from=1-1, to=1-2]
	\arrow["\can"', from=1-1, to=3-1]
	\arrow["\can", from=1-2, to=2-2]
	\arrow["{\tilde{\psi}'}", from=2-2, to=3-2]
	\arrow["{\iota_x}", from=3-1, to=2-2]
	\arrow["\psi"', from=3-1, to=3-2]
\end{tikzcd}\]
is commutative. Then the definition of $\psi_U$ shows that $F(U) \xto{\can} F_x \xto{\tilde{\psi}'} A$ is equal to $\psi_U$. In other words, the diagram
\[\begin{tikzcd}
	{F(U)} & A \\
	{F_x}
	\arrow["{\psi_U}", from=1-1, to=1-2]
	\arrow["\can"', from=1-1, to=2-1]
	\arrow["{\tilde{\psi}'}"', from=2-1, to=1-2]
\end{tikzcd}\]
is commutative. Since this holds for every open neighbourhood $U$ of $x$ in $X$, the definition of $\tilde{\psi}$ shows that $\tilde{\psi}' = \tilde{\psi}$. This completes the proof that the homomorphism $\theta_x : \cat{O}_x \to G_x$ of condensed rings is the $T_x$-coalescence of the condensed $\cat{O}_x$-algebra $\cat{O}_x$.

Let $x \in X$. Since the homomorphism $\theta_x : \cat{O}_x \to G_x$ of condensed rings is the $T_x$-coalescence of the condensed $\cat{O}_x$-algebra $\cat{O}_x$, \cref{lem:extension of continuous valuation along coalescence} shows that there exists a unique continuous valuation $u_x$ on $G_x$ such that $(\theta_x)^{-1}(u_x) = v_x$. Let us define $\cat{U}_x := \{u_x\}$. Then we obtain a family $\cat{U} := (\cat{U}_x)_{x \in X}$. Then the triple $(X, G, \cat{U})$ is an object of $\cat{C}_1$ and $(\id_X, \theta)$ is a morphism $(X, G, \cat{U}) \to (X, \cat{O}, \cat{V})$ in $\cat{C}_1$. By \cref{prop:construction of coreflection in C f}, there exists a coreflection
\begin{equation}
\left( (X', \cat{O}', \cat{V}') ,\,  
(X', \cat{O}', \cat{V}') \xto{(\pi, \pi^{\#})} (X, G, \cat{U})
\right)
\end{equation}
of the object $(X, G, \cat{U}) \in |\cat{C}_1|$ along the inclusion functor $\cat{C}_f \mon \cat{C}_1$ which has the following property.
\begin{enumerate}
\item
The underlying set of $X'$ is equal to the underlying set of $X$.

\item
The map $X' \xto{\pi} X$ is equal to the identity map $x \mapsto x$.

\item
For every $x \in X'$, the homomorphism $\pi^{\#}_x : G_x \to \cat{O}'_x$ is an isomorphism of condensed rings.
\end{enumerate}

We define $(\beta, \beta^{\#}) : (X', \cat{O}', \cat{V}') \to (X, \cat{O}, \cat{V})$ to be the composition 
\begin{equation}
(X', \cat{O}', \cat{V}') \xto{(\pi, \pi^{\#})} (X, G, \cat{U}) \xto{(\id_X, \theta)} (X, \cat{O}, \cat{V})
\end{equation}
in $\cat{C}_1$. We prove that this object $(X', \cat{O}', \cat{V}')$ of $\cat{C}_f$ and the morphism $(\beta, \beta^{\#}) : (X', \cat{O}', \cat{V}') \to (X, \cat{O}, \cat{V})$ in $\cat{C}_f$ satisfy all of the conditions stated in this proposition.

First of all, the underlying set of $X'$ is equal to the underlying set of $X$ by the construction of $(X', \cat{O}', \cat{V}')$. Furthermore, the map $\beta : X' \to X$ is equal to the composition $X' \xto{\pi} X \xto{\id_X} X$. Since the map $X' \xto{\pi} X$ is equal to the identity map $x \mapsto x$ by the construction of $(X', \cat{O}', \cat{V}')$, it follows that the map $\beta : X' \to X$ is equal to the identity map $x \mapsto x$.

For each $x \in X$, the homomorphism $\beta^{\#}_x : \cat{O}_x \to \cat{O}'_x$ is equal to the composition $\cat{O}_x \xto{\theta_x} G_x \xto{\pi^{\#}_x} \cat{O}'_x$. We have already proved that the homomorphism $\theta_x : \cat{O}_x \to G_x$ of condensed rings is the $T_x$-coalescence of the condensed $\cat{O}_x$-algebra $\cat{O}_x$. Moreover, the homomorphism $\pi^{\#}_x : G_x \to \cat{O}'_x$ is an isomorphism of condensed rings by the construction of $(X', \cat{O}', \cat{V}')$. It follows that the homomorphism $\beta^{\#}_x : \cat{O}_x \to \cat{O}'_x$ is the $T_x$-coalescence of the condensed $\cat{O}_x$-algebra $\cat{O}_x$.

Next suppose that $(Y, \cat{O}_Y, \cat{V}_Y)$ is an object of $\cat{C}_c$ and that $(\rho, \rho^{\#}) : (Y, \cat{O}_Y, \cat{V}_Y) \to (X, \cat{O}, \cat{V})$ is a morphism in $\cat{C}_f$. We prove that there exists a unique morphism $(\rho', {\rho'}^{\#}) : (Y, \cat{O}_Y, \cat{V}_Y) \to (X', \cat{O}', \cat{V}')$ in $\cat{C}_f$ such that the following diagram is commutative.
\[\begin{tikzcd}
	{(Y, \cat{O}_Y, \cat{V}_Y)} & {(X', \cat{O}', \cat{V}')} \\
	& {(X, \cat{O}, \cat{V})}
	\arrow["{(\rho', {\rho'}^{\#})}", from=1-1, to=1-2]
	\arrow["{(\rho,\rho^{\#})}"', from=1-1, to=2-2]
	\arrow["{(\beta, \beta^{\#})}", from=1-2, to=2-2]
\end{tikzcd}\]

For each $y \in Y$, let us write $w_y$ for the unique element of $\cat{V}_{Y,y}$. Let $U$ be any open subset of $X$. Let $(f,g) \in T(U)$. Then for each $y \in \rho^{-1}(U)$, we have
\begin{equation}
|f_{\rho(y)}|_{v_{\rho(y)}} \leq |g_{\rho(y)}|_{v_{\rho(y)}} \neq 0 .
\end{equation}
On the other hand, the diagram
\[\begin{tikzcd}
	{\cat{O}(U)} & {\cat{O}_Y(\rho^{-1}(U))} \\
	{\cat{O}_{\rho(y)}} & {\cat{O}_{Y,y}}
	\arrow["{\rho^{\#}_U}", from=1-1, to=1-2]
	\arrow["\can"', from=1-1, to=2-1]
	\arrow["\can", from=1-2, to=2-2]
	\arrow["{\rho^{\#}_y}"', from=2-1, to=2-2]
\end{tikzcd}\]
is commutative. Moreover, we have $( \rho^{\#}_y )^{-1}(w_y) = v_{\rho(y)}$ since $(\rho, \rho^{\#}) : (Y, \cat{O}_Y, \cat{V}_Y) \to (X, \cat{O}, \cat{V})$ is a morphism in $\cat{C}_f$. Therefore 
\begin{equation}
\big| \big( (\rho^{\#}_U)_*(f) \big)_y \big|_{w_y} \leq
\big| \big( (\rho^{\#}_U)_*(g) \big)_y \big|_{w_y} \neq 0 .
\end{equation}
Since $(Y, \cat{O}_Y, \cat{V}_Y)$ is an object of $\cat{C}_c$, it follows that the condensed $\cat{O}_{Y,y}$-algebra $\cat{O}_{Y,y}$ is $\\ \big( (\rho^{\#}_U)_*(f) \big)_y \big/ \big( (\rho^{\#}_U)_*(g) \big)_y$-coalescent. Therefore $\cat{O}_{Y,y}$ is a $\big( (\rho^{\#}_U)_*(f) \big) / \big( (\rho^{\#}_U)_*(g) \big)$-coalescent condensed $\cat{O}_Y(\rho^{-1}(U))$-algebra via $\cat{O}_Y(\rho^{-1}(U)) \xto{\can} \cat{O}_{Y,y}$. Since this holds for every $y \in \rho^{-1}(U)$, \cref{prop:testing coalescence by stalks} shows that the $\cat{O}_Y(\rho^{-1}(U))$-algebra $\cat{O}_Y(\rho^{-1}(U))$ is $\big( (\rho^{\#}_U)_*(f) \big) / \big( (\rho^{\#}_U)_*(g) \big)$-coalescent. Therefore $\cat{O}_Y(\rho^{-1}(U))$ is an $f/g$-coalescent condensed $\cat{O}(U)$-algebra via $\rho^{\#}_U : \cat{O}(U) \to \cat{O}_Y(\rho^{-1}(U))$. Since this holds for every $(f,g) \in T(U)$, we conclude that $\cat{O}_Y(\rho^{-1}(U))$ is a $T(U)$-coalescent condensed $\cat{O}(U)$-algebra via $\rho^{\#}_U : \cat{O}(U) \to \cat{O}_Y(\rho^{-1}(U))$. Then the universality of the $T(U)$-coalescence $\iota_U : \cat{O}(U) \to \cat{O}(U)_{\approx T(U)}$ of the condensed $\cat{O}(U)$-algebra $\cat{O}(U)$ shows that there exists a unique homomorphism $\tilde{\rho}^{\#}_U : \cat{O}(U)_{\approx T(U)} \to \cat{O}_Y(\rho^{-1}(U))$ of condensed rings such that the diagram
\[\begin{tikzcd}
	{\cat{O}_Y(\rho^{-1}(U))} & {\cat{O}(U)_{\approx T(U)}} \\
	& {\cat{O}(U)}
	\arrow["{\tilde{\rho}^{\#}_U}"', from=1-2, to=1-1]
	\arrow["{\rho^{\#}_U}", from=2-2, to=1-1]
	\arrow["{\iota_U}"', from=2-2, to=1-2]
\end{tikzcd}\]
is commutative. If $U,V$ are open subsets of $X$ such that $V \sub U$, then tha following diagram is commutative.
\[\begin{tikzcd}
	& {\cat{O}(U)_{\approx T(U)}} \\
	{\cat{O}_Y(\rho^{-1}(U))} & {\cat{O}(U)} & {\cat{O}(U)_{\approx T(U)}} \\
	{\cat{O}_Y(\rho^{-1}(V))} & {\cat{O}(V)} \\
	& {\cat{O}(V)_{\approx T(V)}}
	\arrow["{\tilde{\rho}^{\#}_U}"', from=1-2, to=2-1]
	\arrow["{\cat{O}_Y(\rho^{-1}(V) \sub \rho^{-1}(U))}"', from=2-1, to=3-1]
	\arrow["{\iota_U}"', from=2-2, to=1-2]
	\arrow["{\rho^{\#}_U}", from=2-2, to=2-1]
	\arrow["{\iota_U}", from=2-2, to=2-3]
	\arrow["{\cat{O}(V \sub U)}", from=2-2, to=3-2]
	\arrow["{\phi_{U,V}}", curve={height=-24pt}, from=2-3, to=4-2]
	\arrow["{\rho^{\#}_V}"', from=3-2, to=3-1]
	\arrow["{\iota_V}", from=3-2, to=4-2]
	\arrow["{\tilde{\rho}^{\#}_V}", from=4-2, to=3-1]
\end{tikzcd}\]
For every $(f,g) \in T(U)$, we have $(f|_V, g|_V) \in T(V)$. Since $\cat{O}_Y(\rho^{-1}(V))$ is a $T(V)$-coalescent condensed $\cat{O}(V)$-algebra via $\rho^{\#}_V : \cat{O}(V) \to \cat{O}_Y(\rho^{-1}(V))$, it follows that $\cat{O}_Y(\rho^{-1}(V))$ is a $T(U)$-coalescent condensed $\cat{O}(U)$-algebra via $\cat{O}(U) \xto{\cat{O}(V \sub U)} \cat{O}(V) \xto{\rho^{\#}_V} \cat{O}_Y(\rho^{-1}(V))$. Then the universality of the $T(U)$-coalescence $\iota_U : \cat{O}(U) \to \cat{O}(U)_{\approx T(U)}$ of the condensed $\cat{O}(U)$-algebra $\cat{O}(U)$ shows that the diagram
\[\begin{tikzcd}
	{\cat{O}_Y(\rho^{-1}(U))} & {\cat{O}(U)_{\approx T(U)}} \\
	{\cat{O}_Y(\rho^{-1}(V))} & {\cat{O}(V)_{\approx T(V)}}
	\arrow["{\cat{O}_Y(\rho^{-1}(V) \sub \rho^{-1}(U))}"', from=1-1, to=2-1]
	\arrow["{\tilde{\rho}^{\#}_U}"', from=1-2, to=1-1]
	\arrow["{\phi_{U,V}}", from=1-2, to=2-2]
	\arrow["{\tilde{\rho}^{\#}_V}", from=2-2, to=2-1]
\end{tikzcd}\]
is commutative. Therefore the family $\tilde{\rho}^{\#} := (\tilde{\rho}^{\#}_U)_{U \sub X \text{ open}}$ is a morphism $F \to \rho_* \cat{O}_Y$ of presheaves of condensed rings on $X$. Note that the diagram
\[\begin{tikzcd}
	{\rho_* \cat{O}_Y} & F \\
	& {\cat{O}}
	\arrow["{\tilde{\rho}^{\#}}"', from=1-2, to=1-1]
	\arrow["{\rho^{\#}}", from=2-2, to=1-1]
	\arrow["\iota"', from=2-2, to=1-2]
\end{tikzcd}\]
is commutative by definition. Since $\rho_* \cat{O}_Y$ is a sheaf of condensed rings on $X$, the universality of the sheafification $\eta : F \to G$ of the presheaf $F$ of condensd rings on $X$ shows that there exists a unique morphism $\hat{\rho}^{\#} : G \to \rho_* \cat{O}_Y$ of sheaves of condensed rings on $X$ such that the diagram
\[\begin{tikzcd}
	{\rho_* \cat{O}_Y} & G \\
	& F
	\arrow["{\hat{\rho}^{\#}}"', from=1-2, to=1-1]
	\arrow["{\tilde{\rho}^{\#}}", from=2-2, to=1-1]
	\arrow["\eta"', from=2-2, to=1-2]
\end{tikzcd}\]
is commutative. Then we have a morphism $(\rho, \hat{\rho}^{\#}) : (Y, \cat{O}_Y) \to (X,G)$ in the category $\cat{D}$, and the commutativity of the diagrams
\[\begin{tikzcd}
	Y & X & {\rho_* \cat{O}_Y} & G \\
	&&& F \\
	& X && {\cat{O}}
	\arrow["\rho", from=1-1, to=1-2]
	\arrow["\rho"', from=1-1, to=3-2]
	\arrow["{\id_X}", from=1-2, to=3-2]
	\arrow["{\hat{\rho}^{\#}}"', from=1-4, to=1-3]
	\arrow["{\tilde{\rho}^{\#}}", from=2-4, to=1-3]
	\arrow["\eta"', from=2-4, to=1-4]
	\arrow["{\rho^{\#}}", curve={height=-12pt}, from=3-4, to=1-3]
	\arrow["\theta"', curve={height=18pt}, from=3-4, to=1-4]
	\arrow["\iota"', from=3-4, to=2-4]
\end{tikzcd}\]
shows that the the diagram
\[\begin{tikzcd}
	{(Y,\cat{O}_Y)} & {(X,G)} \\
	& {(X,\cat{O})}
	\arrow["{(\rho,\hat{\rho}^{\#})}", from=1-1, to=1-2]
	\arrow["{(\rho,\rho^{\#})}"', from=1-1, to=2-2]
	\arrow["{(\id_X,\theta)}", from=1-2, to=2-2]
\end{tikzcd}\]
is commutative in $\cat{D}$. Then for each $y \in Y$, the following diagram is commutative.
\[\begin{tikzcd}
	{\cat{O}_{Y,y}} & {G_{\rho(y)}} \\
	& {\cat{O}_{\rho(y)}}
	\arrow["{\hat{\rho}^{\#}_y}"', from=1-2, to=1-1]
	\arrow["{\rho^{\#}_y}", from=2-2, to=1-1]
	\arrow["{\theta_{\rho(y)}}"', from=2-2, to=1-2]
\end{tikzcd}\]
Therefore
\begin{equation}
(\theta_{\rho(y)})^{-1} \left( ( \hat{\rho}^{\#}_y )^{-1}(w_y) \right) =
(\rho^{\#}_y)^{-1}(w_y).
\end{equation}
On the other hand, since $(\rho, \rho^{\#}) : (Y, \cat{O}_Y, \cat{V}_Y) \to (X, \cat{O}, \cat{V})$ is a morphism in $\cat{C}_f$, we have
\begin{equation}
(\rho^{\#}_y)^{-1}(w_y) = v_{\rho(y)} .
\end{equation}
Consequently, we have
\begin{equation}
(\theta_{\rho(y)})^{-1} \left( ( \hat{\rho}^{\#}_y )^{-1}(w_y) \right) =
(\rho^{\#}_y)^{-1}(w_y) = v_{\rho(y)} .
\end{equation}
Since $u_{\rho(y)}$ is a unique continuous valuation on $G_{\rho(y)}$ such that $(\theta_{\rho(y)})^{-1}(u_{\rho(y)}) = v_{\rho(y)}$, we conclude that
\begin{equation}
( \hat{\rho}^{\#}_y )^{-1}(w_y) = u_{\rho(y)} .
\end{equation}
Since this holds for every $y \in Y$, we conclude that $(\rho, \hat{\rho}^{\#})$ is a morphism $(Y, \cat{O}_Y, \cat{V}_Y) \to (X,G, \cat{U})$ in the category $\cat{C}_1$. Note that the diagram
\[\begin{tikzcd}
	{(Y,\cat{O}_Y,\cat{V}_Y)} & {(X,G,\cat{U})} \\
	& {(X,\cat{O},\cat{V})}
	\arrow["{(\rho,\hat{\rho}^{\#})}", from=1-1, to=1-2]
	\arrow["{(\rho,\rho^{\#})}"', from=1-1, to=2-2]
	\arrow["{(\id_X,\theta)}", from=1-2, to=2-2]
\end{tikzcd}\]
is commutative in $\cat{C}_1$. Since $(Y, \cat{O}_Y, \cat{V}_Y)$ is an object of $\cat{C}_f$ and since
\begin{equation}
\left( (X', \cat{O}', \cat{V}') ,\,  
(X', \cat{O}', \cat{V}') \xto{(\pi, \pi^{\#})} (X, G, \cat{U})
\right)
\end{equation}
is the coreflection of the object $(X, G, \cat{U}) \in |\cat{C}_1|$ along the inclusion functor $\cat{C}_f \mon \cat{C}_1$, it follows that there exists a unique morphism $(\rho', {\rho'}^{\#}) : (Y, \cat{O}_Y, \cat{V}_Y) \to (X', \cat{O}', \cat{V}')$ in $\cat{C}_f$ such that the following diagram is commutative.
\[\begin{tikzcd}
	{(Y,\cat{O}_Y,\cat{V}_Y)} & {(X', \cat{O}', \cat{V}')} \\
	& {(X,G,\cat{U})}
	\arrow["{(\rho', {\rho'}^{\#})}", from=1-1, to=1-2]
	\arrow["{(\rho,\hat{\rho}^{\#})}"', from=1-1, to=2-2]
	\arrow["{(\pi, \pi^{\#})}", from=1-2, to=2-2]
\end{tikzcd}\]
Then the following diagram is commutative.
\[\begin{tikzcd}
	{(Y,\cat{O}_Y,\cat{V}_Y)} & {(X', \cat{O}', \cat{V}')} \\
	& {(X,G,\cat{U})} \\
	& {(X,\cat{O},\cat{V})}
	\arrow["{(\rho', {\rho'}^{\#})}", from=1-1, to=1-2]
	\arrow["{(\rho,\hat{\rho}^{\#})}"', from=1-1, to=2-2]
	\arrow["{(\rho,\rho^{\#})}"', curve={height=24pt}, from=1-1, to=3-2]
	\arrow["{(\pi, \pi^{\#})}"', from=1-2, to=2-2]
	\arrow["{(\beta, \beta^{\#})}", curve={height=-30pt}, from=1-2, to=3-2]
	\arrow["{(\id_X,\theta)}"', from=2-2, to=3-2]
\end{tikzcd}\]

On the other hand, suppose that $(\rho'', {\rho''}^{\#}) : (Y, \cat{O}_Y, \cat{V}_Y) \to (X', \cat{O}', \cat{V}')$ is another morphism in $\cat{C}_f$ such that the following diagram is commutative.
\[\begin{tikzcd}
	{(Y, \cat{O}_Y, \cat{V}_Y)} & {(X', \cat{O}', \cat{V}')} \\
	& {(X, \cat{O}, \cat{V})}
	\arrow["{(\rho'', {\rho''}^{\#})}", from=1-1, to=1-2]
	\arrow["{(\rho,\rho^{\#})}"', from=1-1, to=2-2]
	\arrow["{(\beta, \beta^{\#})}", from=1-2, to=2-2]
\end{tikzcd}\]
We prove that $(\rho'', {\rho''}^{\#}) = (\rho', {\rho'}^{\#})$. First of all, the diagrams
\[\begin{tikzcd}
	Y & {X'} & Y & {X'} \\
	& X && X
	\arrow["{\rho''}", from=1-1, to=1-2]
	\arrow["\rho"', from=1-1, to=2-2]
	\arrow["\beta", from=1-2, to=2-2]
	\arrow["{\rho'}", from=1-3, to=1-4]
	\arrow["\rho"', from=1-3, to=2-4]
	\arrow["\beta", from=1-4, to=2-4]
\end{tikzcd}\]
are commutative and $\beta : X' \to X$ is equal to the identity map $x \mapsto x$. Therefore $\rho'' = \rho'$. Moreover, for each open subset $U$ of $X$, the following diagram is commutative.
\[\begin{tikzcd}
	{\cat{O}_Y(\rho^{-1}(U))} & {\cat{O}'(\beta^{-1}(U))} & {G(U)} \\
	& {\cat{O}(U)} & {F(U)}
	\arrow["{\rho''^{\#}_{\beta^{-1}(U)}}"', from=1-2, to=1-1]
	\arrow["{\pi^{\#}_U}"', from=1-3, to=1-2]
	\arrow["{\rho^{\#}_U}", from=2-2, to=1-1]
	\arrow["{\beta^{\#}_U}"', from=2-2, to=1-2]
	\arrow["{\theta_U}"', from=2-2, to=1-3]
	\arrow["{\iota_U}"', from=2-2, to=2-3]
	\arrow["{\eta_U}"', from=2-3, to=1-3]
\end{tikzcd}\]
Then the defintion of $\tilde{\rho}^{\#}_U$ shows that $\rho''^{\#}_{\beta^{-1}(U)} \of \pi^{\#}_U \of \eta_U = \tilde{\rho}^{\#}_U$. Since this holds for every open subset $U$ of $X$, we conclude that the diagram
\[\begin{tikzcd}
	{\rho_* \cat{O}_Y} & {\beta_* \cat{O}'} & G \\
	&& F
	\arrow["{\beta_* \rho''^{\#}}"', from=1-2, to=1-1]
	\arrow["{\pi^{\#}}"', from=1-3, to=1-2]
	\arrow["{\tilde{\rho}^{\#}_U}", from=2-3, to=1-1]
	\arrow["\eta"', from=2-3, to=1-3]
\end{tikzcd}\]
is commutative. Then the definition of $\hat{\rho}^{\#}$ shows that $\beta_* \rho''^{\#} \of \pi^{\#} = \hat{\rho}^{\#}$. On the other hand, note that $\beta = \id_X \of \pi = \pi$ and therefore $\rho'' \of \pi = \rho'' \of \beta = \rho$. Therefore the diagram
\[\begin{tikzcd}
	{(Y,\cat{O}_Y,\cat{V}_Y)} & {(X', \cat{O}', \cat{V}')} \\
	& {(X,G,\cat{U})}
	\arrow["{(\rho'', {\rho''}^{\#})}", from=1-1, to=1-2]
	\arrow["{(\rho,\hat{\rho}^{\#})}"', from=1-1, to=2-2]
	\arrow["{(\pi, \pi^{\#})}", from=1-2, to=2-2]
\end{tikzcd}\]
is commutative. Then the definition of $(\rho', {\rho'}^{\#})$ shows that $(\rho'', {\rho''}^{\#}) = (\rho', {\rho'}^{\#})$. This completes the proof.
\end{proof}

\begin{proof}[Proof of \cref{prop:construction of coreflection in C c}]
Let $(X, \cat{O}, \cat{V})$ be an object of $\cat{C}_f$. For each $x \in X$, let us write $v_x$ for the valuation on the condensed ring $\cat{O}_x$ which is the unique element of $\cat{V}_x$. Let $\N$ be the set of natural numbers (including 0), ordered by the usual order. We consider $\N$ as a small category. We define a functor $D : \N^{\op} \to \cat{C}_f$ inductively as follows.
\begin{itemize}
\item
We define $D(0) := (X, \cat{O}, \cat{V})$.

\item
Let $n \in \N$ with $n \geq 1$. Suppose we have defined $D(n-1) \in |\cat{C}_f|$. Then we define $D(n)$ and $D(n-1 \leq n) : D(n) \to D(n-1)$ to be the object and morphism in $\cat{C}_f$, respectively, which are obtained by applying \cref{lem:induction step for construction of coreflection in C c} to $D(n-1) \in |\cat{C}_f|$.
\end{itemize}
For each $n \in \N$, let us write $D(n) = (X_n, \cat{O}_n, \cat{V}_n)$. For each $m,n \in \N$ with $m \leq n$, let us write $D(m \leq n) = (f_{n,m}, f^{\#}_{n,m}) : D(n) \to D(m)$. For each $n \in \N$ and each $x \in X_n$, let us write $v_{n,x}$ for the valuation on the condensed ring $\cat{O}_{n,x}$ which is the unique element of $\cat{V}_{n,x}$.

For each $n \in \N$ with $n \geq 1$, \cref{lem:induction step for construction of coreflection in C c} shows that the underlying set of $X_{n-1}$ is equal to the underlying set of $X_n$ and that the map $f_{n,n-1} : X_n \to X_{n-1}$ is equal to the identity map $x \mapsto x$. By induction, we conclude that for every $n \in \N$, the underlying set of $X_n$ is equal to the underlying set of $X$ and that for every $m,n \in \N$ with $m \leq n$, the map $f_{n,m} : X_n \to X_m$ is equal to the identity map $\id_X : X \to X$. For each $n \in \N$, let $\pi_n : X \to X_n$ be the identity map $\id_X : X \to X$. Let $\tilde{X}$ be the set $X$ endowed with the initial topology induced by the family of maps $\pi_n : X \to X_n$ for $n \in \N$. Then
\begin{equation}
(\tilde{X}, (\tilde{X} \xto{\pi_n} X_n))_{n \in |\N^{\op}|}
\end{equation}
is the limit of the functor 
\[\begin{tikzcd}
	{\N^{\op}} & {\cat{C}_f} & {\cat{C}} && {\cat{D}} && {\ub{Top}.}
	\arrow["D", from=1-1, to=1-2]
	\arrow["\inc", hook, from=1-2, to=1-3]
	\arrow["{\text{forgetful functor}}", from=1-3, to=1-5]
	\arrow["{\text{forgetful functor}}", from=1-5, to=1-7]
\end{tikzcd}\]
Using the construction described in \cref{prop:construction of limits in the category D}, we obtain a limit
\begin{equation}
\left( (\tilde{X}, \tilde{\cat{O}}),
\left( (\tilde{X}, \tilde{\cat{O}}) \xto{(\pi_n, \pi^{\#}_n)} (X_n, \cat{O}_n) \right)_{n \in |\N^{\op}|} \right)
\end{equation}
of the functor 
\[\begin{tikzcd}
	{\N^{\op}} & {\cat{C}_f} & {\cat{C}} && {\cat{D}.}
	\arrow["D", from=1-1, to=1-2]
	\arrow["\inc", hook, from=1-2, to=1-3]
	\arrow["{\text{forgetful functor}}", from=1-3, to=1-5]
\end{tikzcd}\]
Using the construction described in \cref{prop:construction of limits in the category C}, we obtain a limit
\begin{equation}
\left( (\tilde{X}, \tilde{\cat{O}}, \tilde{\cat{V}}),
\left( (\tilde{X}, \tilde{\cat{O}}, \tilde{\cat{V}}) \xto{(\pi_n, \pi^{\#}_n)} (X_n, \cat{O}_n, \cat{V}_n) \right)_{n \in |\N^{\op}|} \right)
\end{equation}
of the functor
\[\begin{tikzcd}
	{\N^{\op}} & {\cat{C}_f} & {\cat{C}.}
	\arrow["D", from=1-1, to=1-2]
	\arrow["\inc", hook, from=1-2, to=1-3]
\end{tikzcd}\]
We claim that 
\begin{equation}
\left( (\tilde{X}, \tilde{\cat{O}}, \tilde{\cat{V}}),
\left( (\tilde{X}, \tilde{\cat{O}}, \tilde{\cat{V}}) \xto{(\pi_0, \pi^{\#}_0)} (X, \cat{O}, \cat{V}) \right) \right)
\end{equation}
is a corefletion of $(X, \cat{O}, \cat{V}) \in |\cat{C}_f|$ along the inclusion functor $\cat{C}_c \mon \cat{C}_f$ which has all the properties described in \cref{prop:construction of coreflection in C c}.

First of all, the defintion shows that the underlying set of $\tilde{X}$ is equal to the underlying set of $X$ and the map $\pi_0 : \tilde{X} \to X$ is equal to the identity map $\id_X : X \to X$.

Next let $x \in \tilde{X}$. By (4) of \cref{prop:construction of limits in the category D}, the pair
\begin{equation}
\left( \tilde{\cat{O}}_x , \left( \cat{O}_{n,x} \xto{\pi^{\#}_{n,x}} \tilde{\cat{O}}_x
\right)_{n \in |\N|} \right)
\end{equation}
is the colimit of the functor
\[\begin{tikzcd}[row sep=tiny]
	\N & {\ub{CRing}} \\
	n & {\cat{O}_{n,x}} & {(\text{on objects})} \\
	{(m \leq n)} & {\left( \cat{O}_{m,x} \xto{f^{\#}_{n,m,x}} \cat{O}_{n,x} \right)} & {(\text{on morphisms}) .}
	\arrow["{D_x}", from=1-1, to=1-2]
	\arrow[maps to, from=2-1, to=2-2]
	\arrow[maps to, from=3-1, to=3-2]
\end{tikzcd}\]
This functor is equal to the composition of the functor
\[\begin{tikzcd}[row sep=tiny]
	\N & {\ub{VCRing}} \\
	n & {(\cat{O}_{n,x},v_{n,x})} & {(\text{on objects})} \\
	{(m \leq n)} & {\left( (\cat{O}_{m,x},v_{m,x}) \xto{f^{\#}_{n,m,x}} (\cat{O}_{n,x},v_{n,x}) \right)} & {(\text{on morphisms})}
	\arrow["{\tilde{D}_x}", from=1-1, to=1-2]
	\arrow[maps to, from=2-1, to=2-2]
	\arrow[maps to, from=3-1, to=3-2]
\end{tikzcd}\]
and the forgetful functor $\ub{VCRing} \to \ub{CRing}$. Furthermore, for every $n \in \N$ with $n \geq 1$, \cref{lem:induction step for construction of coreflection in C c} shows that the homomorphism
\begin{equation}
\cat{O}_{n-1,x} \xto{f^{\#}_{n,n-1,x}} \cat{O}_{n,x}
\end{equation}
of condensed rings is the $T_{n-1,x}$-coalescence of the condensed $\cat{O}_{n-1,x}$-algebra $\cat{O}_{n-1,x}$, where
\begin{equation}
T_{n-1,x} := \set{(f,g) \in \cat{O}_{n-1,x}(*) \times \cat{O}_{n-1,x}(*)}
{ |f|_{v_{n-1,x}} \leq |g|_{v_{n-1,x}} \neq 0 } .
\end{equation}
Then \cref{prop:construction of coalescence of vcrings} shows that the following hold.
\begin{itemize}
\item
There exists a unique continuous valuation $\tilde{v}_x$ on the condensed ring $\tilde{\cat{O}}_x$ such that $(\pi^{\#}_{0,x})^{-1}(\tilde{v}_x) = v_x$.

\item
For every $n \in \N$, we have $(\pi^{\#}_{n,x})^{-1}(\tilde{v}_x) = v_{n,x}$.

\item
The pair
\begin{equation}
\left( ( \tilde{\cat{O}}_x , \tilde{v}_x ) , 
( \cat{O}_x , v_x ) \xto{\pi^{\#}_{0,x}} ( \tilde{\cat{O}}_x , \tilde{v}_x ) \right)
\end{equation}
is the coalescence of the valued condensed ring $( \cat{O}_x , v_x )$ in the sense of \cref{df:coalescence of valued condensed ring}.
\end{itemize}
On the other hand, \cref{prop:construction of limits in the category C} shows that the set $\tilde{\cat{V}}_x$ is equal to the set of all continuous valuation $w$ on $\tilde{\cat{O}}_x$ such that $(\pi^{\#}_{n,x})^{-1}(w) = v_{n,x}$ for all $n \in \N$. Therefore we conclude that
\begin{equation}
\tilde{\cat{V}}_x = \{\tilde{v}_x\} . 
\end{equation}

Thus we have proved that the triple $(\tilde{X}, \tilde{\cat{O}}, \tilde{\cat{V}})$ is an object of $\cat{C}_1$ and that for each $x \in \tilde{X}$, the pair
\begin{equation}
\left( ( \tilde{\cat{O}}_x , \tilde{v}_x ) , 
( \cat{O}_x , v_x ) \xto{\pi^{\#}_{0,x}} ( \tilde{\cat{O}}_x , \tilde{v}_x ) \right)
\end{equation}
is the coalescence of the valued condensed ring $( \cat{O}_x , v_x )$ in the sense of \cref{df:coalescence of valued condensed ring}, where $\tilde{v}_x$ is the unique element of $\tilde{\cat{V}}_x$.

Next we show that the triple $(\tilde{X}, \tilde{\cat{O}}, \tilde{\cat{V}})$ is an object of $\cat{C}_f$. Let $U$ be any open subset of $\tilde{X}$. Let $f,g \in \tilde{\cat{O}}(U)(*)$. We prove that the set
\begin{equation}
V := \set{x \in U}{|f_x|_{\tilde{v}_x} \leq |g_x|_{\tilde{v}_x} \neq 0}
\end{equation}
is an open subset of $\tilde{X}$. In order to prove this, it suffices to show that for every $x \in U$, there exists an open neighbourhood $W$ of $x$ in $X$ such that $W \cap V$ is open in $\tilde{X}$. Therefore let us fix $x \in U$ arbitrarily. We already know from (4) of \cref{prop:construction of limits in the category D} that the pair
\begin{equation}
\left( \tilde{\cat{O}}_x , \left( \cat{O}_{n,x} \xto{\pi^{\#}_{n,x}} \tilde{\cat{O}}_x
\right)_{n \in |\N|} \right)
\end{equation}
is the colimit of the functor
\[\begin{tikzcd}[row sep=tiny]
	\N & {\ub{CRing}} \\
	n & {\cat{O}_{n,x}} & {(\text{on objects})} \\
	{(m \leq n)} & {\left( \cat{O}_{m,x} \xto{f^{\#}_{n,m,x}} \cat{O}_{n,x} \right)} & {(\text{on morphisms}) .}
	\arrow["{D_x}", from=1-1, to=1-2]
	\arrow[maps to, from=2-1, to=2-2]
	\arrow[maps to, from=3-1, to=3-2]
\end{tikzcd}\]
Since the category $\N$ is filtered, \cref{prop:evaluation of CRing preserves limits and filtered colimits} shows that the pair
\begin{equation}
\left( \tilde{\cat{O}}_x(*) ,
\left( \cat{O}_{n,x}(*) \xto{(\pi^{\#}_{n,x})_*} \tilde{\cat{O}}_x(*)
\right)_{n \in |\N|} \right)
\end{equation}
is the colimit of the functor
\[\begin{tikzcd}[row sep=tiny]
	\N & {\ub{Ring}} \\
	n & {\cat{O}_{n,x}(*)} & {(\text{on objects})} \\
	{(m \leq n)} & {\left( \cat{O}_{m,x}(*) \xto{(f^{\#}_{n,m,x})_*} \cat{O}_{n,x}(*) \right)} & {(\text{on morphisms}) .}
	\arrow[from=1-1, to=1-2]
	\arrow[maps to, from=2-1, to=2-2]
	\arrow[maps to, from=3-1, to=3-2]
\end{tikzcd}\]
Then the construction of filtered colimits in $\ub{Ring}$ shows that there exist an $n \in \N$ and elements $s,t \in \cat{O}_{n,x}(*)$ such that
\begin{equation}
(\pi^{\#}_{n,x})_*(s) = f_x \quad \text{ and } \quad (\pi^{\#}_{n,x})_*(t) = g_x .
\end{equation}
On the other hand, \cref{prop:compatibility of stalk and evaluation at S} shows that
\begin{equation}
\cat{O}_{n,x}(*) =
\underset{x \in U_n \sub X_n \text{ open}}{\colim} \; \cat{O}_n(U_n)(*) ,
\end{equation}
where the colimit is taken in $\ub{Ring}$. Then the construction of filtered colimits in $\ub{Ring}$ shows that there exist an open neighbourhood $U_n$ of $x$ in $X_n$ and elements $f_n, g_n \in \cat{O}_n(U_n)(*)$ such that
\begin{equation}
(f_n)_x = s \quad \text{ and } \quad (g_n)_x = t .
\end{equation}
Since $(X_n, \cat{O}_n, \cat{V}_n)$ is an object of $\cat{C}_f$, the set 
\begin{equation}
V_n := \set{y \in U_n}{|(f_n)_y|_{v_{n,y}} \leq |(g_n)_y|_{v_{n,y}} \neq 0}
\end{equation}
is an open subset of $X_n$. On the other hand, the following diagram is commutative.
\[\begin{tikzcd}
	{\tilde{\cat{O}}(\pi_n^{-1}(U_n))} & {\cat{O}_n(U_n)} \\
	{\tilde{\cat{O}}_x} & {\cat{O}_{n,x}}
	\arrow["\can"', from=1-1, to=2-1]
	\arrow["{(\pi^{\#}_n)_{U_n}}"', from=1-2, to=1-1]
	\arrow["\can", from=1-2, to=2-2]
	\arrow["{\pi^{\#}_{n,x}}", from=2-2, to=2-1]
\end{tikzcd}\]
Therefore we have
\begin{align}
\Big( ((\pi^{\#}_n)_{U_n})_* (f_n) \Big)_x & = (\pi^{\#}_{n,x})_*((f_n)_x) =
(\pi^{\#}_{n,x})_*(s) = f_x \: ; \\
\Big( ((\pi^{\#}_n)_{U_n})_* (g_n) \Big)_x & = (\pi^{\#}_{n,x})_*((g_n)_x) =
(\pi^{\#}_{n,x})_*(t) = g_x .
\end{align}
Moreover, \cref{prop:compatibility of stalk and evaluation at S} shows that
\begin{equation}
\tilde{\cat{O}}_x (*) = \underset{x \in W \sub \tilde{X} \text{ open}}{\colim} \; \tilde{\cat{O}}(W)(*) ,
\end{equation}
where the colimit is taken in $\ub{Ring}$. Then the construction of filtered colimits in $\ub{Ring}$ shows that there exists an open neighbourhood $W$ of $x$ in $\tilde{X}$ such that $W \sub \pi_n^{-1}(U_n) \cap U$ and 
\begin{equation}
((\pi^{\#}_n)_{U_n})_* (f_n) |_W = f|_W \quad ; \quad
((\pi^{\#}_n)_{U_n})_* (g_n) |_W = g|_W .
\end{equation}
We claim that
\begin{equation}
W \cap V = W \cap \pi_n^{-1}(V_n) .
\end{equation}
Indeed, for $y \in W$, we have
\begin{equation}
f_y = \Big( ((\pi^{\#}_n)_{U_n})_* (f_n) \Big)_y \quad \text{ and } \quad
g_y = \Big( ((\pi^{\#}_n)_{U_n})_* (g_n) \Big)_y
\end{equation}
since $((\pi^{\#}_n)_{U_n})_* (f_n) |_W = f|_W$ and $((\pi^{\#}_n)_{U_n})_* (g_n) |_W = g|_W$. Furthermore, the commutativity of the diagram
\[\begin{tikzcd}
	{\tilde{\cat{O}}(\pi_n^{-1}(U_n))} & {\cat{O}_n(U_n)} \\
	{\tilde{\cat{O}}_y} & {\cat{O}_{n,y}}
	\arrow["\can"', from=1-1, to=2-1]
	\arrow["{(\pi^{\#}_n)_{U_n}}"', from=1-2, to=1-1]
	\arrow["\can", from=1-2, to=2-2]
	\arrow["{\pi^{\#}_{n,y}}", from=2-2, to=2-1]
\end{tikzcd}\]
shows that
\begin{align}
\Big( ((\pi^{\#}_n)_{U_n})_* (f_n) \Big)_y = (\pi^{\#}_{n,y})_*((f_n)_y) \: ; \\
\Big( ((\pi^{\#}_n)_{U_n})_* (g_n) \Big)_y = (\pi^{\#}_{n,y})_*((g_n)_y) .
\end{align}
Therefore $f_y = (\pi^{\#}_{n,y})_*((f_n)_y)$ and $g_y = (\pi^{\#}_{n,y})_*((g_n)_y)$. On the other hand, by the construction of $\tilde{v}_y$, we have
\begin{equation}
(\pi^{\#}_{n,y})^{-1} (\tilde{v}_y) = v_{n,y} .
\end{equation}
Consequently, the condition
\begin{equation}
|f_y|_{\tilde{v}_y} \leq |g_y|_{\tilde{v}_y} \neq 0
\end{equation}
is equivalent to the condition
\begin{equation}
|(f_n)_y|_{v_{n,y}} \leq |(g_n)_y|_{v_{n,y}} \neq 0 .
\end{equation}
Therefore we have $y \in V$ if and only if $\pi_n(y) = y \in V_n$. This shows that $W \cap V = W \cap \pi_n^{-1}(V_n)$. Then, from the fact that $W$ is an open subset of $\tilde{X}$ and that $V_n$ is an open subset of $X_n$, we conclude that $W \cap V$ is an open subset of $\tilde{X}$. This completes the proof that the triple $(\tilde{X}, \tilde{\cat{O}}, \tilde{\cat{V}})$ is an object of $\cat{C}_f$.

Thus we have proved that $(\tilde{X}, \tilde{\cat{O}}, \tilde{\cat{V}})$ is an object of $\cat{C}_f$. Since we have already proved that $( \tilde{\cat{O}}_x , \tilde{v}_x )$ is an object of $\ub{VCRing}_c$ for every $x \in \tilde{X}$, we conclude that the triple $(\tilde{X}, \tilde{\cat{O}}, \tilde{\cat{V}})$ is an object of $\cat{C}_c$.

Next suppose that $(Y, \cat{O}_Y, \cat{V}_Y)$ is an object of $\cat{C}_c$ and that $(\rho, \rho^{\#}) : (Y, \cat{O}_Y, \cat{V}_Y) \to (X, \cat{O}, \cat{V})$ is a morphism in $\cat{C}_f$. We prove that there exists a unique morphism $(\sigma, \sigma^{\#}) : (Y, \cat{O}_Y, \cat{V}_Y) \to (\tilde{X}, \tilde{\cat{O}}, \tilde{\cat{V}})$ in $\cat{C}_c$ such that the following diagram is commutative.
\[\begin{tikzcd}
	{(Y, \cat{O}_Y, \cat{V}_Y)} & {(\tilde{X}, \tilde{\cat{O}}, \tilde{\cat{V}})} \\
	& {(X, \cat{O}, \cat{V})}
	\arrow["{(\sigma,\sigma^{\#})}", from=1-1, to=1-2]
	\arrow["{(\rho,\rho^{\#})}"', from=1-1, to=2-2]
	\arrow["{(\pi_0, \pi^{\#}_0)}", from=1-2, to=2-2]
\end{tikzcd}\]
We construct a cone
\begin{equation}
\left( (Y, \cat{O}_Y, \cat{V}_Y),
\left( (Y, \cat{O}_Y, \cat{V}_Y) \xto{(\rho_n, \rho^{\#}_n)} (X_n, \cat{O}_n, \cat{V}_n) \right)_{n \in |\N^{\op}|} \right)
\end{equation}
on the functor
\[\begin{tikzcd}
	{\N^{\op}} & {\cat{C}_f} & {\cat{C}.}
	\arrow["D", from=1-1, to=1-2]
	\arrow["\inc", hook, from=1-2, to=1-3]
\end{tikzcd}\]
inductively as follows.
\begin{itemize}
\item
We define $(\rho_0, \rho^{\#}_0) := (\rho, \rho^{\#}) : (Y, \cat{O}_Y, \cat{V}_Y) \to (X, \cat{O}, \cat{V}) = (X_0, \cat{O}_0, \cat{V}_0)$.

\item
Let $n \in \N$ with $n \geq 1$. Suppose that we have defined the morphism $(\rho_{n-1}, \rho^{\#}_{n-1}) : (Y, \cat{O}_Y, \cat{V}_Y) \to (X_{n-1}, \cat{O}_{n-1}, \cat{V}_{n-1})$ in $\cat{C}$. Since $(Y, \cat{O}_Y, \cat{V}_Y)$ is an object of $\cat{C}_c$, \cref{lem:induction step for construction of coreflection in C c} shows that there exists a unique morphism $(\rho_n, \rho^{\#}_n) : (Y, \cat{O}_Y, \cat{V}_Y) \to (X_n, \cat{O}_n, \cat{V}_n)$ in $\cat{C}$ such that the diagram
\[\begin{tikzcd}
	{(Y, \cat{O}_Y, \cat{V}_Y)} & {(X_n, \cat{O}_n, \cat{V}_n)} \\
	& {(X_{n-1}, \cat{O}_{n-1}, \cat{V}_{n-1})}
	\arrow["{(\rho_n, \rho^{\#}_n)}", from=1-1, to=1-2]
	\arrow["{(\rho_{n-1}, \rho^{\#}_{n-1})}"', from=1-1, to=2-2]
	\arrow["{(f_{n,n-1},f^{\#}_{n,n-1})}", from=1-2, to=2-2]
\end{tikzcd}\]
is commutative. This is the definition of $(\rho_n, \rho^{\#}_n) : (Y, \cat{O}_Y, \cat{V}_Y) \to (X_n, \cat{O}_n, \cat{V}_n)$.
\end{itemize}
Since the cone
\begin{equation}
\left( (\tilde{X}, \tilde{\cat{O}}, \tilde{\cat{V}}),
\left( (\tilde{X}, \tilde{\cat{O}}, \tilde{\cat{V}}) \xto{(\pi_n, \pi^{\#}_n)} (X_n, \cat{O}_n, \cat{V}_n) \right)_{n \in |\N^{\op}|} \right)
\end{equation}
is the limit of the functor
\[\begin{tikzcd}
	{\N^{\op}} & {\cat{C}_f} & {\cat{C} ,}
	\arrow["D", from=1-1, to=1-2]
	\arrow["\inc", hook, from=1-2, to=1-3]
\end{tikzcd}\]
there exists a unique morphism $(\sigma, \sigma^{\#}) : (Y, \cat{O}_Y, \cat{V}_Y) \to (\tilde{X}, \tilde{\cat{O}}, \tilde{\cat{V}})$ in $\cat{C}$ such that the following diagram is commutative for every $n \in \N$.
\[\begin{tikzcd}
	{(Y, \cat{O}_Y, \cat{V}_Y)} & {(\tilde{X}, \tilde{\cat{O}}, \tilde{\cat{V}})} \\
	& {(X_n, \cat{O}_n, \cat{V}_n)}
	\arrow["{(\sigma, \sigma^{\#})}", from=1-1, to=1-2]
	\arrow["{(\rho_n, \rho^{\#}_n)}"', from=1-1, to=2-2]
	\arrow["{(\pi_n, \pi^{\#}_n)}", from=1-2, to=2-2]
\end{tikzcd}\]
In particular, the following diagram is commutative.
\[\begin{tikzcd}
	{(Y, \cat{O}_Y, \cat{V}_Y)} & {(\tilde{X}, \tilde{\cat{O}}, \tilde{\cat{V}})} \\
	& {(X, \cat{O}, \cat{V})}
	\arrow["{(\sigma, \sigma^{\#})}", from=1-1, to=1-2]
	\arrow["{(\rho, \rho^{\#})}"', from=1-1, to=2-2]
	\arrow["{(\pi_0, \pi^{\#}_0)}", from=1-2, to=2-2]
\end{tikzcd}\]
On the other hand, suppose that $(\tau, \tau^{\#}) : (Y, \cat{O}_Y, \cat{V}_Y) \to (\tilde{X}, \tilde{\cat{O}}, \tilde{\cat{V}})$ is another morphism in $\cat{C}_c$ such that the following diagram is commutative.
\[\begin{tikzcd}
	{(Y, \cat{O}_Y, \cat{V}_Y)} & {(\tilde{X}, \tilde{\cat{O}}, \tilde{\cat{V}})} \\
	& {(X, \cat{O}, \cat{V})}
	\arrow["{(\tau, \tau^{\#})}", from=1-1, to=1-2]
	\arrow["{(\rho,\rho^{\#})}"', from=1-1, to=2-2]
	\arrow["{(\pi_0, \pi^{\#}_0)}", from=1-2, to=2-2]
\end{tikzcd}\]
We prove that $(\tau, \tau^{\#}) = (\sigma,\sigma^{\#})$. By the definition of $(\sigma, \sigma^{\#})$, it suffices to show that the diagram
\[\begin{tikzcd}
	{(Y, \cat{O}_Y, \cat{V}_Y)} & {(\tilde{X}, \tilde{\cat{O}}, \tilde{\cat{V}})} \\
	& {(X_n, \cat{O}_n, \cat{V}_n)}
	\arrow["{(\tau, \tau^{\#})}", from=1-1, to=1-2]
	\arrow["{(\rho_n, \rho^{\#}_n)}"', from=1-1, to=2-2]
	\arrow["{(\pi_n, \pi^{\#}_n)}", from=1-2, to=2-2]
\end{tikzcd}\]
is commutative for every $n \in \N$. For this we use induction on $n$. The case $n=0$ follows directly from the assumption. Let $n \in \N$ with $n \geq 1$ and suppose that the diagram
\[\begin{tikzcd}
	{(Y, \cat{O}_Y, \cat{V}_Y)} & {(\tilde{X}, \tilde{\cat{O}}, \tilde{\cat{V}})} \\
	& {(X_{n-1}, \cat{O}_{n-1}, \cat{V}_{n-1})}
	\arrow["{(\tau, \tau^{\#})}", from=1-1, to=1-2]
	\arrow["{(\rho_{n-1}, \rho^{\#}_{n-1})}"', from=1-1, to=2-2]
	\arrow["{(\pi_{n-1}, \pi^{\#}_{n-1})}", from=1-2, to=2-2]
\end{tikzcd}\]
is commutative. Then the following diagram is commutative.
\[\begin{tikzcd}
	{(Y, \cat{O}_Y, \cat{V}_Y)} & {(\tilde{X}, \tilde{\cat{O}}, \tilde{\cat{V}})} & {(X_n, \cat{O}_n, \cat{V}_n)} \\
	&& {(X_{n-1}, \cat{O}_{n-1}, \cat{V}_{n-1})}
	\arrow["{(\tau, \tau^{\#})}", from=1-1, to=1-2]
	\arrow["{(\rho_{n-1}, \rho^{\#}_{n-1})}"', curve={height=18pt}, from=1-1, to=2-3]
	\arrow["{(\pi_n, \pi^{\#}_n)}", from=1-2, to=1-3]
	\arrow["{(\pi_{n-1}, \pi^{\#}_{n-1})}"'{pos=0.2}, from=1-2, to=2-3]
	\arrow["{(f_{n,n-1},f^{\#}_{n,n-1})}", from=1-3, to=2-3]
\end{tikzcd}\]
Then the definition of $(\rho_n, \rho^{\#}_n)$ shows that $(\pi_n, \pi^{\#}_n) \of (\tau, \tau^{\#}) = (\rho_n, \rho^{\#}_n)$. In other words, the diagram
\[\begin{tikzcd}
	{(Y, \cat{O}_Y, \cat{V}_Y)} & {(\tilde{X}, \tilde{\cat{O}}, \tilde{\cat{V}})} \\
	& {(X_n, \cat{O}_n, \cat{V}_n)}
	\arrow["{(\tau, \tau^{\#})}", from=1-1, to=1-2]
	\arrow["{(\rho_n, \rho^{\#}_n)}"', from=1-1, to=2-2]
	\arrow["{(\pi_n, \pi^{\#}_n)}", from=1-2, to=2-2]
\end{tikzcd}\]
is commutative. This completes the proof.
\end{proof}

\subsection{Relation between $\cat{C}_l$, $\cat{C}_f$ and $\cat{C}_c$}

\begin{prop} \label{prop:C l and C f are coreflective in C 1}
The categories $\cat{C}_l$ and $\cat{C}_f$ are coreflective full subcategories of $\cat{C}_1$. The category $\cat{C}_c$ is a coreflective full subcategories of $\cat{C}_f$.
\end{prop}

\begin{proof}
\cref{prop:construction of coreflection in C l,prop:construction of coreflection in C f} show that every object of $\cat{C}_1$ has a coreflection along the inclusion functors $\cat{C}_l \mon \cat{C}_1$, $\cat{C}_f \mon \cat{C}_1$. \cref{prop:construction of coreflection in C c} shows that every object of $\cat{C}_f$ has a coreflection along the inclusion functor $\cat{C}_c \mon \cat{C}_f$.
\end{proof}

\begin{nt}
The right adjoints of the inclusion functors $\cat{C}_l \mon \cat{C}_1$, $\cat{C}_f \mon \cat{C}_1$, $\cat{C}_c \mon \cat{C}_f$ are denoted by
\[\begin{tikzcd}[row sep=tiny]
	{\cat{C}_1} & {\cat{C}_l ,} & X & {X_l \: ;} \\
	{\cat{C}_1} & {\cat{C}_f ,} & X & {X_f \: ;} \\
	{\cat{C}_f} & {\cat{C}_c ,} & X & {X_c \: ,}
	\arrow[from=1-1, to=1-2]
	\arrow[maps to, from=1-3, to=1-4]
	\arrow[from=2-1, to=2-2]
	\arrow[maps to, from=2-3, to=2-4]
	\arrow[from=3-1, to=3-2]
	\arrow[maps to, from=3-3, to=3-4]
\end{tikzcd}\]
respectively.
\end{nt}

\begin{prop} \label{prop:if X is in C c then X l f is in C c} \,
\begin{enumerate}
\item
For every object $X \in |\cat{C}_l|$, we have $X_f \in |\cat{C}_l|$.

\item
For every object $X \in |\cat{C}_c|$, we have $(X_l)_f\in |\cat{C}_c|$.
\end{enumerate}
\end{prop}

\begin{proof}~
\begin{enumerate}
\item
Let $X = (X, \cat{O}, \cat{V})$ be an object of $\cat{C}_l$. Write $X_f = (\tilde{X}, \tilde{\cat{O}}, \tilde{\cat{V}})$. Let $(\pi, \pi^{\#}) : X_f \to X$ be the canonical morphism. Since coreflections along functors are unique up to canonical isomorphisms, we may assume that the pair $(X_f , X_f \xto{(\pi, \pi^{\#})} X)$ has the properties described in \cref{prop:construction of coreflection in C f}. Then for every $x \in \tilde{X}$, the valued condensed ring $(\tilde{\cat{O}}_x, |\cdot|_x)$ is isomorphic to the valued condensed ring $(\cat{O}_{\pi(x)}, |\cdot|_{\pi(x)})$, which is an object of $\ub{VCRing}_l$ since $X$ is an object of $\cat{C}_l$. Therefore the valued condensed ring $(\tilde{\cat{O}}_x, |\cdot|_x)$ is also an object of $\ub{VCRing}_l$. Since this holds for every $x \in \tilde{X}$, we conclude that $X_f = (\tilde{X}, \tilde{\cat{O}}, \tilde{\cat{V}})$ is an object of $\cat{C}_l$.

\item
Let $X = (X, \cat{O}, \cat{V})$ be an object of $\cat{C}_c$. Write $X_l = (X_l, \cat{O}_l, \cat{V}_l)$ and $(X_l)_f = (X_{lf}, \cat{O}_{lf}, \cat{V}_{lf})$. Let $(\pi, \pi^{\#}) : X_l \to X$ and $(\rho, \rho^{\#}) : (X_l)_f \to X_l$ be the canonical morphisms. Since coreflections along functors are unique up to canonical isomorphisms, we may assume that the pairs $(X_l , X_l \xto{(\pi, \pi^{\#})} X)$ and $((X_l)_f , (X_l)_f \xto{(\rho, \rho^{\#})} X_l)$ have the properties described in \cref{prop:construction of coreflection in C l} and \cref{prop:construction of coreflection in C f}, respectively. 
Then for every $x \in X_{lf}$, the valued condensed ring $(\cat{O}_{l,\rho(x)}, |\cdot|_{\rho(x)})$ is the localization of the valued condensed ring $(\cat{O}_{\pi \of \rho(x)}, |\cdot|_{\pi \of \rho(x)})$ in the sense of \cref{df:localization of valued condensed ring}, and the valued condensed ring $(\cat{O}_{lf,x}, |\cdot|_x)$ is isomorphic to the valued condensed ring $(\cat{O}_{l,\rho(x)}, |\cdot|_{\rho(x)})$. Since $X$ is an object of $\cat{C}_c$, the valued condensed ring $(\cat{O}_{\pi \of \rho(x)}, |\cdot|_{\pi \of \rho(x)})$ is an object of $\ub{VCRing}_c$. Then \cref{prop:localization of objects of vcring c is in vcring c} shows that the valued condensed ring $(\cat{O}_{l,\rho(x)}, |\cdot|_{\rho(x)})$ is also an object of $\ub{VCRing}_c$. Then $(\cat{O}_{lf,x}, |\cdot|_x)$ is also an object of $\ub{VCRing}_c$. Since this holds for every $x \in X_{lf}$, we conclude that $(X_l)_f = (X_{lf}, \cat{O}_{lf}, \cat{V}_{lf})$ is an object of $\cat{C}_c$.
\end{enumerate}
\end{proof}

\begin{prop} \label{prop:construction of coreflection in C l cap C c}
Let $(X, \cat{O}, \cat{V})$ be an object of $\cat{C}$. Then there exists a coreflection
\begin{equation}
\left( (\tilde{X}, \tilde{\cat{O}}, \tilde{\cat{V}}) ,\,  
(\tilde{X}, \tilde{\cat{O}}, \tilde{\cat{V}}) \xto{(\pi, \pi^{\#})} (X, \cat{O}, \cat{V})
\right)
\end{equation}
of $(X, \cat{O}, \cat{V})$ along the inclusion functor $\cat{C}_l \cap \cat{C}_c \mon \cat{C}$ which has the following property.
\begin{enumerate}
\item
The underlying set of $\tilde{X}$ is the set of all pairs $(x,v)$ consisting of a point $x \in X$ and a valuation $v \in \cat{V}_x$.

\item
The map $\tilde{X} \xto{\pi} X$ is equal to the map $(x,v) \mapsto x$.

\item
For every $(x,v) \in \tilde{X}$, the homomorphism $\pi^{\#}_{(x,v)} : \cat{O}_x \to \tilde{\cat{O}}_{(x,v)}$ of condensed rings satisfies $\left( \pi^{\#}_{(x,v)} \right)^{-1}(\tilde{v}) = v$, where $\tilde{v}$ denotes the unique element of $\tilde{\cat{V}}_{(x,v)}$. Moreover, the pair
\begin{equation}
\left( (\tilde{\cat{O}}_{(x,v)} , \tilde{v}) \; , \; (\cat{O}_x , v) \xto{\pi^{\#}_{(x,v)}} (\tilde{\cat{O}}_{(x,v)} , \tilde{v}) \right)
\end{equation}
is the reflection of the valued condensed ring $(\cat{O}_x , v)$ along the inclusion functor $\ub{VCRing}_l \cap \ub{VCRing}_c \mon \ub{VCRing}$.
\end{enumerate}
\end{prop}

\begin{proof}
By \cref{prop:construction of coreflection in C 1}, there exists a coreflection
\begin{equation}
\left( (X_1, \cat{O}_1, \cat{V}_1) ,\,  
(X_1, \cat{O}_1, \cat{V}_1) \xto{(\alpha, \alpha^{\#})} (X, \cat{O}, \cat{V})
\right)
\end{equation}
of $(X, \cat{O}, \cat{V}) \in |\cat{C}|$ along the inclusion functor $\cat{C}_1 \mon \cat{C}$ which has the following property.
\begin{enumerate}
\item
The underlying set of $X_1$ is the set of all pairs $(x,v)$ consisting of a point $x \in X$ and a valuation $v \in \cat{V}_x$.

\item
The map $X_1 \xto{\alpha} X$ is equal to the map $(x,v) \mapsto x$.

\item
For every $(x,v) \in X_1$, the homomorphism $\alpha^{\#}_{(x,v)} : \cat{O}_x \to \cat{O}_{1,(x,v)}$ is an isomorphism of condensed rings. If $|\cdot|_{(x,v)}$ denotes the unique element of $\cat{V}_{1,(x,v)}$, then $\left( \alpha^{\#}_{(x,v)} \right)^{-1}(|\cdot|_{(x,v)}) = v$.
\end{enumerate}
By \cref{prop:construction of coreflection in C f}, there exists a coreflection
\begin{equation}
\left( (X_2, \cat{O}_2, \cat{V}_2) ,\,  
(X_2, \cat{O}_2, \cat{V}_2) \xto{(\beta, \beta^{\#})} (X_1, \cat{O}_1, \cat{V}_1)
\right)
\end{equation}
of $(X_1, \cat{O}_1, \cat{V}_1) \in |\cat{C}_1|$ along the inclusion functor $\cat{C}_f \mon \cat{C}_1$ which has the following property.
\begin{enumerate}
\item
The underlying set of $X_2$ is equal to the underlying set of $X_1$, and the map $X_2 \xto{\beta} X_1$ is equal to the identity map $x \mapsto x$.

\item
For every $x \in X_2$, the homomorphism $\beta^{\#}_x : (\cat{O}_{1,x}, |\cdot|_x) \to (\cat{O}_{2,x}, |\cdot|_x)$ is an isomorphism of valued condensed rings.
\end{enumerate}
By \cref{prop:construction of coreflection in C c}, there exists a coreflection
\begin{equation}
\left( (X_3, \cat{O}_3, \cat{V}_3) ,\,  
(X_3, \cat{O}_3, \cat{V}_3) \xto{(\gamma, \gamma^{\#})} (X_2, \cat{O}_2, \cat{V}_2)
\right)
\end{equation}
of $(X_2, \cat{O}_2, \cat{V}_2) \in |\cat{C}_f|$ along the inclusion functor $\cat{C}_c \mon \cat{C}_f$ which has the following property.
\begin{enumerate}
\item
The underlying set of $X_3$ is equal to the underlying set of $X_2$, and the map $X_3 \xto{\pi} X_2$ is equal to the identity map $x \mapsto x$.

\item
For every $x \in X_3$, the homomorphism $\gamma^{\#}_x : (\cat{O}_{2,x}, |\cdot|_x) \to (\cat{O}_{3,x}, |\cdot|_x)$ of valued condensed rings is the coalescence of the valued condensed ring $(\cat{O}_{2,x}, |\cdot|_x)$ in the sense of \cref{df:coalescence of valued condensed ring}.
\end{enumerate}
By \cref{prop:construction of coreflection in C l}, there exists a coreflection
\begin{equation}
\left( (X_4, \cat{O}_4, \cat{V}_4) ,\,  
(X_4, \cat{O}_4, \cat{V}_4) \xto{(\delta, \delta^{\#})} (X_3, \cat{O}_3, \cat{V}_3)
\right)
\end{equation}
of $(X_3, \cat{O}_3, \cat{V}_3) \in |\cat{C}_1|$ along the inclusion functor $\cat{C}_l \mon \cat{C}_1$ which has the following properties.
\begin{enumerate}
\item
The underlying set of $X_4$ is equal to the underlying set of $X_3$, and the map $X_4 \xto{\pi} X_3$ is equal to the identity map $x \mapsto x$.

\item
For every $x \in X_4$, the homomorphism $\delta^{\#}_x : (\cat{O}_{3,x}, |\cdot|_x) \to (\cat{O}_{4,x}, |\cdot|_x)$ of valued condensed rings is the localization of the valued condensed ring $(\cat{O}_{3,x}, |\cdot|_x)$ in the sense of \cref{df:localization of valued condensed ring}.
\end{enumerate}
By \cref{prop:construction of coreflection in C f}, there exists a coreflection
\begin{equation}
\left( (\tilde{X}, \tilde{\cat{O}}, \tilde{\cat{V}}) ,\,  
(\tilde{X}, \tilde{\cat{O}}, \tilde{\cat{V}}) \xto{(\epsilon, \epsilon^{\#})} (X_4, \cat{O}_4, \cat{V}_4)
\right)
\end{equation}
of $(X_4, \cat{O}_4, \cat{V}_4) \in |\cat{C}_1|$ along the inclusion functor $\cat{C}_f \mon \cat{C}_1$ which has the following property.
\begin{enumerate}
\item
The underlying set of $\tilde{X}$ is equal to the underlying set of $X_4$, and the map $\tilde{X} \xto{\epsilon} X$ is equal to the identity map $x \mapsto x$.

\item
For every $x \in \tilde{X}$, the homomorphism $\epsilon^{\#}_x : (\cat{O}_{4,x}, |\cdot|_x) \to (\tilde{\cat{O}}_x, |\cdot|_x)$ is an isomorphism of valued condensed rings.
\end{enumerate}

Let $(\pi,\pi^{\#}) : (\tilde{X}, \tilde{\cat{O}}, \tilde{\cat{V}}) \to (X, \cat{O}, \cat{V})$ be the composition $(\alpha, \alpha^{\#}) \of (\beta, \beta^{\#}) \of (\gamma, \gamma^{\#}) \of (\delta, \delta^{\#}) \of (\epsilon, \epsilon^{\#})$ in the category $\cat{C}$. We claim that the pair 
\begin{equation}
\left( (\tilde{X}, \tilde{\cat{O}}, \tilde{\cat{V}}) ,\,  
(\tilde{X}, \tilde{\cat{O}}, \tilde{\cat{V}}) \xto{(\pi, \pi^{\#})} (X, \cat{O}, \cat{V})
\right)
\end{equation}
is a coreflection of $(X, \cat{O}, \cat{V}) \in |\cat{C}|$ along the inclusion functor $\cat{C}_l \cap \cat{C}_c \mon \cat{C}$ which has all the following properties stated in this proposition.

First of all, the construction shows that the underlying set of $\tilde{X}$ is equal to the underlying set of $X_1$ and that the map $\beta \of \gamma \of \delta \of \epsilon : \tilde{X} \to X_1$ is equal to the identity map $x \mapsto x$. Moreover, the underlying set of $X_1$ is equal to the set of all pairs $(x,v)$ consisting of a point $x \in X$ and a valuation $v \in \cat{V}_x$, and the map $X_1 \xto{\alpha} X$ is equal to the map $(x,v) \mapsto x$. Consequently, the underlying set of $\tilde{X}$ is also equal to the set of all pairs $(x,v)$ consisting of a point $x \in X$ and a valuation $v \in \cat{V}_x$, and the map $\pi = \alpha \of \beta \of \gamma \of \delta \of \epsilon : \tilde{X} \to X$ is equal to the map $(x,v) \mapsto x$.

Next let $\xi = (x,v) \in \tilde{X}$. The homomomorphism $\pi^{\#}_{\xi} : \cat{O}_x \to \tilde{\cat{O}}_{\xi}$ of condensed rings is equal to the composition
\[\begin{tikzcd}
	{\cat{O}_x} & {\cat{O}_{1,\xi}} & {\cat{O}_{2,\xi}} & {\cat{O}_{3,\xi}} & {\cat{O}_{4,\xi}} & {\tilde{\cat{O}}_{\xi} ,}
	\arrow["{\alpha^{\#}_{\xi}}", from=1-1, to=1-2]
	\arrow["{\beta^{\#}_{\xi}}", from=1-2, to=1-3]
	\arrow["{\gamma^{\#}_{\xi}}", from=1-3, to=1-4]
	\arrow["{\delta^{\#}_{\xi}}", from=1-4, to=1-5]
	\arrow["{\epsilon^{\#}_{\xi}}", from=1-5, to=1-6]
\end{tikzcd}\]
and the construction shows that these are in fact homomorphisms of valued condensed rings
\[\begin{tikzcd}
	{(\cat{O}_x,v)} & {(\cat{O}_{1,\xi},|\cdot|_{\xi})} & {(\cat{O}_{2,\xi},|\cdot|_{\xi})} & {(\cat{O}_{3,\xi},|\cdot|_{\xi})} & {(\cat{O}_{4,\xi},|\cdot|_{\xi})} & {(\tilde{\cat{O}}_{\xi},|\cdot|_{\xi}) .}
	\arrow["{\alpha^{\#}_{\xi}}", from=1-1, to=1-2]
	\arrow["{\beta^{\#}_{\xi}}", from=1-2, to=1-3]
	\arrow["{\gamma^{\#}_{\xi}}", from=1-3, to=1-4]
	\arrow["{\delta^{\#}_{\xi}}", from=1-4, to=1-5]
	\arrow["{\epsilon^{\#}_{\xi}}", from=1-5, to=1-6]
\end{tikzcd}\]
Moreover, the construction shows that the following hold.
\begin{itemize}
\item
The homomorphisms
$\alpha^{\#}_{\xi} : (\cat{O}_x,v) \to (\cat{O}_{1,\xi},|\cdot|_{\xi})$, 
$\beta^{\#}_{\xi} : (\cat{O}_{1,\xi}, |\cdot|_{\xi}) \to (\cat{O}_{2,\xi}, |\cdot|_{\xi})$ and
$\epsilon^{\#}_{\xi} : (\cat{O}_{4,\xi}, |\cdot|_{\xi}) \to (\tilde{\cat{O}}_{\xi}, |\cdot|_{\xi})$ are isomorphisms of valued condensed rings.

\item
The homomorphism $\gamma^{\#}_{\xi} : (\cat{O}_{2,\xi}, |\cdot|_{\xi}) \to (\cat{O}_{3,\xi}, |\cdot|_{\xi})$ of valued condensed rings is the coalescence of the valued condensed ring $(\cat{O}_{2,\xi}, |\cdot|_{\xi})$ in the sense of \cref{df:coalescence of valued condensed ring}.

\item
The homomorphism $\delta^{\#}_{\xi} : (\cat{O}_{3,\xi}, |\cdot|_{\xi}) \to (\cat{O}_{4,\xi}, |\cdot|_{\xi})$ of valued condensed rings is the localization of the valued condensed ring $(\cat{O}_{3,\xi}, |\cdot|_{\xi})$ in the sense of \cref{df:localization of valued condensed ring}.
\end{itemize}
Then \cref{cor:coalescence followed by localization of vcrings} shows that the pair
\begin{equation}
\left( (\tilde{\cat{O}}_{\xi},|\cdot|_{\xi}) \; , \; (\cat{O}_x,v) \xto{\pi^{\#}_{\xi}} (\tilde{\cat{O}}_{\xi},|\cdot|_{\xi}) \right)
\end{equation}
is the reflection of the valued condensed ring $(\cat{O}_x,v)$ along the inclusion functor $\ub{VCRing}_l \cap \ub{VCRing}_c \mon \ub{VCRing}$.

It remains to prove that the pair 
\begin{equation}
\left( (\tilde{X}, \tilde{\cat{O}}, \tilde{\cat{V}}) ,\,  
(\tilde{X}, \tilde{\cat{O}}, \tilde{\cat{V}}) \xto{(\pi, \pi^{\#})} (X, \cat{O}, \cat{V})
\right)
\end{equation}
is a coreflection of $(X, \cat{O}, \cat{V}) \in |\cat{C}|$ along the inclusion functor $\cat{C}_l \cap \cat{C}_c \mon \cat{C}$. For simplicity, let us write
\begin{equation}
X = (X, \cat{O}, \cat{V}) \quad ; \quad
\tilde{X} = (\tilde{X}, \tilde{\cat{O}}, \tilde{\cat{V}}) \quad ; \quad
X_i = (X_i, \cat{O}_i, \cat{V}_i) \quad (i = 1,2,3,4) .
\end{equation}
First of all, $X_3$ is an object of $\cat{C}_c$ by definition. Then \cref{prop:if X is in C c then X l f is in C c} shows that $\tilde{X} = ((X_3)_l)_f$ is an object of $\cat{C}_l \cap \cat{C}_c$. On the other hand, suppose that $Y$ is an object of $\cat{C}_l \cap \cat{C}_c$. Since $Y$ is an object of $\cat{C}_f$, the universality of the coreflection $(\tilde{X}, \tilde{X} \xto{(\epsilon, \epsilon^{\#})} X_4)$ of $X_4$ along the inclusion functor $\cat{C}_f \mon \cat{C}_1$ shows that the map
\[\begin{tikzcd}
	{\cat{C}(Y,\tilde{X})} & {\cat{C}(Y,X_4),} & {(\sigma,\sigma^{\#})} & {(\epsilon, \epsilon^{\#}) \of (\sigma,\sigma^{\#})}
	\arrow[from=1-1, to=1-2]
	\arrow[maps to, from=1-3, to=1-4]
\end{tikzcd}\]
is bijective. Since $Y$ is an object of $\cat{C}_l$, the universality of the coreflection $(X_4, X_4 \xto{(\delta, \delta^{\#})} X_3)$ of $X_3$ along the inclusion functor $\cat{C}_l \mon \cat{C}_1$ shows that the map
\[\begin{tikzcd}
	{\cat{C}(Y,X_4)} & {\cat{C}(Y,X_3),} & {(\sigma,\sigma^{\#})} & {(\delta, \delta^{\#}) \of (\sigma,\sigma^{\#})}
	\arrow[from=1-1, to=1-2]
	\arrow[maps to, from=1-3, to=1-4]
\end{tikzcd}\]
is bijective. Since $Y$ is an object of $\cat{C}_c$, the universality of the coreflection $(X_3, X_3 \xto{(\gamma, \gamma^{\#})} X_2)$ of $X_2$ along the inclusion functor $\cat{C}_c \mon \cat{C}_f$ shows that the map
\[\begin{tikzcd}
	{\cat{C}(Y,X_3)} & {\cat{C}(Y,X_2),} & {(\sigma,\sigma^{\#})} & {(\gamma, \gamma^{\#}) \of (\sigma,\sigma^{\#})}
	\arrow[from=1-1, to=1-2]
	\arrow[maps to, from=1-3, to=1-4]
\end{tikzcd}\]
is bijective. Since $Y$ is an object of $\cat{C}_f$, the universality of the coreflection $(X_2, X_2 \xto{(\beta, \beta^{\#})} X_1)$ of $X_1$ along the inclusion functor $\cat{C}_f \mon \cat{C}_1$ shows that the map
\[\begin{tikzcd}
	{\cat{C}(Y,X_2)} & {\cat{C}(Y,X_1),} & {(\sigma,\sigma^{\#})} & {(\beta, \beta^{\#}) \of (\sigma,\sigma^{\#})}
	\arrow[from=1-1, to=1-2]
	\arrow[maps to, from=1-3, to=1-4]
\end{tikzcd}\]
is bijective. Since $Y$ is an object of $\cat{C}_1$, the universality of the coreflection $(X_1, X_1 \xto{(\alpha, \alpha^{\#})} X)$ of $X$ along the inclusion functor $\cat{C}_1 \mon \cat{C}$ shows that the map
\[\begin{tikzcd}
	{\cat{C}(Y,X_1)} & {\cat{C}(Y,X),} & {(\sigma,\sigma^{\#})} & {(\alpha, \alpha^{\#}) \of (\sigma,\sigma^{\#})}
	\arrow[from=1-1, to=1-2]
	\arrow[maps to, from=1-3, to=1-4]
\end{tikzcd}\]
is bijective. Therefore we conclude that the map
\[\begin{tikzcd}
	{\cat{C}(Y,\tilde{X})} & {\cat{C}(Y,X),} & {(\sigma,\sigma^{\#})} & {(\pi, \pi^{\#}) \of (\sigma,\sigma^{\#})}
	\arrow[from=1-1, to=1-2]
	\arrow[maps to, from=1-3, to=1-4]
\end{tikzcd}\]
is bijective. This completes the proof.
\end{proof}

\begin{prop} \label{prop:intersections of C l C f C c are coreflective in C} \;
The categories $\cat{C}_1$, $\cat{C}_l$, $\cat{C}_f$, $\cat{C}_c$, $\cat{C}_l \cap \cat{C}_f$ and $\cat{C}_l \cap \cat{C}_c$ are all coreflective full subcategories of $\cat{C}$. Moreover, these categories are all complete and cocomplete.
\end{prop}

\begin{proof}

\cref{prop:construction of coreflection in C 1} shows that every object of $\cat{C}$ has a coreflection along the inclusion functor $\cat{C}_1 \mon \cat{C}$. Therefore $\cat{C}_1$ is a coreflective full subcategory of $\cat{C}$. Moreover, we have already proved in \cref{prop:C l and C f are coreflective in C 1} that the categories $\cat{C}_l$ and $\cat{C}_f$ are coreflective full subcategories of $\cat{C}_1$. It follows that the categories $\cat{C}_l$ and $\cat{C}_f$ are coreflective full subcategories of $\cat{C}$. Then the category $\cat{C}_c$ is also a coreflective full subcategory of $\cat{C}$ since it is a coreflective full subcategories of $\cat{C}_f$ by \cref{prop:C l and C f are coreflective in C 1}. 

By \cref{prop:if X is in C c then X l f is in C c}, the functor $\cat{C}_1 \to \cat{C}_f$, $X \mapsto X_f$ restricts to a functor $\cat{C}_l \to \cat{C}_l \cap \cat{C}_f$. This is right adjoint to the inclusion functor $\cat{C}_l \cap \cat{C}_f \mon \cat{C}_l$. Therefore $\cat{C}_l \cap \cat{C}_f$ is a coreflective full subcategory of $\cat{C}_l$. Since $\cat{C}_l$ is a coreflective full subcategory of $\cat{C}$, it follows that $\cat{C}_l \cap \cat{C}_f$ is a coreflective full subcategory of $\cat{C}$. Furthermore, 
\cref{prop:construction of coreflection in C l cap C c} shows that every object of $\cat{C}$ has a coreflection along the inclusion functor $\cat{C}_l \cap \cat{C}_c \mon \cat{C}$. Therefore $\cat{C}_l \cap \cat{C}_c$ is a coreflective full subcategory of $\cat{C}$.

By \cref{cor:completeness of the category C}, the category $\cat{C}$ is complete. By \cref{cor:cocompleteness of the category C}, the category $\cat{C}$ is cocomplete. Since the categories $\cat{C}_1$, $\cat{C}_l$, $\cat{C}_f$, $\cat{C}_c$, $\cat{C}_l \cap \cat{C}_f$ and $\cat{C}_l \cap \cat{C}_c$ are all coreflective full subcategories of $\cat{C}$, they are also complete and cocomplete by the dual of Proposition 3.5.3 and Proposition 3.5.4 of \cite{Borceux:cat1}.
\end{proof}

\section{Relation to adic spectra} \label{sec:Relation to adic spectra}

In this section, we first review some basic definitions of topological algebra, and then discuss some relations between topological algebra and condensed mathematics. After that, we prove the main result of this paper, which compares the adic spectra of Huber pairs, interpreted as objects of $\cat{C}$, with the images of some simple objects under the right adjoint of the inclusion functor $\cat{C}_l \cap \cat{C}_c \mon \cat{C}$.

\subsection{Review of topological algebra}

\begin{nt}
In this paper, a complete uniform space is not necessarily Hausdorff.
\end{nt}

\subsubsection{Non-archimedean abelian groups and rings}

\begin{df} \;
\begin{enumerate}
\item
A \ti{non-archimedean abelian group} is a topological abelian group $M$ whose zero element has a fundamental system of neighbourhoods consisting of subgroups of $M$.

\item (\cite{Henkel:paper}, Definition 2.1, iii ; \cite{Morel:note}, Definition II.1.1.1, (i) ; \cite{Wedhorn:note}, Definition 5.23)
A \ti{non-archimedean ring} is a commutative unital topological ring whose underlying topological abelian group is a non-archimedean abelian group.
\end{enumerate}
\end{df}

\begin{df}
Let $M$ be a topological abelian group. A \ti{zero sequence} in $M$ is a sequence of elements of $M$ which converges to $0$ in $M$.
\end{df}

\begin{prop}
Let $M$ be a complete non-archimedean abelian group. Let $(x_n)_{n \in \N}$ be a zero sequence in $M$. Then $(x_n)_{n \in \N}$ is a summable family of elements of $M$.
\end{prop}

\begin{proof}
By the Cauchy's criterion (\cite{Bourbaki:top1}, Chapter III, \S 5.2, Theorem 1), it suffices to show that, for every neighbourhood $U$ of $0$ in $M$, there exists a finite subset $F$ of $\N$ such that $\sum_{n \in G} x_n \in U$ for every finite subset $G$ of $\N$ satisfying $F \cap G = \ku$. Since $M$ is a non-archimedean abelian group, we may assume that $U$ is a subgroup of $M$. Since $(x_n)_{n \in \N}$ is a zero sequence in $M$, there exists an $N \in \N$ such that $x_n \in U$ for every $n \in \N$ with $n > N$. Then $F := \set{n \in \N}{0 \leq n \leq N}$ is a finite subset of $\N$. Suppose that $G$ is a finite subset of $\N$ satisfying $F \cap G = \ku$. If $n \in G$, then $n > N$ and hence $x_n \in U$. Since $U$ is a subgroup of $M$, we conclude that $\sum_{n \in G} x_n \in U$. This completes the proof.
\end{proof}

\subsubsection{Huber rings}

\begin{df} (\cite{Morel:note}, Definition II.1.1.1, (iii) ; \cite{Wedhorn:note}, Definition 6.1)
\begin{enumerate}
\item
A \ti{Huber ring} is a commutative unital topological ring $A$ with the property that there exist an open subring $A_0$ of $A$ and a finitely generated ideal $I_0$ of $A_0$ such that the subspace topology induced on $A_0$ is equal to the $I_0$-adic topology.

\item
Let $A$ is a Huber ring. A \ti{couple of definition} of $A$ is a pair $(A_0, I_0)$ consisting of an open subring $A_0$ of $A$ and a finitely generated ideal $I_0$ of $A_0$ such that the subspace topology induced on $A_0$ is equal to the $I_0$-adic topology. If $(A_0, I_0)$ is a couple of definition of $A$, then $A_0$ is called a \ti{ring of definition} of $A$.
\end{enumerate}
\end{df}

\begin{df} (\cite{Henkel:paper}, Definition 2.9 ; \cite{Morel:note}, Definition II.1.1.1, (iv) ; \cite{Wedhorn:note}, Definition 6.10)

A \ti{Tate ring} is a Huber ring which has a topologically nilpotent unit.
\end{df}

\begin{prop} \tup{(\cite{Morel:note}, Corollary II.3.1.9, (iv) ; \cite{Wedhorn:note}, Remark 6.8)}

Let $A$ be a Huber ring. Let $\iota : A \to \hat{A}$ be the Hausdorff completion of $A$.
\begin{enumerate}
\item
$\hat{A}$ is a Huber ring.

\item
Suppose that $(A_0, I_0)$ is a couple of definition of $A$. Let $\hat{A_0}$ be the closure of $\iota(A_0)$ in $\hat{A}$. Then $(\hat{A_0}, \, \iota(I_0) \cdot \hat{A_0})$ is a couple of definition of $\hat{A}$.
\end{enumerate}
\end{prop}

\subsubsection{Rings of restricted formal power series}

\begin{prop} \label{prop:well-definedness of Rings of restricted formal power series}
\tup{(\cite{Morel:note}, Proposition II.3.3.2 ; \cite{Wedhorn:note}, Definition 5.48)}

Let $R$ be a non-archimedean ring. Let $n \in \N$. Consider the ring $R[[X_1 ,\ldots, X_n]]$ of all formal power series in $n$ variables over $R$. For $\nu = (\nu_1 ,\ldots, \nu_n) \in {\N}^n$, we write $X^{\nu} := {X_1}^{\nu_1} \cdots {X_n}^{\nu_n}$.
\begin{enumerate}
\item 
Let $R \langle X_1 ,\ldots, X_n \rangle$ be the subset of $R[[X_1 ,\ldots, X_n]]$ consisting of all formal power series $\sum_{\nu \in {\N}^n} r_{\nu} X^{\nu} \in R[[X_1 ,\ldots, X_n]]$ which has the following property: For every neighbourhood $U$ of $0$ in $R$, there exists a finite subset $F$ of ${\N}^n$ such that $r_{\nu} \in U$ for every $\nu \in {\N}^n \sm F$. Then $R \langle X_1 ,\ldots, X_n \rangle$ is a subring of $R[[X_1 ,\ldots, X_n]]$.

\item
There exists a unique topology $\cat{T}$ on $R \langle X_1 ,\ldots, X_n \rangle$ which satisfies the following conditions.
\begin{enumerate}
\item $R \langle X_1 ,\ldots, X_n \rangle$ becomes a non-archimedean ring if it is given the topology $\cat{T}$.
\item For every open additive subgroup $U$ of $R$, let us write $U_{\langle X \rangle}$ for the set of all $\sum_{\nu \in {\N}^n} r_{\nu} X^{\nu} \in R \langle X_1 ,\ldots, X_n \rangle$ such that $r_{\nu} \in U$ for all $\nu \in {\N}^n$. Then the set
\begin{equation}
\set{ U_{\langle X \rangle} }{ U \text{ is an open additive subgroup of } R }
\end{equation}
is a fundamental system of neighbourhoods of $0$ with respect to $\cat{T}$.
\end{enumerate}

\item The canonical inclusion $R \mon R \langle X_1 ,\ldots, X_n \rangle$ is continuous when $R \langle X_1 ,\ldots, X_n \rangle$ is endowed with the topology $\cat{T}$.
\end{enumerate}
\end{prop}

\begin{df}
Let $R$ be a non-archimedean ring. Let $n \in \N$. We use the notation of \cref{prop:well-definedness of Rings of restricted formal power series}. The topological ring $R \langle X_1 ,\ldots, X_n \rangle$, endowed with the topology $\cat{T}$, is called the \ti{ring of restricted formal power series} in $n$ variables over $R$.
\end{df}

\begin{prop} \tup{(\cite{Morel:note}, Proposition II.3.3.3, (iii) ; \cite{Wedhorn:note}, Proposition 5.49, (3))}

Let $R$ be a non-archimedean ring. Let $n \in \N$. If $R$ is complete Hausdorff, then so is $R \langle X_1 ,\ldots, X_n \rangle$.
\end{prop}

\begin{prop} \tup{(\cite{Morel:note}, Proposition II.3.3.6, (ii) ; \cite{Wedhorn:note}, Proposition 6.21, (2))}

Let $A$ be a Huber ring. Let $n \in \N$.
\begin{enumerate}
\item
$A \langle X_1 ,\ldots, X_n \rangle$ is a Huber ring.

\item
If $(A_0, I_0)$ is a couple of definition of $A$, then $( A_0 \langle X_1 ,\ldots, X_n \rangle , I_0 \cdot A_0 \langle X_1 ,\ldots, X_n \rangle )$ is a couple of definition of $A \langle X_1 ,\ldots, X_n \rangle$.
\end{enumerate}
\end{prop}

\begin{df} (\cite{Morel:note}, Definition IV.1.1.1 ; \cite{Wedhorn:note}, Definition 6.36)

Let $A$ be a Tate ring. Let $\hat{A}$ be the Hausdorff completion of $A$. Then $A$ is called \ti{strongly Noetherian} if $\hat{A} \langle X_1 ,\ldots, X_n \rangle$ is a Noetherian ring for all $n \in \N$.
\end{df}

\subsubsection{Localization}

\begin{prop}
Let $A$ be a Huber ring. Let $g \in A$. Let $T$ be a subset of $A$ which generates an open ideal of $A$. Let $\phi : A \to A_g$ be the canonical homomorphism of $A$ into the localization $A_g$ of $A$ by $g$.
\begin{enumerate}
\item
There exists a unique topology $\cat{T}$ on $A_g$ such that the following hold.
\begin{enumerate}
\item
$A_g$ is a non-archimedean ring when endowed with $\cat{T}$. We write $A \left( \frac{T}{g} \right)$ for the topological ring $A_g$ endowed with the topology $\cat{T}$.

\item
The map $\phi : A \to A \left( \frac{T}{g} \right)$ is continuous. The set $\set{\phi(f) \cdot \phi(g)^{-1}}{f \in T}$ is power-bounded in $A \left( \frac{T}{g} \right)$.

\item
Let $B$ be a non-archimedean ring. Let $\psi : A \to B$ be a continuous ring homomorphism. Suppose that $\psi(g)$ is invertible in $B$ and that $\set{\psi(f) \cdot \psi(g)^{-1}}{f \in T}$ is power-bounded in $B$. Then there exists a unique continuous ring homomorphism $\tilde{\psi} : A \left( \frac{T}{g} \right) \to B$ such that the following diagram is commutative.
\[\begin{tikzcd}
	A & B \\
	{A \left( \frac{T}{g} \right)}
	\arrow["\psi", from=1-1, to=1-2]
	\arrow["\phi"', from=1-1, to=2-1]
	\arrow["{\tilde{\psi}}"', from=2-1, to=1-2]
\end{tikzcd}\]
\end{enumerate}

\item
$A \left( \frac{T}{g} \right)$ is a Huber ring.

\item
Suppose that $(A_0, I_0)$ is a couple of definition of $A$. Let $B_0$ be the subring of $A_g$ generated by $\phi(A_0) \cup \set{\phi(f) \cdot \phi(g)^{-1}}{f \in T}$. Then $(B_0, \phi(I_0) \cdot B_0)$ is a couple of definition of $A \left( \frac{T}{g} \right)$.

\end{enumerate}
\end{prop}

\begin{proof}~
\begin{enumerate}
\item
Lemma II.3.3.5 of \cite{Morel:note} or Lemma 6.20 of \cite{Wedhorn:note} shows that, for every $n \in \N$ and every open additive subgroup $U$ of $A$, the additive subgroup $T^n \cdot U$ of $A$ is open in $A$. Then the assertion of (1) follows from Proposition II.3.4.1 of \cite{Morel:note} or Proposition 5.51 of \cite{Wedhorn:note}.

\item
This is Proposition II.3.4.3 of \cite{Morel:note} or Proposition 6.21, (3) of \cite{Wedhorn:note}.

\item
This follows from Remark II.3.4.4 of \cite{Morel:note}. This is also stated in 8.1 of \cite{Wedhorn:note}.
\end{enumerate}
\end{proof}

\begin{df} (\cite{Morel:note}, Definition II.3.4.5 ; \cite{Wedhorn:note}, Definition 5.51)

Let $A$ be a Huber ring. Let $g \in A$. Let $T$ be a subset of $A$ which generates an open ideal of $A$.
\begin{enumerate}
\item
The Hausdorff completion of $A \left( \frac{T}{g} \right)$ is denoted by $A \left \langle \frac{T}{g} \right \rangle$.

\item
Let $\phi : A \to A \left( \frac{T}{g} \right)$ be the canonical homomorphism. Let $\iota : A \left( \frac{T}{g} \right) \to A \left \langle \frac{T}{g} \right \rangle$ be the canonical homomorphism into the Hausdorff completion. The composition $\iota \of \phi : A \to A \left \langle \frac{T}{g} \right \rangle$ is also called the \ti{canonical homomorphism}.
\end{enumerate}
\end{df}

\subsubsection{First countability}

\begin{prop} \label{prop:metrizability and first countability}
Let $M$ be a Hausdorff topological abelian group. The following conditions are equivalent.
\begin{enumerate}
\item
The zero element of $M$ has a countable fundamental system of neighbourhoods in $M$.

\item
$M$ is first countable.

\item
There exists a metric $d$ on $M$ such that the uniformity induced by $d$ on $M$ is equal to the uniform structure on $M$ associated to the topological abelian group structure of $M$.
\end{enumerate}
\end{prop}

\begin{proof}
The implication $(3) \To (2) \To (1)$ is trivial. The implication $(1) \To (3)$ follows from \cite{Bourbaki:top2}, Chapter IX, \S 3.1, Proposition 1.
\end{proof}

\begin{cor} \label{cor:c H f c top ab is Baire}
Let $M$ be a complete Hausdorff first countable topological abelian group. Then $M$ is a Baire space.
\end{cor}

\begin{proof}
By \cref{prop:metrizability and first countability}, there exists a metric $d$ on $M$ with the following property.
\begin{enumerate}
\item $M$ becomes a complete metric space when equipped with the metric $d$.
\item The topology of $M$ is equal to the topology induced by the metric $d$.
\end{enumerate}
Then the result follows from the theorem of Baire (\cite{Bourbaki:top2}, Chapter IX, \S 5.3, Theorem 1).
\end{proof}

\begin{prop} \label{prop:completeness of quotients of f c topological abelian groups}
Let $M$ be a complete Hausdorff first countable topological abelian group. Let $N$ be a closed subgroup of $M$. Then the quotient $M/N$ is complete Hausdorff.
\end{prop}

\begin{proof}
By \cref{prop:metrizability and first countability}, the topology of $M$ is metrizable. Then the result follows from \cite{Bourbaki:top2}, Chapter IX, \S 3.1, Proposition 4.
\end{proof}

\subsection{Topological algebra and condensed mathematics}

\begin{prop} \label{prop:power-boundedness and coalescence}
Let $R$ be a complete Hausdorff non-archimedean ring. Let $f \in R$ be a power-bounded element of $R$. Then $\tu{R}$ is an $f/1$-coalescent condensed $\tu{R}$-module. 
\end{prop}

\begin{proof}
Let $S$ be an arbitrary light profinite set. Let $K$ be the set of all continuous maps $\phi: S \times \N_{\infty} \to R$ such that $f(x,\infty) = 0$ for every $x \in S$. We prove that the map
\[\begin{tikzcd}
	K & {K,} & \phi & {\phi - f \cdot (\phi \of \tau_S)}
	\arrow[from=1-1, to=1-2]
	\arrow[maps to, from=1-3, to=1-4]
\end{tikzcd}\]
is bijective.

Let $\psi \in K$. Let us define $\phi : S \times \N_{\infty} \to R$ as follows. For $x \in S$, define $\phi(x, \infty) := 0$. For each $(x,n) \in S \times \N$, the sequence $(\psi(x,k+n))_{k \in \N}$ is a zero sequence in $R$. Since $f$ is power-bounded in $R$, the sequence $(f^k \psi(x,k+n))_{k \in \N}$ is also a zero sequence in $R$. Since $R$ is complete Hausdorff, its sum exists and is unique. Therefore we can define
\begin{equation}
\phi(x,n) := \sum_{k \in \N} f^k \psi(x,k+n) .
\end{equation}

Next we prove that $\phi : S \times \N_{\infty} \to R$ is continuous. Let $(x_0, n_0) \in S \times \N_{\infty}$. We show that $\phi$ is continuous at $(x_0,n_0)$. Let $U$ be an arbitrary neighbourhood of $0$ in $R$ which is an additive subgroup of $R$. We show that there exists a neighbourhood $W$ of $(x_0, n_0)$ in $S \times \N_{\infty}$ such that $\phi(W) \sub \phi(x_0,n_0) + U$. 

Since $f$ is power-bounded in $R$, there exists a neighbourhood $U'$ of $0$ in $R$ which is an additive subgroup of $R$ and such that $f^k \cdot U' \sub U$ for every $k \in \N$.

Suppose that $n_0 \in \N$. Since $\psi$ is continuous at $(x_0, \infty)$, there exist a neighbourhood $G$ of $x_0$ in $S$ and an $N \in \N$ such that $\psi(x, k+n_0) \in U'$ for every $x \in G$ and every $k \in \N$ with $k \geq N$. For each $k \in \N$ with $0 \leq k < N$, the continuity of $\psi$ at $(x_0 ,k+n_0)$ shows that there exists a neighbourhood $H_k$ of $x_0$ in $S$ such that $\psi(x, k+n_0) - \psi(x_0, k+n_0) \in U'$ for every $x \in H_k$. Then $H := G \cap H_0 \cap \cdots \cap H_{N-1}$ is a neighbourhood of $x_0$ in $S$. Then $W := H \times \{n_0\}$ is a neighbourhood of $(x_0,n_0)$ in $S \times \N_{\infty}$. Let us prove $\phi(W) \sub \phi(x_0,n_0) + U$. Suppose $(x,n) \in W$. Then $n = n_0$. For every $k \in \N$, we claim that $\psi(x, k+n_0) - \psi(x_0, k+n_0) \in U'$. Indeed, if $k \geq N$, then $x, x_0 \in G$ and $k \geq N$ implies that $\psi(x, k+n_0), \psi(x_0, k+n_0) \in U'$. Therefore $\psi(x, k+n_0) - \psi(x_0, k+n_0) \in U'$. On the other hand, if $0 \leq k < N$, then $x \in H_k$ implies that $\psi(x, k+n_0) - \psi(x_0, k+n_0) \in U'$. Thus $\psi(x, k+n_0) - \psi(x_0, k+n_0) \in U'$ for every $k \in \N$. Then $f^k \psi(x, k+n_0) - f^k \psi(x_0, k+n_0) \in U$ for every $k \in \N$. Then
\begin{equation}
\phi(x,n_0) - \phi(x_0, n_0) = \sum_{k \in \N} \Big( f^k \psi(x, k+n_0) - f^k \psi(x_0, k+n_0) \Big)
\: \in U ,
\end{equation}
since $U$ is a closed additive subgroup of $R$. This proves $\phi(W) \sub \phi(x_0,n_0) + U$.

Next suppose that $n_0 = \infty$. Since $\psi$ is continuous at $(x_0, \infty)$, there exist a neighbourhood $G$ of $x_0$ in $S$ and an $N \in \N$ such that $\psi(x, k) \in U'$ for every $x \in G$ and every $k \in \N$ with $k \geq N$. Then $V := \{ \infty \} \cup \set{n \in \N}{n \geq N}$ is a neighbourhood of $\infty$ in $\N_{\infty}$. Therefore $W := G \times V$ is a neighbourhood of $(x_0, \infty)$ in $S \times \N_{\infty}$. Let us prove $\phi(W) \sub U$. If $(x,n) \in W$, then $\psi(x, k+n) \in U'$ for every $k \in \N$. Then $f^k \psi(x, k+n) \in U$ for every $k \in \N$. Since $U$ is a closed additive subgroup of $R$, we have
\begin{equation}
\phi(x,n) = \sum_{k \in \N} f^k \psi(x, k+n) \: \in U.
\end{equation}
This shows that $\phi(W) \sub U$.

Thus we have proved that $\phi : S \times \N_{\infty} \to R$ is continuous. Then $\phi \in K$. Moreover, we claim that $\phi - f \cdot (\phi \of \tau_S) = \psi$. Indeed, for every $(x,n) \in S \times \N$, we have
\begin{align}
\phi(x,n) - f \cdot \phi(x,n+1)
& = \sum_{k \in \N} f^k \psi(x,k+n) - f \cdot \sum_{k \in \N} f^k \psi(x,k+n+1) \\
& = \sum_{k \in \N} f^k \psi(x,k+n) - \sum_{\substack{k \in \N \\ k \geq 1}} f^k \psi(x,k+n) \\
& = \psi(x,n) .
\end{align}
Furthermore, we have $\phi(x,\infty) - f \cdot \phi(x,\infty) = 0 = \psi(x,\infty)$ for every $x \in S$. Thus $\phi - f \cdot (\phi \of \tau_S) = \psi$.

Next suppose that $\rho \in K$ satisfies $\rho - f \cdot (\rho \of \tau_S) = \psi$. We claim that $\rho = \phi$. For $x \in S$, we have $\rho(x, \infty) = 0 = \phi(x, \infty)$. It remains to prove that $\rho(x,n) = \phi(x,n)$ for every $(x,n) \in S \times \N$. Let $(x,n) \in S \times \N$. For $N \in \N$, we have
\begin{align}
\sum_{k=0}^{N} f^k \psi(x,k+n)
& = \sum_{k=0}^{N} f^k \Big( \rho(x,k+n) - f \cdot \rho(x,k+n+1) \Big) \\
& = \sum_{k=0}^{N} \Big( f^k \rho(x,k+n) - f^{k+1} \rho(x,k+1+n) \Big) \\
& = \rho(x,n) - f^{N+1} \rho(x, N+1+n) .
\end{align}
We pass to the limit $N \to \infty$. Then
\begin{equation}
\lim_{N \to \infty} \: \sum_{k=0}^{N} f^k \psi(x,k+n) = \phi(x,n)
\end{equation}
by definition. On the other hand, since $(\rho(x,N+1+n))_{N \in \N}$ is a zero sequence in $R$ and $f$ is power-bounded in $R$, the sequence $(f^{N+1} \rho(x,N+1+n))_{N \in \N}$ is also a zero sequence in $R$. Therefore
\begin{equation}
\lim_{N \to \infty} \: \Big( \rho(x,n) - f^{N+1} \rho(x, N+1+n) \Big) = \rho(x,n).
\end{equation}
Thus $\phi(x,n) = \rho(x,n)$. This completes the proof.
\end{proof}

\begin{prop} \label{prop:lifting compact sets}
Let $M$ be a first countable complete Hausdorff non-archimedean abelian group. Let $N$ be a closed subgroup of $M$. Let $M \xto{\pi} M/N$ be the canonical surjection. If $S$ is a compact subset of $M/N$, then there exists a subspace $\tilde{S}$ of $M$ which is a light profinite set and satisfies $\pi(\tilde{S}) = S$.
\end{prop}

\begin{proof}
Since $M$ is a first countable non-archimedean abelian group, there exists a family $(U_n)_{n \in \N}$ which satisfies the following conditions.
\begin{itemize}
\item The set $\set{U_n}{n \in \N}$ is a fundamental system of neighbourhoods of $0$ in $M$.
\item $U_{n+1} \sub U_n$ for every $n \in \N$.
\item For every $n \in \N$, $U_n$ is a subgroup of $M$.
\end{itemize}
For each $n \in \N$, define $P_n := M/U_n$ and $Q_n := (M/N)/(\pi(U_n))$. These are discrete abelian groups. Write $M \xto{\sigma_n} P_n$ and $M/N \xto{\tau_n} Q_n$ for the canonical surjections. There exists a unique continuous homomorphism $P_n \xto{\pi_n} Q_n$ such that the diagram
\[\begin{tikzcd}
	M & {P_n} \\
	{M/N} & {Q_n}
	\arrow["{\sigma_n}", from=1-1, to=1-2]
	\arrow["\pi"', from=1-1, to=2-1]
	\arrow["{\pi_n}", from=1-2, to=2-2]
	\arrow["{\tau_n}"', from=2-1, to=2-2]
\end{tikzcd}\]
is commutative. Note that $\pi_n$ is surjective. Furthermore, for each $n,m \in \N$ with $n \geq m$, there are unique continous homomorphisms $P_n \xto{\mu_{n,m}} P_m$ and $Q_n \xto{\nu_{n,m}} Q_m$ such that the diagrams 
\[\begin{tikzcd}
	M && {M/N} \\
	{P_n} & {P_m} & {Q_n} & {Q_m}
	\arrow["{\sigma_n}"', from=1-1, to=2-1]
	\arrow["{\sigma_m}", from=1-1, to=2-2]
	\arrow["{\tau_n}"', from=1-3, to=2-3]
	\arrow["{\tau_m}", from=1-3, to=2-4]
	\arrow["{\mu_{n,m}}"', from=2-1, to=2-2]
	\arrow["{\nu_{n,m}}"', from=2-3, to=2-4]
\end{tikzcd}\]
are commutative. Note that $\mu_{n.m}, \nu_{n,m}$ are surjective and the diagram
\[\begin{tikzcd}
	{P_n} & {P_m} \\
	{Q_n} & {Q_m}
	\arrow["{\mu_{n,m}}", from=1-1, to=1-2]
	\arrow["{\pi_n}"', from=1-1, to=2-1]
	\arrow["{\pi_m}", from=1-2, to=2-2]
	\arrow["{\nu_{n,m}}"', from=2-1, to=2-2]
\end{tikzcd}\]
is commutative.

By induction, we define a finite subset $F_n \sub P_n$ for each $n \in \N$ as follows. First let $n = 0$. Since $S$ is a compact subset of $M/N$, the set $\tau_0(S)$ is also a compact subset of $Q_0$. Since $Q_0$ is discrete, we conclude that $\tau_0(S)$ is finite. Since the map $\pi_0 : P_0 \to Q_0$ is surjective, there exists a finite subset $F_0$ of $P_0$ such that $\pi_0(F_0) = \tau_0(S)$. Next suppose $n \geq 1$ and that we have chosen a finite subset $F_{n-1} \sub P_{n-1}$. Suppose $a \in \tau_n(S)$ and $x \in F_{n-1}$ satisfy $\pi_{n-1}(x) = \nu_{n,n-1}(a)$. Since the diagram
\[\begin{tikzcd}
	{P_n} & {P_{n-1}} \\
	{Q_n} & {Q_{n-1}}
	\arrow["{\mu_{n,n-1}}", from=1-1, to=1-2]
	\arrow["{\pi_n}"', from=1-1, to=2-1]
	\arrow["{\pi_{n-1}}", from=1-2, to=2-2]
	\arrow["{\nu_{n,n-1}}"', from=2-1, to=2-2]
\end{tikzcd}\]
is a pullback diagram in $\ub{Ab}$, there exists a unique $\xi_n(a,x) \in P_n$ such that $\mu_{n,n-1}(\xi_n(a,x)) = x$ and $\pi_n(\xi_n(a,x)) = a$. We define
\begin{equation}
F_n := \set{\xi_n(a,x) \in P_n}{ a \in \tau_n(S) \: ; \: x \in F_{n-1} \: ; \: \pi_{n-1}(x) = \nu_{n,n-1}(a)}.
\end{equation}
Since $\tau_n(S)$ is a compact subset of the discrete space $Q_n$, $\tau_n(S)$ is finite. $F_{n-1}$ is also finite. Therefore $F_n$ is a finite subset of $P_n$. This completes the definition of $F_n$.

Consider the inverse system $( (P_n)_{n \in \N} , (\mu_{n,m})_{n \geq m} )$ over $\N$ in $\ub{Top}$. Since $M$ is complete Hausdorff, the cone $(M, (\sigma_n)_{n \in \N} )$ is the limit of this inverse system in $\ub{Top}$. On the other hand, the defintion shows that $\mu_{n,m}(F_n) \sub F_m$ for every $n,m \in \N$ with $n \geq m$. Thus we have an inverse system $( (F_n)_{n \in \N} , (F_n \xto{\mu_{n,m}} F_m)_{n \geq m} )$ over $\N$ in $\ub{Top}$. If we define
\begin{equation}
\tilde{S} := \bigcap_{n \in \N} \sigma_n^{-1}(F_n) ,
\end{equation}
Then the cone $( \tilde{S} , (\tilde{S} \xto{\sigma_n} F_n)_{n \in \N} )$ is the limit of the inverse system $( (F_n)_{n \in \N} , (F_n \xto{\mu_{n,m}} F_m)_{n \geq m} )$ in $\ub{Top}$.

In particular, $\tilde{S}$ is a countable limit in $\ub{Top}$ of the finite discrete spaces $F_n$ $(n \in \N)$. Therefore $\tilde{S}$ is a light profinite set.

Next we prove that $\pi(\tilde{S}) \sub S$. If $\eta \in \tilde{S}$, then $\tau_n( \pi(\eta) ) = \pi_n( \sigma_n(\eta) ) \in \pi_n(F_n) \sub \tau_n(S)$ for every $n \in \N$. Therefore $\pi(\eta) \in \bigcap_{n \in \N} (\pi(U_n) + S)$. However, the set $\set{\pi(U_n)}{n \in \N}$ is a fundamental system of neighbourhoods of $0$ in $M/N$. Therefore $\bigcap_{n \in \N} (\pi(U_n) + S)$ is equal to the closure of $S$ in $M/N$. Moreover, since $N$ is a closed subgroup of $M$, $M/N$ is Hausdorff. Therefore the compact subspace $S$ of $M/N$ is closed in $M/N$. Accordingly, we have $\bigcap_{n \in \N} (\pi(U_n) + S) = S$. Therefore $\pi(\eta) \in \bigcap_{n \in \N} (\pi(U_n) + S) = S$. This shows that $\pi(\tilde{S}) \sub S$.

Next we prove that $S \sub \pi(\tilde{S})$. Suppose $\alpha \in S$. We construct an $\eta \in \tilde{S}$ such that $\pi(\eta) = \alpha$. We define $\eta_n \in F_n$ for $n \in \N$ inductively as follows. Fisrt suppose $n=0$. Since $\pi_0(F_0) = \tau_0(S)$, there exists an $\eta_0 \in F_0$ such that $\pi_0(\eta_0) = \tau_0(\alpha)$. Next suppose $n \geq 1$ and that we have chosen $\eta_{n-1} \in F_{n-1}$ in such a way that $\pi_{n-1}(\eta_{n-1}) = \tau_{n-1}(\alpha)$. Since $\pi_{n-1}(\eta_{n-1}) = \tau_{n-1}(\alpha) = \nu_{n,n-1}(\tau_n(\alpha))$, we can define $\eta_n := \xi_n(\tau_n(\alpha), \eta_{n-1})$. Then we have $\eta_n \in F_n$ and $\pi_n(\eta_n) = \tau_n(\alpha)$. This completes the defintion of $\eta_n$.

For each $n \in \N$, we have $\mu_{n+1,n}(\eta_{n+1}) = \eta_n$ by definition. Therefore $\mu_{n,m}(\eta_n) = \eta_m$ for every $n,m \in \N$ with $n \geq m$. Since the cone $( \tilde{S} , (\tilde{S} \xto{\sigma_n} F_n)_{n \in \N} )$ is the limit of the inverse system $( (F_n)_{n \in \N} , (F_n \xto{\mu_{n,m}} F_m)_{n \geq m} )$ in $\ub{Top}$, there exists a unique $\eta \in \tilde{S}$ such that $\sigma_n(\eta) = \eta_n$ for every $n \in \N$. We prove $\pi(\eta) = \alpha$. For every $n \in \N$, we have $\tau_n( \pi(\eta) )= \pi_n( \sigma_n(\eta) ) = \pi_n(\eta_n) = \tau_n(\alpha)$. Therefore $\pi(\eta) - \alpha \in \pi(U_n)$ for every $n \in \N$. On the other hand, the set $\set{\pi(U_n)}{n \in \N}$ is a fundamental system of neighbourhoods of $0$ in $M/N$. Since $N$ is a closed subgroup of $M$, $M/N$ is Hausdorff. Therefore $\bigcap_{n \in \N} \pi(U_n) = 0$. Thus we conclude that $\pi(\eta) - \alpha = 0$, and hence $\pi(\eta) = \alpha$.
\end{proof}

\begin{cor} \label{cor:compact subsets are light profinite sets}
Let $M$ be a first countable complete Hausdorff non-archimedean abelian group. Then every compact subspace $S$ of $M$ is a light profinite set.
\end{cor}

\begin{proof}
Consider the closed subgroup $\{0\}$ of $M$ and apply \cref{prop:lifting compact sets}.
\end{proof}

\begin{cor} \label{cor:underbar of quotient map is still epic}
Let $M$ be a first countable complete Hausdorff non-archimedean abelian group. Let $N$ be a closed subgroup of $M$.
\begin{enumerate}
\item $\tu{M} \to \tu{M/N}$ is an epimorphism when viewed as a morphism in $\ub{CSet}$. 

\item The sequence
\[\begin{tikzcd}
	0 & {\tu{N}} & {\tu{M}} & {\tu{M/N}} & 0
	\arrow[from=1-1, to=1-2]
	\arrow[from=1-2, to=1-3]
	\arrow[from=1-3, to=1-4]
	\arrow[from=1-4, to=1-5]
\end{tikzcd}\]
is exact in $\ub{CAb}$.
\end{enumerate}
\end{cor}

\begin{proof}~
\begin{enumerate}
\item
Let $M \xto{\pi} M/N$ be the canonical surjection. Suppose that $S$ is a light profinite set and that $\phi \in \tu{M/N}(S)$. Then $\phi(S)$ is a compact subset of $M/N$. By \cref{prop:lifting compact sets}, there exists a subspace $\tilde{S} \sub M$ which is a light profinite set and satisfies $\pi(\tilde{S}) = \phi(S)$. Let $\tilde{S} \xto{j} M$ be the inclusion. We define $S' := S \times_{\phi(S)} \tilde{S}$, where the fiber product is taken in $\ub{Top}$.

Since $M/N$ is a first countable complete Hausdorff non-archimedean abelian group, \cref{cor:compact subsets are light profinite sets} shows that $\phi(S)$ is a light profinite set. Moreover, $S$ and $\tilde{S}$ are also light profinite sets. Then \cref{prop:closure properties of light profinite sets} shows that $S'$ is also a light profinte set.

Let us write $S' \xto{p} S$ and $S' \xto{q} \tilde{S}$ for the projections. The following diagram is commutative, and the left square is a pullback in $\ub{Top}$.
\[\begin{tikzcd}
	{S'} & {\tilde{S}} & M \\
	S & {\phi(S)} & {M/N}
	\arrow["q", from=1-1, to=1-2]
	\arrow["p"', from=1-1, to=2-1]
	\arrow["j", hook, from=1-2, to=1-3]
	\arrow["\pi", two heads, from=1-2, to=2-2]
	\arrow["\pi", two heads, from=1-3, to=2-3]
	\arrow["\phi"', from=2-1, to=2-2]
	\arrow["\inc"', hook, from=2-2, to=2-3]
\end{tikzcd}\]
Since $\tilde{S} \xto{\pi} \phi(S)$ is surjective, the map $S' \xto{p} S$ is also surjective. Moreover we have $j \of q \in \tu{M}(S')$ and $\tu{\pi}_{S'}(j \of q) = \pi \of j \of q = \phi \of p = \tu{M/N}(p)(\phi)$. 

By \cref{prop:monics and epics in CSet}, we conclude that $\tu{\pi} : \tu{M} \to \tu{M/N}$ is an epimorphism when viewed as a morphism in $\ub{CSet}$. 

\item
By definition, it is immediate that the sequence
\[\begin{tikzcd}
	0 & {\tu{N}} & {\tu{M}} & {\tu{M/N}}
	\arrow[from=1-1, to=1-2]
	\arrow[from=1-2, to=1-3]
	\arrow[from=1-3, to=1-4]
\end{tikzcd}\]
is exact in $\ub{CAb}$. Moreover, (1) shows that $\tu{M} \to \tu{M/N}$ is an epimorphism in $\ub{CAb}$.
\end{enumerate}
\end{proof}

\subsection{Restricted power series}

\subsubsection{Settings}

$\\[2mm]$ Throughout this subsection, we fix the following notation.

\begin{nt} \;
\begin{enumerate}
\item $R$ denotes a complete Hausdorff non-archimedean ring.

\item Let $n \in \N$.
\begin{enumerate}
\item $\tu{R}[X_1 , \ldots , X_n]$ denotes the polynomial ring in $n$ variables over $\tu{R}$. $\iota_n : \tu{R} \mon \tu{R}[X_1 , \ldots , X_n]$ denotes the canonical inclusion.

\item $R \langle X_1 , \ldots , X_n \rangle$ denotes the ring of restricted formal power series in $n$ variables over $R$. $j_n : R \mon R \langle X_1 , \ldots , X_n \rangle$ denotes the canonical inclusion.

\item $\phi_n : \tu{R}[X_1, \ldots , X_n] \to \tu{R \langle X_1 , \ldots , X_n \rangle}$ denotes the unique homomorphism of condensed rings such that $\phi_{n,*}(X_i) = X_i$ for $1 \leq i \leq n$ and the following diagram is commutative.
\[\begin{tikzcd}
	{\tu{R}} & {\tu{R}[X_1, \ldots , X_n]} \\
	{\tu{R \langle X_1 , \ldots , X_n \rangle}}
	\arrow["{\iota_n}", hook, from=1-1, to=1-2]
	\arrow["{\tu{j_n}}"', from=1-1, to=2-1]
	\arrow["{\phi_n}", from=1-2, to=2-1]
\end{tikzcd}\]

\item We consider $\tu{R}[X_1, \ldots , X_n]$ and $\tu{R \langle X_1 , \ldots , X_n \rangle}$ as condensed $\tu{R}[X_1 , \ldots , X_n]$-algebras via the identity $\tu{R}[X_1, \ldots , X_n] \to \tu{R}[X_1, \ldots , X_n]$ and $\phi_n : \tu{R}[X_1, \ldots , X_n] \to \tu{R \langle X_1 , \ldots , X_n \rangle}$ respectively.
\end{enumerate}
\end{enumerate}
\end{nt}

\subsubsection{Statement}

\begin{prop} \label{prop:coalescence of polynomial rings}
For every $n \in \N$, we have
\begin{equation}
\tu{R}[X_1 , \ldots , X_n]_{\approx \{ (X_1,1) , \ldots , (X_n,1) \}}
= \tu{R \langle X_1 , \ldots , X_n \rangle}.
\end{equation}
More precisely, $\phi_n : \tu{R}[X_1, \ldots , X_n] \to \tu{R \langle X_1 , \ldots , X_n \rangle}$ is the $\{ (X_1,1) , \ldots , (X_n,1) \}$-coalescence of $\tu{R}[X_1, \ldots , X_n]$ as a condensed $\tu{R}[X_1, \ldots , X_n]$-algebra.
\end{prop}

The proof of this proposition is the content of this subsection.

\subsubsection{The one variable case}

\begin{nt} For simplicity, we use the following notation.
\begin{enumerate}
\item $\tu{R}[X]$ denotes the polynomial ring in one variable over $\tu{R}$. $\iota : \tu{R} \mon \tu{R}[X]$ denotes the canonical inclusion.

\item $R \langle X \rangle$ denotes the ring of restricted formal power series in one variable over $R$. $j : R \mon R \langle X \rangle$ denotes the canonical inclusion.

\item $\phi : \tu{R}[X] \to \tu{R \langle X \rangle}$ denotes the unique homomorphism of condensed rings such that $\phi_{*}(X) = X$ and the following diagram is commutative.
\[\begin{tikzcd}
	{\tu{R}} & {\tu{R}[X]} \\
	{\tu{R \langle X \rangle}}
	\arrow["\iota", hook, from=1-1, to=1-2]
	\arrow["{\tu{j}}"', from=1-1, to=2-1]
	\arrow["\phi", from=1-2, to=2-1]
\end{tikzcd}\]

\item We consider $\tu{R}[X]$ and $\tu{R \langle X \rangle}$ as condensed $\tu{R}[X]$-algebras via the identity $\tu{R}[X] \to \tu{R}[X]$ and $\phi : \tu{R}[X] \to \tu{R \langle X \rangle}$ respectively.

\item For $F \in R \langle X \rangle$ and $n \in \N$, we write $F_n$ for the coefficient of $n$-th degree in $F$.
\end{enumerate}
\end{nt}

In addtion to these notation, we use \cref{nt:N infty,nt:translation and point map of N infty,nt:relative translation and point map of N infty}.

\begin{lem} \label{lem:shift map is bijective}
Let $M$ be a condensed set. Let $S$ be a light profinite set. Then the map
\[\begin{tikzcd}
	{M(S \times \N_{\infty})} & {M(S) \times M(S \times \N_{\infty}),} & t & {\big( M(i_{S,0})(t) , M(\tau_S)(t) \big)}
	\arrow[from=1-1, to=1-2]
	\arrow[maps to, from=1-3, to=1-4]
\end{tikzcd}\]
is bijective.
\end{lem}

\begin{proof}
In $\ub{Prof}$, the maps
\[\begin{tikzcd}[row sep=small]
	S \\
	& {S \times \N_{\infty}} \\
	{S \times \N_{\infty}}
	\arrow["{i_{S,0}}", from=1-1, to=2-2]
	\arrow["{\tau_S}"', from=3-1, to=2-2]
\end{tikzcd}\]
induces a homeomorphism
\[\begin{tikzcd}
	{S \sqcup (S \times \N_{\infty})} & {S \times \N_{\infty}.}
	\arrow["\sim", from=1-1, to=1-2]
\end{tikzcd}\]
Since $M$ is a sheaf on the site $\ub{Prof}$, the maps
\[\begin{tikzcd}[row sep=small]
	{M(S)} \\
	& {M(S \times \N_{\infty})} \\
	{M(S \times \N_{\infty})}
	\arrow["{M(i_{S,0})}"', from=2-2, to=1-1]
	\arrow["{M(\tau_S)}", from=2-2, to=3-1]
\end{tikzcd}\]
induce a bijection
\[\begin{tikzcd}
	{M(S \times \N_{\infty})} & {M(S) \times M(S \times \N_{\infty}),} & t & {\big( M(i_{S,0})(t) , M(\tau_S)(t) \big) .}
	\arrow[from=1-1, to=1-2]
	\arrow[maps to, from=1-3, to=1-4]
\end{tikzcd}\]
\end{proof}

\begin{lem} \label{lem:definition of Gamma S}
Let $S$ be a light profinite set.
\begin{enumerate}
\item
For each $f \in \tu{R \langle X \rangle}(S)$, define $\Gamma_S (f) : S \times \N_{\infty} \to R$ as follows. 
\[\begin{tikzcd}[ampersand replacement=\&]
	{S \times \N_{\infty}} \& {R,} \& {(x,n)} \& {\left\{
	\begin{aligned}
	& f(x)_n & (n \in \N) \\
	& 0 & (n = \infty)
	\end{aligned}
	\right.}
	\arrow["{\Gamma_S(f)}", from=1-1, to=1-2]
	\arrow[maps to, from=1-3, to=1-4]
\end{tikzcd}\]
Then $\Gamma_S (f) \in \ker \tu{R} (i_{S,\infty})$.

\item The map
\[\begin{tikzcd}
	{\tu{R \langle X \rangle}(S)} & {\ker \tu{R}(i_{S,\infty}),} & f & {\Gamma_S (f)}
	\arrow["{\Gamma_S}", from=1-1, to=1-2]
	\arrow[maps to, from=1-3, to=1-4]
\end{tikzcd}\]
is bijective.
\end{enumerate}
\end{lem}

\begin{proof}~
\begin{enumerate}
\item
First we prove that $\Gamma_S (f) \in \tu{R}(S \times \N_{\infty})$, i.e., the map $\Gamma_S (f) : S \times \N_{\infty} \to R$ is continuous. Let $(x_0,n_0) \in S \times \N_{\infty}$. We show that the map $\Gamma_S (f) : S \times \N_{\infty} \to R$ is continuous at $(x_0,n_0)$. 

First suppose that $n_0 \in \N$. Then the map $S \xto{i_{S,n_0}} S \times \N \xto{\Gamma_S (f)} R$ is of the form $x \mapsto f(x)_{n_0}$. Thus it is equal to the composition of the continuous map $S \xto{f} R \langle X \rangle$ and the continuous map $R \langle X \rangle \to R$, $F \mapsto F_{n_0}$. Consequently, the map $S \xto{i_{S,n_0}} S \times \N \xto{\Gamma_S (f)} R$ is continuous. Since $i_{S,n_0} : S \to S \times \N_{\infty}$ defines a homeomorphism of $S$ onto the open neighbourhood $S \times \{n_0\}$ of $(x_0,n_0)$ in $S \times \N_{\infty}$, it follows that $\Gamma_S (f) : S \times \N_{\infty} \to R$ is continuous at $(x_0,n_0)$. 

Next suppose that $n_0 = \infty$. Let $U$ be an arbitrary neighbourhood of $0$ in $R$ which is an additive subgroup of $R$. We show that there exists a neighbourhood $W$ of $(x_0,n_0)$ in $S \times \N_{\infty}$ such that $\Gamma_S (f) (W) \sub U$. Since the map $f : S \to R \langle X \rangle$ is continous, there exists a neighbourhood $G$ of $x_0$ in $S$ such that $f(x)_n - f(x_0)_n \in U$ for every $n \in \N$ and every $x \in G$. On the other hand, since $f(x_0) \in R \langle X \rangle$ is a restricted power series, there exists an $N \in \N$ such that $f(x_0)_n \in U$ for every $n \geq N$. Then $V := \set{n \in \N}{n \geq N} \cup \{ \infty \}$ is a neighbourhood of $\infty = n_0$ in $\N_{\infty}$. Then $W := G \times V$ is a neighbourhood of $(x_0,n_0)$ in $S \times \N_{\infty}$. We show that $\Gamma_S (f) (x,n) \in U$ for every $(x,n) \in W$. Let $(x,n) \in W$. If $n = \infty$, then $\Gamma_S (f) (x,n) = 0 \in U$. Suppose $n \in \N$. Then $x \in G$ and $n \geq N$. By the choice of $G$ and $N$, we have
\begin{equation}
f(x)_n - f(x_0)_n \in U \quad \text{ and } \quad f(x_0)_n \in U .
\end{equation}  
Then
\begin{equation}
\Gamma_S(f)(x,n) = f(x)_n = \big( f(x)_n - f(x_0)_n \big) + f(x_0)_n \, \in U ,
\end{equation}
since $U$ is an additive subgroup of $R$. This completes the proof that the map $\Gamma_S (f) : S \times \N_{\infty} \to R$ is continuous. 

Thus we have proved that $\Gamma_S (f) \in \tu{R}(S \times \N_{\infty})$. Then $\tu{R}(i_{S,\infty}) ( \Gamma_S (f) ) = \Gamma_S (f) \of i_{S,\infty} = 0$ by definition. Therefore $\Gamma_S (f) \in \ker \tu{R}(i_{S,\infty})$.

\item
First we prove the injectivity of the map $\Gamma_S$. If $f,f' \in \tu{R \langle X \rangle}(S)$ satisfy $\Gamma_S(f) = \Gamma_S(f')$, then
\begin{equation}
f(x)_n = \Gamma_S(f)(x,n) = \Gamma_S(f')(x,n) = f'(x)_n
\end{equation}
for every $n \in \N$ and every $x \in S$. It follows that $f(x) = f'(x)$ for every $x \in S$. Therefore $f=f'$. This proves the injectivity of the map $\Gamma_S$. 

Next we prove the surjectivity of the map $\Gamma_S$. Suppose $g \in \ker \tu{R}(i_{S,\infty})$. We show that there exists an $f \in \tu{R \langle X \rangle}(S)$ such that $\Gamma_S(f) = g$. $g : S \times \N_{\infty} \to R$ is a continuous map and satisfies $g \of i_{S,\infty} = 0$. Therefore if $x \in S$, then $( g(x,n) )_{n \in \N}$ is a zero sequence in $R$. Hence 
\begin{equation}
f(x) := \sum_{n \in \N} g(x,n) \cdot X^n
\end{equation}
is an element of $R \langle X \rangle$. Thus we obtain a map $f : S \to R \langle X \rangle$, $x \mapsto f(x)$. We claim that this map is continuous. Let $x_0 \in S$. We show that the map $f : S \to R \langle X \rangle$ is continuous at $x_0$. Let $U$ be an arbitrary neighbourhood of $0$ in $R$ which is an additive subgroup of $R$. We prove that there exists a neighbourhood $H$ of $x_0$ in $S$ such that $f(x)_n - f(x_0)_n \in U$ for every $n \in \N$ and every $x \in H$.

Since the map $g : S \times \N_{\infty} \to R$ is continuous at $(x_0,\infty)$ and since $g(x_0,\infty) = 0$, there exist a neighbourhood $G$ of $x_0$ in $S$ and an $N \in \N$ such that $g(x,n) \in U$ for every $x \in G$ and every $n \in \N$ with $n \geq N$. For each $n \in \N$ with $0 \leq n < N$, the map $g : S \times \N_{\infty} \to R$ is continuous at $(x_0,n)$. Therefore there exists a neighbourhood $G_n$ of $x_0$ in $S$ such that $g(x,n) - g(x_0,n) \in U$ for every $x \in G_n$. Then $H := G \cap G_0 \cap \cdots \cap G_{N-1}$ is a neighbourhood of $x_0$ in $S$. We prove that $f(x)_n - f(x_0)_n \in U$ for every $n \in \N$ and every $x \in H$. Let $x \in H$ and $n \in \N$. If $0 \leq n < N$, then $x \in G_n$. By the choice of $G_n$, we have
\begin{equation}
f(x)_n - f(x_0)_n = g(x,n) - g(x_0,n) \, \in U .
\end{equation}
Suppose $n \geq N$. Since $x, x_0 \in G$ and $n \geq N$, we have $g(x,n), g(x_0,n) \in U$ by the choice of $G$ and $N$. Then
\begin{equation}
f(x)_n - f(x_0)_n = g(x,n) - g(x_0,n) \, \in U ,
\end{equation}
since $U$ is an additive subgroup of $U$. This completes the proof that the map $f : S \to R \langle X \rangle$ is continuous.

Thus $f \in \tu{R \langle X \rangle}(S)$. For every $x \in S$ and every $n \in \N$, we have $\Gamma_S(f)(x,n) = f(x)_n = g(x,n)$ by definition. Moreover, for every $x \in S$, we have $\Gamma_S(f)(x,\infty) = 0 = g(x,\infty)$ since $g \of i_{S,\infty} = 0$. It follows that $\Gamma_S(f) = g$. This completes the proof.
\end{enumerate}
\end{proof}

\begin{lem} \label{lem:functoriality of Gamma S}
If $S' \xto{c} S$ is any morphism in $\ub{Prof}$, then the diagram
\[\begin{tikzcd}
	{\tu{R \langle X \rangle}(S)} & {\ker \tu{R}(i_{S,\infty})} \\
	{\tu{R \langle X \rangle}(S')} & {\ker \tu{R}(i_{S',\infty})}
	\arrow["{\Gamma_S}", from=1-1, to=1-2]
	\arrow["{\tu{R \langle X \rangle}(c)}"', from=1-1, to=2-1]
	\arrow["{\tu{R \langle X \rangle}(c \times \id_{\N_{\infty}})}", from=1-2, to=2-2]
	\arrow["{\Gamma_{S'}}"', from=2-1, to=2-2]
\end{tikzcd}\]
is commutative.
\end{lem}

\begin{proof}
Let $f \in \tu{R \langle X \rangle}(S)$. We prove that
\begin{equation}
\big( \tu{R \langle X \rangle}(c \times \id_{\N_{\infty}}) \of \Gamma_S \big) (f)
= \big( \Gamma_{S'} \of \tu{R \langle X \rangle} (c) \big) (f) .
\end{equation}
For every $x \in S$ and $n \in \N$, we have
\begin{align}
\Big( \big( \tu{R \langle X \rangle}(c \times \id_{\N_{\infty}}) \of \Gamma_S \big) (f) \Big) (x,n)
& = \Big( \big( \Gamma_S (f) \big) \of \big( c \times \id_{\N_{\infty}} \big) \Big) (x,n) \\
& = \Gamma_S (f) (c(x),n) = f(c(x))_n \: ; \\
\Big( \big( \Gamma_{S'} \of \tu{R \langle X \rangle} (c) \big) (f) \Big) (x,n)
& = \Big( \Gamma_{S'} \big( f \of c \big) \Big) (x,n) = f(c(x))_n . 
\end{align}
Thus
\begin{equation}
\Big( \tu{R \langle X \rangle}(c \times \id_{\N_{\infty}}) \of \Gamma_S \, (f) \Big) (x,n)
= f(c(x))_n
= \Big( \Gamma_{S'} \of \tu{R \langle X \rangle} (c) \, (f) \Big) (x,n) .
\end{equation}
For every $x \in S$, we have
\begin{equation}
\Big( \big( \tu{R \langle X \rangle}(c \times \id_{\N_{\infty}}) \of \Gamma_S \big) (f) \Big) (x,\infty)
= 0 = 
\Big( \big( \Gamma_{S'} \of \tu{R \langle X \rangle} (c) \big) (f) \Big) (x,\infty) .
\end{equation}
It follows that
\begin{equation}
\big( \tu{R \langle X \rangle}(c \times \id_{\N_{\infty}}) \of \Gamma_S \big) (f)
= \big( \Gamma_{S'} \of \tu{R \langle X \rangle} (c) \big) (f) .
\end{equation}
\end{proof}

\begin{lem} \label{lem:definition of Theta S}
Let $S$ be a light profinite set. 
\begin{enumerate}
\item
For each $f \in \tu{R \langle X \rangle}(S)$, define $\Theta_S (f) : S \times \N_{\infty} \to R \langle X \rangle$ as follows.
\[\begin{tikzcd}[ampersand replacement=\&]
	{S \times \N_{\infty}} \& {R \langle X \rangle,} \& {(x,n)} \& {\left\{
	\begin{aligned}
	& \sum_{k \in \N} f(x)_{k+n} X^k & (n \in \N) \\
	& 0 & (n = \infty)
	\end{aligned}
	\right.}
	\arrow["{\Theta_S(f)}", from=1-1, to=1-2]
	\arrow[maps to, from=1-3, to=1-4]
\end{tikzcd}\]
Then $\Theta_S (f) \in \ker \tu{R \langle X \rangle} (i_{S,\infty})$. Consequently, we have a well-defined map
\[\begin{tikzcd}
	{\tu{R \langle X \rangle}(S)} & {\ker \tu{R \langle X \rangle}(i_{S,\infty}),} & f & {\Theta_S(f).}
	\arrow["{\Theta_S}", from=1-1, to=1-2]
	\arrow[maps to, from=1-3, to=1-4]
\end{tikzcd}\]

\item
The following diagram is commutative.
\[\begin{tikzcd}
	& {\tu{R \langle X \rangle}(S)} \\
	{\tu{R \langle X \rangle}(S)} & {\ker \tu{R \langle X \rangle}(i_{S,\infty})}
	\arrow["\id", from=2-1, to=1-2]
	\arrow["{\Theta_S}"', from=2-1, to=2-2]
	\arrow["{\tu{R \langle X \rangle}(i_{S,0})}"', from=2-2, to=1-2]
\end{tikzcd}\]
\end{enumerate}
\end{lem}

\begin{proof}~
\begin{enumerate}
\item
Let $f \in \tu{R \langle X \rangle}(S)$. First we prove that the map $\Theta_S (f) : S \times \N_{\infty} \to R \langle X \rangle$ is continuous. Let $(x_0,n_0) \in S \times \N_{\infty}$. We show that the map $\Theta_S (f) : S \times \N_{\infty} \to R \langle X \rangle$ is continuous at $(x_0,n_0)$. Let $U$ be an arbitrary neighbourhood of $0$ in $R$ which is an additive subgroup of $R$. We prove that there exists a neighbourhood $W$ of $(x_0,n_0)$ in $S \times \N_{\infty}$ such that $(\Theta_S(f)(x,n))_m - (\Theta_S(f)(x_0,n_0))_m \in U$ for every $m \in \N$ and every $(x,n) \in W$.

First suppose that $n_0 \in \N$. Since the map $f : S \to R \langle X \rangle$ is continous at $x_0$, there exists a neighbourhood $G$ of $x_0$ in $S$ such that $f(x)_m - f(x_0)_m \in U$ for every $m \in \N$ and every $x \in G$. Then $W := G \times \{ n_0 \}$ is a neighbourhood of $(x_0,n_0)$ in $S \times \N_{\infty}$. We prove that $(\Theta_S(f)(x,n))_m - (\Theta_S(f)(x_0,n_0))_m \in U$ for every $m \in \N$ and every $(x,n) \in W$. Let $(x,n) \in W$ and $m \in \N$. Then $x \in G$ and $n=n_0$. By the choice of $G$, we have $f(x)_{m+n_0} - f(x_0)_{m+n_0} \in U$. Therefore
\begin{align}
(\Theta_S(f)(x,n))_m - (\Theta_S(f)(x_0,n_0))_m
& = f(x)_{m+n} - f(x_0)_{m+n_0} \\
& = f(x)_{m+n_0} - f(x_0)_{m+n_0} \, \in U .
\end{align}

Next suppose that $n_0 = \infty$. Since the map $f : S \to R \langle X \rangle$ is continous at $x_0$, there exists a neighbourhood $G$ of $x_0$ in $S$ such that $f(x)_m - f(x_0)_m \in U$ for every $m \in \N$ and every $x \in G$. On the other hand, since $f(x_0) \in R \langle X \rangle$ is a restricted power series, there exists an $N \in \N$ such that $f(x_0)_m \in U$ for every $m \in \N$ with $m \geq N$. Then $V := \set{n \in \N}{n \geq N} \cup \{ \infty \}$ is a neighbourhood of $\infty = n_0$ in $\N_{\infty}$. Then $W := G \times V$ is a neighbourhood of $(x_0,n_0)$ in $S \times \N_{\infty}$. We prove that $(\Theta_S(f)(x,n))_m - (\Theta_S(f)(x_0,n_0))_m \in U$ for every $m \in \N$ and every $(x,n) \in W$. Since $n_0 = \infty$, we have $(\Theta_S(f)(x_0,n_0))_m = (\Theta_S(f)(x_0,\infty))_m = 0$ for every $m \in \N$. Therefore it suffices to show that $(\Theta_S(f)(x,n))_m \in U$ for every $m \in \N$ and every $(x,n) \in W$. Let $(x,n) \in W$ and $m \in \N$. If $n = \infty$, then $(\Theta_S(f)(x,n))_m = 0 \in U$. Suppose $n \in \N$. Then $n \geq N$. Therefore $m+n \geq N$. By the choice of $N$, we have $f(x_0)_{m+n} \in U$. On the other hand, we have $x \in G$. By the choice of $G$, we have $f(x)_{m+n} - f(x_0)_{m+n} \in U$. It follows that
\begin{equation}
(\Theta_S(f)(x,n))_m = f(x)_{m+n} = \big( f(x)_{m+n} - f(x_0)_{m+n} \big) + f(x_0)_{m+n} \, \in U ,
\end{equation}
since $U$ is an additive subgroup of $R$. This completes the proof that the map $\Theta_S (f) : S \times \N_{\infty} \to R \langle X \rangle$ is continuous. 

Thus $\Theta_S (f) \in \tu{R \langle X \rangle}(S \times \N_{\infty})$. Then we have $\tu{R \langle X \rangle}(i_{S,\infty}) (\Theta_S (f)) = \Theta_S (f) \of i_{S,\infty} = 0$ by definition. Therefore $\Theta_S (f) \in \ker \tu{R \langle X \rangle}(i_{S,\infty})$. This completes the proof of the assertion (1).

\item
Let $f \in \tu{R \langle X \rangle}(S)$. We prove that
\begin{equation}
\big( \tu{R \langle X \rangle}(i_{S,0}) \of \Theta_S \big) (f) = f .
\end{equation}
For every $x \in S$, we have
\begin{align}
\Big( \big( \tu{R \langle X \rangle}(i_{S,0}) \of \Theta_S \big) (f) \Big) (x)
& = \Big( \big( \Theta_S (f) \big) \of i_{S,0} \Big) (x) \\
& = \big( \Theta_S (f) \big) (x,0) \\
& = \sum_{k \in \N} f(x)_k X^k = f(x) .
\end{align}
Therefore we have
\begin{equation}
\big( \tu{R \langle X \rangle}(i_{S,0}) \of \Theta_S \big) (f) = f .
\end{equation}
This completes the proof.
\end{enumerate}
\end{proof}

\begin{lem} \label{lem:commutativity of Gamma S and Theta S}
For every light profinite set $S$, the following diagram is commutative.
\[\begin{tikzcd}
	{\tu{R \langle X \rangle}(S)} & {\ker \tu{R \langle X \rangle}(i_{S,\infty})} \\
	{\ker \tu{R}(i_{S,\infty})} & {\ker \tu{R \langle X \rangle}(i_{S,\infty})}
	\arrow["{\Theta_S}", from=1-1, to=1-2]
	\arrow["{\Gamma_S}"', from=1-1, to=2-1]
	\arrow["{\id - X \cdot (\tu{R \langle X \rangle}(\tau_S))}", from=1-2, to=2-2]
	\arrow["{\tu{j}_{S \times \N_{\infty}}}"', from=2-1, to=2-2]
\end{tikzcd}\]
\end{lem}

\begin{proof}
Let $f \in \tu{R \langle X \rangle}(S)$. We prove that
\begin{equation}
\big( \tu{j}_{S \times \N_{\infty}} \of \Gamma_S \big) (f)
= \Theta_S (f) - X \cdot \big( \tu{R \langle X \rangle}(\tau_S) \of \Theta_S \big) (f) .
\end{equation}
For every $x \in S$ and every $n \in \N$, we have
\begin{align}
\Big( \big( \tu{j}_{S \times \N_{\infty}} \of \Gamma_S \big) (f) \Big) (x,n)
& = \Big( j \of \big( \Gamma_S (f) \big) \Big) (x,n) \\
& = j \big( \Gamma_S (f) (x,n) \big) \\
& = j \big( f(x)_n \big) \\
& = f(x)_n \: ; \\
\Theta_S (f)(x,n) - X \cdot \Big( \big( \tu{R \langle X \rangle}(\tau_S) \of \Theta_S \big) (f) \Big) (x,n) 
& = \Theta_S (f)(x,n) - X \cdot \Big( \big( \Theta_S (f) \big) \of \tau_S \Big) (x,n) \\
& = \Theta_S (f)(x,n) - X \cdot \Theta_S (f) (x,n+1) \\
& = \left( \sum_{k \in \N} f(x)_{k+n} X^k \right) - X \cdot \left( \sum_{k \in \N} f(x)_{k+n+1} X^k \right) \\
& = \left( \sum_{k \in \N} f(x)_{k+n} X^k \right) - \left( \sum_{\substack{k \in \N \\ k \geq 1}} f(x)_{k+n} X^k \right) \\
& = f(x)_n .
\end{align}
Therefore
\begin{equation}
\Big( \big( \tu{j}_{S \times \N_{\infty}} \of \Gamma_S \big) (f) \Big) (x,n)
= f(x)_n
= \Theta_S (f)(x,n) - X \cdot \Big( \big( \tu{R \langle X \rangle}(\tau_S) \of \Theta_S \big) (f) \Big) (x,n) .
\end{equation}
For every $x \in S$, we have
\begin{align}
\Big( \big( \tu{j}_{S \times \N_{\infty}} \of \Gamma_S \big) (f) \Big) (x,\infty)
& = \Big( j \of \big( \Gamma_S (f) \big) \Big) (x,\infty) \\
& = j \big( \Gamma_S (f) (x,\infty) \big) \\
& = j \big( 0 \big) = 0 \: ; \\
\Theta_S (f)(x,\infty) - X \cdot \Big( \big( \tu{R \langle X \rangle}(\tau_S) \of \Theta_S \big) (f) \Big) (x,\infty) 
& = \Theta_S (f)(x,\infty) - X \cdot \Big( \big( \Theta_S (f) \big) \of \tau_S \Big) (x,\infty) \\
& = \Theta_S (f)(x,\infty) - X \cdot \Theta_S (f) (x,\infty) \\
& = 0 - X \cdot 0 = 0 .
\end{align}
Therefore
\begin{equation}
\Big( \big( \tu{j}_{S \times \N_{\infty}} \of \Gamma_S \big) (f) \Big) (x,\infty)
= 0
= \Theta_S (f)(x,\infty) - X \cdot \Big( \big( \tu{R \langle X \rangle}(\tau_S) \of \Theta_S \big) (f) \Big) (x,\infty) .
\end{equation}
It follows that
\begin{equation}
\big( \tu{j}_{S \times \N_{\infty}} \of \Gamma_S \big) (f)
= \Theta_S (f) - X \cdot \big( \tu{R \langle X \rangle}(\tau_S) \of \Theta_S \big) (f) .
\end{equation}
This completes the proof.
\end{proof}

\begin{prop} \label{prop:coalescence of polynomial rings in one variable}
\begin{equation}
\tu{R}[X]_{\approx \{ (X,1) \}} = \tu{R \langle X \rangle} .
\end{equation}
More precisely, $\phi : \tu{R}[X] \to \tu{R \langle X \rangle}$ is the $\{(X,1)\}$-coalescence of $\tu{R}[X]$ as a condensed $\tu{R}[X]$-algebra.
\end{prop}

\begin{proof}
Since $\phi : \tu{R}[X] \to \tu{R \langle X \rangle}$ is already a homomorphism of condensed $\tu{R}[X]$-algebras, it suffices to show that $\phi : \tu{R}[X] \to \tu{R \langle X \rangle}$ is the $\{(X,1)\}$-coalescence of $\tu{R}[X]$ as a condensed $\tu{R}[X]$-module.

Since $R \langle X \rangle$ is a complete Hausdorff non-archimedean ring and since $X \in R \langle X \rangle$ is a power-bounded element of $R \langle X \rangle$, \cref{prop:power-boundedness and coalescence} shows that $\tu{R \langle X \rangle}$ is a $X/1$-coalescent condensed $\tu{R \langle X \rangle}$-module. Since $\phi : \tu{R}[X] \to \tu{R \langle X \rangle}$ satisfies $\phi_{*}(X) = X$, it follows that $\tu{R \langle X \rangle}$ is a $X/1$-coalescent condensed $\tu{R}[X]$-module.

Let $M$ be any condensed $\tu{R}[X]$-module which is $X/1$-coalescent. Let $\beta : \tu{R}[X] \to M$ be any homomorphism of condensed $\tu{R}[X]$-modules. We prove that there exists a unique homomorphism $\gamma : \tu{R \langle X \rangle} \to M$ of condensed $\tu{R}[X]$-modules such that the followng diagram is commutative.
\[\begin{tikzcd}
	{\tu{R}[X]} & M \\
	{\tu{R \langle X \rangle}}
	\arrow["\beta", from=1-1, to=1-2]
	\arrow["\phi"', from=1-1, to=2-1]
	\arrow["\gamma"', from=2-1, to=1-2]
\end{tikzcd}\]

Let $S$ be any light profinite set. Since $M$ is $(X,1)$-coalescent, the map
\[\begin{tikzcd}
	{\delta_S : \ker M(i_{S,\infty})} && {\ker M(i_{S,\infty})}
	\arrow["{\id - X \cdot M(\tau_S)}", from=1-1, to=1-3]
\end{tikzcd}\]
is bijective. Therefore we can define a map $\gamma_S : \tu{R \langle X \rangle}(S) \to M(S)$ by
\begin{equation}
\gamma_S := M(i_{S,0}) \of {\delta_S}^{-1} \of (\beta \of \iota)_{S \times \N_{\infty}} \of \Gamma_S .
\end{equation}
Namely, we define $\gamma_S$ so that the following diagram is commutative.
\[\begin{tikzcd}
	{\tu{R \langle X \rangle}(S)} && {M(S)} \\
	&& {\ker M(i_{S,\infty})} \\
	{\ker \tu{R}(i_{S,\infty})} & {\ker \tu{R}[X](i_{S,\infty})} & {\ker M(i_{S,\infty})}
	\arrow["{\gamma_S}", from=1-1, to=1-3]
	\arrow["{\Gamma_S}"', from=1-1, to=3-1]
	\arrow["{M(i_{S,0})}"', from=2-3, to=1-3]
	\arrow["{\iota_{S \times \N_{\infty}}}"', from=3-1, to=3-2]
	\arrow["{\beta_{S \times \N_{\infty}}}"', from=3-2, to=3-3]
	\arrow["{{\delta_S}^{-1}}"', from=3-3, to=2-3]
\end{tikzcd}\]

We claim that the family $\gamma := (\gamma_S)_{S \in |\ub{Prof}|}$ is a homomorphism $\tu{R \langle X \rangle} \to M$ of condensed $\tu{R}[X]$-modules.

First of all, \cref{lem:functoriality of Gamma S} shows that $\gamma = (\gamma_S)_{S \in |\ub{Prof}|}$ is a map of condensed sets $\tu{R \langle X \rangle} \to M$. Moreover, it is clear that for every light profinite set $S$, the map $\gamma_S$ is a homomorphism of abelian groups. Therefore it remains to prove that, for every light profinite set $S$, the map $\gamma_S$ commutes with the action of $\tu{R}[X](S)$. Since $\tu{R}[X](S) = \tu{R}(S)[X]$, it suffices to show that $\gamma_S$ commutes with the action of $\tu{R}(S)$ and the action of $X$. 

Let $f \in \tu{R \langle X \rangle}(S)$. First we let $r \in \tu{R}(S)$ and we prove that $\gamma_S \big( (j \of r) \cdot f \big) = r \cdot \gamma_S(f)$. Write $p_S : S \times \N_{\infty} \to S$ for the continuous map $(x,n) \mapsto x$. For $(x,n) \in S \times \N$, we have
\begin{align}
\Gamma_S \big( (j \of r) \cdot f \big) (x,n)
& = \big( (j \of r)(x) \cdot f(x) \big)_n = r(x) \cdot (f(x))_n \\
& = \big( (r \of p_S) (x,n) \big) \cdot \big( \Gamma_S(f)(x,n) \big) = \big( (r \of p_S) \cdot \Gamma_S(f) \big) (x,n).
\end{align}
For $x \in S$, we have 
\begin{equation}
\Gamma_S \big( (j \of r) \cdot f \big) (x,\infty) = 0 = \big( (r \of p_S) \cdot \Gamma_S(f) \big) (x,\infty) .
\end{equation}
Therefore $\Gamma_S \big( (j \of r) \cdot f \big) = (r \of p_S) \cdot (\Gamma_S(f))$. Next consider the map
\[\begin{tikzcd}
	{\ker \tu{R}(i_{S,\infty})} & {\ker \tu{R}[X](i_{S,\infty})} & {\ker M(i_{S,\infty})} & {\ker M(i_{S,\infty}) .}
	\arrow["{\iota_{S \times \N_{\infty}}}", from=1-1, to=1-2]
	\arrow["{\beta_{S \times \N_{\infty}}}", from=1-2, to=1-3]
	\arrow["{{\delta_S}^{-1}}", from=1-3, to=1-4]
\end{tikzcd}\]
This is an $\tu{R}(S \times \N_{\infty})$-linear map. Therefore
\begin{align}
& {\delta_S}^{-1} \of (\beta \of \iota)_{S \times \N_{\infty}} \big( (r \of p_S) \cdot (\Gamma_S(f)) \big) \\
= & (r \of p_S) \cdot \Big( {\delta_S}^{-1} \of (\beta \of \iota)_{S \times \N_{\infty}} \of \Gamma_S (f) \Big) .
\end{align}
Since $M$ is a condensed $\tu{R}$-module, we have
\begin{align}
& M(i_{S,0}) \Big( (r \of p_S) \cdot \big( {\delta_S}^{-1} \of (\beta \of \iota)_{S \times \N_{\infty}} \of \Gamma_S (f) \big) \Big) \\
= & \Big( \tu{R}(i_{S,0})(r \of p_S) \Big) \cdot \Big( M(i_{S,0}) \of {\delta_S}^{-1} \of (\beta \of \iota)_{S \times \N_{\infty}} \of \Gamma_S (f) \Big) \\
= & (r \of p_S \of i_{S,0}) \cdot \gamma_S(f) \\
= & r \cdot \gamma_S(f).
\end{align}
Thus we have
\begin{align}
\gamma_S \big( (j \of r) \cdot f \big)
& = M(i_{S,0}) \of {\delta_S}^{-1} \of (\beta \of \iota)_{S \times \N_{\infty}} \of \Gamma_S \, \big( (j \of r) \cdot f \big) \\
& = M(i_{S,0}) \of {\delta_S}^{-1} \of (\beta \of \iota)_{S \times \N_{\infty}} \big( (r \of p_S) \cdot (\Gamma_S(f)) \big) \\
& = M(i_{S,0}) \Big( (r \of p_S) \cdot \big( {\delta_S}^{-1} \of (\beta \of \iota)_{S \times \N_{\infty}} \of \Gamma_S (f) \big) \Big) \\
& = r \cdot \gamma_S(f) .
\end{align}

Next we prove that $\gamma_S(X \cdot f) = X \cdot \gamma_S(f)$. Consider the following diagrams.
\[\begin{tikzcd}
	{\tu{R}(S \times \N_{\infty})} & {M(S \times \N_{\infty})} & {\tu{R}(S \times \N_{\infty})} & {M(S \times \N_{\infty})} \\
	{\tu{R}(S)} & {M(S)} & {\tu{R}(S \times \N_{\infty})} & {M(S \times \N_{\infty})}
	\arrow["{(\beta \of \iota)_{S \times \N_{\infty}}}", from=1-1, to=1-2]
	\arrow["{\tu{R}(i_{S,0})}"', from=1-1, to=2-1]
	\arrow["{M(i_{S,0})}", from=1-2, to=2-2]
	\arrow["{(\beta \of \iota)_{S \times \N_{\infty}}}", from=1-3, to=1-4]
	\arrow["{\tu{R}(\tau_S)}"', from=1-3, to=2-3]
	\arrow["{M(\tau_S)}", from=1-4, to=2-4]
	\arrow["{(\beta \of \iota)_{S}}"', from=2-1, to=2-2]
	\arrow["{(\beta \of \iota)_{S \times \N_{\infty}}}"', from=2-3, to=2-4]
\end{tikzcd}\]
These diagrams are commutative. Therefore
\begin{align}
M(i_{S,0})\big( (\beta \of \iota)_{S \times \N_{\infty}} \of \Gamma_S \, (X \cdot f) \big)
& = (\beta \of \iota)_S \of \tu{R}(i_{S,0}) \, \big( \Gamma_S(X \cdot f) \big) \: ; \\
M(\tau_S)\big( (\beta \of \iota)_{S \times \N_{\infty}} \of \Gamma_S \, (X \cdot f) \big)
& = (\beta \of \iota)_{S \times \N_{\infty}} \of \tu{R}(\tau_S) \, \big( \Gamma_S(X \cdot f) \big) .
\end{align}
For every $(x,n) \in S \times \N$, we have
\begin{equation}
\Gamma_S (X \cdot f) (x,n) = (X \cdot f(x))_n = \left\{
\begin{aligned}
& 0 & (n = 0) \\
& f(x)_{n-1} & (n \geq 1)
\end{aligned}
\right. .
\end{equation}
For every $x \in S$, we have
\begin{equation}
\Gamma_S (X \cdot f) (x,\infty) = 0 .
\end{equation}
Therefore we have
\begin{align}
\tu{R}(i_{S,0}) \big( \Gamma_S(X \cdot f) \big)
& = \Gamma_S(X \cdot f) \of i_{S,0} = 0 \: ; \\
\tu{R}(\tau_S) \big( \Gamma_S(X \cdot f) \big)
& = \Gamma_S(X \cdot f) \of \tau_S = \Gamma_S(f) .
\end{align}
Hence
\begin{align}
M(i_{S,0})\big( (\beta \of \iota)_{S \times \N_{\infty}} \of \Gamma_S \, (X \cdot f) \big)
& = (\beta \of \iota)_S \of \tu{R}(i_{S,0}) \, \big( \Gamma_S(X \cdot f) \big) \\
& = (\beta \of \iota)_S \, (0) = 0 \: ; \\
M(\tau_S)\big( (\beta \of \iota)_{S \times \N_{\infty}} \of \Gamma_S \, (X \cdot f) \big)
& = (\beta \of \iota)_{S \times \N_{\infty}} \of \tu{R}(\tau_S) \, \big( \Gamma_S(X \cdot f) \big) \\
& = (\beta \of \iota)_{S \times \N_{\infty}} \of \Gamma_S \, (f) .
\end{align}
On the other hand, let us write
\begin{equation}
\xi := {\delta_S}^{-1} \of (\beta \of \iota)_{S \times \N_{\infty}} \of \Gamma_S \, (f) \: \in \ker M(i_{S,\infty}) .
\end{equation}
By \cref{lem:shift map is bijective}, there exists a unique $\eta \in M(S \times \N_{\infty})$ such that $M(i_{S,0})(\eta) = X \cdot \gamma_S(f)$ and $M(\tau_S)(\eta)=\xi$. Then we have
\begin{equation}
M(i_{S,\infty})(\eta) = M(\tau_S \of i_{S,\infty})(\eta) = M(i_{S,\infty}) \of M(\tau_S) \, (\eta) = M(i_{S,\infty})(\xi) = 0, 
\end{equation}
which shows that $\eta \in \ker M(i_{S,\infty})$. Moreover, we have
\begin{align}
M(i_{S,0}) \big( \delta_S (\eta) \big)
& = M(i_{S,0})( \eta - X \cdot M(\tau_S)(\eta) ) \\
& = M(i_{S,0})(\eta)- X \cdot \big( M(i_{S,0}) \of M(\tau_S)(\eta) \big) \\
& = X \cdot \gamma_S(f) - X \cdot M(i_{S,0})(\xi) \\
& = X \cdot \gamma_S(f) - X \cdot \gamma_S(f) = 0 \: ; \\
M(\tau_S) \big( \delta_S (\eta) \big)
& = M(\tau_S)( \eta - X \cdot M(\tau_S)(\eta) ) \\
& = M(\tau_S)(\eta) - X \cdot \big( M(\tau_S) \of M(\tau_S)(\eta) \big) \\
& = \xi - X \cdot M(\tau_S)(\xi) \\
& = \delta_S (\xi) \\
& = (\beta \of \iota)_{S \times \N_{\infty}} \of \Gamma_S \, (f) .
\end{align}
Thus we have
\begin{align}
M(i_{S,0})\big( (\beta \of \iota)_{S \times \N_{\infty}} \of \Gamma_S \, (X \cdot f) \big)
& = 0 \\
& = M(i_{S,0}) \big( \delta_S (\eta) \big) \: ; \\
M(\tau_S)\big( (\beta \of \iota)_{S \times \N_{\infty}} \of \Gamma_S \, (X \cdot f) \big)
& = (\beta \of \iota)_{S \times \N_{\infty}} \of \Gamma_S \, (f) \\
& = M(\tau_S) \big( \delta_S (\eta) \big) .
\end{align}
Then \cref{lem:shift map is bijective} shows that
\begin{equation}
(\beta \of \iota)_{S \times \N_{\infty}} \of \Gamma_S \, (X \cdot f)
= \delta_S (\eta) .
\end{equation}
Therefore
\begin{align}
\gamma_S(X \cdot f)
& = M(i_{S,0}) \of {\delta_S}^{-1} \of (\beta \of \iota)_{S \times \N_{\infty}} \of \Gamma_S \, (X \cdot f) \\
& = M(i_{S,0}) \of {\delta_S}^{-1} \of \delta_S \, (\eta) \\
& = M(i_{S,0})(\eta) \\
& = X \cdot \gamma_S(f) .
\end{align}
This completes the proof that $\gamma = (\gamma_S)_{S \in |\ub{Prof}|}$ is a homomorphism of condensed $\tu{R}[X]$-modules $\tu{R \langle X \rangle} \to M$.

Next we prove that $\beta = \gamma \of \phi$. Let $S$ be any light profinite set. We show that $\beta_S = \gamma_S \of \phi_S$. Since both $\beta_S$ and $\gamma_S \of \phi_S$ are homomorphisms of $\tu{R}[X](S)$-modules $\tu{R}[X](S) \to M(S)$, it suffices to show that $\beta_S (1) = \gamma_S \of \phi_S \, (1)$. Since $\phi : \tu{R}[X] \to \tu{R \langle X \rangle}$ is a homomorphism of condensed rings, we have $\phi_S (1) = 1$. Therefore it suffices to show that $\beta_S (1) = \gamma_S (1)$. The following diagrams are commutative.
\[\begin{tikzcd}
	{\tu{R}(S \times \N_{\infty})} & {M(S \times \N_{\infty})} & {\tu{R}(S \times \N_{\infty})} & {M(S \times \N_{\infty})} \\
	{\tu{R}(S)} & {M(S)} & {\tu{R}(S \times \N_{\infty})} & {M(S \times \N_{\infty})}
	\arrow["{(\beta \of \iota)_{S \times \N_{\infty}}}", from=1-1, to=1-2]
	\arrow["{\tu{R}(i_{S,0})}"', from=1-1, to=2-1]
	\arrow["{M(i_{S,0})}", from=1-2, to=2-2]
	\arrow["{(\beta \of \iota)_{S \times \N_{\infty}}}", from=1-3, to=1-4]
	\arrow["{\tu{R}(\tau_S)}"', from=1-3, to=2-3]
	\arrow["{M(\tau_S)}", from=1-4, to=2-4]
	\arrow["{(\beta \of \iota)_{S}}"', from=2-1, to=2-2]
	\arrow["{(\beta \of \iota)_{S \times \N_{\infty}}}"', from=2-3, to=2-4]
\end{tikzcd}\]
Therefore
\begin{align}
M(i_{S,0}) \big( (\beta \of \iota)_{S \times \N_{\infty}} \of \Gamma_S \, (1) \big)
& = (\beta \of \iota)_S \of \tu{R}(i_{S,0}) \, \big( \Gamma_S \, (1) \big) \: ; \\
M(\tau_S) \big( (\beta \of \iota)_{S \times \N_{\infty}} \of \Gamma_S \, (1) \big)
& = (\beta \of \iota)_{S \times \N_{\infty}} \of \tu{R}(\tau_S) \, \big( \Gamma_S \, (1) \big) .
\end{align}
For every $(x,n) \in S \times \N$, we have
\begin{equation}
\Gamma_S (1) (x,n) = (1)_n = \left\{
\begin{aligned}
& 1 & (n = 0) \\
& 0 & (n \geq 1)
\end{aligned}
\right. .
\end{equation}
For every $x \in S$, we have
\begin{equation}
\Gamma_S (1) (x,\infty) = 0 .
\end{equation}
Therefore we have
\begin{align}
\tu{R}(i_{S,0}) \big( \Gamma_S (1) \big)
& = \Gamma_S (1) \of i_{S,0} = 1 \: ; \\
\tu{R}(\tau_S) \big( \Gamma_S (1) \big)
& = \Gamma_S (1) \of \tau_S = 0 .
\end{align}
Hence
\begin{align}
M(i_{S,0})\big( (\beta \of \iota)_{S \times \N_{\infty}} \of \Gamma_S \, (1) \big)
& = (\beta \of \iota)_S \of \tu{R}(i_{S,0}) \, \big( \Gamma_S (1) \big) \\
& = (\beta \of \iota)_S \, (1) \\
& = \beta_S \of \iota_S \, (1) \\ 
& = \beta_S (1) \: ; \\
M(\tau_S) \big( (\beta \of \iota)_{S \times \N_{\infty}} \of \Gamma_S \, (1) \big)
& = (\beta \of \iota)_{S \times \N_{\infty}} \of \tu{R}(\tau_S) \, \big( \Gamma_S (1) \big) \\
& = (\beta \of \iota)_{S \times \N_{\infty}} \, (0) = 0 .
\end{align}
Then we have
\begin{align}
\delta_S \big( (\beta \of \iota)_{S \times \N_{\infty}} \of \Gamma_S \, (1) \big)
& = \big( (\beta \of \iota)_{S \times \N_{\infty}} \of \Gamma_S \, (1) \big) - X \cdot M(\tau_S) \big( (\beta \of \iota)_{S \times \N_{\infty}} \of \Gamma_S \, (1) \big) \\
& = \big( (\beta \of \iota)_{S \times \N_{\infty}} \of \Gamma_S \, (1) \big) - X \cdot 0 \\
& = (\beta \of \iota)_{S \times \N_{\infty}} \of \Gamma_S \, (1) .
\end{align}
Therefore
\begin{align}
\gamma_S (1)
& = M(i_{S,0}) \of {\delta_S}^{-1} \of (\beta \of \iota)_{S \times \N_{\infty}} \of \Gamma_S \, (1) \\
& = M(i_{S,0}) \of {\delta_S}^{-1} \of \delta_S \, \big( (\beta \of \iota)_{S \times \N_{\infty}} \of \Gamma_S \, (1) \big) \\
& = M(i_{S,0}) \, \big( (\beta \of \iota)_{S \times \N_{\infty}} \of \Gamma_S \, (1) \big) \\
& = \beta_S (1).
\end{align}
This completes the proof that $\beta = \gamma \of \phi$.

Next suppose that $\epsilon : \tu{R \langle X \rangle} \to M$ is a homomorphism of condensed $\tu{R}[X]$-modules such that $\beta = \epsilon \of \phi$. We prove that $\epsilon = \gamma$. Let $S$ be any light profinite set. We show that $\epsilon_S = \gamma_S$. \cref{lem:definition of Theta S,lem:commutativity of Gamma S and Theta S} show that the following diagram is commutative.
\[\begin{tikzcd}
	& {\tu{R \langle X \rangle}(S)} && {M(S)} \\
	{\tu{R \langle X \rangle}(S)} & {\ker \tu{R \langle X \rangle}(i_{S,\infty})} && {\ker M(i_{S,\infty})} \\
	{\ker \tu{R}(i_{S,\infty})} & {\ker \tu{R \langle X \rangle}(i_{S,\infty})} && {\ker M(i_{S,\infty})}
	\arrow["{\epsilon_S}", from=1-2, to=1-4]
	\arrow["\id", curve={height=-12pt}, from=2-1, to=1-2]
	\arrow["{\Theta_S}", from=2-1, to=2-2]
	\arrow["{\Gamma_S}"', from=2-1, to=3-1]
	\arrow["{\tu{R \langle X \rangle}(i_{S,0})}"', from=2-2, to=1-2]
	\arrow["{\epsilon_{S \times \N_{\infty}}}", from=2-2, to=2-4]
	\arrow["{\id - X \cdot (\tu{R \langle X \rangle}(\tau_S))}", from=2-2, to=3-2]
	\arrow["{M(i_{S,0})}"', from=2-4, to=1-4]
	\arrow["{\id - X \cdot M(\tau_S) \, = \, \delta_S}", from=2-4, to=3-4]
	\arrow["{\tu{j}_{S \times \N_{\infty}}}"', from=3-1, to=3-2]
	\arrow["{\epsilon_{S \times \N_{\infty}}}"', from=3-2, to=3-4]
\end{tikzcd}\]
Therefore
\begin{align}
\epsilon_S
& = M(i_{S,0}) \of {\delta_S}^{-1} \of \epsilon_{S \times \N_{\infty}} \of \tu{j}_{S \times \N_{\infty}} \of \Gamma_S \\
& = M(i_{S,0}) \of {\delta_S}^{-1} \of \epsilon_{S \times \N_{\infty}} \of (\phi \of \iota)_{S \times \N_{\infty}} \of \Gamma_S \\
& = M(i_{S,0}) \of {\delta_S}^{-1} \of (\epsilon \of \phi \of \iota)_{S \times \N_{\infty}} \of \Gamma_S \\
& = M(i_{S,0}) \of {\delta_S}^{-1} \of (\beta \of \iota)_{S \times \N_{\infty}} \of \Gamma_S \\
& = \gamma_S .
\end{align}
This completes the proof that $\epsilon = \gamma$.
\end{proof}

\subsubsection{The general case}

\begin{lem} \label{lem:another description of coalescence of polynomial ring}
Let $n \in \N$. Let $S$ be a condensed ring. Let $\psi : \tu{R}[X_1, \ldots, X_n] \to S$ be a homomorphism of condensed rings. We consider $S$ as a condensed $\tu{R}[X_1, \ldots, X_n]$-algebra via $\psi$.

Write $g_i := \psi_{*}(X_i) \in S(*)$ for $1 \leq i \leq n$. Then the following conditions are equivalent.
\begin{enumerate}
\item
$\psi : \tu{R}[X_1, \ldots, X_n] \to S$ is the $\{ (X_1,1) , \ldots , (X_n,1) \}$-coalescence of $\tu{R}[X_1, \ldots, X_n]$ as a condensed $\tu{R}[X_1, \ldots, X_n]$-algebra.

\item
\begin{enumerate}
\item $S$ is $\{ (g_1,1), \ldots , (g_n,1) \}$-coalescent when considered as a condensed $S$-algebra via the identity $S \to S$. 

\item
Let $U$ be any condensed ring and let $h_1, \ldots , h_n \in U(*)$. Suppose that $U$ is $\{ (h_1,1), \ldots , (h_n,1) \}$-coalescent when considered as a condensed $U$-algebra via the identity $U \to U$. Then for every homomorphism $\sigma : \tu{R} \to U$ of condensed rings, there exists a unique homomorphism $\tau : S \to U$ of condensed rings such that $\tau_{*}(g_i) = h_i$ for $1 \leq i \leq n$ and the following diagram is commutative.
\[\begin{tikzcd}
	{\tu{R}} & U \\
	{\tu{R}[X_1 , \ldots , X_n]} & S
	\arrow["\sigma", from=1-1, to=1-2]
	\arrow["{\iota_n}"', from=1-1, to=2-1]
	\arrow["\psi"', from=2-1, to=2-2]
	\arrow["\tau"', from=2-2, to=1-2]
\end{tikzcd}\]
\end{enumerate}

\end{enumerate}
\end{lem}

\begin{proof}
First suppose that the condition (1) holds. We show that (2) holds. First of all, $S$ is an $\{ (X_1,1) , \ldots , (X_n,1) \}$-coalescent condensed $\tu{R}[X_1, \ldots, X_n]$-algebra via $\psi : \tu{R}[X_1, \ldots, X_n] \to S$. Since $g_i = \psi_{*}(X_i)$ for $1 \leq i \leq n$, we conclude that $S$ is a $\{ (g_1,1), \ldots , (g_n,1) \}$-coalescent condensed $S$-algebra via the identity $S \to S$.

Next let $U$ be any condensed ring. Let $\sigma : \tu{R} \to U$ be any homomorphism of condensed rings. Let $h_1, \ldots , h_n \in U(*)$ and suppose that $U$ is $\{ (h_1,1), \ldots , (h_n,1) \}$-coalescent when considered as a condensed $U$-algebra via the identity $U \to U$. By \cref{prop:universality of monoid algebras over condensed rings}, there exists a unique homomorphism $\tilde{\sigma} : \tu{R}[X_1 , \ldots , X_n] \to U$ of condensed rings such that $\tilde{\sigma}_*(X_i) = h_i$ for each $1 \leq i \leq n$ and the diagram
\[\begin{tikzcd}
	{\tu{R}} & U \\
	{\tu{R}[X_1 , \ldots , X_n]}
	\arrow["\sigma", from=1-1, to=1-2]
	\arrow["{\iota_n}"', from=1-1, to=2-1]
	\arrow["{\tilde{\sigma}}"', from=2-1, to=1-2]
\end{tikzcd}\]
is commutative. Let us consider $U$ as a condensed $\tu{R}[X_1 , \ldots , X_n]$-algebra via $\tilde{\sigma} : \tu{R}[X_1 , \ldots , X_n] \to U$. The condensed $U$-algebra $U$ is $\{ (h_1,1), \ldots , (h_n,1) \}$-coalescent by assumption. Since $\tilde{\sigma}_*(X_i) = h_i$ for each $1 \leq i \leq n$, we conclude that $U$ is a $\{ (X_1,1) , \ldots , (X_n,1) \}$-coalescent condensed $\tu{R}[X_1 , \ldots , X_n]$-algebra via $\tilde{\sigma} : \tu{R}[X_1 , \ldots , X_n] \to U$. Then the condition (1) implies that there exists a unique homomorphism $\tau : S \to U$ of condensed $\tu{R}[X_1 , \ldots , X_n]$-algebras such that the following diagram is commutative.
\[\begin{tikzcd}
	{\tu{R}[X_1 , \ldots , X_n]} & U \\
	S
	\arrow["{\tilde{\sigma}}", from=1-1, to=1-2]
	\arrow["\psi"', from=1-1, to=2-1]
	\arrow["\tau"', from=2-1, to=1-2]
\end{tikzcd}\]
Then the diagram
\[\begin{tikzcd}
	{\tu{R}} & U \\
	{\tu{R}[X_1 , \ldots , X_n]} & S
	\arrow["\sigma", from=1-1, to=1-2]
	\arrow["{\iota_n}"', from=1-1, to=2-1]
	\arrow["{\tilde{\sigma}}", from=2-1, to=1-2]
	\arrow["\psi"', from=2-1, to=2-2]
	\arrow["\tau"', from=2-2, to=1-2]
\end{tikzcd}\]
is commutative, and we have $\tau_*(g_i) = \tau_* \psi_*(X_i) = \tilde{\sigma}_*(X_i) = h_i$ for each $1 \leq i \leq n$.

On the other hand, suppose that $\tau' : S \to U$ is another homomorphism of condensed rings such that $\tau'_*(g_i) = h_i$ for each $1 \leq i \leq n$ and the diagram
\[\begin{tikzcd}
	{\tu{R}} & U \\
	{\tu{R}[X_1 , \ldots , X_n]} & S
	\arrow["\sigma", from=1-1, to=1-2]
	\arrow["{\iota_n}"', from=1-1, to=2-1]
	\arrow["\psi"', from=2-1, to=2-2]
	\arrow["{\tau'}"', from=2-2, to=1-2]
\end{tikzcd}\]
is commutative. Then we have $h_i = \tau'_*(g_i) = \tau'_* \psi_*(X_i)$ for each $1 \leq i \leq n$, and the diagram
\[\begin{tikzcd}
	{\tu{R}} & U \\
	{\tu{R}[X_1 , \ldots , X_n]}
	\arrow["\sigma", from=1-1, to=1-2]
	\arrow["{\iota_n}"', from=1-1, to=2-1]
	\arrow["{\tau' \of \psi}"', from=2-1, to=1-2]
\end{tikzcd}\]
is commutative. By the definition of $\tilde{\sigma}$, we conclude that $\tau' \of \psi = \tilde{\sigma}$. In other words, the diagram
\[\begin{tikzcd}
	{\tu{R}[X_1 , \ldots , X_n]} & U \\
	S
	\arrow["{\tilde{\sigma}}", from=1-1, to=1-2]
	\arrow["\psi"', from=1-1, to=2-1]
	\arrow["{\tau'}"', from=2-1, to=1-2]
\end{tikzcd}\]
is commutative. Then $\tau'$ is a homomorphism of condensed $\tu{R}[X_1 , \ldots , X_n]$-algebras. By the definition of $\tau$, we have $\tau' = \tau$. This completes the proof that the condition (1) implies (2).

Conversely, suppose that the condition (2) holds. We show that (1) holds. First of all, $S$ is a $\{ (g_1,1), \ldots , (g_n,1) \}$-coalescent condensed $S$-algebra via the identity $S \to S$. Since $g_i = \psi_{*}(X_i)$ for $1 \leq i \leq n$, we conclude that $S$ is an $\{ (X_1,1) , \ldots , (X_n,1) \}$-coalescent condensed $\tu{R}[X_1, \ldots, X_n]$-algebra via $\psi : \tu{R}[X_1, \ldots, X_n] \to S$. 

Next let $U$ be any $\{ (X_1,1) , \ldots , (X_n,1) \}$-coalescent condensed $\tu{R}[X_1, \ldots, X_n]$-algebra. Let $\sigma : \tu{R}[X_1, \ldots, X_n] \to U$ be any homomorphism of condensed $\tu{R}[X_1, \ldots, X_n]$-algebras. We prove that there exists a unique homomorphism $\tau : S \to U$ of condensed $\tu{R}[X_1, \ldots, X_n]$-algebras such that the diagram
\[\begin{tikzcd}
	{\tu{R}[X_1 , \ldots , X_n]} & U \\
	S
	\arrow["\sigma", from=1-1, to=1-2]
	\arrow["\psi"', from=1-1, to=2-1]
	\arrow["\tau"', from=2-1, to=1-2]
\end{tikzcd}\]
is commutative.

Write $h_i := \sigma_*(X_i) \in U(*)$ for each $1 \leq i \leq n$. Since $U$ is an $\{ (X_1,1) , \ldots , (X_n,1) \}$-coalescent condensed $\tu{R}[X_1, \ldots, X_n]$-algebra via $\sigma : \tu{R}[X_1, \ldots, X_n] \to U$, it follows that $U$ is an $\{ (h_1,1), \ldots , (h_n,1) \}$-coalescent condensed $U$-algebra via the identity $U \to U$. Then the condition (2) implies that there exists a unique homomorphism $\tau : S \to U$ of condensed rings such that $\tau_*(g_i) = h_i$ for each $1 \leq i \leq n$ and the diagram
\[\begin{tikzcd}
	{\tu{R}} & {\tu{R}[X_1 , \ldots , X_n]} & U \\
	{\tu{R}[X_1 , \ldots , X_n]} && S
	\arrow["{\iota_n}", from=1-1, to=1-2]
	\arrow["{\iota_n}"', from=1-1, to=2-1]
	\arrow["\sigma", from=1-2, to=1-3]
	\arrow["\psi"', from=2-1, to=2-3]
	\arrow["\tau"', from=2-3, to=1-3]
\end{tikzcd}\]
is commutative. Then we have $\tau_* \psi_*(X_i) = \tau_*(g_i) = h_i = \sigma_*(X_i)$ for each $1 \leq i \leq n$. By \cref{prop:universality of monoid algebras over condensed rings}, we conclude that the diagram
\[\begin{tikzcd}
	{\tu{R}[X_1 , \ldots , X_n]} & U \\
	S
	\arrow["\sigma", from=1-1, to=1-2]
	\arrow["\psi"', from=1-1, to=2-1]
	\arrow["\tau"', from=2-1, to=1-2]
\end{tikzcd}\]
is commutative. Then $\tau : S \to U$ is a homomoprhism of condensed $\tu{R}[X_1, \ldots, X_n]$-algebras.

On the other hand, suppose that $\tau' : S \to U$ is another homomorphism of condensed $\tu{R}[X_1, \ldots, X_n]$-algebras such that the diagram
\[\begin{tikzcd}
	{\tu{R}[X_1 , \ldots , X_n]} & U \\
	S
	\arrow["\sigma", from=1-1, to=1-2]
	\arrow["\psi"', from=1-1, to=2-1]
	\arrow["{\tau'}"', from=2-1, to=1-2]
\end{tikzcd}\]
is commutative. Then $\tau'_*(g_i) = \tau'_* \psi_*(X_i) = \sigma_*(X_i) = h_i$ for each $1 \leq i \leq n$, and the diagram
\[\begin{tikzcd}
	{\tu{R}} & {\tu{R}[X_1 , \ldots , X_n]} & U \\
	{\tu{R}[X_1 , \ldots , X_n]} && S
	\arrow["{\iota_n}", from=1-1, to=1-2]
	\arrow["{\iota_n}"', from=1-1, to=2-1]
	\arrow["\sigma", from=1-2, to=1-3]
	\arrow["\psi"', from=2-1, to=2-3]
	\arrow["{\tau'}"', from=2-3, to=1-3]
\end{tikzcd}\]
is commutative. By the definition of $\tau$, we conclude that $\tau' = \tau$. This completes the proof.
\end{proof}

\begin{proof}[Proof of \cref{prop:coalescence of polynomial rings}]
By \cref{lem:another description of coalescence of polynomial ring}, the assertion of \cref{prop:coalescence of polynomial rings} is equivalent to the following assertions.
\begin{enumerate}
\item The condensed $\tu{R \langle X_1 , \ldots , X_n \rangle}$-algebra $\tu{R \langle X_1 , \ldots , X_n \rangle}$ is $\{ (X_1,1) , \ldots , (X_n,1) \}$-coalescent.

\item
Let $U$ be any condensed ring and let $h_1, \ldots , h_n \in U(*)$. Suppose that $U$ is $\{ (h_1,1), \ldots , (h_n,1) \}$-coalescent when considered as a condensed $U$-algebra via the identity $U \to U$. Then for every homomorphism $\sigma : \tu{R} \to U$ of condensed rings, there exists a unique homomorphism $\tau : \tu{R \langle X_1 , \ldots , X_n \rangle} \to U$ of condensed rings such that $\tau_{*}(X_i) = h_i$ for $1 \leq i \leq n$ and the following diagram is commutative.
\[\begin{tikzcd}
	{\tu{R}} & U \\
	{\tu{R \langle X_1 , \ldots , X_n \rangle}}
	\arrow["\sigma", from=1-1, to=1-2]
	\arrow["{\tu{j_n}}"', from=1-1, to=2-1]
	\arrow["\tau"', from=2-1, to=1-2]
\end{tikzcd}\]
\end{enumerate}
We use the induction on $n$ to prove these assertions. The case $n=0$ is trivial. The case $n=1$ follows from the fact that \cref{prop:coalescence of polynomial rings} holds in the case $n=1$ by \cref{prop:coalescence of polynomial rings in one variable}. 

Suppose $n \geq 2$. First of all, since the elements $X_1, \ldots, X_n \in R \langle X_1 , \ldots , X_n \rangle$ are power-bounded in $R \langle X_1 , \ldots , X_n \rangle$, \cref{prop:power-boundedness and coalescence} shows that the condensed $\tu{R \langle X_1 , \ldots , X_n \rangle}$-algebra $\tu{R \langle X_1 , \ldots , X_n \rangle}$ is $\{ (X_1,1) , \ldots , (X_n,1) \}$-coalescent.

Next let $U$ be any condensed ring. Let $\sigma : \tu{R} \to U$ be any homomorphism of condensed rings. Let $h_1, \ldots, h_n \in U(*)$ and suppose that the condensed $U$-algebra $U$ is $\{ (h_1,1), \ldots, (h_n,1) \}$-coalescent. By the induction hypothesis, there exists a unique homomorphism $\tilde{\sigma} : \tu{R \langle X_1 , \ldots , X_{n-1} \rangle} \to U$ of condensed rings such that $\tilde{\sigma}_*(X_i) = h_i$ for each $1 \leq i \leq n-1$ and the diagram
\[\begin{tikzcd}
	{\tu{R}} & U \\
	{\tu{R \langle X_1 , \ldots , X_{n-1} \rangle}}
	\arrow["\sigma", from=1-1, to=1-2]
	\arrow["{\tu{j_{n-1}}}"', from=1-1, to=2-1]
	\arrow["{\tilde{\sigma}}"', from=2-1, to=1-2]
\end{tikzcd}\]
is commutative. On the other hand, let $j_{n-1,n} : R \langle X_1 , \ldots , X_{n-1} \rangle \to R \langle X_1 , \ldots , X_{n-1} \rangle \langle X_n \rangle = R \langle X_1 , \ldots , X_n \rangle$ be the canonical inclusion. We apply the case $n=1$ to the complete Hausdorff non-archimedean ring $R \langle X_1 , \ldots , X_{n-1} \rangle$. We conclude that there exists a unique homomorphism $\tau : R \langle X_1 , \ldots , X_n \rangle \to U$ of condensed rings such that $\tau_*(X_n) = h_n$ and the diagram
\[\begin{tikzcd}
	{\tu{R \langle X_1 , \ldots , X_{n-1} \rangle}} & U \\
	{\tu{R \langle X_1 , \ldots , X_n \rangle}}
	\arrow["{\tilde{\sigma}}", from=1-1, to=1-2]
	\arrow["{\tu{j_{n-1,n}}}"', from=1-1, to=2-1]
	\arrow["\tau"', from=2-1, to=1-2]
\end{tikzcd}\]
is commutative. Then the diagram
\[\begin{tikzcd}
	{\tu{R}} & U \\
	{\tu{R \langle X_1 , \ldots , X_{n-1} \rangle}} \\
	{\tu{R \langle X_1 , \ldots , X_n \rangle}}
	\arrow["\sigma", from=1-1, to=1-2]
	\arrow["{\tu{j_{n-1}}}", from=1-1, to=2-1]
	\arrow["{\tu{j_n}}"', shift right=5, curve={height=30pt}, from=1-1, to=3-1]
	\arrow["{\tilde{\sigma}}"', from=2-1, to=1-2]
	\arrow["{\tu{j_{n-1,n}}}", from=2-1, to=3-1]
	\arrow["\tau"', curve={height=24pt}, from=3-1, to=1-2]
\end{tikzcd}\]
is commutative, and we have $\tau_*(X_i) = h_i$ for each $1 \leq i \leq n$.

On the other hand, suppose that $\tau' : \tu{R \langle X_1 , \ldots , X_n \rangle} \to U$ is another homomorphism of condensed rings such that $\tau'_*(X_i) = h_i$ for each $1 \leq i \leq n$ and the diagram
\[\begin{tikzcd}
	{\tu{R}} & U \\
	{\tu{R \langle X_1 , \ldots , X_n \rangle}}
	\arrow["\sigma", from=1-1, to=1-2]
	\arrow["{\tu{j_n}}"', from=1-1, to=2-1]
	\arrow["{\tau'}"', from=2-1, to=1-2]
\end{tikzcd}\]
is commutative. The $\tau'_*(j_{n-1,n}(X_i)) = \tau'_*(X_i)= h_i$ for each $1 \leq i \leq n-1$, and the diagram
\[\begin{tikzcd}
	{\tu{R}} & U \\
	{\tu{R \langle X_1 , \ldots , X_{n-1} \rangle}} & {\tu{R \langle X_1 , \ldots , X_n \rangle}}
	\arrow["\sigma", from=1-1, to=1-2]
	\arrow["{\tu{j_{n-1}}}"', from=1-1, to=2-1]
	\arrow["{\tu{j_n}}", from=1-1, to=2-2]
	\arrow["{\tu{j_{n-1,n}}}"', from=2-1, to=2-2]
	\arrow["{\tau'}"', from=2-2, to=1-2]
\end{tikzcd}\]
is commutative. By the definition of $\tilde{\sigma}$, we conclude that $\tau' \of \tu{j_{n-1,n}} = \tilde{\sigma}$. In other words, the diagram
\[\begin{tikzcd}
	{\tu{R \langle X_1 , \ldots , X_{n-1} \rangle}} & U \\
	{\tu{R \langle X_1 , \ldots , X_n \rangle}}
	\arrow["{\tilde{\sigma}}", from=1-1, to=1-2]
	\arrow["{\tu{j_{n-1,n}}}"', from=1-1, to=2-1]
	\arrow["{\tau'}"', from=2-1, to=1-2]
\end{tikzcd}\]
is commutative. Moreover, we have $\tau'(X_n) = h_n$. By the definition of $\tau$, we conclude that $\tau' = \tau$. This completes the proof.
\end{proof}

\subsection{Relation to localization}

\subsubsection{Settings}

$\\[2mm]$ Throughout this subsection, we fix the following notation.

\begin{nt} \;
\begin{enumerate}
\item $A$ denotes a complete Hausdorff Huber ring.

\item $g \in A$ is an element of $A$.

\item $T$ denotes a finite subset of $A$ which generates an open ideal of $A$.

\item $l : A \to A \left \langle \frac{T}{g} \right \rangle$ denotes the canonical homomorphism.

\item We write $T' := \set{(f,g)}{f \in T}$. This is a subset of $A \times A$.
\end{enumerate}
\end{nt}

\subsubsection{Statement}

\begin{prop} \label{prop:comparison of localizations of Huber rings} \;
\begin{enumerate}
\item
There exists a unique homomorphism of condensed rings $\rho : (\tu{A}_g)_{\approx T'} \to \tu{A \left \langle \frac{T}{g} \right \rangle}$ such that the diagram
\[\begin{tikzcd}
	{\tu{A}} & {(\tu{A}_g)_{\approx T'}} \\
	& {\tu{A \left \langle \frac{T}{g} \right \rangle}}
	\arrow["\can", from=1-1, to=1-2]
	\arrow["{\tu{l}}"', from=1-1, to=2-2]
	\arrow["\rho", from=1-2, to=2-2]
\end{tikzcd}\]
is commutative. Here the homomorphism $\tu{A} \xto{\can} (\tu{A}_g)_{\approx T'}$ of condensed rings is the composition of the canonical homomorphism $\tu{A} \to \tu{A}_g$ into the localization $\tu{A}_g$ of the condensed $\tu{A}$-algebra $\tu{A}$ by $g$, and the $T'$-coalescence $\tu{A}_g \to (\tu{A}_g)_{\approx T'}$ of the condensed $\tu{A}$-algebra $\tu{A}_g$.

\item $\rho : (\tu{A}_g)_{\approx T'} \to \tu{A \left \langle \frac{T}{g} \right \rangle}$ is an epimorphism when considered as a morphism in $\ub{CSet}$.

\item
Suppose that one of the following holds.
\begin{enumerate}
\item $A$ is a strongly Noetherian Tate ring.
\item $A$ has a Noetherian ring of definition.
\end{enumerate}
Then $\rho : (\tu{A}_g)_{\approx T'} \to \tu{A \left \langle \frac{T}{g} \right \rangle}$ is an isomorphism of condensed rings.
\end{enumerate}
\end{prop}

The proof of this proposition is the content of this subsection.

\subsubsection{Notations}

\begin{nt} \;
\begin{enumerate}
\item We enumerate all the elements of $T$ and write $T = \{ f_1 , \ldots , f_n \}$, where $n \in \N$.

\item $A \langle X_1 , \ldots , X_n \rangle$ denotes the ring of restricted formal power series in $n$ variables over $A$. We write $j : A \mon A \langle X_1 , \ldots , X_n \rangle$ for the canonical inclusion.

\item
${A \langle X_1 , \ldots , X_n \rangle}^{n}$ denotes the direct sum of $n$ copies of the topological $A \langle X_1 , \ldots , X_n \rangle$-module $A \langle X_1 , \ldots , X_n \rangle$ in $\ub{TopMod}_{A \langle X_1 , \ldots , X_n \rangle}$.

\item $m : {A \langle X_1 , \ldots , X_n \rangle}^{n} \to A \langle X_1 , \ldots , X_n \rangle$ denotes the continuous $A \langle X_1 , \ldots , X_n \rangle$-linear map induced by the continuous $A \langle X_1 , \ldots , X_n \rangle$-linear maps
\[\begin{tikzcd}
	{A \langle X_1 , \ldots , X_n \rangle} && {A \langle X_1 , \ldots , X_n \rangle} & {(1 \leq i \leq n).}
	\arrow["{(f_i-gX_i) \cdot \id}", from=1-1, to=1-3]
\end{tikzcd}\]

\item $J$ denotes the image of $m : {A \langle X_1 , \ldots , X_n \rangle}^{n} \to A \langle X_1 , \ldots , X_n \rangle$, i.e., the ideal of $A \langle X_1 , \ldots , X_n \rangle$ generated by $f_1 - g X_1 , \ldots , f_n - g X_n$. Let $\overline{J \,}$ be the closure of $J$ in $A \langle X_1 , \ldots , X_n \rangle$. We define
\begin{equation}
B := A \langle X_1 , \ldots , X_n \rangle / \overline{J \,} .
\end{equation}
We write $q : \overline{J \,} \mon A \langle X_1 ,\dots, X_n \rangle$ for the inclusion and $p : A \langle X_1 , \ldots , X_n \rangle \epi B$ for the canonical surjection.

\item We write $r : A \to B$ for the composition $p \of j$. We consider $B$ as a commutative unital topological $A$-algebra via $r$.
\end{enumerate}
\end{nt}

\begin{nt} \;
\begin{enumerate}
\item
$\tu{A}[X_1 , \ldots , X_n]$ denotes the polynomial ring in $n$ variables over $\tu{A}$. $\iota : \tu{A} \to \tu{A}[X_1 , \ldots , X_n]$ denotes the canonical inclusion.

\item
${\tu{A}[X_1 , \ldots , X_n]}^{n}$ denotes the direct sum of $n$ copies of the condensed $\tu{A}[X_1 , \ldots , X_n]$-module $\tu{A}[X_1 , \ldots , X_n]$ in $\ub{CMod}_{\tu{A}[X_1 , \ldots , X_n]}$.

\item $\mu : {\tu{A}[X_1 , \ldots , X_n]}^{n} \to \tu{A}[X_1 , \ldots , X_n]$ denotes the homomorphism of condensed $\tu{A}[X_1 , \ldots , X_n]$-modules induced by the homomorphisms of condensed $\tu{A}[X_1 , \ldots , X_n]$-modules
\[\begin{tikzcd}
	{\tu{A}[X_1 , \ldots , X_n]} && {\tu{A}[X_1 , \ldots , X_n]} & {(1 \leq i \leq n).}
	\arrow["{(f_i-gX_i) \cdot \id}", from=1-1, to=1-3]
\end{tikzcd}\]

\item $R$ denotes the cokernel of $\mu : {\tu{A}[X_1 , \ldots , X_n]}^{n} \to \tu{A}[X_1 , \ldots , X_n]$ in $\ub{CMod}_{\tu{A}[X_1 , \ldots , X_n]}$. $\pi : \tu{A}[X_1 , \ldots , X_n] \epi R$ denotes the canonical epimorphism.

\item $R$ has a unique structure of a condensed $\tu{A}[X_1 , \ldots , X_n]$-algebra such that $\pi : \tu{A}[X_1 , \ldots , X_n] \epi R$ is a homomorphism of condensed $\tu{A}[X_1 , \ldots , X_n]$-algebras. Via this structure we consider $R$ as a condensed $\tu{A}[X_1 , \ldots , X_n]$-algebra.
\end{enumerate}
\end{nt}

\begin{nt} \;
\begin{enumerate}
\item We write $T'' := \{ (X_1,1) ,\ldots, (X_n,1) \}$. This is a subset of $\tu{A}[X_1 , \ldots , X_n](*) \times \tu{A}[X_1 , \ldots , X_n](*)$.

\item $\phi : \tu{A}[X_1 , \ldots , X_n] \to \tu{A \langle X_1 , \dots , X_n \rangle}$ denotes the unique homomorphism of condensed rings such that $\phi_{*}(X_i)=X_i$ for $1 \leq i \leq n$ and the following diagram is commutative.
\[\begin{tikzcd}
	{\tu{A}} & {\tu{A}[X_1 ,\dots, X_n]} \\
	{\tu{A \langle X_1 ,\ldots, X_n \rangle}}
	\arrow["\iota", from=1-1, to=1-2]
	\arrow["{\tu{j}}"', from=1-1, to=2-1]
	\arrow["\phi", from=1-2, to=2-1]
\end{tikzcd}\]
We consider $\tu{A \langle X_1 , \dots , X_n \rangle}$ as a condensed $\tu{A}[X_1 , \ldots , X_n]$-algebra via $\phi$.

\item $\psi : R \to R_{\approx T''}$ denotes the $T''$-coalescence of the condensed $\tu{A}[X_1 , \ldots , X_n]$-algebra $R$. 
\end{enumerate}
\end{nt}

\subsubsection{The homomorphism $\rho$}

\begin{proof}[Proof of (1) of \cref{prop:comparison of localizations of Huber rings}]
The element $l(g) \in A \left \langle \frac{T}{g} \right \rangle$ is invertible in $A \left \langle \frac{T}{g} \right \rangle$. Moreover, for every $f \in T$, the element $l(f) \cdot l(g)^{-1} \in A \left \langle \frac{T}{g} \right \rangle$ is power-bouded in $A \left \langle \frac{T}{g} \right \rangle$. Then \cref{prop:power-boundedness and coalescence} shows that the $\tu{A \left \langle \frac{T}{g} \right \rangle}$-algebra $\tu{A \left \langle \frac{T}{g} \right \rangle}$ is $\big( l(f) \cdot l(g)^{-1} \big) / 1$-coalescent. Since the element $l(g) \in A \left \langle \frac{T}{g} \right \rangle$ is invertible in $A \left \langle \frac{T}{g} \right \rangle$, we conclude that the $\tu{A \left \langle \frac{T}{g} \right \rangle}$-algebra $\tu{A \left \langle \frac{T}{g} \right \rangle}$ is $l(f) / l(g)$-coalescent. Therefore the $\tu{A}$-algebra $\tu{A \left \langle \frac{T}{g} \right \rangle}$ is $f/g$-coalescent. Since this holds for every $f \in T$, the $\tu{A}$-algebra $\tu{A \left \langle \frac{T}{g} \right \rangle}$ is $T'$-coalescent. Then \cref{prop:localization and coalescence commute} shows that there exists a unique homomorphism $\rho : (\tu{A}_g)_{\approx T'} \to \tu{A \left \langle \frac{T}{g} \right \rangle}$ of condensed $\tu{A}$-algebras such that the diagram
\[\begin{tikzcd}
	{\tu{A}} & {(\tu{A}_g)_{\approx T'}} \\
	& {\tu{A \left \langle \frac{T}{g} \right \rangle}}
	\arrow["\can", from=1-1, to=1-2]
	\arrow["{\tu{l}}"', from=1-1, to=2-2]
	\arrow["\rho", from=1-2, to=2-2]
\end{tikzcd}\]
is commutative.

On the other hand, let $\rho' : (\tu{A}_g)_{\approx T'} \to \tu{A \left \langle \frac{T}{g} \right \rangle}$ be another homomorphism of condensed rings such that the diagram
\[\begin{tikzcd}
	{\tu{A}} & {(\tu{A}_g)_{\approx T'}} \\
	& {\tu{A \left \langle \frac{T}{g} \right \rangle}}
	\arrow["\can", from=1-1, to=1-2]
	\arrow["{\tu{l}}"', from=1-1, to=2-2]
	\arrow["{\rho'}", from=1-2, to=2-2]
\end{tikzcd}\]
is commutative. Then $\rho'$ is necessarily a homomorphism of condensed $\tu{A}$-algebras $(\tu{A}_g)_{\approx T'} \to \tu{A \left \langle \frac{T}{g} \right \rangle}$. Therefore $\rho' = \rho$.
\end{proof}

\subsubsection{The homomorphisms $\sigma, \sigma'$}

\begin{lem}
We have a canonical isomorphism of condensed $\tu{A \langle X_1 ,\ldots, X_n \rangle}$-modules
\begin{equation}
\tu{{A \langle X_1 ,\ldots, X_n \rangle}^{n}} \simeq \Big( \tu{A \langle X_1 ,\ldots, X_n \rangle} \Big) ^{n} ,
\end{equation}
where $\Big( \tu{A \langle X_1 ,\ldots, X_n \rangle} \Big) ^{n}$ denotes the direct sum of $n$ copies of the condensed $\tu{A \langle X_1 ,\ldots, X_n \rangle}$-module $\tu{A \langle X_1 ,\ldots, X_n \rangle}$.

In the following, we identify $\tu{{A \langle X_1 ,\ldots, X_n \rangle}^{n}}$ with $\Big( \tu{A \langle X_1 ,\ldots, X_n \rangle} \Big) ^{n}$ via this isomorphism.
\end{lem}

\begin{proof}
It follows immediately from the definition that the functor
\[\begin{tikzcd}
	{\ub{TopMod}_{A \langle X_1 ,\ldots, X_n \rangle}} & {\ub{CMod}_{\tu{A \langle X_1 ,\ldots, X_n \rangle}} ,} & M & {\tu{M}}
	\arrow[from=1-1, to=1-2]
	\arrow[maps to, from=1-3, to=1-4]
\end{tikzcd}\]
is an additive functor between additive categories.
\end{proof}

\begin{prop} \label{prop:construction of the homomorphism chi} \;
\begin{enumerate}
\item There exists a unique homomorphism $\chi : \tu{A \langle X_1 ,\ldots, X_n \rangle} \to R_{\approx T''}$ of condensed $\tu{A}[X_1 , \ldots , X_n]$-algebras such that the following diagram is commutative.
\[\begin{tikzcd}
	{\tu{A}[X_1 , \dots , X_n]} & R \\
	{\tu{A \langle X_1 , \dots , X_n \rangle}} & {R_{\approx T''}}
	\arrow["\pi", from=1-1, to=1-2]
	\arrow["\phi"', from=1-1, to=2-1]
	\arrow["\psi", from=1-2, to=2-2]
	\arrow["\chi"', from=2-1, to=2-2]
\end{tikzcd}\]

\item We have the following commutative diagram in $\ub{CMod}_{\tu{A}[X_1 , \ldots , X_n]}$.
\[\begin{tikzcd}
	{{\tu{A}[X_1 , \dots , X_n]}^{n}} & {\tu{A}[X_1 , \dots , X_n]} & R & 0 \\
	{\tu{{A \langle X_1 , \dots , X_n \rangle}^{n}}} & {\tu{A \langle X_1 , \dots , X_n \rangle}} & {R_{\approx T''}} & 0
	\arrow["\mu", from=1-1, to=1-2]
	\arrow["{\phi^{n}}"', from=1-1, to=2-1]
	\arrow["\pi", from=1-2, to=1-3]
	\arrow["\phi"', from=1-2, to=2-2]
	\arrow[from=1-3, to=1-4]
	\arrow["\psi", from=1-3, to=2-3]
	\arrow["{\tu{m}}"', from=2-1, to=2-2]
	\arrow["\chi"', from=2-2, to=2-3]
	\arrow[from=2-3, to=2-4]
\end{tikzcd}\]
Moreover, the rows are exact in $\ub{CMod}_{\tu{A}[X_1 , \ldots , X_n]}$.
\end{enumerate}
\end{prop}

\begin{proof}~
\begin{enumerate}
\item
This follows from \cref{prop:coalescence of polynomial rings}.

\item
By definition, the diagram is commutative and the sequence
\[\begin{tikzcd}
	{{\tu{A}[X_1 , \ldots , X_n]}^{n}} & {\tu{A}[X_1 , \ldots , X_n]} & R & 0
	\arrow["\mu", from=1-1, to=1-2]
	\arrow["\pi", from=1-2, to=1-3]
	\arrow[from=1-3, to=1-4]
\end{tikzcd}\]
is exact in $\ub{CMod}_{\tu{A}[X_1 , \ldots , X_n]}$. We apply the $T''$-coalescence functor $\ub{CMod}_{\tu{A}[X_1 , \ldots , X_n]} \to \ub{CMod}_{\tu{A}[X_1 , \ldots , X_n] \approx T''}$ to this sequence. By \cref{prop:coalescence of polynomial rings}, we obtain a sequence
\[\begin{tikzcd}
	{\tu{{A \langle X_1 , \ldots , X_n \rangle}^{n}}} & {\tu{A \langle X_1 , \ldots , X_n \rangle}} & {R_{\approx T''}} & {0 .}
	\arrow["{\tu{m}}", from=1-1, to=1-2]
	\arrow["\chi", from=1-2, to=1-3]
	\arrow[from=1-3, to=1-4]
\end{tikzcd}\]
Since the $T''$-coalescence functor $\ub{CMod}_{\tu{A}[X_1 , \ldots , X_n]} \to \ub{CMod}_{\tu{A}[X_1 , \ldots , X_n] \approx T''}$ is a left adjoint, it follows that this sequence is exact in $\ub{CMod}_{\tu{A}[X_1 , \ldots , X_n] \approx T''}$. By \cref{cor:CModR approx T is Grothendieck Abelian}, this sequence is also exact in $\ub{CMod}_{\tu{A}[X_1 , \ldots , X_n]}$.
\end{enumerate}
\end{proof}

\begin{prop} \label{prop:construction of the homomorphism sigma} \;
\begin{enumerate}
\item
There exists a unique homomorphism $\sigma : R_{\approx T''} \to \tu{B}$ of condensed rings such that the following diagram is commutative.
\[\begin{tikzcd}
	{\tu{A \langle X_1 , \dots , X_n \rangle}} & {R_{\approx T''}} \\
	{\tu{B}}
	\arrow["\chi", from=1-1, to=1-2]
	\arrow["{\tu{p}}"', from=1-1, to=2-1]
	\arrow["\sigma", from=1-2, to=2-1]
\end{tikzcd}\]

\item
$\sigma : R_{\approx T''} \to \tu{B}$ is an epimorphism when considered as a morphism in $\ub{CSet}$.

\item
We have the following commutative diagram in $\ub{CMod}_{\tu{A}[X_1 , \ldots , X_n]}$.
\[\begin{tikzcd}
	& {\tu{{A \langle X_1 , \dots , X_n \rangle}^{n}}} & {\tu{A \langle X_1 , \dots , X_n \rangle}} & {R_{\approx T''}} & 0 \\
	0 & {\tu{\overline{J \,}}} & {\tu{A \langle X_1 , \dots , X_n \rangle}} & {\tu{B}} & 0
	\arrow["{\tu{m}}", from=1-2, to=1-3]
	\arrow["{\tu{m}}"', from=1-2, to=2-2]
	\arrow["\chi", from=1-3, to=1-4]
	\arrow[equals, from=1-3, to=2-3]
	\arrow[from=1-4, to=1-5]
	\arrow["\sigma", from=1-4, to=2-4]
	\arrow[from=2-1, to=2-2]
	\arrow["{\tu{q}}"', from=2-2, to=2-3]
	\arrow["{\tu{p}}"', from=2-3, to=2-4]
	\arrow[from=2-4, to=2-5]
\end{tikzcd}\]
Moreover, the rows are exact in $\ub{CMod}_{\tu{A}[X_1 , \ldots , X_n]}$.

\item
The following conditions are equivalent.
\begin{enumerate}
\item $m : {A \langle X_1 , \ldots , X_n \rangle}^{n} \to A \langle X_1 , \ldots , X_n \rangle$ is strict.

\item $\tu{m} : \tu{{A \langle X_1 , \dots , X_n \rangle}^{n}} \to \tu{\overline{J \,}}$ is an epimorphism in $\ub{CMod}_{\tu{A}[X_1 , \ldots , X_n]}$.

\item $\sigma : R_{\approx T''} \to \tu{B}$ is an isomorphism of condensed rings.
\end{enumerate}

\end{enumerate}
\end{prop}

\begin{proof}~
\begin{enumerate}
\item
We first claim that $\tu{B}$ is a $T''$-coalescent condensed $\tu{A}[X_1 , \ldots , X_n]$-algebra via $\tu{A}[X_1 , \ldots , X_n] \xto{\phi} \tu{A \langle X_1 , \ldots , X_n \rangle} \xto{\tu{p}} \tu{B}$. Indeed, for $1 \leq i \leq n$, the element $p(X_i) \in B$ is power-bounded in $B$ since the element $X_i \in A \langle X_1 , \ldots , X_n \rangle$ is power-bounded in $A \langle X_1 , \ldots , X_n \rangle$. $B$ is complete Hausdorff by \cref{prop:completeness of quotients of f c topological abelian groups}. Therefore \cref{prop:power-boundedness and coalescence} shows that $\tu{B}$ is a $p(X_i) / 1$-coalescent condensed $\tu{B}$-algebra. Therefore $\tu{B}$ is a $X_i/1$-coalescent condensed $\tu{A}[X_1 , \ldots , X_n]$-algebra via $\tu{A}[X_1 , \ldots , X_n] \xto{\phi} \tu{A \langle X_1 , \ldots , X_n \rangle} \xto{\tu{p}} \tu{B}$. Since this holds for every $1 \leq i \leq n$, we conclude that $\tu{B}$ is a $T''$-coalescent condensed $\tu{A}[X_1 , \ldots , X_n]$-algebra via $\tu{A}[X_1 , \ldots , X_n] \xto{\phi} \tu{A \langle X_1 , \ldots , X_n \rangle} \xto{\tu{p}} \tu{B}$.

From now on let us consider $\tu{B}$ as a condensed $\tu{A}[X_1 , \ldots , X_n]$-algebra via $\tu{A}[X_1 , \ldots , X_n] \xto{\phi} \tu{A \langle X_1 , \ldots , X_n \rangle} \xto{\tu{p}} \tu{B}$. Since the diagram
\[\begin{tikzcd}
	{{\tu{A}[X_1 , \ldots , X_n]}^{n}} & {\tu{A}[X_1 , \ldots , X_n]} \\
	{\tu{{A \langle X_1 , \ldots , X_n \rangle}^{n}}} & {\tu{A \langle X_1 , \dots , X_n \rangle}} & {\tu{B}}
	\arrow["\mu", from=1-1, to=1-2]
	\arrow["{\phi^n}"', from=1-1, to=2-1]
	\arrow["\phi"', from=1-2, to=2-2]
	\arrow["{\tu{m}}", from=2-1, to=2-2]
	\arrow["0"', curve={height=30pt}, from=2-1, to=2-3]
	\arrow["{\tu{p}}", from=2-2, to=2-3]
\end{tikzcd}\]
is commutative, the definition of $R$ shows that there exists a unique homomorphis $\sigma_0 : R \to \tu{B}$ of condensed $\tu{A}[X_1 , \ldots , X_n]$-algebras such that the diagram
\[\begin{tikzcd}
	{\tu{A}[X_1 , \ldots , X_n]} & R \\
	{\tu{A \langle X_1 , \dots , X_n \rangle}} & {\tu{B}}
	\arrow["\pi", from=1-1, to=1-2]
	\arrow["\phi"', from=1-1, to=2-1]
	\arrow["{\sigma_0}", from=1-2, to=2-2]
	\arrow["{\tu{p}}"', from=2-1, to=2-2]
\end{tikzcd}\]
is commutative. Since the condensed $\tu{A}[X_1 , \ldots , X_n]$-algebra $\tu{B}$ is a $T''$-coalescent, there exists a unique homomorphism $\sigma : R_{\approx T''} \to \tu{B}$ of condensed $\tu{A}[X_1 , \ldots , X_n]$-algebras such that the diagram
\[\begin{tikzcd}
	R & {R_{\approx T''}} \\
	{\tu{B}}
	\arrow["\psi", from=1-1, to=1-2]
	\arrow["{\sigma_0}"', from=1-1, to=2-1]
	\arrow["\sigma", from=1-2, to=2-1]
\end{tikzcd}\]
is commutative. Then the diagram
\[\begin{tikzcd}
	{\tu{A \langle X_1 , \dots , X_n \rangle}} \\
	{\tu{A}[X_1 , \ldots , X_n]} & R & {R_{\approx T''}} \\
	{\tu{A \langle X_1 , \dots , X_n \rangle}} & {\tu{B}}
	\arrow["\chi", from=1-1, to=2-3]
	\arrow["\phi", from=2-1, to=1-1]
	\arrow["\pi"', from=2-1, to=2-2]
	\arrow["\phi"', from=2-1, to=3-1]
	\arrow["\psi"', from=2-2, to=2-3]
	\arrow["{\sigma_0}"', from=2-2, to=3-2]
	\arrow["\sigma", from=2-3, to=3-2]
	\arrow["{\tu{p}}"', from=3-1, to=3-2]
\end{tikzcd}\]
is commutative. Since $\tu{B}$ is a $T''$-coalescent condensed $\tu{A}[X_1 , \ldots , X_n]$-algebra, \cref{prop:coalescence of polynomial rings} shows that the diagram
\[\begin{tikzcd}
	{\tu{A \langle X_1 , \dots , X_n \rangle}} & {R_{\approx T''}} \\
	{\tu{B}}
	\arrow["\chi", from=1-1, to=1-2]
	\arrow["{\tu{p}}"', from=1-1, to=2-1]
	\arrow["\sigma", from=1-2, to=2-1]
\end{tikzcd}\]
is commutative.

On the other hand, if $\sigma_1 : R_{\approx T''} \to \tu{B}$ is another homomorphism of condensed rings such that the diagram
\[\begin{tikzcd}
	{\tu{A \langle X_1 , \dots , X_n \rangle}} & {R_{\approx T''}} \\
	{\tu{B}}
	\arrow["\chi", from=1-1, to=1-2]
	\arrow["{\tu{p}}"', from=1-1, to=2-1]
	\arrow["{\sigma_1}", from=1-2, to=2-1]
\end{tikzcd}\]
is commutative, then the diagram
\[\begin{tikzcd}
	{\tu{A}[X_1 , \ldots , X_n]} & R \\
	{\tu{A \langle X_1 , \dots , X_n \rangle}} & {R_{\approx T''}} \\
	{\tu{B}}
	\arrow["\pi", from=1-1, to=1-2]
	\arrow["\phi"', from=1-1, to=2-1]
	\arrow["\psi", from=1-2, to=2-2]
	\arrow["\chi", from=2-1, to=2-2]
	\arrow["{\tu{p}}"', from=2-1, to=3-1]
	\arrow["{\sigma_1}", from=2-2, to=3-1]
\end{tikzcd}\]
is commutative, and $\sigma_1 : R_{\approx T''} \to \tu{B}$ is necessarily a homomorphism of condensed $\tu{A}[X_1 , \ldots , X_n]$-algebras. Then the definition of $\sigma_0$ shows that $\sigma_1 \of \psi = \sigma_0$. The the definition of $\sigma$ shows that $\sigma_1 = \sigma$.

\item
By \cref{cor:underbar of quotient map is still epic}, the map $\tu{A \langle X_1 , \ldots , X_n \rangle} \xto{\tu{p}} \tu{B}$ is an epimorphism when considered as a morphism in $\ub{CSet}$.

\item
By definition, the diagram is commutative. By \cref{prop:construction of the homomorphism chi}, the sequence
\[\begin{tikzcd}
	{\tu{{A \langle X_1 , \ldots , X_n \rangle}^{n}}} & {\tu{A \langle X_1 , \ldots , X_n \rangle}} & {R_{\approx T''}} & {0 .}
	\arrow["{\tu{m}}", from=1-1, to=1-2]
	\arrow["\chi", from=1-2, to=1-3]
	\arrow[from=1-3, to=1-4]
\end{tikzcd}\]
is exact in $\ub{CMod}_{\tu{A}[X_1 , \ldots , X_n]}$. By \cref{cor:underbar of quotient map is still epic}, the sequence
\[\begin{tikzcd}
	0 & {\tu{\overline{J \,}}} & {\tu{A \langle X_1 , \dots , X_n \rangle}} & {\tu{B}} & 0
	\arrow[from=1-1, to=1-2]
	\arrow["{\tu{q}}", from=1-2, to=1-3]
	\arrow["{\tu{p}}", from=1-3, to=1-4]
	\arrow[from=1-4, to=1-5]
\end{tikzcd}\]
is exact in $\ub{CAb}$. By \cref{prop:CMod R to CAb preserves and reflects exact sequences}, it is also exact in $\ub{CMod}_{\tu{A}[X_1 , \ldots , X_n]}$.

\item
The condition (c) is equivalent to the condition that $\sigma : R_{\approx T''} \to \tu{B}$ be an isomorphism in $\ub{CMod}_{\tu{A}[X_1 , \ldots , X_n]}$. Since the homomorphism $\sigma : R_{\approx T''} \to \tu{B}$ is an epimorphism in $\ub{CMod}_{\tu{A}[X_1 , \ldots , X_n]}$ by (2), it is an isomorphism in $\ub{CMod}_{\tu{A}[X_1 , \ldots , X_n]}$ if and only if its kernel in $\ub{CMod}_{\tu{A}[X_1 , \ldots , X_n]}$ is zero. By applying the Snake Lemma to the diagram of (3), we see that the kernel of $\sigma$ in $\ub{CMod}_{\tu{A}[X_1 , \ldots , X_n]}$ is isomorphic to the cokernel of $\tu{m}$ in $\ub{CMod}_{\tu{A}[X_1 , \ldots , X_n]}$. This is zero if and only if the condition (b) holds. Thus we have proved that the condition (c) is equivalent to the condition (b). 

On the other hand, let $J'$ be the ideal $J$ of $A \langle X_1, \ldots, X_n \rangle$ endowed with the quotient topology induced by the map $m : {A \langle X_1, \ldots, X_n \rangle}^{n} \to J$. Then the map $b : J' \to \overline{J \,}$, $x \mapsto x$ is continuous. The following diagram is commutative in $\ub{CMod}_{\tu{A}[X_1 , \ldots , X_n]}$.
\[\begin{tikzcd}
	0 & {\tu{\ker m}} & {\tu{{A \langle X_1 , \ldots , X_n \rangle}^{n}}} & {\tu{J'}} & 0 \\
	0 & {\tu{\ker m}} & {\tu{{A \langle X_1 , \ldots , X_n \rangle}^{n}}} & {\tu{\overline{J \,}}}
	\arrow[from=1-1, to=1-2]
	\arrow["\inc", from=1-2, to=1-3]
	\arrow[equals, from=1-2, to=2-2]
	\arrow["{\tu{m}}", from=1-3, to=1-4]
	\arrow[equals, from=1-3, to=2-3]
	\arrow[from=1-4, to=1-5]
	\arrow["{\tu{b}}", from=1-4, to=2-4]
	\arrow[from=2-1, to=2-2]
	\arrow["\inc"', from=2-2, to=2-3]
	\arrow["{\tu{m}}"', from=2-3, to=2-4]
\end{tikzcd}\]
By \cref{cor:underbar of quotient map is still epic}, the rows are exact in $\ub{CAb}$. By \cref{prop:CMod R to CAb preserves and reflects exact sequences}, the rows are also exact in $\ub{CMod}_{\tu{A}[X_1 , \ldots , X_n]}$. Therefore we conclude that the condition (b) holds if and only if $\tu{b} : \tu{J'} \to \tu{\overline{J \,}}$ is an isomorphism in $\ub{CMod}_{\tu{A}[X_1 , \ldots , X_n]}$. This is the case if and only if the map $\tu{b} : \tu{J'} \to \tu{\overline{J \,}}$ is an isomorphism in $\ub{CSet}$. Since both $J'$ and $\overline{J \,}$ are first countable, Proposition 4.4.2 of \cite{Kedlaya:note} shows that the map $\tu{b} : \tu{J'} \to \tu{\overline{J \,}}$ is an isomorphism in $\ub{CSet}$ if and only if the map $b : J' \to \overline{J \,}$ is a homeomorphism. If this holds, then the topology of $J'$ coincides with the subspace topology induced by the topology of $A \langle X_1, \ldots, X_n \rangle$. Therefore (a) holds. Conversely, if (a) holds, then \cref{prop:completeness of quotients of f c topological abelian groups} shows that $J'=J$ is already complete. Since $A \langle X_1, \ldots, X_n \rangle$ is Hausdorff, we conclude that $J$ is closed in $A \langle X_1, \ldots, X_n \rangle$. Therefore $J = \overline{J \,}$, and the map $b = \id_J : J' = J \to J = \overline{J \,}$ is a homeomorphism. Thus we have proved that the condition (b) is equivalent to the condition (a).
\end{enumerate}
\end{proof}

\begin{prop} \label{prop:construction of the homomorphism sigma prime} \;
\begin{enumerate}
\item 
There exists a unique homomorphism $\sigma' : (R_{\approx T''})_g \to \tu{B}_g$ of condensed rings such that the following diagram is commutative.
\[\begin{tikzcd}
	{R_{\approx T''}} & {(R_{\approx T''})_g} \\
	{\tu{B}} & {\tu{B}_g}
	\arrow["\can", from=1-1, to=1-2]
	\arrow["\sigma"', from=1-1, to=2-1]
	\arrow["{\sigma'}", from=1-2, to=2-2]
	\arrow["\can"', from=2-1, to=2-2]
\end{tikzcd}\]

\item
$\sigma' : (R_{\approx T''})_g \to \tu{B}_g$ is an epimorphism when considered as a morphism in $\ub{CSet}$.

\item If $\sigma : R_{\approx T''} \to \tu{B}$ is an isomorphism of condensed rings, then $\sigma' : (R_{\approx T''})_g \to \tu{B}_g$ is also an isomorphism of condensed rings.
\end{enumerate}
\end{prop}

\begin{proof}
Since $\sigma : R_{\approx T''} \to \tu{B}$ is a homomorphism of condensed $\tu{A}$-algebras and since any homomorphism $\sigma' : (R_{\approx T''})_g \to \tu{B}_g$ of condensed rings such that the diagram
\[\begin{tikzcd}
	{R_{\approx T''}} & {(R_{\approx T''})_g} \\
	{\tu{B}} & {\tu{B}_g}
	\arrow["\can", from=1-1, to=1-2]
	\arrow["\sigma"', from=1-1, to=2-1]
	\arrow["{\sigma'}", from=1-2, to=2-2]
	\arrow["\can"', from=2-1, to=2-2]
\end{tikzcd}\]
is commutative is necessarily a homomorphism of condensed $\tu{A}$-algebras $(R_{\approx T''})_g \to \tu{B}_g$, the assetions (1) and (3) follow from \cref{prop:universality of localization of condensed algebra}.

Let us prove the assertion (2). For this we use \cref{prop:monics and epics in CSet}. Let $S$ be any light profinite set. Let $t \in \tu{B}_g (S)$ be any element. Since $\tu{B}_g (S) = \big( \tu{B}(S) \big)_g$, there exists an element $t_0 \in \tu{B}(S)$ and an $n \in \N$ such that $t = t_0 / g^n$. By \cref{prop:construction of the homomorphism sigma}, the map $\sigma : R_{\approx T''} \to \tu{B}$ is an epimorphism when considered as a morphism in $\ub{CSet}$. Then \cref{prop:monics and epics in CSet} shows that there exists a cover $(S_i \xto{c_i} S)_{i \in I}$ of $S$ in $\ub{Prof}$ such that for every $i \in I$, there exists an element $u_i \in R_{\approx T''}(S_i)$ such that $\sigma_{S_i}(u_i) = \tu{B}(c_i)(t_0)$. Then we have an element $u_i / g^n \in (R_{\approx T''}(S_i))_g = (R_{\approx T''})_g (S_i)$ for each $i \in I$, and 
\begin{equation}
\sigma'_{S_i}(u_i / g^n) = \sigma_{S_i}(u_i) / g^n = \tu{B}(c_i)(t_0) / g^n = \tu{B}(c_i)(t_0 / g^n) = \tu{B}(c_i)(t)
\end{equation}
for each $i \in I$. Then it follows from \cref{prop:monics and epics in CSet} that 
$\sigma' : (R_{\approx T''})_g \to \tu{B}_g$ is an epimorphism when considered as a morphism in $\ub{CSet}$.
\end{proof}

\subsubsection{The homomorphism $\tau$}

\begin{lem}
The subset $r(T)$ of $B$ generates an open ideal of $B$.
\end{lem}

\begin{proof}
Let $(A_0, I_0)$ be a couple of definition of $A$. Since $T$ generates an open ideal of $A$, there exists an $e \in \N$ such that $I_0^e \sub T \cdot A$. Therefore
\begin{equation}
I_0^e \cdot A_0 \langle X_1, \ldots, X_n \rangle \sub T \cdot A \langle X_1, \ldots, X_n \rangle
\end{equation}
in $A \langle X_1, \ldots, X_n \rangle$. Since $(A_0 \langle X_1, \ldots, X_n \rangle, I_0 \cdot A_0 \langle X_1, \ldots, X_n \rangle)$ is a couple of definition of $A \langle X_1, \ldots, X_n \rangle$, we conclude that $T = j(T)$ generates an open ideal of $A \langle X_1, \ldots, X_n \rangle$. Since the map $p : A \langle X_1, \ldots, X_n \rangle \to B$ is surjective and open, it follows that $p(j(T)) = r(T)$ generates an open ideal of $B$.
\end{proof}

\begin{nt} \;
\begin{enumerate}
\item $k_1 : B \to B \left( \frac{r(T)}{r(g)} \right)$ denotes the canonical homomorphism.

\item $k_2 : B \left( \frac{r(T)}{r(g)} \right) \to B \left \langle \frac{r(T)}{r(g)} \right \rangle$ denotes the canonical homomorphism into the completion.

\item We write $k : B \to B \left \langle \frac{r(T)}{r(g)} \right \rangle$ for the composition $k_2 \of k_1$.
\end{enumerate}
\end{nt}

\begin{lem} \label{lem:properties of localizations of B} \;
\begin{enumerate}
\item $k_1 : B \to B \left( \frac{r(T)}{r(g)} \right)$ is an open map.

\item $B \left( \frac{r(T)}{r(g)} \right)$ is complete (but not necessarily Hausdorff).

\item $k_2 : B \left( \frac{r(T)}{r(g)} \right) \to B \left \langle \frac{r(T)}{r(g)} \right \rangle$ is surjective and open. Its kernel is equal to the closure $\overline{\{0\}}$ of $\{0\}$ in $B \left( \frac{r(T)}{r(g)} \right)$.

\item The following conditions are equivalent.
\begin{enumerate}
\item The kernel of $k_1 : B \to B \left( \frac{r(T)}{r(g)} \right)$ is a closed ideal of $B$.
\item $B \left( \frac{r(T)}{r(g)} \right)$ is Hausdorff.
\item $k_2 : B \left( \frac{r(T)}{r(g)} \right) \to B \left \langle \frac{r(T)}{r(g)} \right \rangle$ is an isomorphism of topological rings.
\end{enumerate}

\end{enumerate}
\end{lem}

\begin{proof}~
\begin{enumerate}
\item
Let $(A_0, I_0)$ be a couple of definition of $A$. Then $(A_0 \langle X_1, \ldots, X_n \rangle, j(I_0) \cdot A_0 \langle X_1, \ldots, X_n \rangle)$ is a couple of definition of $A \langle X_1, \ldots, X_n \rangle$. Since the map $p : A \langle X_1, \ldots, X_n \rangle \to B$ is surjective and open, the pair $(p(A_0 \langle X_1, \ldots, X_n \rangle), r(I_0) \cdot p(A_0 \langle X_1, \ldots, X_n \rangle))$ is a couple of definition of $B$. Let us write $B_0 := p(A_0 \langle X_1, \ldots, X_n \rangle)$. Let $B_1$ be the subring of $B \left( \frac{r(T)}{r(g)} \right)$ generated by $k_1(B_0) \cup \set{k_1(r(f_i)) \cdot k_1(r(g))^{-1} }{ 1 \leq i \leq n }$. Then $(B_1, k_1(r(I_0)) \cdot B_1)$ is a couple of definition of $B \left( \frac{r(T)}{r(g)} \right)$. On the other hand, for every $1 \leq i \leq n$, the definition of $B$ shows that $r(f_i) = r(g) \cdot p(X_i)$ in $B$. Therefore $k_1(r(f_i)) \cdot k_1(r(g))^{-1} = k_1(p(X_i))$ in $B \left( \frac{r(T)}{r(g)} \right)$. Since $p(X_i) \in B_0$ by definition, we have $k_1(p(X_i)) \in k_1(B_0)$. Since this holds for every $1 \leq i \leq n$, we concude that $B_1 = k_1(B_0)$. Therefore $( k_1(B_0), k_1(r(I_0) \cdot B_0) )$ is a couple of definition of $B \left( \frac{r(T)}{r(g)} \right)$. This immediately implies that the map $k_1 : B \to B \left( \frac{r(T)}{r(g)} \right)$ is an open map.

\item
By (1), the subring $k_1(B)$ is an open subring of $B \left( \frac{r(T)}{r(g)} \right)$ which is isomorphic to the quotient $B / \ker k_1$ as topological rings. Using the same proof as that of \cref{prop:completeness of quotients of f c topological abelian groups}, one can check that the quotient $B / \ker k_1$ is complete. Therefore the open subring $k_1(B)$ of $B \left( \frac{r(T)}{r(g)} \right)$ is also complete. Then $B \left( \frac{r(T)}{r(g)} \right)$ is itself complete by \cite{Bourbaki:top1}, Chapter III, \S 3.3, Proposition 4.

\item
It suffices to show that the quotient map $k_3 : B \left( \frac{r(T)}{r(g)} \right) \to B \left( \frac{r(T)}{r(g)} \right) \Big/ \, \overline{\{0\}}$ coincides with the Hausdorff completion of $B \left( \frac{r(T)}{r(g)} \right)$. First of all, $B \left( \frac{r(T)}{r(g)} \right) \Big/ \, \overline{\{0\}}$ is Hausdorff, and the map $k_3 : B \left( \frac{r(T)}{r(g)} \right) \to B \left( \frac{r(T)}{r(g)} \right) \Big/ \, \overline{\{0\}}$ has a dense image since it is surjective. Next, from the fact that $B \left( \frac{r(T)}{r(g)} \right)$ is complete by (2), one deduces that $B \left( \frac{r(T)}{r(g)} \right) \Big/ \, \overline{\{0\}}$ is complete. This can be proven in the same way as the proof of \cref{prop:completeness of quotients of f c topological abelian groups}. In addition, for every open subgroup $U$ of $B \left( \frac{r(T)}{r(g)} \right)$, one has $k_3^{-1}(k_3(U)) = U$ since $\overline{\{0\}} \sub U$. This shows that the topology of $B \left( \frac{r(T)}{r(g)} \right)$ is equal to the initial topology induced by $k_3 : B \left( \frac{r(T)}{r(g)} \right) \to B \left( \frac{r(T)}{r(g)} \right) \Big/ \, \overline{\{0\}}$. Therefore we conclude that the map $k_3 : B \left( \frac{r(T)}{r(g)} \right) \to B \left( \frac{r(T)}{r(g)} \right) \Big/ \, \overline{\{0\}}$ coincides with the Hausdorff completion of $B \left( \frac{r(T)}{r(g)} \right)$.

\item
By (3), the condition (c) holds if and only if $\overline{\{0\}} = \{0\}$. This condition is equivalent to (b).

If (b) holds, then (a) holds since the kernel of $k_1$ is equal to the inverse image $k_1^{-1}(\{0\})$ of the closed set $\{0\}$ of $B \left( \frac{r(T)}{r(g)} \right)$.

Conversely, suppose that (a) holds. We use the notation introduced in the proof of (1). Then
\begin{equation}
\overline{\{0\}} = \bigcap_{e \in \N} k_1(r(I_0)^e \cdot B_0)
\end{equation}
since $( k_1(B_0), k_1(r(I_0) \cdot B_0) )$ is a couple of definition of $B \left( \frac{r(T)}{r(g)} \right)$. Then 
\begin{equation}
k_1^{-1} \left( \overline{\{0\}} \right)
= \bigcap_{e \in \N} k_1^{-1} k_1(r(I_0)^e \cdot B_0)
= \bigcap_{e \in \N} \big( (r(I_0)^e \cdot B_0) + \ker k_1 \big) .
\end{equation}
Since $(B_0, r(I_0) \cdot B_0)$ is a couple of definition of $B$, it follows that
\begin{equation}
\bigcap_{e \in \N} \big( (r(I_0)^e \cdot B_0) + \ker k_1 \big) = \overline{\ker k_1} .
\end{equation}
Since we assumed the condition (a), we have
\begin{equation}
\overline{\ker k_1} = \ker k_1 .
\end{equation}
Consequently, 
\begin{equation}
k_1^{-1} \left( \overline{\{0\}} \right) = \ker k_1 .
\end{equation}
Since the map $k_1 : B \to B \left( \frac{r(T)}{r(g)} \right)$ is equal to the canonical map $B \to B_{r(g)}$, we have
\begin{equation}
\overline{\{0\}}
= k_1 \left( k_1^{-1} \left( \overline{\{0\}} \right) \right) \cdot B \left( \frac{r(T)}{r(g)} \right) 
= k_1 \left( \ker k_1 \right) \cdot B \left( \frac{r(T)}{r(g)} \right) = 0 .
\end{equation}
Therefore (b) holds.
\end{enumerate}
\end{proof}

\begin{prop} \label{prop:construction of the homomorphism tau} \;
\begin{enumerate}
\item
There exists a unique homomorphism $\tau : \tu{B}_g \to \tu{B \left \langle \frac{r(T)}{r(g)} \right \rangle}$ of condensed rings such that the following diagram is commutative.
\[\begin{tikzcd}
	{\tu{B}} & {\tu{B}_g} \\
	{\tu{B \left \langle \frac{r(T)}{r(g)} \right \rangle}}
	\arrow["\can", from=1-1, to=1-2]
	\arrow["{\tu{k}}"', from=1-1, to=2-1]
	\arrow["\tau", from=1-2, to=2-1]
\end{tikzcd}\]

\item
$\tau : \tu{B}_g \to \tu{B \left \langle \frac{r(T)}{r(g)} \right \rangle}$ is an epimorphism when considered as a morphism in $\ub{CSet}$.

\item
If the kernel of $k_1 : B \to B \left( \frac{r(T)}{r(g)} \right)$ is a closed ideal of $B$, then $\tau : \tu{B}_g \to \tu{B \left \langle \frac{r(T)}{r(g)} \right \rangle}$ is an isomorphism of condensed rings.
\end{enumerate}
\end{prop}

\begin{proof}~
\begin{enumerate}
\item
Since the element $k(r(g)) \in B \left \langle \frac{r(T)}{r(g)} \right \rangle$ is invertible in $B \left \langle \frac{r(T)}{r(g)} \right \rangle$, \cref{prop:universality of localization of condensed algebra} shows that there exists a unique homomorphism $\tau : \tu{B}_g \to \tu{B \left \langle \frac{r(T)}{r(g)} \right \rangle}$ of condensed $\tu{A}$-algebras such that the following diagram is commutative.
\[\begin{tikzcd}
	{\tu{B}} & {\tu{B}_g} \\
	{\tu{B \left \langle \frac{r(T)}{r(g)} \right \rangle}}
	\arrow["\can", from=1-1, to=1-2]
	\arrow["{\tu{k}}"', from=1-1, to=2-1]
	\arrow["\tau", from=1-2, to=2-1]
\end{tikzcd}\]

On the other hand, suppose that $\tau' : \tu{B}_g \to \tu{B \left \langle \frac{r(T)}{r(g)} \right \rangle}$ is another homomorphism of condensed rings such that the diagram
\[\begin{tikzcd}
	{\tu{B}} & {\tu{B}_g} \\
	{\tu{B \left \langle \frac{r(T)}{r(g)} \right \rangle}}
	\arrow["\can", from=1-1, to=1-2]
	\arrow["{\tu{k}}"', from=1-1, to=2-1]
	\arrow["{\tau'}", from=1-2, to=2-1]
\end{tikzcd}\]
 is commutative. Then $\tau' : \tu{B}_g \to \tu{B \left \langle \frac{r(T)}{r(g)} \right \rangle}$ is necessarily a homomorphism of $\tu{A}$-algebras. By definition of $\tau$, we have $\tau' = \tau$.

\item
We use \cref{prop:monics and epics in CSet}. Let $S$ be any light profinite set. Let $b \in \tu{B \left \langle \frac{r(T)}{r(g)} \right \rangle}(S) = \mathrm{Cont} \left( S, B \left \langle \frac{r(T)}{r(g)} \right \rangle \right)$. Since $S$ is compact, its continuous image $b(S)$ in $B \left \langle \frac{r(T)}{r(g)} \right \rangle$ is also compact. On the other hand, (3) of \cref{lem:properties of localizations of B} shows that
\begin{equation}
B \left \langle \frac{r(T)}{r(g)} \right \rangle =
\bigcup_{e \in \N} k(B) \cdot k(r(g))^{-e}.
\end{equation}
Moreover, (1) and (3) of \cref{lem:properties of localizations of B} show that $k(B) \cdot k(r(g))^{-e}$ is open in $B \left \langle \frac{r(T)}{r(g)} \right \rangle$ for every $e \in \N$. Consequently, there exists an $e \in \N$ such that
\begin{equation}
b(S) \sub k(B) \cdot k(r(g))^{-e} .
\end{equation}
Then $k(r(g))^e \cdot b(S) \sub k(B)$, and $k(r(g))^e \cdot b(S)$ is compact since $b(S)$ is compact and $k(r(g))$ is an invertible element of $B \left \langle \frac{r(T)}{r(g)} \right \rangle$. Moreover, the map $k : B \to k(B)$ is surjective and open by (1) and (3) of \cref{lem:properties of localizations of B}, and has a closed kernel since $B \left \langle \frac{r(T)}{r(g)} \right \rangle$ is Hausdorff. Then \cref{prop:lifting compact sets} shows that there exists a subspace $\tilde{S}$ of $B$ which is a light profinite set and satisfies $k(\tilde{S}) = k(r(g))^e \cdot b(S)$. Let $S' := S \times_{(k(r(g))^e \cdot b(S))} \tilde{S}$ be the fiber product in $\ub{Top}$, with projections $c : S' \to S$ and $d : S' \to \tilde{S}$. More precisely, the left square of the following diagram is a pullback in $\ub{Top}$. 
\[\begin{tikzcd}
	{S'} &&& {\tilde{S}} & B \\
	S & {b(S)} && {k(r(g))^e \cdot b(S)} & {B \left \langle \frac{r(T)}{r(g)} \right \rangle}
	\arrow["d", from=1-1, to=1-4]
	\arrow["c"', from=1-1, to=2-1]
	\arrow["\inc", hook, from=1-4, to=1-5]
	\arrow["k", two heads, from=1-4, to=2-4]
	\arrow["k", from=1-5, to=2-5]
	\arrow["b"', from=2-1, to=2-2]
	\arrow["{x \mapsto k(r(g))^e \cdot x}"', from=2-2, to=2-4]
	\arrow["\inc"', hook, from=2-4, to=2-5]
\end{tikzcd}\]
By \cref{prop:closure properties of light profinite sets}, the topological space $S'$ is a light profinite set. The map $c : S' \to S$ is surjective since the map $k : \tilde{S} \to k(r(g))^e \cdot b(S)$ is surjective. Let $d' : S' \to B$ be the composition of $d :  S' \to \tilde{S}$ and the inclusion $\tilde{S} \mon B$. Then we have $d' \in \mathrm{Cont}(S',B) = \tu{B}(S')$ and $\tu{k}_{S'}(d') = k \of d' = k(r(g))^e \cdot (b \of c)$. Therefore the element $d'/g^e \in \tu{B}(S')_g = \tu{B}_g(S')$ satisfies 
\begin{equation}
\tau_{S'}(d'/g^e) = \tu{k}_{S'}(d') \cdot k(r(g))^{-e} = k(r(g))^e \cdot (b \of c) \cdot k(r(g))^{-e} = b \of c = \tu{B}(c)(b) .
\end{equation}
By \cref{prop:monics and epics in CSet}, we conclude that $\tau : \tu{B}_g \to \tu{B \left \langle \frac{r(T)}{r(g)} \right \rangle}$ is an epimorphism when considered as a morphism in $\ub{CSet}$.

\item
By (2), it suffices to show that $\tau : \tu{B}_g \to \tu{B \left \langle \frac{r(T)}{r(g)} \right \rangle}$ is a monomorphism when considered as a morphism in $\ub{CSet}$. For this we use \cref{prop:monics and epics in CSet}. Let $S$ be any light profinite set. We prove that the map $\tau_S : \tu{B}_g(S) \to \tu{B \left \langle \frac{r(T)}{r(g)} \right \rangle}(S)$ is injective. Since this map is a homomorphism of commutative unital rings, it suffices to show that the kernel of $\tau_S : \tu{B}_g(S) \to \tu{B \left \langle \frac{r(T)}{r(g)} \right \rangle}(S)$ is equal to the zero ideal. Let $b/g^e \in \tu{B}(S)_g = \tu{B}_g(S)$, where $b \in \tu{B}(S) = \mathrm{Cont}(S,B)$ and $e \in \N$. Suppose that $\tau_S(b/g^e) = 0$. We prove that $b/g^e = 0$. We have
\begin{equation}
0 = \tau_S(b/g^e) = \tu{k}_S(b) \cdot k(r(g))^{-e} = (k \of b) \cdot k(r(g))^{-e}.
\end{equation}
Since $k(r(g))^{-e} \in B \left \langle \frac{r(T)}{r(g)} \right \rangle$ is an invertible element of $B \left \langle \frac{r(T)}{r(g)} \right \rangle$, it follows that $k \of b = 0$ in $\mathrm{Cont}\left(S, B \left \langle \frac{r(T)}{r(g)} \right \rangle \right)$. Therefore $b(S)$ is contained in the kernel of $k : B \to B \left \langle \frac{r(T)}{r(g)} \right \rangle$. $\ker k$ is a complete Hausdorff first countable topological abelian group since it is a closed additive subgroup of the complete Hausdorff first countable topological abelian group $B$. Then \cref{cor:c H f c top ab is Baire} shows that $\ker k$ is a Baire space. On the other hand, since the kernel of $k_1 : B \to B \left( \frac{r(T)}{r(g)} \right)$ is assumed to be closed in $B$, (4) of \cref{lem:properties of localizations of B} shows that $k_2 : B \left( \frac{r(T)}{r(g)} \right) \to B \left \langle \frac{r(T)}{r(g)} \right \rangle$ is an isomorphism of topological rings. Since $k = k_2 \of k_1$, it follows that the kernel of $k$ coincides with the kernel of $k_1$. Thus if we write
\begin{equation}
G_d := \set{x \in B}{ r(g)^d \cdot x = 0}
\end{equation}
for $d \in \N$, we have
\begin{equation}
\ker k = \ker k_1 = \bigcup_{d \in \N} G_d .
\end{equation}
Moreover, for every $d \in \N$, the additive subgroup $G_d$ is closed in $B$ since $B$ is Hausdorff and the multiplication of $B$ is continuous. Therefore $G_d$ is also a closed additive subgroup of $\ker k$ for every $d \in \N$. Consequently, there exists a $d_0 \in \N$ such that the interior of $G_{d_0}$ in $\ker k$ is nonempty. From the fact that $\ker k$ is a non-archimedean abelian group, one easily checks that $G_d$ is an open additive subgroup of $\ker k$ for every $d \in \N$ with $d \geq d_0$. Therefore we have an open covering
\begin{equation}
\ker k = \bigcup_{\substack{d \in \N \\ d \geq d_0}} G_d
\end{equation}
of the topological space $\ker k$. Since $b(S)$ is contained in $\ker k$ and since $b(S)$ is compact, we conclude that there exists a $d \in \N$ such that
\begin{equation}
b(S) \sub G_d .
\end{equation}
This means that $r(g)^d \cdot b = 0$ in $\mathrm{Cont}(S,B)$. Consequently, we have
\begin{equation}
b/g^e = (r(g)^d \cdot b) / g^{d+e} = 0
\end{equation}
in $\tu{B}(S)_g = \tu{B}_g(S)$. This completes the proof.
\end{enumerate}
\end{proof}

\subsubsection{Expressing $\rho$ in terms of $\sigma'$ and $\tau$}

\begin{lem} \label{lem:construction of the homomorphism u}
There exists a unique isomorphism $u : B \left \langle \frac{r(T)}{r(g)} \right \rangle \xto{\sim} A \left \langle \frac{T}{g} \right \rangle$ of topological rings such that the following diagram is commutative.
\[\begin{tikzcd}
	A & B \\
	{A \left \langle \frac{T}{g} \right \rangle} & {B \left \langle \frac{r(T)}{r(g)} \right \rangle}
	\arrow["r", from=1-1, to=1-2]
	\arrow["l"', from=1-1, to=2-1]
	\arrow["k", from=1-2, to=2-2]
	\arrow["u", from=2-2, to=2-1]
\end{tikzcd}\]
\end{lem}

\begin{proof}
In order to prove this assertion, it suffices to show that the continuous homomorphism $k \of r : A \to B \left \langle \frac{r(T)}{r(g)} \right \rangle$ has the same universal property as that of $l : A \to A \left \langle \frac{T}{g} \right \rangle$ described in Proposition II.3.4.6 of \cite{Morel:note}. More precisely, it suffices to prove the following assertions. 
\begin{enumerate}
\item
The topological ring $B \left \langle \frac{r(T)}{r(g)} \right \rangle$ is a complete Hausdorff non-archimedean ring. The map $k \of r : A \to B \left \langle \frac{r(T)}{r(g)} \right \rangle$ is a continuous ring homomorphism. The element $k(r(g)) \in B \left \langle \frac{r(T)}{r(g)} \right \rangle$ is invertible in $B \left \langle \frac{r(T)}{r(g)} \right \rangle$ and the set $\set{k(r(f)) \cdot k(r(g))^{-1}}{f \in T}$ is power-bounded in $B \left \langle \frac{r(T)}{r(g)} \right \rangle$.

\item
Let $C$ be any complete Hausdorff non-archimedean ring. Let $d : A \to C$ be any continuous ring homomorphism. Suppose that the element $d(g) \in C$ is invertible in $C$ and that the set $\set{d(f) \cdot d(g)^{-1}}{f \in T}$ is power-bounded in $C$. Then there exists a unique continuous ring homomorphism $e : B \left \langle \frac{r(T)}{r(g)} \right \rangle \to C$ such that the following diagram is commutative.
\[\begin{tikzcd}
	A & C \\
	{B \left \langle \frac{r(T)}{r(g)} \right \rangle}
	\arrow["d", from=1-1, to=1-2]
	\arrow["{k \of r}"', from=1-1, to=2-1]
	\arrow["e"', from=2-1, to=1-2]
\end{tikzcd}\]
\end{enumerate}
The assertion (1) directly follows from the proprerty of $k : B \to B \left \langle \frac{r(T)}{r(g)} \right \rangle$. Let us prove the assertion (2). Let $C$ be any complete Hausdorff non-archimedean ring. Let $d : A \to C$ be any continuous ring homomorphism. Suppose that the element $d(g) \in C$ is invertible in $C$ and that the set $\set{d(f) \cdot d(g)^{-1}}{f \in T}$ is power-bounded in $C$. By Corollary II.3.3.4 of \cite{Morel:note} or Corollary 5.50 of \cite{Wedhorn:note}, there exists a unique continuous ring homomorphism $d_1 : A \langle X_1, \ldots, X_n \rangle \to C$ such that $d_1(X_i) = d(f_i) \cdot d(g)^{-1}$ for each $1 \leq i \leq n$ and the diagram
\[\begin{tikzcd}
	A & C \\
	{A \langle X_1, \ldots, X_n \rangle}
	\arrow["d", from=1-1, to=1-2]
	\arrow["j"', from=1-1, to=2-1]
	\arrow["{d_1}"', from=2-1, to=1-2]
\end{tikzcd}\]
is commutative. Then we have
\begin{equation}
d_1(f_i - g X_i) = d(f_i) - d(g) \cdot d(f_i) \cdot d(g)^{-1} = 0
\end{equation}
for every $1 \leq i \leq n$. Therefore the ideal $J$ of $A \langle X_1, \ldots, X_n \rangle$ is contained in the kernel of $d_1 : A \langle X_1, \ldots, X_n \rangle \to C$. Since $C$ is Hausdorff, the kernel of $d_1 : A \langle X_1, \ldots, X_n \rangle \to C$ is closed in $A \langle X_1, \ldots, X_n \rangle$. Therefore the closure $\overline{J \,}$ of $J$ is also contained in the kernel of $d_1 : A \langle X_1, \ldots, X_n \rangle \to C$. By the definition of $B$, there exists a unique continuous ring homomorphism $d_2 : B \to C$ such that the diagram
\[\begin{tikzcd}
	{A \langle X_1, \ldots, X_n \rangle} & C \\
	B
	\arrow["{d_1}", from=1-1, to=1-2]
	\arrow["p"', from=1-1, to=2-1]
	\arrow["{d_2}"', from=2-1, to=1-2]
\end{tikzcd}\]
is commutative. Then by assumption the element $d_2(r(g)) = d(g) \in C$ is invertible in $C$ and the set $\set{d_2(r(f)) \cdot d_2(r(g))^{-1}}{f \in T} = \set{d(f) \cdot d(g)^{-1}}{f \in T}$ is power-bounded in $C$. Therefore the universal property of $k : B \to B \left \langle \frac{r(T)}{r(g)} \right \rangle$ (\cite{Morel:note}, Proposition II.3.4.6) shows that there exists a unique continuous ring homomorphism $e : B \left \langle \frac{r(T)}{r(g)} \right \rangle \to C$ such that the diagram
\[\begin{tikzcd}
	B & C \\
	{B \left \langle \frac{r(T)}{r(g)} \right \rangle}
	\arrow["{d_2}", from=1-1, to=1-2]
	\arrow["k"', from=1-1, to=2-1]
	\arrow["e"', from=2-1, to=1-2]
\end{tikzcd}\]
is commutative. Then the diagram
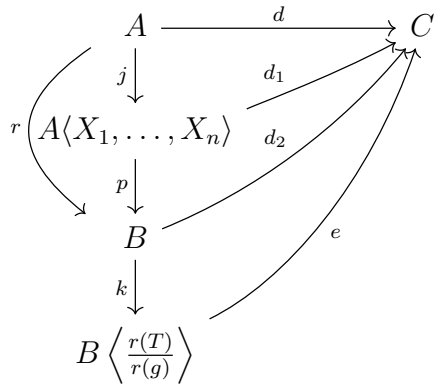
\[\begin{tikzcd}
	A && C \\
	{A \langle X_1, \ldots, X_n \rangle} \\
	B \\
	{B \left \langle \frac{r(T)}{r(g)} \right \rangle}
	\arrow["d", from=1-1, to=1-3]
	\arrow["j"', from=1-1, to=2-1]
	\arrow["r"', shift right=3, curve={height=30pt}, from=1-1, to=3-1]
	\arrow["{d_1}", curve={height=6pt}, from=2-1, to=1-3]
	\arrow["p"', from=2-1, to=3-1]
	\arrow["{d_2}", curve={height=12pt}, from=3-1, to=1-3]
	\arrow["k"', from=3-1, to=4-1]
	\arrow["e"', curve={height=24pt}, from=4-1, to=1-3]
\end{tikzcd}\]
is commutative.

On the other hand, suppose that $e' : B \left \langle \frac{r(T)}{r(g)} \right \rangle \to C$ is another continuous ring homomorphism such that the diagram
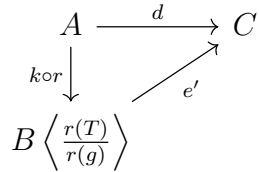
\[\begin{tikzcd}
	A & C \\
	{B \left \langle \frac{r(T)}{r(g)} \right \rangle}
	\arrow["d", from=1-1, to=1-2]
	\arrow["{k \of r}"', from=1-1, to=2-1]
	\arrow["{e'}"', from=2-1, to=1-2]
\end{tikzcd}\]
is commutative. We prove that $e' = e$. The diagram
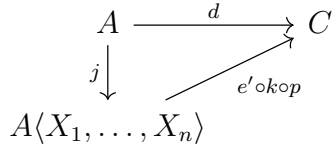
\[\begin{tikzcd}
	A & C \\
	{A \langle X_1, \ldots, X_n \rangle}
	\arrow["d", from=1-1, to=1-2]
	\arrow["j"', from=1-1, to=2-1]
	\arrow["{e' \of k \of p}"', from=2-1, to=1-2]
\end{tikzcd}\]
is commutative. Furthermore, for each $1 \leq i \leq n$, we have $p(f_i - g X_i) = 0$ in $B$ by definition. Therefore $d(f_i) = e' \of k \of p (f_i) = e' \of k \of p (g X_i) = d(g) \cdot \big( e' \of k \of p(X_i) \big)$. Thus $e' \of k \of p(X_i) = d(f_i) \cdot d(g)^{-1}$. Then the definition of $d_1$ shows that $e' \of k \of p = d_1$. In other words, the diagram
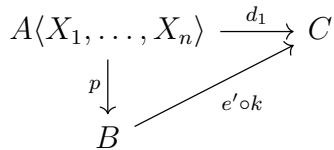
\[\begin{tikzcd}
	{A \langle X_1, \ldots, X_n \rangle} & C \\
	B
	\arrow["{d_1}", from=1-1, to=1-2]
	\arrow["p"', from=1-1, to=2-1]
	\arrow["{e' \of k}"', from=2-1, to=1-2]
\end{tikzcd}\]
is commutative. Then the definition of $d_2$ shows that $e' \of k = d_2$. In other words, the diagram
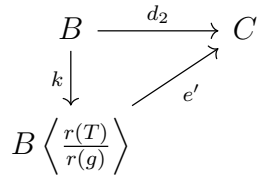
\[\begin{tikzcd}
	B & C \\
	{B \left \langle \frac{r(T)}{r(g)} \right \rangle}
	\arrow["{d_2}", from=1-1, to=1-2]
	\arrow["k"', from=1-1, to=2-1]
	\arrow["{e'}"', from=2-1, to=1-2]
\end{tikzcd}\]
is commutative. Then the definition of $e$ shows that $e' = e$. This completes the proof.
\end{proof}

\begin{lem} \label{lem:construction of the homomorphism omega}
There exists a unique isomorphism $\omega : (\tu{A}_g)_{\approx T'} \xto{\sim} (R_{\approx T''})_g$ of condensed rings such that the following diagram is commutative.
\[\begin{tikzcd}
	{\tu{A}} && {(\tu{A}_g)_{\approx T'}} \\
	{\tu{A \langle X_1 , \dots , X_n \rangle}} & {R_{\approx T''}} & {(R_{\approx T''})_g}
	\arrow["\can", from=1-1, to=1-3]
	\arrow["{\tu{j}}"', from=1-1, to=2-1]
	\arrow["\omega", from=1-3, to=2-3]
	\arrow["\chi"', from=2-1, to=2-2]
	\arrow["\can"', from=2-2, to=2-3]
\end{tikzcd}\]
\end{lem}

\begin{proof}
By \cref{prop:construction of the homomorphism chi}, the following diagram is commutative.
\[\begin{tikzcd}
	{\tu{A}} & {\tu{A}[X_1, \ldots, X_n]} & R \\
	& {\tu{A \langle X_1, \ldots, X_n \rangle}} & {R_{\approx T''}} & {(R_{\approx T''})_g}
	\arrow["\iota", from=1-1, to=1-2]
	\arrow["{\tu{j}}"', from=1-1, to=2-2]
	\arrow["\pi", from=1-2, to=1-3]
	\arrow["\phi", from=1-2, to=2-2]
	\arrow["\psi", from=1-3, to=2-3]
	\arrow["\chi"', from=2-2, to=2-3]
	\arrow["\can"', from=2-3, to=2-4]
\end{tikzcd}\]
Therefore it suffices to prove that there exists a unique isomorphism $\omega : (\tu{A}_g)_{\approx T'} \xto{\sim} (R_{\approx T''})_g$ of condensed rings such that the following diagram is commutative.
\[\begin{tikzcd}
	{\tu{A}} &&& {(\tu{A}_g)_{\approx T'}} \\
	{\tu{A}[X_1, \ldots, X_n]} & R & {R_{\approx T''}} & {(R_{\approx T''})_g}
	\arrow["\can", from=1-1, to=1-4]
	\arrow["\iota"', from=1-1, to=2-1]
	\arrow["\omega", from=1-4, to=2-4]
	\arrow["\pi"', from=2-1, to=2-2]
	\arrow["\psi"', from=2-2, to=2-3]
	\arrow["\can"', from=2-3, to=2-4]
\end{tikzcd}\]
\cref{prop:localization and coalescence commute} shows that the homomorphism $\tu{A} \xto{\can} (\tu{A}_g)_{\approx T'}$ of condensed $\tu{A}$-algebras has a certain universal property. Therefore, in order to prove the assertion above, it suffices to prove that the composition
\[\begin{tikzcd}
	{\xi : \tu{A}} & {\tu{A}[X_1, \ldots, X_n]} & R & {R_{\approx T''}} & {(R_{\approx T''})_g}
	\arrow["\iota", from=1-1, to=1-2]
	\arrow["\pi", from=1-2, to=1-3]
	\arrow["\psi", from=1-3, to=1-4]
	\arrow["\can", from=1-4, to=1-5]
\end{tikzcd}\]
has the same universal property as that of $\tu{A} \xto{\can} (\tu{A}_g)_{\approx T'}$
described in \cref{prop:localization and coalescence commute}. More precisely, it suffices to prove the following assertions.
\begin{enumerate}
\item
The condensed $\tu{A}$-algebra $(R_{\approx T''})_g$ is $T'$-coalescent. The element $\xi_*(g) \in (R_{\approx T''})_g(*)$ is invertible in $(R_{\approx T''})_g(*)$. 

\item
Let $R'$ be any $T'$-coalescent condensed $\tu{A}$-algebra. Let $\xi' : \tu{A} \to R'$ be any homomorphism of condensed $\tu{A}$-algebras such that the element $\xi'_*(g) \in R'(*)$ is invertible in $R'(*)$. Then there exists a unique homomorphism $\eta : (R_{\approx T''})_g \to R'$ of condensed $\tu{A}$-algebras such that the diagram
\[\begin{tikzcd}
	{\tu{A}} & {R'} \\
	{(R_{\approx T''})_g}
	\arrow["{\xi'}", from=1-1, to=1-2]
	\arrow["\xi"', from=1-1, to=2-1]
	\arrow["\eta"', from=2-1, to=1-2]
\end{tikzcd}\]
is commutative.
\end{enumerate}

In this proof, let us write $\lambda : R_{\approx T''} \to (R_{\approx T''})_g$ for the canonical homomorphism into the localization. Since $(R_{\approx T''})_g(*) = ((R_{\approx T''})(*))_g$ by definition, the element $\xi_*(g) \in (R_{\approx T''})_g(*)$ is invertible in $(R_{\approx T''})_g(*)$. Furthermore, by \cref{prop:localization and monoid algebra of coalescent algebras}, the condensed $\tu{A}[X_1, \ldots, X_n]$-algebra $(R_{\approx T''})_g$ is $T''$-coalescent. Therefore for each $1 \leq i \leq n$, the condensed $(R_{\approx T''})_g$-algebra $(R_{\approx T''})_g$ is $(\lambda_* \psi_* \pi_*(X_i) ) / 1$-coalescent. Since the element $\xi_*(g) = \lambda_* \psi_* \pi_*(g) \in (R_{\approx T''})_g(*)$ is invertible in $(R_{\approx T''})_g(*)$, it follows that the condensed $(R_{\approx T''})_g$-algebra $(R_{\approx T''})_g$ is $(\lambda_* \psi_* \pi_*(g \cdot X_i)) / (\lambda_* \psi_* \pi_*(g))$-coalescent. On the other hand, the definition of $R$ shows that $\pi_*(g \cdot X_i) = \pi_*(f_i)$ in $R(*)$. It follows that the condensed $(R_{\approx T''})_g$-algebra $(R_{\approx T''})_g$ is $(\lambda_* \psi_* \pi_*(f_i)) / (\lambda_* \psi_* \pi_*(g))$-coalescent. Therefore the condensed the the condensed $\tu{A}$-algebra $(R_{\approx T''})_g$ is $f_i/g$-coalescent. Since this holds for every $1 \leq i \leq n$, we conclude that the condensed $\tu{A}$-algebra $(R_{\approx T''})_g$ is $T'$-coalescent.

Next suppose that $R'$ is a $T'$-coalescent condensed $\tu{A}$-algebra and that $\xi' : \tu{A} \to R'$ is a homomorphism of condensed $\tu{A}$-algebras such that the element $\xi'_*(g) \in R'(*)$ is invertible in $R'(*)$. By \cref{prop:universality of monoid algebras over condensed rings}, there exists a unique homomorphism of condensed $\tu{A}$-algebras $\eta_1 : \tu{A}[X_1, \ldots, X_n] \to R'$ such that $(\eta_1)_*(X_i) = \xi'_*(f_i) \cdot \xi'_*(g)^{-1}$ for each $1 \leq i \leq n$ and the diagram
\[\begin{tikzcd}
	{\tu{A}} & {R'} \\
	{\tu{A}[X_1 \ldots, X_n]}
	\arrow["{\xi'}", from=1-1, to=1-2]
	\arrow["\iota"', from=1-1, to=2-1]
	\arrow["{\eta_1}"', from=2-1, to=1-2]
\end{tikzcd}\]
is commuative. Let us consider $R'$ as a condensed $\tu{A}[X_1, \ldots, X_n]$-algebra via $\eta_1 : \tu{A}[X_1, \ldots, X_n] \to R'$. We have $(\eta_1)_*(g \cdot X_i) = \xi'_*(g) \cdot (\eta_1)_*(X_i) = \xi'_*(f_i) = (\eta_1)_*(f_i)$ for every $1 \leq i \leq n$. It follows that the composition ${\tu{A}[X_1, \ldots, X_n]}^n \xto{\mu} \tu{A}[X_1, \ldots, X_n] \xto{\eta_1} R'$ is equal to the zero homomorphism. By the definition of $R$, there exists a unique homormophism $\eta_2 : R \to R'$ of condensed $\tu{A}[X_1, \ldots, X_n]$-algebras such that the diagram
\[\begin{tikzcd}
	{\tu{A}[X_1 \ldots, X_n]} & {R'} \\
	R
	\arrow["{\eta_1}", from=1-1, to=1-2]
	\arrow["\pi"', from=1-1, to=2-1]
	\arrow["{\eta_2}"', from=2-1, to=1-2]
\end{tikzcd}\]
is commutative. Next we show that $R'$ is a $T''$-coalescent condensed $\tu{A}[X_1, \ldots, X_n]$-algebra via $\eta_1 : \tu{A}[X_1, \ldots, X_n] \to R'$. Indeed, $R'$ is a $T'$-coalescent condensed $\tu{A}$-algebra via $\xi' : \tu{A} \to R'$ by assumption. Therefore for each $1 \leq i \leq n$, the condensed $R'$-algebra $R'$ is 
$\xi'_*(f_i) / \xi'_*(g)$-coalescent. Since the element $\xi'_*(g) \in R'(*)$ is invertible in $R'(*)$ by assumption, it follows that the condensed $R'$-algebra $R'$ is 
$(\xi'_*(f_i) \cdot \xi'_*(g)^{-1}) /1$-coalescent. Since $(\eta_1)_*(X_i) = \xi'_*(f_i) \cdot \xi'_*(g)^{-1}$, we conclude that $R'$ is an $X_i/1$-coalescent condensed $\tu{A}[X_1, \ldots, X_n]$-algebra via $\eta_1 : \tu{A}[X_1, \ldots, X_n] \to R'$. Since this holds for every $1 \leq i \leq n$, it follows that $R'$ is a $T''$-coalescent condensed $\tu{A}[X_1, \ldots, X_n]$-algebra via $\eta_1 : \tu{A}[X_1, \ldots, X_n] \to R'$. On the other hand, the element $(\eta_1)_*(g) = \xi'_*(g) \in R'(*)$ is invertible in $R'(*)$ by assumption. Then \cref{prop:localization and coalescence commute} shows that there exists a unique homomorphism $\eta : (R_{\approx T''})_g \to R'$ of condensed $\tu{A}[X_1, \ldots, X_n]$-algebras such that the diagram
\[\begin{tikzcd}
	R & {R'} \\
	{(R_{\approx T''})_g}
	\arrow["{\eta_2}", from=1-1, to=1-2]
	\arrow["{\lambda \of \psi}"', from=1-1, to=2-1]
	\arrow["\eta"', from=2-1, to=1-2]
\end{tikzcd}\]
is commutative. Then the diagram
\[\begin{tikzcd}
	{\tu{A}} && {R'} \\
	{\tu{A}[X_1 \ldots, X_n]} \\
	R \\
	{(R_{\approx T''})_g}
	\arrow["{\xi'}", from=1-1, to=1-3]
	\arrow["\iota", from=1-1, to=2-1]
	\arrow["\xi"', shift right=4, curve={height=30pt}, from=1-1, to=4-1]
	\arrow["{\eta_1}", curve={height=6pt}, from=2-1, to=1-3]
	\arrow["\pi", from=2-1, to=3-1]
	\arrow["{\eta_2}"', curve={height=12pt}, from=3-1, to=1-3]
	\arrow["{\lambda \of \psi}", from=3-1, to=4-1]
	\arrow["\eta"', curve={height=30pt}, from=4-1, to=1-3]
\end{tikzcd}\]
is commutative, and $\eta : (R_{\approx T''})_g \to R'$ is a homomorphism of condensed $\tu{A}$-algebras.

On the other hand, suppose that $\eta' : (R_{\approx T''})_g \to R'$ is another homomorphism of condensed $\tu{A}$-algebras such that the diagram
\[\begin{tikzcd}
	{\tu{A}} & {R'} \\
	{(R_{\approx T''})_g}
	\arrow["{\xi'}", from=1-1, to=1-2]
	\arrow["\xi"', from=1-1, to=2-1]
	\arrow["{\eta'}"', from=2-1, to=1-2]
\end{tikzcd}\]
is commutative. We prove $\eta' = \eta$. The diagram
\[\begin{tikzcd}
	{\tu{A}} &&& {R'} \\
	{\tu{A}[X_1 \ldots, X_n]} & R && {(R_{\approx T''})_g}
	\arrow["{\xi'}", from=1-1, to=1-4]
	\arrow["\iota"', from=1-1, to=2-1]
	\arrow["\pi"', from=2-1, to=2-2]
	\arrow["{\lambda \of \psi}"', from=2-2, to=2-4]
	\arrow["{\eta'}"', from=2-4, to=1-4]
\end{tikzcd}\]
is commuative, and $\eta' \of \lambda \of \psi \of \pi : \tu{A}[X_1 \ldots, X_n] \to R'$ is a homomorphism of condensed $\tu{A}$-algebras. Moreover, for each $1 \leq i \leq n$, the definition of $R$ shows that $\pi_*(g \cdot X_i) = \pi_*(f_i)$ in $R(*)$. Therefore $\xi'_*(g) \cdot \big( \eta'_* \lambda_* \psi_* \pi_*(X_i) \big) = \eta'_* \lambda_* \psi_* \pi_*(g \cdot X_i) = \eta'_* \lambda_* \psi_* \pi_*(f_i) = \xi'_*(f_i)$. Since the element $\xi'_*(g) \in R'(*)$ is invertible in $R'(*)$ by assumption, it follows that $\eta'_* \lambda_* \psi_* \pi_*(X_i) = \xi'_*(f_i) \cdot \xi'_*(g)^{-1}$. Then the definition of $\eta_1$ shows that $\eta' \of \lambda \of \psi \of \pi = \eta_1$. In other words, the following diagram is commuative.
\[\begin{tikzcd}
	{\tu{A}[X_1 \ldots, X_n]} & {R'} \\
	R & {(R_{\approx T''})_g}
	\arrow["{\eta_1}", from=1-1, to=1-2]
	\arrow["\pi"', from=1-1, to=2-1]
	\arrow["{\lambda \of \psi}"', from=2-1, to=2-2]
	\arrow["{\eta'}"', from=2-2, to=1-2]
\end{tikzcd}\]
Then $\eta' \of \lambda \of \psi : R \to R'$ is a homomorphism of condensed $\tu{A}[X_1, \ldots, X_n]$-algebras. By the definition of $\eta_2$, we have $\eta' \of \lambda \of \psi = \eta_2$. In other words, the following diagram is commutative.
\[\begin{tikzcd}
	R & {R'} \\
	{(R_{\approx T''})_g}
	\arrow["{\eta_2}", from=1-1, to=1-2]
	\arrow["{\lambda \of \psi}"', from=1-1, to=2-1]
	\arrow["{\eta'}"', from=2-1, to=1-2]
\end{tikzcd}\]
Then $\eta' : (R_{\approx T''})_g \to R'$ is a homomorphism of condensed $\tu{A}[X_1, \ldots, X_n]$-algebras. By the definition of $\eta$, we conclude that $\eta' = \eta$. This completes the proof.
\end{proof}

\begin{prop} \label{prop:expression of rho in terms of sigma and tau} \;
\begin{enumerate}
\item
The following diagram is a commutative diagram in $\ub{CRing}$.
\[\begin{tikzcd}
	{\tu{A}} &&& {(\tu{A}_g)_{\approx T'}} \\
	& {\tu{A \langle X_1 , \dots , X_n \rangle}} & {R_{\approx T''}} & {(R_{\approx T''})_g} \\
	& {\tu{B}} & {\tu{B}_g} \\
	{\tu{A \left \langle \frac{T}{g} \right \rangle}} & {\tu{B \left \langle \frac{r(T)}{r(g)} \right \rangle}}
	\arrow["\can", from=1-1, to=1-4]
	\arrow["{\tu{j}}", from=1-1, to=2-2]
	\arrow["{\tu{l}}"', from=1-1, to=4-1]
	\arrow["\omega", from=1-4, to=2-4]
	\arrow["\chi", from=2-2, to=2-3]
	\arrow["{\tu{p}}"', from=2-2, to=3-2]
	\arrow["\can", from=2-3, to=2-4]
	\arrow["\sigma", from=2-3, to=3-2]
	\arrow["{\sigma'}", from=2-4, to=3-3]
	\arrow["\can"', from=3-2, to=3-3]
	\arrow["{\tu{k}}"', from=3-2, to=4-2]
	\arrow["\tau", from=3-3, to=4-2]
	\arrow["{\tu{u}}", from=4-2, to=4-1]
\end{tikzcd}\]

\item
The equality $\rho = \tu{u} \of \tau \of \sigma' \of \omega$ holds.

\item
$\rho$ is an epimorphism when considered as a morphism in $\ub{CSet}$.

\item
If $\sigma$ and $\tau$ are isomorphisms of condensed rings, then $\rho$ is also an isomorphism of condensed rings.

\end{enumerate}
\end{prop}

\begin{proof}~
\begin{enumerate}
\item
This is an immediate consequence of \cref{prop:construction of the homomorphism chi,prop:construction of the homomorphism sigma,prop:construction of the homomorphism sigma prime,prop:construction of the homomorphism tau,lem:construction of the homomorphism u,lem:construction of the homomorphism omega}.

\item
This follows from (1) and the uniqueness of $\rho$.

\item
By \cref{prop:construction of the homomorphism sigma prime,prop:construction of the homomorphism tau}, the maps $\tau : \tu{B}_g \to \tu{B \left \langle \frac{r(T)}{r(g)} \right \rangle}$ and $\sigma' : (R_{\approx T''})_g \to \tu{B}_g$ are epimorphisms when considered as morphisms in $\ub{CSet}$. By \cref{lem:construction of the homomorphism u,lem:construction of the homomorphism omega}, the maps $\tu{u} : \tu{B \left \langle \frac{r(T)}{r(g)} \right \rangle} \xto{\sim} \tu{A \left \langle \frac{T}{g} \right \rangle}$ and $\omega : (\tu{A}_g)_{\approx T'} \xto{\sim} (R_{\approx T''})_g$ are isomorphisms of condensed sets. Thus the assertion (3) follows from (2).

\item
If $\sigma : R_{\approx T''} \to \tu{B}$ is an isomorphism of condensed rings, then \cref{prop:construction of the homomorphism sigma prime} shows that the homomorphism $\sigma' : (R_{\approx T''})_g \to \tu{B}_g$ is also an isomorphism of condensed rings. By \cref{lem:construction of the homomorphism u,lem:construction of the homomorphism omega}, the homomorphisms $\tu{u} : \tu{B \left \langle \frac{r(T)}{r(g)} \right \rangle} \xto{\sim} \tu{A \left \langle \frac{T}{g} \right \rangle}$ and $\omega : (\tu{A}_g)_{\approx T'} \xto{\sim} (R_{\approx T''})_g$ are isomorphisms of condensed rings. Thus the assertion (4) follows from (2).
\end{enumerate}
\end{proof}

\subsubsection{The case where $A$ is a strongly Noetherian Tate ring}

\begin{lem} \label{lem:m is strict in stronlgy Noetherian Tate case}
Suppose that $A$ is a strongly Noetherian Tate ring. Then $m : {A \langle X_1 , \ldots , X_n \rangle}^{n} \to A \langle X_1 , \ldots , X_n \rangle$ is strict.
\end{lem}

\begin{proof}
If $A$ is a strongly Noetherian Tate ring, then the ring $A \langle X_1 , \ldots , X_n \rangle$ is a complete Hausdorff Noetherian Tate ring. Then the assertion follows form Proposition II.4.2.2 of \cite{Morel:note}.
\end{proof}

\begin{lem} \label{lem:B is already the localization in in stronlgy Noetherian Tate case}
Suppose that $A$ is a Tate ring.
\begin{enumerate}
\item $r(g) \in B$ is invertible in $B$.

\item The canonical homomorphism $\tu{B} \to \tu{B}_g$ is an isomorphism of condensed $\tu{A}$-algebras.

\item $k : B \to B \left \langle \frac{r(T)}{r(g)} \right \rangle$ is an isomorphism of topological rings.
\end{enumerate}
\end{lem}

\begin{proof}~
\begin{enumerate}
\item
Since $A$ has a topologically nilpotent unit, any open ideal of $A$ must be equal to $A$. Since $T$ generates an open ideal of $A$, we conclude that $T$ generates the unit ideal $A$. Therefore there exist elements $a_1, \ldots, a_n \in A$ such that $1 = \sum_{i = 1}^{n} a_i f_i$. Since $p(f_i - g X_i) = 0$ in $B$ for each $1 \leq i \leq n$ by definition, we have
\begin{equation}
r(g) \cdot \left( \sum_{i = 1}^{n} p(a_i X_i) \right) =
\sum_{i = 1}^{n} p(a_i) \cdot p(g X_i) = 
\sum_{i = 1}^{n} p(a_i) \cdot p(f_i) = 
p \left( \sum_{i = 1}^{n} a_i f_i \right) = p(1) = 1
\end{equation}
in $B$. Thus $r(g)$ is invertible in $B$.

\item
This immediately follows from (1) and the definition of $\tu{B}_g$. 

\item
Since $r(g)$ is invertible in $B$ by (1), the map $k_1 : B \to B \left( \frac{r(T)}{r(g)} \right)$ is an isomorphism of commutative unital rings. Furthermore, the map $k_1 : B \to B \left( \frac{r(T)}{r(g)} \right)$ is continuous by definition, and is an open map by (1) of \cref{lem:properties of localizations of B}. It follows that $k_1 : B \to B \left( \frac{r(T)}{r(g)} \right)$ is an isomorphism of topological rings. Then $B \left( \frac{r(T)}{r(g)} \right)$ is Hausdorff since $B$ is Hausdorff. Then (4) of \cref{lem:properties of localizations of B} shows that $k_2 : B \left( \frac{r(T)}{r(g)} \right) \to B \left \langle \frac{r(T)}{r(g)} \right \rangle$ is an isomorphism of topological rings. Consequently, $k = k_2 \of k_1 : B \to B \left \langle \frac{r(T)}{r(g)} \right \rangle$ is an isomorphism of topological rings.
\end{enumerate}
\end{proof}

\begin{prop}
Suppose that $A$ is a strongly Noetherian Tate ring. Then $\rho : (\tu{A}_g)_{\approx T'} \to \tu{A \left \langle \frac{T}{g} \right \rangle}$ is an isomorphism of condensed rings.
\end{prop}

\begin{proof}
By \cref{prop:expression of rho in terms of sigma and tau}, it suffices to show that the homomorphisms $\sigma : R_{\approx T''} \to \tu{B}$ and $\tau : \tu{B}_g \to \tu{B \left \langle \frac{r(T)}{r(g)} \right \rangle}$ are isomorphisms of condensed rings. 

By \cref{lem:m is strict in stronlgy Noetherian Tate case}, the map $m : {A \langle X_1 , \ldots , X_n \rangle}^{n} \to A \langle X_1 , \ldots , X_n \rangle$ is strict. Then \cref{prop:construction of the homomorphism sigma} shows that the homomorphism $\sigma : R_{\approx T''} \to \tu{B}$ is an isomorphism of condensed rings. 

By \cref{prop:construction of the homomorphism tau}, the following diagram is commutative.
\[\begin{tikzcd}
	{\tu{B}} & {\tu{B}_g} \\
	{\tu{B \left \langle \frac{r(T)}{r(g)} \right \rangle}}
	\arrow["\can", from=1-1, to=1-2]
	\arrow["{\tu{k}}"', from=1-1, to=2-1]
	\arrow["\tau", from=1-2, to=2-1]
\end{tikzcd}\]
Moreover, \cref{lem:B is already the localization in in stronlgy Noetherian Tate case} shows that the canonical homomorphism $\tu{B} \to \tu{B}_g$ and the homomorphism $\tu{k} : \tu{B} \to \tu{B \left \langle \frac{r(T)}{r(g)} \right \rangle}$ are isomorphisms of condensed rings. It follows that $\tau : \tu{B}_g \to \tu{B \left \langle \frac{r(T)}{r(g)} \right \rangle}$ is an isomorphism of condensed rings. 
\end{proof}

\subsubsection{The case where $A$ has a Noetherian ring of definition} $\\$

Note that $A \langle X_1 ,\ldots, X_n \rangle$ and $B$ are complete Hausdorff Huber rings.

\begin{lem} \label{lem:Noetherian rings of definition are inherited to A X and B}
Suppose that $A$ has a Noetherian ring of definition. Then both $A \langle X_1 ,\ldots, X_n \rangle$ and $B$ have Noetherian rings of definition.
\end{lem}

\begin{proof}
Let $(A_0,I_0)$ be a couple of definition of $A$ such that $A_0$ is a Noetherian ring. Then $(A_0 \langle X_1 ,\ldots, X_n \rangle, I_0 \cdot A_0 \langle X_1 ,\ldots, X_n \rangle)$ is a couple of definition of $A \langle X_1 ,\ldots, X_n \rangle$. On the other hand, the polynomial ring $A_0[X_1 ,\ldots, X_n]$ is a dense subring of $A_0 \langle X_1 ,\ldots, X_n \rangle$, and $A_0 \langle X_1 ,\ldots, X_n \rangle$ can be seen as the Hausdorff completion of the polynomial ring $A_0[X_1 ,\ldots, X_n]$ endowed with the subspace topology $\cat{T}$ induced by the topology of $A \langle X_1 ,\ldots, X_n \rangle$ (\cite{Morel:note}, Proposition II.3.3.3 ; \cite{Wedhorn:note}, Proposition 5.49). Moreover, the definition shows that the topology $\cat{T}$ coincides with the $I_0 \cdot A_0[X_1 ,\ldots, X_n]$-adic topology on $A_0[X_1 ,\ldots, X_n]$. Consequently, $A_0 \langle X_1 ,\ldots, X_n \rangle$ is equal to the Hausdorff completion of the polynomial ring $A_0[X_1 ,\ldots, X_n]$, which is Noetherian since $A_0$ is Noetherian, with respect to the $I_0 \cdot A_0[X_1 ,\ldots, X_n]$-adic topology. Therefore $A_0 \langle X_1 ,\ldots, X_n \rangle$ is Noetherian (\cite{Bourbaki:comalg}, Chapter III, \S 3.4, Proposition 8). This shows that $A \langle X_1 ,\ldots, X_n \rangle$ has a Noetherian ring of definition. Moreover, since the map $p : A \langle X_1 ,\ldots, X_n \rangle \to B$ is surjective and open, the subring $p(A_0 \langle X_1 ,\ldots, X_n \rangle)$ of $B$ is a ring of definition of $B$. Since $A_0 \langle X_1 ,\ldots, X_n \rangle$ is Noetherian, its quotient $p(A_0 \langle X_1 ,\ldots, X_n \rangle)$ is also Noetherian. This completes the proof.
\end{proof}

\begin{lem} \label{lem:complete Hausdorff Huber rings with Noetherian rings of definition}
Let $C$ be a complete Hausdorff Huber ring which has a Noetherian ring of definition.
\begin{enumerate}
\item Every ideal of $C$ is closed in $C$.

\item Let $\nu \in \N$. Every continuous $C$-linear map $C^{\nu} \to C$ is strict.
\end{enumerate}
\end{lem}

\begin{proof}
Let $(C_0,I_0)$ be a couple of definition of $C$ such that $C_0$ is a Noetherian ring. 
\begin{enumerate}
\item
First we claim that $C_0$ is a Zariski ring in the sense of \cite{Bourbaki:comalg}, Chapter III, \S 3.3, Definition 2. Indeed, since $C_0$ is an open subring of $C$, the subring $C_0$ is also closed in $C$. Since $C$ is complete Hausdorff, it follows that $C_0$ is also complete Hausdorff. Moreover, $C_0$ is Noetherian and has the $I_0$-adic topology. Therefore the topological ring $C_0$ is a Zariski ring (\cite{Bourbaki:comalg}, Chapter III, \S 3.4, Proposition 8).

Let $I$ be any ideal of $C$. Then $C_0 \cap I$ is an ideal of $C_0$. Since $C_0$ is a Zariski ring, the ideal $C_0 \cap I$ is a closed ideal of $C_0$. Since $C_0$ is complete Hausdorff, it follows that $C_0 \cap I$ is also complete Hausdorff. On the other hand, since $C_0$ is open in $C$, the intersection $C_0 \cap I$ is an open subgroup of $I$. Then by \cite{Bourbaki:top1}, Chapter III, \S 3.3, Proposition 4, $I$ is complete. Since $C$ is complete Hausdorff, it follows that $I$ is closed in $C$.

\item
Let $\nu \in \N$. Let $c : C^{\nu} \to C$ be a continuous $C$-linear map. Let $K$ be the image of $c : C^{\nu} \to C$. For each $e \in \N$, let us write $N(e)$ for the direct sum of $\nu$ copies of $I_0^e$. Then the set $\set{N(e)}{e \in \N}$ is a fundamental system of neighbourhoods of $0$ in $C^{\nu}$. We prove that $c(N(e))$ is a neighbourhood of $0$ in $K$ for every $e \in \N$. Let us fix $e \in \N$ arbitrarily. First of all, we note that the intersection $K_0 := C_0 \cap K$ is an ideal of $C_0$. Since $C_0$ is a Noetherian ring, this ideal $K_0$ is finitely generated. Let $x_1, \ldots, x_t \in K_0$ be a finite system of generators of $K_0$. Since the elements $x_1, \ldots, x_t$ belong to the image $K$ of the map $c : C^{\nu} \to C$, there exist elements $y_1, \ldots, y_t \in C^{\nu}$ such that $c(y_i) = x_i$ for every $1 \leq i \leq t$. Then, since the multiplication map of $C$ is continuous, there exists a $d \in \N$ such that $I_0^d \cdot y_i \sub N(e)$ for every $1 \leq i \leq t$. Then we have
\begin{equation}
c(N(e)) \bus c \Big( I_0^d \cdot y_1 + \cdots + I_0^d \cdot y_t \Big)
= I_0^d \cdot x_1 + \cdots + I_0^d \cdot x_t
= I_0^d \cdot K_0 .
\end{equation}
On the other hand, since $C_0$ is Noetherian and has the $I_0$-adic topology, the subspace topology on $K_0$ induced by the topology of $C_0$ is equal to the $I_0$-adic topology on the $C_0$-module $K_0$ (\cite{Bourbaki:comalg}, Chapter III, \S 3.2, Theorem 2). Therefore $I_0^d \cdot K_0$ is a neighbourhood of $0$ in $K_0$. Furthermore, since $C_0$ is open in $C$, the intersection $K_0 = C_0 \cap K$ is open in $K$. Therefore $I_0^d \cdot K_0$ is also a neighbourhood of $0$ in $K$. Since $c(N(e))$ contains $I_0^d \cdot K_0$, we conclude that $c(N(e))$ is also a neighbourhood of $0$ in $K$. This completes the proof.
\end{enumerate}
\end{proof}

\begin{prop}
Suppose that $A$ has a Noetherian ring of definition. Then $\rho : (\tu{A}_g)_{\approx T'} \to \tu{A \left \langle \frac{T}{g} \right \rangle}$ is an isomorphism of condensed rings.
\end{prop}

\begin{proof}
By \cref{prop:expression of rho in terms of sigma and tau}, it suffices to show that the homomorphisms $\sigma : R_{\approx T''} \to \tu{B}$ and $\tau : \tu{B}_g \to \tu{B \left \langle \frac{r(T)}{r(g)} \right \rangle}$ are isomorphisms of condensed rings. 

By \cref{lem:Noetherian rings of definition are inherited to A X and B}, the Huber rings $A \langle X_1 ,\ldots, X_n \rangle$ and $B$ have Noetherian rings of definition. Then \cref{lem:complete Hausdorff Huber rings with Noetherian rings of definition} shows that the map $m : {A \langle X_1 , \ldots , X_n \rangle}^{n} \to A \langle X_1 , \ldots , X_n \rangle$ is strict and that the kernel of $k_1 : B \to B \left( \frac{r(T)}{r(g)} \right)$ is a closed ideal of $B$. By \cref{prop:construction of the homomorphism sigma,prop:construction of the homomorphism tau}, we conclude that the homomorphisms $\sigma : R_{\approx T''} \to \tu{B}$ and $\tau : \tu{B}_g \to \tu{B \left \langle \frac{r(T)}{r(g)} \right \rangle}$ are isomorphisms of condensed rings.
\end{proof}

\subsection{Relation to adic spectra}

\subsubsection{Settings}

$\\[2mm]$
Throughout this subsection, we fix the following notation.

\begin{nt} \,
\begin{enumerate}
\item $(A,A^+)$ denotes a Huber pair (\cite{Morel:note}, Definition III.1.7), or an affinoid ring (\cite{Wedhorn:note}, Definition 7.14), such that $A$ is complete Hausdorff.

\item $X = \mathrm{Spa}(A,A^+)$ is the adic spectrum of $(A,A^+)$ (\cite{Morel:note}, Definition III.2.1 ; \cite{Wedhorn:note}, Definition 7.23).

\item For each element $g \in A$ of $A$ and each finite subset $T \sub A$ of $A$ which generates an open ideal of $A$, we write
\begin{equation}
R\left( \frac{T}{g} \right) :=
\set{x \in X}{|f|_x \leq |g|_x \neq 0 \text{ for every } f \in T} .
\end{equation}
The subsets of $X$ of this form are called \ti{rational subsets} of $X$ (\cite{Morel:note}, III.2 ; \cite{Wedhorn:note}, Definition 7.29). The rational subsets of $X$ form an open basis for the topology of $X$ (\cite{Morel:note}, Corollary III.2.4 ; \cite{Wedhorn:note}, Theorem 7.35, (2)).

\item $\cat{O}$ denotes the structure presheaf on the adic spectrum $X = \mathrm{Spa}(A,A^+)$ of $(A,A^+)$ (\cite{Morel:note}, Definition III.6.2.1 ; \cite{Wedhorn:note}, 8.1).
\end{enumerate}
\end{nt}

\begin{nt} \,
\begin{enumerate}
\item
For each open subset $U$ of $X$, define
\begin{equation}
\tu{\cat{O}}^{\mathrm{pre}}(U) := \tu{\cat{O}(U)}.
\end{equation}
We obtain a presheaf of condensed rings $\tu{\cat{O}}^{\mathrm{pre}}$ on the topological space $X$.

\item
The sheafification of $\tu{\cat{O}}^{\mathrm{pre}}$ is denoted by $\eta : \tu{\cat{O}}^{\mathrm{pre}} \to \tu{\cat{O}}$.
\end{enumerate}
\end{nt}

\begin{nt} \,
\begin{enumerate}
\item
$(\{*\}, F, \cat{U})$ is an object of $\cat{C}$ defined as follows.
\begin{enumerate}
\item
$\{*\}$ denotes the topological space with a unique point $*$, equipped with the unique topology.

\item
$F$ denotes the unique sheaf of condensed rings on $\{*\}$ such that $F(\{*\}) = \tu{A}$ (cf. \cref{prop:taking global section over * is an equivalence}).

\item
$\cat{U}$ denotes the family $(\cat{U}_*)_{* \in \{*\}}$ consisting of a single set $\cat{U}_* = \mathrm{Spa}(A,A^+)$. By \cref{rem:comparison of continuous valuations on R and on underbar R}, this set is a set of continuous valuations on the stalk $F_* = \tu{A}$ of $F$ at the point $* \in \{*\}$.
\end{enumerate}

\item 
$(\tilde{X}, \tilde{\cat{O}}, \tilde{\cat{V}})$ denotes the image of the object $(\{*\}, F, \cat{U}) \in |\cat{C}|$ under the right adjoint of the inclusion $\cat{C}_l \cap \cat{C}_c \mon \cat{C}$ (\cref{prop:intersections of C l C f C c are coreflective in C}). The canonical morphism $(\tilde{X}, \tilde{\cat{O}}, \tilde{\cat{V}}) \to (\{*\}, F, \cat{U})$ is denoted by $(\pi, \pi^{\#})$.
\end{enumerate}
\end{nt}

\subsubsection{Adic spectra as objects of $\cat{C}$}

\begin{prop} \label{prop:valuation on the stalk of underbar O}
Let $x \in X$. Let $c : \tu{A} \to \tu{\cat{O}}_x$ be the composition of the homomorphism $\eta_X : \tu{A} = \tu{\cat{O}}^{\mathrm{pre}}(X) \to \tu{\cat{O}}(X)$ and the canonical homomorphism $\tu{\cat{O}}(X) \to \tu{\cat{O}}_x$. Then there exists a unique continuous valuation $v_x$ on $\tu{\cat{O}}_x$ such that $c^{-1}(v_x) = x$.
\end{prop}

\begin{proof}
By \cref{cor:existence of sheafification functor of condensed presheaves}, the homomorphism $\eta_x : (\tu{\cat{O}}^{\mathrm{pre}})_x \to \tu{\cat{O}}_x$ is an isomorphism of condensed rings. Moreover, the diagram
\[\begin{tikzcd}
	{\tu{\cat{O}}^{\mathrm{pre}}(X)} & {\tu{\cat{O}}(X)} \\
	{(\tu{\cat{O}}^{\mathrm{pre}})_x} & {\tu{\cat{O}}_x}
	\arrow["{\eta_X}", from=1-1, to=1-2]
	\arrow["\can"', from=1-1, to=2-1]
	\arrow["c", from=1-1, to=2-2]
	\arrow["\can", from=1-2, to=2-2]
	\arrow["{\eta_x}"', from=2-1, to=2-2]
\end{tikzcd}\]
is commutative. Therefore it suffices to show that there exists a unique continuous valuation $v_x$ on the condensed ring $(\tu{\cat{O}}^{\mathrm{pre}})_x$ such that $\rho_X^{-1}(v_x) = x$, where we write $\rho_U : \tu{\cat{O}}^{\mathrm{pre}}(U) \to (\tu{\cat{O}}^{\mathrm{pre}})_x$ for the canonical homomorphism for each open neighbourhood $U$ of $x$ in $X$.

Let $U$ be a rational subset of $X$ such that $x \in U$. Let $l_U : A \to \cat{O}(U)$ be the canonical map. Corollary III.4.3.3 of \cite{Morel:note} shows that there exists a unique continuous valuation $v_U$ on the topological ring $\cat{O}(U)$ such that $v_U \of l_U = x$. Moreover, the value group of $v_U$ coincides with the value group of the valuation $x$. On the other hand, since the topological ring $\cat{O}(U)$ is a Huber ring, it is first countable. Then \cref{rem:comparison of continuous valuations on R and on underbar R} shows that the notion of continuous valuations on the topological ring $\cat{O}(U)$ coincides with the notion of continuous valuations on the condensed ring $\tu{\cat{O}(U)}$. Consequently, $v_U$ is a unique continuous valuation on the condensed ring $\tu{\cat{O}(U)}$ such that $(\tu{l_U})^{-1}(v_U) = x$. Moreover, the value group of $v_U$ coincides with the value group of the valuation $x$. It follows that if $U,V$ are a rational subset of $X$ such that $x \in V \sub U$, then the restriction map $r_{U,V} : \cat{O}(U) \to \cat{O}(V)$ satisfies $(\tu{r_{U,V}})^{-1}(v_V) = v_U$ and induces an isomorphism of the value group of $v_U$ onto that of $v_V$. On the other hand, since rational subsets of $X$ form a basis for the topology of $X$, we have
\begin{equation}
(\tu{\cat{O}}^{\mathrm{pre}})_x
= \underset{x \in U \sub X \text{ rational}}{\colim} \: \tu{\cat{O}}^{\mathrm{pre}}(U)
= \underset{x \in U \sub X \text{ rational}}{\colim} \: \tu{\cat{O}(U)} .
\end{equation}
Then \cref{prop:extension of continuous valuation to filtered colimits} shows that there exists a unique continuous valuation $v_x$ on the condensed ring $(\tu{\cat{O}}^{\mathrm{pre}})_x$ such that $\rho_U^{-1}(v_x) = v_U$ for every rational subset $U$ of $X$ with $x \in U$. In particular, we have $\rho_X^{-1}(v_x) = v_X = x$.

On the other hand, suppose that $v'_x$ is another continuous valuation on the condensed ring $(\tu{\cat{O}}^{\mathrm{pre}})_x$ such that $\rho_X^{-1}(v'_x) = x$. Then, for every rational subset $U$ of $X$ with $x \in U$, we have $(\tu{l_U})^{-1} \big( \rho_U^{-1}(v'_x) \big) = \rho_X^{-1}(v'_x) = x$. By the uniqueness of $v_U$, we conclude that $\rho_U^{-1}(v'_x) = v_U$. Then the uniqueness of $v_x$ implies that $v'_x = v_x$. This completes the proof.
\end{proof}

\begin{prop} \label{prop:adic spectra is an object of C l cap C c}
For each $x \in X$, define $\tu{\cat{V}}_x := \{v_x\}$. Write $\tu{\cat{V}} := (\tu{\cat{V}}_x)_{x \in X}$. Then $(X, \tu{\cat{O}}, \tu{\cat{V}})$ is an object of $\cat{C}_l \cap \cat{C}_c$.
\end{prop}

\begin{proof}
Let $x \in X$. The homomorphism $\eta_x : (\tu{\cat{O}}^{\mathrm{pre}})_x \to \tu{\cat{O}}_x$ is an isomorphism of condensed rings by \cref{cor:existence of sheafification functor of condensed presheaves}. Furthermore, \cref{prop:compatibility of stalk and evaluation at S} shows that 
\begin{align}
(\tu{\cat{O}}^{\mathrm{pre}})_x (*)
& = \underset{x \in U \sub X \text{ open}}{\colim} \:
\left( \tu{\cat{O}}^{\mathrm{pre}}(U)(*) \right) \\
& = \underset{x \in U \sub X \text{ open}}{\colim} \:
\left( \tu{\cat{O}(U)}(*) \right) \\
& = \underset{x \in U \sub X \text{ open}}{\colim} \: \cat{O}(U) = \cat{O}_x .
\end{align}
Therefore we have an isomorphism $(\eta_x)_* : \cat{O}_x = (\tu{\cat{O}}^{\mathrm{pre}})_x (*) \xto{\sim} \tu{\cat{O}}_x (*)$ of commutative unital rings. From the construction of $v_x$ given in the proof of \cref{prop:valuation on the stalk of underbar O}, one concludes that the valuation $v_x \of (\eta_x)_*$ on $\cat{O}_x$ is equal to the valuation $|\cdot|_x$ on $\cat{O}_x$ defined in III.6.3 of \cite{Morel:note}. Then (i) of Proposition III.6.3.1 of \cite{Morel:note} shows that $\tu{\cat{O}}_x (*)$ is a local ring whose maximal ideal is equal to the support of $v_x$. Therefore $(X, \tu{\cat{O}}, \tu{\cat{V}})$ is an object of $\cat{C}_l$.

Next we show that $(X, \tu{\cat{O}}, \tu{\cat{V}})$ is an object of $\cat{C}_f$. Let $U$ be an open subset of $X$. Let $f,g \in \tu{\cat{O}}(U)$. We prove that the set
\begin{equation}
V := \set{x \in U}{|f_x|_{v_x} \leq |g_x|_{v_x} \neq 0}
\end{equation}
is an open subset of $X$. For this it suffices to show that for every $x \in U$, there exists an open neighbourhood $W$ of $x$ in $X$ such that $V \cap W$ is open in $X$. Let us fix $x \in U$ arbitrarily. Since $(\eta_x)_* : \cat{O}_x = (\tu{\cat{O}}^{\mathrm{pre}})_x (*) \xto{\sim} \tu{\cat{O}}_x (*)$ is an isomorphism of commutative unital rings, there exist an open neighbourhood $U'$ of $x$ in $X$ and elements $f',g' \in \cat{O}(U)$ such that $(\eta_x)_*(f'_x) = f_x$ and $(\eta_x)_*(g'_x) = g_x$. Then the commutativity of the diagram
\[\begin{tikzcd}
	{\tu{\cat{O}}^{\mathrm{pre}}(U')} & {\tu{\cat{O}}(U')} \\
	{(\tu{\cat{O}}^{\mathrm{pre}})_x} & {\tu{\cat{O}}_x}
	\arrow["{\eta_{U'}}", from=1-1, to=1-2]
	\arrow["\can"', from=1-1, to=2-1]
	\arrow["\can", from=1-2, to=2-2]
	\arrow["{\eta_x}"', from=2-1, to=2-2]
\end{tikzcd}\]
shows that
\begin{align}
\big( (\eta_{U'})_*(f') \big)_x = (\eta_x)_*(f'_x) = f_x \: ; \\
\big( (\eta_{U'})_*(g') \big)_x = (\eta_x)_*(g'_x) = g_x .
\end{align}
On the other hand, \cref{prop:compatibility of stalk and evaluation at S} shows that
\begin{equation}
\tu{\cat{O}}_x (*) = \underset{x \in W \sub X \text{ open}}{\colim} \; 
\tu{\cat{O}}(W)(*) ,
\end{equation}
where the colimit is taken in $\ub{Ring}$. Then the construction of filtered colimits in $\ub{Ring}$ shows that there exists an open neighbourhood $W$ of $x$ in $X$ such that $W \sub U \cap U'$ and
\begin{equation}
( (\eta_{U'})_*(f') )|_W = f|_W \quad \text{ and } \quad
( (\eta_{U'})_*(g') )|_W = g|_W .
\end{equation}
Then for every $y \in W$, the commutativity of the diagram
\[\begin{tikzcd}
	{\tu{\cat{O}}^{\mathrm{pre}}(U')} & {\tu{\cat{O}}(U')} \\
	{(\tu{\cat{O}}^{\mathrm{pre}})_y} & {\tu{\cat{O}}_y}
	\arrow["{\eta_{U'}}", from=1-1, to=1-2]
	\arrow["\can"', from=1-1, to=2-1]
	\arrow["\can", from=1-2, to=2-2]
	\arrow["{\eta_y}"', from=2-1, to=2-2]
\end{tikzcd}\]
and the equalities $( (\eta_{U'})_*(f') )|_W = f|_W$ and $( (\eta_{U'})_*(g') )|_W = g|_W$ show that
\begin{align}
(\eta_y)_*(f'_y) = \big( (\eta_{U'})_*(f') \big)_y = f_y \: ; \\
(\eta_y)_*(g'_y) = \big( (\eta_{U'})_*(g') \big)_y = g_y .
\end{align}
Furthermore, we already know that the valuation $v_y \of (\eta_y)_*$ on $\cat{O}_y$ is equal to the valuation $|\cdot|_y$ on $\cat{O}_y$ defined in III.6.3 of \cite{Morel:note}. Consequently, we have
\begin{equation}
V \cap W = \set{y \in W}{|f'_y|_y \leq |g'_y|_y \neq 0} .
\end{equation}
By Lemma III.6.3.2 of \cite{Morel:note}, this set is open in $X$. This completes the proof that $(X, \tu{\cat{O}}, \tu{\cat{V}})$ is an object of $\cat{C}_f$.

Next we show that $(X, \tu{\cat{O}}, \tu{\cat{V}})$ is an object of $\cat{C}_c$. Let $x \in X$. We prove that the valued condensed ring $(\tu{\cat{O}}_x, v_x)$ is an object of $\ub{VCRing}_c$. We already know that the homomorphism $\eta_x : ( (\tu{\cat{O}}^{\mathrm{pre}})_x , |\cdot|_x ) \xto{\sim} (\tu{\cat{O}}_x, v_x)$ is an isomorphism of valued condensed rings, where $|\cdot|_x$ denotes the valuation on $(\tu{\cat{O}}^{\mathrm{pre}})_x (*) = \cat{O}_x$ defined in III.6.3 of \cite{Morel:note}. Therefore it suffices to show that the valued condensed ring $( (\tu{\cat{O}}^{\mathrm{pre}})_x , |\cdot|_x )$ is an object of $\ub{VCRing}_c$. Let $s,t \in (\tu{\cat{O}}^{\mathrm{pre}})_x (*) = \cat{O}_x$ and suppose that $|s|_x \leq |t|_x \neq 0$. Since the ring $\cat{O}_x$ is a local ring whose maximal ideal is equal to the support of $|\cdot|_x$ (\cite{Morel:note}, Proposition III.6.3.1, (i)), the element $t \in \cat{O}_x$ is invertible in $\cat{O}_x$. Therefore there exist an open neighbourhood $U$ of $x$ in $X$ and elements $f,g \in \cat{O}(U)$ such that $f_x = s$, $g_x = t$ and $g \in \cat{O}(U)$ is invertible in $\cat{O}(U)$. By Lemma III.6.3.2 of \cite{Morel:note}, the set
\begin{equation}
V := \set{y \in U}{|f_y|_y \leq |g_y|_y \neq 0}
\end{equation}
is an open subset of $X$. Moreover, since $f_x = s$, $g_x = t$ and $|s|_x \leq |t|_x \neq 0$, we have $x \in V$. On the other hand, suppose that $W$ is any rational subset of $X$ such that $x \in W \sub V$. The element $g|_W \in \cat{O}(W)$ is invertible in $\cat{O}(W)$ since $g \in \cat{O}(U)$ is invertible in $\cat{O}(U)$. Moreover, for every $y \in W$ we have $|f_y|_y \leq |g_y|_y \neq 0$ and therefore
\begin{equation}
\big| \big( (f|_W) \cdot (g|_W)^{-1} \big)_y \big|_y = |f_y|_y \cdot (|g_y|_y)^{-1} \leq 1 .
\end{equation}
Consequently, we have $(f|_W) \cdot (g|_W)^{-1} \in \cat{O}(W)^{+}$, where $\cat{O}(W)^{+}$ is defined in Definition III.6.2.3 of \cite{Morel:note}. Since $(\cat{O}(W), \cat{O}(W)^{+})$ is a Huber pair by Lemma III.6.2.4 of \cite{Morel:note}, the element $(f|_W) \cdot (g|_W)^{-1} \in \cat{O}(W)$ is power-bounded in the topological ring $\cat{O}(W)$. Then \cref{prop:power-boundedness and coalescence} shows that $\tu{\cat{O}(W)}$ is an $( (f|_W) \cdot (g|_W)^{-1} ) / 1$-coalescent condensed $\tu{\cat{O}(W)}$-algebra. Since $g|_W \in \cat{O}(W)$ is invertible in $\cat{O}(W)$, it follows that $\tu{\cat{O}(W)}$ is an $(f|_W) / (g|_W)$-coalescent condensed $\tu{\cat{O}(W)}$-algebra. Therefore $\tu{\cat{O}(W)}$ is an $f/g$-coalescent condensed $\tu{\cat{O}(U)}$-algebra via the restriction $\tu{\cat{O}(U)} \to \tu{\cat{O}(W)}$. In other words, $\tu{\cat{O}}^{\mathrm{pre}}(W)$ is an $f/g$-coalescent condensed $\tu{\cat{O}}^{\mathrm{pre}}(U)$-algebra via the restriction $\tu{\cat{O}}^{\mathrm{pre}}(U) \to \tu{\cat{O}}^{\mathrm{pre}}(W)$. For each rational subset $W$ of $X$ such that $x \in W \sub V$, let us consider $\tu{\cat{O}}^{\mathrm{pre}}(W)$ as a condensed $\tu{\cat{O}}^{\mathrm{pre}}(U)$-algebra via the restriction $\tu{\cat{O}}^{\mathrm{pre}}(U) \to \tu{\cat{O}}^{\mathrm{pre}}(W)$. Moreover, let us consider $(\tu{\cat{O}}^{\mathrm{pre}})_x$ as a condensed $\tu{\cat{O}}^{\mathrm{pre}}(U)$-algebra via the canonical homomorphism $\tu{\cat{O}}^{\mathrm{pre}}(U) \to (\tu{\cat{O}}^{\mathrm{pre}})_x$. Then we can write
\begin{equation}
(\tu{\cat{O}}^{\mathrm{pre}})_x =
\underset{\substack{x \in W \sub V \\ W \sub X \text{ rational}}}{\colim} \:
\tu{\cat{O}}^{\mathrm{pre}}(W)
\end{equation}
in the category $\ub{CAlg}_{\tu{\cat{O}}^{\mathrm{pre}}(U)}$. By \cref{cor:coalescent modules are closed under limits and colimits}, we conclude that $(\tu{\cat{O}}^{\mathrm{pre}})_x$ is an $f/g$-coalescent condensed $\tu{\cat{O}}^{\mathrm{pre}}(U)$-algebra via the canonical homomorphism $\tu{\cat{O}}^{\mathrm{pre}}(U) \to (\tu{\cat{O}}^{\mathrm{pre}})_x$. Since $f_x = s$ and $g_x = t$, it follows that the condensed $(\tu{\cat{O}}^{\mathrm{pre}})_x$-algebra $(\tu{\cat{O}}^{\mathrm{pre}})_x$ is $s/t$-coalescent. This shows that the valued condensed ring $( (\tu{\cat{O}}^{\mathrm{pre}})_x , |\cdot|_x )$ is an object of $\ub{VCRing}_c$. This completes the proof.
\end{proof}

\begin{prop} \label{prop:sigma sharp inverse v x is equal to x}
Let $(\sigma, \sigma^{\#}) : (X, \tu{\cat{O}}) \to (\{*\}, F)$ be the morphism in $\cat{D}$ defined as follows. 
\begin{enumerate}
\item
The continuous map $\sigma : X \to \{*\}$ is the unique map $X \to \{*\}$.

\item
The morphism $\sigma^{\#} : F \to \sigma_* \tu{\cat{O}}$ of sheaves of condensed rings on $\{*\}$ is the unique one such that the homomorphism $\sigma^{\#}_{\{*\}}: F(\{*\}) \to \tu{\cat{O}}(X)$ of condensed rings is equal to the homomorphism $\eta_X : \tu{A} = \tu{\cat{O}}^{\mathrm{pre}}(X) \to \tu{\cat{O}}(X)$.
\end{enumerate}
For each $x \in X$, consider the homomorphism $\sigma^{\#}_x : \tu{A} = F_* \to \tu{\cat{O}}_x$ of condensed rings. Then the valuation $v_x$ on the condensed ring $\tu{\cat{O}}_x$ is a unique continuous valuation on $\tu{\cat{O}}_x$ such that $(\sigma^{\#}_x)^{-1}(v_x) = x$. In particular, we have
\begin{equation}
(\sigma^{\#}_x)^{-1}(v_x) = x \in \mathrm{Spa}(A,A^+) = \cat{U}_*
\end{equation}
for each $x \in X$, and therefore the pair $(\sigma, \sigma^{\#})$ is a morphism $(X, \tu{\cat{O}}, \tu{\cat{V}}) \to (\{*\}, F, \cat{U})$ in $\cat{C}$.
\end{prop}

\begin{proof}
For each $x \in X$, the homomorphism $\sigma^{\#}_x : F_* \to \tu{\cat{O}}_x$ of condensed rings is equal to the composition $c : \tu{A} \to \tu{\cat{O}}_x$ of the homomorphism $\eta_X : \tu{A} = \tu{\cat{O}}^{\mathrm{pre}}(X) \to \tu{\cat{O}}(X)$ and the canonical homomorphism $\tu{\cat{O}}(X) \to \tu{\cat{O}}_x$. Then the assertion follows from \cref{prop:valuation on the stalk of underbar O}.
\end{proof}

\subsubsection{Statement of the main result}

\begin{thm} \label{thm:comparison of coreflection and adic spectrum} \;
\begin{enumerate}
\item
There exists a unique morphism $(\tau, \tau^{\#}) : (X, \tu{\cat{O}}, \tu{\cat{V}}) \to (\tilde{X}, \tilde{\cat{O}}, \tilde{\cat{V}})$ in $\cat{C}$ such that the diagram
\[\begin{tikzcd}
	{(X, \tu{\cat{O}}, \tu{\cat{V}})} & {(\tilde{X}, \tilde{\cat{O}}, \tilde{\cat{V}})} \\
	& {(\{*\}, F, \cat{U})}
	\arrow["{(\tau, \tau^{\#})}", from=1-1, to=1-2]
	\arrow["{(\sigma, \sigma^{\#})}"', from=1-1, to=2-2]
	\arrow["{(\pi, \pi^{\#})}", from=1-2, to=2-2]
\end{tikzcd}\]
is commutative.

\item The continuous map $\tau : X \to \tilde{X}$ is a homeomorphism of topological spaces.

\item For each $x \in X$, the homomorphism $\tau^{\#}_x : \tilde{\cat{O}}_{\tau(x)} \to \tu{\cat{O}}_x$ of condensed rings is an epimorphism when considered as a morphism in $\ub{CSet}$.

\item
Suppose that one of the following holds.
\begin{enumerate}
\item $A$ is a strongly Noetherian Tate ring.
\item $A$ has a Noetherian ring of definition.
\end{enumerate}
Then for each $x \in X$, the homomorphism $\tau^{\#}_x : \tilde{\cat{O}}_{\tau(x)} \to \tu{\cat{O}}_x$ is an isomorphism of condensed rings. Consequently, the morphism $(\tau, \tau^{\#}) : (X, \tu{\cat{O}}, \tu{\cat{V}}) \to (\tilde{X}, \tilde{\cat{O}}, \tilde{\cat{V}})$ is an isomorphism in $\cat{C}$.
\end{enumerate}
\end{thm}

\subsubsection{Proof of \cref{thm:comparison of coreflection and adic spectrum}}

\begin{as} \label{as:explicit description of tilde X}
Since coreflections along functors are unique up to canonical isomorphisms, we may assume that the pair
\begin{equation}
\left( (\tilde{X}, \tilde{\cat{O}}, \tilde{\cat{V}}) \, , \,
(\tilde{X}, \tilde{\cat{O}}, \tilde{\cat{V}}) \xto{(\pi, \pi^{\#})} (\{*\}, F, \cat{U})
\right)
\end{equation}
has the properties described in \cref{prop:construction of coreflection in C l cap C c}. In other words, we may assume that the following hold.
\begin{enumerate}
\item
The underlying set of $\tilde{X}$ is the set of all pairs $(*,v)$ consisting of the point $* \in \{*\}$ and a valuation $v \in \cat{U}_* = \mathrm{Spa}(A,A^+)$.

\item
For every $(*,v) \in \tilde{X}$, let us write $\tilde{v}_{(*,v)}$ for the unique element of $\tilde{\cat{V}}_{(*,v)}$. Then the homomorphism $\pi^{\#}_{(*,v)} : \tu{A} = F_* \to \tilde{\cat{O}}_{(*,v)}$ of condensed rings satisfies $\left( \pi^{\#}_{(*,v)} \right)^{-1}\left( \tilde{v}_{(*,v)} \right) = v$. Moreover, the pair
\begin{equation}
\left( (\tilde{\cat{O}}_{(*,v)} , \tilde{v}_{(*,v)}) \; , \; (\tu{A} , v) \xto{\pi^{\#}_{(*,v)}} (\tilde{\cat{O}}_{(*,v)} , \tilde{v}_{(*,v)}) \right)
\end{equation}
is the reflection of the valued condensed ring $(\tu{A} , v)$ along the inclusion functor $\ub{VCRing}_l \cap \ub{VCRing}_c \mon \ub{VCRing}$.
\end{enumerate}
Throughout the proof of \cref{thm:comparison of coreflection and adic spectrum}, we assume that these are satisfied.
\end{as}

\begin{proof}[Proof of (1) of \cref{thm:comparison of coreflection and adic spectrum}]
By \cref{prop:adic spectra is an object of C l cap C c}, $(X, \tu{\cat{O}}, \tu{\cat{V}})$ is an object of $\cat{C}_l \cap \cat{C}_c$. Then the assertion (1) of \cref{thm:comparison of coreflection and adic spectrum} is an immediate consequence of the universal property of the coreflection 
\begin{equation}
\left( (\tilde{X}, \tilde{\cat{O}}, \tilde{\cat{V}}) \, , \,
(\tilde{X}, \tilde{\cat{O}}, \tilde{\cat{V}}) \xto{(\pi, \pi^{\#})} (\{*\}, F, \cat{U})
\right)
\end{equation}
of $(\{*\}, F, \cat{U}) \in |\cat{C}|$ along the inclusion functor $\cat{C}_l \cap \cat{C}_c \mon \cat{C}$.
\end{proof}

\begin{lem} \label{lem:identification of the map tau}
Recall from (1) of \cref{as:explicit description of tilde X} that the underlying set of $\tilde{X}$ is equal to the set of all pairs $(*,v)$ consisting of the point $* \in \{*\}$ and a valuation $v \in \cat{U}_* = \mathrm{Spa}(A,A^+) = X$. 
\begin{enumerate}
\item
The map $\tau : X \to \tilde{X}$ is of the form $x \mapsto (*,x)$.

\item
For each $x \in X$, the diagram
\[\begin{tikzcd}
	{(\tu{\cat{O}}_x,v_x)} & {(\tilde{\cat{O}}_{(*,x)},\tilde{v}_{(*,x)})} \\
	& {(\tu{A},x)}
	\arrow["{\tau^{\#}_x}"', from=1-2, to=1-1]
	\arrow["{\sigma^{\#}_x}", from=2-2, to=1-1]
	\arrow["{\pi^{\#}_{(*,x)}}"', from=2-2, to=1-2]
\end{tikzcd}\]
is commutative in $\ub{VCRing}$.
\end{enumerate}
\end{lem}

\begin{proof}
Let $x \in X$. Let us write $\tau(x) = (*,v)$, where $v \in \cat{U}_* = \mathrm{Spa}(A,A^+) = X$. Then the following diagram is commutative.
\[\begin{tikzcd}
	{\tu{\cat{O}}_x} & {\tilde{\cat{O}}_{(*,v)}} \\
	& {F_* = \tu{A}}
	\arrow["{\tau^{\#}_x}"', from=1-2, to=1-1]
	\arrow["{\sigma^{\#}_x}", from=2-2, to=1-1]
	\arrow["{\pi^{\#}_{(*,v)}}"', from=2-2, to=1-2]
\end{tikzcd}\]
Therefore
\begin{equation}
(\sigma^{\#}_x)^{-1}(v_x) =
\big( \pi^{\#}_{(*,v)} \big)^{-1} \big( (\tau^{\#}_x)^{-1}(v_x) \big) .
\end{equation}
Since $(\tau, \tau^{\#}) : (X, \tu{\cat{O}}, \tu{\cat{V}}) \to (\tilde{X}, \tilde{\cat{O}}, \tilde{\cat{V}})$ is a morphism in $\cat{C}_1$, the valuation $(\tau^{\#}_x)^{-1}(v_x)$ on the condensed ring $\tilde{\cat{O}}_{(*,v)}$ must be equal to the unique element $\tilde{v}_{(*,v)}$ of $\tilde{\cat{V}}_{(*,v)}$. Then (2) of \cref{as:explicit description of tilde X} shows that
\begin{equation}
\big( \pi^{\#}_{(*,v)} \big)^{-1} \big( (\tau^{\#}_x)^{-1}(v_x) \big) =
\big( \pi^{\#}_{(*,v)} \big)^{-1} \big( \tilde{v}_{(*,v)} \big) = v .
\end{equation}
On the other hand, \cref{prop:sigma sharp inverse v x is equal to x} shows that
\begin{equation}
(\sigma^{\#}_x)^{-1}(v_x) = x .
\end{equation}
Therefore
\begin{equation}
x = (\sigma^{\#}_x)^{-1}(v_x) =
\big( \pi^{\#}_{(*,v)} \big)^{-1} \big( (\tau^{\#}_x)^{-1}(v_x) \big) = v .
\end{equation}
Consequently, we have $\tau(x) = (*,v) = (*,x)$ and the diagram
\[\begin{tikzcd}
	{(\tu{\cat{O}}_x,v_x)} & {(\tilde{\cat{O}}_{(*,x)},\tilde{v}_{(*,x)})} \\
	& {(\tu{A},x)}
	\arrow["{\tau^{\#}_x}"', from=1-2, to=1-1]
	\arrow["{\sigma^{\#}_x}", from=2-2, to=1-1]
	\arrow["{\pi^{\#}_{(*,x)}}"', from=2-2, to=1-2]
\end{tikzcd}\]
is commutative in $\ub{VCRing}$.
\end{proof}

\begin{proof}[Proof of (2) of \cref{thm:comparison of coreflection and adic spectrum}]
By (1) of \cref{as:explicit description of tilde X}, the underlying set of $\tilde{X}$ is equal to the set of all pairs $(*,x)$ consisting of the point $* \in \{*\}$ and a valuation $x \in \cat{U}_* = \mathrm{Spa}(A,A^+) = X$. Furthermore, (1) of \cref{lem:identification of the map tau} shows that the map $\tau : X \to \tilde{X}$ is of the form $x \mapsto (*,x)$. It follows that the map $\tau : X \to \tilde{X}$ is bijective. In addition, the map $\tau : X \to \tilde{X}$ is continuous by definition. Therefore it remains to prove that the map $\tau : X \to \tilde{X}$ is open. Since the rational subsets of $X$ form a basis of the topology of $X$, it suffices to show that $\tau(U)$ is open in $\tilde{X}$ for every rational subset $U$ of $X$.

Let $U = R \left( \frac{T}{g} \right)$ be any rational subset of $X$, where $g \in A$ is an element of $A$ and $T \sub A$ is a finite subset of $A$ which generates an open ideal of $A$. We show that the set $\tau(U)$ is an open subset of $\tilde{X}$. Note that we have a homomorphism $\pi^{\#}_{\{*\}} : \tu{A} = F(\{*\}) \to \tilde{\cat{O}}(\tilde{X})$ of condensed rings. We claim that
\begin{equation}
\tau(U) = \bigcap_{f \in T} \set{\xi \in \tilde{X}}{
\left| \left( \left( \pi^{\#}_{\{*\}} \right)_*(f) \right)_{\xi} \right|_{\tilde{v}_{\xi}} \leq
\left| \left( \left( \pi^{\#}_{\{*\}} \right)_*(g) \right)_{\xi} \right|_{\tilde{v}_{\xi}} \neq 0}.
\end{equation}
Indeed, for every $\xi = (*,x) \in \tilde{X}$, the following diagram is commutaive.
\[\begin{tikzcd}
	{\tu{A}} & {F(\{*\})} & {\tilde{\cat{O}}(\tilde{X})} \\
	& {F_*} & {\tilde{\cat{O}}_{\xi}}
	\arrow[equals, from=1-1, to=1-2]
	\arrow[equals, from=1-1, to=2-2]
	\arrow["{\pi^{\#}_{\{*\}}}", from=1-2, to=1-3]
	\arrow["\can", from=1-2, to=2-2]
	\arrow["\can", from=1-3, to=2-3]
	\arrow["{\pi^{\#}_{\xi}}"', from=2-2, to=2-3]
\end{tikzcd}\]
Therefore
\begin{align}
& \bigcap_{f \in T} \set{\xi \in \tilde{X}}{
\left| \left( \left( \pi^{\#}_{\{*\}} \right)_*(f) \right)_{\xi} \right|_{\tilde{v}_{\xi}} \leq
\left| \left( \left( \pi^{\#}_{\{*\}} \right)_*(g) \right)_{\xi} \right|_{\tilde{v}_{\xi}} \neq 0} \\
= & \bigcap_{f \in T} \set{\xi \in \tilde{X}}{
\left| \left( \pi^{\#}_{\xi} \right)_*(f) \right|_{\tilde{v}_{\xi}} \leq
\left| \left( \pi^{\#}_{\xi} \right)_*(g) \right|_{\tilde{v}_{\xi}} \neq 0 }.
\end{align}
Moreover, for every $\xi = (*,x) \in \tilde{X}$, (2) of \cref{as:explicit description of tilde X} shows that
\begin{equation}
\left( \pi^{\#}_{\xi} \right)^{-1} \left( \tilde{v}_{\xi} \right) = x .
\end{equation}
Therefore
\begin{align}
& \bigcap_{f \in T} \set{\xi \in \tilde{X}}{
\left| \left( \pi^{\#}_{\xi} \right)_*(f) \right|_{\tilde{v}_{\xi}} \leq
\left| \left( \pi^{\#}_{\xi} \right)_*(g) \right|_{\tilde{v}_{\xi}} \neq 0 } \\
= & \bigcap_{f \in T} \set{(*,x) \in \tilde{X}}{|f|_x \leq |g|_x \neq 0 } .
\end{align}
On the other hand, we have
\begin{equation}
U = R \left( \frac{T}{g} \right) = \bigcap_{f \in T} \set{x \in X}{|f|_x \leq |g|_x \neq 0}
\end{equation}
by definition. Moreover, the map $\tau : X \to \tilde{X}$ is a bijection of the form $x \mapsto (*,x)$ by (1) of \cref{lem:identification of the map tau}. Therefore
\begin{equation}
\bigcap_{f \in T} \set{(*,x) \in \tilde{X}}{|f|_x \leq |g|_x \neq 0 } = \tau(U) . 
\end{equation}
Consequently, we have
\begin{equation}
\tau(U) = \bigcap_{f \in T} \set{\xi \in \tilde{X}}{
\left| \left( \left( \pi^{\#}_{\{*\}} \right)_*(f) \right)_{\xi} \right|_{\tilde{v}_{\xi}} \leq
\left| \left( \left( \pi^{\#}_{\{*\}} \right)_*(g) \right)_{\xi} \right|_{\tilde{v}_{\xi}} \neq 0}.
\end{equation}
This set is open in $\tilde{X}$ since $(\tilde{X}, \tilde{\cat{O}}, \tilde{\cat{V}})$ is an object of $\cat{C}_f$ and since $T$ is a finite set. It follows that $\tau(U)$ is open in $\tilde{X}$. This completes the proof.
\end{proof}

\begin{lem} \label{lem:properties of section over rational subsets}
Let $U = R \left( \frac{T}{g} \right)$ be a rational subset of $X$, where $g \in A$ is an element of $A$ and $T \sub A$ is a nonempty finite subset of $A$ which generates an open ideal of $A$. Let us write $T' := \set{(f,g)}{f \in T}$.
\begin{enumerate}
\item
$\tilde{\cat{O}}(\tau(U))$ is a $T'$-coalescent condensed $\tu{A}$-algebra via $\tu{A} = F(\{*\}) \xto{\pi^{\#}_{\{*\}}} \tilde{\cat{O}}(\tilde{X}) \xto{\res} \tilde{\cat{O}}(\tau(U))$. The element $((\pi^{\#}_{\{*\}})_*(g) )|_{\tau(U)} \in \tilde{\cat{O}}(\tau(U))(*)$ is invertible in $\tilde{\cat{O}}(\tau(U))(*)$. 

\item
$\tu{\cat{O}}(U)$ is a $T'$-coalescent condensed $\tu{A}$-algebra via $\tu{A} = F(\{*\}) \xto{\sigma^{\#}_{\{*\}}} \tu{\cat{O}}(X) \xto{\res} \tu{\cat{O}}(U)$. The element $((\sigma^{\#}_{\{*\}})_*(g)) |_U \in \tu{\cat{O}}(U)(*)$ is invertible in $\tu{\cat{O}}(U)(*)$.
\end{enumerate}
\end{lem}

\begin{proof}
Let $x \in U$ and $f \in T$. Then we have $|f|_x \leq |g|_x \neq 0$ since $x \in U = R \left( \frac{T}{g} \right)$. By \cref{lem:identification of the map tau}, we have the homomoprhisms of valued condensed rings
$\pi^{\#}_{\tau(x)} : (\tu{A},x) \to (\tilde{\cat{O}}_{\tau(x)},\tilde{v}_{\tau(x)})$ and $\sigma^{\#}_x : (\tu{A},x) \to (\tu{\cat{O}}_x,v_x)$. Therefore
\begin{align}
\left| \left( \pi^{\#}_{\tau(x)} \right)_* (f) \right|_{\tilde{v}_{\tau(x)}} & \leq
\left| \left( \pi^{\#}_{\tau(x)} \right)_* (g) \right|_{\tilde{v}_{\tau(x)}} \neq 0 \: ; \\
\left| \left( \sigma^{\#}_x \right)_* (f) \right|_{v_x} & \leq 
\left| \left( \sigma^{\#}_x \right)_* (g) \right|_{v_x} \neq 0 .
\end{align}
On the other hand, the following diagrams are commutative.
\[\begin{tikzcd}
	{\tu{A}} & {F(\{*\})} & {\tilde{\cat{O}}(\tilde{X})} & {\tu{A}} & {F(\{*\})} & {\tu{\cat{O}}(X)} \\
	&& {\tilde{\cat{O}}(\tau(U))} &&& {\tu{\cat{O}}(U)} \\
	{\tu{A}} & {F_*} & {\tilde{\cat{O}}_{\tau(x)}} & {\tu{A}} & {F_*} & {\tu{\cat{O}}_x}
	\arrow[equals, from=1-1, to=1-2]
	\arrow[equals, from=1-1, to=3-1]
	\arrow["{\pi^{\#}_{\{*\}}}", from=1-2, to=1-3]
	\arrow["\can"', from=1-2, to=3-2]
	\arrow["\res", from=1-3, to=2-3]
	\arrow[equals, from=1-4, to=1-5]
	\arrow[equals, from=1-4, to=3-4]
	\arrow["{\sigma^{\#}_{\{*\}}}", from=1-5, to=1-6]
	\arrow["\can"', from=1-5, to=3-5]
	\arrow["\res", from=1-6, to=2-6]
	\arrow["\can", from=2-3, to=3-3]
	\arrow["\can", from=2-6, to=3-6]
	\arrow[equals, from=3-1, to=3-2]
	\arrow["{\pi^{\#}_{\tau(x)}}"', from=3-2, to=3-3]
	\arrow[equals, from=3-4, to=3-5]
	\arrow["{\sigma^{\#}_x}"', from=3-5, to=3-6]
\end{tikzcd}\]
Therefore
\begin{align}
\left| \left( \left( \pi^{\#}_{\{*\}} \right)_* (f) \right)_{\tau(x)} \right|_{\tilde{v}_{\tau(x)}} & \leq
\left| \left( \left( \pi^{\#}_{\{*\}} \right)_* (g) \right)_{\tau(x)} \right|_{\tilde{v}_{\tau(x)}} \neq 0 \: ; \\
\left| \left( \left( \sigma^{\#}_{\{*\}} \right)_* (f) \right)_x \right|_{v_x} & \leq 
\left| \left( \left( \sigma^{\#}_{\{*\}} \right)_* (g) \right)_x \right|_{v_x} \neq 0 .
\end{align}
Since $T \neq \ku$, we may conclude that
\begin{align}
& \left| \left( \left( \pi^{\#}_{\{*\}} \right)_* (g) \right)_{\tau(x)} \right|_{\tilde{v}_{\tau(x)}} \neq 0 \: ; \\
& \left| \left( \left( \sigma^{\#}_{\{*\}} \right)_* (g) \right)_x \right|_{v_x} \neq 0 .
\end{align}
Then since $(X, \tu{\cat{O}}, \tu{\cat{V}})$ and $(\tilde{X}, \tilde{\cat{O}}, \tilde{\cat{V}})$ are objects of $\cat{C}_l \cap \cat{C}_c$, the following hold for every $x \in U$.
\begin{itemize}
\item
$\tilde{\cat{O}}_{\tau(x)}$ is a 
$\left( \left( \pi^{\#}_{\{*\}} \right)_* (f) \right)_{\tau(x)} \Big/
 \left( \left( \pi^{\#}_{\{*\}} \right)_* (g) \right)_{\tau(x)}$-coalescent condensed 
$\tilde{\cat{O}}_{\tau(x)}$-algebra for every $f \in T$.

\item
$\tu{\cat{O}}_x$ is a 
$\left( \left( \sigma^{\#}_{\{*\}} \right)_* (f) \right)_x \Big/
 \left( \left( \sigma^{\#}_{\{*\}} \right)_* (g) \right)_x$-coalescent condensed 
$\tu{\cat{O}}_x$-algebra for every $f \in T$.

\item
The element $\left( \left( \pi^{\#}_{\{*\}} \right)_* (g) \right)_{\tau(x)} 
\in \tilde{\cat{O}}_{\tau(x)}(*)$ is invertible in $\tilde{\cat{O}}_{\tau(x)}(*)$.

\item
The element $\left( \left( \sigma^{\#}_{\{*\}} \right)_* (g) \right)_x 
\in \tu{\cat{O}}_x(*)$ is invertible in $\tu{\cat{O}}_x(*)$.
\end{itemize}
Then \cref{prop:testing invertibility by stalk} shows that the following hold.
\begin{itemize}
\item
The element $\left. \left( \pi^{\#}_{\{*\}} \right)_* (g) \right|_{\tau(U)}
\in \tilde{\cat{O}}(\tau(U))(*)$ is invertible in $\tilde{\cat{O}}(\tau(U))(*)$. 

\item
The element $\left. \left( \sigma^{\#}_{\{*\}} \right)_* (g) \right|_U
\in \tu{\cat{O}}(U)(*)$ is invertible in $\tu{\cat{O}}(U)(*)$.
\end{itemize}
Moreover, the following hold for every $x \in U$ and every $f \in T$.
\begin{itemize}
\item
$\tilde{\cat{O}}_{\tau(x)}$ is a 
$\left. \left( \pi^{\#}_{\{*\}} \right)_* (f) \right|_{\tau(U)} \Big/
 \left. \left( \pi^{\#}_{\{*\}} \right)_* (g) \right|_{\tau(U)}$-coalescent condensed 
$\tilde{\cat{O}}(\tau(U))$-algebra via the canonical homomorphism 
$\tilde{\cat{O}}(\tau(U)) \to \tilde{\cat{O}}_{\tau(x)}$.

\item
$\tu{\cat{O}}_x$ is a 
$\left. \left( \sigma^{\#}_{\{*\}} \right)_* (f) \right|_U \Big/
 \left. \left( \sigma^{\#}_{\{*\}} \right)_* (g) \right|_U$-coalescent condensed 
$\tu{\cat{O}}(U)$-algebra via the canonical homomorphism $\tu{\cat{O}}(U) \to \tu{\cat{O}}_x$.
\end{itemize}
Then \cref{prop:testing coalescence by stalks} shows that the following hold for every $f \in T$.
\begin{itemize}
\item
$\tilde{\cat{O}}(\tau(U))$ is a 
$\left. \left( \pi^{\#}_{\{*\}} \right)_* (f) \right|_{\tau(U)} \Big/
 \left. \left( \pi^{\#}_{\{*\}} \right)_* (g) \right|_{\tau(U)}$-coalescent condensed 
$\tilde{\cat{O}}(\tau(U))$-algebra. 

\item
$\tu{\cat{O}}(U)$ is a 
$\left. \left( \sigma^{\#}_{\{*\}} \right)_* (f) \right|_U \Big/
 \left. \left( \sigma^{\#}_{\{*\}} \right)_* (g) \right|_U$-coalescent condensed 
$\tu{\cat{O}}(U)$-algebra.
\end{itemize}
Therefore the following hold for every $f \in T$.
\begin{itemize}
\item
$\tilde{\cat{O}}(\tau(U))$ is a 
$f/g$-coalescent condensed $\tu{A}$-algebra via 
$\tu{A} = F(\{*\}) \xto{\pi^{\#}_{\{*\}}} \tilde{\cat{O}}(\tilde{X}) \xto{\res} \tilde{\cat{O}}(\tau(U))$.  

\item
$\tu{\cat{O}}(U)$ is a 
$f/g$-coalescent condensed $\tu{A}$-algebra via 
$\tu{A} = F(\{*\}) \xto{\sigma^{\#}_{\{*\}}} \tu{\cat{O}}(X) \xto{\res} \tu{\cat{O}}(U)$.
\end{itemize}
This completes the proof.
\end{proof}

\begin{proof}[Proof of (3) of \cref{thm:comparison of coreflection and adic spectrum}]
Let $x \in X$. Let $M$ be any condensed set. Let $\alpha, \beta : \tu{\cat{O}}_x \to M$ be any maps of condensed sets such that $\alpha \of \tau^{\#}_x = \beta \of \tau^{\#}_x$. We prove that $\alpha = \beta$.

Let $U$ be any rational subset of $X$ such that $x \in U$. Let us write $U = R \left( \frac{T}{g} \right)$, where $g \in A$ is an element of $A$ and $T \sub A$ is a finite subset of $A$ which generates an open ideal of $A$. We may assume $T \neq \ku$. Let us write $T' := \set{(f,g)}{f \in T}$. By \cref{prop:comparison of localizations of Huber rings}, there exists a unique homomorphism $\gamma : (\tu{A}_g)_{\approx T'} \to \tu{A \left \langle \frac{T}{g} \right \rangle} = \tu{\cat{O}}^{\mathrm{pre}}(U)$ of condensed rings such that the diagram
\[\begin{tikzcd}
	{\tu{\cat{O}}^{\mathrm{pre}}(X)} && {\tu{A}} \\
	{\tu{\cat{O}}^{\mathrm{pre}}(U)} & {\tu{A \left \langle \frac{T}{g} \right \rangle}} & {(\tu{A}_g)_{\approx T'}}
	\arrow["\res"', from=1-1, to=2-1]
	\arrow[equals, from=1-3, to=1-1]
	\arrow["\can", from=1-3, to=2-3]
	\arrow[equals, from=2-2, to=2-1]
	\arrow["\gamma", from=2-3, to=2-2]
\end{tikzcd}\]
is commutative. Moreover, $\gamma : (\tu{A}_g)_{\approx T'} \to \tu{\cat{O}}^{\mathrm{pre}}(U)$ is an epimorphism when considered as a morphism in $\ub{CSet}$. On the other hand, by (1) of \cref{lem:properties of section over rational subsets} and \cref{prop:localization and coalescence commute}, there exists a unique homomorphism $\delta : (\tu{A}_g)_{\approx T'} \to \tilde{\cat{O}}(\tau(U))$ of condensed rings such that the diagram
\[\begin{tikzcd}
	{\tilde{\cat{O}}(\tilde{X})} & {\tu{A} = F(\{*\})} \\
	{\tilde{\cat{O}}(\tau(U))} & {(\tu{A}_g)_{\approx T'} }
	\arrow["\res"', from=1-1, to=2-1]
	\arrow["{\pi^{\#}_{\{*\}}}"', from=1-2, to=1-1]
	\arrow["\can", from=1-2, to=2-2]
	\arrow["\delta", from=2-2, to=2-1]
\end{tikzcd}\]
is commutative. Then the following diagram is commutative.
\[\begin{tikzcd}
	&&& {\tilde{\cat{O}}(\tau(U))} & {(\tu{A}_g)_{\approx T'} } \\
	&&& {\tilde{\cat{O}}(\tilde{X})} \\
	{\tu{\cat{O}}(U)} && {\tu{\cat{O}}(X)} && {\tu{A}} \\
	&&& {\tu{\cat{O}}^{\mathrm{pre}}(X)} \\
	&&& {\tu{\cat{O}}^{\mathrm{pre}}(U)} & {(\tu{A}_g)_{\approx T'} }
	\arrow["{\tau^{\#}_{\tau(U)}}"', from=1-4, to=3-1]
	\arrow["\delta"', from=1-5, to=1-4]
	\arrow["\res"', from=2-4, to=1-4]
	\arrow["{\tau^{\#}_{\tilde{X}}}"', from=2-4, to=3-3]
	\arrow["\res"', from=3-3, to=3-1]
	\arrow["\can"', from=3-5, to=1-5]
	\arrow["{\pi^{\#}_{\{*\}}}"', from=3-5, to=2-4]
	\arrow["{\sigma^{\#}_{\{*\}}}"', from=3-5, to=3-3]
	\arrow[equals, from=3-5, to=4-4]
	\arrow["\can", from=3-5, to=5-5]
	\arrow["{\eta_X}", from=4-4, to=3-3]
	\arrow["\res", from=4-4, to=5-4]
	\arrow["{\eta_U}", from=5-4, to=3-1]
	\arrow["\gamma", from=5-5, to=5-4]
\end{tikzcd}\]
By (2) of \cref{lem:properties of section over rational subsets} and \cref{prop:localization and coalescence commute}, we conclude that the following diagarm is commutative.
\[\begin{tikzcd}
	{\tu{\cat{O}}(U)} & {\tilde{\cat{O}}(\tau(U))} \\
	{\tu{\cat{O}}^{\mathrm{pre}}(U)} & {(\tu{A}_g)_{\approx T'} }
	\arrow["{\tau^{\#}_{\tau(U)}}"', from=1-2, to=1-1]
	\arrow["{\eta_U}", from=2-1, to=1-1]
	\arrow["\delta"', from=2-2, to=1-2]
	\arrow["\gamma", from=2-2, to=2-1]
\end{tikzcd}\]
Then the following diagram is commutative.
\[\begin{tikzcd}
	M & {\tu{\cat{O}}_x} && {\tilde{\cat{O}}_{\tau(x)}} \\
	&& {\tu{\cat{O}}(U)} & {\tilde{\cat{O}}(\tau(U))} \\
	& {(\tu{\cat{O}}^{\mathrm{pre}})_x} & {\tu{\cat{O}}^{\mathrm{pre}}(U)} & {(\tu{A}_g)_{\approx T'} }
	\arrow["{\alpha, \beta}"', from=1-2, to=1-1]
	\arrow["{\tau^{\#}_x}"', from=1-4, to=1-2]
	\arrow["\can"', from=2-3, to=1-2]
	\arrow["\can"', from=2-4, to=1-4]
	\arrow["{\tau^{\#}_{\tau(U)}}"', from=2-4, to=2-3]
	\arrow["{\eta_x}", from=3-2, to=1-2]
	\arrow["{\eta_U}", from=3-3, to=2-3]
	\arrow["{\rho_U}", from=3-3, to=3-2]
	\arrow["\delta"', from=3-4, to=2-4]
	\arrow["\gamma", from=3-4, to=3-3]
\end{tikzcd}\]
Here we write $\rho_U : \tu{\cat{O}}^{\mathrm{pre}}(U) \to (\tu{\cat{O}}^{\mathrm{pre}})_x$ for the canonical homomorphism. Since $\alpha \of \tau^{\#}_x = \beta \of \tau^{\#}_x$, we conclude that $\alpha \of \eta_x \of \rho_U \of \gamma = \beta \of \eta_x \of \rho_U \of \gamma$. Since $\gamma : (\tu{A}_g)_{\approx T'} \to \tu{\cat{O}}^{\mathrm{pre}}(U)$ is an epimorphism when considered as a morphism in $\ub{CSet}$, it follows that $\alpha \of \eta_x \of \rho_U = \beta \of \eta_x \of \rho_U$.

Thus we have proved that for every rational subset $U$ of $X$ such that $x \in U$, we have $\alpha \of \eta_x \of \rho_U = \beta \of \eta_x \of \rho_U : (\tu{\cat{O}}^{\mathrm{pre}})_x \to M$. On the other hand, since the forgetful functor $\ub{CRing} \to \ub{CSet}$ preserves filtered colimits, we have
\begin{equation}
(\tu{\cat{O}}^{\mathrm{pre}})_x
= \underset{x \in U \sub X \text{ rational}}{\colim} \: \tu{\cat{O}}^{\mathrm{pre}}(U)
\end{equation}
in the category $\ub{CSet}$. Consequently, we have $\alpha \of \eta_x = \beta \of \eta_x$. Furthermore, by \cref{cor:existence of sheafification functor of condensed presheaves}, the homomorphism $\eta_x : (\tu{\cat{O}}^{\mathrm{pre}})_x \to \tu{\cat{O}}_x$ is an isomorphism of condensed rings, and in particular an isomorphism of condensed sets. It follows that $\alpha = \beta$. This completes the proof.
\end{proof}

\begin{lem} \label{lem:identification of stalk of underbar O}
Suppose that one of the following holds.
\begin{enumerate}
\item $A$ is a strongly Noetherian Tate ring.
\item $A$ has a Noetherian ring of definition.
\end{enumerate}
Then for every $x \in X$, the pair
\begin{equation}
\left( (\tu{\cat{O}}_x,v_x) \; , \; (\tu{A} , x) \xto{\sigma^{\#}_x} (\tu{\cat{O}}_x,v_x) \right)
\end{equation}
is the reflection of the valued condensed ring $(\tu{A} , x)$ along the inclusion functor $\ub{VCRing}_l \cap \ub{VCRing}_c \mon \ub{VCRing}$. 
\end{lem}

\begin{proof}
Let $x \in X$. Then $(\tu{\cat{O}}_x,v_x)$ is an object of $\ub{VCRing}_l \cap \ub{VCRing}_c$ since $(X, \tu{\cat{O}}, \tu{\cat{V}})$ is an object of $\cat{C}_l \cap \cat{C}_c$ by \cref{prop:adic spectra is an object of C l cap C c}. Moreover, \cref{prop:sigma sharp inverse v x is equal to x} shows that $\sigma^{\#}_x$ is a homomorphism $(\tu{A} , x) \to (\tu{\cat{O}}_x,v_x)$ of valued condensed rings.

On the other hand, suppose that $(R,w)$ is an object of $\ub{VCRing}_l \cap \ub{VCRing}_c$ and that $\phi : (\tu{A} , x) \to (R,w)$ is a homomorphism of valued condensed rings. We prove that there exists a unique homomorphism $\psi : (\tu{\cat{O}}_x,v_x) \to (R,w)$ of valued condensed rings such that the diagram
\[\begin{tikzcd}
	{(\tu{A} , x)} & {(R,w)} \\
	{(\tu{\cat{O}}_x,v_x)}
	\arrow["\phi", from=1-1, to=1-2]
	\arrow["{\sigma^{\#}_x}"', from=1-1, to=2-1]
	\arrow["\psi"', from=2-1, to=1-2]
\end{tikzcd}\]
is commutative.

We first claim that for every rational subset $U$ of $X$ such that $x \in U$, there exists a unique homomorphism $\phi_U : \tu{\cat{O}}^{\mathrm{pre}}(U) \to R$ of condensed rings such that the diagram
\[\begin{tikzcd}
	{\tu{\cat{O}}^{\mathrm{pre}}(X)} & {\tu{A}} & {R} \\
	{\tu{\cat{O}}^{\mathrm{pre}}(U)}
	\arrow[equals, from=1-1, to=1-2]
	\arrow["{\tu{\cat{O}}^{\mathrm{pre}}(U \sub X)}"', from=1-1, to=2-1]
	\arrow["\phi", from=1-2, to=1-3]
	\arrow["{\phi_U}"', from=2-1, to=1-3]
\end{tikzcd}\]
is commutative. Let us write $U = R \left( \frac{T}{g} \right)$ where $g \in A$ is an element of $A$ and $T \sub A$ is a finite subset of $A$ which generates an open ideal of $A$. We may assume that $T \neq \ku$. Choose some element $f \in T$. Since $x \in U = R \left( \frac{T}{g} \right)$, we have $|f|_x \leq |g|_x \neq 0$. In paricular $|g|_x \neq 0$. Since $\phi : (\tu{A} , x) \to (R,w)$ is a homomorphism of valued condensed rings, we have $|\phi_*(g)|_w \neq 0$. Since $(R,w)$ is an object of $\ub{VCRing}_l$, it follows that $\phi_*(g) \in R(*)$ is invertible in $R(*)$. Furthermore, for every $f \in T$, we have $|f|_x \leq |g|_x \neq 0$ since $x \in U = R \left( \frac{T}{g} \right)$. Since $\phi : (\tu{A} , x) \to (R,w)$ is a homomorphism of valued condensed rings, we have $|\phi_*(f)|_w \leq |\phi_*(g)|_w \neq 0$. Since $(R,w)$ is an object of $\ub{VCRing}_c$, it follows that $R$ is a $\phi_*(f) / \phi_*(g)$-coalescent condensed $R$-algebra. Therefore $R$ is an $f/g$-coalescent condensed $\tu{A}$-algebra via $\phi : \tu{A} \to R$. Thus we have proved that $\phi_*(g) \in R(*)$ is invertible in $R(*)$ and that $R$ is an $T'$-coalescent condensed $\tu{A}$-algebra via $\phi : \tu{A} \to R$, where $T' = \set{(f,g)}{f \in T}$. Since we assumed that $A$ is a strongly Noetherian Tate ring or has a Noetherian ring of definition, \cref{prop:comparison of localizations of Huber rings} and \cref{prop:localization and coalescence commute} show that there exists a unique homomorphism $\phi_U : \tu{A \left \langle \frac{T}{g} \right \rangle} \to R$ of condensed rings such that the diagram
\[\begin{tikzcd}
	{\tu{A}} & R \\
	{\tu{A \left \langle \frac{T}{g} \right \rangle}}
	\arrow["\phi", from=1-1, to=1-2]
	\arrow["{\tu{l}}"', from=1-1, to=2-1]
	\arrow["{\phi_U}"', from=2-1, to=1-2]
\end{tikzcd}\]
is commutative, where $l : A \to A \left \langle \frac{T}{g} \right \rangle$ denotes the canonical homomorphism. Since $\tu{\cat{O}}^{\mathrm{pre}}(U) = \tu{\cat{O}(U)} = \tu{A \left \langle \frac{T}{g} \right \rangle}$ and the restriction $\tu{\cat{O}}^{\mathrm{pre}}(U \sub X) : \tu{\cat{O}}^{\mathrm{pre}}(X) \to \tu{\cat{O}}^{\mathrm{pre}}(U)$ is equal to the homomorphism $\tu{l} : \tu{A} \to \tu{A \left \langle \frac{T}{g} \right \rangle}$, we conclude that the homomorphism $\phi_U$ is a unique homomorphism $\tu{\cat{O}}^{\mathrm{pre}}(U) \to R$ of condensed rings such that the diagram
\[\begin{tikzcd}
	{\tu{\cat{O}}^{\mathrm{pre}}(X)} & {\tu{A}} & {R} \\
	{\tu{\cat{O}}^{\mathrm{pre}}(U)}
	\arrow[equals, from=1-1, to=1-2]
	\arrow["{\tu{\cat{O}}^{\mathrm{pre}}(U \sub X)}"', from=1-1, to=2-1]
	\arrow["\phi", from=1-2, to=1-3]
	\arrow["{\phi_U}"', from=2-1, to=1-3]
\end{tikzcd}\]
is commutative.

If $U,V$ are rational subsets of $X$ such that $x \in V \sub U$, then the following diagram is commutative.
\[\begin{tikzcd}
	{\tu{\cat{O}}^{\mathrm{pre}}(X)} & {\tu{A}} & {R} \\
	\\
	{\tu{\cat{O}}^{\mathrm{pre}}(U)} && {\tu{\cat{O}}^{\mathrm{pre}}(V)}
	\arrow[equals, from=1-1, to=1-2]
	\arrow["{\tu{\cat{O}}^{\mathrm{pre}}(U \sub X)}"', from=1-1, to=3-1]
	\arrow["{\tu{\cat{O}}^{\mathrm{pre}}(V \sub X)}", from=1-1, to=3-3]
	\arrow["\phi", from=1-2, to=1-3]
	\arrow["{\tu{\cat{O}}^{\mathrm{pre}}(V \sub U)}"', from=3-1, to=3-3]
	\arrow["{\phi_V}"', from=3-3, to=1-3]
\end{tikzcd}\]
Then the uniqueness of $\phi_U$ shows that the composition $\tu{\cat{O}}^{\mathrm{pre}}(U) \xto{\tu{\cat{O}}^{\mathrm{pre}}(V \sub U)} \tu{\cat{O}}^{\mathrm{pre}}(V) \xto{\phi_V} R$ is equal to $\phi_U$. Since we have
\begin{equation}
(\tu{\cat{O}}^{\mathrm{pre}})_x =
\underset{x \in U \sub X \text{ rational}}{\colim} \: \: \tu{\cat{O}}^{\mathrm{pre}}(U) ,
\end{equation}
in the category $\ub{CRing}$, we conclude that there exists a unique homomorphism $\tilde{\phi} : (\tu{\cat{O}}^{\mathrm{pre}})_x \to R$ of condensed rings such that the diagram
\[\begin{tikzcd}
	{\tu{\cat{O}}^{\mathrm{pre}}(U)} & R \\
	{(\tu{\cat{O}}^{\mathrm{pre}})_x}
	\arrow["{\phi_U}", from=1-1, to=1-2]
	\arrow["\can"', from=1-1, to=2-1]
	\arrow["{\tilde{\phi}}"', from=2-1, to=1-2]
\end{tikzcd}\]
is commutative for every rational subset $U$ of $X$ such that $x \in U$.

By \cref{cor:existence of sheafification functor of condensed presheaves}, the homomorphism $\eta_x : (\tu{\cat{O}}^{\mathrm{pre}})_x \to \tu{\cat{O}}_x$ is an isomorphism of condensed rings. Moreover, the following diagram is commutative.
\[\begin{tikzcd}
	{\tu{A}} &&& {\tu{A}} \\
	& {F(\{*\})} & {\tu{\cat{O}}(X)} & {\tu{\cat{O}}^{\mathrm{pre}}(X)} & R \\
	{\tu{A}} & {F_*} & {\tu{\cat{O}}_x} & {(\tu{\cat{O}}^{\mathrm{pre}})_x}
	\arrow[equals, from=1-1, to=2-2]
	\arrow[equals, from=1-1, to=3-1]
	\arrow[equals, from=1-4, to=1-1]
	\arrow[equals, from=1-4, to=2-4]
	\arrow["\phi", curve={height=-12pt}, from=1-4, to=2-5]
	\arrow["{\sigma^{\#}_{\{*\}}}", from=2-2, to=2-3]
	\arrow["\can"', from=2-2, to=3-2]
	\arrow["\can"', from=2-3, to=3-3]
	\arrow["{\eta_X}"', from=2-4, to=2-3]
	\arrow["{\phi_X}"', from=2-4, to=2-5]
	\arrow["\can", from=2-4, to=3-4]
	\arrow[equals, from=3-1, to=3-2]
	\arrow["{\sigma^{\#}_x}"', from=3-2, to=3-3]
	\arrow["{\eta_x^{-1}}"', from=3-3, to=3-4]
	\arrow["{\tilde{\phi}}"', curve={height=12pt}, from=3-4, to=2-5]
\end{tikzcd}\]
Therefore if we define $\psi := \tilde{\phi} \of \eta_x^{-1} : \tu{\cat{O}}_x \to R$, the diagram
\[\begin{tikzcd}
	{\tu{A}} & R \\
	{\tu{\cat{O}}_x}
	\arrow["\phi", from=1-1, to=1-2]
	\arrow["{\sigma^{\#}_x}"', from=1-1, to=2-1]
	\arrow["\psi"', from=2-1, to=1-2]
\end{tikzcd}\]
is commutative. Then we have
\begin{equation}
(\sigma^{\#}_x)^{-1} \big( \psi^{-1}(w) \big) = \phi^{-1}(w) .
\end{equation}
Since $\phi : (\tu{A} , x) \to (R,w)$ is a homomorphism of valued condensed rings, we have $\phi^{-1}(w) = x$. Therefore
\begin{equation}
(\sigma^{\#}_x)^{-1} \big( \psi^{-1}(w) \big) = \phi^{-1}(w) = x.
\end{equation}
On the other hand, \cref{prop:sigma sharp inverse v x is equal to x} shows that the valuation $v_x$ is a unique continuous valuation on $\tu{\cat{O}}_x$ such that $(\sigma^{\#}_x)^{-1}(v_x) = x$. Consequently we have
\begin{equation}
\psi^{-1}(w) = v_x .
\end{equation}
It follows that the diagram
\[\begin{tikzcd}
	{(\tu{A} , x)} & {(R,w)} \\
	{(\tu{\cat{O}}_x,v_x)}
	\arrow["\phi", from=1-1, to=1-2]
	\arrow["{\sigma^{\#}_x}"', from=1-1, to=2-1]
	\arrow["\psi"', from=2-1, to=1-2]
\end{tikzcd}\]
is commutative in $\ub{VCRing}$.

On the other hand, suppose that $\psi' : (\tu{\cat{O}}_x,v_x) \to (R,w)$ is another homomorphism of valued condensed rings such that the diagram
\[\begin{tikzcd}
	{(\tu{A} , x)} & {(R,w)} \\
	{(\tu{\cat{O}}_x,v_x)}
	\arrow["\phi", from=1-1, to=1-2]
	\arrow["{\sigma^{\#}_x}"', from=1-1, to=2-1]
	\arrow["{\psi'}"', from=2-1, to=1-2]
\end{tikzcd}\]
is commutative. We prove that $\psi' = \psi$. For each rational subset $U$ of $X$ such that $x \in U$, the following diagram is commutative.
\[\begin{tikzcd}
	& {\tu{\cat{O}}^{\mathrm{pre}}(X)} &&& {\tu{A}} \\
	&& {\tu{\cat{O}}(X)} & {F(\{*\})} \\
	{\tu{\cat{O}}^{\mathrm{pre}}(U)} & {(\tu{\cat{O}}^{\mathrm{pre}})_x} & {\tu{\cat{O}}_x} & {F_*} & {\tu{A}} & R
	\arrow["{\eta_X}"', from=1-2, to=2-3]
	\arrow["{\tu{\cat{O}}^{\mathrm{pre}}(U \sub X)}"', from=1-2, to=3-1]
	\arrow["\can"', from=1-2, to=3-2]
	\arrow[equals, from=1-5, to=1-2]
	\arrow[equals, from=1-5, to=2-4]
	\arrow[equals, from=1-5, to=3-5]
	\arrow["\phi", curve={height=-6pt}, from=1-5, to=3-6]
	\arrow["\can"', from=2-3, to=3-3]
	\arrow["{\sigma^{\#}_{\{*\}}}"', from=2-4, to=2-3]
	\arrow["\can"', from=2-4, to=3-4]
	\arrow["\can"', from=3-1, to=3-2]
	\arrow["{\eta_x}"', from=3-2, to=3-3]
	\arrow["{\psi'}"', curve={height=24pt}, from=3-3, to=3-6]
	\arrow["{\sigma^{\#}_x}"', from=3-4, to=3-3]
	\arrow[equals, from=3-5, to=3-4]
\end{tikzcd}\]
Then the uniqueness of $\phi_U$ shows that the following diagram is commutative.
\[\begin{tikzcd}
	{\tu{\cat{O}}^{\mathrm{pre}}(U)} & R \\
	{(\tu{\cat{O}}^{\mathrm{pre}})_x} & {\tu{\cat{O}}_x}
	\arrow["{\phi_U}", from=1-1, to=1-2]
	\arrow["\can"', from=1-1, to=2-1]
	\arrow["{\eta_x}"', from=2-1, to=2-2]
	\arrow["{\psi'}"', from=2-2, to=1-2]
\end{tikzcd}\]
This holds for every rational subset $U$ of $X$ such that $x \in U$. Then the uniqueness of $\tilde{\phi}$ shows that $\psi' \of \eta_x = \tilde{\phi}$. Consequently, we have $\psi = \tilde{\phi} \of \eta_x^{-1} = \psi' \of \eta_x \of \eta_x^{-1} = \psi'$. This completes the proof.
\end{proof}

\begin{proof}[Proof of (4) of \cref{thm:comparison of coreflection and adic spectrum}]
Let $x \in X$. \cref{lem:identification of the map tau} shows that $\tau(x) = (*,x)$ and that the diagram
\[\begin{tikzcd}
	{(\tu{\cat{O}}_x,v_x)} & {(\tilde{\cat{O}}_{(*,x)},\tilde{v}_{(*,x)})} \\
	& {(\tu{A},x)}
	\arrow["{\tau^{\#}_x}"', from=1-2, to=1-1]
	\arrow["{\sigma^{\#}_x}", from=2-2, to=1-1]
	\arrow["{\pi^{\#}_{(*,x)}}"', from=2-2, to=1-2]
\end{tikzcd}\]
is commutative in $\ub{VCRing}$. Moreover, (2) of \cref{as:explicit description of tilde X} shows that the pair
\begin{equation}
\left( (\tilde{\cat{O}}_{(*,x)} , \tilde{v}_{(*,x)}) \; , \; (\tu{A} , x) \xto{\pi^{\#}_{(*,x)}} (\tilde{\cat{O}}_{(*,x)} , \tilde{v}_{(*,x)}) \right)
\end{equation}
is a reflection of the valued condensed ring $(\tu{A} , x)$ along the inclusion functor $\ub{VCRing}_l \cap \ub{VCRing}_c \mon \ub{VCRing}$. On the other hand, \cref{lem:identification of stalk of underbar O} shows that the pair
\begin{equation}
\left( (\tu{\cat{O}}_x,v_x) \; , \; (\tu{A} , x) \xto{\sigma^{\#}_x} (\tu{\cat{O}}_x,v_x) \right)
\end{equation}
is also a reflection of the valued condensed ring $(\tu{A} , x)$ along the inclusion functor $\ub{VCRing}_l \cap \ub{VCRing}_c \mon \ub{VCRing}$. Then the uniqueness of reflections along functors implies that the homomorphism $\tau^{\#}_x : (\tilde{\cat{O}}_{(*,x)},\tilde{v}_{(*,x)}) \to (\tu{\cat{O}}_x,v_x)$ must be an isomorphism of valued condensed rings. In particular, the homomorphism $\tau^{\#}_x : \tilde{\cat{O}}_{(*,x)} \to \tu{\cat{O}}_x$ is an isomorphism of condensed rings.

Thus we have proved that the homomorphism $\tau^{\#}_x : \tilde{\cat{O}}_{\tau(x)} \to \tu{\cat{O}}_x$ is an isomorphism of condensed rings for every $x \in X$. On the other hand, recall from (2) of this \cref{thm:comparison of coreflection and adic spectrum} that the map $\tau : X \to \tilde{X}$ is a homeomorphism of topological spaces. Therefore, for every $x \in X$, we have $\tu{\cat{O}}_x = (\tau_* \tu{\cat{O}})_{\tau(x)}$ and the homomorphism $(\tau^{\#})_{\tau(x)} : \tilde{\cat{O}}_{\tau(x)} \to (\tau_* \tu{\cat{O}})_{\tau(x)}$ is equal to the homomorphism $\tau^{\#}_x : \tilde{\cat{O}}_{\tau(x)} \to \tu{\cat{O}}_x$. Consequently, for every $\xi \in \tilde{X}$, the homomorphism $(\tau^{\#})_{\xi} : \tilde{\cat{O}}_{\xi} \to (\tau_* \tu{\cat{O}})_{\xi}$ is an isomorphism of condensed rings. Then \cref{prop:testing sheaf isomorphisms by stalks} shows that the morphism $\tau^{\#} : \tilde{\cat{O}} \to \tau_* \tu{\cat{O}}$ is an isomorphism of sheaves of condensed rings on $\tilde{X}$. In addition, the map $\tau : X \to \tilde{X}$ is a homeomorphism by (2) of this \cref{thm:comparison of coreflection and adic spectrum}. It follows that the morphism $(\tau, \tau^{\#}) : (X, \tu{\cat{O}}, \tu{\cat{V}}) \to (\tilde{X}, \tilde{\cat{O}}, \tilde{\cat{V}})$ is an isomorphism in $\cat{C}$.
\end{proof}

\printbibliography

\end{document}